\documentclass[11pt]{amsart}
\usepackage{amsfonts,amsrefs,latexsym,amsmath,amsthm, amssymb, mathrsfs, verbatim,cancel,slashed}
\usepackage{url,color,fullpage}
\usepackage{pifont}
\usepackage{upgreek}
\usepackage{fancyhdr,hyperref}

\usepackage[percent]{overpic}
\usepackage{caption} 
\usepackage{subcaption}
\usepackage{pict2e}
\usepackage[hypcap=false]{caption}

\newcommand{\R}{\mathbb{R}}

\newcommand{\eps}{\varepsilon}

\newcommand{\norm}[1]{\left\lVert#1\right\rVert}
\newcommand{\mcl}[1]{\mathcal{ #1 }}

\newcommand{\fm}[1]{\begin{align*} #1 \end{align*}}
\newcommand{\eq}[1]{\begin{equation}\begin{aligned} #1 \end{aligned}\end{equation}}
\newcommand{\lra}[1]{\langle #1 \rangle}
\newcommand{\red}[1]{{\color{red} #1 }}
\newcommand{\blue}[1]{{\color{blue} #1 }}
\newcommand{\abs}[1]{\left| #1 \right|}
\newcommand{\kh}[1]{\left( #1 \right)}
\newcommand{\wt}[1]{\widetilde{ #1 }}
\newcommand{\wh}[1]{\widehat{ #1 }}

\DeclareMathOperator{\Ric}{Ric}

\newcommand{\rmd}{\mathrm{d}}

\newtheorem{theorem}{Theorem}[section]
\newtheorem{proposition}[theorem]{Proposition}
\newtheorem{lemma}[theorem]{Lemma}

\newtheorem{conjecture}[theorem]{Conjecture}

\newcommand{\perturbsize}{\mathring{\epsilon}}

\theoremstyle{definition}
\newtheorem{definition}[theorem]{Definition}
\newtheorem{remark}[theorem]{Remark}

\newcommand{\Riemanninvariant}{\mathcal{R}}
\newcommand{\RRiemann}{\mathcal{R}_{(+)}}
\newcommand{\LRiemann}{\mathcal{R}_{(-)}}
\newcommand{\dataRRiemann}{\mathring{\mathcal{R}}_{(+)}}
\newcommand{\dataLRiemann}{\mathring{\mathcal{R}}_{(-)}}
\newcommand{\flatRRiemann}{\mathcal{R}_{(+)}^{(\textsf{flat})}}
\newcommand{\flatLRiemann}{\mathcal{R}_{(-)}^{(\textsf{flat})}}

\newcommand{\profileoverr}{\mathscr{A}}

\newcommand{\globalprofileoverr}{\mathscr{B}}

\newcommand{\BothRiemann}{\mathcal{R}_{(\pm)}}

\newcommand{\Riemannfunction}{F}
\newcommand{\InverseRiemannfunction}{I}

\newcommand{\Transport}{\mathbf{B}}
\newcommand{\Lunit}{L}
\newcommand{\uLunit}{\underline{L}}

\newcommand{\flatLunit}{L_{(\textsf{flat})}}
\newcommand{\flatuLunit}{\underline{L}_{(\textsf{flat})}}

\newcommand{\muX}{\breve{X}}
\newcommand{\muuLunit}{\breve{\uLunit}}

\newcommand{\Speed}{c}

\newcommand{\Datafunctionforeikonal}{\mathfrak{U}}

\newcommand{\datasize}{\mathring{\upalpha}}

\newcommand{\tstar}{T_*}
\newcommand{\ustar}{U_*}
\newcommand{\SBT}{V}

\newcommand{\Tboot}{T_{\textsf{Boot}}}

\newcommand{\Tcrease}{T_{\textsf{Crease}}}
\newcommand{\Rcrease}{R_{\textsf{Crease}}}
\newcommand{\Tblowup}{T_{\textsf{Blowup}}}
\newcommand{\Rblowup}{R _{\textsf{Blowup}}}
\newcommand{\TCH}{T_{\textsf{CH}}}
\newcommand{\RCH}{R_{\textsf{CH}}}
\newcommand{\texistenceintermediateboundary}{T_{\textsf{Intermediate-Boundary}}}
\newcommand{\tboundary}{T_{\textsf{Boundary}}}
\newcommand{\ublowup}{U_{\textsf{Blowup}}}
\newcommand{\uboot}{U_{\textsf{Boot}}}
\newcommand{\ucrease}{U_{\textsf{Crease}}}
\newcommand{\rcrease}{R_{\textsf{Crease}}}
\newcommand{\umin}{U_{\textsf{Min}}}

\newcommand{\RIlinear}{\mathfrak{R}}

\newcommand{\integralcurvelunit}{ \mathfrak{t} }
\newcommand{\integralcurveulunit}{ \underline{\mathfrak{t}}}

\newcommand{\lifespanconstant}{\upsigma}

\newcommand{\globalexistenceslab}[1]{\mathscr{S}^{#1}}
\newcommand{\awayglobalexistenceslab}[1]{\mathscr{S}_{\textsf{Away}}^{#1}}
\newcommand{\nearglobalexistenceslab}[1]{\mathscr{S}_{\textsf{Near}}^{#1}}

\newcommand{\similarprofile}{\psi}
\newcommand{\similarprofileoverr}{\widetilde{\mathscr{A}}} 
\newcommand{\antiderivativesimilarprofile}{\Theta}

\newcommand{\globalsimilarprofileoverr}{\widetilde{\mathscr{B}}}

\newcommand{\crease}{\partial_-\mathcal{B}}
\newcommand{\futurecrease}{\partial_-\mathcal{B}_{\textsf{Future}}}
\newcommand{\pastcrease}{\partial_+\mathcal{B}_{\textsf{Past}}}

\newcommand{\futuresinghyp}{\mathcal{B}_{\textsf{Future}}}
\newcommand{\pastsinghyp}{\mathcal{B}_{\textsf{Past}}}

\newcommand{\futureCauchyhor}{\underline{\mathcal{C}}_{\textsf{Future}}^{\textsf{Max}}}
\newcommand{\pastCauchyhor}{\underline{\mathcal{C}}_{\textsf{Past}}^{\textsf{Max}}}

\newcommand{\mylittleo}{o}

\begin{document}

\title{Two global stability theorems for small-data $3D$ compressible Euler solutions: Unique maximal globally hyperbolic developments with gradient-blowup $\&$ smooth global existence in the entire spacetime}

\author{Leonardo Abbrescia}
\address{School of Mathematics, Georgia Institute of Technology, Atlanta, GA 30332}
\email{abbrescia@math.gatech.edu}

\author{Jared Speck}
\address{Department of Mathematics, Vanderbilt University, Nashville, TN 37235}
\email{jared.speck@vanderbilt.edu}

\author{Dongxiao Yu}
\address{Department of Mathematics, Vanderbilt University, Nashville, TN 37235}
\email{dongxiao.yu@vanderbilt.edu}

\begin{abstract}
	We prove a pair of complementary small-data stability theorems that yield the global structure -- to the future and past -- 
	of classical solutions to the isentropic, spherically symmetric $3D$ compressible
	Euler equations. We allow for any equation of state with positive sound speed, 
	except our shock-formation results do not hold for the Chaplygin gas.
	We expect that our results can be extended away 
	from symmetry with the help of a robust set of techniques developed by Christodoulou, Luk--Speck, and Abbrescia--Speck. 
	
	In our first theorem, we study the Cauchy problem for open sets of smooth data that are perturbations of
	trivial data with vanishing velocity and constant positive density. The perturbed data have finite kinetic energy 
	(the initial velocity is even allowed to be zero), and the perturbed density is equal to a positive constant plus
	an ``asymptotically flat tail'' of size $\datasize/r$ as $r \to \infty$, where $\datasize$ is a small, signed, non-vanishing constant.
	Our main theorem yields a complete description of the existence and uniqueness of the maximal globally hyperbolic development (MGHD) of the data, 
	to the future and the past, for solutions that develop gradient-singularities on hypersurfaces that extend to spatial infinity. 
	This is the first MGHD existence-uniqueness-stability result for shock-forming solutions for any 
	multi-dimensional quasilinear hyperbolic wave-like system. While prior works have yielded portions of ``candidate'' MGHDs,
	it is known that neither the existence nor the uniqueness of an MGHD can be inferred from local considerations;
	in an appendix, we show that for many quasilinear hyperbolic PDEs, including the compressible Euler equations, 
	one must construct the entire MGHD and understand the structure of its boundary to know that it exists and is unique.

	In our second theorem, we consider similar initial data, except the small $1/r$ density tail has the opposite sign.
	We prove global existence in the entire spacetime, to the future and past.
 	The solution's asymptotic behavior towards null infinity is not linear, 
	but rather is distorted by a logarithmic bending of the sound cones away from the flat Minkowski ones, leading to 
	logarithmically enhanced dispersion.
	This is the first small-data global existence result in the entire spacetime 
	for any $3D$ quasilinear hyperbolic wave-like system that fails to satisfy the null condition 
	and the weak null condition.
	The main mechanism of stabilization is a global-in-spacetime, though decaying, rarefaction effect, tied to the $1/r$ tail.
	In both theorems, we use nonlinear geometric optics to show that the nonlinear radiation field has a global sign and
	to show that relative to a system of geometric coordinates tied to an eikonal function,
	the nonlinear terms exhibit good null structure, leading to decay.

\end{abstract}

\maketitle

\setcounter{tocdepth}{1}
\tableofcontents

\section{Introduction}
\label{S:INTRO}
We prove a pair of \emph{small}-data
global stability theorems for classical solutions to the Cauchy problem for the $3D$ isentropic compressible Euler equations in spherical symmetry:
\begin{subequations}
        \begin{align}
            \partial_t v^r + v^r \partial_r v^r = - \frac{\Speed^2}{\varrho} \partial_r \varrho, \label{E:INTROVELOCITYRADIALEQUATION} 
							\\
            \partial_t \varrho + \partial_r (\varrho v^r) = -2 \frac{\varrho v^r}{r}, \label{E:INTRODENSITYRADIALEQUATION}
        \end{align}
\end{subequations}
where $\varrho : \mathbb{R} \times [0,\infty) \rightarrow [0,\infty)$ is the density,
$v^r : \mathbb{R} \times [0,\infty) \rightarrow \mathbb{R}$ is the radial velocity component,
$\Speed = \Speed(\varrho) > 0$ is the speed of sound,
and $(t,r) \in \mathbb{R} \times [0,\infty)$ are standard spherical coordinate functions.
The speed of sound is determined by the equation of state
(see Sect.\,\ref{SS:SPHERICALLYSYMMETRICEULER}).
By ``small-data,'' we mean that we study perturbations -- with a small-amplitude \emph{non-vanishing tail} -- of ``trivial'' 
solutions featuring constant positive density and vanishing velocity; see below.

For open sets of smooth data (at time $t = 0$), our work affirmatively answers the following questions about solutions, 
of fundamental mathematical and physical interest:
\begin{quote}
	\begin{itemize}
		\item Is there \emph{any} largest possible subset $\mathbf{M}^{\textsf{Max}}$ 
			of spacetime on which a \underline{classical} (i.e., $C^1$ or smoother) solution exists
			and is uniquely determined by the initial data?
		\item If so, is the region $\mathbf{M}^{\textsf{Max}}$ itself uniquely determined by the initial data?
		\item Is $\mathbf{M}^{\textsf{Max}}$ the entire spacetime? If not, what happens along $\partial \mathbf{M}^{\textsf{Max}}$,
			i.e., why can't the solution be uniquely extended as a classical solution past it?
	\end{itemize}
\end{quote}

In our first theorem, we show the existence and \underline{uniqueness}
of a classical\footnote{Although our work here concerns only classical solutions, it turns out that the global structure of the classical MGHD is important for
\emph{constructing} a class of weak, piecewise-smooth solutions 
that emanate from the gradient-blowup region; see \cite{abbresciaspecknotices} for further discussion. \label{FN:RELEVANTFORSDP}}  
\emph{maximal globally hyperbolic development} (MGHD) for open sets of smooth data that lead to
shock-forming solutions. In particular, $\mathbf{M}^{\textsf{Max}}$ exists and is uniquely determined by the data, 
but it is not the entire spacetime. 
The notion of an MGHD relies on the fact that the causal structure of the compressible Euler equations is governed by a
\emph{solution-dependent Lorentzian metric}, a fact that is not true for a general hyperbolic PDE; 
see Appendix~\ref{A:GLOBALHYPERBOLICITYANDMGHDS}.
Roughly, an MGHD is a largest possible classical solution + region that is determined by the initial data,
and fascinating phenomena can occur along the boundary. Specifically, for the solutions we study, the MGHD future-boundary consists of a singular portion extending to spatial infinity, along which gradient-blowup (i.e., shock formation) occurs, 
and a pre-compact Cauchy horizon, i.e., a null hypersurface that emanates from the point\footnote{In the ``lifted'' picture in $\mathbb{R}^{1+3}$, this point corresponds to a two-dimensional sphere.}
of first blowup and terminates along 
the axis of symmetry $\lbrace r=0 \rbrace$. 
The two future-boundary portions are continuously -- but not differentiably -- joined at a corner (the aforementioned point of first blowup) that we call the \emph{crease}.
The solution behaves similarly for negative times (see Remark~\ref{R:DISCRETESYMMETRY}) and thus analogous results hold for the MGHD past-boundary; 
see Fig.\,\ref{F:FULLMGHD}.

In our second theorem, for related open sets of small, smooth data, 
we prove classical global existence to the \emph{future and past}, i.e., $\mathbf{M}^{\textsf{Max}}$ is the entire spacetime.
Both theorems hold for any equation of state such that the sound speed is positive when the density is positive, 
except for our shock-formation MGHD theorem, which does not apply to the Chaplygin gas equation of state.\footnote{A Chaplygin gas has an equation of state 
$p = C_1 - \frac{C_2}{\varrho}$,
where $C_1$ and $C_2$ are constants.
Due to remarkable null structure in the fluid equations for the Chaplygin gas 
(e.g., the factor $1
				+
				\Speed' 
				\InverseRiemannfunction'$
on RHS~\eqref{E:MUEVOLUTION} completely vanishes for a Chaplygin gas), shock formation is not expected in this case, though our global existence theorem still applies (albeit with simpler asymptotics compared to other equations of state. \label{FN:INTROCHAPLYGINGASISSPECIAL}}

The novelty of this paper can be summarized as follows:
\begin{quote}
	Both of our theorems yield the first truly global, small-data results 
	of their type for any quasilinear wave-like PDE system in 
	$3D$ that fails to satisfy Klainerman's null condition \cite{sK1984} as well as Lindblad--Rodnianski's weak null condition \cite{hLiR2003}. 
	The \emph{global sign of the radiation field} is the key mechanism driving the behavior of the main terms in the equations, 
	while (sound) wave dispersion and nonlinear geometric optics allow us to control the nonlinear error terms.
\end{quote}

While the existence of at least one MGHD follows from a Zorn's lemma argument (see Theorem~\ref{T:EXISTENCEOFMGHDS}),
it should be appreciated that the uniqueness of an MGHD is a highly nontrivial, global, geometric phenomenon that does not always hold. 
In fact, in principle, it is possible that there is not even a single 
\underline{unique} globally hyperbolic development of the data, let alone a maximal one; see Sect.\,\ref{SSS:NOGHD}.
While the celebrated works of Choquet-Bruhat \cite{CB1952} and Choquet-Bruhat--Geroch \cite{cBgR1969} show that 
a unique MGHD always exists for Einstein's equations, in the remarkable paper \cite{fEhRjS2019},
Eperon--Reall--Sbierski gave an example of a quasilinear wave equation in $1D$ and $C^{\infty}$ data 
that evolves into two distinct globally hyperbolic developments such that the corresponding classical solutions
\emph{disagree} on part of the intersection of the two domains, a horrifying breakdown in classical determinism.
For the shock-forming solutions in our first main theorem, 
we were able to prove uniqueness only because we derived the global structure of the \underline{entire MGHD past- and future- boundaries} 
and proved that in $(t,r)$ coordinates, each of the two boundaries is the graph of a continuous function, unlike the example in \cite{fEhRjS2019}; 
see Appendix~\ref{A:GLOBALHYPERBOLICITYANDMGHDS} for further discussion.

In both theorems, the initial data are \emph{small}, smooth perturbations of data with vanishing velocity and a constant positive density. 
More precisely, the perturbed velocity data decays to $0$ at least as fast as $1/r^2$, leading to  
finite kinetic energy,\footnote{By this, we mean that at time $0$, we have 
$\int_0^{\infty} \varrho |v^r|^2 \, r^2 \mathrm{d}r < \infty$. \label{FN:FINITEKINETICENERGY}} 
and the density data converges to a positive constant at the ``asymptotically flat'' rate of $1/r$.
Importantly, in terms of an expansion of the density data in terms of powers of $1/r$ (for $r$ large), 
the $1/r$ term is \emph{non-zero}, with an arbitrary small but non-vanishing amplitude.
This allows us to exploit various monotonicity properties in the \emph{radiation field} 
generated by the main part of the density perturbation tail, the $1/r$ part. In fact, in both theorems, the radiation field has 
a \underline{global sign}, though it decays at a polynomial rate away from the wave zone; see Remark~\ref{R:SIGNOFRADIATIONFIELD}.

It is of interest to note that in both theorems, the data are \emph{not} close to that of a ``simple wave,'' 
i.e., near time $0$, the solution features both inwards-moving and outwards-moving sound waves, in nearly equal proportions.
More precisely, in terms of
a system of spherical Riemann invariants $\mathcal{R}_{\pm}$ 
(see Def.\,\ref{D:SPHERICALRIEMANNINVARIANTS}), we choose the data for each of $\RRiemann,\LRiemann$ to 
be approximately equal; see e.g., 
\eqref{E:RPLUSDATAISCLOSETOBACKGROUNDATTIME0}--\eqref{E:RMINUSDATAISCLOSETOBACKGROUNDATTIME0}.
This corresponds to  the ``very small'' velocity data highlighted above.
While the dynamics of $\RRiemann$ eventually dominates to the future (and $\LRiemann$ dominates to the past),
the domination does not become evident until \emph{dispersive effects} take over, eventually leading to 
asymptotically distinct behavior for $\RRiemann$ and $\LRiemann$. Indeed, the crux of the PDE analysis is about 
\emph{competition between dispersive effects and resonant nonlinearities}.
 
\subsection{Brief motivation and context}
\label{SS:BRIEFMOTIVATION}
Here we provide brief motivation and context for our results. In Sect.\,\ref{SS:HEURISTICDESCRIPTIONOFMAINTHEOREMS}, 
we provide additional context and further explain some key connections between our work and the vast literature on the structure of
classical solutions to $3D$ compressible Euler equations and related nonlinear wave equations.

Our MGHD results for shock-forming solutions
are especially motivated by Christodoulou's groundbreaking monograph \cite{dC2007},
in which he derived a large portion of a future-MGHD for open sets of shock-forming solutions -- without symmetry assumptions -- 
to the $3D$ irrotational and isentropic relativistic Euler equations in the small-data regime, i.e., small perturbations of data with constant positive density and vanishing velocity. The data-perturbations in \cite{dC2007} 
did not have a $1/r$ tail, but rather were compactly supported and hence, thanks to domain of dependence considerations, 
the solutions were trivial in the exterior of a flat, outgoing sound cone. 
In \cite{dC2007}, Christodoulou developed many important techniques for studying MGHDs, tied to nonlinear geometric optics, 
and offered a vision for the future of the field.
The portion of an MGHD derived in \cite{dC2007}*{Theorem 15.1, Proposition 15.3} was not explicitly constructed, but rather was described as a union of 
globally hyperbolic developments of the data tied to an uncountably infinite family of foliations of a portion of spacetime by sound cones
and flat, spacelike hypersurfaces that terminate on the ``gradient-blowup'' portion of the MGHD boundary. 
Thanks to the properties of the solution induced by 
our $1/r$ data-tails, our description of the MGHD is fully constructive and complete, which in particular allows us to conclude uniqueness.
Another crucial result of Christodoulou's monograph is \cite{dC2007}*{Theorem 13.1}, which shows that in the small-data regime, 
for compactly supported, isentropic, and irrotational data perturbations, 
away from the interior region, 
the only kinds of singularities that can form in principle are shocks due to gradient-blowup caused by the infinite density of outgoing characteristics.
Put differently, 
Christodoulou proved \emph{conditional global existence} in an exterior region, i.e., thanks to dispersive effects,
global existence holds as long as outgoing characteristics never become infinitely dense. To date, his analytical framework is the only one that has been able to prove conditional global existence; 
all subsequent works on multi-dimensional shock formation either used his framework or were restricted to a limited data regime 
that enjoyed additional properties guaranteeing that a shock quickly forms in the flow.
Although we do not explicitly state an analogous result in this paper, 
the proofs of our main theorems show that a conditional global existence result, similar to the one in \cite{dC2007}*{Theorem 13.1},
holds for data perturbations with a small $1/r$ tail 
-- in the entire spacetime --
albeit with more complicated asymptotics induced by the tail.

Our global existence results are motivated by the backwards scattering results derived in third author's paper \cite{dY2025}, 
which \emph{constructed} an interesting class of candidate ``asymptotic profiles'' and used them to construct
non-trivial future-global solutions to the $3D$ irrotational and isentropic compressible Euler equations as well as related nonlinear wave equations. 
We clarify that \cite{dY2025} yielded an existence result -- only to the future -- and the stability of the constructed solutions, from the point of view of the Cauchy problem, was not investigated in that paper. 

Although we have treated only the compressible Euler equations in detail, our methods are robust and applicable to other quasilinear wave equations 
that fail to satisfy the null condition; the novelty of this paper is our discovery that a regime of ``$\frac{1}{r}$ tail'' initial data can be the \emph{cause} of tractable \underline{global, nonlinear} dynamics tied to the radiation field having a sign, rather than the precise details of the equations. 
For this reason, we expect that our results can be generalized to other classes of hyperbolic PDEs besides wave equations. 
Moreover, even though our results are in spherical symmetry, they are compatible with a set of robust multi-dimensional techniques 
developed by Christodoulou \cite{dC2007}, Speck--Luk \cites{jLjS2020a,jS2019c}, 
and Abbrescia--Speck \cites{lAjS2022,abbrescia2025remarkable,abbrescia2026emergencecauchyhorizoncrease}, which we anticipate can be used to generalize our results away from symmetry.
We especially highlight that the techniques developed by Abbrescia--Speck apply in particular to perturbations of symmetric solutions 
(including the ones studied in this paper), 
a class of solutions for which various convexity assumptions used in other approaches to studying shock formation do not hold.
To keep the paper a manageable length, and to help readers navigate the main ideas, we have treated in detail a well-chosen open set of initial data, 
based on perturbing a specific one-parameter family of profiles. In Appendix~\ref{A:EXTENDRESULTSTOOTHERPROFILES}, 
we outline how to generalize our results to a much larger class of data.

\subsection{Heuristic description of the main theorems and context}
\label{SS:HEURISTICDESCRIPTIONOFMAINTHEOREMS}
We now describe our results in more detail, though still at a heuristic level. Readers interested in the detailed theorems can immediately jump to Theorems~\ref{T:MAINMGHDEXISTENCETHEOREM} and \ref{T:MAINGLOBALEXISTENCETHEOREM}. 

\subsubsection{Globally stable, unique maximal globally hyperbolic developments for shock-forming solutions}
\label{SSS:MGHDTHEOREMHEURISTIC}
In our first theorem, using a constructive approach, we provide the first proof of existence \emph{and} uniqueness of the \emph{maximal globally hyperbolic development} (MGHD) of data for shock-forming solutions to the equations. We refer to Appendix~\ref{A:GLOBALHYPERBOLICITYANDMGHDS} for a rigorous definition of an MGHD and
corresponding background material. Here, we will only provide a heuristic description. 
Roughly, for a hyperbolic wave-like system, a \emph{globally hyperbolic development} (GHD) of given smooth data is a region of spacetime $\mathbf{M}\subset \R^{1+3}$ that is equal to the union of the domain of dependence of all of its points, and on which the solution exists classically. 
A GHD is called \emph{maximal} if it is inextendible as a GHD. That is, $\mathbf{M}$ is maximal if for any GHD $\check{\mathbf{M}} \subset \R^{1+3}$ 
such $\mathbf{M} \subset \check{\mathbf{M}}$, it must hold that
$\check{\mathbf{M}} = \mathbf{M}$. One can therefore view an MGHD as a largest spacetime region on which fixed initial data launches a classical solution. In other words, understanding an MGHD is equivalent to one of the most fundamental questions in PDE theory: to what extent can one construct a classical solution from smooth data?  The (possibly empty) \emph{boundary} of an MGHD consists of the spacetime points on which the corresponding classical solution develops a singularity or points beyond which the initial data no longer uniquely determines a classical solution.
We already highlight the following point, fleshed out below:
\begin{quote}
	In principle, there might be two \emph{distinct} maximal classical solutions, defined on regions $\mathbf{M}_1$ and $\mathbf{M}_2$,
	such that the corresponding maximal classical solutions disagree at some points in $\mathbf{M}_1 \cap \mathbf{M}_2$. 
	Such a \underline{breakdown in classical determinism} actually happens in some model quasilinear wave equation 
	solutions launched by smooth data; see 
	Sect.\,\ref{SSS:EPREXAMPLE}.
	In our main theorem, we show that this does not happen, at least for the open sets of data that we study.
\end{quote}

\begin{figure}
		\centering	
			\begin{overpic}[scale=.65, grid = false, tics=3, trim=-.5cm -1cm -1cm -.5cm, clip]{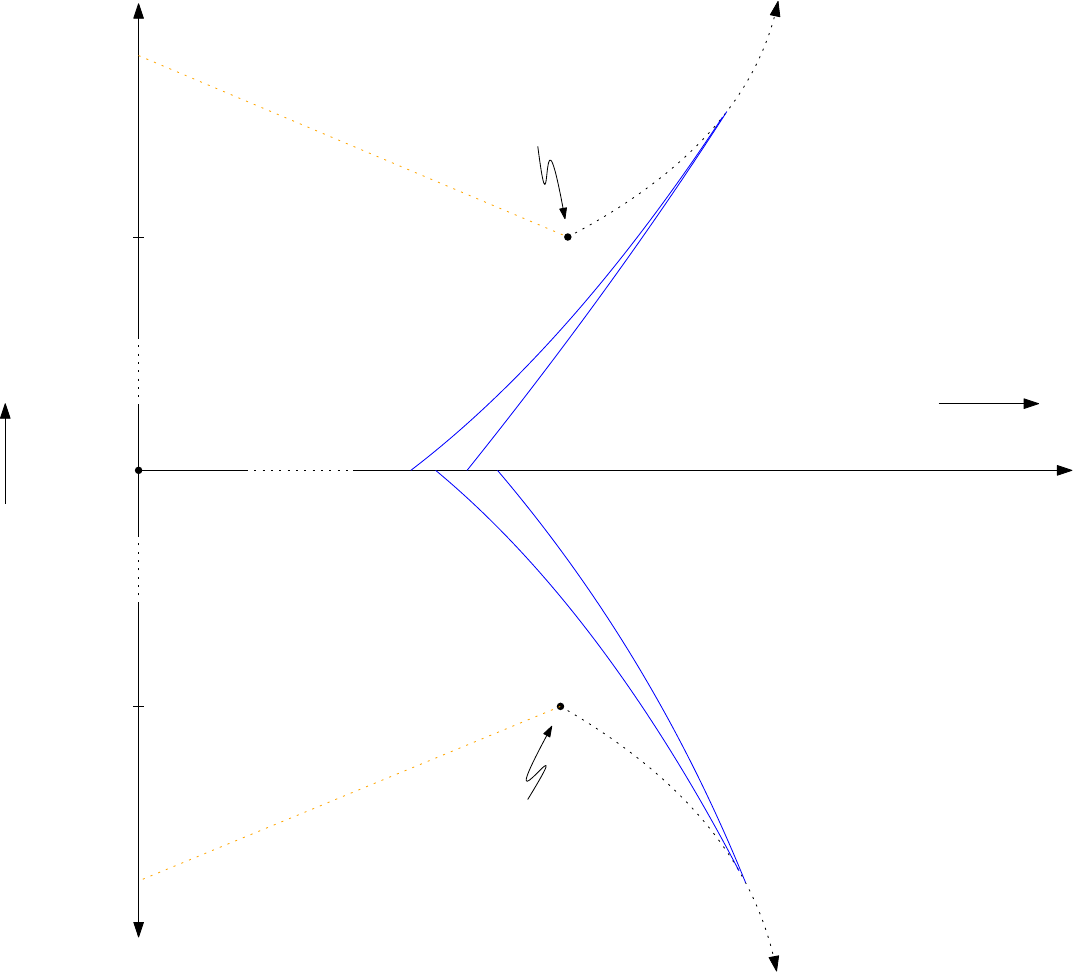} 
				\put (22,83) {\small{$\futureCauchyhor$}}
				\put (69,81) {\small{$\futuresinghyp$}}
				\put (69,9) {\small{$\pastsinghyp$}}
				\put (21,12) {\small{$\pastCauchyhor$}}
				\put (8,48) {\small{$(0,0)$}}
				\put (87,56) {\small{$r$}}
				\put (1,48) {\small{$t$}}
				\put (30,78.5) {\small{$\futurecrease = (\Tcrease^{\textsf{Future}}, \Rcrease^{\textsf{Future}})$}}
				\put (0,68) {\small{$(\Tcrease^{\textsf{Future}},0)$}}
				\put (31,17) {\small{$\pastcrease = (\Tcrease^{\textsf{Past}},\Rcrease^{\textsf{Past}})$}}
				\put (0,27) {\small{$(\Tcrease^{\textsf{Past}},0)$}}
			\end{overpic}
			\caption{The MGHD of a shock-forming solution from Theorem~\ref{T:ROUGHMGHDTHEOREM}
			in the differential structure of standard spherical coordinates $(t,r)$. The solution is nearly symmetric about the $r$-axis. 
			$\mathbf{M}^{\textsf{Max}}$ is the region between $\futureCauchyhor\cup\futuresinghyp\cup\pastCauchyhor\cup\pastsinghyp$. 
			Every point on this figure represents a sphere in the ``lifted'' spacetime $\R^{1+3}$. 
			When viewed as a subset of $\mathbb{R}^{1+3}$, $\mathbf{M}^{\textsf{Max}}$ is open, which means that in $(t,r)$ coordinates,
			$\mathbf{M}^{\textsf{Max}}$ does not contain any of its boundary portion except the one on the axis of symmetry $\lbrace r = 0 \rbrace$.
			If the initial data is of size $\datasize$, then 
			$\Tcrease^{\textsf{Future}},\Rcrease^{\textsf{Future}},|\Tcrease^{\textsf{Past}}|,\Rcrease^{\textsf{Past}} \approx 
			\exp(\frac{C}{\datasize})$. The solid blue curves denote outgoing characteristic curves. 
			These characteristics tangentially graze the future- and past- singular boundaries, along which they have infinite density.}
		\label{F:FULLMGHD}
		\end{figure}

We now state a rough version our main MGHD theorem; see Theorem~\ref{T:MAINMGHDEXISTENCETHEOREM} for the precise statement
and Fig.\,\ref{F:FULLMGHD} for a graphical depiction of the results.

	\begin{theorem}[Rough version of the main MGHD theorem] \label{T:ROUGHMGHDTHEOREM}
		For the spherically symmetric isentropic compressible Euler equations in $3D$, 
		under any $C^4$ equation of state with positive sound speed except the Chaplygin gas,
		there exists an open set of $C^3$, 
	``compressive'' data with small velocity, finite kinetic energy, and with density that is a small perturbation of a positive constant $\overline{\varrho}$
		and that converges to $\overline{\varrho}$ at the asymptotically flat rate of $\frac{\datasize}{r}$ as $r \to \infty$,
		where $\datasize > 0$ is any sufficiently small amplitude.
		Given such data, the following results hold (see Fig.\,\ref{F:FULLMGHD}):
			\begin{enumerate}
				\item (Classical solution). With respect to the differential structure of the standard Cartesian coordinates on $\mathbb{R}^{1+3}$,
					there is a spherically symmetric maximal globally hyperbolic development 
					$\mathbf{M}^{\textsf{Max}} \subsetneq \mathbb{R}^{1+3}$, on which the solution is classical.\footnote{Although we only treat $C^3$ data in 	
					Theorem~\ref{T:ROUGHMGHDTHEOREM},
					standard propagation of regularity arguments could be used to show that 
					more regular data would lead to more regular solutions, i.e., the solutions are as smooth as the data, except
					that our $C^k$-type norms involve $r$-weights at the top order that degenerate near $r=0$.\label{E:MGHDPROPAGATIONOFREGULARITY}}
					By globally hyperbolic, we mean with respect to the acoustical metric $\mathbf{g}$, the solution-dependent Lorentzian metric
					that determines the intrinsic causal structure of the solution; see \eqref{E:ACOUSTICALMETRIC}--\eqref{E:INVERSEACOUSTICALMETRICINSPHERICALSYMMETRY}.
				\item  (MGHD boundary). Viewed as a subset of $(t,r)$-coordinate space, 
					the boundary of $\mathbf{M}^{\textsf{Max}}$ consists of the following sets, depicted in Fig.\,\ref{F:FULLMGHD}:
					\begin{itemize}
						\item A future-\textbf{crease} $\futurecrease$, which is the unique spacetime point 
						(corresponding to a $2D$ sphere in $\R^{1+3}$) whose temporal component is the first positive time for which the gradient of the physical fluid variables 
						(i.e., the velocity and density) blows up, 
						while the undifferentiated fluid variables remain bounded.\footnote{The undifferentiated variables are in fact $C^{1/3}$ at $\futurecrease$ and $\pastcrease$; see, for example, \cite{jLjS2024}*{Corollary~4.5} for a proof. \label{FN:CONETHIRDREGULARITY}} 
						There is also a past crease $\pastcrease$, which is the unique spacetime point whose temporal component is the first negative time for which the gradient of the fluid variables blows up, while the fluid variables remain bounded. 
						The density of the future-outgoing characteristics is infinite at $\futurecrease$, 
						and the density of the past-outgoing characteristics is infinite at $\pastcrease$. 
						The two creases are not part of 
						$\mathbf{M}^{\textsf{Max}}$, but they are part of its topological closure: 
						$\futurecrease, \pastcrease \subset \overline{\mathbf{M}}^{\textsf{Max}}$.
						\item A future-\textbf{singular boundary} $\futuresinghyp$, which is a $C^1$ curve-portion 
						(corresponding to a hypersurface-portion in $\R^{1+3}$) that emanates from the future-crease and flows towards future temporal and spatial infinity, along which the gradient of the physical fluid variables continue to blow-up, while the undifferentiated fluid variables 
						remain bounded.\footnote{The undifferentiated variables can in fact be shown to be $C^{1/2}$ at points in 
						$\futuresinghyp \backslash \futurecrease$, and similarly for points in $\pastsinghyp \backslash \pastcrease$.
						This can be proved by modifying the proof of
						\cite{jLjS2024}*{Corollary~4.5}, which addressed only the $C^{1/3}$-regularity at $\futurecrease$ and $\pastcrease$. \label{FN:CONEHALFREGULARITY}} There is also a past-singular boundary $\pastsinghyp$, which is a $C^1$ curve-portion (i.e., a hypersurface-portion in $\R^{1+3}$) 
						that emanates from the past-crease and flows towards past temporal and spatial infinity, along which the gradient of the physical fluid variables continues to blow-up, while the undifferentiated fluid variables remain bounded. 
						The density of the future-outgoing characteristics is infinite along $\futuresinghyp$, 
						and the density of the past-outgoing characteristics is infinite along $\pastsinghyp$. 
						The two singular boundaries are not part of $\mathbf{M}^{\textsf{Max}}$, 
						but they are part of its topological closure: $\futuresinghyp, \pastsinghyp \subset  \overline{\mathbf{M}}^{\textsf{Max}}$.
						\item A future-\textbf{Cauchy horizon} $\futureCauchyhor$, which is a $C^1$ curve-portion that emanates from $\futurecrease$ and terminates to the future 
						when it intersects the axis of symmetry $\{r=0\}$. There is also a past-Cauchy horizon $\pastCauchyhor$, which is a $C^1$ curve-portion that emanates from $\pastcrease$ towards the past and terminates to the past when it intersects the axis of symmetry $\{r=0\}$. The solution remains smooth along both the future- and past- Cauchy horizons
						(which by definition do not contain their crease boundary points, where gradient-blowup occurs). The two Cauchy horizons are not part of $\mathbf{M}^{\textsf{Max}}$, 
						but they are part of its topological closure: $\futureCauchyhor, \pastCauchyhor \subset  \overline{\mathbf{M}}^{\textsf{Max}}$.
					\end{itemize}
				\item  (Maximality and Uniqueness). 
					$\mathbf{M}^{\textsf{Max}}$ is maximal (i.e., there is no way to uniquely extend the classical solution beyond it), 
					and the region $\mathbf{M}^{\textsf{Max}}$ itself as well as the solution on it are uniquely determined by the data. 
					Moreover, \textbf{any} other globally hyperbolic development of the data is
					defined only on a \textbf{subset} of $\mathbf{M}^{\textsf{Max}}$, and the corresponding solution agrees with the one 
					from Point $1$.
				\item  (Behavior of the outgoing characteristics). 
					There is an eikonal function $u$ defined in a large subset of $\mathbf{M}^{\textsf{Max}}$ away from the origin, 
					extending to spatial infinity, such that for positive times, the 
					level sets $\lbrace u = \mbox{constant} \rbrace$ are outgoing characteristics
					that are distortions of the flat characteristics $\lbrace t-r = \mbox{constant} \rbrace$ by error terms that grow logarithmically in time.
					Relative to the $(t,u)$ coordinates, the nonlinear solution decays like a linear wave towards 
					the trivial solution with vanishing velocity and constant density $\overline{\varrho}$.
					Similarly, for negative times, the level sets of the eikonal function are logarithmic
					distortions of the flat characteristics $\lbrace t+r = \mbox{constant} \rbrace$. 
				\item  (Global compression). The density of the future- and past- outgoing characteristics grows as 
						$|t|$ increase, until it becomes infinite on the two singular boundaries. 
						That is, the solution is globally compressive. 
						The compression is caused by the radiation field having a \underline{global sign}; see Remark~\ref{R:SIGNOFRADIATIONFIELD}.
			\end{enumerate}
		\end{theorem}

\begin{remark}[Simplified notation in the bulk of the paper]
	\label{R:SIMPLIFIEDNOTATION}
	The bulk of the paper concerns PDE analysis for times $t \geq 0$.
	To simplify the notation, in the bulk, we will sometimes drop the label ``Future'' from various symbols, i.e., 
	we denote the time $\Tcrease^{\textsf{Future}}$ from Fig.\,\ref{F:FULLMGHD} by $\Tcrease$,
	we denote the radial value $\rcrease^{\textsf{Future}}$ from Fig.\,\ref{F:FULLMGHD} by $\rcrease$,
	etc.
\end{remark}

At its most primitive level, Theorem~\ref{T:ROUGHMGHDTHEOREM} is a stable shock-formation result for solutions arising from smooth data.
In recent years, there has been an explosion of mathematical activity tied to shock-formation problems,\footnote{Other kinds of compressible Euler singularities, 
such as the celebrated ``implosions'' \cites{fMpRiRjS2022a,fMpRiRjS2022b}, can also arise from smooth data. Unlike shocks, implosions seem to be unstable.} 
especially in multi-dimensions. 
The first multi-dimensional results were Alinhac's important papers \cites{sA1999a,sA1999b} on quasilinear wave equations for ``non-degenerate'' 
shock-forming solutions that enjoy convexity properties. His approach, which relied on a Nash--Moser iteration scheme to close the energy estimates, 
was able to follow the solution until the flat constant-time hypersurface 
$\Sigma_{T_{\textsf{Shock}}}$ of first-blowup, but not to larger regions of spacetime. Thanks to the convexity assumptions, within $\Sigma_{T_{\textsf{Shock}}}$, the singular points were \emph{isolated}, i.e., Alinhac derived a finite number of isolated points on the boundary of an MGHD.

As we highlighted in Sect.\,\ref{SS:BRIEFMOTIVATION},
in his breakthrough monograph \cite{dC2007}, Christodoulou developed techniques much stronger than Alinhac's, 
with a profound geometric perspective, and applied them to study large portions of an MGHD, for small-data shock-forming solutions 
to the wave equations of irrotational and isentropic relativistic fluid mechanics. Crucially, Christodoulou avoided Nash--Moser-type estimates, and instead closed
the energy estimates with the help of techniques that he co-developed with Klainerman in their celebrated proof of the stability of the Minkowski spacetime \cite{dCsK1993};
this is a main reason why he was able to study a much larger MGHD portion compared to Alinhac.
The portions of an MGHD studied in \cite{dC2007} featured a large, though implicit, portion of an MGHD boundary where gradient-blowup occurs.
No Cauchy horizon was derived in \cite{dC2007}, but the ``tilted'' foliations used there yielded a region of existence close to one. 
New formulations of the fluid equations exhibiting remarkable structures were derived in \cites{jLjS2020a,jS2019c},
and those equations were used to prove the first stable shock formation result for $2D$ isentropic compressible fluids with vorticity in \cite{jLjS2018}.
For $2D$ isentropic compressible fluids with vorticity under adiabatic equations of state, with data of large gradient size $\frac{1}{\epsilon}$ that enjoy an Alinhac-type convexity assumption, \cite{shkoller2024geometry} provided a description of a portion of an MGHD that lies between two
parallel constant-time slices separated by distance $\mathcal{O}(\epsilon)$. This included an $\mathcal{O}(\epsilon)$-sized portion of a
singular boundary, where gradient-blowup occurs, and an $\mathcal{O}(\epsilon)$-size portion of a Cauchy horizon, where no singularity develops.
\cite{jLjS2024} proved stable shock formation for nearly plane symmetric solutions to 
the $3D$ compressible Euler equations with vorticity and entropy, which is much more difficult than the $2D$ isentropic case due to the presence of a vorticity stretching term and an entropy-gradient term in the equations, both of which lead to a challenging regularity theory. 
The main results of \cite{jLjS2024} did not require any convexity assumption, but the solutions were studied only until the time of first gradient-blowup, i.e., most of the MGHD boundary was not studied there. Under the same setup in $3D$, \cite{lAjS2022,abbrescia2026emergencecauchyhorizoncrease} derived 
an $\mathcal{O}(1)$-size portion of both the singular boundary and the Cauchy horizon, as well as the full structure of the crease, 
using only a  ``transversal convexity'' assumption, which even in $1D$ is essential to ensure that a continuum of characteristics 
does not collapse into a low-dimensional set; were that to occur, it would be very challenging to study the questions of 
existence and uniqueness of an MGHD. Crucially, the techniques developed in \cite{lAjS2022,abbrescia2026emergencecauchyhorizoncrease} are applicable to perturbations -- away from symmetry --
of the solutions that we study in this paper. The top-order regularity theory in \cite{lAjS2022,abbrescia2026emergencecauchyhorizoncrease}
relied on remarkable energy identities derived in \cite{abbrescia2025remarkable}, 
which allow one to propagate an extra degree worth of regularity for the vorticity and entropy on arbitrary globally hyperbolic domains,
which in turn allows one to treat the flow as a perturbation of an irrotational and isentropic one.

We also mention that for solutions that enjoy an Alinhac-type convexity assumption, there is another approach, initiated in \cite{tBsSvV2019a}, 
to following compressible Euler solutions to the time of first gradient-blowup, based on self-similarity and modulation parameters.
For deriving the full structure of an MGHD, especially in the context of solutions that exhibit a \emph{competition between dispersion and nonlinearities}, 
the framework initiated by Christodoulou in \cite{dC2007} is by far the most powerful one to date. In particular, 
it interacts well with the large array of techniques that have been developed to study the long-time behavior of solutions to 
wave-like equations, and it features many geometric ingredients that are of immense value for following solutions up to the MGHD boundary. 
There are many other works that address various aspects of shock formation; 
for a more detailed discussion of the history of the field as well as additional context,
we refer readers to the surveys \cites{abbresciaspecknotices,abbrescia2023relativistic,gHsKjSwW2016}. 

A surprisingly subtle fact concerning the uniqueness aspect of the MGHD from Theorem\,\ref{T:ROUGHMGHDTHEOREM} is that uniqueness \emph{cannot} be deduced from analysis on localized regions in spacetime or even from knowing \emph{existence}! In fact, the proof of uniqueness of $\mathbf{M}^{\textsf{Max}}$ fundamentally relies on knowing 
the complete global structure and topology of its boundary $\partial \mathbf{M}^{\textsf{Max}}$, 
as is depicted in Fig.\,\ref{F:FULLMGHD}. We discuss these issues in detail in Appendix~\ref{A:GLOBALHYPERBOLICITYANDMGHDS}. 
Here, we highlight that in the breakthrough paper \cite{fEhRjS2019}, Eperon--Reall--Sbierski proved that MGHD-uniqueness holds provided that 
$\mathbf{M}^{\textsf{Max}}$
\emph{everywhere lies on one side of its boundary}; see Theorem~\ref{T:UNIQUEONESIDEDNESS} for the precise statement. Surprisingly, for a given quasilinear wave equation,
this property is \emph{not guaranteed to hold} for all MGHDs of smooth data. More precisely, in that same paper, the authors gave a remarkable constructive example in one spatial dimension of a quasilinear wave equation for which \emph{a fixed} $C^\infty$ initial data evolved classically \textbf{in two distinct globally hyperbolic fashions} on regions of spacetime which overlap. More precisely, they found a smooth initial data set which gave rise to two classical solutions $\Phi_1, \Phi_2$, each respectively defined on regions $\mathbf{M}_1,\mathbf{M}_2 \subset \R^{1+1}$ such that $\mathbf{M}_1 \neq \mathbf{M}_2$ and such that there exists a $p \in \mathbf{M}_1 \cap \mathbf{M}_2$ with $\Phi_1(p) \neq \Phi_2(p)$.
Their result is an example of a strong form of \emph{failure of classical determinism}, and a key conclusion of Theorem~\ref{T:ROUGHMGHDTHEOREM}
is that such failure does \emph{not} occur for the set of compressible Euler data that we treat in this paper. 

The above discussion suggests the following open problems, of intense mathematical and physical interest:
\begin{quote}
	\begin{itemize}
		\item (Non-unique MGHDs). Is there any smooth data for the compressible Euler equations that yields two distinct, 
				intersecting globally hyperbolic developments $\mathbf{M}_1$ and $\mathbf{M}_2$
				such that
				$\mathbf{M}_1 \neq \mathbf{M}_2$,
				and such that the corresponding solutions \emph{disagree} at some point $p \in \mathbf{M}_1 \cap \mathbf{M}_2$?
		\item (Not even a single unique GHD).
				Is there any smooth data for the compressible Euler equations on a spacelike Cauchy hypersurface $\Sigma$ 
				such that \emph{there is no neighborhood $\mathcal{O}$ of $\Sigma$} 
				on which one can guarantee that all classical solutions that are defined on globally hyperbolic subsets
				must agree at every point in $\mathcal{O}$?
		\end{itemize}
\end{quote}
While at first glance, the second question above might seem like it follows from local well-posedness, it in fact does not! 
See Sect.\,\ref{SSS:NOGHD} for further discussion.
An example would show that for the given data, 
there is no neighborhood of the Cauchy hypersurface on which classical determinism holds.

\subsubsection{Global existence}
\label{SSS:GLOBALEXISTENCEHEURISTIC}
In our second main theorem, we consider a related but distinct set of initial data (see Remark~\ref{R:TWOTHEOREMSDATACONNECTION})
and prove global existence to the future and past.\footnote{This is also an MGHD existence and uniqueness theorem. In this case, the MGHD is the solution on the entire spacetime, so there is no boundary. \label{FN:GLOBALISALSOMGHD}} 
Here we give a rough statement of the theorem; see Theorem~\ref{T:MAINGLOBALEXISTENCETHEOREM}
for the precise result. 

	\begin{theorem}[Rough version of the main global existence theorem] \label{T:ROUGHGLOBALEXISTENCETHEOREM}
			For the spherically symmetric isentropic compressible Euler equations in $3D$, 
			under any $C^3$ equation of state with positive sound speed,
			there exists an open set of\footnote{Note that Theorem~\ref{T:ROUGHGLOBALEXISTENCETHEOREM} assumes one fewer derivatives compared to Theorem~\ref{T:ROUGHMGHDTHEOREM}.
The reason is that in the proof of Theorem~\ref{T:ROUGHMGHDTHEOREM}, the extra derivative is used only to derive the structure of the singular boundary 
and Cauchy horizon, which are not present for the global solutions of Theorem~\ref{T:ROUGHGLOBALEXISTENCETHEOREM}.
			\label{FN:EXTRADERIVATIVENOTNEEDEDFORGLOBAL}}
			$C^2$, ``rarefactive'' data with small velocity, finite kinetic energy, and with density that is a small perturbation of a positive constant $\overline{\varrho}$
			and that converges to $\overline{\varrho}$ at the asymptotically flat rate of $\frac{1}{r}$ as $r \to \infty$.
			For the data, the following results hold.
			\begin{enumerate}
				\item (Global existence). The corresponding classical\footnote{Although we only treat $C^2$ data in Theorem~\ref{T:ROUGHGLOBALEXISTENCETHEOREM},
					standard propagation of regularity arguments could be used to show that 
					more regular data would lead to more regular solutions, i.e., the solutions are as smooth as the data, except
					that our $C^k$-type norms involve $r$-weights at the top order that degenerate near $r=0$.\label{E:GLOBALEXISTENCEPROPAGATIONOFREGULARITY}} 
					solution is
					global in space and time to the future and the past.
				\item (Behavior of the outgoing characteristics). There is an eikonal function $u$ defined on $\lbrace r \geq 1 \rbrace$, such that for $t \geq 0$,
					the level sets $\lbrace u = \mbox{constant} \rbrace$ are outgoing characteristics
					that are distortions of the flat characteristics $\lbrace t-r = \mbox{constant} \rbrace$ by error terms that grow\footnote{For the Chaplygin gas 
					equation of state, these error terms grow much more slowly.} 
					logarithmically in time.
					For $t \leq 0$, the
					level sets $\lbrace u = \mbox{constant} \rbrace$ are outgoing characteristics
					that are distortions of the flat characteristics $\lbrace t + r = \mbox{constant} \rbrace$ by error terms that grow logarithmically in time.
					Relative to the $(t,u)$ coordinates, the nonlinear solution decays, to the future and past, 
					like a linear wave, towards the trivial solution with vanishing velocity and constant density $\overline{\varrho}$.
				\item (Global rarefaction).
					The outgoing characteristics separate logarithmically in time as $|t| \to \infty$, that is, the solution is globally rarefactive.
					The rarefaction is caused by the radiation field having a \underline{global sign}; see Remark~\ref{R:SIGNOFRADIATIONFIELD}.
			\end{enumerate}
\end{theorem}
				 
	\begin{remark}[Connection between the data for the two theorems]
		\label{R:TWOTHEOREMSDATACONNECTION}
		Roughly speaking, relative to a system of spherical Riemann invariants (see Def.\,\ref{D:SPHERICALRIEMANNINVARIANTS}), the initial data for
		Theorem~\ref{T:ROUGHGLOBALEXISTENCETHEOREM} has the opposite sign as the data for Theorem~\ref{T:ROUGHMGHDTHEOREM}.
	\end{remark}

A surprising aspect of our global existence theorem is that equations 
\eqref{E:INTROVELOCITYRADIALEQUATION}--\eqref{E:INTRODENSITYRADIALEQUATION}
fail to satisfy the null condition and the weak null condition. 
By this, we roughly mean that the once-differentiated fluid equations always feature Riccati-type resonant interactions,\footnote{For example, 
aside from the case of the Chaplygin gas equation of state, RHS~\eqref{E:RIGHTRIEMANNEQUATIONINTERIORREGION} always features
a dangerous Riccati-term of type $r (\uLunit \RRiemann)^2$, stemming from the commutator term $- r [\uLunit,\Lunit]\RRiemann$, studied in more detail in equation~\eqref{E:COMMUTATOROFLANDULUNIT}.} aside from the case of the Chaplygin gas equation of state. Below, we discuss these resonant terms in more detail. 
In many classically studied \emph{small}-data problems for nonlinear hyperbolic PDEs, 
failure of the null condition and weak null condition has been shown to lead to singularity formation tied to the presence of the resonant terms. 
However, resonant interactions have previously been exploited to show \emph{future}-global existence when the data are sufficiently 
\emph{large and have a ``good'' sign}.
In particular, there are important compressible Euler papers showing that for adiabatic equations of state with initially small density and 
an initially \emph{large, signed} velocity, \emph{future}-global existence holds.
For example, in \cite{mGra1998}, Grassin assumed that the initial density is small and that the initial velocity gradient $\nabla v$, viewed as a matrix,
has spectrum that is bounded from below, strictly away from $0$. He then proved \emph{future}-global existence and dispersive estimates for the solution, due to spreading effects tied to the large velocity gradient. However, his result does not extend to the past direction because the spectrum of the initial $\nabla v$ has the wrong sign for past-global existence. 
In contrast, our initial velocity can be arbitrarily small (even zero) with unsigned gradient tensor, 
and our global existence result holds to the future \emph{and} past.
In a nutshell, for our global existence theorem, 
the initial data evolve to create a radiation field that drives the resonant terms to have a \emph{good rarefactive sign}
(see Remark~\ref{R:SIGNOFRADIATIONFIELD}), to the past and future, which causes the outgoing characteristics to separate and stabilizes the solution, though with distorted asymptotics. 
In our shock-forming MGHD theorem, the radiation field has the opposite sign, leading to compression, and the resonant terms cause gradient-blowup to occur in the solution, all the way out to spatial infinity,
in a form of ``controlled blowup.''

Next, we highlight the important papers \cites{mHjJ2018a,mHjJ2018b}, which proved future-global existence for open sets of \emph{expanding} solutions to the 
compressible Euler equations under various adiabatic equations of state. These are solutions in which the density is compactly supported and satisfies the physical vacuum boundary condition at the fluid-vacuum boundary, and the velocity is approximately radial and size $1$ at the vacuum boundary. While not directly related to our work here, these results provide fascinating examples of how outward velocity motion can future-globally stabilize fluids (those solutions are possibly unstable to the past, due to the past-inward velocity motion).

A simplified model that features \emph{both} global compression and global rarefaction is our recent work \cite{abbresciaBlueSierskiSpeck2025quasilinear}. There, we studied the Cauchy problem for a quasilinear wave equation on $\mathbb{R}^{1+1}$
that fails to satisfy the null condition. The initial data on $\{t=0\}$ \emph{monotonically decreases} from a positive constant $\mathfrak{c}_1 > 0$ 
at $x = -\infty$ to a negative constant $\mathfrak{c}_2 < 0$ at $x = +\infty$. This monotonicity forced the Riccati terms in the equation to have 
a global \emph{compressive} sign towards the future and a global \emph{rarefactive} sign towards the past. In particular, the solutions studied in \cite{abbresciaBlueSierskiSpeck2025quasilinear} developed a future-crease where a shock forms and from which a future-Cauchy horizon and future singular boundary emanate. The Cauchy horizon extends all the way to negative spatial infinity (hence, it fails to be pre-compact, which is different than in this paper) and the singular boundary extends all the way to positive spatial infinity. 
On the other hand, for the same data, as $t \downarrow - \infty$, the characteristics smoothly and monotonically expand and the solution therefore
enjoys smooth past-global existence. As a whole, the boundary of this MGHD also enjoyed the global structure that 
is sufficient to conclude MGHD-uniqueness. 
However, we stress that the results of \cite{abbresciaBlueSierskiSpeck2025quasilinear} were much easier to derive compared to those of this paper due to the lack of dispersive effects. We also stress the difference in the topology on the ambient spacetime of $\R^{1+1}$ for \cite{abbresciaBlueSierskiSpeck2025quasilinear} compared to the ambient spherically symmetric spacetime $\R\times[0,\infty)$ in the present paper, which in $(t,r)$ coordinates has a boundary $\{r=0\}$. Morally, it is the \emph{reflection} of waves off of $\{r = 0\}$ that preserves the sign of the radiation field towards the future and past, leading to the nearly-time-symmetric results of Theorems~\,\ref{T:ROUGHMGHDTHEOREM} and \ref{T:ROUGHGLOBALEXISTENCETHEOREM}. There is no reflection of waves in $\R^{1+1}$, and hence the solution is compressive to the future and rarefactive to the past.

We also mention the interesting paper \cite{gCfeKyHyS2025}, in which the authors constructed some future-global, spherically symmetric solutions
to the compressible Euler equations for adiabatic equations of state $p = \varrho^{\upgamma}$, $1 < \upgamma < 3$,
in a \emph{supersonic regime}, quite different than the small-data regime we treat here. 
For the solutions in \cite{gCfeKyHyS2025}, both characteristic directions are outgoing to the future,
e.g., the radial variable $r$ increases along the future-directed integral curves of both characteristic vectorfields. Consequently, for $t \geq 0$,
no sound waves travel inwards towards the axis of symmetry $\lbrace r = 0 \rbrace$, 
and the authors were able to exploit various ``domain-invariance'' properties\footnote{This means that there are subsets of solution-space such that the solution remains in them for all future-times if it starts in them.}
to globally control the solution. In contrast, our small-data solutions feature inwards and outwards moving sound waves, for positive and negative times, 
and for the general equations of state that we treat, 
our only means to control the solution is dispersive estimates and the sign of the radiation field.

Finally, we highlight Wang's paper \cite{qW2025}, in which she proved exterior future-global existence 
for irrotational compressible fluids without symmetry,
for a large class of ``signed'' initial data 
(a ``future-rarefaction'' assumption was made on a term involving the density and its future-ingoing null derivative).
The data were non-trivial only in an annulus and hence domain of dependence considerations imply that 
the solution is completely trivial in the exterior of the flat sound cone emanating from the boundary of the annulus. 
Because of the sign assumption, the data have future-rarefactive properties, but are not required to be small. 
The author proved that rarefactive-type solutions exist for non-negative times
in the non-trivial region between an future-outgoing inner sound cone and the flat future-outgoing outer sound cone.
The dynamics were not studied in the interior region or past-direction, but it is likely that for at least some of the subsets of data, 
shock singularities will develop to the past, due to the signed data causing ``past-compression.''
It is also interesting to note that in the exterior vacuum region, 
the solutions in \cite{qW2025} are \emph{unstable} under perturbations of the data 
due to small-data shock formation, 
as was shown by Christodoulou--Miao in \cite{dCsM2014}, using the techniques that Christodoulou pioneered in his monograph \cite{dC2007} on the relativistic Euler equations.

As we have mentioned, the main novelty of Theorem~\ref{T:ROUGHGLOBALEXISTENCETHEOREM} is that the compressible Euler equations 
\eqref{E:INTROVELOCITYRADIALEQUATION}--\eqref{E:INTRODENSITYRADIALEQUATION}
\emph{fail to satisfy the null condition and the weak null condition} (aside from the Chaplygin gas equation of state)
and hence there are resonant interactions. Resonant terms appear, for example, in the evolution equation
\eqref{E:RIGHTRIEMANNEQUATIONINTERIORREGION} for the $\uLunit$ derivative of the spherical Riemann invariant $\RRiemann$,
where $\uLunit = \partial_t + (v^r - \Speed) \partial_r$ is an ingoing null vectorfield and $\Lunit = \partial_t + (v^r + \Speed) \partial_r$ is the complementary outgoing null vectorfield; see Def.\,\ref{D:NULLVECTORFIELDS}.
The equation takes\footnote{Throughout the paper, if $\mathbf{X}$ is a vectorfield and $f$ is a scalar function, 
then $\mathbf{X} f := \mathbf{X}^{\alpha} \partial_{\alpha} f$ denotes the derivative of $f$ in the direction $\mathbf{X}$. \label{FN:VECTORFIELDDERIVATIVESOFSCALARS}} 
the schematic form: 
\begin{align} \label{E:INTROSCHEMATICRICCATIFLUIDEQUATION}
	\Lunit (r \uLunit \RRiemann) 
	& = - r [\uLunit,\Lunit]\RRiemann + \cdots, 
\end{align}	
and aside from the case of the Chaplygin gas equation of state, the commutator term 
$- r [\uLunit,\Lunit]\RRiemann$ on RHS~\eqref{E:INTROSCHEMATICRICCATIFLUIDEQUATION}
contains a dangerous, resonant Riccati-type term proportional to $r (\uLunit \RRiemann)^2$. 

For $3D$ wave equations that \emph{do} satisfy the null condition, 
small-data global existence in $3D$ is a classic result proved independently by Klainerman and Christodoulou
\cites{sK1984,dC1986a}. Those results were immensely influential and have since been generalized in an incredible number of directions, far too many for us to exhaustively discuss here.
A notable extension of them is Lindblad--Rodnianski's proof of the stability of Minkowski spacetime in wave coordinates \cite{hLiR2010},
which remarkably simplified the original proof \cite{dCsK1993}. In \cite{hLiR2010}, the authors proved small-data global existence for a system of PDEs without the null condition,
but rather with a generalization of it coined the ``weak null condition'' \cite{hLiR2003}. 
The weak null condition applies to systems of hyperbolic PDEs in which dangerous quadratic nonlinearities are present, 
but they are \emph{not} Riccati-type self-interactions. That is, dangerous derivative quadratic terms are allowed to appear, but \emph{never} on the right-hand side of an equation 
whose principal left-hand part features the same unknown. In contrast, equation \eqref{E:INTROSCHEMATICRICCATIFLUIDEQUATION} fails to satisfy even the weak null condition.

We highlight that like the data in our main results, 
the initial data in \cite{hLiR2010} featured a $1/r$ tail, a Schwarzschild tail, which by the positive mass theorem had to be present.
Before \cite{hLiR2010}, almost all global existence results for nonlinear wave equations concerned compactly supported data or data that decayed faster than $1/r$.
The fact that the authors were able to prove global existence in the presence of the weak null condition with a $1/r$ tail helped motivate us to 
investigate how the tail affects equations for which even the weak null condition fails. It is interesting to note that in \cite{hLiR2010},
the $1/r$ tail was a ``nuisance'' in the PDE estimates that required careful care to handle. In contrast, in this paper, the $1/r$ tail yields monotonic terms that allow us to understand the global structure of the classical solution.

Generally throughout the literature, the presence of Riccati terms has been viewed as  
a \emph{threat} to small-data global existence. The reason is simply that the ODE $\dot{y} = y^2, y(0) = y_0$, features blowup at the time
$T_{(\textsf{Singular})} = \frac{1}{y_0}$. A shortcoming of this simple ODE is that it does not capture $3D$ wave dispersion, since $3D$ linear waves decay at the rate $1/(1 + |t|)$.
A standard ODE model for $3D$ nonlinear wave equations that fail to satisfy the null condition is the following nonlinear ODE, which mimics the wave decay in $3D$:
$\dot{y} = \frac{y^2}{1+|t|}, y(0) = y_0$. The term $\frac{1}{1+|t|}$, though decaying, is not integrable in time, and hence this model 
still exhibits blowup, for example at the positive time $T_{(\textsf{Singular})} = \exp\left(\frac{1}{y_0} \right) - 1$ when $y_0 > 0$. 
In particular, \emph{all} solutions 
aside from $y \equiv 0$ blow up, and the expected future-lifespan is approximately the exponential of $1$ over the size of the data. 
In \cite{fJsK1984}, John--Klainerman showed that even if the null condition fails, the lifespan of a nonlinear $3D$ wave equation is at least 
$\exp\left(\frac{1}{C \datasize} \right)$, where $\datasize$ is the size of the data and $C > 0$ is a constant. 
This shows that the ODE model $\dot{y} = \frac{y^2}{1+|t|}, y(0) = y_0$ from above gives a reasonably good lower bound for the 
wave equation solution's classical lifespan, at least when the wave data are smooth and compactly supported.

None of the above discussion addressed whether singularities actually form in the nonlinear PDE solutions when the null condition and weak null condition both fail.
Fritz John's work \cite{fJ1981} showed via a proof by contradiction that for a large class of $3D$ nonlinear wave equations that fail to satisfy the null condition
and weak null condition, solutions launched by small, smooth, compactly supported data \emph{always} lead to finite-time singularity formation. 
Sideris later proved a related result for the $3D$ compressible Euler equations in his influential paper \cite{tS1985}, in which he
considered smooth initial data that are compactly supported perturbations of trivial data with vanishing velocity 
and entropy and constant positive density. He showed via a proof by contradiction that when certain integrals of such data have a sign (roughly, a compressive sign),
the solution must blow up to the future in finite time. His work did not reveal the nature of the singularity or
address what happens for negative times (or equivalently, when the integrals have the opposite sign).
Moreover, in \cite{fJ1985}, John showed that for $3D$ spherically symmetric quasilinear wave equations of type $\partial_t^2 \phi = \Speed^2(\partial_t\phi) \Delta \phi$,
with $\Speed(0) > 0$ and $\Speed'(0) \neq 0$, \emph{all} small, compactly supported, non-trivial data lead to finite-time blowup in the solution
at a time approximately equal to $\exp\left(\frac{1}{\widetilde{C} \datasize} \right)$, where $\widetilde{C} > 0$ is a constant. 
Our small-data global existence theorem could also be proved for the equation $\partial_t^2 \phi = \Speed^2(\partial_t\phi) \Delta \phi$ for initial data with an appropriate $1/r$ tail.
Similarly, in the recent work \cite{klainerman2026inevitable}, the authors prove that given arbitrary non-trivial, smooth 
compactly supported initial data for the $3D$ irrotational and isentropic compressible Euler equations, 
if the data are \emph{re-scaled} by a sufficiently small constant, 
then the resulting solution forms a shock in finite time.
Therefore, a key takeaway of Theorem~\ref{T:ROUGHGLOBALEXISTENCETHEOREM} is that without the compact support assumption on the data, the results of \cites{fJ1985,klainerman2026inevitable} are not true in general.

Roughly, in Theorem~\ref{T:ROUGHGLOBALEXISTENCETHEOREM}, our integrals\footnote{We never need to consider these integrals in our analysis.} 
have the opposite sign compared to Sideris' paper \cite{tS1985},
which causes the resonant Riccati-type terms to have a decay-inducing effect, rather than causing blowup.
There does not seem to be a simple way to capture this with an ODE model, since all nontrivial data for
$\dot{y} = \frac{y^2}{1+|t|}, y(0) = y_0$ lead to blowup in the future or past, depending on the sign of $y_0$. 
We will now briefly describe the global existence dynamics for the compressible fluids that we study in Theorem~\ref{T:ROUGHGLOBALEXISTENCETHEOREM},
focusing for simplicity on the spherical velocity component $v^r$; 
see Sect.\,\ref{SS:BREIFOVERVIEWOFANALYSIS} for a more detailed overview of the analysis.
Standard linear wave theory\footnote{Linearizing \eqref{E:INTROVELOCITYRADIALEQUATION}--\eqref{E:INTRODENSITYRADIALEQUATION} leads to linear wave equations for $\varrho$ and $v^r$.} in $3D$ 
predicts that the Minkowskian ingoing null derivative of $v^r$, namely $|(\partial_t - \partial_r) v^r|$,
should decay to the future like $\frac{1}{(1 + t + |u_{\textsf{(flat)}}|)(1 + |u_{\textsf{(flat)}}|)}$, where $u_{\textsf{(flat)}} := t-r$
is a linear outgoing eikonal function (i.e., its level sets are characteristic for the linear wave equation).
In contrast, in the context of Theorem~\ref{T:ROUGHGLOBALEXISTENCETHEOREM}, 
we prove that $|\uLunit v^r|$ decays to the future like $\frac{1}{\upmu} \frac{1}{(1 + t + |u|)(1 + |u)}$,
where $\uLunit = \partial_t + (v^r - \Speed) \partial_r$ is the ingoing null vectorfield,
$u$ is an outgoing eikonal function (it satisfies $\Lunit u = 0$, where $\Lunit = \partial_t + (v^r + \Speed) \partial_r$),
and $\upmu$ \emph{grows} logarithmically in time, to the future; see 
\eqref{E:MUWEIGHTEDVECTORFIELDS},
\eqref{E:GLOBALEXISTENCEFARFROMTIMEAXISPOINTWISEAPRIORIMULBARRPLUS},
\eqref{E:GLOBALEXISTENCEFARFROMTIMEAXISPOINTWISEAPRIORILBARRMINUS},
and \eqref{E:GLOBALEXISTENCEMUPOINTWISESTIMATE}.
In particular, the time decay rate of $\uLunit v^r$ is \emph{logarithmically faster} compared to the linear wave equation.
Moreover, due to our very small velocity data,
our solutions behave very similarly as $t \downarrow - \infty$, (see Remark~\ref{R:DISCRETESYMMETRY}) 
and in particular, analogous time-decay results hold for $t < 0$.
The factor of $\frac{1}{\upmu}$ has a crucial geometric interpretation: 
it is the density of the $\Lunit$-characteristics, and the logarithmic-in-time growth of $\upmu$ is
a quantified manifestation of rarefaction, i.e., of the separation of the outgoing characteristics.
The factor of $\upmu$ contains the ``information'' of the Riccati terms in the flow, and the separation of the characteristics is tantamount to the Riccati terms having a \emph{good sign}.

\subsection{A brief overview of the analysis}
\label{SS:BREIFOVERVIEWOFANALYSIS}
In this section, we provide a brief overview of the analysis we use to prove our main results.
To capture wave dispersion (i.e., the decay of nonlinear fluid sound waves), we find it convenient to use spherical Riemann invariants 
$\mathcal{R}_{\pm}$, 
defined below in Def.\,\ref{D:SPHERICALRIEMANNINVARIANTS}. 
The density $\varrho$ and spherical velocity component $v^r$
can be written as smooth functions of $\mathcal{R}_{\pm}$, so to control the fluid, it suffices to control $\mathcal{R}_{\pm}$.
The main advantage of the Riemann invariants is that they satisfy a quasilinear, diagonal transport system, featuring factors of $1/r$ that generate \emph{wave dispersion};
see \eqref{E:OUTGOINGRIEMANNINVARIANTEVOLUTION}--\eqref{E:INGOINGRIEMANNINVARIANTEVOLUTION}.
By our definitions, $\RRiemann = \LRiemann \equiv 0$ corresponds to a solution with constant density $\overline{\varrho} > 0$ and vanishing velocity.
The solutions in our main theorems are perturbations of such a trivial solution. Our perturbations are such that the initial velocity is 
\emph{very small}, with a tiny amplitude and decaying at least like $1/r^2$ at spatial infinity (vanishing velocity data is allowed); 
see Remark~\ref{R:TINYVELOCITYDATA}. 
Relatedly, in both of our theorems, 
the past- and future-dynamics are almost mirror images of each other. We further this explain fact in the following remark.

\begin{remark}[It suffices to understand the solution for $t \geq 0$]
\label{R:DISCRETESYMMETRY}
The compressible Euler equations
\eqref{E:VELOCITYRADIALEQUATION}--\eqref{E:DENSITYRADIALEQUATION} are invariant under the discrete symmetry
$t \rightarrow - t$ and $v^r \rightarrow - v^r$, which leaves the density initial data invariant and
preserves our close-to-the-profile data-assumptions 
\eqref{E:RPLUSDATAISCLOSETOBACKGROUNDATTIME0}--\eqref{E:RMINUSDATAISCLOSETOBACKGROUNDATTIME0} in the first theorem
as well as our close-to-the-profile data-assumptions \eqref{E:GLOBALEXISTENCERPLUSDATAISCLOSETOBACKGROUNDATTIME0}--\eqref{E:GLOBAEXISTENCERMINUSDATAISCLOSETOBACKGROUNDATTIME0} in the second theorem, 
and hence preserves the velocity-smallness assumptions clarified in Remark~\ref{R:TINYVELOCITYDATA}.
In particular, the past-dynamics are nearly a mirror image of the future-dynamics.
Moreover, this discrete symmetry interchanges $\RRiemann$ and $\LRiemann$ as well as the null vectorfields
$\Lunit$ and $\uLunit$ from Def.\,\ref{D:NULLVECTORFIELDS}.
\end{remark}

The following point will guide all the subsequent discussion:
\begin{quote}
	In the small-data regime for the nonlinear equations, the most dangerous term to the future is the future-ingoing null derivative $\uLunit \RRiemann$, where
	$\uLunit = \partial_t + (v^r - \Speed) \partial_r$. All other derivative-terms in the nonlinear equations, including $\Lunit \RRiemann$,
	$\uLunit \LRiemann$, and $\Lunit \LRiemann$,
	where $\Lunit := \partial_t + (v^r + \Speed) \partial_r$,
	decay much faster to the future (at an integrable-in-time rate)
	and are harmless in the context of global existence. They are also harmless in the context of shock formation in that they don't blow up.
	Similarly, the most dangerous term to the past is the past-ingoing null derivative
	$\Lunit \LRiemann$. 
\end{quote}

\subsubsection{Motivation via the linearized system}
\label{SSS:MOTIVATIONVIALINEARIZEDSYSTEM}
To help motivate our results as well as the forthcoming long-time nonlinear analysis, 
here we derive formulas for solutions to a linearized version of  \eqref{E:OUTGOINGRIEMANNINVARIANTEVOLUTION}--\eqref{E:INGOINGRIEMANNINVARIANTEVOLUTION}.
We denote the linear solution variables by $(\flatRRiemann,\flatLRiemann)$.
Using \eqref{E:SPEEDOFSOUNDEXPANSION}--\eqref{E:SPEEDOFSOUNDERRORFUNCTIONVANISHESATORIGIN},
it is straightforward to compute that the linearized (around the solution $(\RRiemann,\LRiemann) \equiv (0,0)$) 
system corresponding to \eqref{E:OUTGOINGRIEMANNINVARIANTEVOLUTION}--\eqref{E:INGOINGRIEMANNINVARIANTEVOLUTION} is:
\begin{subequations}
\begin{align}
	\flatLunit \flatRRiemann
	& = \frac{1}{r} (\flatLRiemann - \flatRRiemann),
		\label{E:LINEARIZEDOUTGOINGRIEMANNINVARIANTEVOLUTION} 
			\\
	\flatuLunit \flatLRiemann
	& = \frac{1}{r} (\flatLRiemann - \flatRRiemann),
	\label{E:LLINEARIZEDINGOINGRIEMANNINVARIANTEVOLUTION} 
\end{align}
\end{subequations}
where:
\begin{subequations}
\begin{align}
	\flatLunit 
	& = \partial_t + \partial_r,
			\label{E:FLATLUNIT} 
				\\
	\flatuLunit 
	& = \partial_t - \partial_r.
		\label{E:FLATULUNIT}
\end{align}

To motivate our main shock-formation/MGHD results, we choose the data for the linearized system to be the following ``background profile'':
	\begin{align} \label{E:BACKGROUNDDATAFORLINEARSYSYSTEM}
	\flatRRiemann(0,r)
	& = \flatLRiemann(0,r)
	:= - \frac{\datasize \arctan(r)}{r},
\end{align}
\end{subequations}
where $\datasize > 0$.
In Theorem~\ref{T:MAINMGHDEXISTENCETHEOREM}, the data will be a perturbation of the background profile in \eqref{E:BACKGROUNDDATAFORLINEARSYSYSTEM}.
In Appendix~\ref{A:EXTENDRESULTSTOOTHERPROFILES}, we outline how to treat a much larger class of data, based on perturbing other profiles.
The main advantage of the profiles is that they lead to radiation fields with a global sign; see below.

It is straightforward to check\footnote{One can check that $\flatRRiemann$ and $\flatLRiemann$ satisfy $\flatLunit( r\flatuLunit \flatRRiemann - 2 \flatRRiemann) = \flatuLunit(r\flatLunit\flatLRiemann + 2 \flatLRiemann) = 0$. From this, one easily obtains \eqref{E:LINEARIZEDSOLUTIONSPECIALRPLUSCOMBINATIONFORMULA}--\eqref{E:LINEARIZEDSOLUTIONSPECIALRMINUSCOMBINATIONFORMULA}. Equations \eqref{E:LINEARRPLUSSOLUTIONFORMULA}--\eqref{E:LINEARRMINUSSOLUTIONFORMULA} then follow from solving \eqref{E:LINEARIZEDSOLUTIONSPECIALRPLUSCOMBINATIONFORMULA}--\eqref{E:LINEARIZEDSOLUTIONSPECIALRMINUSCOMBINATIONFORMULA} by the method of characteristics.} 
that the linear solution launched by the data \eqref{E:BACKGROUNDDATAFORLINEARSYSYSTEM} is:
\begin{subequations}
\begin{align} \label{E:LINEARRPLUSSOLUTIONFORMULA}
	\flatRRiemann 
	& =
	\frac{\datasize}{2}
				\frac{1}{r^2}
				\left\lbrace
				 (t + r) \arctan (t - r)
					- 
					(t + r) \arctan (t + r)
				-
				\frac{1}{2} 
				\ln \left(\frac{1 + (t - r)^2}{1 + (t+ r)^2} \right)
				\right\rbrace,
				\\
	\flatLRiemann  \label{E:LINEARRMINUSSOLUTIONFORMULA}
	& =
	\frac{\datasize}{2}
				\frac{1}{r^2}
				\left\lbrace
					(t - r) \arctan ((t + r))
					- 
					(t - r) \arctan ((t - r))
					-
					\frac{1}{2} 
				\ln \left(\frac{1 + (t + r)^2}{1 + (t - r)^2} \right)
				\right\rbrace,
\end{align}				
\end{subequations}
and that the linear solution satisfies the following identities:
\begin{subequations}
\begin{align} \label{E:LINEARIZEDSOLUTIONSPECIALRPLUSCOMBINATIONFORMULA} 
	r
	\flatuLunit
	\flatRRiemann
	- 
	2 \flatRRiemann
	& =  
			\frac{2 \datasize}{1 + (t-r)^2},
				\\
	r
	\flatLunit
	\flatLRiemann
	+ 
	2 \flatLRiemann
	& =  
		\frac{- 2 \datasize}{1 + (t+r)^2},
		\label{E:LINEARIZEDSOLUTIONSPECIALRMINUSCOMBINATIONFORMULA} 
\end{align}
\end{subequations}

\begin{align} \label{E:BACKGROUNDTIMEDERIVATIVEDATAFORLINEARSYSYSTEM}
	\partial_t \flatRRiemann(0,r)
	& = - \partial_t \flatLRiemann(0,r)
	= 
	\datasize
	\partial_r \left(\frac{\arctan(r)}{r} \right).
\end{align}

From \eqref{E:LINEARRPLUSSOLUTIONFORMULA}--\eqref{E:LINEARRMINUSSOLUTIONFORMULA}, it is straightforward to check that
$\flatRRiemann$ and $\flatLRiemann$ decay at the expected rate of a linear wave, even though their data has a $1/r$ tail: 
$|\flatRRiemann| + |\flatLRiemann| \lesssim \frac{\datasize}{1 + t + r}$.
Moreover, the crucial identity \eqref{E:LINEARIZEDSOLUTIONSPECIALRPLUSCOMBINATIONFORMULA}
captures the \emph{future radiation field} of the linear solution.
In particular, since $\datasize > 0$, 
$r
	\flatuLunit
	\flatRRiemann
	- 
	2 \flatRRiemann$ 
has a positive sign, everywhere in spacetime. 
Because of the decay of 
$\flatRRiemann$ and $\flatLRiemann$, the main part of this expression is
$r
	\flatuLunit
	\flatRRiemann$, 
i.e., for sufficiently large times, 
$r
	\flatuLunit
	\flatRRiemann$ has a positive sign.
In the context of Theorem~\ref{T:ROUGHMGHDTHEOREM},
a similar radiation field estimate drives the resonant ``null-condition-failing'' terms and causes the blowup of $\uLunit \RRiemann$, all the way out to spatial infinity,
along the future singular boundary curve $\futuresinghyp$ in Fig.\,\ref{F:FULLMGHD}.
Similarly, \eqref{E:LINEARIZEDSOLUTIONSPECIALRMINUSCOMBINATIONFORMULA} captures the \emph{past radiation field} of the linear solution,
and in the context of Theorem~\ref{T:ROUGHMGHDTHEOREM},
a similar radiation field estimate drives the resonant ``null-condition-failing'' terms and causes the blowup of $\Lunit \LRiemann$, all the way out to spatial infinity,
along the past singular boundary curve $\pastsinghyp$ in Fig.\,\ref{F:FULLMGHD}.

To motivate our global existence results, we choose the data for the linearized system to have the opposite sign compared to \eqref{E:BACKGROUNDDATAFORLINEARSYSYSTEM}:
\begin{subequations}
\begin{align} \label{E:INTROGLOBALEXISTENCEBACKGROUNDDATAFORLINEARSYSYSTEM}
	\flatRRiemann(0,r)
	& = \flatLRiemann(0,r)
	:= \frac{\datasize \arctan(r)}{r},
\end{align}
\end{subequations}
where $\datasize > 0$.
The corresponding linear solution has the opposite sign compared to \eqref{E:LINEARRPLUSSOLUTIONFORMULA}--\eqref{E:LINEARRMINUSSOLUTIONFORMULA}:
\begin{subequations}
\begin{align} \label{E:INTROGLOBALEXISTENCELINEARRPLUSSOLUTIONFORMULA}
	\flatRRiemann 
	& =
	-
	\frac{\datasize}{2}
				\frac{1}{r^2}
				\left\lbrace
				 (t + r) \arctan (t - r)
					- 
					(t + r) \arctan (t + r)
				-
				\frac{1}{2} 
				\ln \left(\frac{1 + (t - r)^2}{1 + (t+ r)^2} \right)
				\right\rbrace,
				\\
\flatLRiemann  \label{E:INTROGLOBALEXISTENCELINEARRMINUSSOLUTIONFORMULA}
	& =
	-
	\frac{\datasize}{2}
				\frac{1}{r^2}
				\left\lbrace
					(t - r) \arctan ((t + r))
					- 
					(t - r) \arctan ((t - r))
					-
					\frac{1}{2} 
				\ln \left(\frac{1 + (t + r)^2}{1 + (t - r)^2} \right)
				\right\rbrace,
\end{align}				
\end{subequations}
and the past-and future- radiation fields are as follows:
\begin{subequations}
\begin{align} \label{E:INTROGLOBALEXISTENCELINEARIZEDSOLUTIONSPECIALRPLUSCOMBINATIONFORMULA} 
	r
	\flatuLunit
	\flatRRiemann
	- 
	2 \flatRRiemann
	& =  
			\frac{- 2 \datasize}{1 + (t-r)^2},
				\\
	r
	\flatLunit
	\flatLRiemann
	+ 
	2 \flatLRiemann
	& =  
		\frac{2 \datasize}{1 + (t+r)^2}.
		\label{E:INTROGLOBALEXISTENCELINEARIZEDSOLUTIONSPECIALRMINUSCOMBINATIONFORMULA} 
\end{align}
\end{subequations}
Once again, the key point is the \emph{sign} of the future-radiation field \eqref{E:INTROGLOBALEXISTENCELINEARIZEDSOLUTIONSPECIALRPLUSCOMBINATIONFORMULA} and 
the past-radiation field \eqref{E:INTROGLOBALEXISTENCELINEARIZEDSOLUTIONSPECIALRMINUSCOMBINATIONFORMULA}.
For example, up to the decaying term $\flatRRiemann$, 
\eqref{E:INTROGLOBALEXISTENCELINEARIZEDSOLUTIONSPECIALRPLUSCOMBINATIONFORMULA} shows that
$r \flatuLunit \flatRRiemann$ is everywhere \emph{negative}. 
Hence, in the nonlinear problem, the corresponding Riccati-type equation 
$\Lunit (r \uLunit \RRiemann)  \sim r (\uLunit \RRiemann)^2 + \cdots$ actually \emph{enhances} the future-decay of 
$|\uLunit \RRiemann|$.
A simple ODE model for this effect (with $\uLunit \RRiemann$ in the role of $y$) is the ODE $\dot{y} = \frac{y^2}{1+t}$ with \emph{negative} initial data $y(0) = y_0 < 0$;
the corresponding solution $y$ exists globally to the future and in fact decays in amplitude towards $0$.
A similar effect happens to the past due to the sign of the radiation field \eqref{E:INTROGLOBALEXISTENCELINEARIZEDSOLUTIONSPECIALRMINUSCOMBINATIONFORMULA},
which, in the nonlinear Riccati equation for $r \Lunit \RRiemann$, 
enhances the decay of the dangerous term $|\Lunit \RRiemann|$ as $t \downarrow - \infty$.

\subsubsection{Nonlinear geometric optics}
\label{SSS:INTRONONLINEARGEOMETRICOPTICS}
In proving Theorems~\ref{T:ROUGHMGHDTHEOREM} and \ref{T:ROUGHGLOBALEXISTENCETHEOREM}, we fundamentally rely on nonlinear geometric optics,
which we implement with the help of an outgoing eikonal function $u$, that is, for $t \geq 0$, a solution $u$ to:
\begin{align} \label{E:INTROEIKONAL}
	\Lunit u 
	& = 0,
\end{align}
where $\Lunit = \partial_t + (v^r + \Speed) \partial_r$ is the outgoing null vectorfield (see Def.\,\ref{D:NULLVECTORFIELDS}).
For $t \geq 0$, we then study the fluid relative to the \emph{geometric coordinate system} $(t,u)$,
whose corresponding partial derivative vectorfields we denote by $\lbrace \frac{\partial}{\partial t}, \frac{\partial}{\partial u} \rbrace$.
The level sets of $u$ are outgoing characteristic curve portions. 
Note that since $\Lunit t = 1$ and $\Lunit u = 0$, we have $\Lunit = \frac{\partial}{\partial t}$.
In both theorems, we avoid constructing $u$ near $r= 0$ because coordinate degeneracies in the equations near $r=0$ make it difficult to control $u$ there. That is not important in the sense that the fluid exhibits a lot of time-decay near the symmetry axis $\lbrace r=0 \rbrace$, and hence
the fluid is easy to control there directly with respect to the $(t,r)$ coordinates. We choose different data for $u$ in the two theorems. The details are not important for the present discussion, but in each theorem, $u$ behaves like $t - r$ plus correction terms that can grow logarithmically in time.

The key advantage of the geometric coordinates is the following collection of ideas, which has roots from many places, 
notably papers of Alinhac \cites{sA1999a,sA1999b}, Christodoulou \cites{dC2007}, and second author \cites{jLjS2020a,jS2016b,jS2019c}:
\begin{quote}
	Relative to the $(t,u)$ coordinates, $\mathcal{R}_{\pm}$ satisfy a system of wave equations in which the nonlinear terms \emph{do} 
	satisfy the null condition, i.e., no resonant terms are present in the equations; 
	see Prop.\,\ref{P:GEOMETRICWAVELIKEEVOLUTIONEQUATIONS}. Consequently, in $(t,u)$ coordinates, 
	one can study the dynamics using techniques developed for more standard long-time existence problems.
	Moreover, if a shock forms in the solution with respect to the original $(t,r)$ coordinate system,
	this must be due to a degeneracy between the two coordinate systems.
\end{quote}
The surprising fact that the null condition holds in $(t,u)$ coordinates
does \emph{not} contradict the fact that 
in the original $(t,r)$ coordinate system, the null condition fails, i.e., that 
there are resonant Riccati terms in the original equation $\Lunit (r \uLunit \RRiemann)  \sim r (\uLunit \RRiemann)^2 + \cdots$.
We will now explain why there is no contradiction.

When one changes variables from $(t,r)$ to $(t,u)$ coordinates, this introduces a change of variables factor $\frac{\partial}{\partial u} r$. One can show that schematically, we have:
\begin{align} \label{E:INTROSCHEMATICIDENTITYFORUDERIVATIVEOFR}
	\frac{\partial}{\partial u} r \sim - \upmu,
\end{align}
where $\upmu > 0$, the \emph{inverse foliation density of the characteristics},
satisfies:
\begin{align} \label{E:INTROMUIDENTITY}
	\upmu
	& = \frac{2}{\uLunit u}.
\end{align}
Geometrically,
$\upmu = 0$ means that the level sets of $u$ have become infinitely dense, signifying shock formation, while
$\upmu$ growing large signifies the separation of the characteristics (rarefaction).
When one changes variables to geometric coordinates in the wave equation for $\RRiemann$,
this generates a change-of-variables term of the form $\Lunit \upmu \cdot \uLunit \RRiemann$.
The key point is that $\upmu$ satisfies an evolution equation of the following schematic form (see \eqref{E:MUEVOLUTION}):
\begin{align} \label{E:INTROMUEVOLUTION}
	\Lunit \upmu
	& \sim - \upmu \uLunit \RRiemann + \cdots,
\end{align}
and it turns out that the change-of-variables term $\Lunit \upmu \cdot \uLunit \RRiemann$
leads to the \emph{complete cancellation of the dangerous Riccati term} $r (\uLunit \RRiemann)^2 + \cdots$
in the wave equation for $\RRiemann$; see \eqref{E:REVAMPEDRPLUSWAVEEQUATION}.
The upshot is the following:
\begin{quote}
	In all our theorems, relative to the $(t,u)$ coordinates, the nonlinear terms satisfy the null condition, and the corresponding nonlinear solution
	behaves like the linear solutions from Sect.\,\ref{SSS:MOTIVATIONVIALINEARIZEDSYSTEM},
	where for positive times, the role of $t-r$ is played by $u$.
\end{quote}
Hence, relative to the $(t,u)$ coordinates, we can approach the PDE analysis like a small-data global existence problem for a wave system that satisfies the null condition, the new difficulty being that the data have a $1/r$ tail, which makes some of the decay-in-$u$ rates worse compared to the case of compactly supported data. 
We must also control $\upmu$ by carefully studying the transport equation \eqref{E:INTROMUEVOLUTION}. 

In the shock-formation/MGHD problem, to estimate $\upmu$, we first use the wave equation for $\RRiemann$ and our data assumptions to
prove a bound of the form $\upmu \uLunit \RRiemann \sim \frac{f(u)}{t}$ for some 
smooth function $f > 0$ (see \eqref{E:EXTERIORREGIONPOINTWISEESTIMATEFORMULBARRPLUS} for the details);
this bound relies on sign of the radiation field, as discussed in Sect.\,\ref{SSS:MOTIVATIONVIALINEARIZEDSYSTEM}. 
We then insert this bound into equation \eqref{E:INTROMUEVOLUTION} and integrate in time in the $(t,u)$ coordinate system 
(recalling that $\Lunit = \frac{\partial}{\partial t}$), thereby concluding (thanks to the non-integrability in $t$ of $\frac{1}{t}$ near $\infty$) 
that $\upmu$ vanishes in finite time along the singular boundary curve $\futuresinghyp$ in Fig.\,\ref{F:MGHD}.
Note that since $\upmu \uLunit \RRiemann \sim \frac{f(u)}{t}$, the vanishing of $\upmu$ causes the blowup of $\uLunit \RRiemann$, i.e., blowup of a first derivative of $\RRiemann$ in the $(t,r)$ differential structure. On the other hand, simple calculations yield that in $(t,u)$ coordinates,
we have $\upmu \uLunit = \upmu \frac{\partial}{\partial t} + 2 \frac{\partial}{\partial u}$ (see Lemma~\ref{L:VECTORFIELDSINTERMSOFGEOMETRICCOORDINATES}),
and hence $\upmu \frac{\partial}{\partial t} \RRiemann + 2 \frac{\partial}{\partial u} \RRiemann$ remains bounded. 
Since $\frac{\partial}{\partial t} \RRiemann$ also remains bounded (see \eqref{E:EXTERIORREGIONPOINTWISEESTIMATEFORRPLUSANDLDERIVATIVES}),
it follows that relative to the geometric coordinates, no first derivative of $\RRiemann$ blows up, even when $\upmu$ vanishes.
Our estimates (see Props.\,\ref{P:RIEMANNINVARIANTSAPRIORIEXTERIORREGIONESTIMATES} and \ref{P:APRIORIESTIMATESININTERIORREGION}) also show 
that \emph{no first derivative of $\LRiemann$ blows up to the future, even in $(t,r)$ coordinates}; see \eqref{E:OUTGOINGRIEMANNINVARIANTEVOLUTION}--\eqref{E:INGOINGRIEMANNINVARIANTEVOLUTION}, \eqref{E:EXTERIORREGIONPOINTWISEESTIMATEFORRPLUSANDLDERIVATIVES}, and \eqref{E:EXTERIORREGIONPOINTWISEESTIMATEFORRMINUSANDLDERIVATIVES}.

For the global existence problem, the analysis proceeds in a similar fashion. The main difference is that the radiation field has the opposite sign,
and thus the evolution equation \eqref{E:INTROMUEVOLUTION} leads to logarithmic-in-time growth of $\upmu$; see 
\eqref{E:GLOBALEXISTENCEMUPOINTWISESTIMATE}.
Hence, the outgoing characteristics separate, and no singularity ever forms.


\subsection{The spherically symmetric isentropic compressible Euler equations} \label{SS:SPHERICALLYSYMMETRICEULER}
In this section, we provide a more detailed introduction to the equations.
The isentropic Euler equations in three spatial dimensions can be expressed as a quasilinear hyperbolic system for the \emph{velocity} $v : \R^{1+3} \to \R^3$ and the \emph{density} $\varrho: \R^{1+3} \to \R_{\ge 0}$ of a compressible fluid. In the standard Cartesian differential structure of $\R^{1+3}$ with coordinates $(t,x) =(t,x^1,x^2,x^3)$, the equations of motion are:\footnote{On RHS~\eqref{E:INTROTRANSPORTDENSITY} and throughout, we are using Einstein's summation convention, i.e., repeated adjacent lowercase Latin indices are summed from $1$ to $3$. \label{FN:EINSTEINSUMMATION}}
\begin{subequations} \label{E:INTRO:COMPRESSIBLEEULER}
\begin{align}
	\Transport v^i 
	& = 
	- \frac{\partial_i p}{\varrho},
	&& (i=1,2,3),
	\label{E:INTROTRANSPORTVI}
		\\
		 \label{E:INTROTRANSPORTDENSITY}
	\Transport \varrho
	& = - \varrho \partial_a v^a,
	&&
\end{align}
\end{subequations}
where:
\begin{align}\label{E:MATERIALDERIVATIVEVECOTRFIELD}
\Transport 
& := \partial_t + v^a \partial_a
\end{align}
denotes the \emph{material derivative vectorfield} and $p$ denotes the \emph{pressure}. \eqref{E:INTRO:COMPRESSIBLEEULER} is an under-determined system of because there are not enough equations. 
To close the system, we assume that $p = p(\varrho)$ is given as a function of $\varrho$, which is called the equation of state. Our results hold 
\emph{for any $C^4$ equation of state except for that of the Chaplygin gas}, provided that the \emph{speed of sound} $\Speed$ is positive whenever $\varrho$ is positive: 
\begin{align} \label{E:INTRO:SPEEDOFSOUND}
	\Speed(\varrho) := \sqrt{\frac{d p}{d \varrho}(\varrho)} > 0, \qquad \qquad \text{when $\varrho > 0$}.
\end{align}
		
In this article, to avoid fluid-vacuum free boundaries,
we study only solutions with strictly positive density. 
More precisely, we will study perturbations of ``background solutions'' 
in which the velocity completely vanishes and
such that the density is everywhere equal to $\bar{\varrho}$, where $\bar{\varrho} > 0$ 
is a constant. By rescaling the time coordinate $t$ by a positive constant if necessary, we can assume the following 
normalization condition:
\begin{align} \label{E:BACKGROUNDSPEEDOFSOUNDISUNITYNORMALIZATION}
	\overline{\Speed}
	& := \Speed(\varrho=\bar{\varrho})
	= 1.
\end{align} 

We are concerned with \emph{spherically symmetric} flows, i.e., the spatial dependence of solutions to \eqref{E:INTROTRANSPORTVI}--\eqref{E:INTROTRANSPORTDENSITY} is restricted to $r = \sqrt{(x^1)^2+(x^2)^2+(x^3)^2}$ and the velocity field is purely radial: 
\begin{align} \label{E:PURELYRADIALVELOCITY}
	v^i(t,x) = v^r(t,r) \frac{x^i}{r}.
\end{align}

Under these assumptions, the compressible Euler equations \eqref{E:INTROTRANSPORTVI}--\eqref{E:INTROTRANSPORTDENSITY} are equivalent to the following PDE system
in $(\varrho,v^r)$:
\begin{subequations}
        \begin{align}
            \partial_t v^r + v^r\partial_r v^r = - \frac{\Speed^2}{\varrho} \partial_r \varrho, \label{E:VELOCITYRADIALEQUATION} 
							\\
            \partial_t \varrho + \partial_r (\varrho v^r) = -2 \frac{\varrho v^r}{r}. \label{E:DENSITYRADIALEQUATION}
        \end{align}
\end{subequations}

For smooth spherically symmetric solutions, $\frac{v^r}{r}$ and $\varrho$ extend to functions of $r \in - (\infty,\infty)$ that are smooth and even.
More precisely, for any non-negative integer $k$ such that $v^r \in C^{2k}$ and $\varrho \in C^{2k + 1}$, 
we have the following ``boundary conditions,'' which we assume throughout the article:
\begin{subequations}
\begin{align}
	\partial_r^{2k} v^r \restriction_{\lbrace r=0 \rbrace} & =0,
		\label{E:EVENDERIVATIVESOFVRVANISHATTHEORIGIN} \\
	\partial_r^{2k+1} \varrho \restriction_{\lbrace r=0 \rbrace} & = 0.
	\label{E:ODDDERIVATIVESOFDENSITYVANISHATTHEORIGIN}
\end{align}
\end{subequations}

\subsection*{Acknowledgements}

L. Abbrescia gratefully acknowledges support from a Travel Support for Mathematicians grant from the Simons Foundation. J. Speck gratefully acknowledges support from NSF grants DMS-2349575 and DMS-2054184. D. Yu gratefully acknowledges support from  a VandyGRAF Fellowship from Vanderbilt University. We are grateful to The Fields Institute for Research in Mathematical Sciences for funding the conference ``Singularity Formation and Propagation in Gas Dynamics'', where some of the ideas in this paper were developed.
    
\section{Riemann invariants and the characteristic vectorfields}
\label{S:RIEMANNINVARIANTSANDCHARACTERISTICVECTORFIELDS}
In this section, we introduce the spherical Riemann invariants $\mathcal{R}_{\pm}$ as well as the corresponding null (i.e., characteristic) vectorfields
$\Lunit$ and $\uLunit$. We also derive, 
as a consequence of \eqref{E:VELOCITYRADIALEQUATION}--\eqref{E:DENSITYRADIALEQUATION},
various evolution equations satisfied by $\mathcal{R}_{\pm}$.
All of our forthcoming PDE estimates are for $\mathcal{R}_{\pm}$.

\subsection{Definitions and identities}
\label{SS:RIEMANNINVARIANTSANDNULLVECTORFIELDSDEFSANDIDS}

\begin{definition}[Spherical Riemann invariants]
\label{D:SPHERICALRIEMANNINVARIANTS}
Recall that $\overline{\varrho} > 0$ denotes a fixed positive ``background'' density, an that we will be studying solutions whose density is close to
$\overline{\varrho}$. Given any density $\varrho$ and any radial velocity component $v^r$, we define the corresponding \emph{spherical Riemann invariants}\footnote{Some papers use a slightly different definition for the Riemann invariants, such as $\Riemanninvariant_{(\pm)} = v^r \pm F(\varrho)$. We chose \eqref{E:SPHERICALRIEMANNINVARIANTS} merely for convenience, in particular to simplify the numerical coefficients in various equations.} to be: 
\begin{align} \label{E:SPHERICALRIEMANNINVARIANTS}
\Riemanninvariant_{(\pm)} 
& := 
\frac{1}{2}
\left\lbrace
\Riemannfunction(\varrho)
\pm 
v^r
\right\rbrace,
& 
\text{where \, \, } 
\Riemannfunction(\varrho)
&
:= 
\int_{\overline{\varrho}}^{\varrho} 
	\frac{\Speed(\varrho')}{\varrho'}
\, \rmd \varrho'.
\end{align}
\end{definition}
Note that given $\RRiemann$ and $\LRiemann$, we can uniquely recover the original fluid variables via the following formulas:
\begin{align} \label{E:ORIGINALVARIABLESINTERMSOFRIEMANNINVARIANTS}
	v^r
	& = 
		\RRiemann
		-
		\LRiemann,	
	& 
	\varrho
	& = 
	\InverseRiemannfunction 
	\circ
	\left\lbrace
	  \RRiemann
		+
		\LRiemann
	\right\rbrace,
\end{align}
where $\InverseRiemannfunction$ is the inverse function of $\Riemannfunction$, i.e., 
\begin{align} \label{E:INVERSERIEMANNFUNCTION}
	\InverseRiemannfunction 
	& := \Riemannfunction^{-1}.
\end{align}
In particular, $\Riemannfunction(\overline{\varrho}) = 0$, and hence:
\begin{align} \label{E:0RIEMANNINVARIANTSMAPPEDTOPOSITIVEDENSITYCONSTANT}
	\InverseRiemannfunction(0)
	& = \overline{\varrho} > 0.
\end{align}
Importantly, the chain rule, \eqref{E:SPHERICALRIEMANNINVARIANTS}, and the fundamental theorem of calculus imply that:
    \begin{align} \label{E:DERIVATIVEOFINVERSERIEMANNFUNCTION}
        I'(z)
				&:=
				\frac{d}{d z} I(z) = \frac{I(z)}{\Speed \circ I(z)}.
    \end{align}
Moreover, for use later on, we define:
\begin{align} \label{E:DERIVATIVEOFSPEEDOFSOUND}
	\Speed'
	:=
	\frac{d}{d\varrho} \Speed(\varrho).
\end{align}

We also note that in view of \eqref{E:BACKGROUNDSPEEDOFSOUNDISUNITYNORMALIZATION} and  
\eqref{E:ORIGINALVARIABLESINTERMSOFRIEMANNINVARIANTS}--\eqref{E:DERIVATIVEOFINVERSERIEMANNFUNCTION},
we can express:
\begin{subequations}
\begin{align} \label{E:SPEEDOFSOUNDEXPANSION}
	\Speed 
	& = 1 + f(\RRiemann,\LRiemann),
\end{align}
where $f$ is a $C^3$ function satisfying:
\begin{align} \label{E:SPEEDOFSOUNDERRORFUNCTIONVANISHESATORIGIN}
	f(0,0) & = 0.
\end{align}
\end{subequations}
In particular, in the regime that we study in this paper, in which the Riemann invariants are close to $0$, the speed of sound will be a small perturbation of $1$. We will often silently use this basic fact throughout the paper.

The Riemann invariants ``diagonalize'' the Euler equations relative to the \emph{characteristic vectorfields} given in the following definition.

\begin{definition}[Characteristic vectorfields] 
\label{D:NULLVECTORFIELDS}
We define the future-outgoing characteristic vectorfield (also known as a ``null vectorfield'') $\Lunit$ and 
the future-ingoing characteristic vectorfield $\uLunit$ as follows:
\begin{align} \label{E:NULLVECTORFIELDS}
	\Lunit
	& 
	:= \partial_t + (v^r + \Speed) \partial_r,
	&
	\uLunit
	& := \partial_t + (v^r - \Speed) \partial_r. 
\end{align}
\end{definition}

For use throughout the paper, we note the following identities,
which follow from \eqref{E:MATERIALDERIVATIVEVECOTRFIELD}, \eqref{E:NULLVECTORFIELDS},
and straightforward computations: 
\begin{subequations}
\begin{align} \label{E:PARTIALRINTERMSOFLUNITANDULUNIT}
	\partial_r
	& = 
			\frac{1}{2 \Speed}
			\left\lbrace
				\Lunit 
				-
				\uLunit
			\right\rbrace, 
				\\
	\Transport
	& = \frac{1}{2}
			\left\lbrace
				\Lunit
				+ 
				\uLunit
			\right\rbrace,
				\label{E:MATERIALDERIVATIVEINTERMSOFNULLVECTORFIELDS}
				\\
\partial_t
& = 
			\frac{1}{2}
			\left\lbrace
				1 
				-
				\frac{v^r}{\Speed}
			\right\rbrace
			\Lunit
			+
			\frac{1}{2}
			\left\lbrace
				1 
				+  
				\frac{v^r}{\Speed}
			\right\rbrace 
			\uLunit.
		\label{E:PARTIALTINTERMSOFLANDLBAR}
\end{align}
\end{subequations}

We also note that based on \eqref{E:NULLVECTORFIELDS} and \eqref{E:PARTIALRINTERMSOFLUNITANDULUNIT}, 
the commutator of $\Lunit$ and $\uLunit$ follows from straightforward computations:
\begin{align} 
\begin{split} \label{E:COMMUTATOROFLANDULUNIT}
	[\Lunit,\uLunit]
	& = 
			\frac{1}{2 \Speed}
			\left\lbrace
				\Lunit v^r
				-
				\Speed' 
				\Lunit \varrho
			\right\rbrace
			\left\lbrace
				\Lunit
				-
				\uLunit
			\right\rbrace	
			-
			\frac{1}{2 \Speed}
			\left\lbrace
				\uLunit v^r
				+
				\Speed' 
				\uLunit \varrho
			\right\rbrace
			\left\lbrace
				\Lunit
				-
				\uLunit
			\right\rbrace	
			\\
			& =
			\frac{1}{2 \Speed}
			\left\lbrace
				\Lunit v^r
				-
				\Speed' 
				\Lunit \varrho
				-
				\uLunit v^r
				-
				\Speed' 
				\uLunit \varrho
			\right\rbrace
			\Lunit
			+
			\frac{1}{2 \Speed}
			\left\lbrace
				-
				\Lunit v^r
				+
				\Speed' 
				\Lunit \varrho
				+
				\uLunit v^r
				+
				\Speed' 
				\uLunit \varrho
			\right\rbrace
			\uLunit.
	\end{split}
\end{align}

\subsection{Diagonalizing the fluid flow}
\label{SS:DIAGONALIZINGTHEFLOW}
In the next lemma, we derive the aforementioned diagonalization of the Euler equations in spherical symmetry.

\begin{lemma}[Equivalent formulation of the evolution equations in terms the Riemann invariants]
\label{L:EVOLUTIONFORRIEMANNINVARIANTS}
Let $\Riemanninvariant_{(\pm)}$ be as in \eqref{E:SPHERICALRIEMANNINVARIANTS},
and let $\Lunit,\uLunit$ be the null vectorfields defined in \eqref{E:NULLVECTORFIELDS}.
If $v^r$ and $\varrho$ are $C^1$ functions, then
equations \eqref{E:VELOCITYRADIALEQUATION}--\eqref{E:DENSITYRADIALEQUATION} hold if and only if
the following evolution equations hold (see Footnote~\ref{FN:VECTORFIELDDERIVATIVESOFSCALARS} regarding the notation):
\begin{subequations}
\begin{align}
	\Lunit \RRiemann
	& = \frac{\Speed}{r} (\LRiemann - \RRiemann),
		\label{E:OUTGOINGRIEMANNINVARIANTEVOLUTION} 
			\\
	\uLunit \LRiemann
	& = \frac{\Speed}{r} (\LRiemann - \RRiemann).
	\label{E:INGOINGRIEMANNINVARIANTEVOLUTION} 
\end{align}
\end{subequations}

\end{lemma}

\begin{proof}
	We will only derive \eqref{E:OUTGOINGRIEMANNINVARIANTEVOLUTION}--\eqref{E:INGOINGRIEMANNINVARIANTEVOLUTION}
	as a consequence of \eqref{E:VELOCITYRADIALEQUATION}--\eqref{E:DENSITYRADIALEQUATION}; the reverse derivation
	can be derived by reversing the arguments.
	
	Expanding $2 \times$ LHS~\eqref{E:OUTGOINGRIEMANNINVARIANTEVOLUTION} with the help of
	\eqref{E:SPHERICALRIEMANNINVARIANTS} and \eqref{E:NULLVECTORFIELDS}, 
	we deduce that
	$2 \Lunit \RRiemann 
	= \partial_t v^r 
	+ \frac{\Speed}{\varrho} \partial_t \varrho 
	+ (v^r + \Speed) \partial_r v^r 
	+ \frac{\Speed}{\varrho}  (v^r + \Speed) \partial_r \varrho$.
	Equations \eqref{E:VELOCITYRADIALEQUATION} and \eqref{E:DENSITYRADIALEQUATION}
	imply that the RHS of this identity is equal to $- 2 \Speed \frac{v^r}{r}$ which, in view of
	\eqref{E:ORIGINALVARIABLESINTERMSOFRIEMANNINVARIANTS}, is equal to $2 \times$ RHS~\eqref{E:OUTGOINGRIEMANNINVARIANTEVOLUTION}.
	We have therefore proved \eqref{E:OUTGOINGRIEMANNINVARIANTEVOLUTION}.
	
	Equation \eqref{E:INGOINGRIEMANNINVARIANTEVOLUTION} follows from a similar argument, and we omit the details.
		
\end{proof}

To derive estimates for solutions to \eqref{E:OUTGOINGRIEMANNINVARIANTEVOLUTION}--\eqref{E:INGOINGRIEMANNINVARIANTEVOLUTION} 
all the way down to $r=0$, we will use the following boundary conditions, 
which follow from Definition~\ref{D:SPHERICALRIEMANNINVARIANTS}
and \eqref{E:EVENDERIVATIVESOFVRVANISHATTHEORIGIN}--\eqref{E:ODDDERIVATIVESOFDENSITYVANISHATTHEORIGIN}:
\begin{align} \label{E:RIEMANNINVARIANTMATCHINGCONDITION}
	\partial_r^{2k} 
	\RRiemann \restriction_{\lbrace r=0 \rbrace}
	& 
	= 
	\partial_r^{2k} \LRiemann \restriction_{\lbrace r=0 \rbrace},
		\\
	\label{E:RIEMANNINVARIANTRDERIVATIVEMATCHINGCONDITION}
	\partial_r^{2k+1} \RRiemann \restriction_{\lbrace r=0 \rbrace}
	& 
	= 
	- \partial_r^{2k+1} \LRiemann \restriction_{\lbrace r=0 \rbrace}.
\end{align}

\subsection{Wave-like evolution equations for the Interior Region}
\label{SS:WAVEEQUATIONSFORINTERIOR}
To derive decay estimates in the Interior Region (i.e., close to the time axis), 
we will rely on the following wave-like equations,
which are straightforward consequences of
\eqref{E:OUTGOINGRIEMANNINVARIANTEVOLUTION}--\eqref{E:INGOINGRIEMANNINVARIANTEVOLUTION}.

\begin{proposition}[Wave-like evolution equations for the Interior Region]
   \label{P:WAVELIKEEQUATIONSFORTHERIEMANNINVARIANTS}
Let $(\RRiemann,\LRiemann)$ be a 
$C^2$ solution to equations \eqref{E:OUTGOINGRIEMANNINVARIANTEVOLUTION}--\eqref{E:INGOINGRIEMANNINVARIANTEVOLUTION}.
Then the following equations are also satisfied:
    \begin{subequations}
        \begin{align}
            \begin{split} \label{E:RIGHTRIEMANNEQUATIONINTERIORREGION}
                \Lunit \Big( r \uLunit \RRiemann - 2 \Speed \RRiemann\Big) & = v^r \uLunit \RRiemann - v^r \Lunit \RRiemann - 2 \RRiemann \Lunit \Speed \\
                &  \ \ - r [\uLunit,\Lunit]\RRiemann + \frac{r}{\Speed} (\uLunit \Speed) \Lunit \RRiemann
									\\
								& = 
								- \boxed{2 \RRiemann \Lunit \Speed}
								-\frac{r}{\Speed} (\Lunit \RRiemann) \uLunit \RRiemann 
								+
								\frac{r}{\Speed} (\Lunit \RRiemann) \Lunit \RRiemann  \\
                &  \ \ - r [\uLunit,\Lunit]\RRiemann + \frac{r}{\Speed} (\uLunit \Speed) \Lunit \RRiemann,
									\\
            \end{split} \\
            \begin{split} \label{E:LEFTRIEMANNEQUATIONINTERIORREGION}
                \uLunit \Big( r \Lunit \LRiemann +  2\Speed\LRiemann\Big)  
								& = v^r \Lunit \LRiemann  - v^r \uLunit \LRiemann  + 2 \LRiemann   \uLunit \Speed \\
                & \ \ - r [\Lunit,\uLunit]\LRiemann + \frac{r}{\Speed} (\Lunit \Speed)\uLunit \LRiemann
									\\
								& = 
								\boxed{2 \LRiemann \uLunit \Speed}
								- \frac{r}{\Speed} (\Lunit \LRiemann)^2  + \frac{r}{\Speed} (\Lunit \RRiemann) \uLunit \LRiemann  
								+\\
                & \ \ - r [\Lunit,\uLunit]\LRiemann + \frac{r}{\Speed} (\Lunit \Speed) \uLunit \LRiemann.
            \end{split} 
        \end{align}
    \end{subequations}
\end{proposition}

\begin{remark}[Boxed terms]
\label{R:BOXEDTERMS}
In some estimates
we will have to exploit cancellations tied to the boxed terms on 
RHSs~\eqref{E:RIGHTRIEMANNEQUATIONINTERIORREGION}--\eqref{E:LEFTRIEMANNEQUATIONINTERIORREGION},
which can be used to partially cancel corresponding terms on the LHSs of the equations.
Specifically, such cancellations are important for some of the top-order estimates
in Sect.\,\ref{S:GLOBALEXISTENCE}, in the region where $r$ is small.
Similar remarks apply to the boxed terms on
RHSs~\eqref{E:REVAMPEDRPLUSWAVEEQUATION}--\eqref{E:REVAMPEDRMINUSWAVEEQUATION}.
\end{remark}

\begin{proof}[Proof of Prop.\,\ref{P:WAVELIKEEQUATIONSFORTHERIEMANNINVARIANTS}]
    We prove \eqref{E:LEFTRIEMANNEQUATIONINTERIORREGION}. \eqref{E:RIGHTRIEMANNEQUATIONINTERIORREGION} 
		can be proved using the same arguments. In fact, the proof is the same after one 
		interchanges the roles of $\Lunit,\uLunit$ and $\RRiemann, \LRiemann$. 
    
    We begin by multiplying both sides of \eqref{E:INGOINGRIEMANNINVARIANTEVOLUTION} by $\frac{r}{\Speed}$ to deduce that 
		$\frac{r}{\Speed} \uLunit \LRiemann = \LRiemann - \RRiemann$. Differentiating this identity 
		with respect to $\Lunit$ and using the fact that $\Lunit \RRiemann = \uLunit \LRiemann$ 
		(see \eqref{E:OUTGOINGRIEMANNINVARIANTEVOLUTION}--\eqref{E:INGOINGRIEMANNINVARIANTEVOLUTION}), we have:
    \begin{align} \label{E:LEFTRIEMANNEQUATIONINTERIORREGIONSTEP1}
        \frac{r}{\Speed} \Lunit \uLunit \LRiemann + \frac{\Lunit r}{\Speed} \uLunit \LRiemann -   \frac{r}{\Speed^2} (\Lunit \Speed) \uLunit \LRiemann = \Lunit \LRiemann - \uLunit \LRiemann.
    \end{align}
    Commuting and differentiating by parts, we have that \eqref{E:LEFTRIEMANNEQUATIONINTERIORREGIONSTEP1} is equivalent to: 
    \begin{align}
        \begin{split} \label{E:LEFTRIEMANNEQUATIONINTERIORREGIONSTEP2}
            & \frac{1}{\Speed} \uLunit \left( r \Lunit \LRiemann \right) + \frac{r}{\Speed} [\Lunit,\uLunit]\LRiemann - \frac{\uLunit r}{\Speed} \Lunit\LRiemann \\
            & \ \ + \frac{\Lunit r}{\Speed} \uLunit \LRiemann -   \frac{r}{\Speed^2} (\Lunit \Speed) \uLunit \LRiemann  =  \Lunit \LRiemann - \uLunit \LRiemann.
        \end{split}
    \end{align}
    Bringing the $-\uLunit \LRiemann$ from the RHS \eqref{E:LEFTRIEMANNEQUATIONINTERIORREGIONSTEP2} over to the left, and putting all of the remaining terms on the LHS over to the right \emph{except} for the linear term $\frac{\Lunit r}{\Speed} \uLunit \LRiemann$, we have:
    \begin{align}
        \begin{split} \label{E:LEFTRIEMANNEQUATIONINTERIORREGIONSTEP3}
            \frac{1}{\Speed}\uLunit \left(r \Lunit \LRiemann \right)  + \left(1+\frac{\Lunit r}{\Speed}\right) \uLunit \LRiemann 
            & = \left(1 + \frac{\uLunit r}{\Speed}\right) \Lunit \LRiemann - \frac{r}{\Speed} [\Lunit,\uLunit] \LRiemann  + \frac{r}{\Speed^2} (\Lunit \Speed) \uLunit \LRiemann. 
        \end{split}
    \end{align}
    Multiplying \eqref{E:LEFTRIEMANNEQUATIONINTERIORREGIONSTEP3} by $\Speed$, differentiating by parts, and using
		the identities $\Lunit r = v^r+\Speed$ and $\uLunit r = v^r - \Speed$ (which follow from \eqref{E:NULLVECTORFIELDS}), we have: 
	\begin{align}
        \begin{split} \label{E:LEFTRIEMANNEQUATIONINTERIORREGIONSTEP4}
          \uLunit \Big( r \Lunit \LRiemann +  2\Speed\LRiemann\Big)  & = v^r \Lunit \LRiemann  - v^r \uLunit \LRiemann  + 2 \LRiemann   \uLunit \Speed - r [\Lunit,\uLunit]\LRiemann + \frac{r}{\Speed} (\Lunit \Speed)\uLunit \LRiemann,
        \end{split}
    \end{align}
		which yields the first equality in \eqref{E:LEFTRIEMANNEQUATIONINTERIORREGION}.
		The second equality follows from using equations 
		\eqref{E:ORIGINALVARIABLESINTERMSOFRIEMANNINVARIANTS} and
		\eqref{E:OUTGOINGRIEMANNINVARIANTEVOLUTION} to algebraically substitute the factors of $v^r$ in
		the first equality with $-\frac{r}{\Speed} (\Lunit \RRiemann)$.
\end{proof}

\section{Acoustic geometry and nonlinear geometric optics}
\label{S:ACOUSTICGEOMETRY}
In this section, we construct the acoustic geometry, i.e, the characteristics corresponding to the outgoing null vectorfield $\Lunit$, as well as some other related geometric quantities. This is known as implementing \emph{nonlinear geometric optics}.
They are essential for our sharp analysis of the gradient-blowup.

\subsection{Construction of the eikonal function}
\label{SS:CONSTRUCTIONOFEIKONAL}
Let $\tstar \geq 0$ be a fixed time. Later, we will set $\tstar=\datasize^{-1}$ for a sufficiently small data-parameter $\datasize$, 
which we introduce in Sect.\,\ref{S:DATA}.

By definition, the eikonal function $u$ is the solution to the following transport equation initial value problem:
\begin{align} \label{E:ODEFOREXTERIOREIKONALFUNCTION}
	\Lunit u 
	& = 0,
	&
	u \restriction_{\Sigma_{\tstar} \cap \lbrace r \geq \frac{\tstar}{2} \rbrace} 
	& 
	= \tstar - r.
\end{align}
Note that the data for $u$ are posed at time $\tstar$ and only in the region $r \geq \frac{\tstar}{2}$,
which corresponds to $u \leq \frac{\tstar}{2}$. 
This setup has two main advantages:

\begin{itemize}
	\item $u$ is defined only away from $r = 0$, where the equations degenerate due to singular factors of $r$
		(see \eqref{E:OUTGOINGRIEMANNINVARIANTEVOLUTION}--\eqref{E:INGOINGRIEMANNINVARIANTEVOLUTION}).
	\item The gradient-blowup will occur approximately in a subset of $\lbrace u \leq 0 \rbrace$. Thus, when $\tstar=\datasize^{-1}$, the gradient-blowup 
		occurs well within the region where $u$ is defined. This is important in that we fundamentally rely on $u$ to derive the
		singular boundary, i.e., the entire curve of gradient-blowup.
\end{itemize}

\subsection{Geometric coordinates and partial derivatives}
\label{SS:GEOMETRICCOORDINATES}
Most of our analysis, especially the part concerning gradient-blowup, 
will take place with respect to the differential structure of the \emph{geometric coordinates}, which we now define.

\begin{definition}[Geometric coordinates and partial derivatives]
Let $u$ be the solution to \eqref{E:ODEFOREXTERIOREIKONALFUNCTION}.
We refer to $(t,u)$ as the geometric coordinates.\footnote{In Prop.\,\ref{P:CHOVDIFFEOMORPHISM}, we prove that $(t,u)$ are a viable coordinate system. \label{FN:GEOMETRICCOORDINATESARECOORIDNATES}}
We denote the corresponding partial derivative
vectorfields by\footnote{Note that $\partial_t \neq \frac{\partial}{\partial t}$. \label{FN:TWOPARTIALTSARENOTEQUAL}} 
$\lbrace \frac{\partial}{\partial t}, \frac{\partial}{\partial u} \rbrace$.
\end{definition}

\subsection{Inverse foliation density}
\label{SS:INVERSEFOLIATIONDENSITY}
As is standard in shock-formation problems, the gradient-blowup of the fluid is tied to the gradient-blowup of $u$, which signifies the \emph{infinite density} 
of the characteristics. Rather than studying the density of the characteristics, we prefer to study its reciprocal $\upmu$, which we now rigorously define.

\begin{definition} [Inverse foliation density]
\label{D:INVERSEFOLIATIONDENSITY}
Let $\Transport$ and $\uLunit$ be the vectorfields defined in \eqref{E:MATERIALDERIVATIVEVECOTRFIELD} and \eqref{E:NULLVECTORFIELDS}.
We define the inverse foliation density to be the following scalar function $\upmu$:
\begin{align} \label{E:DEFOFINVERSEFOLIATIONDENSITY}
\upmu
& := \frac{1}{\Transport u}
	= \frac{2}{\uLunit u}
\end{align}
where the second equality follows from \eqref{E:MATERIALDERIVATIVEINTERMSOFNULLVECTORFIELDS} and the fact
that $\Lunit u = 0$.
\end{definition}

\begin{remark}[The behavior of $\upmu$ in our main results]
\label{R:BEHAVIOROFMUINMAINRESULTS}
In our main gradient-blowup results, 
we will have $\upmu > 0$, except along the gradient-blowup curve, where $\upmu = 0$. In particular, the vanishing of $\upmu$
signifies the infinite density of the characteristics and the blowup of the fluid.
In our global existence results, we will prove that $\upmu$ is strictly positive in all of spacetime.
\end{remark}

For future use, we use 
\eqref{E:PARTIALRINTERMSOFLUNITANDULUNIT}--\eqref{E:PARTIALTINTERMSOFLANDLBAR}, 
\eqref{E:ODEFOREXTERIOREIKONALFUNCTION},
and
\eqref{E:DEFOFINVERSEFOLIATIONDENSITY}
to derive the following identity for the data of $\upmu$ (at time $\tstar$):
\begin{align} \label{E:MGHDIDENTITYFORMUATTIMEZERO}	
	\upmu \restriction_{\Sigma_{\tstar} \cap\lbrace r \geq \frac{\tstar}{2} \rbrace}
	& = \frac{1}{\Speed} \restriction_{\Sigma_{\tstar} \cap \lbrace  r \geq \frac{\tstar}{2} \rbrace}.
\end{align}

\subsection{\texorpdfstring{$\upmu$}{inverse foliation density}-weighted vectorfields}
\label{SS:WEIGHTEDVECTORFIELDS}
Our PDE analysis fundamentally relies on $\upmu$-weighted vectorfields.

\begin{definition}[$\upmu$-weighted vectorfields]
\label{D:MUWEIGHTEDVECTORFIELDS}
Let $\upmu$ be the inverse foliation density from \eqref{E:DEFOFINVERSEFOLIATIONDENSITY}.
We define the following $\upmu$-weighted vectorfields.
\begin{align} \label{E:MUWEIGHTEDVECTORFIELDS}
	\muX 
	& := - \Speed \upmu \partial_r,
	&
	\muuLunit
	& := \upmu \uLunit.
\end{align}
\end{definition}

In the next lemma, we record some basic but important properties of the above vectorfields.

\begin{lemma}[Properties of the $\upmu$-weighted vectorfields] 
\label{L:VECTORFIELDSINTERMSOFGEOMETRICCOORDINATES}
The following identities hold:
\begin{align} \label{E:VECTORFIELDSINTERMSOFGEOMETRICCOORDINATES}
	\Lunit
	& = \frac{\partial}{\partial t},
		\\
	\muX  
	& = \frac{\partial}{\partial u},
	&
	\muuLunit
	& = \upmu 
		\frac{\partial}{\partial t}
		+
		2 \frac{\partial}{\partial u}.
\end{align}	

In particular, we have the following commutation identities:
\begin{subequations}
\begin{align}
	[\Lunit, \muX]
	& = 0,
		\label{E:VANISHINGCOMMUTATOROFLANDMUX} 
		\\
	[\muuLunit,\muX]
	& = - (\muX \upmu) \Lunit,
		\label{E:COMMUATATOROFMUULUNITANDMUX} 
			\\
[\Lunit,\muuLunit]
	& = (\Lunit \upmu)
			\Lunit.
\label{E:COMMUTATOROFLANDMUULUNIT}
\end{align}
\end{subequations}

\end{lemma}

\begin{proof}
		From \eqref{E:NULLVECTORFIELDS} and \eqref{E:ODEFOREXTERIOREIKONALFUNCTION}, we have that
		$\Lunit t = 1$ and $\Lunit u = 0$. Hence, in geometric coordinates, we must have
		$\Lunit = \frac{\partial}{\partial t}$.
		Moreover, by \eqref{E:MUWEIGHTEDVECTORFIELDS}, 
		\eqref{E:PARTIALRINTERMSOFLUNITANDULUNIT},
		and
		\eqref{E:ODEFOREXTERIOREIKONALFUNCTION},
		we have $\muX t = 0$ and 
		$\muX u = - \Speed \upmu \partial_r u =  \Speed \upmu \frac{\uLunit u}{2 \Speed} = 1$.
		Hence, in geometric coordinates, we must have
		$\muX = \frac{\partial}{\partial u}$.
		The identity 
		$\muuLunit
		= \upmu 
		\frac{\partial}{\partial t}
		+
		2 \frac{\partial}{\partial u}
		$
		follows from the other two identities and \eqref{E:PARTIALRINTERMSOFLUNITANDULUNIT}.

		\eqref{E:VANISHINGCOMMUTATOROFLANDMUX}--\eqref{E:COMMUTATOROFLANDMUULUNIT} follow easily from \eqref{E:VECTORFIELDSINTERMSOFGEOMETRICCOORDINATES}
		and Def.\,\ref{D:MUWEIGHTEDVECTORFIELDS}.
\end{proof}

\subsection{The change of variables map}
\label{SS:CHOVMAP}
Our analysis relies on a careful study of the change of variables map from geometric to Cartesian coordinates.
In this short section, we define it and reveal its basic properties.

\begin{definition}[Change of variables map]
\label{D:CHOVFROMTUTOTRCOORDINATES}
We define the change of variables map from geometric to Cartesian coordinates as follows:
\begin{align} \label{E:CHOVFROMTUTOTRCOORDINATES}
	\Upsilon(t,u) := (t,r),
\end{align}
where on RHS~\eqref{E:CHOVFROMTUTOTRCOORDINATES}, we are viewing $r = r(t,u)$.
\end{definition}

\begin{lemma}[Basic properties of $\Upsilon$]
\label{L:BASICPROPERTIESOFCHOVMAP}
	The Jacobian matrix of $\Upsilon$ is as follows:
	\begin{align} \label{E:CHOVJACOBIAN}
		\mathrm{d} \Upsilon(t,u)
		& := 
			\begin{pmatrix}
					\frac{\partial}{\partial t} t &  \frac{\partial}{\partial u} t
						\\
				 \frac{\partial}{\partial t} r & \frac{\partial}{\partial u} r
			\end{pmatrix}
			= 
			\begin{pmatrix}
					1 & 0 
						\\
				 \Lunit r & \upmu X r
				\end{pmatrix}
			=
			\begin{pmatrix}
					1 & 0 
						\\
					v^r + c & - \Speed \upmu
			\end{pmatrix}.
	\end{align}

Moreover, its Jacobian determinant is as follows:
\begin{align} \label{E:MAINSTATEMENTCHOVJACOBIANDTERMINANT}
		\textnormal{det} (\mathrm{d} \Upsilon(t,u))
		& = - \Speed \upmu.
	\end{align}

\end{lemma}

\begin{proof}
	\eqref{E:CHOVJACOBIAN} follows easily from \eqref{E:NULLVECTORFIELDS}
	and Lemma~\ref{L:VECTORFIELDSINTERMSOFGEOMETRICCOORDINATES}.
	
	\eqref{E:MAINSTATEMENTCHOVJACOBIANDTERMINANT} then follows easily from~\eqref{E:CHOVJACOBIAN}.
\end{proof}

\begin{remark}[Notational convention regarding the two different coordinate systems]
	\label{R:NOTATIONALCONVENTIONCOORDIANTEDEPENDENCE}
    As we explained in Sect.\,\ref{SSS:INTRONONLINEARGEOMETRICOPTICS} and Sect.\,\ref{SS:CONSTRUCTIONOFEIKONAL}, we carry out our PDE analysis in geometric coordinates $(t,u)$ \emph{only in exterior regions}, i.e., regions away from $\{r=0\}$; see Def.\,\ref{D:GLOBALLHYPERBOLICEXTERIORBOOTSTRAPREGION}. In contrast, we carry out our PDE analysis in interior regions containing $\{r=0\}$ in the standard spherical coordinate $(t,r)$; see Def.\,\ref{D:GHINTERIORBOOTSTRAP}. As a matter of convenience, in most of the article, 
		we do \emph{not} make any notational distinction for functions (such as $\RRiemann$) defined on exterior regions in geometric coordinates versus interior regions in spherical coordinates. That is, we denote $\RRiemann(t,u)$ for the former and $\RRiemann(t,r)$ for the latter, or even use notation without reference to any coordinate system, such as $\RRiemann$. This is consistent with the conventions of differential geometry.
		Similar statements hold for $\LRiemann$, as well other geometric objects that we use, i.e., we denote $\Lunit = \frac{\partial}{\partial t}$ as well as $\Lunit = \partial_t + (v^r + \Speed)\partial_r$. This is not an issue away from the MGHD boundary because the two coordinate systems are related by the change of variables 
		$\Upsilon(t,u) = (t,r)$ that is shown (in Prop.\,\ref{P:EXISTENCEUPTOCREASEANDSINGULARBOUNDARYANDPORTIONOFCAUCHYHORIZON}) 
		to be as smooth as the solution away from the MGHD boundary.

    One exception is that in certain parts of Prop.\,\ref{P:EXISTENCEUPTOCREASEANDSINGULARBOUNDARYANDPORTIONOFCAUCHYHORIZON}, 
		we prove results in an exterior region in $(t,u)$-space that contains part of the MGHD boundary, 
		as well as other results concerning the behavior of the solution in $(t,r)$-coordinates
		on the corresponding image of the region under $\Upsilon$,
		and the two coordinate systems are continuously -- but not smoothly -- related on the MGHD boundary.
		Hence, different conclusions hold in the two coordinate systems when the MGHD boundary is taken into account 
		(see, e.g., Points~2, 5, and 8 in Prop.\,\ref{P:EXISTENCEUPTOCREASEANDSINGULARBOUNDARYANDPORTIONOFCAUCHYHORIZON}), 
		and there we use notation such as $\RRiemann\circ \Upsilon^{-1}$, 
		$\LRiemann \circ \Upsilon^{-1}$ to explicitly indicate that we are viewing the Riemann invariants as a function of $(t,r)$.
    
\end{remark}



\section{Geometric evolution equations in the Exterior Region}
\label{S:GEOMETRICEVOLUTIONEQUATIONSINEXTERIOR}
In this section, we derive the evolution equations that we use to study the solution in the ``Exterior Region'' away from $r=0$, where the most interesting analysis occurs.

We start with a lemma that yields a transport equation satisfied by the inverse foliation density.

\begin{lemma}[Transport equation equation for $\upmu$]
\label{L:MUEVOLUTION}
The inverse foliation density defined in \eqref{E:DEFOFINVERSEFOLIATIONDENSITY}
satisfies the following transport equation,
where $\Speed' = \Speed'(\varrho) = \Speed' \circ I (\RRiemann + \RRiemann)$ is the derivative of the speed of sound with respect to the density, and
$I = I(\RRiemann + \LRiemann)$ and $I'= I'(\RRiemann + \LRiemann)$ are the functions from \eqref{E:ORIGINALVARIABLESINTERMSOFRIEMANNINVARIANTS}--\eqref{E:DERIVATIVEOFINVERSERIEMANNFUNCTION}:
\begin{align} 
\begin{split} \label{E:MUEVOLUTION}
	\Lunit \upmu
	& = 
		-
			\frac{1}{2 \Speed}
			\left\lbrace
				\muuLunit v^r
				+
				\Speed' 
				\muuLunit \varrho
			\right\rbrace
			+
			\frac{\upmu }{2 \Speed}
			\left\lbrace
				\Lunit v^r
				-
				\Speed' 
				\Lunit \varrho
			\right\rbrace
				\\
			& = 
			-
			\frac{1}{2 \Speed}
			\left\lbrace
				1
				+
				\Speed' 
				\InverseRiemannfunction'
			\right\rbrace
			\muuLunit 
			\RRiemann
			+
			\frac{1}{2 \Speed}
			\left\lbrace
				1
				-
				\Speed' 
				\InverseRiemannfunction'
			\right\rbrace
			\muuLunit \LRiemann
				\\
		& \ \
			+
			\frac{\upmu }{2 \Speed}
			\left\lbrace
				1
				-
				\Speed' 
				\InverseRiemannfunction'
			\right\rbrace
			\Lunit 
			\RRiemann
			-
			\frac{\upmu }{2 \Speed}
			\left\lbrace
				1
				+
				\Speed' 
				\InverseRiemannfunction'
			\right\rbrace
			\Lunit 
			\LRiemann.
\end{split}
\end{align}
\end{lemma}

\begin{proof}
	Recall that $\Lunit u = 0$ by construction. Hence using
	\eqref{E:DERIVATIVEOFSPEEDOFSOUND},
	\eqref{E:COMMUTATOROFLANDULUNIT}, 
	and
	\eqref{E:DEFOFINVERSEFOLIATIONDENSITY},
	we see that
	$\Lunit \upmu = -\frac{2}{(\uLunit u)^2} \Lunit \uLunit u 
	=  
	-\frac{2}{(\uLunit u)^2}
	[\Lunit, \uLunit] u
	=-
	\frac{1}{\uLunit u}
	\frac{1}{\Speed}
			\left\lbrace
				-
				\Lunit v^r
				+
				\Speed' 
				\Lunit \varrho
				+
				\uLunit v^r
				+
				\Speed' 
				\uLunit \varrho
			\right\rbrace
	$,
	which, in view of 
	\eqref{E:ORIGINALVARIABLESINTERMSOFRIEMANNINVARIANTS}--\eqref{E:INVERSERIEMANNFUNCTION},
	\eqref{E:DEFOFINVERSEFOLIATIONDENSITY},
	and
	\eqref{E:MUWEIGHTEDVECTORFIELDS}, 
	implies \eqref{E:MUEVOLUTION}.
\end{proof}

\begin{remark}[Main term in the analysis]
	\label{R:MAINTERMINANALYSIS}
	In Theorem~\ref{T:MAINMGHDEXISTENCETHEOREM}, the main term that will drive the vanishing of $\upmu$ (and hence the gradient-blowup of $\RRiemann$)
	is the first product 
	$
	-
			\frac{1}{2 \Speed}
			\left\lbrace
				1
				+
				\Speed' 
				\InverseRiemannfunction'
			\right\rbrace
			\muuLunit 
			\RRiemann
			+
			\frac{1}{2 \Speed}
			\left\lbrace
				1
				-
				\Speed' 
				\InverseRiemannfunction'
			\right\rbrace
			\muuLunit \LRiemann
		$
		on RHS~\eqref{E:MUEVOLUTION}. It turns out that the details of some of the factors multiplying $\muuLunit \LRiemann$ 
		in this expression are not important. 
		For this reason, in the next definition, we define the null condition failure constant, which
		is obtained by Taylor expanding the factors to order $0$ around the background solution $(\RRiemann,\LRiemann) \equiv (0,0)$,
		thereby capturing the ``important part.''
\end{remark}

\begin{definition}[Null condition failure constant]
	\label{D:NULLCONDITIONFAILURECONSTANT}
	Recall that $\Speed'$ is the derivative of the speed of sound with respect to the density and that
	$I'$ is the function from \eqref{E:DERIVATIVEOFINVERSERIEMANNFUNCTION}.
	We define $\lifespanconstant$ to be the following constant:
	\begin{align} \label{E:NULLCONDITIONFAILURECONSTANT}
		\lifespanconstant
		& := 
				1
				+
				\Speed' 
				\InverseRiemannfunction'|_{(\RRiemann,\LRiemann) = (0,0)}.
	\end{align}
\end{definition}

Failure of the null condition near the background solution $(\RRiemann,\LRiemann) = (0,0)$ is equivalent 
to $\lifespanconstant \neq 0$. For almost all equations of state and background densities $\overline{\varrho}$, 
we have $\lifespanconstant \neq 0$. The one exception is 
a Chaplygin gas, a i.e., a fluid for which the equation of state is 
$p = C_1 - \frac{C_2}{\varrho}$,
where $C_1$ and $C_2$ are constants.
For a Chaplygin gas, the expression on RHS~\eqref{E:NULLCONDITIONFAILURECONSTANT} identically vanishes, and shocks are not expected to form. 
However, we can still prove global existence in the Chaplygin gas case. Unless we state otherwise, we assume that 
the equation of state is \emph{not} the Chaplygin gas one, and that $\lifespanconstant \neq 0$. For definiteness,
without loss of generality,\footnote{To treat the case $\lifespanconstant < 0$, we could just change the sign of the parameter $\datasize$ featured below in our assumptions on the data.}
we assume that:
\begin{align} \label{E:LIFESPANCONSTANT}
	\lifespanconstant
	& > 0.
\end{align}


We now derive geometric wave equations satisfied by the Riemann invariants in the region where the eikonal function is defined.
We will use the equations in conjunction with \eqref{E:MUEVOLUTION} to derive sharp estimates in the region
where the gradient-blowup occurs.
These equations feature important cancellations tied to $\upmu$, which are not present in the equations of
Prop.\,\ref{P:WAVELIKEEQUATIONSFORTHERIEMANNINVARIANTS}.
In particular, the equations of Prop.\,\ref{P:GEOMETRICWAVELIKEEVOLUTIONEQUATIONS} have the following key structural property 
(which is a version of the null condition): the dangerous quadratic term $(\muuLunit \RRiemann)^2$
does not appear on the RHS of either of the wave equations.

\begin{proposition}[Geometric wave-like evolution equations]
\label{P:GEOMETRICWAVELIKEEVOLUTIONEQUATIONS}
Let $(\RRiemann,\LRiemann)$ be a 
$C^2$ solution to equations \eqref{E:OUTGOINGRIEMANNINVARIANTEVOLUTION}--\eqref{E:INGOINGRIEMANNINVARIANTEVOLUTION}.
Then the following equations are also satisfied (see Remark~\ref{R:BOXEDTERMS} concerning the boxed terms):
\begin{subequations}
\begin{align}
	\begin{split} \label{E:REVAMPEDRPLUSWAVEEQUATION}
		\Lunit 
		\left\lbrace
			r
			\muuLunit \RRiemann
			-
			2 \upmu \Speed \RRiemann
		\right\rbrace
		& = 
			- \boxed{2 \RRiemann \Lunit (\upmu \Speed)} 
			- 
			\frac{r}{\Speed} (\Lunit \RRiemann) \muuLunit \RRiemann
			+
		  \upmu \frac{r}{\Speed} (\Lunit \RRiemann) \Lunit \RRiemann
				\\
		& \ \
			+ r (\Lunit \upmu) \Lunit \RRiemann
			+ \frac{r}{\Speed} (\muuLunit \Speed) \Lunit \RRiemann,
 \end{split}
			\\
	\begin{split} \label{E:REVAMPEDRMINUSWAVEEQUATION}
	\muuLunit
	\left\lbrace
		r \Lunit \LRiemann
		+
		2 \Speed \LRiemann
	\right\rbrace
	& =	
		\boxed{2 \LRiemann \muuLunit \Speed}
		+
		r (\Lunit \upmu) \Lunit \RRiemann
		-
		r (\Lunit \upmu) \Lunit \LRiemann
			\\
	& + \frac{r}{\Speed} (\Lunit \RRiemann) \muuLunit \LRiemann
	+ \upmu \frac{r}{\Speed} (\Lunit \Speed) \Lunit \RRiemann
	- \upmu \frac{r}{\Speed} (\Lunit \RRiemann) \Lunit \RRiemann.
\end{split}
\end{align}	
\end{subequations}
\end{proposition}

\begin{proof}[Proof of Prop.\,\ref{P:GEOMETRICWAVELIKEEVOLUTIONEQUATIONS}]
To prove \eqref{E:REVAMPEDRPLUSWAVEEQUATION}, we start by multiplying both sides of \eqref{E:OUTGOINGRIEMANNINVARIANTEVOLUTION} by $r$ to 
deduce that
$r \Lunit \RRiemann = \Speed(\LRiemann - \RRiemann)$. 
Differentiating this equation with respect to $\muuLunit$ and using the identity $\uLunit r = v^r - \Speed$
as well as the fact that $\Lunit \RRiemann = \uLunit \LRiemann$ (see \eqref{E:OUTGOINGRIEMANNINVARIANTEVOLUTION}--\eqref{E:INGOINGRIEMANNINVARIANTEVOLUTION}),
 we further deduce that:
    \begin{align} \label{E:RPLUSEQUATIONEXTERIORREGIONSTEP1}
      r \muuLunit \Lunit \RRiemann 
			& = (\muuLunit \Speed)(\LRiemann - \RRiemann)
				+
				\upmu \Speed \Lunit \RRiemann 
				+
				\upmu (\Speed - v^r) \Lunit \RRiemann
				-
				\Speed \muuLunit \RRiemann. 
    \end{align}
		Commuting the vectorfields $\muuLunit$ and $\Lunit$ on LHS~\eqref{E:RPLUSEQUATIONEXTERIORREGIONSTEP1} 
		with the help of \eqref{E:COMMUTATOROFLANDMUULUNIT},
		we find that:
    \begin{align} \label{E:RPLUSEQUATIONEXTERIORREGIONSTEP2}
       r \Lunit \muuLunit \RRiemann 
				& = 
				r (\Lunit \upmu) \Lunit \RRiemann 
				+
				(\muuLunit \Speed)(\LRiemann - \RRiemann)
				+
				\upmu \Speed \Lunit \RRiemann 
				+
				\upmu (\Speed - v^r) \Lunit \RRiemann
				-
				\Speed \muuLunit \RRiemann.
    \end{align}
		From \eqref{E:RPLUSEQUATIONEXTERIORREGIONSTEP2},
		the identity $\Lunit r = v^r + \Speed$, 
		and \eqref{E:ORIGINALVARIABLESINTERMSOFRIEMANNINVARIANTS},
		we deduce that:
		  \begin{align} \label{E:RPLUSEQUATIONEXTERIORREGIONSTEP3}
       \Lunit 
			\left\lbrace
					r  \muuLunit \RRiemann 
			\right\rbrace
				& = 
			 	r (\Lunit \upmu) \Lunit \RRiemann 
				-
				v^r
				\muuLunit \Speed
				+
				\upmu \Speed \Lunit \RRiemann 
				+
				\upmu (\Speed - v^r) \Lunit \RRiemann
				+
				v^r \muuLunit \RRiemann.
    \end{align}
		Using equations 
		\eqref{E:ORIGINALVARIABLESINTERMSOFRIEMANNINVARIANTS} and
		\eqref{E:OUTGOINGRIEMANNINVARIANTEVOLUTION} to algebraically substitute the factors of $v^r$ on RHS~\eqref{E:RPLUSEQUATIONEXTERIORREGIONSTEP3}
		with $-\frac{r}{\Speed} (\Lunit \RRiemann)$,
		and carrying out straightforward computations, we conclude \eqref{E:REVAMPEDRPLUSWAVEEQUATION}.
		
		Equation \eqref{E:REVAMPEDRMINUSWAVEEQUATION} follows from multiplying
		both sides of \eqref{E:INGOINGRIEMANNINVARIANTEVOLUTION} by $\upmu r$ to 
		deduce that
		$r \muuLunit \LRiemann = \upmu \Speed(\LRiemann - \RRiemann)$,
		and then taking an $\Lunit$ derivative of this identity and using arguments similar to the ones given above; we omit the details.

	\end{proof}

\section{Data and conventions for constants}
\label{S:DATA}
Our main goal in this section is to construct the open sets of Riemann invariant initial data 
for which our main MGHD existence and uniqueness results hold.
We first introduce our conventions for explicit and implicit constants;
we will use them throughout the paper in our PDE analysis.

\subsection{Conventions for constants}
\label{SS:CONVENTIONSFORCONSTANTS}
\begin{itemize}
    \item Explicit constants $C > 0$ are allowed to vary from line to line. They are allowed to depend on the equation of state, but can be chosen to be \emph{independent} of the amplitude of the data $\datasize$ (see \eqref{E:BACKGROUNDPROFILEWITHSMALLAMPLITUDEFACTOR} and the bootstrap parameter 
		$\eps$ (see \eqref{E:BOUNDONBOOTSTRAPPARAMETERSIZE},\eqref{E:EPSILONSMALLNESSASSUMPTIONININTERIOREGION},\eqref{E:GLOBALEXISTENCEBOUNDONBOOTSTRAPPARAMETERSIZE}), provided that $\datasize$ and $\eps$ are sufficiently small.
    \item By $A \lesssim B$, we mean that there exists a $C > 0$ satisfying the properties above such that $A \le C B$.
    \item By $A \approx B$, we mean that $A \lesssim B$ and $B \lesssim A$. 
    \item By $A = \mathcal{O}(B)$, we mean that there exists a $C>0$ with the properties above such that $|A| \le C B$.
    \item For any $\upalpha < 1$, $\upalpha^{1^+}$ denotes any quantity $Q$ such that there is a number $\updelta > 0$, 
		independent of all sufficiently small $\upalpha$, with $|Q| \leq \alpha^{1+ \updelta}$.
    \item We will also use the following standard ``little o'' notation:
\begin{align} \label{E:LITTLEONOTATION}
	\mylittleo(\datasize)
\end{align}
to denote any $\datasize$-dependent quantity $Q$ such $|Q| \to 0$ as $\datasize \downarrow 0$.
\end{itemize}

\subsection{Data at time \texorpdfstring{$0$}{zero}}
\label{SS:DATAATTIME0}
We denote the data for the Riemann invariants at time $0$ as follows:
\begin{align} \label{E:RIEMANNINVARIANTDATAATTIME0}
	\dataRRiemann(r)
	& := \RRiemann(0,r),
	&
	\dataLRiemann(r)
	& := \LRiemann(0,r).
\end{align}

For reasons explained in the discussion surrounding \eqref{E:RIEMANNINVARIANTMATCHINGCONDITION}--\eqref{E:RIEMANNINVARIANTRDERIVATIVEMATCHINGCONDITION},
we assume that the data satisfy the following boundary conditions:
\begin{subequations}
\begin{align} \label{E:VELOCITYTIMEAXISCOMPATABILITYCONDITIONFORNONLINEARDATAATTIME0}
	\partial_r^{2k} \dataRRiemann(0)
	& = \partial_r^{2k} \dataLRiemann(0),
	&& k = 0,1,
		\\
	\partial_r^{2k+1} \dataRRiemann(0)
	& = - \partial_r^{2k+1} \dataLRiemann(0),
	& & k = 0,1.
	\label{E:DENSITYTIMEAXISCOMPATABILITYCONDITIONFORNONLINEARDATAATTIME0}
\end{align}
\end{subequations}

Our data for $\mathcal{R}_{\pm}$ will be a \emph{perturbation} of the ``background'' $\profileoverr(r)$,
where:
\begin{align} \label{E:BACKGROUNDPROFILEWITHSMALLAMPLITUDEFACTOR}
	\profileoverr(r)
	& := - \datasize \frac{\arctan(r)}{r},
\end{align}
and the number $\datasize > 0$ is a sufficiently small amplitude.

\begin{remark}[Specific profile only for convenience]
		\label{R:SPECIFICPROFILE}
	It is only for convenience/ to help guide the reader that we have focused on the background \eqref{E:BACKGROUNDPROFILEWITHSMALLAMPLITUDEFACTOR}.
	Our main results can be extended to a much larger class of background ``profiles'' that enjoy suitable properties. 
	We outline the extension in
	Appendix~\ref{A:EXTENDRESULTSTOOTHERPROFILES}.
\end{remark}	

Our analysis will show that the minus sign on RHS~\eqref{E:BACKGROUNDPROFILEWITHSMALLAMPLITUDEFACTOR} is
fundamental for ensuring that $\uLunit \RRiemann$ blows up in finite time.
Note that the background satisfies the following estimates:
\begin{subequations}
\begin{align}
		|\partial_r^K \profileoverr(r)|
		& \lesssim \frac{\datasize}{(1 + r)^{1+K}},
		& K & \geq 0,	\label{E:DATAISASYMPTOTICALLYFLAT} \\
		|\partial_r^K (r \profileoverr(r))|
		& \lesssim \frac{\datasize}{(1 + r)^{1+K}},
		& K & \geq 1. \label{E:RTIMESDATAISBETTER}
\end{align}
\end{subequations}
The estimates in \eqref{E:DATAISASYMPTOTICALLYFLAT} capture that the data are asymptotically flat.
The estimates in \eqref{E:RTIMESDATAISBETTER} capture that the data is in fact asymptotically flat in a strong sense.
\eqref{E:RTIMESDATAISBETTER} is fundamentally important for our analysis and is not implied by
\eqref{E:DATAISASYMPTOTICALLYFLAT}.

In our main MGHD theorem, we will consider data that are a perturbation of $\profileoverr$, where the undifferentiated perturbation has a faster decay rate at infinity compared to the background, and a smaller amplitude as well.  More precisely, to measure the size of the data perturbations, we will use the following weighted $C^3$ norm:
\begin{align} \label{E:PERTURBATIONTIME0NORM}
	\| \phi \|_{C_*^3(\Sigma_0)}
	& := 
		\sum_{K=0}^3
		\sup_{r \geq 0} 
		|(1 + r)^{2+K} \partial_r^K \phi(r)|.
\end{align}

To prove our main results, we will assume that $(\dataRRiemann,\dataLRiemann)$
is close to $(\profileoverr,\profileoverr)$ in the following sense:
\begin{subequations}
\begin{align}
		\| \dataRRiemann - \profileoverr  \|_{C_*^3(\Sigma_0)}
		& \leq \datasize^{1^+},
			 \label{E:RPLUSDATAISCLOSETOBACKGROUNDATTIME0} \\
		\| \dataLRiemann - \profileoverr  \|_{C_*^3(\Sigma_0)}
		& \leq \datasize^{1^+}.
		\label{E:RMINUSDATAISCLOSETOBACKGROUNDATTIME0}
\end{align}
\end{subequations}

\begin{remark}[Tiny velocity data and dominant density data]
Note that \label{R:TINYVELOCITYDATA}
\eqref{E:ORIGINALVARIABLESINTERMSOFRIEMANNINVARIANTS} and
\eqref{E:RPLUSDATAISCLOSETOBACKGROUNDATTIME0}--\eqref{E:RMINUSDATAISCLOSETOBACKGROUNDATTIME0} imply that
$\| v^r  \|_{C_*^3(\Sigma_0)} \lesssim \datasize^{1^+}$.
In particular, $|v^r \restriction_{t=0}|$ decays at least as fast as $\frac{\datasize^{1^+}}{(1+r)^2}$,
as we highlighted in the introduction.
On the other hand, 
\eqref{E:SPHERICALRIEMANNINVARIANTS}--\eqref{E:SPEEDOFSOUNDERRORFUNCTIONVANISHESATORIGIN},
\eqref{E:BACKGROUNDPROFILEWITHSMALLAMPLITUDEFACTOR},
and \eqref{E:RPLUSDATAISCLOSETOBACKGROUNDATTIME0}--\eqref{E:RMINUSDATAISCLOSETOBACKGROUNDATTIME0}
imply that $\varrho\restriction_{t=0} = \overline{\varrho} - 2 \datasize \frac{\arctan(r)}{r} + \mathcal{O}(\datasize^{1^+})\frac{1}{(1 + r)^2}$,
and hence the dominant term in the initial data is the density.
\end{remark}

\subsubsection{Estimates for the Riemann invariants at time $0$}
\label{SSS:ESTIMATESFORRIEMANNINVARIANTSATTIME0}
In this section, we use our assumptions on the data to derive estimates for the Riemann invariants at time $0$.

\begin{lemma}[Estimates for the Riemann invariants at time $0$]
\label{L:RIEMANNINVARIANTESTIMATESATTIME0}
Assume that the Riemann invariant data functions $(\dataRRiemann,\dataLRiemann)$ satisfy
\eqref{E:RPLUSDATAISCLOSETOBACKGROUNDATTIME0}--\eqref{E:RMINUSDATAISCLOSETOBACKGROUNDATTIME0}.
If $\datasize$ is sufficiently small, then the solution
$(\RRiemann,\LRiemann)$ to
\eqref{E:OUTGOINGRIEMANNINVARIANTEVOLUTION}--\eqref{E:INGOINGRIEMANNINVARIANTEVOLUTION}
obeys the following estimates at time $0$:

\noindent \underline{\textbf{Crude estimates for the Riemann invariants}}.
For $0 \leq J + K \leq 3$, we have:
\begin{subequations}
\begin{align} \label{E:TIME0CRUDEPOINTWISEBOUNDFORDATAOFRPLUSUPTO3DERIVATIVES}
	|\Lunit^J \uLunit^K \RRiemann(0,r)| 
	& \lesssim \frac{\datasize}{(1 + r)^{1+J+K}},
		\\
	|\partial_t \RRiemann(0,r)| 
	& \lesssim \frac{\datasize}{(1 + r)^2}.
	\label{E:TIME0CRUDEPOINTWISEBOUNDFORDATAOFPARTIALTRPLUS}
\end{align}
\end{subequations}

For $0 \leq J \leq 3$ and $0 \leq J + K' \leq 2$, we have:
\begin{subequations}
\begin{align} \label{E:TIME0CRUDEPOINTWISEBOUNDFORDATAOFRMINUSUPTO3LDERIVATIVES}
	|\Lunit^J \LRiemann(0,r)| 
	& \lesssim \frac{\datasize}{(1 + r)^{1+J}},
		\\
	|\Lunit^J \uLunit^{K'} \uLunit \LRiemann(0,r)| 
	& \lesssim \frac{\datasize}{(1 + r)^{2+J+K'}},
		\label{E:TIME0CRUDEPOINTWISEBOUNDFORDATAOFRMINUSUPTO3DERIVATIVESATLEASTONELBAR} \\
	|\partial_t \LRiemann(0,r)| 
	& \lesssim \frac{\datasize}{(1 + r)^2}.
	\label{E:TIME0CRUDEPOINTWISEBOUNDFORDATAOFPARTIALTRMINUS}
\end{align}
\end{subequations}

\noindent \underline{\textbf{Sharp estimates for special combinations involving the Riemann invariants}}.
For $0 \leq J + K \leq 2$, we have:
\begin{subequations}
\begin{align} 
\begin{split} \label{E:TIME0TRCOORDINATESALLDERIVATIVESTRANSPORTEDMODIFIEDULUNITRPLUSDATA}
	\Lunit^J
	\uLunit^K
	\left\lbrace
			r
			\uLunit
			\RRiemann(0,r)
			- 2 \RRiemann(0,r)
		\right\rbrace
		&	= 
			2 \datasize 
				\left\lbrace
				\left(\partial_t + \partial_r \right)^J
				\left( \partial_t - \partial_r \right)^K
				\frac{1}{1 + (t-r)^2}
				\right\rbrace \restriction_{t = 0}
				\\
		& \ \
			+
			\mathcal{O}(\datasize^{1^+}) \frac{1}{(1 + r)^{2+J+K}},
\end{split}
	\\
\begin{split}  \label{E:TIME0TRCOORDINATESPARTIALTDERIVATIVETRANSPORTEDMODIFIEDULUNITRPLUSDATA}
		\left\lbrace
			r
			\uLunit
			\partial_t \RRiemann(0,r)
			- 2 \partial_t \RRiemann(0,r)
		\right\rbrace
		&	= 
			2 \datasize 
				\left\lbrace
				\partial_t
				\left( \frac{1}{1 + (t-r)^2} \right)
				\right\rbrace \restriction_{t= 0}
				\\
		& \ \
			+
			\mathcal{O}(\datasize^{1^+}) \frac{1}{(1 + r)^3}.
\end{split}
\end{align}
\end{subequations}

For $J=0,1,2$, we have:
\begin{subequations}
\begin{align} 
\begin{split} \label{E:TIME0TRCOORDINATESPURELDERIVATIVESTRANSPORTEDMODIFIEDLUNITRMINUSDATA}
	\Lunit^J
	\left\lbrace
			r
			\Lunit \LRiemann(0,r)
			+ 
			2 \LRiemann(0,r)
		\right\rbrace 
		& = 
			- 2 [\datasize + \mathcal{O}(\datasize^{1^+})]
				\left\lbrace
				\left(\partial_t + \partial_r \right)^J
				\frac{1}{1 + (t + r)^2} 
				\right\rbrace
				\restriction_{t=0}.
\end{split}
\end{align}

For $0 \leq J + K' \leq 2$, we have:
\begin{align} 
\begin{split} \label{E:TIME0TRCOORDINATESATLEASTONELBARDERIVATIVETRANSPORTEDMODIFIEDLUNITRMINUSDATA}
	\Lunit^J
	\uLunit^{K'}
	\uLunit
	\left\lbrace
			r
			\Lunit \LRiemann(0,r)
			+ 
			2 \LRiemann(0,r)
		\right\rbrace 
		& = 
			\mathcal{O}(\datasize^{1^+}) \frac{1}{(1 + r)^{3+J+K'}}.
\end{split}
\end{align}
\end{subequations}

Finally, we have:
\begin{align}
\begin{split} \label{E:TIME0TRCOORDINATESPARTIALDDERIVATIVETRANSPORTEDMODIFIEDLUNITRMINUSDATA}
	\left\lbrace
			r
			\Lunit \partial_t \LRiemann(0,r)
			+ 
			2 \partial_t \LRiemann(0,r)
		\right\rbrace 
		& = 
			- 2 
				[\datasize + \mathcal{O}(\datasize^{1^+})]
				\left\lbrace
					\partial_t
				 \left( \frac{1}{1 + (t + r)^2} \right)
				\right\rbrace
				\restriction|_{t=0}
					\\
		& \ \
		+ 
		\mathcal{O}(\datasize^{1^+}) \frac{1}{(1 + r)^3}.
\end{split}
\end{align}

\end{lemma}

\begin{proof}
Throughout, we will silently use \eqref{E:SPEEDOFSOUNDEXPANSION}--\eqref{E:SPEEDOFSOUNDERRORFUNCTIONVANISHESATORIGIN},
which effectively allow us to replace $\Speed$ with $1$ up to harmless error terms.

\eqref{E:TIME0CRUDEPOINTWISEBOUNDFORDATAOFRPLUSUPTO3DERIVATIVES}--\eqref{E:TIME0CRUDEPOINTWISEBOUNDFORDATAOFPARTIALTRMINUS}
are straightforward consequences of
equations \eqref{E:OUTGOINGRIEMANNINVARIANTEVOLUTION}--\eqref{E:INGOINGRIEMANNINVARIANTEVOLUTION}
and the data-assumptions \eqref{E:RPLUSDATAISCLOSETOBACKGROUNDATTIME0}--\eqref{E:RMINUSDATAISCLOSETOBACKGROUNDATTIME0}. 

The remaining estimates in the lemma are more subtle in that they rely on more properties of the profile. 
We refer the reader to Sect.\,\ref{SSS:MOTIVATIONVIALINEARIZEDSYSTEM} for motivation for these estimates.
All the proofs are similar, so we will only prove 
\eqref{E:TIME0TRCOORDINATESALLDERIVATIVESTRANSPORTEDMODIFIEDULUNITRPLUSDATA} in the case $J=K=0$.
To proceed, we use equation \eqref{E:NULLVECTORFIELDS} to deduce:
\begin{align} \label{E:ULUNITRRIEMANNINTERMSOFPARTIALRRRIEMANNANDLUNITRRIEMANN}
	\uLunit \RRiemann
	& =
		- 2 \partial_r \RRiemann
		- 2 (\Speed - 1) \partial_r \RRiemann
		+
		\Lunit \RRiemann.
\end{align}
Using 
\eqref{E:BACKGROUNDSPEEDOFSOUNDISUNITYNORMALIZATION},
equation \eqref{E:OUTGOINGRIEMANNINVARIANTEVOLUTION},
and the assumptions \eqref{E:RPLUSDATAISCLOSETOBACKGROUNDATTIME0}--\eqref{E:RMINUSDATAISCLOSETOBACKGROUNDATTIME0}
to control the terms on $r \times \mbox{RHS~\eqref{E:ULUNITRRIEMANNINTERMSOFPARTIALRRRIEMANNANDLUNITRRIEMANN}}$,
we deduce the following pointwise estimate at time $0$:
\begin{align} 
	\begin{split} \label{E:TIME0ESTIMATEULUNITRRIEMANNINTERMSOFPARTIALRRRIEMANNANDLUNITRRIEMANN}
	r \uLunit \RRiemann(0,r)
	& =
		- 2 r \partial_r \RRiemann(0,r)
		+
		\mathcal{O}(\datasize^{1^+}) \frac{1}{1 + r^2}
			\\
		&
		= - 2 \partial_r (r \RRiemann(0,r))
		+ 
		2 \RRiemann(0,r) 
		+
		\mathcal{O}(\datasize^{1^+}) \frac{1}{1 + r^2}.
\end{split}
\end{align}
By \eqref{E:RPLUSDATAISCLOSETOBACKGROUNDATTIME0}, we have that 
$\partial_r (r \RRiemann(0,r)) = \partial_r (r \profileoverr (r)) + \mathcal{O}(\datasize^{1^+}) \frac{1}{1 + r^2}
= - \partial_r (\datasize \arctan(r)) + \mathcal{O}(\datasize^{1^+}) \frac{1}{1 + r^2}
= - \frac{\datasize}{1 + r^2} + \mathcal{O}(\datasize^{1^+}) \frac{1}{1 + r^2}
$.
Inserting this estimate into the first term on RHS~\eqref{E:TIME0ESTIMATEULUNITRRIEMANNINTERMSOFPARTIALRRRIEMANNANDLUNITRRIEMANN},
we conclude
\eqref{E:TIME0TRCOORDINATESALLDERIVATIVESTRANSPORTEDMODIFIEDULUNITRPLUSDATA} in the case $J=K=0$.

\end{proof}

\subsection{``Data'' induced by the solution at time \texorpdfstring{$\tstar$}{t-star}}
\label{SS:DATAATTIMETSTAR}
For the rest of the paper, $\datasize > 0$ is a parameter that controls the small size of the amplitude of the data as well as the size of the neighborhood
around the background data for which our main results hold.
Our main results will hold whenever $\datasize$ is sufficiently small.
The interesting part of the analysis occurs for times much later than $\tstar$, where:
\begin{align} \label{E:TSTARDEF}
	\tstar
	& := \frac{1}{\datasize}.
\end{align}

In the rest of the paper, 
\begin{align} \label{E:SIGMATDEF}
	\Sigma_t 
\end{align}
denotes the flat hypersurface of constant Cartesian time $t$ in $(t,r)$-coordinate space, 
intersected with the physical region $\lbrace r \geq 0 \rbrace$.

\subsubsection{Estimates relative to the $(t,r)$ coordinates}
\label{SSS:DATAATTIMETSTARRELATIVETOTRCOORDINATES}
In the next lemma, we provide estimates relative to the $(t,r)$ coordinates for the state of the solution at time $\tstar$,
long before any singularity has formed. 

\begin{lemma}[Estimates at time $\tstar$ with respect to the $(t,r)$ coordinates]
\label{L:TIMETSTARDATAESTIMATESINTRCOORDINATES}
Under the assumptions 
\eqref{E:RPLUSDATAISCLOSETOBACKGROUNDATTIME0}--\eqref{E:RMINUSDATAISCLOSETOBACKGROUNDATTIME0},
if $\datasize$ is sufficiently small, then 
the solution $(\RRiemann,\LRiemann)$ to 
\eqref{E:OUTGOINGRIEMANNINVARIANTEVOLUTION}--\eqref{E:INGOINGRIEMANNINVARIANTEVOLUTION}
obeys the following estimates a time $\tstar := \frac{1}{\datasize}$:

\medskip

\noindent \underline{\textbf{Crude estimates for the Riemann invariants}}.
For $0 \leq J + K \leq 3$, we have:
\begin{subequations}
\begin{align} \label{E:CRUDEPOINTWISEBOUNDFORDATAOFRPLUSUPTO3DERIVATIVES}
	|\Lunit^J \uLunit^K \RRiemann(\tstar,r)| 
	& \lesssim \frac{\datasize}{(1 + \tstar + r)^{1+J}(1 + |\tstar - r|)^K},
		\\
	|\partial_t \RRiemann(\tstar,r)| 
	& \lesssim \frac{\datasize}{(1 + \tstar + r)(1 + |\tstar - r|)}.
	\label{E:CRUDEPOINTWISEBOUNDFORDATAOFPARTIALTRPLUS}
\end{align}
\end{subequations}

For $0 \leq J \leq 3$ and $0 \leq J + K' \leq 2$, we have:
\begin{subequations}
\begin{align} \label{E:CRUDEPOINTWISEBOUNDFORDATAOFRMINUSUPTO3LDERIVATIVES}
	|\Lunit^J \LRiemann(\tstar,r)| 
	& \lesssim \frac{\datasize}{(1 + \tstar + r)^{1+J}},
		\\
	|\Lunit^J \uLunit^{K'} \uLunit \LRiemann(\tstar,r)| 
	& \lesssim \frac{\datasize}{(1 + \tstar + r)^{2+J}(1 + |\tstar - r|)^{K'}},
		\label{E:CRUDEPOINTWISEBOUNDFORDATAOFRMINUSUPTO3DERIVATIVESATLEASTONELBAR} \\
	|\partial_t \LRiemann(\tstar,r)| 
	& \lesssim \frac{\datasize}{(1 + \tstar + r)^2}.
	\label{E:CRUDEPOINTWISEBOUNDFORDATAOFPARTIALTRMINUS}
\end{align}
\end{subequations}

\noindent \underline{\textbf{Sharp estimates for special combinations involving the Riemann invariants}}.
For $0 \leq J + K \leq 2$, we have:
\begin{subequations}
\begin{align} 
\begin{split} \label{E:TRCOORDINATESALLDERIVATIVESTRANSPORTEDMODIFIEDULUNITRPLUSTSTARDATA}
	\Lunit^J
	\uLunit^K
	\left\lbrace
			r
			\uLunit
			\RRiemann
			- 2 \RRiemann
		\right\rbrace \restriction_{\Sigma_{\tstar}}
		&	= 
			2 \datasize 
				\left\lbrace
				\left(\partial_t + \partial_r \right)^J
				\left( \partial_t - \partial_r \right)^K
				\frac{1}{1 + (t-r)^2}
				\right\rbrace \restriction_{t= \tstar}
				\\
		& \ \
			+
			\mathcal{O}\left(\frac{\datasize^{1^+}}{(1 + \tstar + r)^J(1 + |\tstar - r|)^{2+K}} \right),
\end{split}
	\\
\begin{split}  \label{E:TRCOORDINATESPARTIALTDERIVATIVETRANSPORTEDMODIFIEDULUNITRPLUSTSTARDATA}
		\left\lbrace
			r
			\uLunit
			\partial_t \RRiemann
			- 2 \partial_t \RRiemann
		\right\rbrace \restriction_{\Sigma_{\tstar}}
		&	= 
			2 \datasize 
				\left\lbrace
				\partial_t
				\left( \frac{1}{1 + (t-r)^2} \right)
				\right\rbrace \restriction_{t= \tstar}
				\\
		& \ \
			+
			\mathcal{O}\left(\frac{\datasize^{1^+}}{(1 + |\tstar - r|)^3} \right).
\end{split}
\end{align}
\end{subequations}

For $J=0,1,2$, we have:
\begin{subequations}
\begin{align} 
\begin{split} \label{E:TRCOORDINATESPURELDERIVATIVESTRANSPORTEDMODIFIEDLUNITRMINUSTSTARDATA}
	\Lunit^J
	\left\lbrace
			r
			\Lunit \LRiemann
			+ 
			2 \LRiemann
		\right\rbrace \restriction_{\Sigma_{\tstar}} 
		& = 
			- 2 [\datasize + \mathcal{O}(\datasize^{1^+})]
				\left\lbrace
				\left(\partial_t + \partial_r \right)^J
				\frac{1}{1 + (t + r)^2} 
				\right\rbrace
				\restriction|_{t=\tstar}.
\end{split}
\end{align}

For $0 \leq J + K' \leq 2$, we have:
\begin{align} 
\begin{split} \label{E:TRCOORDINATESATLEASTONELBARDERIVATIVETRANSPORTEDMODIFIEDLUNITRMINUSTSTARDATA}
	\Lunit^J
	\uLunit^{K'}
	\uLunit
	\left\lbrace
			r
			\Lunit \LRiemann
			+ 
			2 \LRiemann
		\right\rbrace \restriction_{\Sigma_{\tstar}} 
		& = 
			\mathcal{O}\left(\frac{\datasize^2}{(1 + \tstar + r)^{2+J}(1 + |\tstar - r|)^{1+K'}} \right).
\end{split}
\end{align}
\end{subequations}

Finally, we have:
\begin{align}
\begin{split} \label{E:TRCOORDINATESPARTIALDDERIVATIVETRANSPORTEDMODIFIEDLUNITRMINUSTSTARDATA}
	\left\lbrace
			r
			\Lunit \partial_t \LRiemann
			+ 
			2 \partial_t \LRiemann
		\right\rbrace \restriction_{\Sigma_{\tstar}} 
		& = 
			- 2 
				[\datasize + \mathcal{O}(\datasize^{1^+})]
				\left\lbrace
					\partial_t
				 \left( \frac{1}{1 + (t + r)^2} \right)
				\right\rbrace
				\restriction|_{t=\tstar}
					\\
		& \ \
		+ 
		\mathcal{O}\left(\frac{\datasize^2}{(1 + \tstar + r)^2(1 + |\tstar - r|)} \right).
\end{split}
\end{align}

\end{lemma}

\begin{proof}[Discussion of proof]
For the main ideas on how to prove the lemma,
we refer the reader to the proofs of
Props.\,\ref{P:GLOBALEXISTENCERIEMANNINVARIANTSAPRIORIEXTERIORREGIONESTIMATES},
\ref{P:RIEMANNINVARIANTSAPRIORIEXTERIORREGIONESTIMATES}, 
\ref{P:APRIORIESTIMATESININTERIORREGION},
and \ref{P:GLOBALEXISTENCESHARPESTIMATESFORMU}
where similar but much more difficult estimates are derived.
Those propositions are much more difficult to prove because there we have to control the solution up to its classical lifespan,
which is of size $\gtrsim \exp(\frac{C}{\datasize})$ in the case of
Props.\,\ref{P:GLOBALEXISTENCERIEMANNINVARIANTSAPRIORIEXTERIORREGIONESTIMATES}
and
\ref{P:RIEMANNINVARIANTSAPRIORIEXTERIORREGIONESTIMATES}, 
and infinite in the case of Props.\,\ref{P:APRIORIESTIMATESININTERIORREGION},
and \ref{P:GLOBALEXISTENCESHARPESTIMATESFORMU}.
In contrast, for Lemma~\ref{L:TIMETSTARDATAESTIMATESINTRCOORDINATES},
we only need to control the solution for $\tstar = \frac{1}{\datasize}$ amounts of time,
a tiny portion of the lifespan in which the nonlinearities are negligible.
Finally, we remark that as before, \eqref{E:SPEEDOFSOUNDEXPANSION}--\eqref{E:SPEEDOFSOUNDERRORFUNCTIONVANISHESATORIGIN}
allow us to replace $\Speed$ with $1$ in the estimates up to harmless error terms.
\end{proof}

\begin{remark}
Lemma~\ref{L:TIMETSTARDATAESTIMATESINTRCOORDINATES} could also be proved by introducing a fluid potential, which would satisfy a nonlinear wave equation, and 
then deriving weighted energy estimates for the wave equation using the method of commuting vectorfields. 
Such a proof would require different (i.e., Sobolev-type) regularity assumptions on the data.
\end{remark}

\subsubsection{Estimates relative to the $(t,u)$ coordinates}
\label{SSS:DATAATTIMETSTARRELATIVETOTUCOORDINATES}
The next lemma is an analog of 
Lemma~\ref{L:TIMETSTARDATAESTIMATESINTRCOORDINATES} 
in which we provide estimates relative to the $(t,u)$ coordinates for the state of the solution at time $\tstar$.

\begin{lemma}[Estimates in $\Sigma_{\tstar} \cap \lbrace u \leq \frac{\tstar}{2} \rbrace$ with respect to the $(t,u)$-coordinates]
\label{L:EXTERIORREGIONDATAESTIMATESTIMETSTAR}
Under the assumptions 
\eqref{E:RPLUSDATAISCLOSETOBACKGROUNDATTIME0}--\eqref{E:RMINUSDATAISCLOSETOBACKGROUNDATTIME0},
if $\datasize$ is sufficiently small, then  
the solution $(\RRiemann,\LRiemann)$ to 
\eqref{E:OUTGOINGRIEMANNINVARIANTEVOLUTION}--\eqref{E:INGOINGRIEMANNINVARIANTEVOLUTION}
obeys the following estimates in $\Sigma_{\tstar} \cap \lbrace r \geq \frac{\tstar}{2} \rbrace = 
\Sigma_{\tstar} \cap \lbrace u \leq \frac{\tstar}{2} \rbrace$,
where $\tstar := \frac{1}{\datasize}$:

\medskip

\noindent \underline{\textbf{Crude estimates for the Riemann invariants}}.
For $0 \leq J + K \leq 3$, we have:
\begin{align} \label{E:ALLDERIVATIVESRPLUSTSTARDATA}
	\left|
	\Lunit^J
	\muX^K
	\RRiemann 
	\right|
	\restriction_{\Sigma_{\tstar} \cap \lbrace u \leq \frac{\tstar}{2} \rbrace},
		\,
		\left|
	\Lunit^J
	\muuLunit^K
	\RRiemann 
	\right|
	& \leq C \frac{\datasize}{(1 + \tstar + |u|)^{1+J}(1+|u|)^K}.
\end{align}

For $0 \leq J \leq 3$, we have:
\begin{align} \label{E:ALLLUNITDERIVATIVESRMINUSTSTARDATA}
	\left|
	\Lunit^J
	\LRiemann 
	\right|
	\restriction_{\Sigma_{\tstar} \cap \lbrace u \leq \frac{\tstar}{2} \rbrace}
	& \leq C \frac{\datasize}{(1 + \tstar + |u|)^{1+J}}.
\end{align}

For $0 \leq J + K' \leq 2$, we have:
\begin{align} \label{E:ATLEASTONEMUXDERIVATIVERMINUSTSTARDATA}
	\left|
	\Lunit^J
	\muX^{K'}
	\muX
	\LRiemann 
	\right|
	\restriction_{\Sigma_{\tstar} \cap \lbrace u \leq \frac{\tstar}{2} \rbrace},
		\,
	\left|
	\Lunit^J
	\muuLunit^{K'}
	\muuLunit
	\LRiemann 
	\right|
	\restriction_{\Sigma_{\tstar} \cap \lbrace u \leq \frac{\tstar}{2} \rbrace}
	& \leq C \frac{\datasize}{(1 + \tstar + |u|)^{2+J}(1+|u|)^{K'}}.
\end{align}

\medskip

\noindent \underline{\textbf{Estimates for special combinations involving the Riemann invariants}}.
For $0 \leq J + K \leq 2$, we have:
\begin{align} 
\begin{split} \label{E:ALLDERIVATIVESTRANSPORTEDMODIFIEDULUNITRPLUSTSTARDATA}
	\Lunit^J
	\muX^K
	\left\lbrace
			r
			\uLunit
			\RRiemann
			- 2 \RRiemann
		\right\rbrace \restriction_{\Sigma_{\tstar} \cap \lbrace u \leq \frac{\tstar}{2} \rbrace}
		&	= 
			2 \datasize 
				\left( \frac{\partial}{\partial t} \right)^J
				\left( \frac{\partial}{\partial u} \right)^K
				\frac{1}{1 + u^2}
				\\
		& \ \
			+
			\mathcal{O}\left(\frac{\datasize^{1^+}}{(1 + \tstar + |u|)^J(1 + |u|)^{2+K}} \right),
\end{split}
\end{align}

\begin{align} 
\begin{split} \label{E:ALLDERIVATIVESTRANSPORTEDMODIFIEDLUNITRMINUSTSTARDATA}
	\Lunit^J
	\muX^K
	\left\lbrace
			r
			\Lunit \LRiemann
			+ 
			2 \LRiemann
		\right\rbrace \restriction_{\Sigma_{\tstar} \cap \lbrace u \leq \frac{\tstar}{2} \rbrace}
		& = - 
				2 
				[\datasize + \mathcal{O}(\datasize^{1^+})]
				\left\lbrace
				\left( \frac{\partial}{\partial t} \right)^J
				\left( \frac{\partial}{\partial u} \right)^K
				\frac{1}{1 + (2t - u)^2} 
				\right\rbrace
				\restriction|_{t=\tstar}
					\\
		& \ \
		+ 
		\mathcal{O}\left(\frac{\datasize^{1^+}}{(1 + \tstar + |u|)^{2+J}(1 + |u|)^K} \right).
\end{split}
\end{align}

\medskip

\noindent \underline{\textbf{Data estimates for the inverse foliation density}}.

\begin{align} \label{E:MUPOINTWISEATTIMETSTAR}
	\upmu \restriction_{\Sigma_{\tstar} \cap \lbrace u \leq \frac{\tstar}{2} \rbrace}
	 &= \frac{2}{-2 \Speed \partial_r u}=\Speed^{-1}
	=1+\mathcal{O}\left(\frac{\datasize}{(1+\tstar+|u|)} \right),
		\\
	|\muX \upmu|
	\restriction_{\Sigma_{\tstar} \cap \lbrace u \leq \frac{\tstar}{2} \rbrace}
	& \lesssim \frac{\datasize}{(1+\tstar+|u|)(1 + |u|)},
		\label{E:XBREVEMUPOINTWISEATTIMETSTAR} \\
	|\muX \muX \upmu|
	\restriction_{\Sigma_{\tstar} \cap \lbrace u \leq \frac{\tstar}{2} \rbrace}
	&	\lesssim \frac{\datasize}{(1+\tstar+|u|)(1 + |u|)^2}.
		\label{E:TWOXBREVEMUPOINTWISEATTIMETSTAR} 
\end{align}

\end{lemma}

\begin{proof}[Proof outline]
These estimates are straightforward to derive based on
Lemma~\ref{L:TIMETSTARDATAESTIMATESINTRCOORDINATES}.
One only needs to estimate $u$ and $\upmu$ at time $\tstar$
via the data assumption in \eqref{E:ODEFOREXTERIOREIKONALFUNCTION},
\eqref{E:MGHDIDENTITYFORMUATTIMEZERO},
the normalization condition \eqref{E:BACKGROUNDSPEEDOFSOUNDISUNITYNORMALIZATION},
and the evolution equations \eqref{E:ODEFOREXTERIOREIKONALFUNCTION} and
\eqref{E:MUEVOLUTION},
and to translate statements in $(t,r)$ coordinates into statements in $(t,u)$ coordinates
using the change of variables map \eqref{E:CHOVFROMTUTOTRCOORDINATES}, whose Jacobian matrix
is \eqref{E:CHOVJACOBIAN}.
Finally, we remark that as before, \eqref{E:SPEEDOFSOUNDEXPANSION}--\eqref{E:SPEEDOFSOUNDERRORFUNCTIONVANISHESATORIGIN}
allow us to replace $\Speed$ with $1$ in the estimates up to harmless error terms.

\end{proof}

\section{Bootstrap assumptions on a globally hyperbolic exterior region}
\label{S:EXTERIORGLOBALLYHYPERBOLICBOOTSTRAPASSUMPTIONS}

\subsection{Globally hyperbolic exterior bootstrap region}
We will prove our main Exterior Region results with the help of a bootstrap argument on globally hyperbolic subsets of spacetime.
The next definition captures the basic structure of the globally hyperbolic subsets that we encounter in our analysis.
See the proof of Prop.\,\ref{P:EXISTENCEUPTOCREASEANDSINGULARBOUNDARYANDPORTIONOFCAUCHYHORIZON}
for the specific globally hyperbolic subsets that play a role in our bootstrap argument.

\begin{definition}[Globally hyperbolic exterior bootstrap region]
	\label{D:GLOBALLHYPERBOLICEXTERIORBOOTSTRAPREGION}
	A subset $\mathbf{GH}_{\textnormal{Boot}}^{{\textnormal{Ext}}}$ of spacetime is said to be a \textbf{globally hyperbolic exterior bootstrap region} if:
	\begin{align} \label{E:TANDUBOUNDSFOREXTERIORBOOTSTRAPREGION}
		\mathbf{GH}_{\textnormal{Boot}}^{{\textnormal{Ext}}} 
		& \subset 
		\lbrace t \geq \tstar \rbrace
		\cap
		\lbrace u \leq \tstar/2 \rbrace,
	\end{align}
	and if for every point $q \in \mathbf{GH}_{\textnormal{Boot}}^{{\textnormal{Ext}}}$, 
	the past-directed integral curves of $\Lunit$ and $\uLunit$
	emanating from $q$ remain in the interior of $\mathbf{GH}_{\textnormal{Boot}}^{{\textnormal{Ext}}}$ 
	until they intersect a point $q'$ satisfying
	$q' \in \mathbf{GH}_{\textnormal{Boot}}^{{\textnormal{Ext}}} \bigcap \left( \Sigma_{\tstar} \cap \lbrace u \leq \tstar/2 \rbrace \right)$. 
\end{definition}

\begin{remark}[Connection to the standard notion of global hyperbolicity]
	\label{R:CONNECTIONTOSTANDARDGLOBALHYPERBOLICITY}
	Our definitions of ``globally hyperbolic'' are not standard but can be shown to be equivalent to standard ones from Lorentzian geometry.
	We have adopted our definitions because they are simple to state and sufficient for our PDE analysis.
	In Appendix~\ref{A:GLOBALHYPERBOLICITYANDMGHDS}, we prove that our maximal globally hyperbolic development
	is globally hyperbolic in the standard sense of Lorentzian geometry; see Prop.\,\ref{P:EQUIVALENCEOFGLOBALHYPERBOLICITYMMAX}.
\end{remark}

\subsection{Bootstrap assumptions on \texorpdfstring{$\mathbf{GH}_{\textnormal{Boot}}^{{\textnormal{Ext}}}$}{an Exterior Region}}
\label{SS:BA}
As is standard for evolution quasilinear problems, we find it convenient to derive estimates via a bootstrap argument. 
To this end, we fix a globally hyperbolic exterior bootstrap region $\mathbf{GH}_{\textnormal{Boot}}^{{\textnormal{Ext}}}$,  
and we assume that the following inequalities hold for some number $\eps$ satisfying:
\begin{align} \label{E:BOUNDONBOOTSTRAPPARAMETERSIZE}
	0 & < \eps \leq \datasize^{3/4}.
\end{align}

\medskip 
\noindent \underline{\textbf{Bootstrap assumptions for $\RRiemann$}}.
For $0 \leq J \leq 3$, we assume:
\begin{align} \label{E:BAEXTERIORREGIONPOINTWISEESTIMATEFORRPLUSANDDERIVATIVES}
	|\Lunit^J \RRiemann|
	& 
	\leq \frac{\eps}{(1 + t + |u|)^{1+J}}.
\end{align}

For $0 \leq J + K' \leq 2$, we assume:
\begin{align} \label{E:BAEXTERIORREGIONPOINTWISEESTIMATEFORRPLUSATLEASONEXDERIVATIVE}
	|\Lunit^J \muX^{K'} \muX \RRiemann|
	& 
	\leq 
	\frac{\eps}{(1 + t + |u|)^{2+J}(1 + |u|)^{K'}} 
	+ 
	\frac{\eps}{(1 + t + |u|)^{1+J}(1 + |u|)^{2+K'}}.
\end{align}

\medskip 
\noindent \underline{\textbf{Bootstrap assumptions for $\LRiemann$}}.
For $0 \leq J \leq 3$, we assume:
\begin{subequations}
\begin{align} \label{E:BAEXTERIORREGIONPOINTWISEESTIMATEFORRMINUSANDLDERIVATIVES}
	|\Lunit^J \LRiemann|
	& 
	\leq \frac{\eps}{(1 + t + |u|)^{1+J}}.
\end{align}

For $0 \leq J + K' \leq 2$, we assume:
\begin{align} \label{E:BAEXTERIORREGIONPOINTWISEESTIMATEFORRMINUSATLEASTONEXDERIVATIVE}
	|\Lunit^J \muX^{K'} \muX \LRiemann|
	& 
	\leq \frac{\eps}{(1 + t + |u|)^{2+J}(1 + |u|)^{K'}}.
\end{align}
\end{subequations}

\medskip 
\noindent \underline{\textbf{Bootstrap assumptions for $\upmu$}}.
We assume:

\begin{subequations}
\begin{align} \label{E:BAMUNOTTOOBIG}
	&
	0
	<
	\upmu
	\leq 2,
		\\
	|\muX \upmu|
	& \leq \frac{10}{1 + |u|},
		\label{E:BAMUXMUNOTTOOBIG} \\
	|\muX \muX \upmu|
	& \leq \frac{10}{(1 + |u|)^2}.
		\label{E:BATWICEMUXMUNOTTOOBIG}
\end{align}
\end{subequations}

\section{Geometric a priori estimates in the exterior region based on the bootstrap assumptions}
\label{S:GEOMETRICEXTERIORREGIONESTIMATES}
In this section, we use the bootstrap assumptions to derive our main a priori estimates in the Exterior Region.

\subsection{Sharp estimates for the inverse foliation density}
\label{SS:SHARPESTIMATESFORMU}
We start with the following proposition, which yields sharp estimates for $\upmu$. These results will be essential for controlling the shape 
of the future-boundary of the MGHD of the data. 
Such sharp information is fundamentally important for our proof that the MGHD exists and is unique.

\begin{proposition}[Sharp estimates for the eikonal function and the inverse foliation density]
\label{P:SHARPESTIMATESFORMU}
Under the data assumptions of Section~\ref{S:DATA} and the bootstrap assumptions of
Section~\ref{S:EXTERIORGLOBALLYHYPERBOLICBOOTSTRAPASSUMPTIONS}, 
if $\datasize$ is sufficiently small, then the following estimates hold on $\mathbf{GH}_{\textnormal{Boot}}^{{\textnormal{Ext}}}$:

\noindent \underline{\textbf{Estimates for the eikonal function}}.
The eikonal function $u$, which is the solution to the transport initial value problem \eqref{E:ODEFOREXTERIOREIKONALFUNCTION},
satisfies the following estimates:
\begin{align} \label{E:SHARPPOINWISECOMPARISONBETWEENUANDTMINUSR}
	u 
	& = t - r 
	+
	\mathcal{O}(\eps)
	\ln \left( \frac{1 + t + |u|}{1 + \tstar + |u|} \right).
\end{align}

\medskip
\noindent \underline{\textbf{Comparison estimates for various coordinate functions}}.
The Cartesian coordinate functions obey the following estimates:
\begin{subequations}
	\begin{align} \label{E:RISAPPROXIMATELYTMINUSU}
		r & = [1 + \mathcal{O}(\eps)](t - u),
			\\
		1 + t
		& \approx
		\frac{t}{2} 
		\leq t - u 
		\leq t + |u|
		\approx 1 + t + |u|
		\approx 2t - u.
		\label{E:COMPARISONBETWEENTOVERUTMINUSUANDTPLUSMODU}
	\end{align}
\end{subequations}

\medskip

\noindent  \underline{\textbf{Estimates for the Riemann invariants}}.
The Riemann invariants satisfy the following estimates:

\begin{subequations}
\begin{align} \label{E:EXTERIORREGIONPOINTWISEESTIMATEFORMODIFIEDLBARRPLUSDERIVATIVE}
		r
		\muuLunit \RRiemann
		- 
		2 
		\upmu \RRiemann
		& =
		\frac{2 \datasize + \mathcal{O}(\datasize^{1^+})}{1 + u^2},	
			\\
		\muX
		\left\lbrace
		r
		\muuLunit \RRiemann
		- 
		2 
		\upmu \RRiemann
		\right\rbrace
		& =
		-
		\frac{4 \datasize u}{(1 + u^2)^2}
		+
		\frac{\mathcal{O}(\datasize^{1^+}) }{(1 + |u|)^3},
			\label{E:MUXDIFFERENTIATEDEXTERIORREGIONPOINTWISEESTIMATEFORMODIFIEDLBARRPLUSDERIVATIVE}
				\\
	\muX \muX
		\left\lbrace
		r
		\muuLunit \RRiemann
		- 
		2 
		\upmu \RRiemann
		\right\rbrace
		& =
		4 
		\frac{\datasize (3u^2 - 1)}{(1 + u^2)^3}
		+
		\frac{\mathcal{O}(\datasize^{1^+}) }{(1 + |u|)^4},
			\label{E:TWICEMUXDIFFERENTIATEDEXTERIORREGIONPOINTWISEESTIMATEFORMODIFIEDLBARRPLUSDERIVATIVE}
\end{align}
\end{subequations}

\begin{subequations}
\begin{align} \label{E:EXTERIORREGIONPOINTWISEESTIMATEFORMULBARRPLUS}
		\muuLunit \RRiemann
		& =
		\frac{2 \datasize + \mathcal{O}(\datasize^{1^+})}{(t-u) (1 + u^2)}
		+
		\mathcal{O}(\eps)
		\frac{1}{(1 + t + |u|)^2},
			\\
		\muX \RRiemann
		& =
		\frac{\datasize + \mathcal{O}(\datasize^{1^+})}{(t-u) (1 + u^2)}
		+
		\mathcal{O}(\eps)
		\frac{1}{(1 + t + |u|)^2},
		\label{E:EXTERIORREGIONPOINTWISEESTIMATEFORMUXRPLUS}
			\\
		\muX \muuLunit \RRiemann
		& =
		-
		\frac{4 \datasize u}{ (t-u)(1 + u^2)^2}
		+
		\frac{\mathcal{O}(\datasize^{1^+}) }{(t-u)(1 + |u|)^3}
		+
		\frac{\mathcal{O}(\eps)}{(1 + t + |u|)^2 (1 + |u|)},
			\label{E:EXTERIORREGIONPOINTWISEESTIMATEFORMUXMULBARRPLUS} \\
		\muX \muX \RRiemann
		& =
		-
		\frac{2 \datasize u}{ (t-u)(1 + u^2)^2}
		+
		\frac{\mathcal{O}(\datasize^{1^+}) }{(t-u)(1 + |u|)^3}
		+
		\frac{\mathcal{O}(\eps)}{(1 + t + |u|)^2 (1 + |u|)},
		\label{E:EXTERIORREGIONPOINTWISEESTIMATEFORMUXMUXRPLUS}
			\\
		\muX \muX \muuLunit \RRiemann
		& =
		4 
		\frac{\datasize (3u^2 - 1)}{(t-u)(1 + u^2)^3}
		+
		\frac{\mathcal{O}(\datasize^{1^+})}{(t-u)(1 + |u|)^4}
		+
		\frac{\mathcal{O}(\eps)}{(1 + t + |u|)^2 (1 + |u|)^2},
			\label{E:EXTERIORREGIONPOINTWISEESTIMATEFORMUXMUXMULBARRPLUS}  \\
	  \muX \muX \muX \RRiemann
		& =
		2 
		\frac{\datasize (3u^2 - 1)}{ (t-u)(1 + u^2)^3}
		+
		\frac{\mathcal{O}(\datasize^{1^+}) }{(t-u)(1 + |u|)^4}
		+
		\frac{\mathcal{O}(\eps)}{(1 + t + |u|)^2 (1 + |u|)^2}.
		\label{E:EXTERIORREGIONPOINTWISEESTIMATEFORMUXMUXMUXRPLUS}
\end{align}
\end{subequations}

\noindent \underline{\textbf{Estimates for the inverse foliation density}}.
Let $\lifespanconstant > 0$ be the constant from \eqref{E:NULLCONDITIONFAILURECONSTANT}--\eqref{E:LIFESPANCONSTANT}.
Then the following estimates hold:

\begin{subequations}
\begin{align} \label{E:LUNITMUPOINTWISEESTIMATEEXTERIOR}
	\Lunit \upmu
	& =
		- \frac{\lifespanconstant \datasize + \mathcal{O}(\datasize^{1^+})}{(t-u)(1 + u^2)}
		+ 
		\frac{\mathcal{O}(\eps)}{(1 + t + |u|)^2},
			\\
	\Lunit \muX \upmu
	& =
		\frac{2 \lifespanconstant \datasize u}{(t-u)(1 + u^2)^2}
		+
		\frac{\mathcal{O}(\datasize^{1^+}) }{(t-u)(1 + |u|)^3}
		+ 
		\frac{\mathcal{O}(\eps)}{(1 + t + |u|)^2 (1 + |u|)},
			\label{E:LUNITMUXMUPOINTWISEESTIMATEEXTERIOR} \\
	\Lunit \muX \muX \upmu
	& =
			2 
			\lifespanconstant \datasize
			\left\lbrace
			\frac{1 - 3u^2}{(t-u)(1 + u^2)^3} 
			\right\rbrace
		+
		\frac{\mathcal{O}(\datasize^{1^+}) }{(t-u)(1 + |u|)^4}
		+ 
		\frac{\mathcal{O}(\eps)}{(1 + t + |u|)^2 (1 + |u|)^2},
		\label{E:LUNITTWICEMUXMUPOINTWISEESTIMATEEXTERIOR} 
\end{align}
\end{subequations}

\begin{subequations}
\begin{align} \label{E:MUPOINTWISEESTIMATEEXTERIOR}
	\upmu
	& = 1 
			+
			\frac{\mathcal{O}(\eps)}{(1 + \tstar + |u|)}
			-
			\frac{\lifespanconstant \datasize + \mathcal{O}(\datasize^{1^+})}{(1 + u^2)} \ln \left( \frac{t-u}{\tstar - u} \right),
				\\
	\muX \upmu
	& = \frac{\mathcal{O}(\eps)}{(1 + \tstar + |u|)(1+|u|)}
			+
			\frac{2 \lifespanconstant \datasize u + \mathcal{O}(\datasize^{1^+})(1 + |u|)}{(1 + u^2)^2} \ln \left( \frac{t-u}{\tstar - u} \right),
			 \label{E:XBREVEMUPOINTWISEESTIMATEEXTERIOR}
				\\
\muX \muX \upmu
	& = \frac{\mathcal{O}(\eps)}{(1 + \tstar + |u|)(1+|u|)^2}
			+
			\left\lbrace
				\frac{2 \lifespanconstant \datasize(1 - 3u^2) + \mathcal{O}(\datasize^{1^+})(1 + u^2)}{(1 + u^2)^3} 
			\right\rbrace
			\ln \left( \frac{t-u}{\tstar - u} \right).
			 \label{E:TWOXBREVEMUPOINTWISEESTIMATEEXTERIOR}
	\end{align}
\end{subequations}

Moreover,
\begin{subequations}
\begin{align}
	\upmu  \label{E:MUSIMPLERUPPERBOUND}
	& \leq 1 
			+
			C \frac{\eps}{(1 + \tstar + |u|)},
				\\
	|\muX \upmu|
	& \leq  C \frac{\eps}{(1 + \tstar + |u|)(1+|u|)}
				+
				2
				\frac{|u|}{1 + u^2}
				+
				C \frac{\datasize^{1^+}}{\datasize}
				\frac{1}{1 + |u|},
	\label{E:MUXMUSIMPLERUPPERBOUND}		
		\\
	|\muX \muX \upmu|
	& \leq  C \frac{\eps}{(1 + \tstar + |u|)(1+|u|)^2}
				+
				2
				\frac{|1 - 3u^2|}{(1 + u^2)^2}
				+
				C \frac{\datasize^{1^+}}{\datasize}
				\frac{1}{(1 + |u|)^2}.
	\label{E:MUXMUXMUSIMPLERUPPERBOUND}	
\end{align}
\end{subequations}

\medskip

\noindent \underline{\textbf{Sharpened estimates in the region where the inverse foliation density is small}}.
In $\lbrace 0 < \upmu \leq \frac{1}{200} \rbrace \cap \mathbf{GH}_{\textnormal{Boot}}^{{\textnormal{Ext}}}$, we have:
\begin{align} \label{E:QUANTITATIVEESTIMATEFORTINREGIONWHEREMUISSMALL}
	\exp
	\left(
	\frac{\frac{99}{100} (1 + u^2)}{\lifespanconstant \datasize}
	\right)
	\leq
	\frac{t-u}{\tstar - u} 
	&
	\leq
	\exp
	\left(
	\frac{\frac{101}{100}(1 + u^2)}{\lifespanconstant \datasize}
	\right),
		\\
	\frac{99}{100} 
	& \leq
	\frac{\lifespanconstant \datasize}{(1 + u^2)}
	 \ln \left( \frac{t-u}{\tstar - u} \right)
	\leq \frac{101}{100}, 
	\label{E:QUANTITATIVEESTIMATEFORLOGTINREGIONWHEREMUISSMALL}
\end{align}

\begin{align} \label{E:QUANTITATIVELUNITMUNEGATIVEWHENMUISSMALL}
	-
	\frac{\left\lbrace \lifespanconstant \datasize + \mathcal{O}(\datasize^{1^+}) \right\rbrace }{(1 + u^2)(\tstar-u)}
	\exp
	\left(
	- \frac{\frac{99}{100}(1 + u^2)}{\lifespanconstant \datasize}
	\right)
	& 
	\leq 
	\Lunit \upmu
	\leq
	-
	\frac{\left\lbrace \lifespanconstant \datasize + \mathcal{O}(\datasize^{1^+}) \right\rbrace }{(1 + u^2)(\tstar-u)}
	\exp
	\left(
	- \frac{\frac{101}{100}(1 + u^2)}{\lifespanconstant \datasize}
	\right),
\end{align}
and:
\begin{align} \label{E:QUANTITATIVEMUXRPLUSUPPERANDLOWERBOUNDSWHENMUISSMALL}
	\frac{\left\lbrace \datasize + \mathcal{O}(\datasize^{1^+}) \right\rbrace }{(1 + u^2)(\tstar-u)}
	\exp
	\left(
	- \frac{\frac{101}{100}(1 + u^2)}{\lifespanconstant \datasize}
	\right)
	& 
	\leq 
	\muX \RRiemann
	\leq
	\frac{\left\lbrace \datasize + \mathcal{O}(\datasize^{1^+}) \right\rbrace }{(1 + u^2)(\tstar-u)}
	\exp
	\left(
	- \frac{\frac{99}{100}(1 + u^2)}{\lifespanconstant \datasize}
	\right).
\end{align}

In $\lbrace 0 < \upmu \leq \frac{1}{200} \rbrace \cap \lbrace u \leq 0 \rbrace \cap \mathbf{GH}_{\textnormal{Boot}}^{{\textnormal{Ext}}}$, we have:
\begin{subequations}
\begin{align} \label{E:MUXMUSHARPBOUNDWHENMUISSMALLANDUISNEGATIVE}
		  \frac{101}{50}
			\frac{u + \mylittleo(\datasize)}{(1 + u^2)} 
		&
		\leq
		\muX \upmu
		\leq \frac{99}{50}
			\frac{u + \mylittleo(\datasize)}{(1 + u^2)},
				\\
			\frac{101}{25}
			\frac{u + \mylittleo(\datasize)}{(1 + u^2)}
		&
		\leq
		\muuLunit \upmu
		\leq \frac{99}{25}
			\frac{u + \mylittleo(\datasize)}{(1 + u^2)}. 
			\label{E:MUULUNITMUSHARPBOUNDWHENMUISSMALLANDUISNEGATIVE}
\end{align}
\end{subequations}

In $\lbrace 0 < \upmu \leq \frac{1}{200} \rbrace \cap \lbrace 0 \leq u \leq \frac{\tstar}{2} \rbrace \cap \mathbf{GH}_{\textnormal{Boot}}^{{\textnormal{Ext}}}$, we have:
\begin{subequations}
\begin{align} \label{E:MUXMUSHARPBOUNDWHENMUISSMALLANDUISPOSITIVE}
		  \frac{99}{50}
			\frac{u + \mylittleo(\datasize)}{(1 + u^2)} 
		&
		\leq
		\muX \upmu
		\leq 
			\frac{101}{50}
			\frac{u + \mylittleo(\datasize)}{(1 + u^2)},
				\\
	 \frac{99}{25}
			\frac{u + \mylittleo(\datasize)}{(1 + u^2)}
		&
		\leq
		\muuLunit \upmu
		\leq 
			\frac{101}{25}
			\frac{u + \mylittleo(\datasize)}{(1 + u^2)}.
			\label{E:MUULUNITMUSHARPBOUNDWHENMUISSMALLANDUISPOSITIVE}
\end{align}
\end{subequations}

In $\lbrace 0 < \upmu \leq \frac{1}{200} \rbrace \cap \lbrace u^2 \leq \frac{1}{3} \rbrace \cap \mathbf{GH}_{\textnormal{Boot}}^{{\textnormal{Ext}}}$, we have:
\begin{subequations}
\begin{align}
		  \frac{99}{50}
			\frac{1 - 3u^2}{(1 + u^2)^2}
			+
			\frac{\mylittleo(\datasize)}{(1 + \tstar + |u|)(1+|u|)^2}
		&
		\leq
		\muX \muX \upmu
		\leq
			\frac{101}{50}
			\frac{1 - 3u^2}{(1 + u^2)^2}
			+
			\frac{\mylittleo(\datasize)}{(1 + \tstar + |u|)(1+|u|)^2},
				\label{E:BEHAVIOROFMUXMUXMUWHENMUISSMALLANDUSQUAREDISSMALL} \\
	\frac{99}{25}
			\frac{1 - 3u^2}{(1 + u^2)^2}
			+
			\frac{\mylittleo(\datasize)}{(1 + \tstar + |u|)(1+|u|)^2}
		&
		\leq
		\muX \muuLunit \upmu
		\leq
			\frac{101}{25}
			\frac{1 - 3u^2}{(1 + u^2)^2}
			+
			\frac{\mylittleo(\datasize)}{(1 + \tstar + |u|)(1+|u|)^2},
				\label{E:BEHAVIOROFMUXMUULUNITMUWHENMUISSMALLANDUSQUAREDISSMALL} \\
			7
			\frac{1 - 3u^2}{(1 + u^2)^2}
			+
			\frac{\mylittleo(\datasize)}{(1 + \tstar + |u|)(1+|u|)^2}
		&
		\leq
		\muuLunit \muuLunit \upmu
		\leq
			8
			\frac{1 - 3u^2}{(1 + u^2)^2}
			+
			\frac{\mylittleo(\datasize)}{(1 + \tstar + |u|)(1+|u|)^2}.
			\label{E:BEHAVIOROFMUULUNITMUULUNITMUWHENMUISSMALLANDUSQUAREDISSMALL}
\end{align}
\end{subequations}

In $\lbrace 0 < \upmu \leq \frac{1}{200} \rbrace \cap \lbrace u^2 \geq \frac{1}{3} \rbrace \cap \mathbf{GH}_{\textnormal{Boot}}^{{\textnormal{Ext}}}$, we have:
\begin{subequations}
\begin{align} \label{E:BEHAVIOROFMUXMUXMUWHENMUISSMALLANDUSQUAREDISLARGE}
		  \frac{101}{50}
			\frac{1 - 3u^2}{(1 + u^2)^2}
			+
			\frac{\mylittleo(\datasize)}{(1 + \tstar + |u|)(1+|u|)^2}
		&
		\leq
		\muX \muX \upmu
		\leq 
		\frac{99}{50}
		\frac{1 - 3u^2}{(1 + u^2)^2} 
		+
		\frac{\mylittleo(\datasize)}{(1 + \tstar + |u|)(1+|u|)^2},
			\\
			  \frac{101}{25}
			\frac{1 - 3u^2}{(1 + u^2)^2}
			+
			\frac{\mylittleo(\datasize)}{(1 + \tstar + |u|)(1+|u|)^2}
		&
		\leq
		\muX \muuLunit \upmu
		\leq 
		\frac{99}{25}
		\frac{1 - 3u^2}{(1 + u^2)^2} 
		+
		\frac{\mylittleo(\datasize)}{(1 + \tstar + |u|)(1+|u|)^2},
			 \label{E:BEHAVIOROFMUULUNITMUXMUWHENMUISSMALLANDUSQUAREDISLARGE} \\
		8
		\frac{1 - 3u^2}{(1 + u^2)^2}
			+
			\frac{\mylittleo(\datasize)}{(1 + \tstar + |u|)(1+|u|)^2}
		&
		\leq
		\muuLunit \muuLunit \upmu
		\leq 
		7
		\frac{1 - 3u^2}{(1 + u^2)^2} 
		+
		\frac{\mylittleo(\datasize)}{(1 + \tstar + |u|)(1+|u|)^2}.
		 \label{E:BEHAVIOROFMUULUNITMUULUNITMUWHENMUISSMALLANDUSQUAREDISLARGE}
\end{align}
\end{subequations}

\end{proposition}

\begin{remark}[Sign of the radiation field]
	\label{R:SIGNOFRADIATIONFIELD}
	LHS~\eqref{E:EXTERIORREGIONPOINTWISEESTIMATEFORMODIFIEDLBARRPLUSDERIVATIVE} is the future-radiation field of the solution.
	Note that it is everywhere \emph{positive} and decays away from the wave zone. This is the main mechanism driving our shock-forming MGHD results.
	Similar results also hold in the context of our global existence result, Theorem~\ref{T:MAINGLOBALEXISTENCETHEOREM}.
	Specifically, LHS~\eqref{E:GLOBALEXISTENCEPOINTWISESPECIALCOMBINATIONRMUULUNITRPLUS} is the future-radiation field, and it is everywhere \textbf{negative},
	leading to global compression.
\end{remark}	

\begin{proof}[Proof of Prop.\,\ref{P:SHARPESTIMATESFORMU}]
\ \\

\noindent \textbf{Proof of \eqref{E:SHARPPOINWISECOMPARISONBETWEENUANDTMINUSR}}:
The identity $\Lunit r = \Speed + v^r$ and the bootstrap assumptions imply that
$\Lunit r = 1 + \frac{\mathcal{O}(\eps)}{(1 + t + |u|)}$. Also using that
$\Lunit t = 1$ and $\Lunit u = 0$, we deduce that:
\begin{align} \label{E:POINTWISEESTIMATEFORLDERIVATIVEOFUMINUSTMINUSR}
	\Lunit [u - (t-r)] 
	& = \frac{\mathcal{O}(\eps)}{(1 + t + |u|)}.
\end{align}
Recalling that $\Lunit = \frac{\partial}{\partial t}$ in geometric coordinates,
we integrate 
	\eqref{E:EXTERIORREGIONPOINTWISEESTIMATEFORLDERIVATIMEOFMODIFIEDLBARRPLUSDERIVATIVE}
	from time $\tstar$ to time $t$ and use the data assumption
in \eqref{E:ODEFOREXTERIOREIKONALFUNCTION}
to conclude \eqref{E:SHARPPOINWISECOMPARISONBETWEENUANDTMINUSR}.

\medskip \noindent \textbf{Proof of \eqref{E:RISAPPROXIMATELYTMINUSU}--\eqref{E:COMPARISONBETWEENTOVERUTMINUSUANDTPLUSMODU}}:
These two inequalities follow easily from
\eqref{E:SHARPPOINWISECOMPARISONBETWEENUANDTMINUSR} and our assumption 
\eqref{E:TANDUBOUNDSFOREXTERIORBOOTSTRAPREGION}
that $t \geq \tstar$ and $u \leq \frac{\tstar}{2}$

In the rest of the proof, we will silently use \eqref{E:RISAPPROXIMATELYTMINUSU}--\eqref{E:COMPARISONBETWEENTOVERUTMINUSUANDTPLUSMODU}.

\medskip

\noindent \textbf{Proof of \eqref{E:EXTERIORREGIONPOINTWISEESTIMATEFORMODIFIEDLBARRPLUSDERIVATIVE}--\eqref{E:TWICEMUXDIFFERENTIATEDEXTERIORREGIONPOINTWISEESTIMATEFORMODIFIEDLBARRPLUSDERIVATIVE}}:
We first use the wave equation \eqref{E:REVAMPEDRPLUSWAVEEQUATION},
Lemma~\ref{L:VECTORFIELDSINTERMSOFGEOMETRICCOORDINATES},
the transport equation \eqref{E:MUEVOLUTION},
	and the bootstrap assumptions
	to deduce that for $K=0,1,2$, we have:
  \begin{align} \label{E:EXTERIORREGIONPOINTWISEESTIMATEFORLDERIVATIMEOFMODIFIEDLBARRPLUSDERIVATIVE}
		\left|
		\Lunit 
		\muX^K
		\left\lbrace
			r
			\muuLunit \RRiemann
			- 2 \upmu \RRiemann
		\right\rbrace
		\right|
		& \lesssim 
		\eps^2 \frac{1}{(1 + t + |u|)^2(1 + |u|)^{2+K}}
		+ 
		\eps^2 \frac{1}{(1 + t + |u|)^3(1 + |u|)^K}.
	\end{align}
	We clarify that the most difficult term on RHS~\eqref{E:REVAMPEDRPLUSWAVEEQUATION} is the boxed one,
	which by equation \eqref{E:MUEVOLUTION} leads to a term of type $\RRiemann \cdot \muuLunit \RRiemann$
	that we bound via \eqref{E:BAEXTERIORREGIONPOINTWISEESTIMATEFORRPLUSANDDERIVATIVES} 
	and \eqref{E:BAEXTERIORREGIONPOINTWISEESTIMATEFORRPLUSATLEASONEXDERIVATIVE}.
	Recalling that $\Lunit = \frac{\partial}{\partial t}$ in geometric coordinates, we integrate 
	\eqref{E:EXTERIORREGIONPOINTWISEESTIMATEFORLDERIVATIMEOFMODIFIEDLBARRPLUSDERIVATIVE}
	from time $\tstar$ to time $t$ and use the data bounds 
	in \eqref{E:ALLDERIVATIVESTRANSPORTEDMODIFIEDULUNITRPLUSTSTARDATA}
	to conclude
	\eqref{E:EXTERIORREGIONPOINTWISEESTIMATEFORMODIFIEDLBARRPLUSDERIVATIVE}--\eqref{E:TWICEMUXDIFFERENTIATEDEXTERIORREGIONPOINTWISEESTIMATEFORMODIFIEDLBARRPLUSDERIVATIVE}.

\medskip

\noindent \textbf{Proof of \eqref{E:EXTERIORREGIONPOINTWISEESTIMATEFORMULBARRPLUS}--\eqref{E:EXTERIORREGIONPOINTWISEESTIMATEFORMUXMUXMUXRPLUS}}:
These estimates follow from
\eqref{E:EXTERIORREGIONPOINTWISEESTIMATEFORMODIFIEDLBARRPLUSDERIVATIVE}--\eqref{E:TWICEMUXDIFFERENTIATEDEXTERIORREGIONPOINTWISEESTIMATEFORMODIFIEDLBARRPLUSDERIVATIVE},
\eqref{E:SHARPPOINWISECOMPARISONBETWEENUANDTMINUSR}, 
the identity $\muX =\frac{1}{2}(\muuLunit - \upmu \Lunit)$ (which follows from Lemma~\ref{L:VECTORFIELDSINTERMSOFGEOMETRICCOORDINATES}), 
the identities $\Lunit r = v^r+\Speed$ and $\uLunit r = v^r - \Speed$ (which follow from \eqref{E:NULLVECTORFIELDS})
and the bootstrap assumptions.

\medskip
\noindent \textbf{Proof of \eqref{E:LUNITMUPOINTWISEESTIMATEEXTERIOR}--\eqref{E:LUNITTWICEMUXMUPOINTWISEESTIMATEEXTERIOR}}:
We use the bootstrap assumptions and the already proven estimates
	\eqref{E:EXTERIORREGIONPOINTWISEESTIMATEFORMODIFIEDLBARRPLUSDERIVATIVE}--\eqref{E:EXTERIORREGIONPOINTWISEESTIMATEFORMUXMUXMUXRPLUS} 
	to bound the terms on the RHS
	of the transport equation \eqref{E:MUEVOLUTION} for $\upmu$ as well as its first and second derivatives with respect to $\muX$.
	We clarify that the main term in the transport equation is the one 
	$-
			\frac{1}{2 \Speed}
			\left\lbrace
				1
				+
				\Speed' 
				\InverseRiemannfunction'
			\right\rbrace
			\muuLunit 
			\RRiemann
	$
	that we highlighted in Remark~\ref{R:MAINTERMINANALYSIS},
	and that by Taylor expanding the nonlinearities and using
	\eqref{E:SPEEDOFSOUNDEXPANSION}--\eqref{E:SPEEDOFSOUNDERRORFUNCTIONVANISHESATORIGIN},
	definition \eqref{E:NULLCONDITIONFAILURECONSTANT},
	and the bootstrap assumptions, we can write this main term as
	$
	- \frac{1}{2} \lifespanconstant 
	\muuLunit 
	\RRiemann
	$
	plus a small error term. 
	The term 
	$
	- \frac{1}{2} \lifespanconstant 
	\muuLunit 
	\RRiemann
	$
	and its derivatives lead to the $\lifespanconstant$-dependent terms on
	RHSs~\eqref{E:LUNITMUPOINTWISEESTIMATEEXTERIOR}--\eqref{E:LUNITTWICEMUXMUPOINTWISEESTIMATEEXTERIOR},
	while the small error terms are part of the $\mathcal{O}(\datasize^{1^+})$-multiplied
	error terms on RHSs~\eqref{E:LUNITMUPOINTWISEESTIMATEEXTERIOR}--\eqref{E:LUNITTWICEMUXMUPOINTWISEESTIMATEEXTERIOR}.
	
\medskip
\noindent \textbf{Proof of \eqref{E:MUPOINTWISEESTIMATEEXTERIOR}--\eqref{E:TWOXBREVEMUPOINTWISEESTIMATEEXTERIOR}}:
	Recall that $\Lunit = \frac{\partial}{\partial t}$ in geometric coordinates.
	Hence, to prove \eqref{E:MUPOINTWISEESTIMATEEXTERIOR}, we can integrate 
	\eqref{E:LUNITMUPOINTWISEESTIMATEEXTERIOR}
	from time $\tstar$ to time $t$
	and use the data bound \eqref{E:MUPOINTWISEATTIMETSTAR}.
	
	The estimates \eqref{E:XBREVEMUPOINTWISEESTIMATEEXTERIOR}--\eqref{E:TWOXBREVEMUPOINTWISEESTIMATEEXTERIOR}
	follow similarly from \eqref{E:LUNITMUXMUPOINTWISEESTIMATEEXTERIOR}--\eqref{E:LUNITTWICEMUXMUPOINTWISEESTIMATEEXTERIOR}
	and the data bounds \eqref{E:XBREVEMUPOINTWISEATTIMETSTAR}--\eqref{E:TWOXBREVEMUPOINTWISEATTIMETSTAR}.

\medskip
\noindent \textbf{Proof of \eqref{E:MUSIMPLERUPPERBOUND}--\eqref{E:MUXMUXMUSIMPLERUPPERBOUND}}:
\eqref{E:MUSIMPLERUPPERBOUND} follows directly from the fact that the 
last product $-
			\frac{\lifespanconstant \datasize + \mathcal{O}(\datasize^{1^+})}{(1 + u^2)} \ln \left( \frac{t-u}{\tstar - u} \right)$
 on RHS~\eqref{E:MUPOINTWISEESTIMATEEXTERIOR} is negative.

To prove \eqref{E:MUXMUSIMPLERUPPERBOUND}, we first use \eqref{E:MUPOINTWISEESTIMATEEXTERIOR} and our bootstrap assumption that $\upmu > 0$
to deduce that everywhere in $\mathbf{GH}_{\textnormal{Boot}}^{{\textnormal{Ext}}}$, we have:
\begin{align} \label{E:TIMEBOUNDINREGIONWHEREMUISPOSITIVE}
	\frac{\lifespanconstant \datasize + \mathcal{O}(\datasize^{1^+})}{(1 + u^2)} \ln \left( \frac{t-u}{\tstar - u} \right)
	& \leq 1 + C \eps.
\end{align}
Using the bound \eqref{E:TIMEBOUNDINREGIONWHEREMUISPOSITIVE} to control the last product on RHS~\eqref{E:XBREVEMUPOINTWISEESTIMATEEXTERIOR}, 
we arrive at \eqref{E:MUXMUSIMPLERUPPERBOUND}. The estimate \eqref{E:MUXMUXMUSIMPLERUPPERBOUND} follows from similar reasoning based on
\eqref{E:TIMEBOUNDINREGIONWHEREMUISPOSITIVE} and \eqref{E:TWOXBREVEMUPOINTWISEESTIMATEEXTERIOR}.

\medskip
\noindent \textbf{Proof of \eqref{E:QUANTITATIVEESTIMATEFORTINREGIONWHEREMUISSMALL}--\eqref{E:QUANTITATIVEESTIMATEFORLOGTINREGIONWHEREMUISSMALL}}:
These estimates are algebraically equivalent, so it suffices to prove
\eqref{E:QUANTITATIVEESTIMATEFORLOGTINREGIONWHEREMUISSMALL}.
This bound follows easily from the estimate \eqref{E:MUPOINTWISEESTIMATEEXTERIOR}
and the assumption that $0 < \upmu \leq \frac{1}{200}$.

\medskip
\noindent \textbf{Proof of \eqref{E:QUANTITATIVELUNITMUNEGATIVEWHENMUISSMALL} and \eqref{E:QUANTITATIVEMUXRPLUSUPPERANDLOWERBOUNDSWHENMUISSMALL}}:
\eqref{E:QUANTITATIVELUNITMUNEGATIVEWHENMUISSMALL} follows from \eqref{E:LUNITMUPOINTWISEESTIMATEEXTERIOR} and \eqref{E:QUANTITATIVEESTIMATEFORTINREGIONWHEREMUISSMALL}, where we use that the latter estimate and \eqref{E:TANDUBOUNDSFOREXTERIORBOOTSTRAPREGION} 
imply that $t-u\geq (\tstar-u)\exp\kh{\frac{\frac{99}{100}(1+u^2)}{\lifespanconstant\datasize}}\geq \frac{1}{2}\tstar \frac{\frac{99}{100}(1+u^2)}{\lifespanconstant\datasize}$ and hence $\frac{\eps}{(1+t+|u|)^2}\leq C\frac{\eps\datasize^2}{(t-u)(1+u^2)}$.
Similarly, \eqref{E:QUANTITATIVEMUXRPLUSUPPERANDLOWERBOUNDSWHENMUISSMALL} follows from 
\eqref{E:EXTERIORREGIONPOINTWISEESTIMATEFORMUXRPLUS} and 
\eqref{E:QUANTITATIVEESTIMATEFORTINREGIONWHEREMUISSMALL}.


\medskip
\noindent \textbf{Proof of \eqref{E:MUXMUSHARPBOUNDWHENMUISSMALLANDUISNEGATIVE}--\eqref{E:BEHAVIOROFMUULUNITMUULUNITMUWHENMUISSMALLANDUSQUAREDISLARGE}}:
\eqref{E:MUXMUSHARPBOUNDWHENMUISSMALLANDUISNEGATIVE}
and \eqref{E:MUXMUSHARPBOUNDWHENMUISSMALLANDUISPOSITIVE}
follow from 
\eqref{E:XBREVEMUPOINTWISEESTIMATEEXTERIOR}
and \eqref{E:QUANTITATIVEESTIMATEFORLOGTINREGIONWHEREMUISSMALL}.

\eqref{E:MUULUNITMUSHARPBOUNDWHENMUISSMALLANDUISNEGATIVE} and \eqref{E:MUULUNITMUSHARPBOUNDWHENMUISSMALLANDUISPOSITIVE}
then follow from 	
\eqref{E:MUXMUSHARPBOUNDWHENMUISSMALLANDUISNEGATIVE}
and \eqref{E:MUXMUSHARPBOUNDWHENMUISSMALLANDUISPOSITIVE},	
\eqref{E:QUANTITATIVELUNITMUNEGATIVEWHENMUISSMALL},
and the identity $\muuLunit = \upmu \Lunit + 2 \muX$ 
(which follows from Lemma~\ref{L:VECTORFIELDSINTERMSOFGEOMETRICCOORDINATES}).	

The estimates 
\eqref{E:BEHAVIOROFMUXMUXMUWHENMUISSMALLANDUSQUAREDISSMALL}--\eqref{E:BEHAVIOROFMUULUNITMUULUNITMUWHENMUISSMALLANDUSQUAREDISLARGE}
follow from a similar argument that also relies on \eqref{E:TWOXBREVEMUPOINTWISEESTIMATEEXTERIOR}.	
	
\end{proof}

\subsection{Inequalities along the integral curves of \texorpdfstring{$\muuLunit$}{the ingoing null vectorfield}}
\label{SS:INEQUALITIESALONGINTEGRALCURVES}
To prove our main a priori estimates in the Exterior Region, 
we will rely on some general estimates along the integral curves of 
$\muuLunit$, which we provide in the next lemma.

\begin{lemma}[Inequalities along the integral curves of $\muuLunit$]
	\label{L:INEQUALITIESALONGINTEGRALCURVES}
	Let $u \rightarrow (\mathfrak{t}_{u_0}(u),u)$ be the $u$-parameterized integral curve of 
	$\muuLunit = \upmu \frac{\partial}{\partial t} + 2 \frac{\partial}{\partial u}$ emanating from the point
	$(\tstar,u_0) \in \Sigma_{\tstar}$, i.e., 
	in view of Lemma~\ref{L:VECTORFIELDSINTERMSOFGEOMETRICCOORDINATES},
	$\mathfrak{t}$ is the solution to the following initial value problem:
	\begin{align} \label{E:IVPFORINTEGRALCURVESOFMUULUNITPARAMETERIZEDBYU}
		\frac{d}{d u}
		\mathfrak{t}_{u_0}(u)
		& = \frac{1}{2} \upmu(\mathfrak{t}(u),u),
		&&
		\mathfrak{t}_{u_0}(u_0) 
		= \tstar. 
	\end{align}
	Let $u_1 \geq u_0$, let $t_1 := \mathfrak{t}_{u_0}(u_1)$, let $r_0 := \tstar - u_0$ be the radial value 
	corresponding to the point with geometric coordinates $(\tstar,u_0)$ (see \eqref{E:ODEFOREXTERIOREIKONALFUNCTION}), 
	and let $r_1 := r(t_1,u_1)$ be the corresponding
	radial value at the point with geometric coordinates $(t_1,u_1)$.
	Under the data assumptions of Section~\ref{S:DATA} and the bootstrap assumptions of
	Section~\ref{S:EXTERIORGLOBALLYHYPERBOLICBOOTSTRAPASSUMPTIONS}, 
	if $\datasize$ is sufficiently small, then
	the following estimates hold on $\mathbf{GH}_{\textnormal{Boot}}^{{\textnormal{Ext}}}$.
	
\medskip
	
\noindent \underline{\textbf{Coordinate function estimates along the integral curves}}.
	
	\begin{align} \label{E:SMALLCHANGEINTPLUSRALONGINTEGRALCURVESOFULUNIT}
		t_1 + r_1
		& = \tstar + r_0 + \mathcal{O}(\eps) \ln \left( \frac{1 + t_1}{1 + \tstar} \right),
	\end{align}
	
	\begin{align} 
	\begin{split} \label{E:TWOTMINUSUALONGINTEGRALCURVESOFULUNIT}
		2t_1 - u_1
		& 
		= 2 \tstar - u_0 
		+
		\mathcal{O}(\eps) \ln \left( \frac{1 + t_1}{1 + \tstar} \right)
		+
		\mathcal{O}(\eps)
		\ln\left( \frac{1 + t_1 + |u_1|}{1 + \tstar + |u_1|} \right)
			\\
		& 
		= 2 \tstar - u_0 
		+
		\mathcal{O}(\eps) \ln \left( \frac{1 + t_1}{1 + \tstar} \right).
	\end{split}
	\end{align}

	\noindent \underline{\textbf{Estimates of integrals involving the geometric coordinates}}.
	Let $A \geq 0$ be a constant, 
	and let $\ln_+$ be the function defined by:
	\begin{align} \label{E:LOGPLUS}
		\ln_+(z) 
		& := \ln(e + z).
	\end{align}
	Then there exists a $C > 0$, depending on $A$, such that the following estimates hold:
	\begin{align} \label{E:INTEGRALESTIMATEALONGMUULUNITONEOVERTPLUSUSQUAREDTIMESONEOVER1PLUSU}
	\int_{u_0}^{u_1}
		\frac{1}{(1 + \mathfrak{t}_{u_0}(u) + |u|)^A (1 + |u|)}
	\, \mathrm{d} u
	& \leq 
			C
			\frac{\ln_+(\tstar + |u_0|)}{[1 + (2 \tstar - u_0)]^A},
				\\
	\int_{u_0}^{u_1}
		\frac{\ln_+(t + |u|) }{(1 + \mathfrak{t}_{u_0}(u) + |u|)^A (1 + |u|)^2}
	\, \mathrm{d} u
	& \leq 
			C
			\frac{\ln_+(\tstar + |u_0|)}{[1 + (2 \tstar - u_0)]^A}.
			 \label{E:INTEGRALESTIMATEALONGMUULUNITONEOVERTPLUSUSQUAREDTIMESLOGTPLUSMODUTIMESONEOVER1PLUSUSQUARED}
	\end{align}

\end{lemma}

\begin{proof}
Throughout the proof, we will silently use \eqref{E:RISAPPROXIMATELYTMINUSU}--\eqref{E:COMPARISONBETWEENTOVERUTMINUSUANDTPLUSMODU}.

\medskip

\noindent \textbf{Proof of \eqref{E:SMALLCHANGEINTPLUSRALONGINTEGRALCURVESOFULUNIT}--\eqref{E:TWOTMINUSUALONGINTEGRALCURVESOFULUNIT}}:
To prove \eqref{E:SMALLCHANGEINTPLUSRALONGINTEGRALCURVESOFULUNIT}, we first use \eqref{E:NULLVECTORFIELDS} to compute that
$\uLunit (t+r) = 1 + (v^r - \Speed)$.
From this identity, 
\eqref{E:SPEEDOFSOUNDEXPANSION}--\eqref{E:SPEEDOFSOUNDERRORFUNCTIONVANISHESATORIGIN},
and the bootstrap assumptions, 
we deduce:
\begin{align} \label{E:POINTWISEESTIMATEFORULUNITDERIVATIVEOFTPLUSR}
	|\uLunit (t + r)| 
	& \lesssim \frac{\eps}{1+t}.
\end{align}
Note that $\uLunit t = 1$, i.e., along the integral curves of $\uLunit$, we have $\uLunit = \frac{d}{dt}$.
Hence, we can integrate \eqref{E:POINTWISEESTIMATEFORULUNITDERIVATIVEOFTPLUSR} along the integral curves of $\uLunit$ 
from time $\tstar$ to time $t_1$ 
to conclude \eqref{E:SMALLCHANGEINTPLUSRALONGINTEGRALCURVESOFULUNIT}.

\eqref{E:TWOTMINUSUALONGINTEGRALCURVESOFULUNIT} then follows from	
\eqref{E:SMALLCHANGEINTPLUSRALONGINTEGRALCURVESOFULUNIT} and
\eqref{E:SHARPPOINWISECOMPARISONBETWEENUANDTMINUSR}.	

\medskip
\noindent \textbf{Proof of \eqref{E:INTEGRALESTIMATEALONGMUULUNITONEOVERTPLUSUSQUAREDTIMESONEOVER1PLUSU}--\eqref{E:INTEGRALESTIMATEALONGMUULUNITONEOVERTPLUSUSQUAREDTIMESLOGTPLUSMODUTIMESONEOVER1PLUSUSQUARED}}:
	To prove \eqref{E:INTEGRALESTIMATEALONGMUULUNITONEOVERTPLUSUSQUAREDTIMESONEOVER1PLUSU},
	we use
	\eqref{E:TWOTMINUSUALONGINTEGRALCURVESOFULUNIT},
	the fact that $u$ increases along the integral curves of $\muuLunit$ (since $\muuLunit u = 2$),
	and our assumption that $u \leq \frac{\tstar}{2}$ (see \eqref{E:TANDUBOUNDSFOREXTERIORBOOTSTRAPREGION})
	to deduce that:
	\begin{align} 
	\begin{split} \label{E:FIRSTBOUNDINTEGRALESTIMATEALONGMUULUNITONEOVERTPLUSUSQUAREDTIMESONEOVER1PLUSU}
	\int_{u_0}^{u_1}
		\frac{1}{(1 + \mathfrak{t}_{u_0}(u) + |u|)^A (1 + |u|)}
	\, \mathrm{d} u
	& \lesssim
		\int_{u_0}^{u_1}
			\frac{1}{(1 + 3 \mathfrak{t}_{u_0}(u) - u)^A (1 + |u|)}
		\, \mathrm{d} u
			\\
	& \lesssim
		\frac{1}{(1 + 2 \tstar - u_0)^A}
		\int_{u_0}^{u_1}
			\frac{1}{(1 + |u|)}
		\, \mathrm{d} u
		\\
	& \lesssim 
		\frac{\ln_+(|u_1|) + \ln_+(|u_0|)}{(1 + 2 \tstar - u_0)^A},
		\\
		& 
		\lesssim 
		\frac{\ln_+(\tstar + |u_0|)}{(1 + 2 \tstar - u_0)^A},
\end{split}
\end{align}
as is desired.

The estimate \eqref{E:INTEGRALESTIMATEALONGMUULUNITONEOVERTPLUSUSQUAREDTIMESLOGTPLUSMODUTIMESONEOVER1PLUSUSQUARED} 
can be proved via a similar argument, and we omit the details.

\end{proof}

\subsection{The main a priori estimates for the Riemann invariants}
\label{SS:MAINAPRIORIFORRIEMANNINVARIANTS}
In the next proposition, we derive our main a priori estimates for the Riemann invariants in the exterior region.

\begin{proposition}[A priori estimates and strict improvement of the bootstrap assumptions for the Riemann invariants in the exterior region]
\label{P:RIEMANNINVARIANTSAPRIORIEXTERIORREGIONESTIMATES}
Under the data assumptions of Section~\ref{S:DATA} and the bootstrap assumptions of
Section~\ref{S:EXTERIORGLOBALLYHYPERBOLICBOOTSTRAPASSUMPTIONS}, 
if $\datasize$ is sufficiently small, then
the following estimates hold on $\mathbf{GH}_{\textnormal{Boot}}^{{\textnormal{Ext}}}$:

\medskip 
For $0 \leq J \leq 3$, we have:
\begin{align} \label{E:EXTERIORREGIONPOINTWISEESTIMATEFORRPLUSANDLDERIVATIVES}
	|\Lunit^J \RRiemann|
	& 
	\leq C \frac{\datasize}{(1 + t + |u|)^{1+J}}.
\end{align}

\medskip 
For $0 \leq J + K' \leq 2$, we have:
\begin{align} \label{E:EXTERIORREGIONPOINTWISEESTIMATEFORRPLUSATLEASTONEXDERIVATIVE}
	|\Lunit^J \muX^{K'} \muX  \RRiemann|
	&  
	\leq 
	C \frac{\datasize}{(1 + t + |u|)^{1+J}(1 + |u|)^{2+K'}}
	+
	C \frac{\datasize}{(1 + t + |u|)^{2+J}(1 + |u|)^{K'}}.
\end{align}

For $0 \leq J \leq 3$, we have:
\begin{align} \label{E:EXTERIORREGIONPOINTWISEESTIMATEFORRMINUSANDLDERIVATIVES}
	|\Lunit^J \LRiemann|
	& 
	\leq C \frac{\datasize}{(1 + t + |u|)^{1+J}}.
\end{align}

For $0 \leq J + K' \leq 2$, we have:
\begin{align} \label{E:EXTERIORREGIONPOINTWISEESTIMATEFORRMINUSATLEASTONEXDERIVATIVE}
	|\Lunit^J \muX^{K'} \muX \LRiemann|
	& 
	\leq C \frac{\datasize}{(1 + t + |u|)^{2+J}(1 + |u|)^{K'}}.
\end{align}

\end{proposition}

\begin{remark}[Strict improvement of the bootstrap assumptions, except possibly $\upmu > 0$]
	\label{R:STRICTIMPROVEMENTOFBOOTSTRAP}
	Recall that we assumed only that $0 < \eps \leq \datasize^{3/4}$
	in the bootstrap assumptions on $\mathbf{GH}_{\textnormal{Boot}}^{{\textnormal{Ext}}}$ that we made
	in Section~\ref{S:EXTERIORGLOBALLYHYPERBOLICBOOTSTRAPASSUMPTIONS}; see \eqref{E:BOUNDONBOOTSTRAPPARAMETERSIZE}.
	It is therefore easy to see that when $\datasize$ is sufficiently small, 
	the estimates of Props.\ \ref{P:SHARPESTIMATESFORMU}
	and \ref{P:RIEMANNINVARIANTSAPRIORIEXTERIORREGIONESTIMATES} 
	collectively imply strict improvements of
	the bootstrap assumptions, except it is
	possible that $\upmu$ vanishes on the closure of $\mathbf{GH}_{\textnormal{Boot}}^{{\textnormal{Ext}}}$.
	In the two propositions imply that the $\eps$-involving bootstrap assumptions in fact hold with $C \datasize$
	in place of $\eps$, which is a strict improvement if $\eps := \datasize^{3/4}$
	and $C \datasize \leq \datasize^{3/4}$.
	Until Sect.\,\ref{S:GLOBALEXISTENCE}, 
	we will often silently use this fact in the following manner: 
	all of the $\eps$-involving bootstrap assumptions we have made and the estimates we have proved
	hold on $\mathbf{GH}_{\textnormal{Boot}}^{{\textnormal{Ext}}}$ with 
	$\eps$ replaced by $C \datasize$.
\end{remark}

\begin{remark}
\label{Remark:improveboundsR+-}
Some of the decay rates stated in Prop.\,\ref{P:RIEMANNINVARIANTSAPRIORIEXTERIORREGIONESTIMATES} are not optimal with respect to decay in $u$;
we chose to prove the simplest estimates that allow us to close the bootstrap argument.
To illustrate an example of an estimate that can be improved, 
we consider the solutions $\flatRRiemann$ and $\flatLRiemann$ to the linearized equations from
\eqref{E:LINEARRPLUSSOLUTIONFORMULA}--\eqref{E:LINEARRMINUSSOLUTIONFORMULA}. If one expresses these solutions 
in the ``flat geometric coordinates'' $(t,u^{(\textsf{flat})}): = (t,t-r)$, then from
these solution formulas, it is not too difficult to show the following improved bounds:
\begin{subequations}
\begin{align}
    |\flatRRiemann(t,u^{(\textsf{flat})})| & \leq C\frac{\datasize}{(1+t+|u^{(\textsf{flat})}|)(1+(u^{(\textsf{flat})})_+)},\\
		|\flatLRiemann(t,u^{(\textsf{flat})})|&\leq C\frac{\datasize\kh{1+\ln\frac{1+t+|u^{(\textsf{flat})}|}{1+|u^{(\textsf{flat})}|}+(u^{(\textsf{flat})})_-}}{(1+t+|u^{(\textsf{flat})}|)^2}.
 \end{align}
\end{subequations}
%
Here $u_+=\max\{u,0\}$ and $u_-=\max\{-u,0\}$. Similar improved estimates hold for the derivatives of $(\flatRRiemann,\flatLRiemann)$.
With a bit of additional effort, we could have proved that our nonlinear Riemann invariants $\BothRiemann$ and their derivatives satisfy the same pointwise bounds with respect to $(t,u)$ as the linear solutions and their derivatives do with respect to $(t,u^{(\textsf{flat})})$.
Similar remarks apply to Props.\,\ref{P:APRIORIESTIMATESININTERIORREGION} 
and \ref{P:GLOBALEXISTENCERIEMANNINVARIANTSAPRIORIEXTERIORREGIONESTIMATES}; see 
Remarks~\ref{Remark:improveboundsR+-:interior} and \ref{Remark:improveboundsR+-:global}.
\end{remark}

\begin{proof}[Proof of Proposition~\ref{P:RIEMANNINVARIANTSAPRIORIEXTERIORREGIONESTIMATES}]
Throughout the proof, we will silently use 
\eqref{E:TANDUBOUNDSFOREXTERIORBOOTSTRAPREGION},
\eqref{E:SHARPPOINWISECOMPARISONBETWEENUANDTMINUSR},
and
\eqref{E:RISAPPROXIMATELYTMINUSU}--\eqref{E:COMPARISONBETWEENTOVERUTMINUSUANDTPLUSMODU}.
We will also silently use \eqref{E:SPEEDOFSOUNDEXPANSION}--\eqref{E:SPEEDOFSOUNDERRORFUNCTIONVANISHESATORIGIN},
which allow us to replace $\Speed$ with $1$ in various estimates, up to harmless error terms.

\medskip

\noindent \textbf{Preliminary estimates}:
	We first use the wave equation \eqref{E:REVAMPEDRPLUSWAVEEQUATION}, 
	Lemma~\ref{L:VECTORFIELDSINTERMSOFGEOMETRICCOORDINATES},
	the transport equation \eqref{E:MUEVOLUTION},
	and the bootstrap assumptions
	to deduce that for $0 \leq J + K \leq 2$, we have:
  \begin{align} \label{E:EXTERIORREGIONPOINTWISEESTIMATEFORCOMMUTEDVERSIONLDERIVATIMEOFMODIFIEDLBARRPLUSDERIVATIVE}
		\left|
		\Lunit 
		\Lunit^J
		\muX^K
		\left\lbrace
			r
			\muuLunit \RRiemann
			- 
			2 \upmu \RRiemann
		\right\rbrace
		\right|
		& \lesssim 
		\eps^2 \frac{1}{(1 + t + |u|)^{2+J}(1 + |u|)^{2+K}}
		+ 
		\eps^2 \frac{1}{(1 + t + |u|)^{3+J}(1 + |u|)^K}.
	\end{align}
	We clarify that the most difficult term on RHS~\eqref{E:REVAMPEDRPLUSWAVEEQUATION} is the boxed one,
	which by equation \eqref{E:MUEVOLUTION} leads to a term of type $\RRiemann \cdot \muuLunit \RRiemann$
	that we bound via \eqref{E:BAEXTERIORREGIONPOINTWISEESTIMATEFORRPLUSANDDERIVATIVES} 
	and \eqref{E:BAEXTERIORREGIONPOINTWISEESTIMATEFORRPLUSATLEASONEXDERIVATIVE}.
	We now restrict to the case $J=0$ in \eqref{E:EXTERIORREGIONPOINTWISEESTIMATEFORCOMMUTEDVERSIONLDERIVATIMEOFMODIFIEDLBARRPLUSDERIVATIVE}.
	Recalling that $\Lunit = \frac{\partial}{\partial t}$ in geometric coordinates, we integrate 
	\eqref{E:EXTERIORREGIONPOINTWISEESTIMATEFORCOMMUTEDVERSIONLDERIVATIMEOFMODIFIEDLBARRPLUSDERIVATIVE}
	from time $\tstar$ to time $t$ and use the data bounds in
	\eqref{E:ALLDERIVATIVESTRANSPORTEDMODIFIEDULUNITRPLUSTSTARDATA}
	to deduce that for $0 \leq  K \leq 2$, we have the following preliminary estimate:
	\begin{align} \label{E:POINTWISEBOUNDALLMUXDERIVATIVESTRANSPORTEDMODIFIEDULUNITRPLUSTSTARDATA}
	\begin{split}	
		\left|
		\muX^K
		\left\lbrace
			r
			\muuLunit \RRiemann
			- 
			2 \upmu \RRiemann
		\right\rbrace
		-
		2 \datasize 
				\left( \frac{\partial}{\partial u} \right)^K
				\frac{1}{1 + u^2}
		\right|
		& \lesssim
			\datasize^{1^+} \frac{1}{(1 + |u|)^{2+K}}.
	\end{split}
	\end{align}
	
	Similarly, we use the wave equation \eqref{E:REVAMPEDRMINUSWAVEEQUATION}, 
	Lemma~\ref{L:VECTORFIELDSINTERMSOFGEOMETRICCOORDINATES},
	the transport equation \eqref{E:MUEVOLUTION},
	the bootstrap assumptions,
	and \eqref{E:XBREVEMUPOINTWISEESTIMATEEXTERIOR}--\eqref{E:TWOXBREVEMUPOINTWISEESTIMATEEXTERIOR}
	to deduce that for $0 \leq J + K \leq 2$, we have:
	\begin{align} \label{E:WAVEEQUATIONFORALLDERIVATIVESRMINUSINHOMOGENEOUSTERMBOUND}
	\left|
	\muuLunit
	\Lunit^J
	\muX^K
	\left\lbrace
		r \Lunit \LRiemann
		+
		2 \LRiemann
	\right\rbrace
	\right|
	& \lesssim \eps^2 \frac{1 + [\ln_+(t + |u|)]^K}{(1 + t + |u|)^{2+J}(1+|u|)^{1+K}}.
	\end{align}
	We now restrict to the case $K=0$ in \eqref{E:WAVEEQUATIONFORALLDERIVATIVESRMINUSINHOMOGENEOUSTERMBOUND}.
	Integrating \eqref{E:WAVEEQUATIONFORALLDERIVATIVESRMINUSINHOMOGENEOUSTERMBOUND} along the integral curve of $\muuLunit$ that connects $(t,u)$ to 
	the ``data point''
	$(\tstar,u_0)$, recalling that $\muuLunit  = \upmu \frac{\partial}{\partial t} + 2 \frac{\partial}{\partial u}$,
	using the integral estimate
	\eqref{E:INTEGRALESTIMATEALONGMUULUNITONEOVERTPLUSUSQUAREDTIMESONEOVER1PLUSU},
	the data bound \eqref{E:ALLDERIVATIVESTRANSPORTEDMODIFIEDLUNITRMINUSTSTARDATA},
	\eqref{E:TWOTMINUSUALONGINTEGRALCURVESOFULUNIT},
	and \eqref{E:SHARPPOINWISECOMPARISONBETWEENUANDTMINUSR},
	we find that for $0 \leq J \leq 2$, we have the following preliminary estimate:
	\begin{align} \label{E:EXTERIORREGIONPOINTWISEESTIMATEFORALLLUNITDERIVATIVESMODIFIEDLRMINUSDERIVATIVE}
	\begin{split}
		\Lunit^J
		[r \Lunit \LRiemann
		+
		2 \LRiemann]
		& 
		=
		-
		2 \datasize
		\left\lbrace
			\left( \frac{\partial}{\partial t} \right)^J
			\frac{1}{1 + (2t - u_0)^2} 
			\right\rbrace
		\restriction|_{t=\tstar}
			\\
		& \ \
		+
		\mathcal{O}\left(\frac{\datasize^{1^+}}{(1 + \tstar + |u_0|)^{2+J}} \right)
		+
		\mathcal{O}(\eps^2) \frac{\ln_+(\tstar + |u_0|)}{[1 + (2 \tstar - u_0)]^{2+J}}
			\\
	& =
		-
		2 \datasize 
		\left( \frac{\partial}{\partial t} \right)^J
		\frac{1}{1 + (2t - u)^2}
		+
		\mathcal{O}(\datasize^{1^+}) \frac{\ln_+(t + |u|)}{(1 + t + |u|)^{2+J}}.
	\end{split}
	\end{align}

\medskip

\noindent \textbf{Proof of \eqref{E:EXTERIORREGIONPOINTWISEESTIMATEFORRPLUSANDLDERIVATIVES} in the case $J=0$}:
	We consider equation \eqref{E:POINTWISEBOUNDALLMUXDERIVATIVESTRANSPORTEDMODIFIEDULUNITRPLUSTSTARDATA} in the case $K=0$.
	We multiply this equation by $r$ and use the identity
	$\muuLunit r = \upmu \uLunit r = \upmu (v^r - \Speed)$, 
	\eqref{E:SPEEDOFSOUNDEXPANSION}--\eqref{E:SPEEDOFSOUNDERRORFUNCTIONVANISHESATORIGIN},
	\eqref{E:RISAPPROXIMATELYTMINUSU}, 
	and the bootstrap assumptions to deduce:
	\begin{align}
	\begin{split} \label{E:POINTWISEBOUNDEXTERIORMUULINITDERIVATIVEOFRADIALCOORDINATESQUAREDTIMESRPLUS}
		\muuLunit [r^2 \RRiemann]
		& 
		=
			2  
			\datasize
			\frac{r}{1 + u^2}
		+
		\mathcal{O}\left(\datasize^{1^+} \frac{t - u}{(1 + |u|)^2} \right)
			\\
		& \ \
		+
		\mathcal{O}\left(\frac{\eps (t-u)}{1 + t + |u|} \right) |\RRiemann|
			\\
		& 
		=
			2  
			\left\lbrace
				\datasize
				+
				\mathcal{O}(\datasize^{1^+})
			\right\rbrace
			\frac{(t-u)}{1 + u^2}.
		\end{split}
	\end{align}
	Recalling that $\muuLunit = 2 \frac{\mathrm{d}}{\mathrm{d} u}$ along the integral curves of $\muuLunit$, we 
integrate \eqref{E:POINTWISEBOUNDEXTERIORMUULINITDERIVATIVEOFRADIALCOORDINATESQUAREDTIMESRPLUS}
with respect to $u$ starting from the data point $(\tstar,u_0)$ and terminating at the final point $(t_1,u_1)$,
use that $t$ increases along the integral curves of $\muuLunit$,
and use the data estimate \eqref{E:CRUDEPOINTWISEBOUNDFORDATAOFRPLUSUPTO3DERIVATIVES} and \eqref{E:SHARPPOINWISECOMPARISONBETWEENUANDTMINUSR}
to deduce:
\begin{align} 
\begin{split} \label{E:ALMOSTFINALPOINTWISEESTIMATEFORRADIALCOORDINATESQUAREDTIMESRPLUS}
	|[r^2 \RRiemann](t_1,u_1)|
	& = |[r^2 \RRiemann](\tstar,u_0)|
	+
	\mathcal{O}(\datasize) t_1
	+
	\mathcal{O}(\datasize)
	\ln_+(|u_0| + |u_1|)
		\\
& = \mathcal{O}(\datasize) (\tstar-u_0)
	+
	\mathcal{O}(\datasize) t_1
	+
	\mathcal{O}(\datasize)
	\ln_+(|u_0| + |u_1|).
\end{split}
\end{align}
Dividing \eqref{E:ALMOSTFINALPOINTWISEESTIMATEFORRADIALCOORDINATESQUAREDTIMESRPLUS} by $r^2(t_1,u_1)$
and using
\eqref{E:SHARPPOINWISECOMPARISONBETWEENUANDTMINUSR},
\eqref{E:RISAPPROXIMATELYTMINUSU}--\eqref{E:COMPARISONBETWEENTOVERUTMINUSUANDTPLUSMODU},
and \eqref{E:TWOTMINUSUALONGINTEGRALCURVESOFULUNIT},
we conclude \eqref{E:EXTERIORREGIONPOINTWISEESTIMATEFORRPLUSANDLDERIVATIVES} in the case $J=0$.

\medskip
	\noindent \textbf{Proof of \eqref{E:EXTERIORREGIONPOINTWISEESTIMATEFORRMINUSANDLDERIVATIVES}}:
	We consider the case $J=0$ in \eqref{E:EXTERIORREGIONPOINTWISEESTIMATEFORALLLUNITDERIVATIVESMODIFIEDLRMINUSDERIVATIVE}.
	Using \eqref{E:RISAPPROXIMATELYTMINUSU} and the bootstrap assumptions to replace $r$ in 
	\eqref{E:EXTERIORREGIONPOINTWISEESTIMATEFORALLLUNITDERIVATIVESMODIFIEDLRMINUSDERIVATIVE} by $t-u$ up to a small error, 
	multiplying the resulting equation by $(t-u)$, 
	and using that $\Lunit (t-u) = 1$,
	we find that:
	\begin{align}
	\begin{split} \label{E:POINTWISEBOUNDFORLUNITTMINUSUSQUAREDRMINUSSQUARED}
		\Lunit [(t-u)^2 \LRiemann]
		& =
		-
		2 \frac{\datasize(t-u)}{1 + (2t - u)^2}
		+
		\mathcal{O}(\datasize^{1^+}) \frac{\ln_+(t + |u|)}{(1 + t + |u|)}.
	\end{split}
	\end{align}
	Recalling that $\Lunit = \frac{\partial}{\partial t}$ in geometric coordinates, we integrate 
	\eqref{E:POINTWISEBOUNDFORLUNITTMINUSUSQUAREDRMINUSSQUARED}
	with respect to time from $\tstar$ to $t$ and use the data estimate
	\eqref{E:ALLLUNITDERIVATIVESRMINUSTSTARDATA}
	to deduce:
	\begin{align}
	\begin{split}	\label{E:EXTERIORREGIONMAINPOINTWISEESTIMATEFORTMINUSSQUAREDTIMESRMINUS}
		|(t-u)^2 \LRiemann(t,u)|
		& \lesssim 
				\datasize
				(1 + t + |u|)
				+
				\datasize^{1^+} \frac{\tstar - u}{1 + |u|}
				+
				\datasize^{1^+}[\ln_+\left( 1 + t + |u| \right)]^2.
	\end{split}
	\end{align}
	The desired bound \eqref{E:EXTERIORREGIONPOINTWISEESTIMATEFORRMINUSANDLDERIVATIVES} in the case $J=0$
	now follows from
	\eqref{E:EXTERIORREGIONMAINPOINTWISEESTIMATEFORTMINUSSQUAREDTIMESRMINUS}.
	
	The bound \eqref{E:EXTERIORREGIONPOINTWISEESTIMATEFORRMINUSANDLDERIVATIVES} in the remaining cases $J=1,2,3$ 
	then follows inductively
	from \eqref{E:EXTERIORREGIONPOINTWISEESTIMATEFORALLLUNITDERIVATIVESMODIFIEDLRMINUSDERIVATIVE},
	the bootstrap assumptions,
	and the identity $\Lunit r = v^r + \Speed = 1 + \frac{\mathcal{O}(\eps)}{1 + t + |u|}$.

\medskip

\noindent \textbf{Proof of \eqref{E:EXTERIORREGIONPOINTWISEESTIMATEFORRPLUSANDLDERIVATIVES} in the cases $J=1,2,3$}:
These estimates follow from induction in $J$
via the evolution equation
\eqref{E:OUTGOINGRIEMANNINVARIANTEVOLUTION},
the bound $\Speed \lesssim 1$ (which follows easily from the bootstrap assumptions),
the already proven estimate
\eqref{E:EXTERIORREGIONPOINTWISEESTIMATEFORRPLUSANDLDERIVATIVES} in the case $J=0$,
and the already proven estimates in \eqref{E:EXTERIORREGIONPOINTWISEESTIMATEFORRMINUSANDLDERIVATIVES}.

\medskip 

\noindent \textbf{Proof of the remaining estimates}:
The remaining estimates stated in the proposition can be inductively derived by commuting
the evolution equations
\eqref{E:OUTGOINGRIEMANNINVARIANTEVOLUTION}--\eqref{E:INGOINGRIEMANNINVARIANTEVOLUTION}
with factors of $\Lunit$ and $\muX$ and using
Lemma~\ref{L:VECTORFIELDSINTERMSOFGEOMETRICCOORDINATES},
the bootstrap assumptions,
the already proven estimates,
\eqref{E:SHARPPOINWISECOMPARISONBETWEENUANDTMINUSR},
and the identity $\muX = \frac{1}{2}(\upmu \Lunit - \muuLunit)$
(which follows from Lemma~\ref{L:VECTORFIELDSINTERMSOFGEOMETRICCOORDINATES}).
All estimates for $\RRiemann$ can be derived 
via equations \eqref{E:EXTERIORREGIONPOINTWISEESTIMATEFORCOMMUTEDVERSIONLDERIVATIMEOFMODIFIEDLBARRPLUSDERIVATIVE}
and \eqref{E:POINTWISEBOUNDALLMUXDERIVATIVESTRANSPORTEDMODIFIEDULUNITRPLUSTSTARDATA}
without the need to derive any further estimates for
$\LRiemann$. In contrast, some of the estimates for 
$\LRiemann$ depend on having already derived estimates for $\RRiemann$.
By ``inductively,'' we mean that it is convenient to bound quantities in the following order:
$\muX \RRiemann$,
$\Lunit \muX \RRiemann$,
$\Lunit \Lunit \muX \RRiemann$,
$\muX \muX \RRiemann$,
$\Lunit \muX \muX \RRiemann$,
$\muX \muX \muX \RRiemann$,
$\muX \LRiemann$,
$\Lunit \muX \LRiemann$,
$\Lunit \Lunit \muX \LRiemann$,
$\muX \muX \LRiemann$,
$\Lunit \muX \muX \LRiemann$,
$\muX \muX \muX \LRiemann$.

For example, 
the estimate 
$|\muX \RRiemann| \leq 
C \frac{\datasize}{(1 + t + |u|)(1 + |u|)^2}
+
C \frac{\datasize}{(1 + t + |u|)^2} 
$
featured in \eqref{E:EXTERIORREGIONPOINTWISEESTIMATEFORRPLUSATLEASTONEXDERIVATIVE}
follows from using
\eqref{E:POINTWISEBOUNDALLMUXDERIVATIVESTRANSPORTEDMODIFIEDULUNITRPLUSTSTARDATA}
in the case $K=0$,
the identity $\muX = \frac{1}{2}(\upmu \Lunit - \muuLunit)$,
the bootstrap assumptions,
and \eqref{E:RISAPPROXIMATELYTMINUSU}--\eqref{E:COMPARISONBETWEENTOVERUTMINUSUANDTPLUSMODU}
to algebraically express $\muX \RRiemann$
as a sum of terms that have already been shown to be $\lesssim \frac{\datasize}{(1 + t + |u|)^2} 
+ 
\frac{\datasize}{(1 + t + |u|)(1 + |u|)^2}$.
As a second example, we note that the estimate $|\muX \LRiemann| \leq C \frac{\datasize}{(1 + t + |u|)^2}$
follows from using equation \eqref{E:INGOINGRIEMANNINVARIANTEVOLUTION},
the identity $\muX = \frac{1}{2}(\upmu \Lunit - \muuLunit)$,
the bootstrap assumptions,
and  \eqref{E:RISAPPROXIMATELYTMINUSU}--\eqref{E:COMPARISONBETWEENTOVERUTMINUSUANDTPLUSMODU}
to algebraically express $\muX \LRiemann$
as a sum of terms that have already been shown to be $\lesssim \frac{\datasize}{(1 + t + |u|)^2}$.

\end{proof}

\section{Diffeomorphism and homeomorphism properties of the change of variables map}
\label{S:CHOVDIFFEOMORPHISM}
In this section, we prove a proposition that exhibits various diffeomorphism and homeomorphism properties of the change of variables map
$\Upsilon$ from Sect.\,\ref{SS:CHOVMAP}. When proving our main MGHD theorem, we will use the proposition 
to transition back and forth between the two coordinate systems.

\begin{proposition}[Diffeomorphism and homeomorphism properties of the change of variables map]
	\label{P:CHOVDIFFEOMORPHISM}
	Let $\mathbf{GH}_{\textnormal{Boot}}^{{\textnormal{Ext}}}$
	be a globally hyperbolic exterior bootstrap region in the sense of Def.\,\ref{D:GLOBALLHYPERBOLICEXTERIORBOOTSTRAPREGION} that,
	in addition, is foliated by intervals of constant $u$ 
	and has a piecewise $C^1$ future-boundary 
	i.e.,
	either:
	\begin{align} \label{E:GHFOLIATEDBYULINESOFCONSTANTU}
		\mathbf{GH}_{\textnormal{Boot}}^{{\textnormal{Ext}}}
		& = \cup_{t \in [\tstar,\Tboot)} \lbrace (t,u) \ | \ A(t) < u < B(t) \rbrace,
	\end{align}
	or:
	\begin{align} \label{E:ALTGHFOLIATEDBYULINESOFCONSTANTU}
		\mathbf{GH}_{\textnormal{Boot}}^{{\textnormal{Ext}}}
		& = \cup_{t \in [\tstar,\Tboot)} \lbrace (t,u) \ | \ - \infty < u < B(t) \rbrace,
	\end{align}
	where $A$ and $B$ are continuous, piecewise $C^1$ functions on $[\tstar,\Tboot)$ that continuously extend to the closed interval $[\tstar,\Tboot]$.
	Under the data assumptions of Section~\ref{S:DATA} and the bootstrap assumptions of
Section~\ref{S:EXTERIORGLOBALLYHYPERBOLICBOOTSTRAPASSUMPTIONS}, 
if $\datasize$ is sufficiently small, then
the change of variables map $\Upsilon$ defined by
\eqref{E:CHOVFROMTUTOTRCOORDINATES} satisfies \eqref{E:MAINSTATEMENTCHOVJACOBIANDTERMINANT}
and enjoys the following properties:

\begin{enumerate}
	\item (Diffeomorphism on the bootstrap region). $\Upsilon$ is a global $C^3$ diffeomorphism on
	$\mathbf{GH}_{\textnormal{Boot}}^{{\textnormal{Ext}}}$.
	\item (Homeomorphism on the closure).
	$\Upsilon(t,u)$ extends to a $C^3$ homeomorphism on the closure of 
	$\mathbf{GH}_{\textnormal{Boot}}^{{\textnormal{Ext}}}$,
	which in the case of \eqref{E:GHFOLIATEDBYULINESOFCONSTANTU} is: 
\begin{align} \label{E:CLOSUREGHFOLIATEDBYULINESOFCONSTANTU}
	\overline{\mathbf{GH}_{\textnormal{Boot}}^{{\textnormal{Ext}}}} = \cup_{t \in [\tstar,\Tboot]} \lbrace (t,u) \ | \ A(t) \leq u \leq B(t) \rbrace,
\end{align}
and in the case of \eqref{E:ALTGHFOLIATEDBYULINESOFCONSTANTU} is:
\begin{align} \label{E:ALTCLOSUREGHFOLIATEDBYULINESOFCONSTANTU}
	\overline{\mathbf{GH}_{\textnormal{Boot}}^{{\textnormal{Ext}}}} = \cup_{t \in [\tstar,\Tboot]} \lbrace (t,u) \ | \ - \infty < u \leq B(t) \rbrace.
\end{align}
\end{enumerate}

\end{proposition}

\begin{proof}
	We assume that $\overline{\mathbf{GH}_{\textnormal{Boot}}^{{\textnormal{Ext}}}}$ is of the form \eqref{E:CLOSUREGHFOLIATEDBYULINESOFCONSTANTU}; the case 
	in which $\overline{\mathbf{GH}_{\textnormal{Boot}}^{{\textnormal{Ext}}}}$ is of the form \eqref{E:ALTGHFOLIATEDBYULINESOFCONSTANTU}
	can be handled using similar arguments.
	We recall that the Jacobian matrix of the change of variables map $\Upsilon(t,u)$ is given by the formula \eqref{E:CHOVJACOBIAN}.
	Since RHS~\eqref{E:CHOVJACOBIAN} is $C^2$, it follows that $\Upsilon$ is $C^3$ on $\mathbf{GH}_{\textnormal{Boot}}^{{\textnormal{Ext}}}$.
	Moreover,  the bootstrap assumptions and the transport equation \eqref{E:MUEVOLUTION}
	imply that $\Lunit \mathrm{d} \Upsilon = \frac{\partial}{\partial t} \mathrm{d} \Upsilon$ is uniformly bounded 
	in $C_{t,u}^2(\mathbf{GH}_{\textnormal{Boot}}^{{\textnormal{Ext}}})$,
	and thus $\Upsilon$ extends to a $C_{t,u}^3$ function on
	$\overline{\mathbf{GH}_{\textnormal{Boot}}^{{\textnormal{Ext}}}}$.
	
	Next, using \eqref{E:CHOVJACOBIAN} and the bootstrap assumptions, we find that:
	\begin{align} \label{E:CHOVJACOBIANDTERMINANT}
		\mbox{det} (\mathrm{d} \Upsilon(t,u))
		& = - \Speed \upmu < 0, && \mbox{on } \mathbf{GH}_{\textnormal{Boot}}^{{\textnormal{Ext}}}.
	\end{align}
	Hence, by the inverse function theorem, $\Upsilon$ is a local diffeomorphism on $\mathbf{GH}_{\textnormal{Boot}}^{{\textnormal{Ext}}}$.
	Moreover, since \eqref{E:CHOVJACOBIAN} implies:
	\begin{align} \label{E:RADIALCOORDINATEWEAKLYDECREASINGINU}
		\frac{\partial}{\partial u} r = - \Speed \upmu,
	\end{align}
	it follows that at each fixed $t \in [\tstar,\Tboot)$,
	the map $u \rightarrow r(t,u)$ is one-to-one for $u \in (A(t), B(t))$ and thus
	$u \rightarrow r(t,u)$ extends to a $C^3$, one-to-one function of $u$ on the closed interval $[(A(t), B(t)]$.
	It follows that the extended map $\Upsilon$ is injective on
	$\cup_{t \in [\tstar,\Tboot)} \lbrace (t,u) \ | \ A(t) \leq u  \leq B(t) \rbrace$.
	The main remaining task is to show that the extended function $r(t,u)$ is such that
	the function $u \rightarrow r(\Tboot,u)$ is one-to-one for $u \in [A(\Tboot), B(\Tboot)]$.
	This implies that $\Upsilon$ is injective on $\overline{\mathbf{GH}_{\textnormal{Boot}}^{{\textnormal{Ext}}}}$,
	which, by compactness, implies that $\Upsilon$ is a homeomorphism on 
	$\overline{\mathbf{GH}_{\textnormal{Boot}}^{{\textnormal{Ext}}}}$.
	
	To show that $u \rightarrow r(\Tboot,u)$ is injective for $u \in [A(\Tboot), B(\Tboot)]$,
	it suffices to prove the following claim:
	that $\upmu(\Tboot,u) > 0$ must hold for all $u \in [A(\Tboot), B(\Tboot)]$
	except at possibly one point $u_*$. The injectiveness of the map $u \rightarrow r(\Tboot,u)$ 
	then follows from the identity \eqref{E:RADIALCOORDINATEWEAKLYDECREASINGINU} 
	(which by continuity continues to hold at time $\Tboot$)
	and the mean value theorem.
	
	To prove the claim from the previous paragraph, we first restrict our attention to the open interval $(A(\Tboot), B(\Tboot))$.
	Since $\upmu$ is positive on $\mathbf{GH}_{\textnormal{Boot}}^{{\textnormal{Ext}}}$,
	at any point $u_* \in (A(\Tboot), B(\Tboot))$ 
	at which the function $u \rightarrow \upmu(\Tboot,u)$ vanishes, $\upmu$ achieves its minimum value of $0$, 
	i.e.,  $\upmu(\Tboot,u_*) = \muX \upmu(\Tboot,u_*) = 0$, $\muX \muX \upmu(\Tboot,u_*) \geq 0$.
	The $\upmu$-estimates of Prop.\,\ref{P:SHARPESTIMATESFORMU} imply that these three conditions are simultaneously possible only for at most a single point $u_*$
	satisfying $u_* = \mylittleo(\datasize)$, at which we must have $\muX \muX \upmu(\Tboot,u_*) > \frac{3}{2}$. That is, 
	the map $u \rightarrow \upmu(\Tboot,u)$ can have at most a single $0$ at an interior point $u_* \in (A(\Tboot), B(\Tboot))$.
	Moreover, if $u \rightarrow \upmu(\Tboot,u)$ does have such an interior zero at $u_*$, then the $\upmu$-estimates of 
	Prop.\,\ref{P:SHARPESTIMATESFORMU} imply that $\muX \upmu(\Tboot,u) > 0$ for all points $u > u_*$ satisfying $\upmu(\Tboot,u) < 1/200$ and
	$\muX \upmu(\Tboot,u) < 0$ for all points $u < u_*$ satisfying $\upmu(\Tboot,u) < 1/200$.
	This implies, in particular, that either $\upmu(\Tboot,A(\Tboot)) \geq 1/200$ or $\muX \upmu(\Tboot,A(\Tboot)) > 0$.
	That is, either $\upmu(\Tboot,A(\Tboot)) \geq 1/200$, or $\upmu(\Tboot,A(\Tboot))$ is more positive than nearby points
	$\upmu(\Tboot,u)$ with $u \in (A(\Tboot),B(\Tboot))$. Either way, we must have $\upmu(\Tboot,A(\Tboot)) > 0$.
	Similar logic implies that $\upmu(\Tboot,B(\Tboot)) > 0$. We have therefore shown that if $u \rightarrow \upmu(\Tboot,u)$
	has a zero on $(A(\Tboot), B(\Tboot))$, then it must be the only such zero, and $\upmu(\Tboot,A(\Tboot)) > 0$, $\upmu(\Tboot,B(\Tboot)) > 0$.
	
	To finish the proof of the claim, we will show via a contradiction argument that on the closed interval $[A(\Tboot), B(\Tboot)]$,
	it is impossible for $u \rightarrow \upmu(\Tboot,u)$ to vanish at both endpoints $A(\Tboot)$ and $B(\Tboot)$ while being
	positive in the interior. If this were to occur, then $\upmu(\Tboot,A(\Tboot)) = 0$,
	$\muX \upmu(\Tboot,A(\Tboot)) \geq 0$,
	$\upmu(\Tboot,B(\Tboot)) = 0$,
	and $\muX \upmu(\Tboot,B(\Tboot)) \leq 0$,
	and the estimates of
	Prop.\,\ref{P:SHARPESTIMATESFORMU} imply that $A(\Tboot) \geq - |\mylittleo(\datasize)|$
	and 
	$B(\Tboot) \leq |\mylittleo(\datasize)|$. Since we also have $A(\Tboot) < B(\Tboot)$,
	it follows that $[A(\Tboot),B(\Tboot)] \subset [-|\mylittleo(\datasize)|, |\mylittleo(\datasize)|]$,
	and the estimates of Prop.\,\ref{P:SHARPESTIMATESFORMU} then imply that
	$\frac{3}{2} < \muX \muX \upmu(\Tboot,u) < 3$ for all $u \in [A(\Tboot),B(\Tboot)]$.
	That is, $u \rightarrow \muX \upmu(\Tboot,u)$ is strictly increasing on 
	$[A(\Tboot),B(\Tboot)]$, making it impossible to have both 
	$\muX \upmu(\Tboot,A(\Tboot)) \geq 0$ and $\muX \upmu(\Tboot,B(\Tboot)) \leq 0$.
	We have therefore demonstrated a contradiction.
	In total, we have proved the claim from the previous paragraph, which completes the proof of the proposition.
	
\end{proof}

\section{The behavior of the solution up to the crease, the singular boundary, and a portion of the Cauchy horizon}
\label{S:EXISTENCEUPTOCREASEANDSINGULARBOUNDARYANDPORTIONOFCAUCHYHORIZON}
In Prop.\,\ref{P:EXISTENCEUPTOCREASEANDSINGULARBOUNDARYANDPORTIONOFCAUCHYHORIZON}, 
we derive the structure of the part of the MGHD that lies in the exterior region.
This marks the culmination of the most difficult and interesting analysis needed for our first main theorem,
though we delay the full statement of our MGHD results until Theorem~\ref{T:MAINMGHDEXISTENCETHEOREM}.

\begin{proposition}[The Exterior Region: Behavior of the solution up to the crease, the singular boundary, and a portion of the Cauchy horizon]
\label{P:EXISTENCEUPTOCREASEANDSINGULARBOUNDARYANDPORTIONOFCAUCHYHORIZON}
Under the data-assumptions 
\eqref{E:RPLUSDATAISCLOSETOBACKGROUNDATTIME0}--\eqref{E:RMINUSDATAISCLOSETOBACKGROUNDATTIME0},
if $\datasize$ is sufficiently small, then 
the solution $(\RRiemann,\LRiemann)$ to 
\eqref{E:OUTGOINGRIEMANNINVARIANTEVOLUTION}--\eqref{E:INGOINGRIEMANNINVARIANTEVOLUTION}
exhibits the following behavior in the ``Exterior Region'' portion of Fig.\,\ref{F:MGHD}
(see Remark~\ref{R:NOTATIONALCONVENTIONCOORDIANTEDEPENDENCE} regarding our notation).

\begin{enumerate}
\item There exist numbers $\ucrease$ and $\Tcrease$ and a corresponding point $(\Tcrease,\ucrease)$, called the crease, 
such that the following estimates hold, 
where $\lifespanconstant > 0$ is the constant from \eqref{E:NULLCONDITIONFAILURECONSTANT} and $\tstar = \frac{1}{\datasize}$:
\begin{align}  \label{E:UCREASEISSMALL}
\ucrease & = \mylittleo(\datasize),
	\\
\label{E:BLOWUPTIMEUPPERANDLOWERBOUNDS}
\exp
	\left(
	\frac{\frac{99}{100} (1 + \ucrease^2)}{\lifespanconstant \datasize}
	\right)
	\leq
	\frac{\Tcrease-\ucrease}{\tstar - \ucrease} 
	&
	\leq
	\exp
	\left(
	\frac{\frac{101}{100}(1 + \ucrease^2)}{\lifespanconstant \datasize}
	\right).
\end{align}
Moreover, there exist functions $\tboundary(u)$, $\Tblowup(u)$, $\TCH(u)$
such that:
\begin{itemize}
	\item $\Tblowup(u)$ is decreasing and $C^2$ on $(-\infty,\ucrease]$
		and satisfies:
	\begin{align} \label{E:EXISTENCETIMEDEPENDINGONUSINGULARBOUNDARYMAINSTATEMENT}
	\exp
	\left(
	\frac{\frac{99}{100} (1 + u^2)}{\lifespanconstant \datasize}
	\right)
	\leq
	\frac{\Tblowup(u)-u}{\tstar - u} 
	&
	\leq
	\exp
	\left(
	\frac{\frac{101}{100}(1 + u^2)}{\lifespanconstant \datasize}
	\right).
\end{align}		
	\item $\TCH(u)$ is $C^3$ and increasing on $[\ucrease,\tstar/2]$.
	\item
\begin{align} \label{E:TASAFUNCTIONOFUBOUNDARYPIECEWISEDEFINITION}
	\tboundary(u)
	& = 
		\begin{cases}
		\Tblowup(u), & u \in  (-\infty,\ucrease],
		\\
			\TCH(u), & u \in (\ucrease,\tstar/2],
	\end{cases}
\end{align}
\end{itemize}
and $\tboundary$ is $C^1$ across $\ucrease$, with:
\begin{subequations}
\begin{align} \label{E:BLOWUPTIMERELATIONS}
	\Tblowup(\ucrease) 
	& = \TCH(\ucrease) = \Tcrease,
		\\
	\frac{d}{du} 
	\Tblowup(\ucrease) 
	& = \frac{d}{du} \TCH(\ucrease) = 0.
\end{align}
\end{subequations}
Furthermore, 
\begin{align} \label{E:SHARPUPPERANDLOWERBOUNDSONVALUEOFTATTOPPOINTOFUEQUALSTSTAROVER2}
	\Tcrease & < \TCH\left(\frac{\tstar}{2} \right) 
	< \Tcrease + \frac{\tstar}{4} + C \datasize.
\end{align}
\item
The solution exists classically with respect to the geometric coordinates $(t,u)$ on the following subset of
geometric coordinate space:
\begin{align} \label{E:GLOBALLYHYPERBOLICEXTERIOREGIONPLUSABITOFCAUCHYHORIZON}
	\mathbf{MGHD}^{\textsf{Ext}}
	& = \lbrace (t,u) \ | \ u \in (-\infty,\tstar/2], \, t \in [\tstar, \tboundary(u)) \rbrace,
\end{align}
which is a globally hyperbolic bootstrap region in the sense of Definition~\ref{D:GLOBALLHYPERBOLICEXTERIORBOOTSTRAPREGION}.
Moreover, with respect to the $(t,u)$ coordinates, the solution uniquely extends to a classical solution on the closure of $\mathbf{MGHD}^{\textsf{Ext}}$,
which is:
\begin{align}
	\overline{\mathbf{MGHD}^{\textsf{Ext}}}
	& = \lbrace (t,u) \ | \ u \in (-\infty,\tstar/2], \, t \in [\tstar, \tboundary(u)] \rbrace.
\end{align}
\item On $\overline{\mathbf{MGHD}^{\textsf{Ext}}}$, $\upmu$ is everywhere positive, except on the singular boundary
	$\futuresinghyp$ defined by:
	\begin{align} \label{E:SINGULARBOUNDARY}
		\futuresinghyp := \lbrace (\Tblowup(u),u) \ | \ u \in (-\infty,\ucrease] \rbrace,
	\end{align}	
	along which $\upmu$ vanishes. Furthermore, with $\futurecrease = (\Tcrease,\ucrease)$ denoting the past boundary of $\futuresinghyp$
	(which is equal to the crease), 
	we have that $\muX \upmu \restriction_{\futuresinghyp \backslash \futurecrease} < 0$ and  $\muX \upmu \restriction_{\futurecrease} = 0$.
	In particular, $\upmu$ is positive on $\underline{\mathcal{C}}^{\tstar/2}$, the Cauchy Horizon portion of the boundary 
	defined by:
	\begin{align} \label{E:CAUCHYHORIZONUPTOTSTAROVER2}
		\underline{\mathcal{C}}^{\tstar/2}
		& := \lbrace (\TCH(u),u) \ | \ u \in (\ucrease,\tstar/2] \rbrace.
	\end{align}
Moreover, $\underline{\mathcal{C}}^{\tstar/2}$ is a future-portion of the integral curve of $\muuLunit$
emanating from the crease.
\item The change of variables map $\Upsilon(t,u) := (t,r)$
is a global $C^3$ diffeomorphism on the set
$\mathbf{MGHD}^{\textsf{Ext}} \cup \underline{\mathcal{C}}^{\tstar/2}$,
on which $\upmu > 0$.
Moreover, $\Upsilon(t,u)$ is a $C^3$ homeomorphism on $\overline{\mathbf{MGHD}^{\textsf{Ext}}}$.
\item The Riemann invariants $\RRiemann \circ \Upsilon^{-1},\LRiemann \circ \Upsilon^{-1}$ exist as classical solutions with respect to the standard $(t,r)$-coordinates on $\Upsilon(\mathbf{MGHD}^{\textsf{Ext}} \cup \underline{\mathcal{C}}^{\tstar/2})$. In addition, they also extend as continuous functions with respect to the standard 	
	coordinates $(t,r)$ on $\Upsilon(\overline{\mathbf{MGHD}^{\textsf{Ext}}})$. 
	In particular, the extended functions are continuous on $\Upsilon(\futuresinghyp)$.
\item $\Upsilon(\futuresinghyp)$ is a $C^1$ curve-portion in $(t,r)$-space such that at each point of $\Upsilon(\futuresinghyp)$, 
$\Lunit = \partial_t + ((v^r + \Speed)\circ \Upsilon^{-1}) \partial_r$ spans the tangent space.
\item $\Upsilon(\underline{\mathcal{C}}^{\tstar/2})$ is a $C^1$- curve-portion in $(t,r)$-space such that at each point of 
$\Upsilon(\underline{\mathcal{C}}^{\tstar/2})$, 
	$\uLunit = \partial_t + ((v^r - \Speed)\circ \Upsilon^{-1})\partial_r$ is a tangent vector. 
\item  $\Upsilon(\futuresinghyp)$ and $\Upsilon(\underline{\mathcal{C}}^{\tstar/2})$ have
		$\Upsilon(\Tcrease,\ucrease)$ (the image of the crease under $\Upsilon$)
		as a common boundary point, and at this point, 
		the curve $\Upsilon(\futuresinghyp) \cup \Upsilon(\underline{\mathcal{C}}^{\tstar/2})$
		has a corner, i.e., the curve is continuous at $\Upsilon(\Tcrease,\ucrease)$, but not differentiable in the $(t,r)$-coordinate differential structure
		because at $\Upsilon(\Tcrease,\ucrease)$, its tangent from the left is $\uLunit(\Tcrease,\ucrease)$,
		while its tangent from the right is $\Lunit(\Tcrease,\ucrease)$.
\item  If $(t_{\star},r_{\star}) \in \Upsilon(\futuresinghyp)$ and $\lbrace (t_n,r_n) \rbrace_{n \in \mathbb{N}}$ is a sequence of points in
	$\Upsilon(\mathbf{MGHD}^{\textsf{Ext}} \cup \underline{\mathcal{C}}^{\tstar/2})$ such that $\lim_{n \to \infty} (t_n,r_n) = (t_{\star},r_{\star})$,
	then $\limsup_{n \to \infty} \left|\partial_r [\RRiemann \circ \Upsilon^{-1}](t_n,r_n) \right| = \infty$, i.e., the $\partial_r$-derivative of 
	$\RRiemann$ blows up on $\Upsilon(\futuresinghyp)$.
\item
For $t \leq 0$, the solution obeys analogous results (see also Remark~\ref{R:DISCRETESYMMETRY}). That is, 
there exists a unique classical solution on a past globally hyperbolic bootstrap region
$\mathbf{MGHD}^{\textsf{Ext}}_{\textsf{Past}}$ contained in $\lbrace t \leq 0 \rbrace$, the past-boundary of 
$\mathbf{MGHD}^{\textsf{Ext}}_{\textsf{Past}}$ is the union of a closed, unbounded past-singular boundary $\pastsinghyp$ that emanates from the past-crease $\pastcrease$
and a pre-compact portion of the past-Cauchy horizon $\underline{\mathcal{C}}_{\textsf{Past}}^{\tstar/2}$ that emanates from $\pastcrease$ towards the past and terminates on $u = \tstar/2$. The change of variables $\Upsilon$ is a global $C^3$ diffeomorphism on $\mathbf{MGHD}^{\textsf{Ext}}_{\textsf{Past}} \cup \underline{\mathcal{C}}_{\textsf{Past}}^{\tstar/2}$ and extends to a global $C^3$ homeomorphism on 
$\overline{\mathbf{MGHD}^{\textsf{Ext}}_{\textsf{Past}}}$. The undifferentiated Riemann invariants extend as continuous functions to 
$\Upsilon(\pastsinghyp)$, and therefore $\uLunit$ also extends as a continuous vectorfield to $\pastsinghyp$, along which it is tangent. Finally, $\partial_r \LRiemann$ blows up when approaching $\pastsinghyp$ from within $\mathbf{MGHD}^{\textsf{Ext}}_{\textsf{Past}}$.

\end{enumerate}
	
\end{proposition}

\begin{proof}
		Throughout the proof, we will silently use the observations made Remark~\ref{R:STRICTIMPROVEMENTOFBOOTSTRAP}.
		
		\medskip
		
		\noindent \underline{\textbf{Controlling the solution up to the crease}}.
		We consider the bootstrap region $\mathbf{GH}_{\Tboot}^{{\textnormal{Ext}}}$, defined as follows in $(t,u)$-coordinates:
		\begin{align} \label{E:BOOTSTRPAREGIONUPTOCREASE}
			\mathbf{GH}_{\Tboot}^{{\textnormal{Ext}}} := [\tstar,\Tboot) \times (-\infty, \tstar/2].
		\end{align}
		If $\datasize$ is sufficiently small and
		$\Tboot$ is near $\tstar$, then by standard local well-posedness,
		on the globally hyperbolic region $\mathbf{GH}_{\Tboot}^{{\textnormal{Ext}}}$,
		$\RRiemann$,
		$\LRiemann$,
		and 
		$\upmu$ are classical solutions to the quasilinear transport equation system
		\eqref{E:OUTGOINGRIEMANNINVARIANTEVOLUTION}--\eqref{E:INGOINGRIEMANNINVARIANTEVOLUTION} 
		+
		\eqref{E:MUEVOLUTION}
		such that $\RRiemann$ and $\LRiemann$ are $C^3$ with respect to the geometric coordinates $(t,u)$,
		such that $\upmu$ is $C^2$ with respect to $(t,u)$,
		and such that the bootstrap assumptions of
		Section~\ref{S:EXTERIORGLOBALLYHYPERBOLICBOOTSTRAPASSUMPTIONS}
		hold with $\eps := \datasize^{3/4}$.
		By Props.\,\ref{P:SHARPESTIMATESFORMU} and \ref{P:RIEMANNINVARIANTSAPRIORIEXTERIORREGIONESTIMATES},
		no bootstrap assumptions are saturated on $\mathbf{GH}_{\Tboot}^{{\textnormal{Ext}}}$,
		except possibly $\inf_{\mathbf{GH}_{\Tboot}^{{\textnormal{Ext}}}} \upmu = 0$; see Remark~\ref{R:STRICTIMPROVEMENTOFBOOTSTRAP}.
		Regardless of whether the latter scenario occurs,
		the up-to-order $3$ derivatives of $\RRiemann$ and $\LRiemann$ with respect to $(t,u)$
		are uniformly bounded on $\mathbf{GH}_{\Tboot}^{{\textnormal{Ext}}}$,
		and similarly for the up-to-order $2$ derivatives of $\upmu$.
		Moreover, if $\inf_{\mathbf{GH}_{\Tboot}^{{\textnormal{Ext}}}} \upmu > 0$, then no bootstrap assumptions are saturated, and
		standard continuation criteria 
		imply that for some $\Delta > 0$, the solution extends classically in $(t,u)$ coordinates to the
		larger globally hyperbolic region
		$\mathbf{GH}_{\Tboot + \Delta}^{{\textnormal{Ext}}} := [\tstar,\Tboot + \Delta) \times (-\infty, \tstar/2]$,
		on which the bootstrap assumptions hold with $\eps$ replaced by $C \datasize$; 
		see Remark~\ref{R:STRICTIMPROVEMENTOFBOOTSTRAP}.
		It therefore follows that $\RRiemann$, $\LRiemann$, and $\upmu$
		exist classically with respect to the $(t,u)$ coordinate system
		on the closure of $\mathbf{GH}_{\Tcrease}^{{\textnormal{Ext}}}$,
		where $\Tcrease$ is the time at which $\upmu(t,u)$ first vanishes at some $u$-value belonging to the interval $(-\infty,\tstar/2]$,
		and the closure is:
		\begin{align} \label{E:CLOSUREBOOTSTRPAREGIONUPTOCREASE}
			\overline{\mathbf{GH}_{\Tcrease}^{{\textnormal{Ext}}}} & = [\tstar,\Tcrease] \times (-\infty, \tstar/2].
		\end{align}
		Moreover, since the estimates of Props.\,\ref{P:SHARPESTIMATESFORMU} and \ref{P:RIEMANNINVARIANTSAPRIORIEXTERIORREGIONESTIMATES} hold on
		$\mathbf{GH}_{\Tcrease}^{{\textnormal{Ext}}}$,
		the evolution equations \eqref{E:OUTGOINGRIEMANNINVARIANTEVOLUTION},
		\eqref{E:INGOINGRIEMANNINVARIANTEVOLUTION},
		and
		\eqref{E:MUEVOLUTION}
		imply that the estimates of Props.\,\ref{P:SHARPESTIMATESFORMU} and \ref{P:RIEMANNINVARIANTSAPRIORIEXTERIORREGIONESTIMATES}
		also hold on the closure
		$\overline{\mathbf{GH}_{\Tcrease}^{{\textnormal{Ext}}}}$.
		In particular, $\RRiemann$ and $\LRiemann$
		extend to $\overline{\mathbf{GH}_{\Tcrease}^{{\textnormal{Ext}}}}$ 
		as solutions having $C^3$ regularity with respect to $(t,u)$, 
		and $\upmu$ extends to $\overline{\mathbf{GH}_{\Tcrease}^{{\textnormal{Ext}}}}$ as a solution with $C^2$ regularity.
		In addition, the $\upmu$-estimates of Prop.\,\ref{P:SHARPESTIMATESFORMU}
		imply that the map $u \rightarrow \upmu(\Tcrease,u)$ has a unique non-degenerate zero $\ucrease$
		satisfying $\ucrease = \mylittleo(\datasize) $, i.e., $\upmu(\Tcrease,\ucrease) = \muX \upmu(\Tcrease,\ucrease) = 0$,
		$1 < \muX \muX \upmu(\Tcrease,\ucrease) < 3$.
		Furthermore,
		by Prop.\,\ref{P:CHOVDIFFEOMORPHISM}, the solution exists classically with respect to the $(t,r)$-coordinate system on 
		$\Upsilon(\mathbf{GH}_{\Tcrease}^{{\textnormal{Ext}}})$.
		For future use, we note that the identities
		$\upmu(\Tcrease,\ucrease) = \muX \upmu(\Tcrease,\ucrease) = 0$ and Lemma~\ref{L:VECTORFIELDSINTERMSOFGEOMETRICCOORDINATES} 
		imply that:
		\begin{align} \label{E:MUULUNITUPMUVANISHESATCREASE}
			\muuLunit \upmu(\Tcrease,\ucrease) = 0.
		\end{align}
		
		\medskip
		
		\noindent \underline{\textbf{Controlling the solution up to the singular boundary}}.
		Similar arguments imply that for any fixed $\hat{u} \leq \ucrease$, 
		$\RRiemann$,
		$\LRiemann$,
		and 
		$\upmu$ exist classically (with the same regularity stated above) with respect to the $(t,u)$ coordinates and satisfy the bootstrap assumptions
		on the following globally hyperbolic region:
			\begin{align} \label{E:BOOTSTRAPREGIONFORCONTROLLINGSOLUTIONUPTOSINGULARBOUNDARYATFIXEDLEFTUCHARACTERISTIC}
				\mathbf{GH}_{\Tblowup(\hat{u})}^{{\textnormal{Ext}}} := [\tstar,\Tblowup(\hat{u})) \times (-\infty, \hat{u}],
			\end{align}
		where $\Tblowup(\hat{u})$ is the time at which $\upmu(t,u)$ first vanishes when $u$ is restricted to the domain $(-\infty,\hat{u}]$.
		The same logic implies that the solution exists classically 
		with respect to the $(t,u)$ coordinates and satisfies all bootstrap assumptions aside from $\upmu > 0$, 
		on the closure of $\mathbf{GH}_{\Tblowup(\hat{u})}^{{\textnormal{Ext}}}$, which is:
			\begin{align} \label{E:CLOSUREDBOOTSTRAPREGIONFORCONTROLLINGSOLUTIONUPTOSINGULARBOUNDARYATFIXEDLEFTUCHARACTERISTIC}
				\overline{\mathbf{GH}_{\Tblowup(\hat{u})}^{{\textnormal{Ext}}}} & = [\tstar,\Tblowup(\hat{u})] \times (-\infty, \hat{u}].
			\end{align}
		Moreover, the same arguments used in the proof of
		Prop.\,\ref{P:CHOVDIFFEOMORPHISM} and in the previous paragraph 
		imply that when $\hat{u} \leq \ucrease$,
		the map $u \rightarrow \upmu(\Tblowup(\hat{u}),u)$ with domain $(-\infty,\hat{u}]$
		vanishes only when $u = \hat{u}$, and that when 
		$\hat{u} < \ucrease$, we have $\muX \upmu(\Tblowup(\hat{u}),\hat{u}) < 0$.
		In particular, 
		$\Tblowup(\ucrease) = \Tcrease$,
		and in the previous paragraph, we showed that $\muX \upmu(\Tblowup(\ucrease),\ucrease) = 0$.
		Note that the map 
		$\hat{u} \rightarrow \Tblowup(\hat{u})$ 
		is a decreasing function of $\hat{u}$ (because the domain  $(-\infty, \hat{u}]$ increases as $\hat{u}$ increases).
		We will now give a more precise description of the behavior of this map.
		By \eqref{E:QUANTITATIVELUNITMUNEGATIVEWHENMUISSMALL} and the implicit function theorem, the map
		$\hat{u} \rightarrow \Tblowup(\hat{u})$ is $C^2$.
		By differentiating the equation $\upmu\left(\Tblowup(\hat{u}),\hat{u}\right) = 0$, we further deduce that:
		\begin{align} \label{E:IMPLICITFUNCTIONTHEOREMRELATIONFORUDERIVATIVEOFBLOWUPCURVE}
			\frac{d}{du} \Tblowup(\hat{u}) 
			& =  - \frac{[\muX \upmu]\left(\Tblowup(\hat{u}),\hat{u}\right)}{[\Lunit \upmu]\left(\Tblowup(\hat{u}),\hat{u}\right)}.
		\end{align}
		The estimates of Prop.\,\ref{P:SHARPESTIMATESFORMU} and the arguments given above
		imply that RHS~\eqref{E:IMPLICITFUNCTIONTHEOREMRELATIONFORUDERIVATIVEOFBLOWUPCURVE} vanishes when $\hat{u} = \ucrease$
		and that for $u < \hat{u}$, RHS~\eqref{E:IMPLICITFUNCTIONTHEOREMRELATIONFORUDERIVATIVEOFBLOWUPCURVE}
		is strictly negative. In turn, this implies that on $(-\infty,\ucrease]$, the map $\hat{u} \rightarrow \Tblowup(\hat{u})$ is strictly decreasing.
		
		Moreover, 
		\eqref{E:QUANTITATIVEESTIMATEFORTINREGIONWHEREMUISSMALL} implies that the function $\hat{u} \rightarrow \Tblowup(\hat{u})$ satisfies:
	\begin{align} \label{E:EXISTENCETIMEDEPENDINGONUSINGULARBOUNDARYPROOF}
	\exp
	\left(
	\frac{\frac{99}{100} (1 + \hat{u}^2)}{\lifespanconstant \datasize}
	\right)
	\leq
	\frac{\Tblowup(\hat{u})-\hat{u}}{\tstar - \hat{u}} 
	&
	\leq
	\exp
	\left(
	\frac{\frac{101}{100}(1 + \hat{u}^2)}{\lifespanconstant \datasize}
	\right).
\end{align}		
\eqref{E:EXISTENCETIMEDEPENDINGONUSINGULARBOUNDARYPROOF} shows in particular that
as $\hat{u} \rightarrow - \infty$, $\Tblowup(\hat{u}) \rightarrow \infty$.
The desired estimate \eqref{E:BLOWUPTIMEUPPERANDLOWERBOUNDS} follows from
\eqref{E:EXISTENCETIMEDEPENDINGONUSINGULARBOUNDARYPROOF} and the bound $\ucrease = \mylittleo(\datasize)$ proved above.
We have therefore shown that the solution exists classically with respect to the $(t,u)$ coordinates on the following globally hyperbolic region:
\begin{align} \label{E:UPARAMETERIZEDEXTERIORREGION}
	\widetilde{\mathbf{GH}}^{{\textnormal{Ext}}}
	& := \lbrace (t,\hat{u}) \ | \ \hat{u} \in (-\infty,\ucrease], \, t \in [\tstar, \Tblowup(\hat{u})) \rbrace,
\end{align}		
and that in $(t,u)$ coordinates, the solution extends as a classical solution to its closure:
\begin{align} \label{E:CLOSUREDUPARAMETERIZEDEXTERIORREGION}
	\overline{\widetilde{\mathbf{GH}}^{{\textnormal{Ext}}}}
	& = \lbrace (t,\hat{u}) \ | \ \hat{u} \in (-\infty,\ucrease], \, t \in [\tstar, \Tblowup(\hat{u})] \rbrace.
\end{align}	
Since $\Tblowup$ is a strictly decreasing function on $(-\infty,\ucrease]$, we can alternatively express:
\begin{align} \label{E:TPARAMETERIZEDEXTERIORREGION}
	\widetilde{\mathbf{GH}}^{{\textnormal{Ext}}}
	& = \lbrace (t,\hat{u}) \ | \ t \in [\tstar, \infty), \, \hat{u} \in (-\infty,\ublowup(t)),
\end{align}
where $\ublowup$ is the inverse function of $\hat{u}\rightarrow \Tblowup(\hat{u})$.
Moreover, we define the smooth function $\texistenceintermediateboundary(u)$ by:
\begin{align}
	\texistenceintermediateboundary(u)
	& := \begin{cases}
				\Tcrease, & \mbox{for } u \in (\ucrease,\tstar/2),
					\\
				\Tblowup(u), & \mbox{for } u \in (-\infty,\ucrease].
			 \end{cases}
\end{align}
By construction, $u \rightarrow \texistenceintermediateboundary(u)$ is continuous and piecewise-$C^2$ on the domain 
$(-\infty,\ucrease]$.

Thus far, we have shown that the solution exists classically on the following globally hyperbolic region, 
depicted with respect to the $(t,u)$ coordinates in Fig.\,\ref{F:EXTERIORSINGULARGHD}:
\begin{align} \label{E:UPARAMETERIZEDEXTERIORREGIONEXTENDEDPASTCREASE}
	\mathbf{GH}^{{\textnormal{Singular;Ext}}}
	& := \mathbf{GH}_{\Tcrease}^{{\textnormal{Ext}}} \cup \widetilde{\mathbf{GH}}^{{\textnormal{Ext}}}
	= \lbrace (t,\hat{u}) \ | \ \hat{u} \in (-\infty,\ucrease], \, t \in [\tstar, \texistenceintermediateboundary(\hat{u})) \rbrace.
\end{align}
It is straightforward to see that $\mathbf{GH}^{{\textnormal{Singular;Ext}}}$ is foliated
by intervals of constant $u$.
Hence, we can apply
Prop.\,\ref{P:CHOVDIFFEOMORPHISM} to conclude that the change of variables map $\Upsilon$ is a $C^3$ diffeomorphism on 
$\mathbf{GH}^{{\textnormal{Singular;Ext}}}$	
that extends to a $C^3$ homeomorphism on the closure, which by \eqref{E:UPARAMETERIZEDEXTERIORREGIONEXTENDEDPASTCREASE} is:
\begin{align} \label{E:CLOSUREUPARAMETERIZEDEXTERIORREGION}
	\overline{\mathbf{GH}^{{\textnormal{Singular;Ext}}}}
	& = \lbrace (t,\hat{u}) \ | \ \hat{u} \in (-\infty,\tfrac{\tstar}{2}], \, t \in [\tstar, \texistenceintermediateboundary(\hat{u})] \rbrace.
\end{align}	

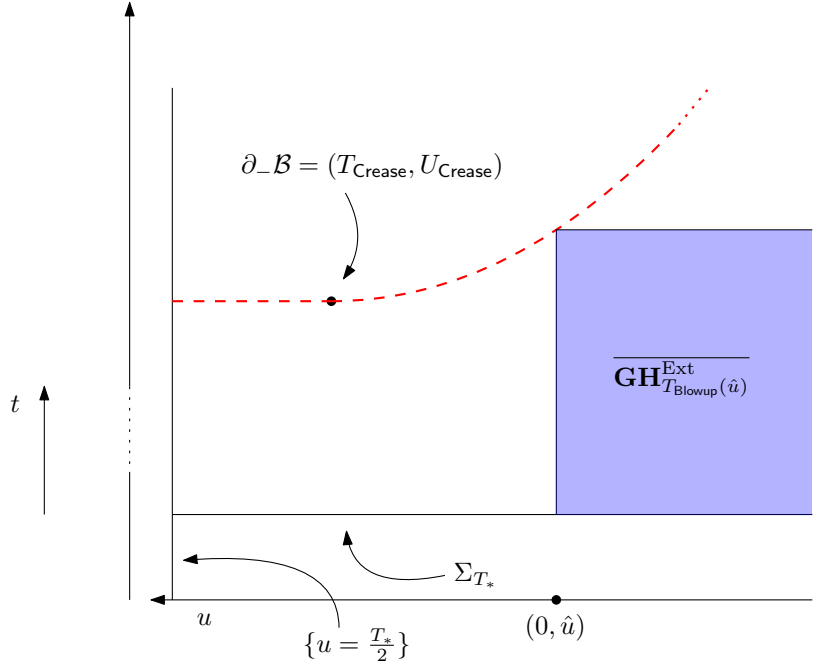
\begin{figure} 
	\centering	
		\begin{overpic}[scale=1, grid = false, tics=3, trim=-.5cm -1cm -1cm -.5cm, clip]{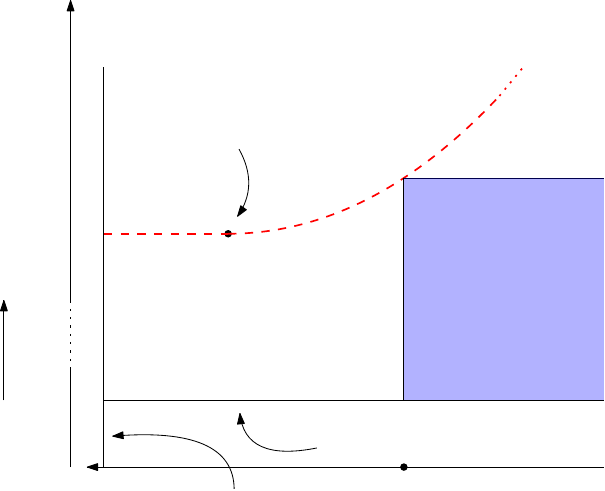} 
				\put (34,6) {\small{$\{ u = \tfrac{\tstar}{2}\}$}}
				\put (51,14) {\small{$\Sigma_{\tstar}$}}
				\put (27,60) {\small{$\crease = (\Tcrease,\ucrease)$}}
				\put (22,9) {\small{$u$}}
				\put (1,33) {\small{$t$}}
				\put (69,36) {\small{$\overline{\mathbf{GH}_{\Tblowup(\hat{u})}^{{\textnormal{Ext}}}}$}}
				\put (59,8) {\small{$(0,\hat{u})$}}
			\end{overpic}
			\caption{The dashed red curve in the figure is the parametrization $\{(\texistenceintermediateboundary(u),u) \, | \, u \in (-\infty,\tfrac{\tstar}{2}]\}$, 
			where the singular boundary is the curved portion on the right of the crease $\crease$. We displayed a dotted portion of the $t$-axis to indicate that 
			$\Tcrease  \approx\exp(C \mathring \upalpha^{-1})$ is significantly larger than $\tstar = \mathring \upalpha^{-1}$. The set $\mathbf{GH}^{{\textnormal{Singular;Ext}}}$ is the entire region below -- and not including -- the dashed red curve and above $\Sigma_{\tstar}$. For a given $\hat{u} < \ucrease$, the infinite, closed rectangular region $\overline{\mathbf{GH}_{\Tblowup(\hat{u})}^{{\textnormal{Ext}}}}$ is the portion shaded in blue, where the top left corner (which lies on the dashed red curve) and horizontal half-line emanating from it as $u \to -\infty$ are included in $\overline{\mathbf{GH}_{\Tblowup(\hat{u})}^{{\textnormal{Ext}}}}$.}
			\label{F:EXTERIORSINGULARGHD}
		\end{figure} 
\medskip
		
		\noindent \underline{\textbf{The behavior of the solution up to a portion of the Cauchy horizon}}.
	Next, given a number $\uboot \in (\ucrease,\tstar/2)$, we consider the following approximately-triangular-shaped 
	bootstrap region $\mathbf{GH}_{\uboot}$, 
	depicted in Fig.\,\ref{F:TRIANGLEBOOTSTRAPREGION}, and
	defined as follows in $(t,u)$-coordinates:
	\begin{align} \label{E:BOOTSTRAPREGIONSMALLPIECEOFCAUCHYHORIZON}
		\mathbf{GH}_{\uboot}
		& := \lbrace (t,u) \ | \ t \in [\tstar,\mathfrak{t}(\uboot)), \, u \in (\mathfrak{u}_{\uboot}(t),\tstar/2) \rbrace,
	\end{align}
		where the function $t \rightarrow \mathfrak{u}_{\uboot}(t)$ is the value of $u$ along the $t$-parameterized integral curve of $\uLunit$
		that emanates from the point with geometric coordinates $(\Tcrease,\uboot)$, 
		and $\mathfrak{t}(\uboot)$ is the value of $t$ along that integral curve when it intersects the characteristic $\lbrace u = \tstar/2 \rbrace$.
		Note that since $\uboot > \ucrease$, the vectorfield $\uLunit = \frac{\partial}{\partial t} + 2 \upmu^{-1} \frac{\partial}{\partial u}$ 
		is well defined, 
		and so are the functions $\mathfrak{t}(\uboot)$ and $\mathfrak{u}_{\uboot}(t)$.
		In particular, 
		$\mathfrak{t}(\uboot) > \Tcrease$,
		$\mathfrak{u}_{\uboot}(\Tcrease) = \uboot$,
		and the sets $\mathbf{GH}_{\uboot}$ are decreasing with respect to the parameter $\uboot$.
		The left lateral boundary of $\mathbf{GH}_{\uboot}$ is $\lbrace u = \tstar/2 \rbrace \cap \lbrace \tstar \leq t \leq \mathfrak{t}(\uboot) \rbrace$,
		the right lateral boundary of $\mathbf{GH}_{\uboot}$ is the portion of the aforementioned 
		integral curve of $\uLunit$ on which $\tstar \leq t \leq \mathfrak{t}(\uboot)$,
		and the bottom boundary of $\mathbf{GH}_{\uboot}$ is 
		$\lbrace t = \tstar \rbrace \cap \lbrace \mathfrak{u}_{\uboot}(\tstar) \leq u \leq \tstar/2 \rbrace$.
		If $\datasize$ is sufficiently small and
		$\uboot$ is near $\tstar/2$, then by standard local well-posedness
		and the fact that we have already shown that the solution exists classically on 
		$\mathbf{GH}^{{\textnormal{Singular;Ext}}}$,
		the solution also
		exists classically and satisfies the bootstrap assumptions on the region $\mathbf{GH}_{\uboot}$, 
		which is globally hyperbolic 
		(its left lateral boundary is a portion of an integral curve of $\Lunit$ 
		and its right lateral boundary is a portion of an integral curve of $\uLunit$).
		By Remark~\ref{R:STRICTIMPROVEMENTOFBOOTSTRAP}, in the case $\eps := \datasize^{3/4}$,
		we see that no bootstrap assumptions are saturated on the bootstrap region $\mathbf{GH}_{\uboot}$,
		except possibly $\inf_{\mathbf{GH}_{\uboot}} \upmu = 0$.
		If the latter scenario does not occur, 
		then by
		Prop.\,\ref{P:CHOVDIFFEOMORPHISM}, the solution and its up-to-order $3$ derivatives with respect to the original coordinate system $(t,r)$
		are uniformly bounded on $\mathbf{GH}_{\uboot}$, and then standard continuation criteria 
		imply that for some $\Delta > 0$, the solution extends classically to the
		larger globally hyperbolic region
		$\mathbf{GH}_{\uboot - \Delta} := \lbrace (t,u) \ | \ t \in [\tstar,\mathfrak{t}(\uboot - \Delta)], 
		\, u \in (\mathfrak{u}_{\uboot - \Delta}(t),\tstar/2) \rbrace$,
		which strictly contains $\mathbf{GH}_{\uboot}$, and
		on which the bootstrap assumptions hold.
		It follows that the solution exists classically with respect to both coordinate systems
		and satisfies the estimates stated in Remark~\ref{R:STRICTIMPROVEMENTOFBOOTSTRAP}
		on the open-at-the-right-side region $\mathbf{GH}_{\umin}^{{\textnormal{Ext}}}$,
		where $\umin = \inf \lbrace u \in (\ucrease,\tstar/2) \ | \ \inf_{\mathbf{GH}_u} \upmu > 0 \rbrace$.
		Note that because $\upmu(\Tcrease,\ucrease) = 0$,
		we have $\umin \geq \ucrease$ and $\inf_{\mathbf{GH}_{\umin}} \upmu = 0$.
		Note also that by the definition of $\umin$,
		for any $u \in (\umin,\tstar/2)$, we have $\inf_{\mathbf{GH}_u} \upmu > 0$.
		For the same reasons given earlier in the proof,
		with respect to the $(t,u)$-coordinates, 
		$\RRiemann$ and $\LRiemann$
		extend to the closure of
		$\mathbf{GH}_{\umin}$ as solutions having $C^3$ regularity with respect to $(t,u)$, 
		$\upmu$ extends to the closure of
		$\mathbf{GH}_{\umin}$ as a solution with $C^2$ regularity.
		In the next part of the proof, we will show that $\umin = \ucrease$.

      \begin{figure} 
			\centering
			\begin{overpic}[scale=1, grid = false, tics=3, trim=-.5cm -1cm -1cm -.5cm, clip]{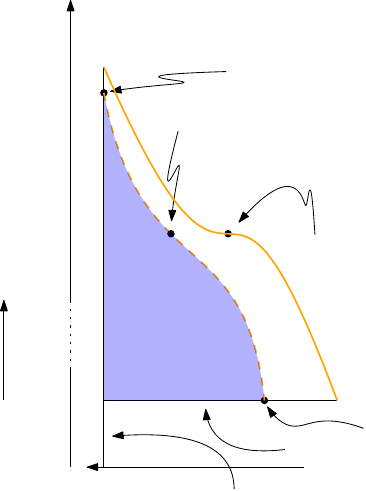} 
				\put (40,7) {\small{$\{ u = \tfrac{\tstar}{2}\}$}}
				\put (56,16) {\small{$\Sigma_{\tstar}$}}
				\put (69,20) {\small{$(\tstar,\mathfrak{u}_{\uboot}(\tstar))$}}
				\put (52,51) {\small{$\crease = (\Tcrease,\ucrease)$}}
				\put (28,74) {\small{$(\Tcrease,\uboot)$}}
				\put (26,11) {\small{$u$}}
				\put (1,33) {\small{$t$}}
				\put (46,83) {\small{$(\mathfrak{t}(\uboot),\tfrac{\tstar}{2})$}}
			\end{overpic}
			\caption{The region $\mathbf{GH}_{\uboot}$, which is the curved triangular open-at-the-right region, where the right lateral boundary is the integral curve of $\uLunit$ passing through the point $(\Tcrease,\uboot)$. The upper left vertex is the point $(t= \mathfrak{t}(\uboot),u =\tfrac{\tstar}{2})$, and the bottom right vertex is the point $(t=\tstar,u =\mathfrak{u}_{\uboot}(\tstar))$. We displayed a dotted portion of the $t$-axis to indicate that $\Tcrease  \approx\exp(C \mathring \upalpha^{-1})$ is significantly larger than $\tstar = \mathring \upalpha^{-1}$. Note that this curved triangle is rather narrow, especially in view of the omitted dotted portion of the $t$-axis, since $\ucrease = \mylittleo(\datasize)$. 
			Finally, the solid orange curve passing through the crease, which is the parametrized curve $\{(t,\mathfrak{u}_{\ucrease}(t)) \, | \, t \in [\tstar,\mathfrak{t}(\ucrease)]\}$, contains a portion of the Cauchy horizon and is derived through continuation criteria; 
			see the discussion leading up to \eqref{E:CAUCHYHORIZONPIECEDEVELOPMENT}. 
			We note that estimates \eqref{E:CAUCHYHORIZONTIMEFUNCTIONINITIALVALUEPROBLEM}--\eqref{E:THIRDDERIVATIVECAUCHYHORIZONTIMEFUNCTIONINITIALVALUEPROBLEM} imply that this orange curve is a cubic graph over $u$ with an inflection point at 
			$(\Tcrease,\ucrease)$.}
		\end{figure}
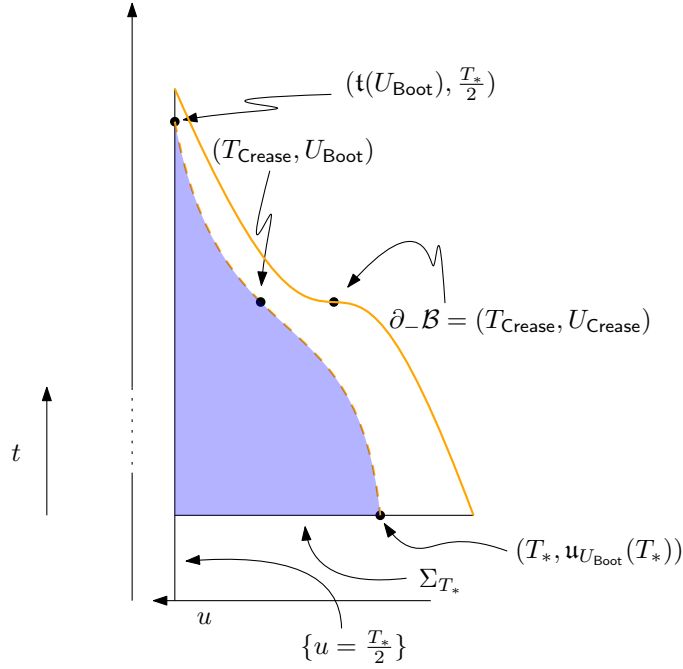 \label{F:TRIANGLEBOOTSTRAPREGION}

		\medskip
		
		\noindent \underline{\textbf{The structure of the Cauchy horizon and the parameterization of the boundary}}.
		Note that as subsets of spacetime, the integral curves of $\muuLunit = \upmu \uLunit$ and $\uLunit$ are the same,
		and that $\muuLunit = \upmu \frac{\partial}{\partial t} + 2 \frac{\partial}{\partial u}$
		is $C^2$ with respect to $(t,u)$ on the closure of $\mathbf{GH}_{\umin}$.
		Hence, by construction in the previous paragraph, the right boundary of
		the closure of $\mathbf{GH}_{\umin}$
		is a portion of the integral curve of $\muuLunit$ that emanates from the point with geometric coordinates $(\Tcrease,\umin)$,
		and $\upmu$ must vanish somewhere on this right boundary.
		Let us denote this integral curve portion by $\breve{\underline{\upgamma}}_{(\textsf{Min})}$.
		$\upmu$ cannot vanish at the upper endpoint $(\mathfrak{t}(\umin),\tstar/2)$ of $\breve{\underline{\upgamma}}_{(\textsf{Min})}$, 
		for otherwise, the $\upmu$-estimates of Prop.\,\ref{P:SHARPESTIMATESFORMU} would imply that
		$\muuLunit \upmu$ is strictly positive at $(\mathfrak{t}(\umin),\tstar/2)$, and thus $\upmu$ would turn 
		strictly negative along $\breve{\underline{\upgamma}}_{(\textsf{Min})}$ to the past of 
		$(\mathfrak{t}(\umin),\tstar/2)$.
		By continuity, this would imply that there is a number $U > \umin$ 
		such that $\upmu$ vanishes in $\mathbf{GH}_U$, a contradiction.
		$\upmu$ cannot vanish at the lower endpoint $(\tstar,\mathfrak{u}_{\min}(\tstar))$ of $\breve{\underline{\upgamma}}_{(\textsf{Min})}$,
		since this point lies in the data hypersurface, where $\upmu \approx 1$.
		Hence, the zero of $\upmu$ along  $\breve{\underline{\upgamma}}_{(\textsf{Min})}$ 
		must occur at an interior point $q$
		of the integral curve.
		Let us denote the geometric coordinates 
		of $q$ by $(t',u')$. Then $\upmu$ must have a minimum at $q$, i.e., $\upmu(t',u') = \muuLunit \upmu(t',u') = 0$, and 
		$\muuLunit \muuLunit \upmu(t',u') \geq 0$. Moreover, we must have $t' \geq \Tcrease$ and $u' \geq \ucrease$ since otherwise, 
		we would have $q \in \mathbf{GH}^{{\textnormal{Singular;Ext}}}$, 
		which is impossible because $\upmu > 0$ holds everywhere in $\mathbf{GH}^{{\textnormal{Singular;Ext}}}$
		($\upmu$ vanishes only along part of its top boundary, specifically the curve $t \rightarrow (t,\ublowup(t))$,
		which is not part of $\mathbf{GH}^{{\textnormal{Singular;Ext}}}$).
		The estimates of Prop.\,\ref{P:SHARPESTIMATESFORMU} imply that $u' = \mylittleo(\datasize)$ and
		$6 < \muuLunit \muuLunit \upmu(t',u') < 7$.
		By construction, $\breve{\underline{\upgamma}}_{(\textsf{Min})}$ intersects $\lbrace t = \Tcrease \rbrace$ at the point
		with geometric coordinates $(\Tcrease,\umin)$, where $\ucrease \leq \umin \leq u'$.
		Since $u' = \mylittleo(\datasize)$ and $\ucrease = \mylittleo(\datasize)$, we also have $\umin = \mylittleo(\datasize)$.
		We will now argue by contradiction that $\ucrease < \umin$ is impossible, and thus $\umin = \ucrease$.
		If $\ucrease < \umin$, then on the one hand, we must have $\muuLunit \upmu(\Tcrease,\umin) < 0$,
		in view of the fact that $\muuLunit \upmu(t',u') = 0$, the fundamental theorem of calculus on $\breve{\underline{\upgamma}}_{(\textsf{Min})}$,
		and the fact that the estimates of Prop.\,\ref{P:SHARPESTIMATESFORMU}
		imply that $6 < \muuLunit \muuLunit \upmu < 7$ at all points on 
		$\breve{\underline{\upgamma}}_{(\textsf{Min})}$ in between $(t',u')$ and $(\Tcrease,\umin)$.
		On the other hand, 
		we must have $\muuLunit \upmu(\Tcrease,\umin) > 0$,
		in view of the fact that $\muuLunit \upmu(\Tcrease,\ucrease) = 0$ (see \eqref{E:MUULUNITUPMUVANISHESATCREASE}), 
		the fundamental theorem of calculus along the integral curve of $\muX = \frac{\partial}{\partial u}$ that is tangent to $\lbrace t = \Tcrease \rbrace$,
		and the fact that the estimates of Prop.\,\ref{P:SHARPESTIMATESFORMU}
		imply that $3 < \muX \muuLunit \upmu(\Tcrease,u) < 5$ for all $u \in [\ucrease,u']$.
		We have demonstrated a contradiction, thereby showing that $\umin = \ucrease$.
		
		We now parameterize the portion of the integral curve of $\breve{\underline{\upgamma}}_{(\textsf{Min})}$ that lies to the future of
		$(\Tcrease,\ucrease)$ by $u$, and we let
		$\TCH(u)$ denote the value of $t$ along that integral curve. Since 
		$\muuLunit = \upmu \frac{\partial}{\partial t} + 2 \frac{\partial}{\partial u}$,
		$\TCH$ is the solution to the following initial value problem on the domain $[\ucrease,\tstar/2]$:
		\begin{align} \label{E:CAUCHYHORIZONTIMEFUNCTIONINITIALVALUEPROBLEM}
			\frac{d}{d u} \TCH(u)
			& = \frac{1}{2} \upmu(\TCH(u),u),
			& 
			\TCH(\ucrease) & = \Tcrease.
		\end{align}
		The solution $\TCH(u)$ to \eqref{E:CAUCHYHORIZONTIMEFUNCTIONINITIALVALUEPROBLEM} is $C^3$ because
		is $\upmu \in C^2$. Moreover, by the chain rule, we have:
		\begin{align} 
		\begin{split} \label{E:SECONDDERIVATIVECAUCHYHORIZONTIMEFUNCTIONINITIALVALUEPROBLEM}
			\frac{d^2}{d u^2} \TCH(u)
			& = \frac{1}{2} [\frac{\partial}{\partial t} \upmu](\TCH(u),u) 
					\frac{d}{d u} \TCH(u)
					+
					\frac{1}{2}
					[\frac{\partial}{\partial u} \upmu](\TCH(u),u) 
						\\
				& = \frac{1}{4} \upmu(\TCH(u),u) [\frac{\partial}{\partial t} \upmu](\TCH(u),u) 
					+ 
					\frac{1}{2}
					[\frac{\partial}{\partial u} \upmu](\TCH(u),u) 
					= \frac{1}{4} [\muuLunit \upmu](\TCH(u),u).
		\end{split}
		\end{align}
		The arguments we gave in the previous paragraph show that on the domain $[\ucrease,\tstar/2]$,
		$\upmu(\TCH(u),u)$ vanishes only when $u = \ucrease$.
		This shows that $\upmu$ is non-zero on the Cauchy horizon portion $\underline{\mathcal{C}}^{\tstar/2}$
		defined in \eqref{E:CAUCHYHORIZONUPTOTSTAROVER2}, since this portion does not include the point
		$(\Tcrease,\ucrease)$.
		The properties of the functions
		$\Tblowup$,
		$\TCH$,
		and $\tboundary$
		stated in the proposition therefore follow
		in a straightforward fashion from the results we have derived,
		including \eqref{E:SECONDDERIVATIVECAUCHYHORIZONTIMEFUNCTIONINITIALVALUEPROBLEM} 
		and the identities
		$\upmu(\Tcrease,\ucrease) = \muuLunit \upmu(\Tcrease,\ucrease) = 0$.
		In particular, \eqref{E:SHARPUPPERANDLOWERBOUNDSONVALUEOFTATTOPPOINTOFUEQUALSTSTAROVER2} follows from
		integrating \eqref{E:CAUCHYHORIZONTIMEFUNCTIONINITIALVALUEPROBLEM} and using 
		the bound $\upmu \leq 1 + C \datasize$ afforded by
		\eqref{E:MUPOINTWISEESTIMATEEXTERIOR} and Remark~\ref{R:STRICTIMPROVEMENTOFBOOTSTRAP}.
		
		For future use, we differentiate \eqref{E:SECONDDERIVATIVECAUCHYHORIZONTIMEFUNCTIONINITIALVALUEPROBLEM} and compute that:
		\begin{align} 
		\begin{split} \label{E:THIRDDERIVATIVECAUCHYHORIZONTIMEFUNCTIONINITIALVALUEPROBLEM}
			\frac{d^3}{d u^3} \TCH(u)
			& = \frac{1}{8} [\muuLunit \muuLunit \upmu](\TCH(u),u).
		\end{split}
		\end{align}
		
		Also for future use, we set: 
		\begin{subequations}
		\begin{align}
			\mathbf{GH}^{\textsf{CH;Ext}} 
			&:= \mathbf{GH}_{\ucrease},
			 \label{E:CAUCHYHORIZONPIECEDEVELOPMENT}	\\
			\mathbf{MGHD}^{\textnormal{Ext}} 
			& := \mathbf{GH}^{{\textnormal{Singular;Ext}}} \cup \mathbf{GH}^{\textnormal{CH;Ext}}.
			\label{E:EXTERIORDEVELOPMENTNOBOUNDARY}
		\end{align}
		\end{subequations}
		
			\begin{figure} 
			\centering
			\begin{overpic}[scale=1, grid = false, tics=3, trim=-.5cm -1cm -1cm -.5cm, clip]{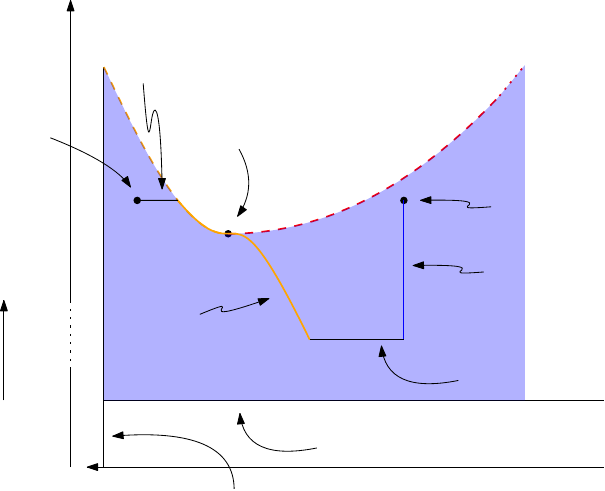} 
				\put (34,6) {\small{$\{ u = \tfrac{\tstar}{2}\}$}}
				\put (51,14) {\small{$\Sigma_{\tstar}$}}
				\put (29,60) {\small{$\crease = (\Tcrease,\ucrease)$}}
				\put (22,9) {\small{$u$}}
				\put (1,33) {\small{$t$}}
                \put (4,60) {\small{$(t,u_1)$}}
                \put (24,68) {\small{$h_1$}}
                \put (71,23) {\small{$h_2$}}
                \put (75,39) {\small{$\ell$}}
                \put (31,33) {\small{$\underline{\ell}$}}
                \put (76,49) {\small{$(t,u_2)$}}
			\end{overpic}
			\caption{The region $\mathbf{MGHD}^{\textnormal{Ext}}$ is the shaded blue region depicted in $(t,u)$-coordinates. The portion of the Cauchy horizon in this region is the dashed orange line, while the singular boundary is the dashed red line and extends out to $u$-spatial negative infinity.  These dashed curves, as well as the crease, are \emph{not} included in  $\mathbf{MGHD}^{\textnormal{Ext}}$, though they are part of its closure $\overline{\mathbf{MGHD}^{\textnormal{Ext}}}$.
            }
			\label{F:FULLEXTERIORGHD}
		\end{figure}

		$\Upsilon(\mathbf{MGHD}^{\textnormal{Ext}})$ 
		(i.e., $\mathbf{MGHD}^{\textnormal{Ext}}$ in n $(t,r)$ coordinates)
		is the portion of the ``Exterior'' region
		in Fig.\,\ref{F:MGHD} that lies above $\Sigma_{\tstar}$.
		
		The arguments we have given thus far show that the solution exists classically with respect to the $(t,u)$ coordinates on the closure of
		$\mathbf{MGHD}^{\textnormal{Ext}}$, which we denote by:
		\begin{align} \label{E:CLOSUREEXTERIORDEVELOPMENT}
			\overline{\mathbf{MGHD}^{\textnormal{Ext}}},
		\end{align}
		whose future-boundary is the $C^1$ curve $\futuresinghyp \cup \underline{\mathcal{C}}^{\tstar/2}$,
		where $\futuresinghyp$ is the singular boundary defined in \eqref{E:SINGULARBOUNDARY}, and
		$\underline{\mathcal{C}}^{\tstar/2}$ is the Cauchy horizon portion defined in \eqref{E:CAUCHYHORIZONUPTOTSTAROVER2}.
		We have also shown that $\upmu$ is positive everywhere on 
		$\overline{\mathbf{MGHD}^{\textnormal{Ext}}}$, except on $\futuresinghyp$, where it vanishes.
		
		\medskip

		\noindent \underline{\textbf{The behavior of $\Upsilon$ on $\mathbf{MGHD}^{\textsf{Ext}}$ and $\overline{\mathbf{MGHD}^{\textsf{Ext}}}$}}.
		Prop.\,\ref{P:CHOVDIFFEOMORPHISM} shows that $\Upsilon \restriction_{\overline{\mathbf{GH}^{{\textnormal{Singular;Ext}}}}}$
		is a homeomorphism and $\Upsilon \restriction_{\overline{\mathbf{GH}^{\textsf{CH;Ext}}}}$ is a homeomorphism.
		We will show that $\Upsilon$ is injective on the closure of the union of these two regions, 
		i.e., the set $\overline{\mathbf{MGHD}^{\textnormal{Ext}}}$ from
		\eqref{E:CLOSUREEXTERIORDEVELOPMENT}. This will
		in particular imply that $\Upsilon$ is a diffeomorphism on the subset 
		$\mathbf{MGHD}^{\textsf{Ext}} \cup \underline{\mathcal{C}}^{\tstar/2}$ of $\overline{\mathbf{MGHD}^{\textnormal{Ext}}}$, 
		on which $\upmu > 0$.
		We have already shown that $\Upsilon$ extends to 
		$\overline{\mathbf{MGHD}^{\textnormal{Ext}}}$ as a $C^3$ function of $(t,u)$. 
		To show the desired injectivity, we
		assume that $(t,u_2) \in \overline{\mathbf{GH}^{{\textnormal{Singular;Ext}}}}\backslash \overline{\mathbf{GH}^{\textsf{CH;Ext}}}$ and 
		$(t,u_1) \in \overline{\mathbf{GH}^{\textsf{CH;Ext}}}$.
		We claim that $u_1 > u_2$. To see this, denote the set of $u$-values of the curve $\breve{\underline{\upgamma}}_{(\textsf{Min})}$ defined above by $S$. Since this curve is the right boundary of  $\overline{\mathbf{GH}^{\textsf{CH;Ext}}}$, it follows that $u_1 \ge u$ for all $u \in S$. On the other hand, $\breve{\underline{\upgamma}}_{(\textsf{Min})}$ is the left boundary of (but not contained in)  $\overline{\mathbf{GH}^{{\textnormal{Singular;Ext}}}}\backslash \overline{\mathbf{GH}^{\textsf{CH;Ext}}}$ and so $u_2 < u$ for all $u \in S$, which was the desired claim. 
		Fix an arbitrary $u_* \in (u_2,\ucrease]$. We can then join $(t,u_1)$ to $(t,u_2)$
		by a continuous, piecewise smooth curve in $\overline{\mathbf{MGHD}^{\textnormal{Ext}}}$
		that starts at $(t,u_1)$ and consists of the following segments:\footnote{Clearly, there are infinitely many piecewise smooth paths one could take from $(t,u_1)$ to $(t,u_2)$. Our order of $h_1 \to \underline{\ell} \to h_2 \to \ell$ is merely a convenient way to show that $r(t,u_1) < r(t,u_2)$.}
		\begin{itemize}
			\item A horizontal line segment $h_1$, along which $t$ is constant as $u$ decreases.
			\item A (possibly empty) portion $\underline{\ell}$ of the integral curve portion $\breve{\underline{\upgamma}}_{(\textsf{Min})}$ 
				from earlier in the proof (which strictly contains the Cauchy horizon portion $\underline{\mathcal{C}}^{\tstar/2}$), 
				along which $t$ decreases as $u$ decreases,
				until stopping at the unique point on $\breve{\underline{\upgamma}}_{(\textsf{Min})}$  whose $u$-value is the fixed $u_* \in (u_2,\ucrease]$
            \item Another (possibly empty) horizontal line $h_2$, along which $t$ is constant as $u$ decreases from $u_*$ until reaching $u_2$.
			\item A (possibly empty) portion $\ell$ of an integral curve of $\Lunit$, along which $u$ is constantly $u_2$ while $t$ increases up to its original
			value from the starting point of $(t,u_1)$.
		\end{itemize}
		Once we show that the coordinate function $r$ strictly increases along all $3$ kinds of segments as one moves from $(t,u_1)$ towards $(t,u_2)$
		(i.e., as we decrease $u$),
		we will have shown that $\Upsilon(t,u_1) \neq \Upsilon(t,u_2)$, which yields the desired injectivity.
		The increasing behavior of $r$ along $h_1$ or $h_2$ (as $u$ decreases) was shown in the proof of Prop.\,\ref{P:CHOVDIFFEOMORPHISM}.
		The fact that $r$ increases along $\underline{\ell}$ as $u$ decreases follows from the
		relations $\muuLunit r = \upmu (v^r - \Speed) \partial_r r = \upmu (v^r - \Speed) \approx - \upmu$ and $\muuLunit u = 2$, and the fact that
		$\upmu > 0$ on $\underline{\ell}$ except at a single point, namely the crease $(\Tcrease,\ucrease)$.
		The fact that $r$ increases along $\ell$ as $t$ increases follows easily from the relation
		$\Lunit r = \Speed + v^r \approx 1$.

        \medskip

        \underline{\textbf{Continuous extension of $\RRiemann \circ \Upsilon^{-1}$ and $\LRiemann \circ \Upsilon^{-1}$ to 
				$\Upsilon(\overline{\mathbf{MGHD}^{\textsf{Ext}}})$}}. Above we proved that relative to the $(t,u)$ coordinates, 
				$(\RRiemann,\LRiemann)$ extend as $C^3$ solutions to $\overline{\mathbf{MGHD}^{\textsf{Ext}}}$. It then follows immediately that 
				$(\RRiemann \circ \Upsilon^{-1},\LRiemann \circ \Upsilon^{-1})$
				extend as continuous functions to $\Upsilon(\overline{\mathbf{MGHD}^{\textsf{Ext}}})$ in $(t,r)$-coordinates, which includes $\Upsilon(\futuresinghyp)$, since $\Upsilon^{-1}: \Upsilon(\overline{\mathbf{MGHD}^{\textsf{Ext}}}) \to \overline{\mathbf{MGHD}^{\textsf{Ext}}}$
				is a homeomorphism.

		\medskip
		
		\noindent \underline{\textbf{The blowup of $\partial_r [\RRiemann \circ \Upsilon^{-1}]$ along $\Upsilon(\futuresinghyp)$}}.
		Let $(t_{\star},r_{\star}) \in \Upsilon(\futuresinghyp)$, and let $\lbrace (t_n,r_n) \rbrace_{n \in \mathbb{N}}$ be a sequence of points in
	$\Upsilon(\mathbf{MGHD}^{\textsf{Ext}} \cup \underline{\mathcal{C}}^{\tstar/2})$ such that $\lim_{n \to \infty} (t_n,r_n) = (t_{\star},r_{\star})$.
	Since $\Upsilon$ is a homeomorphism on $\overline{\mathbf{MGHD}^{\textsf{Ext}}}$ and
	$\upmu$ vanishes along $\futuresinghyp$, it follows that:
	\begin{align} \label{E:SEQUENCEOFMUVALUESHASVANISHINGLIMITASSINGULARBOUNDARYISAPPROACHED}
		\lim_{n \to \infty} \upmu \circ \Upsilon^{-1}(t_n,r_n)
		& = \upmu \circ \Upsilon^{-1}(t_{\star},r_{\star}) = 0.
	\end{align}
	Moreover, using that $\Speed = 1 + \mathcal{O}(|\RRiemann| + |\LRiemann|)$ and using the pointwise 
	estimates $|\RRiemann| \lesssim \frac{\datasize}{1 + t + |u|}$ and $|\LRiemann| \lesssim \frac{\datasize}{1 + t + |u|}$
	afforded by Prop.\,\ref{P:RIEMANNINVARIANTSAPRIORIEXTERIORREGIONESTIMATES},
	we find that:
	\begin{align} \label{E:SEQUENCEOFSPEEDOFSOUNDVALUESHASNONZEROLIMITASSINGULARBOUNDARYISAPPROACHED}
		\lim_{n \to \infty} \Speed \circ \Upsilon^{-1}(t_n,r_n)
		& = \Speed \circ \Upsilon^{-1}(t_{\star},r_{\star}) > 0.
	\end{align}
	On the other hand, \eqref{E:TANDUBOUNDSFOREXTERIORBOOTSTRAPREGION} and the estimate \eqref{E:QUANTITATIVEMUXRPLUSUPPERANDLOWERBOUNDSWHENMUISSMALL} 
	show that:
	\begin{align} \label{E:SEQUENCEOFMUXRPLUSVALUESNONZEROLIMITASSINGULARBOUNDARYISAPPROACHED}
		\lim_{n \to \infty} [\muX \RRiemann] \circ \Upsilon^{-1}(t_n,r_n)
		& > 0.
	\end{align}
	Combining \eqref{E:SEQUENCEOFMUVALUESHASVANISHINGLIMITASSINGULARBOUNDARYISAPPROACHED}--\eqref{E:SEQUENCEOFMUXRPLUSVALUESNONZEROLIMITASSINGULARBOUNDARYISAPPROACHED} and appealing to the identity 
	$\muX = - \upmu \Speed \partial_r$ (see \eqref{E:MUWEIGHTEDVECTORFIELDS}), 
	we conclude that
	$\limsup_{n \to \infty} \left|\partial_r [\RRiemann \circ \Upsilon^{-1}(t_n,r_n)] \right| = \infty$,
	which is the desired blowup-result.
	
	\medskip 
	\noindent \underline{\textbf{The structure of $\Upsilon(\futuresinghyp)$ and the continuous tangential extension of $\Lunit$ to $\Upsilon(\futuresinghyp)$}}.
	The vectorfield $V := \frac{(\muX \upmu)}{\Lunit \upmu} \Lunit - \muX$ satisfies $V u = - 1$ and $V \upmu = 0$. In particular, $V$ is tangent to
	$\futuresinghyp$. Moreover, the arguments given just below \eqref{E:CLOSUREDBOOTSTRAPREGIONFORCONTROLLINGSOLUTIONUPTOSINGULARBOUNDARYATFIXEDLEFTUCHARACTERISTIC} show that
	along $\futuresinghyp$, we have $\Lunit \upmu < 0$ and 
	$\muX \upmu < 0$, except at the crease, where $\muX \upmu = 0$. We now fix any point $(t,r) \in \Upsilon(\futuresinghyp)$. Observe that the 
	pushforward of $V$ by $\Upsilon$ at $(t,r)$, i.e., $\mathrm{d} \Upsilon_{\Upsilon^{-1}(t,r)} V_{\Upsilon^{-1}(t,r)}$, is well-defined because 
	$\Upsilon$ is $C^3$ with a continuous inverse. Moreover, this pushforward-vector is tangent to $\Upsilon(\futuresinghyp)$. 
	Using \eqref{E:CHOVJACOBIAN}, we compute $d \Upsilon_{\Upsilon^{-1}(t,r)} V_{\Upsilon^{-1}(t,r)} = \frac{\muX \upmu}{\Lunit \upmu} \circ \Upsilon^{-1}(t,r) \left(\partial_t + ((v^r + \Speed - \upmu \Speed)\circ \Upsilon^{-1}(t,r))\partial_r\right) =\frac{\muX \upmu}{\Lunit \upmu} \circ \Upsilon^{-1}(t,r)\left( \partial_t + ((v^r + \Speed)\circ \Upsilon^{-1}(t,r))\partial_r\right)$, where the second identity follows from the fact that $\upmu = 0$ on $\futuresinghyp$. 
	In particular, everywhere along $\Upsilon(\futuresinghyp)$, the tangent vector to $\Upsilon(\futuresinghyp)$ is a multiple of $\Lunit \circ \Upsilon^{-1}= \partial_t + ( v^r + \Speed)\circ \Upsilon^{-1} \partial_r$, as is desired.
	
	To show that $\Upsilon(\futuresinghyp)$ is a $C^1$ curve-portion in $(t,r)$-space, we will first
	describe $\Upsilon(\futuresinghyp)$ as a $u$-parameterized curve $u \rightarrow (\Tblowup(u),\Rblowup(u))$ in $(t,r)$ space 
	for $u \in (-\infty,\ucrease]$, where we recall that $\Tblowup(u)$ is the Cartesian time of the unique point where 
	$\upmu \restriction_{[u,-\infty)}$ first vanishes, and $\Rblowup(u)$ is the value of $r$ at $\Upsilon((\Tblowup(u),u))$.
	In particular, $(\Tblowup(\ucrease),\Rblowup(\ucrease)) = (\Tcrease,\Rcrease)$, where $(\Tcrease,\Rcrease) = \Upsilon((\Tblowup(\ucrease),\ucrease))$.
	The map $u \rightarrow (\Tblowup(u),\Rblowup(u))$ is the image under $\Upsilon$ of the $C^2$ curve $\lbrace (t,u) \ | \ \upmu(t,u) = 0 \rbrace$,
	and thus $u \rightarrow (\Tblowup(u),\Rblowup(u))$ is a $C^2$ function of $u$ with a tangent vector $V$ satisfying $V u = -1$.
	We aim to show that along this curve, $r$ is a $C^1$ function of $t$. It suffices to show that:
	\begin{subequations}
	\begin{align}
		\frac{d}{du} \Tblowup(u)
		& = 0, && u = \ucrease,
			\label{E:DERIVATIVEOFTBLOWUPZEROATCREASE} \\
		\frac{d}{du} \Tblowup(u)
		& < 0, && u \in (-\infty,\ucrease),	
		\label{E:DERIVATIVEOFTBLOWUPNEGATIVEAWAYFROMCREASE} \\
		\frac{d^2}{du^2} \Tblowup(u)
		& > 0, && u = \ucrease,
		\label{E:SECONDDERIVATIVEOFTBLOWUPPOSITIVEATCREASE}
	\end{align}
	\end{subequations}
	
	\begin{subequations}
	\begin{align}
		\frac{d}{du} \Rblowup(u)
		& = 0, && u = \ucrease,
		\label{E:DERIVATIVEOFRBLOWUPZEROATCREASE}
			\\
		\frac{d}{du} \Rblowup(u)
		& < 0, && u \in (-\infty,\ucrease),	
		\label{E:DERIVATIVEOFRBLOWUPNEGATIVEAWAYFROMCREASE} \\
		\frac{d^2}{du^2} \Rblowup(u)
		& > 0, && u = \ucrease.
		\label{E:SECONDDERIVATIVEOFRBLOWUPPOSITIVEATCREASE}
	\end{align}
	\end{subequations}
	Note that $\frac{d}{du} \Tblowup(u) = V t$, $\frac{d^2}{du^2} \Tblowup(u) = V V t$, 
	$\frac{d}{du} \Rblowup(u) = V r$, and $\frac{d^2}{du^2} \Rblowup(u) = V V r$. 
	Hence, to prove \eqref{E:DERIVATIVEOFTBLOWUPZEROATCREASE}--\eqref{E:SECONDDERIVATIVEOFRBLOWUPPOSITIVEATCREASE},
	we start by computing:
	\begin{subequations}
	\begin{align}
		V t &
		= \frac{(\muX \upmu)}{\Lunit \upmu},
			\label{E:VTFORMULA} \\
		V V t
		& = \frac{(\muX \upmu)}{\Lunit \upmu} \Lunit \left( \frac{(\muX \upmu)}{\Lunit \upmu} \right)
		- \muX \left( \frac{(\muX \upmu)}{\Lunit \upmu} \right),
		\label{E:TWOVTFORMULA}
	\end{align}
	\end{subequations}
	and:
	\begin{subequations}
	\begin{align}
		V r &
		= \frac{(\muX \upmu)}{\Lunit \upmu} (v^r + \Speed) 
		+ \upmu \Speed,
			\label{E:VRFORMULA} \\
		V V r
		& = \frac{(\muX \upmu)}{\Lunit \upmu} \Lunit \left( \frac{(\muX \upmu)}{\Lunit \upmu} (v^r + \Speed)  \right)
		- 
		\muX \left(  \frac{(\muX \upmu)}{\Lunit \upmu} (v^r + \Speed)  \right)
		+ 
		\frac{(\muX \upmu)}{\Lunit \upmu} \Lunit \left( \upmu \Speed \right)
		- \muX ( \upmu \Speed).
		\label{E:TWOVRFORMULA}
	\end{align}
	\end{subequations}
	In our above analysis, we showed that:
	\begin{subequations}
	\begin{align}
		\upmu  \restriction_{\futuresinghyp} & = 0,
			\label{E:MUVANISHESONSINGULARBOUNDARY} \\
		\Lunit \upmu \restriction_{\futuresinghyp} & < 0,
			\\
		\muX \upmu \restriction_{\futurecrease} & = 0,
			\\
		\muX \upmu \restriction_{\futuresinghyp \backslash \futurecrease} & < 0,
			\\
		\muX \muX \upmu \restriction_{\futurecrease} & > 0.
		\label{E:TRANSVERSALCONVEXITYATCREASE}
	\end{align}
	\end{subequations}
	
	Using \eqref{E:VTFORMULA}--\eqref{E:TRANSVERSALCONVEXITYATCREASE}, the estimates of Prop.\,\ref{P:RIEMANNINVARIANTSAPRIORIEXTERIORREGIONESTIMATES}, 
	and the simple estimate
	$v^r + \Speed = 1 + \frac{\mathcal{O}(\datasize)}{1 + t + r}$,
	we compute that:
	\begin{align}
		V t \restriction_{\futurecrease}  & = 0,
			\label{E:VTZEROATCREASE} \\
		V t 	\restriction_{\futuresinghyp \backslash \futurecrease}  & > 0,
		\label{E:VTPOSITIVEAWAYFROMCREASE}
	\end{align}
	
	\begin{align}
		V V t \restriction_{\futurecrease}  & > 0,
		\label{E:TWICEVTPOSITIVEATCREASE}
	\end{align}
	
	\begin{align}
		V r \restriction_{\futurecrease}  & = 0,
			\label{E:VRZEROATCREASE} \\
		V r 	\restriction_{\futuresinghyp \backslash \futurecrease}  & > 0,
		\label{E:VRPOSITIVEAWAYFROMCREASE}
	\end{align}
	
	\begin{align}
		V V r \restriction_{\futurecrease}  & > 0,
		\label{E:TWICEVRPOSITIVEATCREASE}
	\end{align}
	which yields the desired results \eqref{E:DERIVATIVEOFTBLOWUPZEROATCREASE}--\eqref{E:SECONDDERIVATIVEOFRBLOWUPPOSITIVEATCREASE}.
	
	\medskip

	\noindent \underline{\textbf{The structure of $\Upsilon(\underline{\mathcal{C}}^{\tstar/2})$}}.
	The analysis here is similar to the analysis for $\Upsilon(\futuresinghyp)$, so we will be terse.
	We let $\overline{\underline{\mathcal{C}}^{\tstar/2}}$ denote the closure of $\underline{\mathcal{C}}^{\tstar/2}$.
	In particular, $\overline{\underline{\mathcal{C}}^{\tstar/2}}$ contains its past boundary point $\futurecrease$.
	We have already shown that $\muuLunit$ is tangent to $\overline{\underline{\mathcal{C}}^{\tstar/2}}$ and that $\upmu$ is everywhere positive on 
	$\underline{\mathcal{C}}^{\tstar/2}$. Hence,
	the pushforward vectorfield $d \Upsilon \cdot \uLunit$ is tangent to $\Upsilon(\underline{\mathcal{C}}^{\tstar/2})$ in $(t,r)$-space.
	To show that $\Upsilon(\underline{\mathcal{C}}^{\tstar/2})$ is a $C^1$ curve-portion in $(t,r)$-space, we will 
	describe $\Upsilon(\overline{\underline{\mathcal{C}}^{\tstar/2}})$ as a $u$-parameterized curve $u \rightarrow (\TCH(u),\RCH(u))$ in $(t,r)$ space 
	for $u \in [\ucrease,\frac{\tstar}{2}]$, where $\TCH$ is the function 
	from earlier in the proof, and $\RCH(u)$ is the value of $r$ at $\Upsilon((\TCH(u),u))$.
	It suffices to show that:
	\begin{subequations}
	\begin{align}
		\frac{d}{du} \TCH(u)
		& = 0, && u = \ucrease,
			\label{E:DERIVATIVEOFTCHZEROATCREASE} 
				\\
		\frac{d}{du} \TCH(u)
		& > 0, && u \in (\ucrease,\frac{\tstar}{2}],
			\label{E:DERIVATIVEOFTCHPOSITIVEAWAYFROMCREASE} 
				\\
	 \frac{d^2}{du^2} \TCH(u)
		& = 0, && u = \ucrease,
		\label{E:SECONDDERIVATIVEOFTCHZEROATCREASE}
			\\
		\frac{d^2}{du^2} \TCH(u)
		& > 0, && u = \ucrease,
		\label{E:THIRDDERIVATIVEOFTCHPOSITIVEATCREASE}
	\end{align}
	\end{subequations}
	
	\begin{subequations}
	\begin{align}
		\frac{d}{du} \RCH(u)
		& = 0, && u = \ucrease,
			\label{E:DERIVATIVEOFRCHZEROATCREASE} 
			\\
		\frac{d}{du} \RCH(u)
		& < 0, && u \in (\ucrease,\frac{\tstar}{2}],
			\label{E:DERIVATIVEOFRCHNEGATIVEAWAYFROMCREASE}
				\\
	 \frac{d^2}{du^2} \RCH(u)
		& > 0, && u = \ucrease,
		\label{E:SECONDDERIVATIVEOFRCHZEROATCREASE}
			\\
		\frac{d^3}{du^3} \RCH(u)
		& > 0, && u = \ucrease,
		\label{E:THIRDDERIVATIVEOFRCHPOSITIVEATCREASE}
	\end{align}
	\end{subequations}
	To prove \eqref{E:DERIVATIVEOFTCHZEROATCREASE}--\eqref{E:THIRDDERIVATIVEOFRCHPOSITIVEATCREASE},
	we start by combining 
	\eqref{E:CAUCHYHORIZONTIMEFUNCTIONINITIALVALUEPROBLEM}--\eqref{E:THIRDDERIVATIVECAUCHYHORIZONTIMEFUNCTIONINITIALVALUEPROBLEM}
	with the following companion formulas, which can be derived through similar arguments based on the fact that
	$\frac{1}{2} \muuLunit u = 1$ and $\frac{1}{2} \muuLunit r = \frac{1}{2} \upmu (v^r - \Speed)$:
	\begin{subequations}
	\begin{align} \label{E:CAUCHYHORIZONRADIALFUNCTIONINITIALVALUEPROBLEM}
			\frac{d}{d u} \RCH(u)
			& = \frac{1}{2} \upmu (v^r - \Speed) \restriction_{(\TCH(u),u)},
			& 
			\RCH(\ucrease) & = \Rcrease,
				\\
			\frac{d^2}{d u^2} \RCH(u)
			& = \frac{1}{4} \muuLunit \left(\upmu(\TCH(u),u) [v^r - \Speed] \right) \restriction_{(\TCH(u),u)},
				 \label{E:SECONDDERIVATIVECAUCHYHORIZONRADIALFUNCTIONINITIALVALUEPROBLEM} \\
			\frac{d^3}{d u^3} \RCH(u)
			& = \frac{1}{8} \muuLunit \muuLunit \left(\upmu(\TCH(u),u) [v^r - \Speed] \right) \restriction_{(\TCH(u),u)}.
				\label{E:THIRDDERIVATIVECAUCHYHORIZONRADIALFUNCTIONINITIALVALUEPROBLEM}
		\end{align}
		\end{subequations}

	The desired results \eqref{E:DERIVATIVEOFTCHZEROATCREASE}--\eqref{E:THIRDDERIVATIVEOFRCHPOSITIVEATCREASE}
	follow from  \eqref{E:DERIVATIVEOFTCHZEROATCREASE}--\eqref{E:THIRDDERIVATIVEOFRCHPOSITIVEATCREASE},
	\eqref{E:CAUCHYHORIZONRADIALFUNCTIONINITIALVALUEPROBLEM}--\eqref{E:THIRDDERIVATIVECAUCHYHORIZONRADIALFUNCTIONINITIALVALUEPROBLEM},
	and the following results, obtained earlier in the proof:
	\begin{subequations}
	\begin{align}
		\upmu  \restriction_{\futurecrease} & = 0,
			\label{E:MUVANISHESONCREASE} \\
			\upmu  \restriction_{\underline{\mathcal{C}}^{\tstar/2}} & > 0,
			\label{E:MUPOSITIVEONCAUCHYHORIZON} 
			\\
		\muuLunit \upmu \restriction_{\futurecrease} & = 0,
			\label{E:MUULUNITMUVANISHESONCREASE} 
				\\
		\muuLunit \muuLunit \upmu \restriction_{\futurecrease} & > 0.
		\label{E:MUULUNITTRANSVERSALCONVEXITYATCREASE}
	\end{align}
	\end{subequations}
	
	\medskip
	
	\noindent \underline{\textbf{Proof that $\Upsilon(\futuresinghyp) \cup \Upsilon(\underline{\mathcal{C}}^{\tstar/2})$ has a corner}}.
	Recall that $\futurecrease \in \futuresinghyp$ and that $\futurecrease$ is a past limit point of $\underline{\mathcal{C}}^{\tstar/2}$,
	i.e., its closure satisfies $\overline{\underline{\mathcal{C}}^{\tstar/2}} = \underline{\mathcal{C}}^{\tstar/2} \cup \lbrace \futurecrease \rbrace$.
	Earlier in the proof, we showed that the tangent vector to the curve $\Upsilon(\futuresinghyp)$
	at $\Upsilon(\futurecrease)$ is proportional to $\Lunit \restriction_{\Upsilon(\futurecrease)}$,
	and that the tangent vector to the curve $\Upsilon(\overline{\underline{\mathcal{C}}^{\tstar/2}})$
	at $\Upsilon(\futurecrease)$ is proportional to $\uLunit \restriction_{\Upsilon(\futurecrease)}$.
	The formulas \eqref{E:NULLVECTORFIELDS} and the
	estimates of Prop.\,\ref{P:RIEMANNINVARIANTSAPRIORIEXTERIORREGIONESTIMATES} together imply that $\Lunit$ and $\uLunit$ are everywhere transversal,
	which shows that the curve  $\Upsilon(\futuresinghyp) \cup \Upsilon(\underline{\mathcal{C}}^{\tstar/2})$ has a corner at $\futurecrease$.

	\medskip 
	\noindent \underline{\textbf{Global hyperbolicity}}.
	Finally, we show that the region
	$\mathbf{MGHD}^{\textsf{Ext}}$ from \eqref {E:GLOBALLYHYPERBOLICEXTERIOREGIONPLUSABITOFCAUCHYHORIZON}
	is globally hyperbolic in the sense of
	Definition~\ref{D:GLOBALLHYPERBOLICEXTERIORBOOTSTRAPREGION}.
	We have shown that $\mathbf{MGHD}^{\textsf{Ext}}$ 
	is foliated by integral curves of $\Lunit$ (i.e., portions of level sets of $u$) that connect points in $\Sigma_{\tstar}$ to points
	on its future-boundary. Moreover, given any past-directed integral curve of $\uLunit$
	emanating from a point in $\mathbf{MGHD}^{\textsf{Ext}}$, both $t$ and $u$ decrease towards $\Sigma_{\tstar}$ (since $\uLunit t = 1$ and $\muuLunit u = 2$).
	By ODE uniqueness, such an integral curve can never intersect the Cauchy horizon boundary $\underline{\mathcal{C}}^{\tstar/2}$
	(which does not belong to $\mathbf{MGHD}^{\textsf{Ext}}$),
	which is itself an integral curve of $\muuLunit$ where the solution is smooth in both coordinate systems.
	Hence, this integral curve terminates along $\Sigma_{\tstar}$ at some $u$-value that is more negative than its starting value.
\end{proof}

\section{The MGHD existence and uniqueness theorem for shock-forming solutions, including the Interior Region}
\label{S:MAINEXISTENCEANDUNIQUENESSTHEOREMINCLUDINGINTERIOR}
In this section, we prove our main theorem on the existence, uniqueness, and stability of the MGHD
for open sets of data that lead to shock-forming solutions. 
The most difficult aspects of the theorem concern the Exterior Region,
which we already treated in Prop.\,\ref{P:EXISTENCEUPTOCREASEANDSINGULARBOUNDARYANDPORTIONOFCAUCHYHORIZON}.
Hence, most of the effort in this section is dedicated towards controlling the solution in the Interior Region.
To that end, we will exploit the decay of the solution with respect to $t-r$ in the interior; despite the failure of the
null condition in the equations, the decay in $t-r$ will allow us to control the solution for an amount of time
$T_{\textsf{Max}}$ that is sufficient for capturing the entire MGHD.

\subsection{Statement of the MGHD existence and uniqueness theorem}
\label{SS:STATEMENTOFMAINEXISTENCETHEOREM}
In this section, we state our main theorem on MGHD existence and uniqueness. 
We provide the proof in Sect.\,\ref{SS:PROOFOFT:MAINEXISTENCETHEOREM}.

The full statement of the theorem refers to the \emph{acoustical metric} on $\mathbb{R}^{1+3}$, which is the following Lorentzian metric:
\begin{align} \label{E:MAINTHEOREMACOUSTICALMETRIC}
	\mathbf{g} 
		& := 
		- \mathrm{d} t \otimes \mathrm{d}t
			+ 
			\Speed^{-2} 
			\sum_{a=1}^3(\mathrm{d}x^a - v^a \mathrm{d}t) \otimes (\mathrm{d}x^a - v^a \mathrm{d}t).
\end{align}
The tensor $\mathbf{g}$ can be used to construct the characteristics.
More precisely, for the $3D$ compressible Euler equations, the characteristic subset in co-tangent space
includes the sound cones $\lbrace \upxi \ | \ (\mathbf{g}^{-1})^{\alpha \beta} \upxi_{\alpha} \upxi_{\beta} = 0 \rbrace$.
The causal structure of the solution and the correct intrinsic notion of global hyperbolicity are tied to $\mathbf{g}$.
We refer to Appendix~\ref{A:GLOBALHYPERBOLICITYANDMGHDS} for further discussion.

\begin{figure}  
\centering
\begin{overpic}[scale=.7, grid = false, tics=3, trim=-.5cm -1cm -1cm -.5cm, clip]{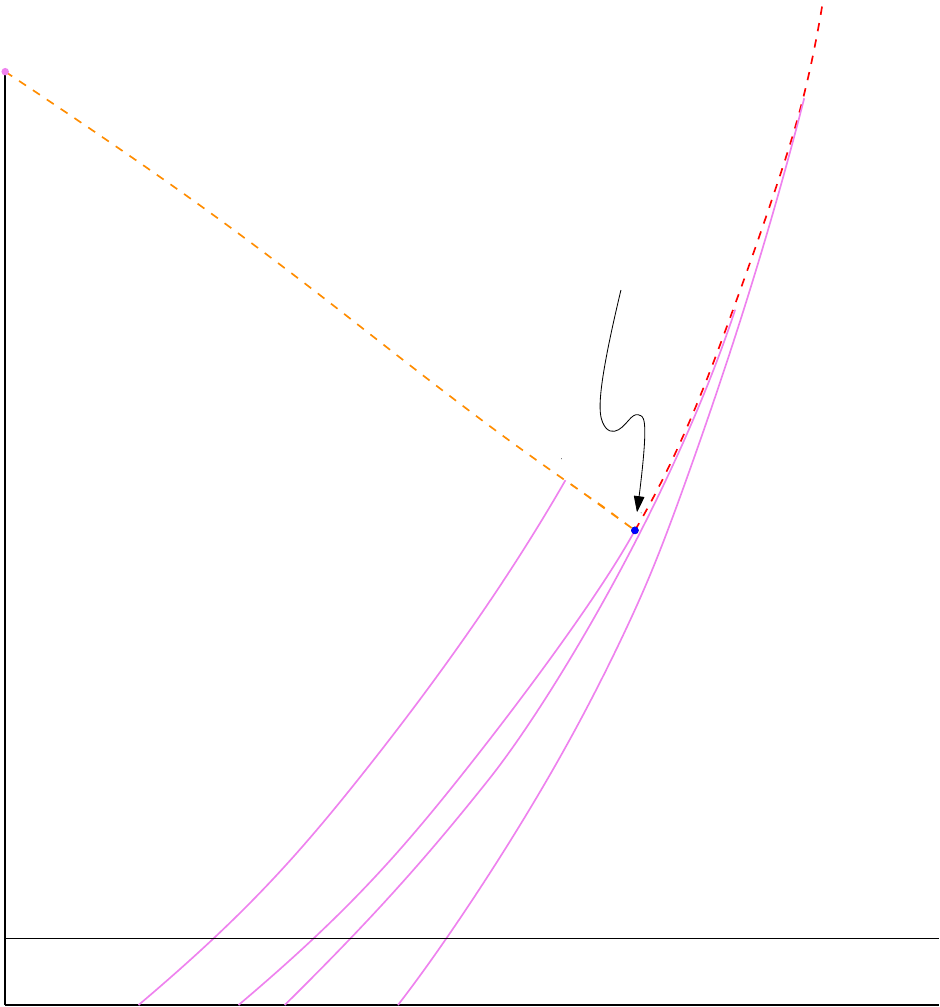}
	\put (0,93) {\small $(T_{\textsf{Max}},0)$}
	\put (35,69) {\small $\futureCauchyhor$}
	\put (73,91) {\small $\Upsilon(\futuresinghyp)$}
	\put (49,81) {\small $(\Tcrease,\rcrease)$}
	\put (56,77) {\small \rotatebox{90}{$=$}}
	\put (50,74) {\small $\Upsilon(\futurecrease)$}
	\put (20,50) {\small \textnormal{Interior}}
	\put (19,20) {\small $\lbrace u = \frac{\tstar}{2} \rbrace$}
	\put (48,40) {\small \textnormal{Exterior}}
	\put(5,3) 		{\vector(1,0){10}}
	\put(7.5,3.5)     {\small $r$}
	\put(2,9) 		{\vector(0,1){10}}
	\put(0,11)     {\small $t$}
	\put(49,13)      {\small $\Sigma_{\tstar}$}
	\put(49,6.5)      {\small $\Sigma_0$}
\end{overpic}
\caption{$\mathbf{M}^{\textsf{Future}}$ in $(t,r)$-coordinates. By definition, ``Exterior'' is the infinite region lying
above $\Sigma_{\tstar}$, to the right of the characteristic curve $\lbrace u = \frac{\tstar}{2} \rbrace$, and
below $\futureCauchyhor \cup \Upsilon(\futuresinghyp)$. 
It was constructed with respect to the $(t,u)$-coordinates in Prop.\,\ref{P:EXISTENCEUPTOCREASEANDSINGULARBOUNDARYANDPORTIONOFCAUCHYHORIZON}.
More precisely, the Exterior is equal to $\Upsilon(\mathbf{MGHD}^{\textsf{Ext}})$, where 
$\mathbf{MGHD}^{\textsf{Ext}}$ is the subset of $(t,u)$-coordinate space defined in \eqref{E:GLOBALLYHYPERBOLICEXTERIOREGIONPLUSABITOFCAUCHYHORIZON}.}
\label{F:MGHD}
  \end{figure}
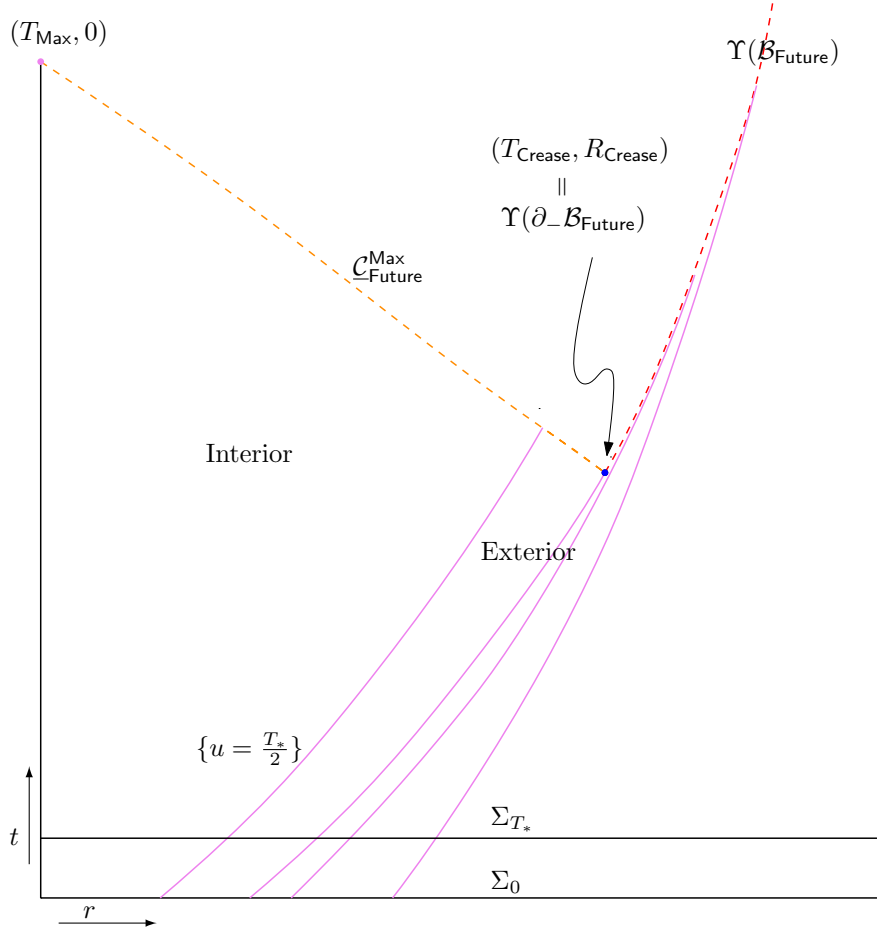

\begin{theorem}[Existence and uniqueness of MGHDs for shock-forming solutions]
\label{T:MAINMGHDEXISTENCETHEOREM}
Consider the compressible Euler equations \eqref{E:OUTGOINGRIEMANNINVARIANTEVOLUTION}--\eqref{E:INGOINGRIEMANNINVARIANTEVOLUTION}
under any $C^4$ equation of state with positive sound speed except that of a Chaplygin gas.
Without loss of generality, assume the positivity condition\footnote{If $\lifespanconstant < 0$, 
then we can just change the sign of $\datasize$, and the theorem still holds. \label{E:CHANGESIGNOFDATAFOROPPOSITESIGNEDNONLINEARITIES}} 
\eqref{E:LIFESPANCONSTANT} on the constant $\lifespanconstant$.
Let $\mathbf{M}^{\textsf{Future}}$ be the subset of $(t,r)$-coordinate space in Fig.\,\ref{F:MGHD}
that includes $\Sigma_0$, lies above $\Sigma_0$ and lies strictly below $\futureCauchyhor \cup \Upsilon(\futuresinghyp)$,
where $\Upsilon(t,u) := (t,r)$ is the change of variables map from geometric coordinates to $(t,r)$ coordinates.
In particular,
$\mathbf{M}^{\textsf{Future}}$ is
the union of three sub-regions depicted in Fig.\ref{F:MGHD}: 
i) the ``flat'' region in between $\Sigma_0$ and $\Sigma_{\tstar}$.
ii) the Interior Region, which contains some points lying to the future of $\Sigma_{\tstar}$; 
and iii) the Exterior Region, which contains other points lying to the future of $\Sigma_{\tstar}$ and which is 
precisely the region where the eikonal function $u$ is defined.\footnote{Put differently, the Exterior Region in Fig.\,\ref{F:MGHD} is
the image of the change of variables map $\Upsilon$ from definition \eqref{E:CHOVFROMTUTOTRCOORDINATES}.}
The intersection of the Interior and Exterior Regions is equal to 
the portion of the characteristic $\lbrace u = \frac{\tstar}{2} \rbrace$ depicted in the figure.

Under the data assumptions stated in Section~\ref{S:DATA} (see Appendix~\ref{A:EXTENDRESULTSTOOTHERPROFILES} for a much larger class 
of data to which our results can be extended), 
if $\datasize$ is sufficiently small,
then the following results hold, where $\tstar = \frac{1}{\datasize}$.

\begin{enumerate}
	\item (\textbf{Existence and uniqueness of the future-solution}).
		The data on $\Sigma_0$ launch a unique classical solution $(\RRiemann,\LRiemann)$ 
		to equations \eqref{E:OUTGOINGRIEMANNINVARIANTEVOLUTION}--\eqref{E:INGOINGRIEMANNINVARIANTEVOLUTION} 
		on $\mathbf{M}^{\textsf{Future}}$ that
		satisfies the boundary conditions \eqref{E:RIEMANNINVARIANTMATCHINGCONDITION}--\eqref{E:RIEMANNINVARIANTRDERIVATIVEMATCHINGCONDITION}.
	\item (\textbf{Estimates}).
		The state of the solution at the ``early'' time $\tstar = \frac{1}{\datasize}$ is provided by
		Lemma~\ref{L:TIMETSTARDATAESTIMATESINTRCOORDINATES}.
		For later times, the solution satisfies the Exterior Region estimates\footnote{We 
        derived the Exterior Region estimates relative to geometric coordinates $(t,u)$, where we constructed the eikonal function $u$ in Sect.\,\ref{S:ACOUSTICGEOMETRY}. By ``satisfies the Exterior Region estimates'' in the theorem, we mean the following. First, given any point  $(t,r_*) \in \Upsilon(\mathbf{MGHD}^{\textnormal{Ext}})$, where 
        $\mathbf{MGHD}^{\textnormal{Ext}}$ is defined in \eqref{E:GLOBALLYHYPERBOLICEXTERIOREGIONPLUSABITOFCAUCHYHORIZON},
        there exists a unique point $(t,u_*) \in \mathbf{MGHD}^{\textnormal{Ext}}$ such that $\Upsilon(t,u_*) = (t,r_*)$; see Point~4 of Prop.\,\ref{P:EXISTENCEUPTOCREASEANDSINGULARBOUNDARYANDPORTIONOFCAUCHYHORIZON}.
        Then the following holds: $|\RRiemann(t,r_*)| \le C \frac{\datasize}{1+t+|u_*|}$, and similarly for the rest of the estimates in Prop.\,\ref{P:RIEMANNINVARIANTSAPRIORIEXTERIORREGIONESTIMATES}.} of Prop.\,\ref{P:RIEMANNINVARIANTSAPRIORIEXTERIORREGIONESTIMATES} and the Interior Region estimates featured in 
		Prop.\,\ref{P:APRIORIESTIMATESININTERIORREGION}.
	\item (\textbf{Structure of the future-boundary}).
		The future-boundary of $\mathbf{M}^{\textsf{Future}}$ is the union of two curves which, by definition,
		do not belong to $\mathbf{M}^{\textsf{Future}}$.
		The first curve is the future-singular boundary $\Upsilon(\futuresinghyp)$, which is closed,
		and which was constructed in Prop.\,\ref{P:EXISTENCEUPTOCREASEANDSINGULARBOUNDARYANDPORTIONOFCAUCHYHORIZON}. The past-boundary
		of $\Upsilon(\futuresinghyp)$ is a single point, the crease, denoted by $\Upsilon(\futurecrease)$ in Fig.\ref{F:MGHD}.
		The second curve is the Cauchy horizon $\futureCauchyhor$, an integral curve of
		$\uLunit$ that does not contain $\futurecrease$ but features it as a past-endpoint. 
		Moreover, $\futureCauchyhor$ 
		terminates to the future at a point $(T_{\textsf{Max}},0)$ on the axis of symmetry
		satisfying the following estimate, where $\Tcrease$ is the value of $t$ at the crease (see \eqref{E:BLOWUPTIMEUPPERANDLOWERBOUNDS}):
	\begin{align} \label{E:MAINTHEOREMESTIMATEFORTMAXININTERIORREGION}
	T_{\textsf{Max}}
	& \leq 2 \Tcrease + \mathcal{O}(1).
	\end{align}
\item (\textbf{Behavior of the solution on the future-boundary}).
	As was shown in Prop.\,\ref{P:EXISTENCEUPTOCREASEANDSINGULARBOUNDARYANDPORTIONOFCAUCHYHORIZON}, 
	the solution does not extend classically to $\Upsilon(\futuresinghyp)$ because $|\partial_r \RRiemann(q)| \to \infty$ 
	as $q$ approaches (from the Exterior Region) any point on
	$\Upsilon(\futuresinghyp)$.
	In contrast, the solution uniquely extends as a $C^2$ solution\footnote{$(\RRiemann,\LRiemann)$ could actually 
	be shown to be $C^3$, aside from the issue that there would be an $r$-weight in front of the third-order derivatives.
	In fact, in the Exterior Region estimates of Prop.\,\ref{P:RIEMANNINVARIANTSAPRIORIEXTERIORREGIONESTIMATES},
	we proved that the solution is $C^3$. In the Interior Region estimates of 
	Prop.\,\ref{P:APRIORIESTIMATESININTERIORREGION}, we only derived estimates for the up-to-second order derivatives of
	$(\RRiemann,\LRiemann)$; we did not try to control the third derivatives because we did not need such control to close the proof.} to 
	$\underline{\mathcal{C}}^{\textsf{max}}$, i.e., the solution exists classically up to the Cauchy horizon
	 (which does not contain its past-endpoint $\Upsilon(\futurecrease)$).
\item (\textbf{Global hyperbolicity}).
	$\mathbf{M}^{\textsf{Future}}$ is globally hyperbolic
	and contains the Cauchy hypersurface $\Sigma_0$. 
	That is, it satisfies the following:
	\begin{itemize}
		\item For every point $q \in \mathbf{M}^{\textsf{Future}}$, the past-directed integral curve of $\uLunit$
		emanating from $q$ remains in the interior of $\mathbf{M}^{\textsf{Future}}$ until it intersects $\Sigma_0$.
	\item For every point $q \in \mathbf{M}^{\textsf{Future}}$, the 
		past-directed integral curve of $\Lunit$
		emanating from $q$ remains in the interior of $\mathbf{M}^{\textsf{Future}}$ 
		until it intersects either the axis of symmetry $\lbrace r = 0 \rbrace$ or $\Sigma_0$.
	\end{itemize}
\item (\textbf{Similar results to the past}). 
For $t \leq 0$, the solution obeys analogous results, i.e., there exists a unique classical solution on a globally hyperbolic set
$\mathbf{M}^{\textsf{Past}}$ contained in $\lbrace t \leq 0 \rbrace$, the past-boundary of 
$\mathbf{M}^{\textsf{Past}}$ is the union of a closed, unbounded past-singular boundary $\pastsinghyp$ that emanates from the past-crease $\pastcrease$
and a pre-compact past-Cauchy horizon $\pastCauchyhor$ that emanates from $\pastcrease$ and terminates at the axis of symmetry.
Moreover, the solution obeys
estimates corresponding to
Props.\,\ref{P:RIEMANNINVARIANTSAPRIORIEXTERIORREGIONESTIMATES} 
and \ref{P:APRIORIESTIMATESININTERIORREGION}
with the following three changes: 
i) replacing $t$ with $-t = |t|$; ii) interchanging the roles of $\Lunit$ and $\uLunit$; iii) interchanging $\RRiemann$ and $\LRiemann$.
In particular, with $\pastsinghyp$ denoting the past singular boundary,
$|\partial_r \LRiemann|$ blows up as $\Upsilon(\pastsinghyp)$ 
is approached from any sequence of points lying in the past Interior Region.
Moreover $\mathbf{M}^{\textsf{Past}}$ is globally hyperbolic in the sense of Point~5, but with the roles of past and future interchanged
and with the roles of $\Lunit$ and $\uLunit$ interchanged.
\item (\textbf{Global hyperbolicity of} $\mathbf{M}^{\textsf{Max}}$ - \textbf{Part $A$}).
	Let $\mathbf{M}^{\textsf{Max}}$ denote the union
	$\mathbf{M}^{\textsf{Future}} \cup \mathbf{M}^{\textsf{Past}}$, viewed as a subset of $(t,r)$ coordinate space.
	Then $\mathbf{M}^{\textsf{Max}}$ is globally hyperbolic with $\Sigma_0$ as a Cauchy hypersurface.
	That is, it satisfies the following:
	\begin{itemize}
		\item For every point $q \in \mathbf{M}^{\textsf{Max}}$, the maximally extended 
			integral curve of $\Lunit$ passing through $q$ intersects either $\Sigma_0$ or $\lbrace r = 0 \rbrace$.
		\item For every point $q' \in \mathbf{M}^{\textsf{Max}}$, the maximally extended 
			integral curve of $\uLunit$ passing through $q$ intersects either $\Sigma_0$ or $\lbrace r = 0 \rbrace$.
	\end{itemize}
\item (\textbf{Global hyperbolicity of} $\mathbf{M}^{\textsf{Max}}$ - \textbf{Part $B$}). 
		If we now view $\mathbf{M}^{\textsf{Max}}$ to be
		a spherically symmetric subset\footnote{In Appendix~\ref{A:GLOBALHYPERBOLICITYANDMGHDS},
		the corresponding set in $\mathbb{R}^{1+3}$ is precisely described as $S(\mathbf{M}^{\textsf{Max}} \times \mathbb{S}^2)$, 
		where $S$ is the map defined in \eqref{E:MAPFROMSPHERICALCOORDINATESTOCARTESIANCOORDINATES}. \label{FN:MOREPRECISEDEFINITIONOFEMBEDDINGINTOR1PLUS3}} 
		of $\mathbb{R}^{1+3}$, 
		then
		with $\mathbf{g}$ denoting the acoustical metric on $\mathbb{R}^{1+3}$ defined in \eqref{E:MAINTHEOREMACOUSTICALMETRIC},
		$\mathbf{M}^{\textsf{Max}}$ is an open set, and
		the Lorentzian manifold $(\mathbf{M}^{\textsf{Max}},\mathbf{g})$
		is globally hyperbolic in the ``standard'' sense of Lorentzian geometry
		(see Def.\,\ref{D:CAUCHYHYPERSURFACESANDGLOBALLYHYPERBOLICREGIONS}), where $\lbrace t = 0 \rbrace$ is a Cauchy hypersurface.
\item (\textbf{Maximality of} $\mathbf{M}^{\textsf{Max}}$) 
Let $\mathbf{M}^{\textsf{Max}}$ be as in the previous point, in particular 
viewed as a spherically symmetric subset of $\mathbb{R}^{1+3}$.
Let $(\RRiemann,\LRiemann)$ be the corresponding classical solution on $\mathbf{M}^{\textsf{Max}}$
guaranteed by the above points.
	Assume that $\widetilde{\mathbf{M}}$ and $(\widetilde{\varrho},\widetilde{v}^1,\widetilde{v}^2,\widetilde{v}^3)$ 
		satisfy the following conditions:
		\begin{itemize}
			\item $\mathbf{M}^{\textsf{Max}} \subset \widetilde{\mathbf{M}} \subset \mathbb{R}^{1+3}$.
			\item $(\widetilde{\varrho},\widetilde{v}^1,\widetilde{v}^2,\widetilde{v}^3)$ is a $C^2$ solution to the compressible
				Euler equations \eqref{E:INTROTRANSPORTVI}--\eqref{E:INTROTRANSPORTDENSITY}
				on $\widetilde{\mathbf{M}}$ that agrees with the spherically symmetric solution
				on $\mathbf{M}^{\textsf{Max}}$ constructed above. Note that while on $\mathbf{M}^{\textsf{Max}}$,
				$(\widetilde{\varrho},\widetilde{v}^1,\widetilde{v}^2,\widetilde{v}^3)$ is spherically symmetric
				and has corresponding spherical Riemann invariants $(\widetilde{\mathcal{R}}_+,\widetilde{\mathcal{R}}_-)$
				satisfying $(\widetilde{\mathcal{R}}_+,\widetilde{\mathcal{R}}_-) \restriction_{\widetilde{\mathbf{M}}} = (\RRiemann,\LRiemann)$,
				$(\widetilde{\varrho},\widetilde{v}^1,\widetilde{v}^2,\widetilde{v}^3)$ is not assumed to be 
				spherically symmetric on $\widetilde{\mathbf{M}} \backslash \mathbf{M}^{\textsf{Max}}$,
				nor is the region $\widetilde{\mathbf{M}} \backslash \mathbf{M}^{\textsf{Max}}$ assumed to be spherically symmetric.
			\item Let $\widetilde{\mathbf{g}}$ denote the acoustical metric \eqref{E:MAINTHEOREMACOUSTICALMETRIC}, evaluated at
				the solution $(\widetilde{\varrho},\widetilde{v}^1,\widetilde{v}^2,\widetilde{v}^3)$.
			We assume that $\widetilde{\mathbf{M}}$ is open 
			and that the Lorentzian manifold $(\widetilde{\mathbf{M}},\widetilde{\mathbf{g}})$ 
			is globally hyperbolic in the sense of Def.\,\ref{D:CAUCHYHYPERSURFACESANDGLOBALLYHYPERBOLICREGIONS}
			and contains $\lbrace t = 0 \rbrace$ as a Cauchy hypersurface.
		\end{itemize} Then $\widetilde{\mathbf{M}} = \mathbf{M}^{\textsf{Max}}$. That is, there is no globally hyperbolic 
		development of the data that strictly contains $\mathbf{M}^{\textsf{Max}}$.
\item (\textbf{Uniqueness of the MGHD}) 		
	Let $(\RRiemann,\LRiemann)$ be the classical solution on $\mathbf{M}^{\textsf{Max}}$
	guaranteed by the above points.
	Assume that $\widehat{\mathbf{M}}$
	and $(\widehat{\varrho},\widehat{v}^1,\widehat{v}^2,\widehat{v}^3)$ satisfy the following conditions:
		\begin{itemize}
			\item $\lbrace t = 0 \rbrace \subset \widehat{\mathbf{M}} \subset \mathbb{R}^{1+3}$. 
			\item $(\widehat{\varrho},\widehat{v}^1,\widehat{v}^2,\widehat{v}^3)$ is a $C^2$ 
				solution to the compressible Euler equations \eqref{E:INTROTRANSPORTVI}--\eqref{E:INTROTRANSPORTDENSITY}
				on $\widehat{\mathbf{M}}$ that has the same initial data\footnote{By ``same initial data,'' we mean that 
				the data for $(\widehat{\varrho},\widehat{v}^1,\widehat{v}^2,\widehat{v}^3)$ are the unique spherically symmetric data
				corresponding to $(\RRiemann,\LRiemann) \restriction_{\lbrace t = 0 \rbrace}$.} 
				on $\lbrace t = 0 \rbrace$ as $(\RRiemann,\LRiemann)$.
			\item Let $\widehat{\mathbf{g}}$ denote the acoustical metric \eqref{E:MAINTHEOREMACOUSTICALMETRIC}, evaluated at
				the solution $(\widehat{\varrho},\widehat{v}^1,\widehat{v}^2,\widehat{v}^3)$.
			We assume that $\widehat{\mathbf{M}}$ is open, and the Lorentzian manifold
			$(\widehat{\mathbf{M}},\widehat{\mathbf{g}})$ is globally hyperbolic in the sense of Def.\,\ref{D:CAUCHYHYPERSURFACESANDGLOBALLYHYPERBOLICREGIONS}, 
			and contains $\lbrace t = 0 \rbrace$ as a Cauchy hypersurface.
			Note that we are not assuming that $(\widehat{\varrho},\widehat{v}^1,\widehat{v}^2,\widehat{v}^3)$ or
			$\widehat{\mathbf{M}}$ are spherically symmetric.
		\end{itemize}
		Then $\widehat{\mathbf{M}} \subset \mathbf{M}^{\textsf{Max}}$, 
		$(\widehat{\varrho},\widehat{v}^1,\widehat{v}^2,\widehat{v}^3)$ is spherically symmetric on 
		$\widehat{\mathbf{M}}$ with corresponding spherical Riemann invariants
		$(\widehat{\mathcal{R}}_+,\widehat{\mathcal{R}}_-)$,
		and
		$(\RRiemann,\LRiemann)\restriction_{\widehat{\mathbf{M}}} = (\widehat{\mathcal{R}}_+,\widehat{\mathcal{R}}_-)$.
		That is, $(\mathbf{M}^{\textsf{Max}},(\RRiemann,\LRiemann))$
		is an extension of \underline{all} other globally hyperbolic developments 
		(see Def.\,\ref{D:GHD})
		of the data on $\lbrace t = 0 \rbrace$
		and thus is the unique MGHD of the data.
	\end{enumerate}

\end{theorem}

\begin{remark}[Carefully distinguish maximality from uniqueness]
\label{R:MAXIMAILITYVSUNIQUENESS}
Point~$9$ says that $(\mathbf{M}^{\textsf{Max}},(\RRiemann,\LRiemann))$ is \underline{a} largest possible globally hyperbolic classical solution, i.e., 
maximality. A subtlety is that Point~$9$ by itself does not preclude the existence of a distinct maximal GHD $\widetilde{\mathbf{M}}^{\textnormal{Max}} \neq \mathbf{M}^{\textnormal{Max}}$ satisfying the same assumptions and conclusions of Point~9. However,
Point~$10$ rules out the possibility of a distinct $\widetilde{\mathbf{M}}^{\textnormal{Max}}$ because
$(\mathbf{M}^{\textsf{Max}},(\RRiemann,\LRiemann))$ \emph{contains} \underline{all} other globally hyperbolic classical solutions;
this ``containment'' is the key ingredient behind MGHD uniqueness. 
In contrast, it is known that for a different hyperbolic PDE, there exist $C^{\infty}$ data such 
that no solution exists with properties as in Point~$10$, even though multiple distinct solutions 
with properties as in Point~$9$ do exist; see Appendix~\ref{A:GLOBALHYPERBOLICITYANDMGHDS}.
\end{remark}

\subsection{State of the solution in the Interior Region}
\label{SS:SOLUTIONININTERIOR}
To control the solution in the Interior Region, we will 
use some of the ``data'' induced by the Exterior Region solution 
(yielded by Prop.\,\ref{P:EXISTENCEUPTOCREASEANDSINGULARBOUNDARYANDPORTIONOFCAUCHYHORIZON})
on the outgoing characteristic portion
$\lbrace u = \frac{\tstar}{2} \rbrace \cap \lbrace \tstar \leq t \leq \TCH\left(\frac{\tstar}{2} \right) \rbrace$.
In the next lemma, we derive the relevant estimates for that ``data.''
All the estimates in the lemma are easy consequences of results we have already derived and their proofs.

\begin{lemma}[State of the solution on $\lbrace u = \frac{\tstar}{2} \rbrace \cap \lbrace \tstar \leq t \leq \TCH\left(\frac{\tstar}{2} \right) \rbrace$]
\label{L:RIEMANNINVARIANTSDATAFORINTERIORPROBLEM}
Under assumptions and conclusions of Prop.\,\ref{P:EXISTENCEUPTOCREASEANDSINGULARBOUNDARYANDPORTIONOFCAUCHYHORIZON},
if $\datasize$ is sufficiently small,
then the following results hold, where $\TCH = \TCH(u)$ is the function from the proposition.

\medskip

\noindent \underline{\textbf{Behavior of $t$ and $r$}}.
The following comparison estimates hold on $\lbrace u = \frac{\tstar}{2} \rbrace \cap \lbrace \tstar \leq t \leq \TCH\left(\frac{\tstar}{2} \right) \rbrace$:
\begin{align} \label{E:COMPARISONBETWEENUANDMINUSRONCREASECHARACTERISTIC}
	u 
	&
	=
	\frac{\tstar}{2} 
	= t - r 
	+
	\mathcal{O}(\datasize)
	\ln \left( \frac{1 + t + \frac{\tstar}{2}}{1 + \frac{3\tstar}{2} } \right)
	= t - r 
	+
	\mathcal{O}(1).
\end{align}

\medskip

\noindent \underline{\textbf{Estimates for the acoustic geometry}}.
The following pointwise estimates hold on $\lbrace u = \frac{\tstar}{2} \rbrace \cap \lbrace \tstar \leq t \leq \TCH\left(\frac{\tstar}{2} \right) \rbrace$:
\begin{subequations}
\begin{align} \label{E:MUPOINTWISEESTIMATECREASECHARACTERISTIC}
	\upmu \restriction_{\lbrace u = \frac{\tstar}{2} \rbrace \cap \lbrace \tstar \leq t \leq \TCH\left(\frac{\tstar}{2} \right) \rbrace}
	& = 1 + \mathcal{O}(\datasize^2),
		\\
	\muX \upmu \restriction_{\lbrace u = \frac{\tstar}{2} \rbrace \cap \lbrace \tstar \leq t \leq \TCH\left(\frac{\tstar}{2} \right) \rbrace}
	& = \mathcal{O}(\datasize^3),
		\label{E:MUXMUPOINTWISEESTIMATECREASECHARACTERISTIC} 
		\\
	\muX \muX \upmu \restriction_{\lbrace u = \frac{\tstar}{2} \rbrace \cap \lbrace \tstar \leq t \leq \TCH\left(\frac{\tstar}{2} \right) \rbrace}
	& = \mathcal{O}(\datasize^4).
	\label{E:MUXTWICEMUPOINTWISEESTIMATECREASECHARACTERISTIC}
\end{align}
\end{subequations}

\medskip 
\noindent \underline{\textbf{Estimates for the Riemann invariants}}.
The following pointwise estimates hold on $\lbrace u = \frac{\tstar}{2} \rbrace \cap \lbrace \tstar \leq t \leq \TCH\left(\frac{\tstar}{2} \right) \rbrace$.

For $0 \leq J + K \leq 3$, we have:
\begin{align} \label{E:INTERIORREGIONCREASECHARACTERISTICPOINTWISEESTIMATEFORRPLUSANDDERIVATIVES}
	|\Lunit^J \uLunit^K  \RRiemann \restriction_{\lbrace u = \frac{\tstar}{2} \rbrace \cap \lbrace \tstar \leq t \leq \TCH\left(\frac{\tstar}{2} \right) \rbrace}|
	& 
	\leq C \frac{\datasize}{(1 + t + r)^{1+J}(1 + |t-r|)^K}.
\end{align}

For $0 \leq J \leq 3$, we have:
\begin{align} \label{E:INTERIORREGIONCREASECHARACTERISTICPOINTWISEESTIMATEFORRMINUSANDLDERIVATIVES}
	|\Lunit^J \LRiemann \restriction_{\lbrace u = \frac{\tstar}{2} \rbrace \cap \lbrace \tstar \leq t \leq \TCH\left(\frac{\tstar}{2} \right) \rbrace}|
	& 
	\leq C \frac{\datasize}{(1 + t + r)^{1+J}}.
\end{align}

For $0 \leq J + K' \leq 2$, we have:
\begin{align} \label{E:INTERIORREGIONCREASECHARACTERISTICPOINTWISEESTIMATEFORRMINUSATLEASTONEXDERIVATIVE}
	|\Lunit^J \uLunit^{K'} \uLunit \LRiemann 
	\restriction_{\lbrace u = \frac{\tstar}{2} \rbrace \cap \lbrace \tstar \leq t \leq \TCH\left(\frac{\tstar}{2} \right) \rbrace}|
	& 
	\leq C \frac{\datasize}{(1 + t + r)^{2+J}(1 + |t-r|)^{K'}}.
\end{align}

\medskip 
\noindent \underline{\textbf{Estimates for special combinations involving the Riemann invariants}}.
The following pointwise estimates hold on $\lbrace u = \frac{\tstar}{2} \rbrace \cap \lbrace \tstar \leq t \leq \TCH\left(\frac{\tstar}{2} \right) \rbrace$.

For $0 \leq J + K \leq 3$, we have:
\begin{align} \label{E:INTERIORREGIONCREASECHARACTERISTICPOINTWISEESTIMATEFORCOMMUTEDVERSIONLDERIVATIMEOFMODIFIEDLBARRPLUSDERIVATIVE}
		\left|
		\Lunit 
		\Lunit^J
		\muX^K
		\left\lbrace
			r
			\muuLunit \RRiemann
			- 
			2 \upmu \RRiemann
		\right\rbrace
		\restriction_{\lbrace u = \frac{\tstar}{2} \rbrace \cap \lbrace \tstar \leq t \leq \TCH\left(\frac{\tstar}{2} \right) \rbrace}
		\right|
		& \lesssim 
		\datasize^{1^+} \frac{1}{(1 + t + r)^{2+J}(1 + |t-r|)^{1+K}}.
	\end{align}
	
	For $0 \leq K \leq 3$, we have:
	\begin{align} \label{E:CREASECHARACTERISTICPOINTWISEBOUNDALLMUXDERIVATIVESTRANSPORTEDMODIFIEDULUNITRPLUSTSTARDATA}
	\begin{split}	
		\muX^K
		\left\lbrace
			r
			\muuLunit \RRiemann
			- 
			2 \upmu \RRiemann
		\right\rbrace
		\restriction_{\lbrace u = \frac{\tstar}{2} \rbrace \cap \lbrace \tstar \leq t \leq \TCH\left(\frac{\tstar}{2} \right) \rbrace}
		& 
		=
				2 \datasize 
				\left( \frac{\partial}{\partial u} \right)^K
				\frac{1}{1 + u^2} \restriction_{\lbrace u = \frac{\tstar}{2} \rbrace \cap \lbrace \tstar \leq t \leq \TCH\left(\frac{\tstar}{2} \right) \rbrace}
					\\
	& \ \
		+
			\frac{\mathcal{O}(\datasize^{1^+})}{(1 + \tstar + r)^J(1 + |t-r|)^{2+K}}.
	\end{split}
	\end{align}

	For $0 \leq J \leq 2$, we have:
	\begin{align} \label{E:CREASECHARACTERISTICPOINTWISEESTIMATEFORALLLUNITDERIVATIVESMODIFIEDLRMINUSDERIVATIVE}
	\begin{split}
		&
		\left\lbrace
		\Lunit^J
		[r \Lunit \LRiemann
		+
		2 \LRiemann]
		\right\rbrace
		\restriction_{\lbrace u = \frac{\tstar}{2} \rbrace \cap \lbrace \tstar \leq t \leq \TCH\left(\frac{\tstar}{2} \right) \rbrace}
			\\
		& =
		\left\lbrace
		-
		2 \datasize 
		\left( \frac{\partial}{\partial t} \right)^J
		\frac{1}{1 + (2t - u)^2}
		+
		\mathcal{O}(\datasize^{1^+}) \frac{\ln_+(t + |u|)}{(1 + t + |u|)^{2+J}}
		\right\rbrace \restriction_{u = \frac{\tstar}{2}}.
	\end{split}
	\end{align}
	
	Finally, we have:
	\begin{align} \label{E:CREASECHARACTERISTICPOINTWISEESTIMATEFORPARTIALTMODIFIEDLRMINUSDERIVATIVE}
	\begin{split}
		\left\lbrace
		[r \Lunit \partial_t \LRiemann
		+
		2 \Speed \partial_t \LRiemann]
		\right\rbrace\restriction_{\lbrace u = \frac{\tstar}{2} \rbrace \cap \lbrace \tstar \leq t \leq \TCH\left(\frac{\tstar}{2} \right) \rbrace}
		& =
		\frac{\mathcal{O}(\datasize)}{(1 + t + r)^3}
		+
		\frac{\mathcal{O}(\datasize^2) \ln(t+r)}{(1 + t + r)^3}.
	\end{split}
	\end{align}
	
\end{lemma}

\begin{proof}
The estimates of Prop.\,\ref{P:RIEMANNINVARIANTSAPRIORIEXTERIORREGIONESTIMATES}
hold on the full exterior region
$\mathbf{MGHD}^{\textsf{Ext}}$
from Prop.\,\ref{P:EXISTENCEUPTOCREASEANDSINGULARBOUNDARYANDPORTIONOFCAUCHYHORIZON}.
Similarly, by Remark~\ref{R:STRICTIMPROVEMENTOFBOOTSTRAP},
on $\mathbf{MGHD}^{\textsf{Ext}}$, the estimates of Prop.\,\ref{P:SHARPESTIMATESFORMU} hold with $\eps$ replaced by 
$C \datasize$. Throughout the proof, we will use these facts.
We will also silently use
\eqref{E:SHARPPOINWISECOMPARISONBETWEENUANDTMINUSR}
and
\eqref{E:RISAPPROXIMATELYTMINUSU}--\eqref{E:COMPARISONBETWEENTOVERUTMINUSUANDTPLUSMODU}.

\medskip

\noindent \underline{\textbf{Proof of \eqref{E:COMPARISONBETWEENUANDMINUSRONCREASECHARACTERISTIC}}}.
This estimate follows from \eqref{E:SHARPPOINWISECOMPARISONBETWEENUANDTMINUSR},
\eqref{E:UCREASEISSMALL}--\eqref{E:BLOWUPTIMEUPPERANDLOWERBOUNDS},
and
\eqref{E:SHARPUPPERANDLOWERBOUNDSONVALUEOFTATTOPPOINTOFUEQUALSTSTAROVER2}.

\medskip

\noindent \underline{\textbf{Proof of \eqref{E:MUPOINTWISEESTIMATECREASECHARACTERISTIC}--\eqref{E:MUXTWICEMUPOINTWISEESTIMATECREASECHARACTERISTIC}}}.
These estimates follow from \eqref{E:MUPOINTWISEESTIMATEEXTERIOR}--\eqref{E:TWOXBREVEMUPOINTWISEESTIMATEEXTERIOR}
and the upper bound on $\TCH\left(\frac{\tstar}{2} \right)$ implied by 
\eqref{E:UCREASEISSMALL}--\eqref{E:BLOWUPTIMEUPPERANDLOWERBOUNDS}
and
\eqref{E:SHARPUPPERANDLOWERBOUNDSONVALUEOFTATTOPPOINTOFUEQUALSTSTAROVER2}.

\medskip

\noindent \underline{\textbf{Proof of \eqref{E:INTERIORREGIONCREASECHARACTERISTICPOINTWISEESTIMATEFORRPLUSANDDERIVATIVES}--\eqref{E:INTERIORREGIONCREASECHARACTERISTICPOINTWISEESTIMATEFORRMINUSATLEASTONEXDERIVATIVE}}}.
These estimates were proved in Prop.\,\ref{P:RIEMANNINVARIANTSAPRIORIEXTERIORREGIONESTIMATES}.

\medskip

\noindent \underline{\textbf{Proof of \eqref{E:INTERIORREGIONCREASECHARACTERISTICPOINTWISEESTIMATEFORCOMMUTEDVERSIONLDERIVATIMEOFMODIFIEDLBARRPLUSDERIVATIVE}--\eqref{E:CREASECHARACTERISTICPOINTWISEBOUNDALLMUXDERIVATIVESTRANSPORTEDMODIFIEDULUNITRPLUSTSTARDATA}}}. 
These estimates follow from \eqref{E:EXTERIORREGIONPOINTWISEESTIMATEFORCOMMUTEDVERSIONLDERIVATIMEOFMODIFIEDLBARRPLUSDERIVATIVE}--\eqref{E:POINTWISEBOUNDALLMUXDERIVATIVESTRANSPORTEDMODIFIEDULUNITRPLUSTSTARDATA} and \eqref{E:SHARPPOINWISECOMPARISONBETWEENUANDTMINUSR}.

\medskip 
\noindent \underline{\textbf{Proof of \eqref{E:CREASECHARACTERISTICPOINTWISEESTIMATEFORALLLUNITDERIVATIVESMODIFIEDLRMINUSDERIVATIVE}}}.
This estimate was proved as \eqref{E:EXTERIORREGIONPOINTWISEESTIMATEFORALLLUNITDERIVATIVESMODIFIEDLRMINUSDERIVATIVE}.

\medskip 
\noindent \underline{\textbf{Proof of \eqref{E:CREASECHARACTERISTICPOINTWISEESTIMATEFORPARTIALTMODIFIEDLRMINUSDERIVATIVE}}}.
First, using \eqref{E:PARTIALTINTERMSOFLANDLBAR}, \eqref{E:MUPOINTWISEESTIMATECREASECHARACTERISTIC}, 
and Prop.\,\ref{P:RIEMANNINVARIANTSAPRIORIEXTERIORREGIONESTIMATES},
we deduce that on $\lbrace u = \frac{\tstar}{2} \rbrace \cap \lbrace \tstar \leq t \leq \TCH\left(\frac{\tstar}{2} \right) \rbrace$,
we have:
\begin{align} \label{E:PARTIALTINTERMSOFLANDMULBARALONGCREASECHARACTERISTICS}
\partial_t
& = 	\mathcal{O}(1)
			\Lunit
			+
			\mathcal{O}(1)
			\muuLunit.
\end{align}	
Next, using \eqref{E:PARTIALTINTERMSOFLANDMULBARALONGCREASECHARACTERISTICS}, 
the estimates of Prop.\,\ref{P:RIEMANNINVARIANTSAPRIORIEXTERIORREGIONESTIMATES}, 
and \eqref{E:SPEEDOFSOUNDEXPANSION}--\eqref{E:SPEEDOFSOUNDERRORFUNCTIONVANISHESATORIGIN},
we further deduce that
on $\lbrace u = \frac{\tstar}{2} \rbrace \cap \lbrace \tstar \leq t \leq \TCH\left(\frac{\tstar}{2} \right) \rbrace$, we have:
\begin{align} \label{E:ZEROPROOFSTEPCREASECHARACTERISTICPOINTWISEESTIMATEFORPARTIALTMODIFIEDLRMINUSDERIVATIVE}
	\mbox{LHS~\eqref{E:CREASECHARACTERISTICPOINTWISEESTIMATEFORPARTIALTMODIFIEDLRMINUSDERIVATIVE}}
	& = 
		\partial_t
		[r \Lunit \LRiemann
		+
		2 \LRiemann]
		+ 
		\frac{\mathcal{O}(\datasize^2)}{(1 + t + |u|)^2(1 + |u|)}.
\end{align}
Using
\eqref{E:PARTIALTINTERMSOFLANDMULBARALONGCREASECHARACTERISTICS},
\eqref{E:WAVEEQUATIONFORALLDERIVATIVESRMINUSINHOMOGENEOUSTERMBOUND},
and 
\eqref{E:EXTERIORREGIONPOINTWISEESTIMATEFORALLLUNITDERIVATIVESMODIFIEDLRMINUSDERIVATIVE}
to bound the first term on RHS~\eqref{E:ZEROPROOFSTEPCREASECHARACTERISTICPOINTWISEESTIMATEFORPARTIALTMODIFIEDLRMINUSDERIVATIVE},
we further deduce that:
\begin{align} \label{E:FIRSTPROOFSTEPCREASECHARACTERISTICPOINTWISEESTIMATEFORPARTIALTMODIFIEDLRMINUSDERIVATIVE}
	\mbox{LHS~\eqref{E:CREASECHARACTERISTICPOINTWISEESTIMATEFORPARTIALTMODIFIEDLRMINUSDERIVATIVE}}
	& = 
		\mathcal{O}(\datasize)
		\frac{1}{[1 + (2t-u)]^3}
		+ 
		\frac{\mathcal{O}(\datasize^2)\ln(t + |u|)}{(1 + t + |u|)^3}.
\end{align}
From \eqref{E:FIRSTPROOFSTEPCREASECHARACTERISTICPOINTWISEESTIMATEFORPARTIALTMODIFIEDLRMINUSDERIVATIVE}
and \eqref{E:COMPARISONBETWEENUANDMINUSRONCREASECHARACTERISTIC},
we conclude \eqref{E:CREASECHARACTERISTICPOINTWISEESTIMATEFORPARTIALTMODIFIEDLRMINUSDERIVATIVE}.
\end{proof}

\subsection{Bootstrap assumptions on a globally hyperbolic Interior Region}
\label{SS:BOOTSTRAPASSUMPTIONSONINTERIORREGION}

\subsubsection{Globally hyperbolic interior bootstrap region}
\label{SSS:GLOBALLYHYPERBOLICINTERIOR}
We will carry out our Interior Region analysis on bootstrap regions
of the type featured in the following definition.

\begin{definition}[Globally hyperbolic interior bootstrap region]
	\label{D:GHINTERIORBOOTSTRAP}
	A subset $\mathbf{GH}_{\textnormal{Boot}}^{{\textnormal{INT}}}$ of spacetime is said to be\footnote{Remark~\ref{R:CONNECTIONTOSTANDARDGLOBALHYPERBOLICITY} applies here as well.} 
	a \textbf{globally hyperbolic interior bootstrap region} 
	if the following conditions hold: 
	\begin{itemize}
	\item $(\RRiemann,\LRiemann)$ is a $C^2$ solution to \eqref{E:OUTGOINGRIEMANNINVARIANTEVOLUTION}--\eqref{E:INGOINGRIEMANNINVARIANTEVOLUTION} 
		on $\mathbf{GH}_{\textnormal{Boot}}^{{\textnormal{INT}}}$.
	\item 
		\begin{align} \label{E:BOUNDSTHATDEFINEINTERIORGLOBALLYHYPEROLICBOOTSTRAPREGION}
			\mathbf{GH}_{\textnormal{Boot}}^{{\textnormal{INT}}} \subset \lbrace \tstar \leq t \leq 4 \Tcrease \rbrace,
		\end{align}
	where $\Tcrease$ is the value of $t$ at the crease (see \eqref{E:BLOWUPTIMEUPPERANDLOWERBOUNDS}).
		In particular, by assumption, the following holds everywhere in $\mathbf{GH}_{\textnormal{Boot}}^{{\textnormal{INT}}}$:
	\begin{align} \label{E:INTERIORTIMEBOUNDASSUMPTION}
		\tstar & \leq t \leq 4 \Tcrease.
	\end{align}
	\item For every point $q \in \mathbf{GH}_{\textnormal{Boot}}^{{\textnormal{INT}}}$, the past-directed integral curve of $\uLunit$
	emanating from $q$ remains in the interior of $\mathbf{GH}_{\textnormal{Boot}}^{{\textnormal{INT}}}$ until it intersects
	a point $q'$ satisfying
	$q' \in \mathbf{GH}_{\textnormal{Boot}}^{{\textnormal{INT}}} 
	\bigcap
	\Sigma^{\textnormal{INT;Data}}
	$,
	where
	$
	\Sigma^{\textnormal{INT;Data}}
	:=
	\left(\Sigma_{\tstar} \cap \lbrace r \leq \tstar/2 \rbrace \right) 
	\bigcup 
	\left( \lbrace u = \frac{\tstar}{2} \rbrace
	\cap 
	\lbrace \tstar \leq t < 4 \Tcrease \rbrace
	\right)
	$
	is the union of the two ``data hypersurfaces'' for the Interior Region, i.e., a spacelike portion joined to a characteristic portion.
	\item For every point $q \in \mathbf{GH}_{\textnormal{Boot}}^{{\textnormal{INT}}}$, the 
	past-directed integral curve of $\Lunit$
	emanating from $q$ remains in the interior of $\mathbf{GH}_{\textnormal{Boot}}^{{\textnormal{INT}}}$
	until it intersects a point $q'$
	satisfying 
	$$q' 
	\in 
	\mathbf{GH}_{\textnormal{Boot}}^{{\textnormal{Ext}}} 
	\bigcap 
	\left(
	\left(\Sigma_{\tstar} \cap \lbrace r \leq \tstar/2 \rbrace \right) 
	\bigcup \lbrace r = 0 \rbrace
	\right).$$
 \end{itemize}	
\end{definition}


\subsubsection{Interior bootstrap assumptions}
\label{SSS:INTERIORBOOTSTRAPASSUMPTIONS}
As in our study of the Exterior Region, to control the solution in the Interior Region, 
we find it convenient to derive estimates via a bootstrap argument. 
To this end, we fix a globally hyperbolic interior bootstrap region $\mathbf{GH}_{\textnormal{Boot}}^{{\textnormal{INT}}}$,  
and we assume that the following inequalities hold for some number $\eps$ satisfying:
\begin{align} \label{E:EPSILONSMALLNESSASSUMPTIONININTERIOREGION}
	0 & < \eps \leq \datasize^{3/4}.
\end{align}

\noindent \underline{\textbf{Bootstrap assumptions for the Riemann invariants}}.

\begin{subequations}
\begin{align} \label{E:POINTWISEBAINTERIORRPLUS}
	|\RRiemann|
	& 
	\leq \frac{\eps}{1 + t + r},
		\\
|\partial_t \RRiemann|
		& 
	\leq \frac{\eps}{(1 + t + r)(1 + |t-r|)},
		\label{E:POINTWISEBAINTERIORPARTIALTRPLUS}
			\\
	|\uLunit \RRiemann|
	& 
	\leq \frac{\eps}{(1 + t + r)(1 + |t-r|)},
		\label{E:POINTWISEBAINTERIORLBARRPLUS}
		\\
|r \uLunit \RRiemann|
	& 
	\leq
	\frac{\eps}{1 + t + r}
	+
	\frac{\eps \ln(t-r)}{(1 + |t - r|)^2},
		\label{E:POINTWISEBAINTERIORRTIMESLBARRPLUS}
			\\
|\Lunit \RRiemann|
	& 
	\leq 
		\frac{\eps}{(1 + t + r)^2},
		\label{E:POINTWISEBAINTERIORLRPLUS}
			\\
|r \Lunit \RRiemann|
	& 
	\leq 
		\frac{\eps}{1 + t + r},
		\label{E:POINTWISEBAINTERIORRTIMESLRPLUS}
			\\
|r \uLunit \partial_t \RRiemann|
		& 
	\leq  
		\frac{\eps}{(1 + t + r)(1 + |t-r|)}
		+
		\frac{\eps}{(1 + |t-r|)^3}
		+
		\frac{\eps^{1+}}{(1 + |t-r|)^2}
		+
		\frac{\eps^{1+} \ln(t-r)}{(1 + |t-r|)^3},
		\label{E:POINTWISEBAINTERIORRTIMESLBARPARTIALTRPLUS}
			\\
|r \Lunit \uLunit \RRiemann|
	& 
	\leq 
	\frac{\eps}{(1 + t + r)(1 + |t-r|)},
	\label{E:POINTWISEBAINTERIORTIMESLUNITLBARRPLUS}
		\\
|r \uLunit \uLunit \RRiemann|
	& 
	\leq 
		\frac{\eps}{(1 + t + r)(1 + |t-r|)}
		+
		\frac{\eps}{(1 + |t-r|)^3}
		+
		\frac{\eps^{1+}}{(1 + |t-r|)^2}
		+
		\frac{\eps^{1+} \ln(t-r)}{(1 + |t-r|)^3},
	\label{E:POINTWISEBAINTERIORTIMESLBARLBARRPLUS}
	\\
|r \Lunit \Lunit \RRiemann|
	& 
	\leq 
	\frac{\eps}{(1 + t + r)^2},
	\label{E:POINTWISEBAINTERIORTIMESLUNITLUNITRPLUS}
\end{align}
\end{subequations}

\begin{subequations}
\begin{align} \label{E:POINTWISEBAINTERIORRMINUS}
	|\LRiemann|
	& 
	\leq \frac{\eps}{1 + t + r},
		\\
|\partial_t \LRiemann|
		& 
	\leq \frac{\eps}{(1 + t + r)^2},
		\label{E:POINTWISEBAINTERIORPARTIALTRMINUS}
			\\
|\uLunit \LRiemann|
		& 
	\leq \frac{\eps}{(1 + t + r)^2},
		\label{E:POINTWISEBAINTERIORLBARRMINUS}
		\\
	|r \uLunit \LRiemann|
	& 
	\leq\frac{\eps}{1 + t + r},
		\label{E:POINTWISEBAINTERIORRTIMESLBARRMINUS}
		\\
|\Lunit \LRiemann|
	&  
	\leq 
	\frac{\eps}{(1 + t + r)^2},
		\label{E:POINTWISEBAINTERIORLRMINUS}
			\\
|r \Lunit \LRiemann|
	& 
	\leq 
	\frac{\eps}{1 + t + r},
		\label{E:POINTWISEBAINTERIORRTIMESLRMINUS}
			\\
	|r \Lunit \Lunit \LRiemann|
	& 
	\leq 
	\frac{\eps}{(1 + t + r)^2},
	\label{E:POINTWISEBAINTERIORRTIMESLUNITLUNITRMINUS}
		\\
	|r \Lunit \partial_t \LRiemann|
		& 
	\leq 
	\frac{\eps}{(1 + t + r)^2},
		\label{E:POINTWISEBAINTERIORRTIMESLUNITPARTIALTRMINUS}
			\\
	|r \uLunit \Lunit \LRiemann|
	& 
	\leq 
	\frac{\eps}{(1 + t + r)^2},
	\label{E:POINTWISEBAINTERIORRTIMESLBARLUNITRMINUS}
		\\
	|r \uLunit \uLunit \LRiemann|
		& 
	\leq 
	\frac{\eps}{(1 + t + r)(1 + |t-r|)}.
		\label{E:POINTWISEBAINTERIORRTIMESLBARLBARRMINUS}
\end{align}
\end{subequations}

\noindent \underline{\textbf{Bootstrap assumption for the coordinates}}.
We assume that:
\begin{align} \label{E:INTERIORREGIONTMINUSRLARGEBOOTSTRAPASSUMPTION}
	\frac{\tstar}{4}
	& \leq
	t - r. 
\end{align}

Since $t \geq \tstar$, we deduce from \eqref{E:INTERIORREGIONTMINUSRLARGEBOOTSTRAPASSUMPTION} that:
\begin{align} \label{E:TPLUSRANDTEQUIVALENTININTERIORREGION}
	t+r
	\leq
	2t - \frac{\tstar}{4}
	\leq 
	\frac{7}{4} t \leq 2 (t+r).
\end{align}	

We will silently use \eqref{E:INTERIORREGIONTMINUSRLARGEBOOTSTRAPASSUMPTION}--\eqref{E:TPLUSRANDTEQUIVALENTININTERIORREGION} 
throughout our analysis of the Interior Region.

\medskip

\noindent \underline{\textbf{Bootstrap assumption for $t-r$ along the integral curves of $\Lunit$}}.
Let $(t_0,r_0) \in \mathbf{GH}_{\textnormal{Boot}}^{{\textnormal{INT}}}$, and let 
$(t_1,r_1) \in \mathbf{GH}_{\textnormal{Boot}}^{{\textnormal{INT}}}$ be any point on the future-directed integral curve of $\Lunit$ emanating from
$(t_0,r_0)$. We assume that:
\begin{align} \label{E:BAINTERIORTMINUSRCHANGEALONGINTEGRALCURVESOFLUNIT}
	\left| (t_1 - r_1) - (t_0 - r_0) \right|
	& \leq \ln \left(\frac{1}{\datasize} \right).
\end{align}

In particular, by \eqref{E:INTERIORREGIONTMINUSRLARGEBOOTSTRAPASSUMPTION},
we have:
\begin{align}
	|t_0 - r_0| & = (1 + \mathcal{O}(\datasize^{1/2})) |t_1-r_1|.
\end{align}

\subsection{Preliminary estimates}
\label{SS:PRELIMINARYESTIMATES}
We start by deriving some simple estimates for various time values in $\mathbf{GH}_{\textnormal{Boot}}^{{\textnormal{INT}}}$.

\begin{lemma}[Some simple estimates for time values in $\mathbf{GH}_{\textnormal{Boot}}^{{\textnormal{INT}}}$]
\label{E:SIMPLEESITMATESFORTIMEVALUES}
Under the data-assumptions stated in Sect.\,\ref{S:DATA}
and the bootstrap assumptions of Sect.\,\ref{SSS:INTERIORBOOTSTRAPASSUMPTIONS},
if $\datasize$ is sufficiently small, then the following results hold.

\medskip

\noindent \underline{\textbf{Bounds for} $\Tcrease$}.
The time $\Tcrease$ from Prop.\,\ref{P:EXISTENCEUPTOCREASEANDSINGULARBOUNDARYANDPORTIONOFCAUCHYHORIZON}
satisfies the following upper and lower bounds,
where $\lifespanconstant > 0$ is the constant from \eqref{E:NULLCONDITIONFAILURECONSTANT}:
\begin{align} \label{E:UPPERANDLOWERBOUNDSFORBLOWUPTIME}
		\tstar
	\exp
	\left(
	\frac{\frac{49}{50}}{\lifespanconstant \datasize}
	\right)
	\leq
	\Tcrease
	&
	\leq
	\tstar
	\exp
	\left(
	\frac{\frac{51}{50}}{\lifespanconstant \datasize}
	\right).
\end{align}

\medskip

\noindent \underline{\textbf{Bound for the logarithmic ratio of any two time values}}.
If $(t_1,r_1), (t_2,r_2) \in \mathbf{GH}_{\textnormal{Boot}}^{{\textnormal{INT}}}$ and $t_1 \leq t_2$,
then the following estimate holds:
\begin{align} \label{E:LOGOFTIMERATIOSBOUNDEDBYCOVERDATASIZEININTERIORREGION}
			\ln \left( \frac{1 + t_1}{1 + t_0} \right)
		& \leq \frac{\frac{51}{50}}{\lifespanconstant \datasize}.
\end{align}

\end{lemma}

\begin{proof}
\eqref{E:UPPERANDLOWERBOUNDSFORBLOWUPTIME} follows from the fact that $\tstar := \frac{1}{\datasize}$,
the bound $\ucrease = \mylittleo(\datasize)$ proved in \eqref{E:UCREASEISSMALL},
and \eqref{E:BLOWUPTIMEUPPERANDLOWERBOUNDS}.

\eqref{E:LOGOFTIMERATIOSBOUNDEDBYCOVERDATASIZEININTERIORREGION} follows from
\eqref{E:UPPERANDLOWERBOUNDSFORBLOWUPTIME} and our 
assumptions that $\tstar \leq t_1 \leq t_2 \leq 4 \Tcrease$.

\end{proof}

In the next lemma, we use the bootstrap assumptions to derive some preliminary estimates along the 
integral curves of $\Lunit$ and $\uLunit$. We will use these estimates in our proof of
Prop.\,\ref{P:APRIORIESTIMATESININTERIORREGION}, which yields our main Interior Region a priori estimates.

\begin{lemma}[Estimates along integral curves of $\Lunit$ and $\uLunit$]
\label{L:INTERIORREGIONESTIMATESALONGINTEGRALCURVES}
Under the bootstrap assumptions of Sect.\,\ref{SSS:INTERIORBOOTSTRAPASSUMPTIONS},
if $\datasize$ is sufficiently small,
then the following results hold.

\medskip

\noindent \underline{\textbf{Behavior of coordinate functions along integral curves of $\uLunit$}}.
Let $(t_1,r_1) \in \mathbf{GH}_{\textnormal{Boot}}^{{\textnormal{INT}}}$ 
be a point on the future-directed integral curve of $\uLunit$ emanating from
$(t_0,r_0) \in \mathbf{GH}_{\textnormal{Boot}}^{{\textnormal{INT}}}$. 
Then the following estimates hold:
\begin{subequations}
\begin{align} \label{E:INTERIORREGIONPOINTWISEESTIMATECHANGEINTPLUSRALONGINTEGRALCURVESOFLBAR}
	t_1 + r_1 
	& = t_0 + r_0
	+ \mathcal{O}(\eps) \frac{t_1 - t_0}{(t_1 + r_1)},
		\\
	t_1 - r_1
	& = t_0 - r_0
		+ 
		2(t_1 - t_0)
		+
		\mathcal{O}(\eps) \frac{t_1 - t_0}{(t_1 + r_1)},
		\label{E:INTERIORREGIONPOINTWISEESTIMATECHANGEINTMINUSRALONGINTEGRALCURVESOFLBAR}
			\\
\frac{1}{t_1 + r_1}
	& = \frac{1}{t_0 + r_0}
	+ 
	\mathcal{O}(\eps) \frac{t_1 - t_0}{(t_1 + r_1)^2(t_0 + r_0)},
	\label{E:RECIPROCALTPLUSRNEARLYCONSTANTALONGINTEGRALCURVESOFULUNIT}
		\\
	\frac{1}{2(1 + t_0 + r_0)}
	& \leq	
	\frac{1}{1 + t_1 + r_1}
	\leq \frac{2}{1 + t_0 + r_0}.
	\label{E:RECIPROCALONEPLUSTPLUSRNEARLYCONSTANTALONGINTEGRALCURVESOFULUNIT}
\end{align}
\end{subequations}

\medskip

\noindent \underline{\textbf{Estimates for integrals along integral curves of $\uLunit$}}.
Let $t \rightarrow (t,\underline{\mathfrak{r}}(t))$ be the integral curve of $\uLunit$ in
$\mathbf{GH}_{\textnormal{Boot}}^{{\textnormal{INT}}}$ joining the point $(t_1,r_1)$ to the point
$(t_2,r_2)$ lying to its future.

Let $B \geq 0$ be a constant. Then the following estimates hold, where the implicit constants
depend on $B$:
\begin{subequations}
\begin{align} \label{E:INTERIORFIRSTINTEGRALESTIMATEALONGINTEGRALCURVESOFLBAR}
	\int_{t_1}^{t_2} \frac{1}{[1 + t + \underline{\mathfrak{r}}(t)]^B
		(1 + |t - \underline{\mathfrak{r}}(t)|)} 
	\, \mathrm{d} t
	& \lesssim 
		\frac{\ln\left(\frac{1 + |t_2 - r_2|}{1 + |t_1 - r_1|} \right)}{(1 + t_2 + r_2)^B},
			\\
	\int_{t_1}^{t_2} \frac{1}{[1 + t + \underline{\mathfrak{r}}(t)]^B
		(1 + |t - \underline{\mathfrak{r}}(t)|)^2} 
	\, \mathrm{d} t
	& \lesssim 
		\frac{1}{(1 + t_2 + r_2)^B (1 + |t_1 - r_1|)}.
		\label{E:INTERIORFIRSTINTEGRALESTIMATEONEOVERONEPLUSTPLUSRTOATIMESONEOVERTMINUSRSQUAREDALONGINTEGRALCURVESOFLBAR}
\end{align}
\end{subequations}

In addition to the above hypotheses, assume that $(t_1,r_1)$ satisfies $|t_1-r_1|^2 \geq t_1 \ln(t_1-r_1)$, and assume that 
$(t_2,r_2) = (t_2,0)$, i.e., that $r_2 = 0$ and thus
$(t_2,r_2)$ lies on the time axis. Then the following estimates hold:
\begin{subequations}
\begin{align} \label{E:THIRDINTEGRALESTIMATEALONGINTEGRALCURVESOFLBAR}
	\int_{t_1}^{t_2} \frac{\underline{\mathfrak{r}}(t) \ln(t - \underline{\mathfrak{r}}(t))}{(1 + |t - \underline{\mathfrak{r}}(t)|)^2} \, \mathrm{d} t
	& \lesssim
		\frac{r_1^2}{1 + t_1 + r_1},
\end{align}

\begin{align} \label{E:FOURTHINTEGRALESTIMATEALONGINTEGRALCURVESOFLBAR}
	\int_{t_1}^{t_2} \frac{\underline{\mathfrak{r}}(t) \ln(t - \underline{\mathfrak{r}}(t))}{(1 + |t - \underline{\mathfrak{r}}(t)|)^3} \, \mathrm{d} t
	& \lesssim
		\frac{r_1^2}{(1 + t_1 + r_1)(1 + |t_1 - r_1|)}.
\end{align}
\end{subequations}

Let $t \rightarrow (t,\underline{\mathfrak{r}}(t))$ be the integral curve of $\uLunit$ in
$\mathbf{GH}_{\textnormal{Boot}}^{{\textnormal{INT}}}$ joining the point $(t_1,r_1)$ to the point
$(t_0,r_0)$ lying to its past.  Let $\updelta > 0$ be a constant.
Then the following estimate holds, where the implicit constant depends on $\updelta$:
\begin{align} \label{E:NEARWAVEZONEINTEGRALESTIMATELOGTMINUSROVER1PLUSTPLUSRTIMES1OVERTMINUSRTSQUAREDALONGINTEGRALCURVESOFLBAR}
\int_{t_0}^{t_1} 
	\frac{\underline{\mathfrak{r}}(t)}{(1 + |t - \underline{\mathfrak{r}}(t)|)^{1 + \updelta}} \, \mathrm{d} t
	& 
	\lesssim 
	r_0.
\end{align}
Moreover, if the two points also satisfy
$(t_0,r_0) \in \left(\Sigma_{\tstar} \cap \lbrace r \leq \frac{\tstar}{2} \rbrace \right) \bigcup \left( \lbrace u = \frac{\tstar}{2} \rbrace \cap \lbrace \tstar \leq t \leq \TCH\left(\frac{\tstar}{2} \right) \rbrace \right)$
and $|t_1-r_1|^2 \leq t_1 \ln(t_1-r_1)$,
then the following estimate holds:
\begin{align} \label{E:COMPARISONOFROAND1PLUSTPLUSR1INWAVEZONE}
	r_0
	& \approx 
	r_1
	\approx
	1 + t_1 + r_1.
\end{align}

\medskip

\noindent \underline{\textbf{Behavior of coordinate functions along integral curves of $\Lunit$}}.
Let $(t_1,r_1) \in \mathbf{GH}_{\textnormal{Boot}}^{{\textnormal{INT}}}$ be a point on the future-directed integral curve of $\Lunit$ emanating from
$(t_0,r_0) \in \mathbf{GH}_{\textnormal{Boot}}^{{\textnormal{INT}}}$. 
Then the following estimates hold on $\mathbf{GH}_{\textnormal{Boot}}^{{\textnormal{INT}}}$:

\begin{align} 
\begin{split} \label{E:BOUNDFORTPLUSRALONGINTEGRALCURVESOFL}
	t_1 + r_1 
	& = 
	t_0 + r_0
	+ 
	2(t_1 - t_0)
	+
	\mathcal{O}(\eps)
	\ln \left( \frac{1 + t_1}{1 + t_0} \right)
		\\
	& = t_0 + r_0
			+
			(1 + \mathcal{O}(\datasize^{3/4})) 2(t_1 - t_0),
	\end{split}
	\\
	\begin{split} \label{E:BOUNDFORTMINUSRALONGINTEGRALCURVESOFL}
	t_1 - r_1
	& = t_0 - r_0
		+
		\mathcal{O}(\eps)
		\ln \left( \frac{1 + t_1}{1 + t_0} \right)
			\\
	& =  (1 + \mathcal{O}(\datasize^{3/4})) (t_0 - r_0).
	\end{split}
\end{align}

\medskip 

\noindent \underline{\textbf{Estimates for integrals along integral curves of $\Lunit$}}.
Let $t \rightarrow (t,\mathfrak{r}(t))$ be the integral curve of $\Lunit$ in
$\mathbf{GH}_{\textnormal{Boot}}^{{\textnormal{INT}}}$ joining the point $(t_0,0)$ on the time axis to the point
$(t_1,r_1)$ lying to its future. Assume that $r_1 \leq t_0$.
Then for any constants $A \geq 0$ and $B \geq 0$, the following
estimates hold, where the implicit constants depend on $A$ and $B$:
\begin{subequations}
\begin{align} \label{E:NEARTIMEAXISINTEGRALESTIMATEROVER1PLUSTPLUSRPOWERAALONGINTEGRALCURVESOFL}
	\int_{t_0}^{t_1} \frac{\mathfrak{r}^A(t)}{[1 + t + \mathfrak{r}(t)]^B} \, \mathrm{d} t
	& \lesssim 
	\frac{r_1^{1+A}}{(1 + t_1 + r_1)^B},
		\\
	\int_{t_0}^{t_1} \frac{\mathfrak{r}^A(t) \ln(t - \mathfrak{r}(t))}{[1 + t + \mathfrak{r}(t)]^B} \, \mathrm{d} t
	& \lesssim 
	\frac{r_1^{1+A} \ln(t_1 - r_1)}{(1 + t_1 + r_1)^B}.
	\label{E:NEARTIMEAXISLOGINVOLVINGINTEGRALESTIMATEROVER1PLUSTPLUSRPOWERAALONGINTEGRALCURVESOFL}
\end{align}
\end{subequations}

Let $t \rightarrow (t,\mathfrak{r}(t))$ be the integral curve of $\Lunit$ in
$\mathbf{GH}_{\textnormal{Boot}}^{{\textnormal{INT}}}$ joining the point $(t_1,r_1)$ to the point
$(t_2,r_2)$ lying to its future. 
Then for any constant $B \geq 0$, the following
estimates hold, where the implicit constants depend on $B$:
\begin{subequations}
\begin{align} \label{E:FARFROMTIMEAXISINTEGRALESTIMATEROVER1PLUSTPLUSRSQUAREDALONGINTEGRALCURVESOFL}
	\int_{t_1}^{t_2} \frac{\mathfrak{r}(t)}{(1 + t + \mathfrak{r}(t))^2} \, \mathrm{d} t
	& \lesssim 
		\ln \left( \frac{1 + r_2}{1 + r_1} \right),
		\\
	\int_{t_1}^{t_2} \frac{\mathfrak{r}(t) \ln(t - \mathfrak{r}(t))}{(1 + t + \mathfrak{r}(t))^2} \, \mathrm{d} t
	& \lesssim 
	\left[ \ln(1 + r_2) \right]^2
	-
		\left[ \ln(1 + r_1) \right]^2,
	\label{E:FARFROMTIMEAXISLOGINVOLVINGINTEGRALESTIMATEROVER1PLUSTPLUSRSQUAREDALONGINTEGRALCURVESOFL}
		\\
		\int_{t_1}^{t_2} \frac{\mathfrak{r}(t)}{(1 + t + \mathfrak{r}(t))^3} \, \mathrm{d} t
	& \lesssim 
		\frac{1}{1 + t_1 + r_1},
		\label{E:FARFROMTIMEAXISINTEGRALESTIMATEROVER1PLUSTPLUSRCUBEDALONGINTEGRALCURVESOFL}
\end{align}
\end{subequations}

\begin{subequations}
\begin{align} \label{E:INTEGRALESTIMATE1OVER1PLUSTPLUSRSQUAREDTIMES1OVERTMINUSRALONGINTEGRALCURVESOFL}
	\int_{t_1}^{t_2} \frac{1}{(1 + t + \mathfrak{r}(t))^2(1 + |t - \mathfrak{r}(t)|)^B} \, \mathrm{d} t
	& \lesssim 
		\frac{1}{(1 + t_1 + r_1)(1 + |t_2 - r_2|)^B},
			\\
	\int_{t_1}^{t_2} \frac{\ln(t - \mathfrak{r}(t))}{(1 + t + \mathfrak{r}(t))(1 + |t - \mathfrak{r}(t)|)^B} \, \mathrm{d} t
	& \lesssim 
		\frac{\ln(t_2 - r_2) \ln\left( \frac{1 + t_2}{1 + t_1} \right)}{(1 + |t_2 - r_2|)^B},
		 \label{E:INTEGRALESTIMATELOGTMINUSROVER1PLUSTPLUSRTIMES1OVERTMINUSRTPOWERAALONGINTEGRALCURVESOFL}
			\\
		\int_{t_1}^{t_2} \frac{1}{(1 + t + \mathfrak{r}(t))(1 + |t - \mathfrak{r}(t)|)^B} \, \mathrm{d} t
	& \lesssim 
		\frac{\ln\left(\frac{1 + t_2}{1 + t_1} \right)}{(1 + |t_2 - r_2|)^B}.
		 \label{E:INTEGRALESTIMATE1OVER1PLUSTPLUSRTIMES1OVERTMINUSRTPOWERBALONGINTEGRALCURVESOFL}
	\end{align}
\end{subequations}

\end{lemma}

\begin{proof}
\noindent \underline{\textbf{Proof of \eqref{E:INTERIORREGIONPOINTWISEESTIMATECHANGEINTPLUSRALONGINTEGRALCURVESOFLBAR}--\eqref{E:RECIPROCALONEPLUSTPLUSRNEARLYCONSTANTALONGINTEGRALCURVESOFULUNIT}}}.
We first use \eqref{E:SPEEDOFSOUNDEXPANSION}--\eqref{E:SPEEDOFSOUNDERRORFUNCTIONVANISHESATORIGIN} and the bootstrap assumptions to deduce that:
\begin{align} \label{E:SPEEDPOINTWISEINTERIOR}
	\Speed & = 1 + \frac{\mathcal{O}(\eps)}{1 + t + r}.
\end{align}
Also using \eqref{E:NULLVECTORFIELDS},
we further deduce that:
\begin{align} \label{E:ULUNITRISMINUSONEPLUSSMALLERROR}
	\uLunit r 
	& = v^r - \Speed
		= - 1 + \frac{\mathcal{O}(\eps)}{1 + t + r}.
\end{align}
Since $\uLunit t = 1$, we also deduce that:
\begin{subequations}
\begin{align} \label{E:ULUNITTPLUSRISSMALLERROR}
	\uLunit (t + r) 
	& = \frac{\mathcal{O}(\eps)}{1 + t + r},
		\\
\uLunit (t - r) 
	& = 2  + \frac{\mathcal{O}(\eps)}{1 + t + r}
	\leq 2 +  \frac{\mathcal{O}(\eps)}{1 + t}.
	\label{E:ULUNITTMINUSRISTWOPLUSSMALLERROR}
\end{align}
\end{subequations}
Setting $f := 1 + t + r$, we deduce from \eqref{E:ULUNITTPLUSRISSMALLERROR} that
$|\uLunit (f^2)| = \mathcal{O}(\eps)$. Integrating this inequality along 
the integral curve of $\uLunit$ connecting $(t_0,r_0)$ to $(t_1,r_1)$
and using the fundamental theorem of calculus, 
we find that:
\begin{align} \label{E:NONLINEARINEQUALITYFORONEPLUSTPLUSRSQUAREDALONGINTEGRALCURVESOFULUNIT}
	(1 + t_0 + r_0)^2
	& = 
		(1 + t_1 + r_1)^2
		+ 
		\mathcal{O}(\eps) (t_1 - t_0)
	= (1 + t_1 + r_1)^2
		\left\lbrace
			1 + \mathcal{O}(\eps) \frac{(t_1 - t_0)}{(1 + t_1 + r_1)^2}
		\right\rbrace.
\end{align}
Taking the square root of \eqref{E:NONLINEARINEQUALITYFORONEPLUSTPLUSRSQUAREDALONGINTEGRALCURVESOFULUNIT}
and Taylor expanding
$
	\left\lbrace
			1 + \mathcal{O}(\eps) \frac{(t_1 - t_0)}{(1 + t_1 + r_1)^2}
		\right\rbrace^{1/2}
$
in the small quantity $\mathcal{O}(\eps) \frac{(t_1 - t_0)}{(1 + t_1 + r_1)^2}$,
we conclude \eqref{E:INTERIORREGIONPOINTWISEESTIMATECHANGEINTPLUSRALONGINTEGRALCURVESOFLBAR} and
\eqref{E:RECIPROCALONEPLUSTPLUSRNEARLYCONSTANTALONGINTEGRALCURVESOFULUNIT}.

\eqref{E:RECIPROCALTPLUSRNEARLYCONSTANTALONGINTEGRALCURVESOFULUNIT} follows from
multiplying \eqref{E:INTERIORREGIONPOINTWISEESTIMATECHANGEINTPLUSRALONGINTEGRALCURVESOFLBAR} by
$\frac{1}{(t_0 + r_0)(t_1 + r_1)}$.

Inequality \eqref{E:INTERIORREGIONPOINTWISEESTIMATECHANGEINTMINUSRALONGINTEGRALCURVESOFLBAR} 
is just an algebraic rearrangement of \eqref{E:INTERIORREGIONPOINTWISEESTIMATECHANGEINTPLUSRALONGINTEGRALCURVESOFLBAR}.

\medskip

\noindent \underline{\textbf{Proof of \eqref{E:INTERIORFIRSTINTEGRALESTIMATEALONGINTEGRALCURVESOFLBAR}--\eqref{E:INTERIORFIRSTINTEGRALESTIMATEONEOVERONEPLUSTPLUSRTOATIMESONEOVERTMINUSRSQUAREDALONGINTEGRALCURVESOFLBAR}}}.
Using \eqref{E:RECIPROCALONEPLUSTPLUSRNEARLYCONSTANTALONGINTEGRALCURVESOFULUNIT},
\eqref{E:ULUNITTMINUSRISTWOPLUSSMALLERROR}, and the fundamental theorem of calculus,
we see that LHS~\eqref{E:INTERIORFIRSTINTEGRALESTIMATEALONGINTEGRALCURVESOFLBAR} is:
\begin{align} 
\begin{split} \label{E:FIRSTPROOFSTEPINTERIORFIRSTINTEGRALESTIMATEALONGINTEGRALCURVESOFLBAR}
& 
\lesssim
\frac{1}{(1 + t_2 + r_2)^B}
\int_{t_1}^{t_2} 
		\frac{1}{(1 + |t - \underline{\mathfrak{r}}(t)|)} 
	\, \mathrm{d} t
	\\
& \lesssim
\frac{1}{(1 + t_2 + r_2)^B}
\int_{t_1}^{t_2} 
		\uLunit \ln(1 + |t - \underline{\mathfrak{r}}(t)|)
	\, \mathrm{d} t
		\\
& =
	\frac{1}{(1 + t_2 + r_2)^B}
	\ln\left(\frac{1 + |t_2 - r_2|}{1 + |t_1 - r_1|} \right),
\end{split}
\end{align}
which yields the desired bound \eqref{E:INTERIORFIRSTINTEGRALESTIMATEALONGINTEGRALCURVESOFLBAR} easily follows.
\eqref{E:INTERIORFIRSTINTEGRALESTIMATEONEOVERONEPLUSTPLUSRTOATIMESONEOVERTMINUSRSQUAREDALONGINTEGRALCURVESOFLBAR}
can be proved via similar arguments, and we omit the details.

\medskip
\noindent \underline{\textbf{Proof of \eqref{E:THIRDINTEGRALESTIMATEALONGINTEGRALCURVESOFLBAR}--\eqref{E:FOURTHINTEGRALESTIMATEALONGINTEGRALCURVESOFLBAR}}}.
Let 
\begin{align} \label{E:FUNCTIONTODIVIDEINTERIORREGIONFORRPLUSESTIMATES}
	f(t,r) : = (t-r)^2 - t\ln(t-r).
\end{align}
We compute that:
$\uLunit f 
= 2\frac{\uLunit(t-r)}{t-r} f 
+ 
2 t \uLunit(t-r) \frac{\ln(t-r) - \frac{1}{2}}{t-r}
- 
\ln(t-r)
$.
From this identity,
the bootstrap assumption $t- r \geq \frac{\tstar}{4}$, and \eqref{E:ULUNITTMINUSRISTWOPLUSSMALLERROR}, 
we deduce that: 
\begin{align} \label{E:CLOSEDINEQUALITYFORDETERMININGNEARTIMEAXISREGION}
	\uLunit f \geq \frac{2}{t-r} f.
\end{align}
By assumption, $f(t_1,r_1) \geq 0$. From this inequality and \eqref{E:CLOSEDINEQUALITYFORDETERMININGNEARTIMEAXISREGION},
it follows that for any point $(t,r)$ on the future-directed integral curve of $\uLunit$ emanating from $(t_1,r_1)$
(in particular, $t \geq t_1$ and $r \leq r_1$),
we have $f(t,r) \geq 0$. This implies that for such points, we have:
\begin{align} \label{E:BOUNDFORRPLUSEVOLUTIONEQUATIONERRORTERMSNEARTIMEAXIS}
	\frac{\ln(t-r)}{(t-r)^2} 
	& \leq \frac{1}{t}
	\leq \frac{1}{t_1}.
\end{align}
Using \eqref{E:BOUNDFORRPLUSEVOLUTIONEQUATIONERRORTERMSNEARTIMEAXIS},
we deduce that LHS~\eqref{E:THIRDINTEGRALESTIMATEALONGINTEGRALCURVESOFLBAR} is:
\begin{align} \label{E:FIRSTPROOFSTEPTHIRDINTEGRALESTIMATEALONGINTEGRALCURVESOFLBAR}
	& \lesssim 
	\frac{1}{t_1}
	\int_{t_1}^{t_2} \underline{\mathfrak{r}}(t)  \, \mathrm{d} t.
\end{align}
Using \eqref{E:ULUNITRISMINUSONEPLUSSMALLERROR}, we can change variables in the integral and use that $\underline{\mathfrak{r}}(t_2) = 0$
to deduce that RHS~\eqref{E:FIRSTPROOFSTEPTHIRDINTEGRALESTIMATEALONGINTEGRALCURVESOFLBAR} is:
\begin{align} \label{E:SECONDPROOFSTEPTHIRDINTEGRALESTIMATEALONGINTEGRALCURVESOFLBAR}
	& \lesssim 
	\frac{1}{t_1}
	\int_0^{t_1} r  \, \mathrm{d} r
	= \frac{r_1^2}{t_1}
	\lesssim
	\frac{r_1^2}{1 + t_1 + r_1},
\end{align}
which yields the desired bound \eqref{E:THIRDINTEGRALESTIMATEALONGINTEGRALCURVESOFLBAR}.

The bound \eqref{E:FOURTHINTEGRALESTIMATEALONGINTEGRALCURVESOFLBAR} can be proved through similar arguments,
and we omit the details.

\medskip
\noindent \underline{\textbf{Proof of \eqref{E:NEARWAVEZONEINTEGRALESTIMATELOGTMINUSROVER1PLUSTPLUSRTIMES1OVERTMINUSRTSQUAREDALONGINTEGRALCURVESOFLBAR}
and \eqref{E:COMPARISONOFROAND1PLUSTPLUSR1INWAVEZONE}}}.
By \eqref{E:ULUNITRISMINUSONEPLUSSMALLERROR}, the function $\underline{\mathfrak{r}}(t)$ is decreasing in $t$.
It follows that LHS~\eqref{E:NEARWAVEZONEINTEGRALESTIMATELOGTMINUSROVER1PLUSTPLUSRTIMES1OVERTMINUSRTSQUAREDALONGINTEGRALCURVESOFLBAR} is:
\begin{align} \label{E:FIRSTPROOFSTEPNEARWAVEZONEINTEGRALESTIMATELOGTMINUSROVER1PLUSTPLUSRTIMES1OVERTMINUSRTSQUAREDALONGINTEGRALCURVESOFLBAR}
& \lesssim
r_0
\int_{t_0}^{t_1} 
	\frac{1}{(1 + |t - \underline{\mathfrak{r}}(t)|)^{1 + \updelta}} \, \mathrm{d} t.
\end{align}
We can bound the integral in 
\eqref{E:FIRSTPROOFSTEPNEARWAVEZONEINTEGRALESTIMATELOGTMINUSROVER1PLUSTPLUSRTIMES1OVERTMINUSRTSQUAREDALONGINTEGRALCURVESOFLBAR}
via the change of variables $w(t) := t - \underline{\mathfrak{r}}(t)$ which, by 
\eqref{E:ULUNITTMINUSRISTWOPLUSSMALLERROR}, satisfies $\mathrm{d} w = (2 + \mathcal{O}(\eps)) \mathrm{d} t$.
Let $w_0 := w(t_0)$ and $w_1 := w(t_1)$.
Using \eqref{E:INTERIORTIMEBOUNDASSUMPTION} and \eqref{E:INTERIORREGIONTMINUSRLARGEBOOTSTRAPASSUMPTION}, we see that
$w_0 \geq \frac{\tstar}{4} > 0$ and $w_1 \leq 4 \Tcrease$.
It follows that RHS~\eqref{E:FIRSTPROOFSTEPNEARWAVEZONEINTEGRALESTIMATELOGTMINUSROVER1PLUSTPLUSRTIMES1OVERTMINUSRTSQUAREDALONGINTEGRALCURVESOFLBAR}
is:
\begin{align} \label{E:SECONDPROOFSTEPNEARWAVEZONEINTEGRALESTIMATELOGTMINUSROVER1PLUSTPLUSRTIMES1OVERTMINUSRTSQUAREDALONGINTEGRALCURVESOFLBAR}
& \lesssim
r_0
\int_{w_0}^{w_1} 
	\frac{1}{(1 + w)^{1 + \updelta}} \, \mathrm{d} w
\lesssim 
r_0
\int_{\frac{\tstar}{4 }}^{4 \Tcrease} 
	\frac{1}{(1 + w)^{1 + \updelta}} \, \mathrm{d} w
\lesssim r_0.
\end{align}
We have therefore proved \eqref{E:NEARWAVEZONEINTEGRALESTIMATELOGTMINUSROVER1PLUSTPLUSRTIMES1OVERTMINUSRTSQUAREDALONGINTEGRALCURVESOFLBAR}.

To prove \eqref{E:COMPARISONOFROAND1PLUSTPLUSR1INWAVEZONE}, we first use 
\eqref{E:COMPARISONBETWEENUANDMINUSRONCREASECHARACTERISTIC} and the fact that $\tstar = \frac{1}{\datasize}$
to deduce that for any point
$(t_0,r_0) \in \left(\Sigma_{\tstar} \cap \lbrace r \leq \frac{\tstar}{2} \rbrace \right) \bigcup 
\left( \lbrace u = \frac{\tstar}{2} \rbrace \cap \lbrace \tstar \leq t \leq \TCH\left(\frac{\tstar}{2} \right) \rbrace \right)$, 
we have $\frac{1}{2 \datasize} + \mathcal{O}(1) \leq t_0-r_0 \leq \frac{1}{\datasize}$.
Thus, by \eqref{E:ULUNITTMINUSRISTWOPLUSSMALLERROR} and our assumption that $|t_1-r_1|^2 \leq t_1 \ln(t_1-r_1)$ 
the difference $\Delta$ defined by:
\begin{align} \label{E:DEFINITIONCHANGEINTMINUSRALONGULUNITINTEGRALCURVE}
	\Delta & := (t_1- r_1) - (t_0 - r_0)
\end{align}
satisfies the following inequalities: 
\begin{align} \label{E:SMALLCHANGEINTMINUSRINWAVEZONE}
	0
	\leq
	\Delta 
	& \leq t_1^{1/2} \ln(t_1-r_1)^{1/2} - \frac{1}{2 \datasize} - \mathcal{O}(1).
\end{align}
From 
\eqref{E:INTERIORREGIONPOINTWISEESTIMATECHANGEINTPLUSRALONGINTEGRALCURVESOFLBAR}--\eqref{E:INTERIORREGIONPOINTWISEESTIMATECHANGEINTMINUSRALONGINTEGRALCURVESOFLBAR} and \eqref{E:SMALLCHANGEINTMINUSRINWAVEZONE},
we deduce that:
\begin{align} \label{E:RMINUSR0ISAPPROXIMATELYUNITYTIMEST0MINUST1}
		- 
		\Delta
		- 
		(t_0 - t_1)
		& =
		r_1 - r_0
		=
		(1 + \mathcal{O}(\eps))
		(t_0 - t_1)
		= 
		- 
		\left\lbrace \frac{1}{2} + \mathcal{O}(\eps) \right \rbrace \Delta.
\end{align}
From \eqref{E:RMINUSR0ISAPPROXIMATELYUNITYTIMEST0MINUST1}, it follows that:
\begin{subequations}
\begin{align}
	r_1 + \frac{\Delta}{3} 
	& \leq r_0 \leq r_1 + \Delta,
		\label{E:R1ISCLOSETOR0} \\
	t_0 + \frac{\Delta}{3} 
	& \leq t_1 \leq t_0 + \Delta.
	\label{E:T1ISCLOSETOT0}
\end{align}
\end{subequations}
From our assumption that $|t_1-r_1|^2 \leq t_1 \ln(t_1-r_1)$
and \eqref{E:INTERIORREGIONTMINUSRLARGEBOOTSTRAPASSUMPTION},
it follows that:
\begin{align} \label{E:T1NEARR1QUANTIFIEDUPPERANDLOWERBOUNDS}
	.9 t_1
	\leq
	t_1 + \mathcal{O}(t_1^{1/2} \ln^{1/2}(t_1))
	& 
	=
	r_1 
	\leq t_1.
\end{align}
Combining \eqref{E:SMALLCHANGEINTMINUSRINWAVEZONE},
\eqref{E:R1ISCLOSETOR0}--\eqref{E:T1ISCLOSETOT0},
and \eqref{E:T1NEARR1QUANTIFIEDUPPERANDLOWERBOUNDS}, we conclude 
\eqref{E:COMPARISONOFROAND1PLUSTPLUSR1INWAVEZONE}.

\medskip

\noindent \underline{\textbf{Proof of \eqref{E:BOUNDFORTPLUSRALONGINTEGRALCURVESOFL}--\eqref{E:BOUNDFORTMINUSRALONGINTEGRALCURVESOFL}}}.
These estimates can be proved via arguments similar to the ones we used to prove
\eqref{E:INTERIORREGIONPOINTWISEESTIMATECHANGEINTPLUSRALONGINTEGRALCURVESOFLBAR}--\eqref{E:RECIPROCALONEPLUSTPLUSRNEARLYCONSTANTALONGINTEGRALCURVESOFULUNIT},
and we omit the details, aside from noting that the last line of \eqref{E:BOUNDFORTMINUSRALONGINTEGRALCURVESOFL}
follows from 
\eqref{E:EPSILONSMALLNESSASSUMPTIONININTERIOREGION},
\eqref{E:INTERIORREGIONTMINUSRLARGEBOOTSTRAPASSUMPTION}, 
and \eqref{E:LOGOFTIMERATIOSBOUNDEDBYCOVERDATASIZEININTERIORREGION}.

\medskip

\noindent \underline{\textbf{Proof of \eqref{E:NEARTIMEAXISINTEGRALESTIMATEROVER1PLUSTPLUSRPOWERAALONGINTEGRALCURVESOFL}--\eqref{E:NEARTIMEAXISLOGINVOLVINGINTEGRALESTIMATEROVER1PLUSTPLUSRPOWERAALONGINTEGRALCURVESOFL}}}.
To prove \eqref{E:NEARTIMEAXISINTEGRALESTIMATEROVER1PLUSTPLUSRPOWERAALONGINTEGRALCURVESOFL}, we first use	
\eqref{E:BOUNDFORTMINUSRALONGINTEGRALCURVESOFL}
and our assumption that $r_1 \leq t_0$ 
to deduce that for all $t \in [t_0,t_1]$, we have
$1 + t + \mathfrak{r}(t) \approx 1 + t_0 \approx 1 + t_1 + r_1$.
This allows us to bound LHS~\eqref{E:NEARTIMEAXISINTEGRALESTIMATEROVER1PLUSTPLUSRPOWERAALONGINTEGRALCURVESOFL} by:
\begin{align} \label{E:FIRSTPROOFSTEPNEARTIMEAXISINTEGRALESTIMATEROVER1PLUSTPLUSRPOWERAALONGINTEGRALCURVESOFL}
	&
	\lesssim
	\frac{1}{(1 + t_1 + r_1)^B} 
	\int_{t_0}^{t_1} \mathfrak{r}^A(t) 
	\, \mathrm{d} t.
\end{align}
Next, we argue as in \eqref{E:SPEEDPOINTWISEINTERIOR}--\eqref{E:ULUNITRISMINUSONEPLUSSMALLERROR}
to deduce that:
\begin{align} \label{E:LUNITRISONEPLUSSMALLERROR}
	\Lunit r 
	& = v^r + \Speed
		= 1 + \frac{\mathcal{O}(\eps)}{1 + t + r}.
\end{align}
We can therefore change variables to $r = \mathfrak{r}(t)$ 
in the integral in \eqref{E:FIRSTPROOFSTEPNEARTIMEAXISINTEGRALESTIMATEROVER1PLUSTPLUSRPOWERAALONGINTEGRALCURVESOFL},
using \eqref{E:LUNITRISONEPLUSSMALLERROR} to deduce that $\mathrm{d} r = (1 + \mathcal{O}(\eps)) \mathrm{d} t$,
and recalling that $\mathfrak{r}(t_0) = 0$,
thereby bounding RHS~\eqref{E:FIRSTPROOFSTEPNEARTIMEAXISINTEGRALESTIMATEROVER1PLUSTPLUSRPOWERAALONGINTEGRALCURVESOFL} by: 
\begin{align} \label{E:SECONDPROOFSTEPNEARTIMEAXISINTEGRALESTIMATEROVER1PLUSTPLUSRPOWERAALONGINTEGRALCURVESOFL}
	&
	\lesssim
	\frac{1}{(1 + t_1 + r_1)^B} 
	\int_0^{r_1} r^A
	\, \mathrm{d} r
	\lesssim 
	\frac{r_1^{A+1}}{(1 + t_1 + r_1)^B}.
\end{align}
We have therefore proved \eqref{E:NEARTIMEAXISINTEGRALESTIMATEROVER1PLUSTPLUSRPOWERAALONGINTEGRALCURVESOFL}.

To prove \eqref{E:NEARTIMEAXISLOGINVOLVINGINTEGRALESTIMATEROVER1PLUSTPLUSRPOWERAALONGINTEGRALCURVESOFL},
we use \eqref{E:BOUNDFORTMINUSRALONGINTEGRALCURVESOFL} to bound the integrand factor
$\ln(t - \mathfrak{r}(t))$ by $\lesssim \ln(t_1 - r_1)$
and then appeal to the estimate \eqref{E:NEARTIMEAXISINTEGRALESTIMATEROVER1PLUSTPLUSRPOWERAALONGINTEGRALCURVESOFL}.

\medskip

\noindent \underline{\textbf{Proof of \eqref{E:FARFROMTIMEAXISINTEGRALESTIMATEROVER1PLUSTPLUSRSQUAREDALONGINTEGRALCURVESOFL}--\eqref{E:INTEGRALESTIMATE1OVER1PLUSTPLUSRTIMES1OVERTMINUSRTPOWERBALONGINTEGRALCURVESOFL}}}.
These estimates can all be proved via similar arguments, so we will only prove a representative one, 
namely \eqref{E:INTEGRALESTIMATELOGTMINUSROVER1PLUSTPLUSRTIMES1OVERTMINUSRTPOWERAALONGINTEGRALCURVESOFL}.
We first use \eqref{E:BOUNDFORTMINUSRALONGINTEGRALCURVESOFL} to bound 
the integrand factor $\frac{\ln(t - \mathfrak{r}(t))}{(1 + |t - \mathfrak{r}(t)|)^B}$
on LHS~\eqref{E:INTEGRALESTIMATELOGTMINUSROVER1PLUSTPLUSRTIMES1OVERTMINUSRTPOWERAALONGINTEGRALCURVESOFL}
by $\lesssim \frac{\ln(t_2 - r_2)}{(1 + |t_2 - r_2|)^B}$. It follows that  
LHS~\eqref{E:INTEGRALESTIMATELOGTMINUSROVER1PLUSTPLUSRTIMES1OVERTMINUSRTPOWERAALONGINTEGRALCURVESOFL}
is bounded by:
\begin{align} \label{E:FIRSTPROOFSTEPINTEGRALESTIMATELOGTMINUSROVER1PLUSTPLUSRTIMES1OVERTMINUSRTPOWERAALONGINTEGRALCURVESOFL}
& \lesssim
 \frac{\ln(t_2 - r_2)}{(1 + |t_2 - r_2|)^B}
\int_{t_1}^{t_2} 
	\frac{1}{1 + t} 
\, \mathrm{d} t,
\end{align}
from which the desired estimate \eqref{E:INTEGRALESTIMATELOGTMINUSROVER1PLUSTPLUSRTIMES1OVERTMINUSRTPOWERAALONGINTEGRALCURVESOFL} readily follows.
	
\end{proof}

\subsection{The main a priori estimates in the Interior Region}
\label{SS:MAINAPRIORIESTIMATESININTERIORREGION}
In this section, we prove Prop.\,\ref{P:APRIORIESTIMATESININTERIORREGION}, which yields
our main a priori estimates in the Interior Region. 


%



\begin{proposition}[The main a priori estimates in the Interior Region]
\label{P:APRIORIESTIMATESININTERIORREGION}
Under the data assumptions stated in Section~\ref{S:DATA}
	and the bootstrap assumptions of Sect.\,\ref{SSS:INTERIORBOOTSTRAPASSUMPTIONS},
if $\datasize$ is sufficiently small,
then the following pointwise estimates hold on $\mathbf{GH}_{\textnormal{Boot}}^{{\textnormal{INT}}}$.

\noindent \underline{\textbf{Estimates for the Riemann invariants}}.
\begin{subequations}
\begin{align} \label{E:POINTWISEAPRIORIINTERIORRPLUS}
	|\RRiemann|
	& 
	\leq C \frac{\datasize}{1 + t + r},
		\\
|\partial_t \RRiemann|
		& 
	\leq C \frac{\datasize}{(1 + t + r)(1 + |t-r|)},
		\label{E:POINTWISEAPRIORIINTERIORPARTIALTRPLUS}
			\\
	|\uLunit \RRiemann|
	& 
	\leq C \frac{\datasize}{(1 + t + r)(1 + |t-r|)},
		\label{E:POINTWISEAPRIORIINTERIORLBARRPLUS}
		\\
|r \uLunit \RRiemann|
	& 
	\leq
	C
	\frac{\datasize}{1 + t + r}
	+
	C \frac{\datasize \ln(t-r)}{(1 + |t - r|)^2},
		\label{E:POINTWISEAPRIORIINTERIORRTIMESLBARRPLUS}
			\\
|\Lunit \RRiemann|
	& 
	\leq 
		C \frac{\datasize}{(1 + t + r)^2},
		\label{E:POINTWISEAPRIORIINTERIORLRPLUS}
			\\
|r \Lunit \RRiemann|
	& 
	\leq 
		C \frac{\datasize}{1 + t + r},
		\label{E:POINTWISEAPRIORIINTERIORRTIMESLRPLUS}
			\\
|r \uLunit \partial_t \RRiemann|
		& 
	\leq  
	C \frac{\datasize}{(1 + t + r)(1 + |t-r|)}
		+
		C \frac{\datasize}{(1 + |t-r|)^3}
		+
		\frac{\datasize^{1+}}{(1 + |t-r|)^2}
		+
		C \frac{\datasize^{1+} \ln(t-r)}{(1 + |t-r|)^3},
		\label{E:POINTWISEAPRIORIINTERIORRTIMESLBARPARTIALTRPLUS}
			\\
|r \Lunit \uLunit \RRiemann|
	& 
	\leq 
	C \frac{\datasize}{(1 + t + r)(1 + |t-r|)},
	\label{E:POINTWISEAPRIORIINTERIORTIMESLUNITLBARRPLUS}
		\\
|r \uLunit \uLunit \RRiemann|
	& 
	\leq 
	C \frac{\datasize}{(1 + t + r)(1 + |t-r|)}
		+
		C \frac{\datasize}{(1 + |t-r|)^3}
		+
		C \frac{\datasize^{1+}}{(1 + |t-r|)^2}
		+
		C \frac{\datasize^{1+} \ln(t-r)}{(1 + |t-r|)^3},
	\label{E:POINTWISEAPRIORIINTERIORTIMESLBARLBARRPLUS}
	\\
|r \Lunit \Lunit \RRiemann|
	& 
	\leq 
	C \frac{\datasize}{(1 + t + r)^2},
	\label{E:POINTWISEAPRIORIINTERIORTIMESLUNITLUNITRPLUS}
\end{align}
\end{subequations}

\begin{subequations}
\begin{align} \label{E:POINTWISEAPRIORIINTERIORRMINUS}
	|\LRiemann|
	& 
	\leq C \frac{\datasize}{1 + t + r},
		\\
|\partial_t \LRiemann|
		& 
	\leq C \frac{\datasize}{(1 + t + r)^2},
		\label{E:POINTWISEAPRIORIINTERIORPARTIALTRMINUS}
			\\
|\uLunit \LRiemann|
		& 
	\leq C \frac{\datasize}{(1 + t + r)^2},
		\label{E:POINTWISEAPRIORIINTERIORLBARRMINUS}
		\\
	|r \uLunit \LRiemann|
	& 
	\leq C \frac{\datasize}{1 + t + r},
		\label{E:POINTWISEAPRIORIINTERIORRTIMESLBARRMINUS}
		\\
|\Lunit \LRiemann|
	& 
	\leq 
		C \frac{\datasize}{(1 + t + r)^2},
		\label{E:POINTWISEAPRIORIINTERIORLRMINUS}
			\\
|r \Lunit \LRiemann|
	& 
	\leq 
		C \frac{\datasize}{1 + t + r},
		\label{E:POINTWISEAPRIORIINTERIORRTIMESLRMINUS}
			\\
	|r \Lunit \Lunit \LRiemann|
	& 
	\leq 
	C \frac{\datasize}{(1 + t + r)^2},
	\label{E:POINTWISEAPRIORIINTERIORRTIMESLUNITLUNITRMINUS}
		\\
	|r \Lunit \partial_t \LRiemann|
		& 
	\leq 
	C \frac{\datasize}{(1 + t + r)^2},
		\label{E:POINTWISEAPRIORIINTERIORRTIMESLUNITPARTIALTRMINUS}
			\\
	|r \uLunit \Lunit \LRiemann|
	& 
	\leq 
	C \frac{\datasize}{(1 + t + r)^2},
	\label{E:POINTWISEAPRIORIINTERIORRTIMESLBARLUNITRMINUS}
		\\
	|r \uLunit \uLunit \LRiemann|
		& 
	\leq 
	C \frac{\datasize}{(1 + t + r)(1 + |t-r|)}.
		\label{E:POINTWISEAPRIORIINTERIORRTIMESLBARLBARRMINUS}
\end{align}
\end{subequations}

\medskip

\noindent \underline{\textbf{Estimates for the coordinates}}.
 \ \\
There exists a $C > 0$ such that for all points $(t,r) \in \mathbf{GH}_{\textnormal{Boot}}^{{\textnormal{INT}}}$, we have
the following estimates, where $\Tcrease$ is the value of $t$ at the crease (see \eqref{E:BLOWUPTIMEUPPERANDLOWERBOUNDS}):
\begin{subequations}
\begin{align} \label{E:INTERIORREGIONTMINUSRLARGEIMPROVEDESTIMATE}
	\frac{\tstar}{2}
	- C
	& 
	\leq
	t - r,	
		\\
	\tstar \leq t \leq 2 \Tcrease + C.
	 \label{E:INTERIORREGIONIMPROVEDMAXIMUMBOUNDONT}
\end{align}
\end{subequations}

\end{proposition}

\begin{remark}[Strict improvement of the Interior bootstrap assumptions]
	\label{R:INTERIORSTRICTIMPROVEMENTOFBOOTSTRAP}
	Recall that we assumed only that $0 < \eps \leq \datasize^{3/4}$
	in our bootstrap assumptions on $\mathbf{GH}_{\textnormal{Boot}}^{{\textnormal{INT}}}$ 
	(see \eqref{E:EPSILONSMALLNESSASSUMPTIONININTERIOREGION}).
	It is therefore easy to see that when $\datasize$ is sufficiently small, the estimates of Prop.\,\ref{P:APRIORIESTIMATESININTERIORREGION}
	yield strict improvements of all the bootstrap assumptions.
\end{remark}

\begin{remark}\label{Remark:improveboundsR+-:interior}
As in Remarks \ref{Remark:improveboundsR+-} and \ref{Remark:improveboundsR+-:global},
some of the estimates in Prop.\,\ref{P:APRIORIESTIMATESININTERIORREGION} could be improved, but we don't need such improvements
to close our bootstrap argument.
For example, in the Interior Region, the linear Riemann invariants could be shown to satisfy:
\begin{align}
|\flatRRiemann|&\leq C\frac{\datasize}{(1+t+r)(1+t-r)},
	\\
|\flatLRiemann|&\leq C\frac{\datasize\kh{1+\ln\frac{1+t+r}{1+t-r}}}{(1+t+r)^2},
\end{align}
and the nonlinear solution could also be shown to satisfy the same estimates.
\end{remark}

\begin{proof}[Proof of Prop.\,\ref{P:APRIORIESTIMATESININTERIORREGION}]
To help guide the reader, we start by noting that the ``data'' for the estimates in this proof are on the spacelike-characteristic surface portion
$\left(\Sigma_{\tstar} \cap \lbrace r \leq \frac{\tstar}{2} \rbrace \right) \bigcup \left( \lbrace u = \frac{\tstar}{2} \rbrace \cap \lbrace \tstar \leq t \leq \TCH\left(\frac{\tstar}{2} \right) \rbrace \right)$; see Fig.\,\ref{F:EXISTENCEANDCONTINUATION}. 
For some estimates, we will also rely on the boundary condition \eqref{E:RIEMANNINVARIANTMATCHINGCONDITION}.
Throughout the proof, we will silently use \eqref{E:SPEEDOFSOUNDEXPANSION}--\eqref{E:SPEEDOFSOUNDERRORFUNCTIONVANISHESATORIGIN},
which allow us to replace $\Speed$ with $1$ in various estimates, up to harmless error terms.

\medskip
\noindent \underline{\textbf{Proof of \eqref{E:POINTWISEAPRIORIINTERIORRMINUS} and \eqref{E:POINTWISEAPRIORIINTERIORRTIMESLRMINUS}}}.
We first use the wave equation \eqref{E:LEFTRIEMANNEQUATIONINTERIORREGION}, the bootstrap assumptions,
and the pointwise estimate \eqref{E:SPEEDPOINTWISEINTERIOR}
to deduce that:
\begin{align} \label{E:INTERIORREGIONWAVEEQUATIONINHOMOGENEOUSTERMPOINTWISEESTIMATEFORRMINUS}
	 \left|\uLunit \Big( r \Lunit \LRiemann +  2 \LRiemann\Big) \right|
	& \lesssim \frac{\eps^2}{(1 + t + r)^2(1 + |t-r|)}.
\end{align}
From \eqref{E:INTERIORREGIONWAVEEQUATIONINHOMOGENEOUSTERMPOINTWISEESTIMATEFORRMINUS},
the integral estimate \eqref{E:INTERIORFIRSTINTEGRALESTIMATEALONGINTEGRALCURVESOFLBAR}, 
the data bounds 
\eqref{E:TRCOORDINATESPURELDERIVATIVESTRANSPORTEDMODIFIEDLUNITRMINUSTSTARDATA} 
and
\eqref{E:CREASECHARACTERISTICPOINTWISEESTIMATEFORALLLUNITDERIVATIVESMODIFIEDLRMINUSDERIVATIVE},
\eqref{E:COMPARISONBETWEENUANDMINUSRONCREASECHARACTERISTIC},
and \eqref{E:INTERIORREGIONPOINTWISEESTIMATECHANGEINTPLUSRALONGINTEGRALCURVESOFLBAR},
we deduce that:
\begin{align} \label{E:INTERIORREGIONPPOINWISEESTIMATEFORRTIMESLUNITRIMINUSPLUS2RMINUS}
	r \Lunit \LRiemann +  2 \LRiemann 
	& =
		-
		2 \datasize 
		\frac{1}{1 + (t + r)^2}
		+
		\frac{\mathcal{O}(\datasize^{1^+}) \ln(|t-r|)}{(1 + t + r)^2}.
\end{align}
In particular, setting $r=0$ in \eqref{E:INTERIORREGIONPPOINWISEESTIMATEFORRTIMESLUNITRIMINUSPLUS2RMINUS}, we deduce that:
\begin{align} \label{E:RMINUSBOUNDONTIMEAXIS}
	\LRiemann \restriction_{\lbrace r=0 \rbrace}
	& = 
		-
		\datasize 
		\frac{1}{1 + t^2}
		+
		\frac{\mathcal{O}(\datasize^{1^+})\ln(t)}{(1 + t)^2}.
\end{align}
Next, we multiply \eqref{E:INTERIORREGIONPPOINWISEESTIMATEFORRTIMESLUNITRIMINUSPLUS2RMINUS} by $r$
and use the identity $\Lunit r = v^r + \Speed$ (see \eqref{E:NULLVECTORFIELDS}), 
the bootstrap assumptions, and \eqref{E:SPEEDPOINTWISEINTERIOR} to deduce the following pointwise bound:
\begin{align} \label{E:INTERIORREGIONPPOINWISEESTIMATEFORRTIMESLUNITRSQUAREDRMINUS}
	\Lunit(r^2 \LRiemann)
	& =
		-
		2 \datasize 
		\frac{r}{1 + (t + r)^2}
		+
		\mcl{O}(\datasize^{1^+}) \frac{\ln_+(t-r) r}{(1 + t + r)^2}.
\end{align}

The remainder of the proof of \eqref{E:POINTWISEAPRIORIINTERIORRMINUS} splits into two cases.

\noindent \textit{Case $1$: Estimates along integral curves of $\Lunit$ that emanate from the time axis}.
Let $(t_0,0)$ be a point on the time axis with $t_0 \geq \tstar$, 
and let $(t_1,r_1)$ be a point lying on the integral curve $t \rightarrow \upgamma_{t_0}(t)$ of $\Lunit$ emanating from
$(t_0,0)$ with $t_1 \geq t_0$. 

\noindent \textit{Sub-case $1a$: $r_1 \leq t_0$}. 
Integrating \eqref{E:INTERIORREGIONPPOINWISEESTIMATEFORRTIMESLUNITRSQUAREDRMINUS} along $\upgamma_{t_0}$ and using
the fundamental theorem of calculus, 
and 
\eqref{E:NEARTIMEAXISINTEGRALESTIMATEROVER1PLUSTPLUSRPOWERAALONGINTEGRALCURVESOFL}--\eqref{E:NEARTIMEAXISLOGINVOLVINGINTEGRALESTIMATEROVER1PLUSTPLUSRPOWERAALONGINTEGRALCURVESOFL},
we deduce that:
\begin{align} \label{E:NEARTIMEAXISPOINTWISEESTIMATESFORRSQUAREDTIMESRMINUS}
	|r_1^2 \LRiemann(t_1,r_1)|
	& \lesssim 
		\datasize \frac{r_1^2}{(1 + t_1 + r_1)^2}
		+
		\datasize^{1^+} r_1^2
		\frac{\ln(t_1 + r_1)}{(1 + t_1 + r_1)^2}.
\end{align}
Dividing \eqref{E:NEARTIMEAXISPOINTWISEESTIMATESFORRSQUAREDTIMESRMINUS} by $r_1^2$,
we conclude that:
\begin{align} \label{E:NEARTIMEAXISPOINTWISEESTIMATESFORRMINUS}
	|\LRiemann(t_1,r_1)|
	& \lesssim 
		 \frac{\datasize}{(1 + t_1 + r_1)^2}
		+
		\mathcal{O}(\datasize^{1^+})
		\frac{\ln(t_1+r_1)}{(1 + t_1 + r_1)^2},
\end{align}
which in particular implies a stronger estimate than \eqref{E:POINTWISEAPRIORIINTERIORRMINUS}.

\noindent \textit{Sub-case $1b$: $r_2 \geq t_0$}. 
We now consider the remaining sub-case, in which the point $(t_2,r_2)$ on $\upgamma_{t_0}$ 
satisfies $r_2 \geq t_0$. We fix $(t_1,r_1)$ to be the point on $\upgamma_{t_0}$ in which $r_1 = t_0$; the point $(t_1,r_1)$ was handled in Sub-case $1a$.
Using \eqref{E:BOUNDFORTMINUSRALONGINTEGRALCURVESOFL} with $(t_2,r_2)$ in the role of $(t_1,r_1)$ and $0$ in the role of $r_0$, we deduce that
$1.1 t_0 \geq t_2 - r_2 \geq .9 t_0$. Also using our assumption $r_2 \geq t_0$, we deduce that:
\begin{align}
	1.9 t_0 & \leq .9 t_0 + r_2 \leq t_2 \leq r_2 + 1.1 t_0 \leq 2.1 r_2.
\end{align}
Also using \eqref{E:INTERIORREGIONTMINUSRLARGEBOOTSTRAPASSUMPTION}--\eqref{E:TPLUSRANDTEQUIVALENTININTERIORREGION}
we deduce that:
\begin{align} \label{E:TANDRCOMPARABLEININTERIORREGIONAWAYFROMTIMEAXIS}
	r_2 
	& \approx 
	t_2
	\approx
	1 + t_2 + r_2.
\end{align} 
To proceed, we integrate \eqref{E:INTERIORREGIONPPOINWISEESTIMATEFORRTIMESLUNITRSQUAREDRMINUS} along $\upgamma_{t_0}$ from
from time $t_1$ to time $t_2$
and use the fundamental theorem of calculus,  
the already proven bound \eqref{E:NEARTIMEAXISPOINTWISEESTIMATESFORRSQUAREDTIMESRMINUS} at $(t_1,r_1)$,
and \eqref{E:FARFROMTIMEAXISINTEGRALESTIMATEROVER1PLUSTPLUSRSQUAREDALONGINTEGRALCURVESOFL}--\eqref{E:FARFROMTIMEAXISLOGINVOLVINGINTEGRALESTIMATEROVER1PLUSTPLUSRSQUAREDALONGINTEGRALCURVESOFL}
to deduce that:
\begin{align} \label{E:FARFROMTIMEAXISPOINTWISEESTIMATESFORRSQUAREDTIMESRMINUS}
	\begin{split}
	|r_2^2 \LRiemann(t_2,r_2)|
	&
	\leq
	|r_1^2 \LRiemann(t_1,r_1)|
	+
		C
		\datasize
		\ln \left( \frac{1 + r_2}{1 + r_1} \right)
		+
		C
		\datasize^{1^+}
		[\ln(1 + r_2)]^2
		\\
	& \lesssim 
		\datasize
		\ln(1 + r_2)
		+
		\datasize^{1^+}
		[\ln(1 + r_2)]^2.
\end{split}
\end{align}
Dividing \eqref{E:FARFROMTIMEAXISPOINTWISEESTIMATESFORRSQUAREDTIMESRMINUS} by $r_2^2$
and using \eqref{E:TANDRCOMPARABLEININTERIORREGIONAWAYFROMTIMEAXIS},
we deduce that:
\begin{align} \label{E:FARFROMTIMEAXISPOINTWISEESTIMATESFORRMINUS}
	\begin{split}
	|\LRiemann(t_2,r_2)|
	& \lesssim 
		\datasize
		\frac{\ln(t_2 + r_2)}{(1 + t_2 + r_2)^2}
		+
		\datasize^{1^+}
		\frac{[\ln(t_2 + r_2)]^2}{(1 + t_2 + r_2)^2},
\end{split}
\end{align}
which in particular implies a stronger estimate than \eqref{E:POINTWISEAPRIORIINTERIORRMINUS}.
We have therefore proved \eqref{E:POINTWISEAPRIORIINTERIORRMINUS} in Case $1$. 

\medskip

\noindent \textit{Case 2: Estimates along integral curves of $\Lunit$ that emanate from $\Sigma_{\tstar} \cap \lbrace r \leq \frac{\tstar}{2} \rbrace$}.
Let $(\tstar,r_0) \in \Sigma_{\tstar} \cap \lbrace r \leq \frac{\tstar}{2} \rbrace$,
and let $(t_1,r_1)$ be a point lying on the integral curve $t \rightarrow \upgamma_{r_0}(t)$ of $\Lunit$ emanating from
$(\tstar,r_0)$ with $t_1 \geq \tstar$. 

\noindent \textit{Sub-case 2a: $r_1 \leq \tstar$}. 
We argue as in the proof of \eqref{E:NEARTIMEAXISPOINTWISEESTIMATESFORRSQUAREDTIMESRMINUS}, 
but this time we incur a non-vanishing data term $r_0^2 \LRiemann(\tstar,r_0)$.
This yields the following bound:
\begin{align} \label{E:INTERIORREGIONNEARSIGMATSTARPOINTWISEESTIMATESFORRSQUAREDTIMESRMINUS}
	|r_1^2 \LRiemann(t_1,r_1)|
	& \lesssim 
		|r_0^2 \LRiemann(\tstar,r_0)|
		+
		\datasize \frac{r_1^2}{(1 + t_1 + r_1)^2}
		+
		\datasize^{1^+} r_1^2
		\frac{\ln(t_1 + r_1)}{(1 + t_1 + r_1)^2}.
\end{align}
Dividing \eqref{E:INTERIORREGIONNEARSIGMATSTARPOINTWISEESTIMATESFORRSQUAREDTIMESRMINUS} by $r_1^2$,
we find that:
\begin{align} \label{E:FIRSTVERSIONINTERIORREGIONNEARSIGMATSTARPOINTWISEESTIMATESFORRMINUS}
	|\LRiemann(t_1,r_1)|
	& \lesssim 
		\frac{r_0^2}{r_1^2} |\LRiemann(\tstar,r_0)| 
		+
		\datasize \frac{1}{(1 + t_1 + r_1)^2}
		+
		\datasize^{1^+}
		\frac{\ln(t_1 + r_1)}{(1 + t_1 + r_1)^2}.
\end{align}
Inequalities \eqref{E:BOUNDFORTPLUSRALONGINTEGRALCURVESOFL}--\eqref{E:BOUNDFORTMINUSRALONGINTEGRALCURVESOFL} and our assumptions that
$r_0 \leq \frac{\tstar}{2}$ and $r_1 \leq \tstar$ collectively imply that $r_1 \geq r_0$ and
$\tstar \approx 1 + \tstar + r_0 \approx 1 + t_1 + r_1$.
Combining these bounds with \eqref{E:FIRSTVERSIONINTERIORREGIONNEARSIGMATSTARPOINTWISEESTIMATESFORRMINUS} and the data bound 
\eqref{E:CRUDEPOINTWISEBOUNDFORDATAOFRMINUSUPTO3LDERIVATIVES} (in the case $J=0$),
we conclude that:
\begin{align} 
	\begin{split} \label{E:INTERIORREGIONNEARSIGMATSTARPOINTWISEESTIMATESFORRMINUS}
	|\LRiemann(t_1,r_1)|
	& \lesssim 
		\frac{r_0^2}{r_1^2} |\LRiemann(\tstar,r_0)| 
		+
		\datasize \frac{1}{(1 + t_1 + r_1)^2}
		+
		\datasize^{1^+}
		\frac{\ln(t_1 + r_1)}{(1 + t_1 + r_1)^2}
			\\
	& \lesssim
		\datasize \frac{1}{1 + t_1 + r_1}
		+
		\datasize^{1^+}
		\frac{\ln(t_1 + r_1)}{(1 + t_1 + r_1)^2},
\end{split}
\end{align}
which in particular implies a stronger estimate than \eqref{E:POINTWISEAPRIORIINTERIORRMINUS}.

\noindent \textit{Sub-case 2b: $r_2 \geq \tstar$}. 
We now consider the final sub-case in which the point $(t_2,r_2)$ on $\upgamma_{r_0}$ 
satisfies $r_2 \geq \tstar$. We fix $(t_1,r_1)$ to be the point on $\upgamma_{r_0}$ in which $r_1 = \tstar$; the point $(t_1,r_1)$ was handled in Sub-case 2a.
The proof of the first inequality in \eqref{E:FARFROMTIMEAXISPOINTWISEESTIMATESFORRSQUAREDTIMESRMINUS} goes through verbatim,
which yields:
\begin{align} \label{E:FARFROMSIGMATSTARPOINTWISEESTIMATESFORRSQUAREDTIMESRMINUS}
	|r_2^2 \LRiemann(t_2,r_2)|
	&
	\leq
	|r_1^2 \LRiemann(t_1,r_1)|
	+
		C
		\datasize
		\ln \left( \frac{1 + r_2}{1 + r_1} \right)
		+
		C
		\datasize^{1^+}
		[\ln(1 + r_2)]^2.
\end{align}
Multiplying \eqref{E:INTERIORREGIONNEARSIGMATSTARPOINTWISEESTIMATESFORRMINUS} by $r_1^2$ and using the resulting bound  
to control the first term on RHS~\eqref{E:FARFROMSIGMATSTARPOINTWISEESTIMATESFORRSQUAREDTIMESRMINUS},
we further deduce that:
\begin{align} \label{E:SECONDVERSIONFARFROMSIGMATSTARPOINTWISEESTIMATESFORRSQUAREDTIMESRMINUS}
	|r_2^2 \LRiemann(t_2,r_2)|
	&
	\lesssim
		\datasize r_1
		+
		\datasize
		\ln \left( \frac{1 + r_2}{1 + r_1} \right)
		+
		C
		\datasize^{1^+}
		[\ln(1 + r_2)]^2.
\end{align}
Using \eqref{E:BOUNDFORTMINUSRALONGINTEGRALCURVESOFL} with $(t_2,r_2)$ in the role of $(t_1,r_1)$ and $\tstar$ in the role of $t_0$, 
and our assumption that $r_0 \leq \frac{\tstar}{2}$, we deduce that
$t_2 - r_2 \leq \frac{3}{4} \tstar$. Since $r_2 \geq \tstar$ and $r_2 \leq t_2$, we further deduce that:
\begin{align} \label{E:TANDRCOMPARABLEININTERIORREGIONAWAYFROMSIGMATSTAR}
	r_2 
	& \approx 
	t_2
	\approx
	1 + t_2 + r_2.
\end{align}
Moreover, using \eqref{E:BOUNDFORTPLUSRALONGINTEGRALCURVESOFL} with $(t_2,r_2)$ in the role of $(t_1,r_1)$ and $\tstar$ in the role of $t_0$, 
and our assumption that $r_0 \leq \frac{\tstar}{2}$, we deduce the simple fact that $r_2 \geq r_1$.
Using this bound, \eqref{E:SECONDVERSIONFARFROMSIGMATSTARPOINTWISEESTIMATESFORRSQUAREDTIMESRMINUS},
and \eqref{E:TANDRCOMPARABLEININTERIORREGIONAWAYFROMSIGMATSTAR}, we divide
\eqref{E:SECONDVERSIONFARFROMSIGMATSTARPOINTWISEESTIMATESFORRSQUAREDTIMESRMINUS} by $r_2^2$ to conclude that:
\begin{align} 
\begin{split} \label{E:FARFROMSIGMATSTARPOINTWISEESTIMATESFORRMINUS}
	|\LRiemann(t_2,r_2)|
	&
	\lesssim
		\datasize \frac{r_1}{r_2^2}
		+
		\datasize
		\frac{1}{r_2^2}
		\ln \left( \frac{1 + r_2}{1 + r_1} \right)
		+
		C
		\datasize^{1^+}
		\frac{1}{r_2^2}
		[\ln(1 + r_2)]^2
			\\
	& \lesssim
	\frac{\datasize}{1 + t_2 + r_2}.
\end{split}
\end{align}
We have therefore proved \eqref{E:POINTWISEAPRIORIINTERIORRMINUS} in Case $2$, which finishes the proof of \eqref{E:POINTWISEAPRIORIINTERIORRMINUS}. 
From this estimate and \eqref{E:INTERIORREGIONPPOINWISEESTIMATEFORRTIMESLUNITRIMINUSPLUS2RMINUS}, we also deduce \eqref{E:POINTWISEAPRIORIINTERIORRTIMESLRMINUS}.

\medskip
\noindent \underline{\textbf{Proof of \eqref{E:POINTWISEAPRIORIINTERIORPARTIALTRMINUS} and \eqref{E:POINTWISEAPRIORIINTERIORRTIMESLUNITPARTIALTRMINUS}}}.
We commute the wave equation \eqref{E:LEFTRIEMANNEQUATIONINTERIORREGION} with $\partial_t$ 
and use \eqref{E:NULLVECTORFIELDS}, 
\eqref{E:PARTIALRINTERMSOFLUNITANDULUNIT}--\eqref{E:PARTIALTINTERMSOFLANDLBAR},
and the bootstrap assumptions, 
carefully noting that the products $2 \LRiemann \partial_t \uLunit \Speed$ 
(which are ``dangerous'' in that they involve a second derivative that lacks an $r$-weight)
cancel from both sides of the resulting equation,
thereby deducing that:
\begin{align} \label{E:INTERIORREGIONRMINUSWAVEEQUATIONCOMMUTEDWITHPARTIALTINHOMOGENEOUSTERMS}
	 \left|\uLunit \Big(r \Lunit \partial_t \LRiemann +  2 \Speed \partial_t \LRiemann \Big) \right|
	& \lesssim 
	\frac{\eps^2}{(1 + t + r)^2(1 + |t-r|)^2}.
\end{align}
Integrating \eqref{E:INTERIORREGIONRMINUSWAVEEQUATIONCOMMUTEDWITHPARTIALTINHOMOGENEOUSTERMS} and arguing
as in the proof of \eqref{E:INTERIORREGIONPPOINWISEESTIMATEFORRTIMESLUNITRIMINUSPLUS2RMINUS}, 
but this time using the integral estimate 
\eqref{E:INTERIORFIRSTINTEGRALESTIMATEONEOVERONEPLUSTPLUSRTOATIMESONEOVERTMINUSRSQUAREDALONGINTEGRALCURVESOFLBAR},
\eqref{E:INTERIORREGIONTMINUSRLARGEBOOTSTRAPASSUMPTION},
\eqref{E:LOGOFTIMERATIOSBOUNDEDBYCOVERDATASIZEININTERIORREGION},
\eqref{E:RECIPROCALONEPLUSTPLUSRNEARLYCONSTANTALONGINTEGRALCURVESOFULUNIT}, 
and the data bounds \eqref{E:TRCOORDINATESPARTIALDDERIVATIVETRANSPORTEDMODIFIEDLUNITRMINUSTSTARDATA} and 
\eqref{E:CREASECHARACTERISTICPOINTWISEESTIMATEFORPARTIALTMODIFIEDLRMINUSDERIVATIVE},
we deduce the following pointwise estimate:
\begin{align} \label{E:INTERIORREGIONPPOINWISEESTIMATEFORRTIMESLUNITPARTIALTRIMINUSQUANTITY}
	\begin{split}
	|r \Lunit \partial_t \LRiemann +  2 \Speed \partial_t \LRiemann|
	& = 
		\frac{\mathcal{O}(\datasize)}{(1 + t + r)^3}
		+
		\frac{\mathcal{O}(\eps^2 \datasize)}{(1 + t + r)^2}.
\end{split}
\end{align}
For future use, we set $r=0$ in \eqref{E:INTERIORREGIONPPOINWISEESTIMATEFORRTIMESLUNITPARTIALTRIMINUSQUANTITY} to deduce:
\begin{align} \label{E:POINTWISEBOUNDFORPARTIALTRMINUSONTIMEAXIS}
	\partial_t \LRiemann \restriction_{\lbrace r=0 \rbrace}
	& = \frac{\mathcal{O}(\datasize)}{(1 + t)^3}
			+
			\frac{\mathcal{O}(\eps^2 \datasize)}{(1 + t)^2}.
\end{align}

Next, we multiply \eqref{E:INTERIORREGIONPPOINWISEESTIMATEFORRTIMESLUNITPARTIALTRIMINUSQUANTITY} by $r$ 
and argue as in the proof of \eqref{E:INTERIORREGIONPPOINWISEESTIMATEFORRTIMESLUNITRSQUAREDRMINUS}
to deduce the following pointwise bound:
\begin{align} \label{E:INTERIORREGIONPPOINWISEESTIMATEFORRTIMESLUNITPARTIALTRSQUAREDTIMESRIMINUS}
	\begin{split}
	|\Lunit (r^2 \partial_t \LRiemann)|
	& \lesssim
	\datasize \frac{r}{(1 + t + r)^3}
	+
	\eps^2 \datasize \frac{r}{(1 + t + r)^2}.
\end{split}
\end{align}

The remainder of the proof of \eqref{E:POINTWISEAPRIORIINTERIORPARTIALTRMINUS} splits into two cases.

\noindent \textit{Case 1: Estimates along integral curves of $\Lunit$ that emanate from the time axis}.
Let $(t_0,0)$ be a point on the time axis with $t_0 \geq \tstar$, 
and let $(t_1,r_1)$ be a point lying on the $t-$parameterized integral curve $t \rightarrow \upgamma_{t_0}(t)$ of $\Lunit$ emanating from
$(t_0,0)$ with $t_1 \geq t_0$. 

\noindent \textit{Sub-case $1a$: $r_1 \leq t_0$}. 
Integrating \eqref{E:INTERIORREGIONPPOINWISEESTIMATEFORRTIMESLUNITPARTIALTRSQUAREDTIMESRIMINUS}
and arguing as in the proof of \eqref{E:NEARTIMEAXISPOINTWISEESTIMATESFORRSQUAREDTIMESRMINUS}, 
we deduce that:
\begin{align} \label{E:NEARTIMEAXISPOINTWISEESTIMATESFORRSQUAREDTIMESPARTIALTRMINUS}
	|r_1^2 \partial_t \LRiemann(t_1,r_1)|
	& \lesssim 
		\datasize \frac{r_1^2}{(1 + t_0)^3}
		+
		\eps^2 \datasize \frac{r_1^2}{(1 + t_0)^2}.
\end{align}
Dividing \eqref{E:NEARTIMEAXISPOINTWISEESTIMATESFORRSQUAREDTIMESPARTIALTRMINUS} by $r_1^2$ and using 
\eqref{E:RECIPROCALONEPLUSTPLUSRNEARLYCONSTANTALONGINTEGRALCURVESOFULUNIT}, we further deduce that:
\begin{align} \label{E:NEARTIMEAXISPOINTWISEESTIMATESFORPARTIALTRMINUS}
	|\partial_t \LRiemann(r_1,t_1)|
	& \lesssim 
		\datasize \frac{1}{(1 + t_1 + r_1)^3}
		+
		\eps^2 \datasize \frac{1}{(1 + t_1 + r_1)^2}.
\end{align}
We have therefore proved \eqref{E:POINTWISEAPRIORIINTERIORPARTIALTRMINUS} in Sub-case $1a$.

\noindent \textit{Sub-case $1b$: $r_2 \geq t_0$}. 
We now consider the remaining sub-case in which the point $(t_2,r_2)$ on $\upgamma_{t_0}$ 
satisfies $r_2 \geq t_0$. We fix $(t_1,r_1)$ to be the point on $\upgamma_{t_0}$ in which $r_1 = t_0$; the point $(t_1,r_1)$ was handled in Sub-case $1a$.

Integrating \eqref{E:INTERIORREGIONPPOINWISEESTIMATEFORRTIMESLUNITPARTIALTRSQUAREDTIMESRIMINUS}
along $\upgamma_{t_0}$ from time $t_1$ to time $t_2$ and using
the fundamental theorem of calculus, 
and \eqref{E:FARFROMTIMEAXISINTEGRALESTIMATEROVER1PLUSTPLUSRSQUAREDALONGINTEGRALCURVESOFL}--\eqref{E:FARFROMTIMEAXISINTEGRALESTIMATEROVER1PLUSTPLUSRCUBEDALONGINTEGRALCURVESOFL},
we deduce that:
\begin{align} \label{E:FARFROMTIMEAXISPOINTWISEESTIMATESFORRSQUAREDTIMESPARTIALTRMINUS}
	r_2^2 \partial_t \LRiemann(t_2,r_2)
	=
	r_1^2 \partial_t \LRiemann(t_1,r_1)
	+
	\frac{\mathcal{O}(\datasize)}{1 + t_1 + r_1}
	+
	\mcl{O}(\eps^2 \datasize) \ln\left( \frac{1 + r_2}{1 + r_1} \right).
\end{align}
Dividing \eqref{E:FARFROMTIMEAXISPOINTWISEESTIMATESFORRSQUAREDTIMESPARTIALTRMINUS} by $r_2^2$,
arguing as in the proof of \eqref{E:FARFROMSIGMATSTARPOINTWISEESTIMATESFORRMINUS}, 
and using the already proven bound \eqref{E:NEARTIMEAXISPOINTWISEESTIMATESFORPARTIALTRMINUS} at $(t_1,r_1)$
as well as \eqref{E:LOGOFTIMERATIOSBOUNDEDBYCOVERDATASIZEININTERIORREGION},
we deduce that:
\begin{align} \label{E:FARFROMTIMEAXISPOINTWISEESTIMATESFORPARTIALTRMINUS}
	\begin{split}
	|\partial_t \LRiemann(t_2,r_2)|
	& \lesssim 
		\frac{r_1^2}{r_2^2} |\partial_t \LRiemann(t_1,r_1)|
		+
		\frac{\datasize}{(1 + t_1 + r_1) r_2^2}
		+
		\frac{\eps^2 \datasize \left( \frac{1 + r_2}{1 + r_1} \right)}{r_2^2}
			\\
		& 
		\lesssim 
		\frac{\datasize}{(1 + t_2 + r_2)^2},
\end{split}
\end{align}
which in particular implies \eqref{E:POINTWISEAPRIORIINTERIORPARTIALTRMINUS} in Sub-case $1b$.
We have therefore proved \eqref{E:POINTWISEAPRIORIINTERIORPARTIALTRMINUS} in Case $1$.

\medskip

\noindent \textit{Case 2: Estimates along integral curves of $\Lunit$ that emanate from $\Sigma_{\tstar} \cap \lbrace r \leq \frac{\tstar}{2} \rbrace$}.
Let $(\tstar,r_0) \in \Sigma_{\tstar} \cap \lbrace r \leq \frac{\tstar}{2} \rbrace$,
and let $(t_1,r_1)$ be a point lying on the integral curve $\upgamma_{r_0}$ of $\Lunit$ emanating from
$(\tstar,r_0)$ with $t_1 \geq \tstar$. 

\noindent \textit{Sub-case 2a: $r_1 \leq \tstar$}. 
We argue as in the proof of 
\eqref{E:NEARTIMEAXISPOINTWISEESTIMATESFORRSQUAREDTIMESPARTIALTRMINUS}--\eqref{E:NEARTIMEAXISPOINTWISEESTIMATESFORPARTIALTRMINUS}, 
but this time we incur a non-vanishing data term $r_0^2 \partial_t\LRiemann(\tstar,r_0)$.
This yields the following bound:
\begin{align} \label{E:INTERIORREGIONNEARSIGMATSTARPOINTWISEESTIMATESFORRSQUAREDTIMESPARTIALTRMINUS}
	|r_1^2 \partial_t \LRiemann(t_1,r_1)|
	& \lesssim 
		|r_0^2 \partial_t \LRiemann(\tstar,r_0)|
		+
		\datasize \frac{r_1^2}{(1 + t_1 + r_1)^3}
		+
		\eps^2 \datasize \frac{r_1^2}{(1 + t_1 + r_1)^2}.
\end{align}
Dividing \eqref{E:INTERIORREGIONNEARSIGMATSTARPOINTWISEESTIMATESFORRSQUAREDTIMESPARTIALTRMINUS} by $r_1^2$
and using the data bound \eqref{E:CRUDEPOINTWISEBOUNDFORDATAOFPARTIALTRMINUS} as well as 
our assumption that $r_0 \leq \frac{\tstar}{2}$,
we find that:
\begin{align} \label{E:FIRSTVERSIONINTERIORREGIONNEARSIGMATSTARPOINTWISEESTIMATESFORPARTIALTRMINUS}
	|\partial_t \LRiemann(t_1,r_1)|
	& \lesssim 
		\frac{\datasize}{r_1^2} 
		+
		\datasize \frac{1}{(1 + t_1 + r_1)^3}
		+
		\eps^2 \datasize \frac{1}{(1 + t_1 + r_1)^2}.
\end{align}
From \eqref{E:FIRSTVERSIONINTERIORREGIONNEARSIGMATSTARPOINTWISEESTIMATESFORPARTIALTRMINUS}
and the estimates
$r_1 \geq r_0$ and
$\tstar \approx 1 + \tstar + r_0 \approx 1 + t_1 + r_1$ shown just above \eqref{E:INTERIORREGIONNEARSIGMATSTARPOINTWISEESTIMATESFORRMINUS},
we conclude that:
\begin{align} 
	\begin{split} \label{E:INTERIORREGIONNEARSIGMATSTARPOINTWISEESTIMATESFORPARTIALTRMINUS}
	|\partial_t \LRiemann(t_1,r_1)|
	& \lesssim 
		\frac{\datasize}{(1 + t_1 + r_1)^2},
\end{split}
\end{align}
which in particular implies \eqref{E:POINTWISEAPRIORIINTERIORPARTIALTRMINUS} in Sub-case $2a$.

\noindent \textit{Sub-case 2b: $r_2 \geq \tstar$}. 
We now consider the final sub-case in which the point $(t_2,r_2)$ on $\upgamma_{r_0}$ 
satisfies $r_2 \geq \tstar$. We fix $(t_1,r_1)$ to be the point on $\upgamma_{r_0}$ in which $r_1 = \tstar$; the point $(t_1,r_1)$ was handled in Sub-case 2a.
The proof of \eqref{E:FARFROMTIMEAXISPOINTWISEESTIMATESFORRSQUAREDTIMESPARTIALTRMINUS} goes through verbatim,
which yields:
\begin{align} \label{E:FARFROMSIGMATSTARPOINTWISEESTIMATESFORRSQUAREDTIMESPARTIALTRMINUS}
	|r_2^2 \partial_t \LRiemann(t_2,r_2)|
	&
	\leq
	|r_1^2 \partial_t \LRiemann(t_1,r_1)|
	+
		C
		\frac{\datasize}{1 + t_1 + r_1}
	+
	C
	\eps^2 \datasize \ln\left( \frac{1 + r_2}{1 + r_1} \right).
\end{align}
Multiplying \eqref{E:INTERIORREGIONNEARSIGMATSTARPOINTWISEESTIMATESFORPARTIALTRMINUS} by $r_1^2$,
inserting the resulting bound into 
RHS~\eqref{E:FARFROMSIGMATSTARPOINTWISEESTIMATESFORRSQUAREDTIMESPARTIALTRMINUS},
and also using \eqref{E:LOGOFTIMERATIOSBOUNDEDBYCOVERDATASIZEININTERIORREGION},
we further deduce that:
\begin{align} \label{E:SECONDVERSIONFARFROMSIGMATSTARPOINTWISEESTIMATESFORRSQUAREDTIMESPARTIALTRMINUS}
	|r_2^2 \partial_t \LRiemann(t_2,r_2)|
	&
	\lesssim
		\datasize.
\end{align}
The estimate \eqref{E:TANDRCOMPARABLEININTERIORREGIONAWAYFROMSIGMATSTAR} holds in the present context, and thus
we can divide \eqref{E:SECONDVERSIONFARFROMSIGMATSTARPOINTWISEESTIMATESFORRSQUAREDTIMESPARTIALTRMINUS} by $r_2^2$
to conclude that:
\begin{align} 
\begin{split} \label{E:FARFROMSIGMATSTARPOINTWISEESTIMATESFORPARTIALTRMINUS}
	|\partial_t \LRiemann(t_2,r_2)|
	&
	\lesssim
		 \frac{\datasize}{(1 + t_2 + r_2)^2}.
\end{split}
\end{align}
We have therefore proved \eqref{E:POINTWISEAPRIORIINTERIORPARTIALTRMINUS} in Case $2b$, 
which finishes the proof of \eqref{E:POINTWISEAPRIORIINTERIORPARTIALTRMINUS}. 

\eqref{E:POINTWISEAPRIORIINTERIORRTIMESLUNITPARTIALTRMINUS} then follows from
\eqref{E:POINTWISEAPRIORIINTERIORPARTIALTRMINUS}, \eqref{E:INTERIORREGIONPPOINWISEESTIMATEFORRTIMESLUNITPARTIALTRIMINUSQUANTITY},
and \eqref{E:SPEEDPOINTWISEINTERIOR}.

\medskip
\noindent \underline{\textbf{Proof of \eqref{E:POINTWISEAPRIORIINTERIORLRMINUS} and \eqref{E:POINTWISEAPRIORIINTERIORRTIMESLUNITLUNITRMINUS}}}.
We commute the wave equation \eqref{E:LEFTRIEMANNEQUATIONINTERIORREGION} with $\Lunit$
and use the identity $\Lunit r = v^r + \Speed$ (see \eqref{E:NULLVECTORFIELDS}), 
\eqref{E:COMMUTATOROFLANDULUNIT},
\eqref{E:SPEEDPOINTWISEINTERIOR},
and the bootstrap assumptions, 
and carefully note that the products $2 \LRiemann \Lunit \uLunit \Speed$ 
(which are ``dangerous'' in that they involve a second derivative that lacks an $r$-weight)
cancel from both sides of the commuted equation,
thereby deducing the following pointwise bound:
\begin{align} \label{E:INTERIORREGIONRMINUSWAVEEQUATIONCOMMUTEDWITHLINHOMOGENEOUSTERMS}
	 \left|\uLunit \Big(r \Lunit \Lunit \LRiemann +  3 \Speed \Lunit \LRiemann \Big) \right|
	& \lesssim 
			\frac{\eps^2}{(1 + t + r)^3(1 + |t-r|)}.
\end{align}
Integrating \eqref{E:INTERIORREGIONRMINUSWAVEEQUATIONCOMMUTEDWITHLINHOMOGENEOUSTERMS} and using the integral estimate
\eqref{E:INTERIORFIRSTINTEGRALESTIMATEALONGINTEGRALCURVESOFLBAR}, the data-bounds
\eqref{E:TRCOORDINATESPURELDERIVATIVESTRANSPORTEDMODIFIEDLUNITRMINUSTSTARDATA}
and \eqref{E:CREASECHARACTERISTICPOINTWISEESTIMATEFORALLLUNITDERIVATIVESMODIFIEDLRMINUSDERIVATIVE},
\eqref{E:COMPARISONBETWEENUANDMINUSRONCREASECHARACTERISTIC},
the identity $\Lunit r = v^r + \Speed$, 
\eqref{E:SPEEDPOINTWISEINTERIOR}, and the bootstrap assumptions,
we deduce the following pointwise bound:
\begin{align} \label{E:INTERIORREGIONPPOINWISEESTIMATEFORRTIMESLUNITLCOMMUTEDRIMINUSQUANTITY}
	\begin{split}
	|r \Lunit \Lunit \LRiemann +  3 \Lunit \LRiemann|
	& = 
		\frac{\mathcal{O}(\datasize)}{(1 + t + r)^3}
		+
		\frac{\mathcal{O}(\datasize^{1+}) \ln(|t-r|)}{(1 + t + r)^3}.
\end{split}
\end{align}

For future use, we set $r=0$ in \eqref{E:INTERIORREGIONPPOINWISEESTIMATEFORRTIMESLUNITLCOMMUTEDRIMINUSQUANTITY} to deduce:
\begin{align} \label{E:LDERIVATIVEOFRMINUSBOUNDONTIMEAXIS}
	\Lunit \LRiemann \restriction_{\lbrace r=0 \rbrace}
	& = \frac{\mathcal{O}(\datasize)}{(1 + t)^3}
		+
		\frac{\mathcal{O}(\datasize^{1+}) \ln(t)}{(1 + t)^3}.
\end{align}

Next, we multiply \eqref{E:INTERIORREGIONPPOINWISEESTIMATEFORRTIMESLUNITLCOMMUTEDRIMINUSQUANTITY} by $r^2$
and use the identity $\Lunit r = v^r + \Speed$, 
\eqref{E:SPEEDPOINTWISEINTERIOR}, and the bootstrap assumptions to deduce the following pointwise bound:
\begin{align} \label{E:INTERIORREGIONPPOINWISEESTIMATEFORLRCUBEDLRMINUS}
	\begin{split}
	|\Lunit (r^3 \Lunit \LRiemann)|
	& \lesssim
	\frac{\mathcal{O}(\datasize)r^2}{(1 + t + r)^3}
		+
		\frac{\mathcal{O}(\datasize^{1+}) \ln(|t-r|)r^2}{(1 + t + r)^3}.
\end{split}
\end{align}

The remainder of the proof of \eqref{E:POINTWISEAPRIORIINTERIORLRMINUS} splits into two cases.

\noindent \textit{Case $1$: Estimates along integral curves of $\Lunit$ that emanate from the time axis}.
Let $(t_0,0)$ be a point on the time axis with $t_0 \geq \tstar$, 
and let $(t_1,r_1)$ be a point lying on the $t$-parameterized integral curve $t \rightarrow \upgamma_{t_0}(t)$ of $\Lunit$ emanating from
$(t_0,0)$ with $t_1 \geq t_0$.

\noindent \textit{Sub-case $1a$: $r_1 \leq t_0$}. 
Starting from \eqref{E:INTERIORREGIONPPOINWISEESTIMATEFORLRCUBEDLRMINUS}, we can emulate the proof of
\eqref{E:NEARTIMEAXISPOINTWISEESTIMATESFORRMINUS} to deduce that:
\begin{align} \label{E:NEARTIMEAXISPOINTWISEESTIMATESFORLRMINUS}
	|\Lunit \LRiemann(t_1,r_1)|
	& \lesssim 
		\frac{\datasize}{(1 + t_1 + r_1)^3}
		+
		\frac{\datasize^{1+} \ln(t_1)}{(1 + t_1 + r_1)^3},
\end{align}
which in particular implies \eqref{E:POINTWISEAPRIORIINTERIORLRMINUS} in Sub-case $1a$.

\noindent \textit{Sub-case $1b$: $r_2 \geq t_0$}. 
We now consider the remaining sub-case in which the point $(t_2,r_2)$ on $\upgamma_{t_0}$, at which we would like to derive the pointwise estimate,
satisfies $r_2 \geq t_0$. We fix $(t_1,r_1)$ to be the point on $\upgamma_{t_0}$ in which $r_1 = t_0$; the point $(t_1,r_1)$ was handled in Sub-case $1a$.
Starting from \eqref{E:INTERIORREGIONPPOINWISEESTIMATEFORLRCUBEDLRMINUS}, and using the already proven bound
\eqref{E:NEARTIMEAXISPOINTWISEESTIMATESFORLRMINUS},
we can emulate the proof of \eqref{E:FARFROMTIMEAXISPOINTWISEESTIMATESFORRSQUAREDTIMESRMINUS}--\eqref{E:FARFROMTIMEAXISPOINTWISEESTIMATESFORRMINUS}
to deduce that:
\begin{align}
	r_2^3 \Lunit \LRiemann(t_2,r_2)
	=
	r_1^3 \Lunit \LRiemann(t_1,r_1)
	+
	\mathcal{O}(\datasize) \ln\left(\frac{t_2}{t_1} \right)
	+
	\mathcal{O}(\datasize^{1+}[\ln(t_2)]^2,
\end{align}
and thus, in view of 
\eqref{E:TANDRCOMPARABLEININTERIORREGIONAWAYFROMTIMEAXIS} 
and
\eqref{E:NEARTIMEAXISPOINTWISEESTIMATESFORLRMINUS}, 
that:
\begin{align} \label{E:FARFROMTIMEAXISPOINTWISEESTIMATESFORLRMINUS}
	\begin{split}
	|\Lunit \LRiemann(t_2,r_2)|
	&
		\lesssim 
		\frac{\datasize \ln(t_2+r_2)}{(1 + t_2 + r_2)^3}	
		+
		\frac{\datasize^{1^+} [\ln(t_2+r_2)]^2}{(1 + t_2 + r_2)^3}.
\end{split}
\end{align}
We have therefore proved \eqref{E:POINTWISEAPRIORIINTERIORLRMINUS} in Case $1b$, thereby finishing the proof in Case $1$.

\noindent \textit{Case 2: Estimates along integral curves of $\Lunit$ that emanate from $\Sigma_{\tstar} \cap \lbrace r \leq \frac{\tstar}{2} \rbrace$}.
Let $(\tstar,r_0) \in \Sigma_{\tstar} \cap \lbrace r \leq \frac{\tstar}{2} \rbrace$,
and let $(t_1,r_1)$ be a point lying on the $t$-parameterized integral curve $t \rightarrow \upgamma_{r_0}(t)$ of $\Lunit$ emanating from
$(\tstar,r_0)$ with $t_1 \geq \tstar$. 

\noindent \textit{Sub-case 2a: $r_1 \leq \tstar$}. 
Starting from \eqref{E:INTERIORREGIONPPOINWISEESTIMATEFORLRCUBEDLRMINUS},
we can emulate the proof of \eqref{E:INTERIORREGIONNEARSIGMATSTARPOINTWISEESTIMATESFORRSQUAREDTIMESRMINUS}, 
thereby deducing that:
\begin{align} \label{E:INTERIORREGIONNEARSIGMATSTARPOINTWISEESTIMATESFORRCUBEDTIMESLRMINUS}
	|r_1^3 \Lunit \LRiemann(t_1,r_1)|
	& \lesssim 
		|r_0^3 \Lunit \LRiemann(\tstar,r_0)| 
		+	
		\frac{\datasize r_1^3}{(1 + t_1 + r_1)^3}
		+
		\frac{\datasize^{1+} \ln(t_1) r_1^3}{(1 + t_1 + r_1)^3}.
\end{align}
Using the data bound \eqref{E:CRUDEPOINTWISEBOUNDFORDATAOFRMINUSUPTO3LDERIVATIVES}
to control the first term on RHS~\eqref{E:INTERIORREGIONNEARSIGMATSTARPOINTWISEESTIMATESFORRCUBEDTIMESLRMINUS},
dividing \eqref{E:INTERIORREGIONNEARSIGMATSTARPOINTWISEESTIMATESFORRCUBEDTIMESLRMINUS} by $r_1^3$, 
using our assumption that 
$r_1 \leq t_0$, and arguing as in the proof of \eqref{E:INTERIORREGIONNEARSIGMATSTARPOINTWISEESTIMATESFORRMINUS},
we conclude that:
\begin{align} \label{E:INTERIORREGIONNEARSIGMATSTARPOINTWISEESTIMATESFORLRMINUS}
	|\Lunit \LRiemann(t_1,r_1)|
	& \lesssim 
		\datasize \frac{1}{(1 + t_1 + r_1)^2}
		+
		\datasize^{1^+}
		\frac{\ln(t_1 + r_1)}{(1 + t_1 + r_1)^3}
	\lesssim 
	\datasize \frac{1}{(1 + t_1 + r_1)^2},
\end{align}
which yields \eqref{E:POINTWISEAPRIORIINTERIORLRMINUS} in Sub-case $2a$.

\noindent \textit{Sub-case 2b: $r_2 \geq \tstar$}. 
The comparison result $r_2 \approx t_2 \approx 1 + t_2 + r_2$ proved in \eqref{E:TANDRCOMPARABLEININTERIORREGIONAWAYFROMSIGMATSTAR}
also holds in the present context.
The desired bound \eqref{E:POINTWISEAPRIORIINTERIORLRMINUS} in Sub-case 2b therefore follows from the already proven bound
\eqref{E:POINTWISEAPRIORIINTERIORRTIMESLRMINUS}.

We have therefore proved \eqref{E:POINTWISEAPRIORIINTERIORLRMINUS}. 
\eqref{E:POINTWISEAPRIORIINTERIORRTIMESLUNITLUNITRMINUS} then follows from
\eqref{E:POINTWISEAPRIORIINTERIORLRMINUS} and \eqref{E:INTERIORREGIONPPOINWISEESTIMATEFORRTIMESLUNITLCOMMUTEDRIMINUSQUANTITY}.

\medskip
\noindent \underline{\textbf{Proof of \eqref{E:POINTWISEAPRIORIINTERIORLRPLUS} and \eqref{E:POINTWISEAPRIORIINTERIORLBARRMINUS}}}.
Using \eqref{E:PARTIALTINTERMSOFLANDLBAR} and the bootstrap assumptions,
we deduce that $|\uLunit \LRiemann| \lesssim |\partial_t \LRiemann| + |\Lunit \LRiemann|$.
From this estimate,
\eqref{E:POINTWISEAPRIORIINTERIORPARTIALTRMINUS},
and \eqref{E:POINTWISEAPRIORIINTERIORLRMINUS},
we conclude \eqref{E:POINTWISEAPRIORIINTERIORLBARRMINUS}.

\eqref{E:POINTWISEAPRIORIINTERIORLRPLUS} then follows from
\eqref{E:POINTWISEAPRIORIINTERIORLBARRMINUS} and the identity
$\Lunit \RRiemann = \uLunit \LRiemann$,
which follows from \eqref{E:OUTGOINGRIEMANNINVARIANTEVOLUTION}--\eqref{E:INGOINGRIEMANNINVARIANTEVOLUTION}.

\medskip
\noindent \underline{\textbf{Proof of \eqref{E:POINTWISEAPRIORIINTERIORRPLUS}}}.
Equation \eqref{E:OUTGOINGRIEMANNINVARIANTEVOLUTION} and the bootstrap assumptions imply that
$|\RRiemann| \lesssim |r \Lunit \RRiemann| + |\LRiemann|$.
From this estimate,
\eqref{E:POINTWISEAPRIORIINTERIORLRPLUS},
and
\eqref{E:POINTWISEAPRIORIINTERIORRMINUS},
we conclude the desired bound \eqref{E:POINTWISEAPRIORIINTERIORRPLUS}.

\medskip
\noindent \underline{\textbf{Proof of \eqref{E:POINTWISEAPRIORIINTERIORRTIMESLBARRPLUS}}}.
 We will prove the following pointwise estimate:
\begin{align} 
\begin{split}
\label{E:EVERYWHEREINTERIORREGIONPOINTWISEESTIMATEFORRTIMESULUNITRPLUSMINUSTWOSPEEDRPLUS}
	 \left|r \uLunit \RRiemann - 2 \RRiemann \right|
	& \lesssim 
	\frac{\datasize}{(1 + |t - r|)^2}
	+
	\frac{\datasize^{1^+} \ln(t-r)}{(1 + |t - r|)^2}.
\end{split}
\end{align}
The desired bound \eqref{E:POINTWISEAPRIORIINTERIORRTIMESLBARRPLUS} then follows from
\eqref{E:POINTWISEAPRIORIINTERIORRPLUS}
and \eqref{E:EVERYWHEREINTERIORREGIONPOINTWISEESTIMATEFORRTIMESULUNITRPLUSMINUSTWOSPEEDRPLUS}.

To prove \eqref{E:EVERYWHEREINTERIORREGIONPOINTWISEESTIMATEFORRTIMESULUNITRPLUSMINUSTWOSPEEDRPLUS}, 
we first use the wave equation \eqref{E:RIGHTRIEMANNEQUATIONINTERIORREGION}
\eqref{E:COMMUTATOROFLANDULUNIT},
\eqref{E:SPEEDPOINTWISEINTERIOR}, and the bootstrap assumptions to deduce the following bound,
where we highlight that a null-condition-failing product of type $r \uLunit \RRiemann \cdot \uLunit \RRiemann$
arises from the commutator term $r[\uLunit,\Lunit]\RRiemann$ on RHS~\eqref{E:RIGHTRIEMANNEQUATIONINTERIORREGION}:
\begin{align} \label{E:INTERIORREGIONWAVEEQUATIONINHOMOGENEOUSTERMPOINTWISEESTIMATEFORRPLUS}
	\left|\Lunit \Big( r \uLunit \RRiemann - 2 \RRiemann\Big) \right|
	& 
	\lesssim \frac{\eps^2}{(1 + t + r)^2(1 + |t-r|)}
	+
	\frac{\eps^2 \ln(t-r)}{(1 + t + r)(1 + |t-r|)^3}.
\end{align}
For future use, we note that the boundary condition \eqref{E:RIEMANNINVARIANTMATCHINGCONDITION} and \eqref{E:RMINUSBOUNDONTIMEAXIS} imply that:
\begin{align} \label{E:RPLUSBOUNDONTIMEAXIS}
		|\RRiemann \restriction_{r=0}|
		& \lesssim 
			\datasize
			\frac{1}{(1 + t)^2}
		+
			\frac{\mathcal{O}(\datasize^{1^+})\ln(t)}{(1 + t)^2}.
\end{align}

The remainder of the proof of \eqref{E:EVERYWHEREINTERIORREGIONPOINTWISEESTIMATEFORRTIMESULUNITRPLUSMINUSTWOSPEEDRPLUS}
splits into two cases.

\noindent \textit{Case $I$: Estimates along integral curves of $\Lunit$ that emanate from the time axis}.
Let $(t_0,0)$ be a point on the time axis with $t_0 \geq \tstar$, 
and let $(t_1,r_1)$ be a point lying on the integral curve $t \rightarrow \upgamma_{t_0}(t)$ of $\Lunit$ emanating from
$(t_0,0)$ with $t_1 \geq t_0$.

Integrating \eqref{E:INTERIORREGIONWAVEEQUATIONINHOMOGENEOUSTERMPOINTWISEESTIMATEFORRPLUS} from time $t_0$ to time $t_1$
and using 
\eqref{E:INTERIORREGIONTMINUSRLARGEBOOTSTRAPASSUMPTION},
\eqref{E:INTEGRALESTIMATE1OVER1PLUSTPLUSRSQUAREDTIMES1OVERTMINUSRALONGINTEGRALCURVESOFL}--\eqref{E:INTEGRALESTIMATELOGTMINUSROVER1PLUSTPLUSRTIMES1OVERTMINUSRTPOWERAALONGINTEGRALCURVESOFL},
\eqref{E:BOUNDFORTMINUSRALONGINTEGRALCURVESOFL},
\eqref{E:RPLUSBOUNDONTIMEAXIS},
\eqref{E:EPSILONSMALLNESSASSUMPTIONININTERIOREGION},
and \eqref{E:LOGOFTIMERATIOSBOUNDEDBYCOVERDATASIZEININTERIORREGION},
we find that:
\begin{align} \label{E:INTERIORRTIMESULUNITRIEMANNMINUS2RIEMANNPOINTWISEMANATINGFROMTIMEAXIS}
	\begin{split}
	|r_1 \uLunit \RRiemann(t_1,r_1) - 2 \RRiemann(t_1,r_1)|
	& = 
		|- 2 \RRiemann(t_0,0)|
		+ 
		\frac{\mathcal{O}(\eps^2)}{(1 + t_0)(1 + |t_1-r_1|)}
			\\
& \ \
		+
		\frac{\mathcal{O}(\eps^2) \ln\left(\frac{1 + t_1}{1+t_0} \right)\ln(t_1-r_1)}{(1 + |t_1-r_1|)^3}
			\\
		& 
		\lesssim 
		\frac{\datasize}{(1 + t_0)^2}
		+
		\frac{\datasize^{1^+}\ln(t_0)}{(1 + t_0)^2}
			\\
	& \ \
		+
		\frac{\mathcal{O}(\eps^2)}{(1 + |t_1-r_1|)^2}
		+
		\frac{\mathcal{O}(\eps^2) \ln\left(\frac{1 + t_1}{1+t_0} \right)\ln(t_1-r_1)}{(1 + |t_1-r_1|)^3}
			\\
	& 	
		\lesssim 
		\frac{\datasize}{1 + (t_1-r_1)^2}
		+
		\frac{\datasize^{1^+}\ln(|t_1-r_1|)}{(1 + |t_1-r_1|)^2}
		+
		\frac{\eps^2 \datasize \ln\left(\frac{1 + t_1}{1+t_0} \right)\ln(t_1-r_1)}{(1 + |t_1-r_1|)^2}
			\\
	& \lesssim 
		\frac{\datasize}{1 + (t_1-r_1)^2}
		+
		\frac{\datasize^{1^+}\ln(|t_1-r_1|)}{(1 + |t_1-r_1|)^2}.
\end{split}
\end{align}
We have therefore proved \eqref{E:EVERYWHEREINTERIORREGIONPOINTWISEESTIMATEFORRTIMESULUNITRPLUSMINUSTWOSPEEDRPLUS} in Case $I$.

\noindent \textit{Case $II$: Estimates along integral curves of $\Lunit$ that emanate from $\Sigma_{\tstar} \cap \lbrace r \leq \frac{\tstar}{2} \rbrace$}.
Let $(\tstar,r_0) \in \Sigma_{\tstar} \cap \lbrace r \leq \frac{\tstar}{2} \rbrace$, 
and let $(t_1,r_1)$ be a point lying on the integral curve $t \rightarrow \upgamma_{r_0}(t)$ of $\Lunit$ emanating from
$(\tstar,r_0)$ with $t_1 \geq t_0$.

Again integrating \eqref{E:INTERIORREGIONWAVEEQUATIONINHOMOGENEOUSTERMPOINTWISEESTIMATEFORRPLUS}
and using 
\eqref{E:INTERIORREGIONTMINUSRLARGEBOOTSTRAPASSUMPTION},
\eqref{E:INTEGRALESTIMATE1OVER1PLUSTPLUSRSQUAREDTIMES1OVERTMINUSRALONGINTEGRALCURVESOFL}--\eqref{E:INTEGRALESTIMATELOGTMINUSROVER1PLUSTPLUSRTIMES1OVERTMINUSRTPOWERAALONGINTEGRALCURVESOFL},
\eqref{E:BOUNDFORTMINUSRALONGINTEGRALCURVESOFL},
\eqref{E:TRCOORDINATESALLDERIVATIVESTRANSPORTEDMODIFIEDULUNITRPLUSTSTARDATA},
\eqref{E:LOGOFTIMERATIOSBOUNDEDBYCOVERDATASIZEININTERIORREGION},
and our assumption that $r_0 \leq \frac{\tstar}{2}$,
we deduce that:
\begin{align} 
\begin{split} \label{E:INTERIORREGIONPOINTWISEESTIMATEFORRTIMESULUNITRPLUSMINUSTWORPLUS}
	 \left|r_1 \uLunit \RRiemann(t_1,r_1) - 2 \RRiemann(t_1,r_1) \right|
	& \lesssim
		\frac{\datasize}{1 + (\tstar - r_0)^2}
			\\
	& \ \
		+
		\frac{\eps^2}{(1 + |t_1-r_1|)^2}
		+
		\frac{\eps^2 \ln\left(\frac{1 + t_1}{1+\tstar} \right)\ln(t_1-r_1)}{(1 + |t_1-r_1|)^3}
				\\
	& \lesssim
		\frac{\datasize}{1 + (t_1 - r_1)^2}
		+ 
		\frac{\mathcal{O}(\datasize^{1^+})\ln(|t_1-r_1|)}{(1 + |t_1-r_1|)^2}.
\end{split}
\end{align}
We have therefore proved \eqref{E:EVERYWHEREINTERIORREGIONPOINTWISEESTIMATEFORRTIMESULUNITRPLUSMINUSTWOSPEEDRPLUS} in Case $II$,
thereby finishing the proof of \eqref{E:EVERYWHEREINTERIORREGIONPOINTWISEESTIMATEFORRTIMESULUNITRPLUSMINUSTWOSPEEDRPLUS}.

\medskip 

\noindent \underline{\textbf{Proof of \eqref{E:POINTWISEAPRIORIINTERIORRTIMESLRPLUS} and \eqref{E:POINTWISEAPRIORIINTERIORRTIMESLBARRMINUS}}}.
These estimates follows from multiplying equations
\eqref{E:OUTGOINGRIEMANNINVARIANTEVOLUTION}--\eqref{E:INGOINGRIEMANNINVARIANTEVOLUTION} 
by $r$ and using
\eqref{E:POINTWISEAPRIORIINTERIORRPLUS}
and
\eqref{E:POINTWISEAPRIORIINTERIORRMINUS}.

\medskip 

\noindent \underline{\textbf{Proof of \eqref{E:POINTWISEAPRIORIINTERIORLBARRPLUS}}}.

\noindent \textit{Case $1$: Estimates along integral curves of $\Lunit$ that emanate from the time axis}.
Let $(t_0,0)$ be a point on the time axis with $t_0 \geq \tstar$, 
and let $(t_1,r_1)$ be a point lying on the integral curve $t \rightarrow \upgamma_{t_0}(t)$ of $\Lunit$ emanating from
$(t_0,0)$ with $t_1 \geq t_0$.
From equation \eqref{E:INTERIORREGIONWAVEEQUATIONINHOMOGENEOUSTERMPOINTWISEESTIMATEFORRPLUS} and the already proven estimate
\eqref{E:POINTWISEAPRIORIINTERIORLRPLUS}, we deduce:
\begin{align} \label{E:INTERIORREGIONLUNITRULUNITRPLUS}
	\left|\Lunit \Big( r \uLunit \RRiemann \Big) \right|
	& 
	\lesssim 
	\frac{\datasize}{(1 + t + r)^2}
	+
	\frac{\eps^2 \ln(t-r)}{(1 + t + r)(1 + |t-r|)^3}.
\end{align}

\noindent \textit{Sub-case $1a$: $r_1 \leq t_0$}. 
We integrate \eqref{E:INTERIORREGIONLUNITRULUNITRPLUS} along $\upgamma_{t_0}$ 
from time $t_0$ to time $t_1$
and use the integral estimates
\eqref{E:NEARTIMEAXISINTEGRALESTIMATEROVER1PLUSTPLUSRPOWERAALONGINTEGRALCURVESOFL}
and \eqref{E:NEARTIMEAXISLOGINVOLVINGINTEGRALESTIMATEROVER1PLUSTPLUSRPOWERAALONGINTEGRALCURVESOFL}
to deduce that:
\begin{align} \label{E:NEARTIMEAXISPOINTWISEESTIMATESFORRTIMESLRMINUS}
	|r_1 \uLunit \RRiemann(t_1,r_1)|
	& \lesssim 
		\frac{\datasize r_1}{(1 + t_1 + r_1)^2}
		+
		\frac{\datasize^{1+} r_1 \ln(t_1-r_1)}{(1 + t_1 + r_1)(1 + |t_1-r_1|)^3},
\end{align}
which in particular implies \eqref{E:POINTWISEAPRIORIINTERIORLBARRPLUS} in Sub-case $1a$.

\noindent \textit{Sub-case $1b$: $r_2 \geq t_0$}. 
We have already proved \eqref{E:POINTWISEAPRIORIINTERIORRTIMESLBARRPLUS}, i.e., that
$|r_2 \uLunit \RRiemann(t_2,r_2)|
 \leq
	C
	\frac{\datasize}{1 + t_2 + r_2}
	+
	C \frac{\datasize \ln(t_2-r_2)}{(1 + |t_2 - r_2|)^2}$.
Dividing this estimate by $r_2$ and using the bound \eqref{E:TANDRCOMPARABLEININTERIORREGIONAWAYFROMTIMEAXIS} 
(which holds in the present context), 
we conclude \eqref{E:POINTWISEAPRIORIINTERIORLBARRPLUS} in Case $1b$, thereby finishing the proof in Case $1$.

\medskip 	
\noindent \textit{Case 2: Estimates along integral curves of $\Lunit$ that emanate from $\Sigma_{\tstar} \cap \lbrace r \leq \frac{\tstar}{2} \rbrace$}.
Let $(\tstar,r_0) \in \Sigma_{\tstar} \cap \lbrace r \leq \frac{\tstar}{2} \rbrace$, 
and let $(t_1,r_1)$ be a point lying on the integral curve $t \rightarrow \upgamma_{r_0}(t)$ of $\Lunit$ emanating from
$(\tstar,r_0)$ with $t_1 \geq \tstar$. 

\medskip 
\noindent \textit{Sub-case 2a: $r_1 \leq \tstar$}. 
We integrate \eqref{E:INTERIORREGIONLUNITRULUNITRPLUS} and argue as in the proof of \eqref{E:NEARTIMEAXISPOINTWISEESTIMATESFORRTIMESLRMINUS},
but this time we incur a non-vanishing data term $r_0 \uLunit \RRiemann(\tstar,r_0)$, which we
bound with \eqref{E:CRUDEPOINTWISEBOUNDFORDATAOFRPLUSUPTO3DERIVATIVES}. This yields the following bound:
\begin{align}
	\begin{split}  \label{E:FIRSTPROOFBOUNDEMANATINGFORMSIGMATSTARPOINTWISEESTIMATESFORRTIMESLRMINUS}
	|r_1 \uLunit \RRiemann(t_1,r_1)|
	& \lesssim 
		\frac{\datasize r_0}{(1 + \tstar + r_0)(1 + |\tstar - r_0|)}
		+
		\frac{\datasize r_1}{(1 + t_1 + r_1)^2}
		+
		\frac{\datasize^{1+} r_1 \ln(t_1-r_1)}{(1 + t_1 + r_1)(1 + |t_1-r_1|)^3}
			\\
	& \lesssim
		\frac{\datasize}{(1 + \tstar + r_0)}
		+
		\frac{\datasize r_1}{(1 + t_1 + r_1)^2}
		+
		\frac{\datasize^{1+} r_1 \ln(t_1-r_1)}{(1 + t_1 + r_1)(1 + |t_1-r_1|)^3}.
	\end{split}
\end{align}
Dividing \eqref{E:FIRSTPROOFBOUNDEMANATINGFORMSIGMATSTARPOINTWISEESTIMATESFORRTIMESLRMINUS} by $r_1$ and using
the estimates
$r_1 \geq r_0$ and
$\tstar \approx 1 + \tstar + r_0 \approx 1 + t_1 + r_1$ shown just above \eqref{E:INTERIORREGIONNEARSIGMATSTARPOINTWISEESTIMATESFORRMINUS},
we find that:
\begin{align} \label{E:SECONDPROOFBOUNDEMANATINGFORMSIGMATSTARPOINTWISEESTIMATESFORRTIMESLRMINUS}
	|\uLunit \RRiemann(t_1,r_1)|
	& \lesssim
		\frac{\datasize}{(1 + t_1 + r_1)^2}
		+
		\frac{\datasize^{1+} \ln(t_1-r_1)}{(1 + t_1 + r_1)(1 + |t_1-r_1|)^3},
\end{align}
which yields \eqref{E:POINTWISEAPRIORIINTERIORLBARRPLUS} in Case $2a$.

\noindent \textit{Sub-case 2b: $r_2 \geq \tstar$}. 
The comparison result $r_2 \approx t_2 \approx 1 + t_2 + r_2$ proved in \eqref{E:TANDRCOMPARABLEININTERIORREGIONAWAYFROMSIGMATSTAR}
also holds in the present context.
The desired bound \eqref{E:POINTWISEAPRIORIINTERIORLBARRPLUS} in Sub-case 2b therefore follows from the already proven bound
\eqref{E:POINTWISEAPRIORIINTERIORRTIMESLBARRPLUS}. We have therefore proved \eqref{E:POINTWISEAPRIORIINTERIORLBARRPLUS}.

\medskip

\noindent \underline{\textbf{Proof of \eqref{E:POINTWISEAPRIORIINTERIORTIMESLUNITLUNITRPLUS} and \eqref{E:POINTWISEAPRIORIINTERIORRTIMESLBARLBARRMINUS}}}.
\eqref{E:POINTWISEAPRIORIINTERIORTIMESLUNITLUNITRPLUS} follows from multiplying equation \eqref{E:OUTGOINGRIEMANNINVARIANTEVOLUTION} by $r$,
taking a $\Lunit$ derivative,
and using the identity $\Lunit r = v^r + \Speed$ (see \eqref{E:NULLVECTORFIELDS}), the bootstrap assumptions,
and the already proven estimates
\eqref{E:POINTWISEAPRIORIINTERIORLRPLUS}
and \eqref{E:POINTWISEAPRIORIINTERIORLRMINUS}.

\eqref{E:POINTWISEAPRIORIINTERIORRTIMESLBARLBARRMINUS} follows from a similar argument based on
equation \eqref{E:INGOINGRIEMANNINVARIANTEVOLUTION}
and the already proven estimates
\eqref{E:POINTWISEAPRIORIINTERIORLBARRPLUS} and \eqref{E:POINTWISEAPRIORIINTERIORLBARRMINUS}.

\medskip

\noindent \underline{\textbf{Proof of \eqref{E:POINTWISEAPRIORIINTERIORTIMESLUNITLBARRPLUS} and \eqref{E:POINTWISEAPRIORIINTERIORRTIMESLBARLUNITRMINUS}}}.
To prove \eqref{E:POINTWISEAPRIORIINTERIORRTIMESLBARLUNITRMINUS}, we use use \eqref{E:INTERIORREGIONWAVEEQUATIONINHOMOGENEOUSTERMPOINTWISEESTIMATEFORRMINUS},
the identity $\uLunit r = v^r - \Speed$ (see \eqref{E:NULLVECTORFIELDS}),
the bootstrap assumptions, and the already proven estimates 
\eqref{E:POINTWISEAPRIORIINTERIORLBARRMINUS} and \eqref{E:POINTWISEAPRIORIINTERIORLRMINUS}
to deduce:
\begin{align} \label{E:INTERIORREGIONWAVEEQUATIONRTIMELBARLRMINUSPOINTWISE}
	\left| r \uLunit \Lunit \LRiemann \right|
	& = \left|\uLunit \Big( r \Lunit \LRiemann +  2 \LRiemann \Big) \right|
		+ \frac{\mathcal{O}(\datasize)}{(1 + t + r)^2}
	= \frac{\mathcal{O}(\datasize)}{(1 + t + r)^2},
\end{align}
as is desired.

\eqref{E:POINTWISEAPRIORIINTERIORTIMESLUNITLBARRPLUS} follows from a similar argument based on 
\eqref{E:INTERIORREGIONWAVEEQUATIONINHOMOGENEOUSTERMPOINTWISEESTIMATEFORRPLUS} and the already proven estimates
\eqref{E:POINTWISEAPRIORIINTERIORLBARRPLUS} and \eqref{E:POINTWISEAPRIORIINTERIORLRPLUS}.

\medskip
\noindent \underline{\textbf{Proof of \eqref{E:POINTWISEAPRIORIINTERIORRTIMESLUNITLUNITRMINUS}}}.
This is the only estimate for $\LRiemann$ in the proposition that has not yet been proved. 
This estimate follows easily from the identity \eqref{E:PARTIALTINTERMSOFLANDLBAR} and the already proven estimates in the proposition.

\medskip \noindent \underline{\textbf{Proof of \eqref{E:POINTWISEAPRIORIINTERIORPARTIALTRPLUS}}}.
The identity \eqref{E:PARTIALTINTERMSOFLANDLBAR} and the bootstrap assumptions imply that:
\begin{align} \label{E:PARTIALTINTERMSOFLANDLBARWITHOOF1COEFFICIENTS}
\partial_t
& = 	\mathcal{O}(1)
			\Lunit
			+
			\mathcal{O}(1)
			\uLunit.
\end{align}	
From \eqref{E:PARTIALTINTERMSOFLANDLBARWITHOOF1COEFFICIENTS} and the already proven bounds
\eqref{E:POINTWISEAPRIORIINTERIORLBARRPLUS} and \eqref{E:POINTWISEAPRIORIINTERIORLRPLUS}, we conclude
\eqref{E:POINTWISEAPRIORIINTERIORPARTIALTRPLUS}.

\medskip
\noindent \underline{\textbf{Proof of \eqref{E:POINTWISEAPRIORIINTERIORRTIMESLBARPARTIALTRPLUS}}}.
We will prove the following pointwise estimate:
\begin{align} \label{E:COMPLETENTERIORRTIMESULUNITPARTIALTRIEMANNMINUS2PARTIALRIEMANNPOINTWISE}
	\begin{split}
	r \uLunit \partial_t \RRiemann - 2 \partial_t \RRiemann
	& =
		\frac{\mathcal{O}(\datasize)}{(1 + |t-r|)^3}
		+
		\frac{\mathcal{O}(\eps^2 \datasize)}{(1 + |t-r|)^2}
		+
		\frac{\mathcal{O}(\eps^2) \ln(t-r)}{(1 + |t-r|)^3}.
\end{split}
\end{align}
Given \eqref{E:COMPLETENTERIORRTIMESULUNITPARTIALTRIEMANNMINUS2PARTIALRIEMANNPOINTWISE}, the
desired bound \eqref{E:POINTWISEAPRIORIINTERIORRTIMESLBARPARTIALTRPLUS}
follows from \eqref{E:COMPLETENTERIORRTIMESULUNITPARTIALTRIEMANNMINUS2PARTIALRIEMANNPOINTWISE} and the already proven bound
\eqref{E:POINTWISEAPRIORIINTERIORPARTIALTRPLUS}.

To initiate the proof of \eqref{E:COMPLETENTERIORRTIMESULUNITPARTIALTRIEMANNMINUS2PARTIALRIEMANNPOINTWISE}, we commute the wave equation \eqref{E:RIGHTRIEMANNEQUATIONINTERIORREGION} with $\partial_t$
and use \eqref{E:NULLVECTORFIELDS}, 
\eqref{E:PARTIALRINTERMSOFLUNITANDULUNIT}--\eqref{E:PARTIALTINTERMSOFLANDLBAR},
\eqref{E:SPEEDPOINTWISEINTERIOR},
and the bootstrap assumptions,
carefully noting that the products $- 2 \RRiemann \partial_t \Lunit \Speed$
(which are ``dangerous'' in that they involve a second derivative that lacks an $r$-weight)
cancel from both sides of the resulting equation,
arises from the differentiating the commutator term $r[\uLunit,\Lunit]\RRiemann$ on RHS~\eqref{E:RIGHTRIEMANNEQUATIONINTERIORREGION}
and noting that a null-condition-failing product of type $r \uLunit \partial_t \RRiemann \cdot \uLunit \RRiemann$
\begin{align} \label{E:INTERIORREGIONRPLUSWAVEEQUATIONCOMMUTEDWITHPARTIALTINHOMOGENEOUSTERMS}
	\left|\Lunit \Big( r \uLunit \partial_t \RRiemann - 2 \Speed \partial_t \RRiemann\Big) \right|
	& 
	\lesssim \frac{\eps^2}{(1 + t + r)(1 + |t-r|)^3}.
\end{align}

For future use, we note that the boundary condition \eqref{E:RIEMANNINVARIANTMATCHINGCONDITION},
\eqref{E:POINTWISEBOUNDFORPARTIALTRMINUSONTIMEAXIS}, and \eqref{E:SPEEDPOINTWISEINTERIOR} imply that:
\begin{align} \label{E:PARTIALTRPLUSBOUNDONTIMEAXIS}
		|\partial_t \RRiemann \restriction_{r=0}|
		& \lesssim 
			\frac{\mathcal{O}(\datasize)}{(1 + t)^3}
		+
		\frac{\mathcal{O}(\eps^2 \datasize)}{(1 + t)^2}.
\end{align}

The remainder of the proof of \eqref{E:COMPLETENTERIORRTIMESULUNITPARTIALTRIEMANNMINUS2PARTIALRIEMANNPOINTWISE} splits into two cases.

\noindent \textit{Case $I$: Estimates along integral curves of $\Lunit$ that emanate from the time axis}.
Let $(t_0,0)$ be a point on the time axis with $t_0 \geq \tstar$, 
and let $(t_1,r_1)$ be a point lying on the $t$-parameterized integral curve $t \rightarrow \upgamma_{t_0}(t)$ of $\Lunit$ emanating from
$(t_0,0)$ with $t_1 \geq t_0$.

Starting from \eqref{E:INTERIORREGIONRPLUSWAVEEQUATIONCOMMUTEDWITHPARTIALTINHOMOGENEOUSTERMS},
we argue as in the proof of \eqref{E:INTERIORRTIMESULUNITRIEMANNMINUS2RIEMANNPOINTWISEMANATINGFROMTIMEAXIS}
with the help of \eqref{E:INTEGRALESTIMATE1OVER1PLUSTPLUSRTIMES1OVERTMINUSRTPOWERBALONGINTEGRALCURVESOFL}
and \eqref{E:PARTIALTRPLUSBOUNDONTIMEAXIS}, thereby deducing that:
\begin{align} \label{E:INTERIORRTIMESULUNITPARTIALTRIEMANNMINUS2PARTIALRIEMANNPOINTWISEMANATINGFROMTIMEAXIS}
	\begin{split}
	r \uLunit \partial_t \RRiemann(t_1,r_1) - 2 \partial_t \RRiemann(t_1,r_1)
	& = - 2 \partial_t \RRiemann(t_0,0)
		+ 
		\frac{\mathcal{O}(\eps^2) \ln\left( \frac{1+t_1}{1+t_0} \right)}{(1 + |t_1-r_1|)^3}
			\\
		& =
		\frac{\mathcal{O}(\datasize)}{(1 + |t_1-r_1|)^3}
		+
		\frac{\mathcal{O}(\eps^2 \datasize)}{(1 + |t_1-r_1|)^2}
		+
		\frac{\mathcal{O}(\eps^2) \ln(t_1-r_1)}{(1 + |t_1-r_1|)^3},
\end{split}
\end{align}
which implies \eqref{E:COMPLETENTERIORRTIMESULUNITPARTIALTRIEMANNMINUS2PARTIALRIEMANNPOINTWISE} in Case $I$.

\noindent \textit{Case $II$: Estimates along integral curves of $\Lunit$ that emanate from $\Sigma_{\tstar} \cap \lbrace r \leq \frac{\tstar}{2} \rbrace$}.
Let $(\tstar,r_0) \in \Sigma_{\tstar} \cap \lbrace r \leq \frac{\tstar}{2} \rbrace$, 
and let $(t_1,r_1)$ be a point lying on the $t$-parameterized integral curve $t \rightarrow \upgamma_{r_0}(t)$ of $\Lunit$ emanating from
$(\tstar,r_0)$ with $t_1 \geq t_0$.

Staring from \eqref{E:INTERIORREGIONRPLUSWAVEEQUATIONCOMMUTEDWITHPARTIALTINHOMOGENEOUSTERMS},
we argue as in the proof of \eqref{E:INTERIORREGIONPOINTWISEESTIMATEFORRTIMESULUNITRPLUSMINUSTWORPLUS}
with the help of \eqref{E:INTEGRALESTIMATE1OVER1PLUSTPLUSRTIMES1OVERTMINUSRTPOWERBALONGINTEGRALCURVESOFL}
and the data estimate \eqref{E:TRCOORDINATESPARTIALTDERIVATIVETRANSPORTEDMODIFIEDULUNITRPLUSTSTARDATA},
thereby deducing that:
\begin{align} \label{E:INTERIORRTIMESULUNITPARTIALTRIEMANNMINUS2PARTIALRIEMANNPOINTWISEMANATINGFROMSIGMATSART}
\begin{split}
	[r \uLunit \partial_t \RRiemann - 2 \partial_t \RRiemann](t_1,r_1)
	& =
		\frac{\mathcal{O}(\datasize)}{(1 + \tstar + r_0)^3}
		+
		\frac{\mathcal{O}(\eps^2) \ln\left( \frac{1+t_1}{1+t_0} \right)}{(1 + |t-r|)^3}
		+
		\frac{\mathcal{O}(\eps^2)}{(1 + t_1 + r_1)^2(1 + |t_1 - r_1|)}
			\\
	& = 
		\frac{\mathcal{O}(\datasize)}{(1 + |t_1-r_1|)^3}
		+
		\frac{\mathcal{O}(\eps^2) \ln(t_1-r_1)}{(1 + |t_1-r_1|)^3}
		+
		\frac{\mathcal{O}(\eps^2 \datasize)}{(1 + |t_1-r_1|)^2}.
\end{split}
\end{align}
We have therefore proved \eqref{E:COMPLETENTERIORRTIMESULUNITPARTIALTRIEMANNMINUS2PARTIALRIEMANNPOINTWISE} in Case $II$,
which finishes the proof of \eqref{E:COMPLETENTERIORRTIMESULUNITPARTIALTRIEMANNMINUS2PARTIALRIEMANNPOINTWISE}.

\medskip
\noindent \underline{\textbf{Proof of \eqref{E:POINTWISEAPRIORIINTERIORTIMESLBARLBARRPLUS}}}.
This is the only estimate for $\RRiemann$ in the proposition that has not yet been proved. 
This estimate follows easily from the identity \eqref{E:PARTIALTINTERMSOFLANDLBAR} and the already proven estimates in the proposition.

\medskip \noindent \underline{\textbf{Proof of \eqref{E:INTERIORREGIONTMINUSRLARGEIMPROVEDESTIMATE} and \eqref{E:INTERIORREGIONIMPROVEDMAXIMUMBOUNDONT}}}.
Let $(t_1,r_1)$ be any point in $\mathbf{GH}_{\textnormal{Boot}}^{{\textnormal{INT}}}$. Since 
$\mathbf{GH}_{\textnormal{Boot}}^{{\textnormal{INT}}}$ is a globally hyperbolic interior bootstrap region by assumption,
the past-directed integral curve of $\uLunit$ through $(t_1,r_1)$ 
intersects a point $(t_0,r_0) \in \left(\Sigma_{\tstar} \cap \lbrace r \leq \tstar/2 \rbrace  \right)
\bigcup 
\left(
\lbrace u = \frac{\tstar}{2} 	\rbrace \cap 
\lbrace \tstar \leq t \leq \TCH\left(\frac{\tstar}{2} \right) \rbrace \right)$. 
With the help of \eqref{E:COMPARISONBETWEENUANDMINUSRONCREASECHARACTERISTIC}, we see that
$t_0 - r_0 \geq \frac{\tstar}{2} - \mathcal{O}(1)$.
From this bound and \eqref{E:INTERIORREGIONPOINTWISEESTIMATECHANGEINTMINUSRALONGINTEGRALCURVESOFLBAR},
we find that $t_1 - r_1 \geq t_0 - r_0$, which yields \eqref{E:INTERIORREGIONTMINUSRLARGEIMPROVEDESTIMATE}.

To prove \eqref{E:INTERIORREGIONIMPROVEDMAXIMUMBOUNDONT}, we first note that we have shown that all the $\eps$-involving bootstrap assumptions 
from Sect.\,\ref{SSS:INTERIORBOOTSTRAPASSUMPTIONS}
hold with $\eps$ replaced by $C \datasize$. We can therefore repeat the proof 
\eqref{E:INTERIORREGIONPOINTWISEESTIMATECHANGEINTMINUSRALONGINTEGRALCURVESOFLBAR}
to show that the estimate holds with $C \datasize$ in place of $\eps$.
Also using \eqref{E:LOGOFTIMERATIOSBOUNDEDBYCOVERDATASIZEININTERIORREGION}, 
we see that	
$t_1 - r_1
= t_0 - r_0
	+ 
	2(t_1 - t_0)
+ \mathcal{O}(1)$,
which is equivalent to:
\begin{align} \label{E:MAXSIZEOFTALONGINTGOINGINTEGRALCURVEOFULUNITINTERMSOFRADIALDIFFERENCES}
	t_1 
	& = t_0 + r_0 - r_1 + \mathcal{O}(1).
\end{align}
The maximum possible value of
$t_0 + r_0$ on RHS~\eqref{E:MAXSIZEOFTALONGINTGOINGINTEGRALCURVEOFULUNITINTERMSOFRADIALDIFFERENCES} occurs at the top of the 
characteristic portion $\lbrace (t,r) \ | \ u(t,r) = \frac{\tstar}{2} \rbrace$, and thus, 
by \eqref{E:COMPARISONBETWEENUANDMINUSRONCREASECHARACTERISTIC},
we have that $t_0 + r_0  \leq 2 \TCH\left(\frac{\tstar}{2} \right) - \frac{\tstar}{2} + \mathcal{O}(1)$. 
In view of \eqref{E:SHARPUPPERANDLOWERBOUNDSONVALUEOFTATTOPPOINTOFUEQUALSTSTAROVER2}, we further deduce that the maximum possible value
of $t_0 + r_0$ is bounded by $2 \Tcrease + \mathcal{O}(1)$. From this bound and
\eqref{E:MAXSIZEOFTALONGINTGOINGINTEGRALCURVEOFULUNITINTERMSOFRADIALDIFFERENCES}, we conclude \eqref{E:INTERIORREGIONIMPROVEDMAXIMUMBOUNDONT}.

\end{proof}

\subsection{Proof of Theorem~\ref{T:MAINMGHDEXISTENCETHEOREM}}
\label{SS:PROOFOFT:MAINEXISTENCETHEOREM}
\ \\

\noindent \underline{\textbf{Opening Remarks}}:
It is a standard result that a unique classical solution exists in the ``small-time'' region in between
$\Sigma_0$ and $\Sigma_{\tstar}$, which, when $\datasize$ is small, is 
a tiny portion of the classical lifespan, during which the nonlinearities are negligible. We refer readers to 
Sect.\,\ref{S:GLOBALEXISTENCE} for detailed proofs of how to control the solution for $0 \leq t \leq \tstar$,
where the data in that section has the opposite sign, but the same strategy for controlling the solution in the small-time region
could be applied; we omit the details here.

Next, we recall that in Prop.\,\ref{P:EXISTENCEUPTOCREASEANDSINGULARBOUNDARYANDPORTIONOFCAUCHYHORIZON}, we constructed the solution 
in the Exterior Region and in particular proved the blowup of $|\partial_r \RRiemann|$ on 
$\Upsilon(\futuresinghyp)$  and of $|\partial_r \LRiemann|$ on $ \Upsilon(\pastsinghyp)$; see Fig.\,\ref{F:MGHD}. 
Therefore, to prove Theorem~\ref{T:MAINMGHDEXISTENCETHEOREM}, 
the main remaining tasks we have to accomplish are to construct the solution in the Interior Region and to prove
uniqueness of the solution and the maximality and uniqueness of the MGHD.

\begin{figure}  
\centering
\begin{overpic}[scale=.7, grid = false, tics=3, trim=-.5cm -1cm -1cm -.5cm, clip]{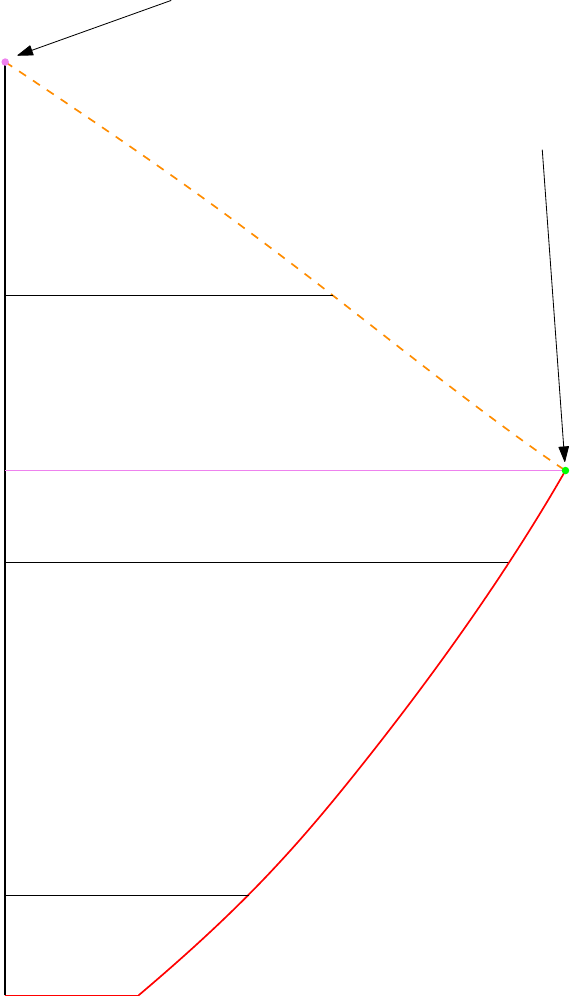}
	\put (6,7) {\tiny $\Sigma_{\tstar} \cap \lbrace  0 \leq r \leq \mathfrak{r}(\tstar) \rbrace$}
	\put (45,35) {\small $\lbrace u = \frac{\tstar}{2} \rbrace$}
	\put (53,21) {\small Exterior Region (from Prop.\,\ref{P:EXISTENCEUPTOCREASEANDSINGULARBOUNDARYANDPORTIONOFCAUCHYHORIZON})}
	\put (22,31) {\small $I$}
	\put (27,47) {\small $\Sigma_T$}
	\put (11,16) {\small $\Sigma_{T_{\textsf{Local}}}$}
	\put (14,56.5) {\tiny $\Sigma_{\TCH({\frac{\tstar}{2}})} \cap \lbrace r \leq \mathfrak{R}(\TCH({\frac{\tstar}{2}})) \rbrace$}
	\put (15,76) {\small $II$}
	\put (15,71.5) {\small $\Sigma_{T'}$}
	\put (39,66.5) {\tiny $\underline{\upgamma}_{[\TCH\left(\frac{\tstar}{2} \right),T']}$}
	\put (39,86) {\tiny $(\TCH({\frac{\tstar}{2}}),\mathfrak{R}(\TCH({\frac{\tstar}{2}})))$}
	\put (15,99) {\small $(T_{\textsf{Max}},0)$}
	\put(5,3) 		{\vector(1,0){10}}
	\put(7.5,3.5)     {\small $r$}
	\put(2,9) 		{\vector(0,1){10}}
	\put(0,11)     {\small $t$}
\end{overpic} 
\caption{Existence in the Interior Region via a priori estimates and continuation}
\label{F:EXISTENCEANDCONTINUATION}                                                                              
\end{figure}
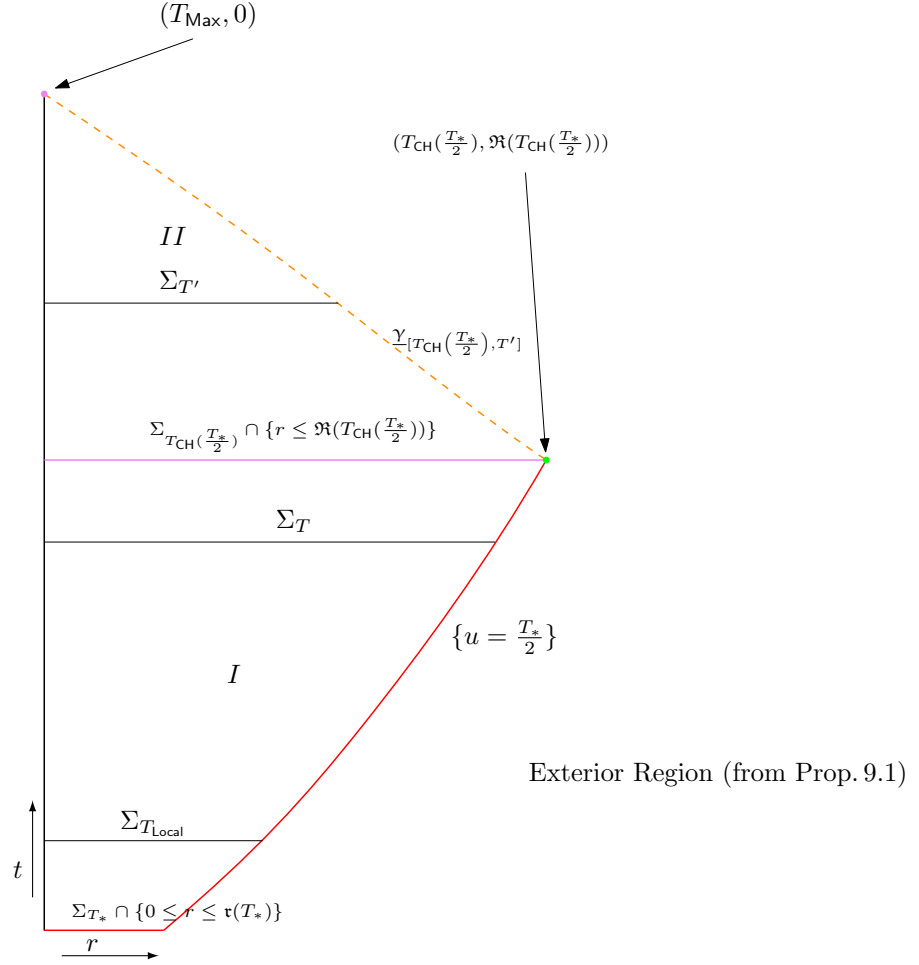

Throughout this proof, $\TCH$ denotes the function yielded by Prop.\,\ref{P:EXISTENCEUPTOCREASEANDSINGULARBOUNDARYANDPORTIONOFCAUCHYHORIZON}.

As a final opening remark, we remind the reader that for almost all results in the theorem, 
it suffices to understand the solution for times $t \geq 0$, 
i.e., to understand the solution on $\mathbf{M}^{\textsf{Future}}$.
The past-dynamics are almost a mirror image of the future-dynamics; see Remark~\ref{R:DISCRETESYMMETRY}.

\medskip

\noindent \underline{\textbf{Separation of the Interior Region into two pieces}}:
In Fig.\,\ref{F:EXISTENCEANDCONTINUATION}, we have divided the Interior Region into two subsets, $I$ and $II$,
separated by the constant-time hypersurface portion $\Sigma_{\TCH({\frac{\tstar}{2}})} \cap \lbrace r \leq \mathfrak{R}(\TCH({\frac{\tstar}{2}})) \rbrace$,
where the radial value $\mathfrak{R}(\TCH({\frac{\tstar}{2}}))$ is discussed below.
Our goal is to prove classical existence in $I \cup II$, i.e., on $\mathbf{M}^{\textsf{Future}}$.

\medskip
\noindent \underline{\textbf{Existence in region $I$}}:
We will first prove existence in region $I$ in Fig.\,\ref{F:EXISTENCEANDCONTINUATION},
which is bounded on the right by the portion of the characteristic $\lbrace (t,r) \ | \ u(t,r) = \frac{\tstar}{2} \rbrace$
that lies between $\Sigma_{\tstar}$ and the Cauchy horizon.
To proceed, we first note that results of Prop.\,\ref{P:EXISTENCEUPTOCREASEANDSINGULARBOUNDARYANDPORTIONOFCAUCHYHORIZON} concerning the change of variables map $\Upsilon$ imply that $\lbrace (t,r) \ | \ u(t,r) = \frac{\tstar}{2} \rbrace$
can be viewed as a $t$-parameterized curve, i.e., 
$t \rightarrow (t,\mathfrak{r}(t))$, where $\mathfrak{r}$ is a $C^3$ function of $t$ on the domain $[\tstar,\TCH\left(\frac{\tstar}{2} \right)]$.
Given any $T \in [\tstar,\TCH\left(\frac{\tstar}{2} \right)]$, we define the following bootstrap region, which is a subset of $I$:
\begin{align} \label{E:TYPEAINTERIORREGIONOFEXISTENCE}
	\mathbf{GH}_{\textnormal{Boot};T}^{{\textnormal{INT}}}
	& := \lbrace (t,r) \ | \ \tstar \leq t < T \rbrace \cap \lbrace 0 \leq r \leq \mathfrak{r}(t) \rbrace.
\end{align}
By standard local well-posedness for quasilinear transport equations,
if $\datasize > 0$ is sufficiently small, then there exists a time $T_{\textsf{Local}} > \tstar$ 
such that a unique classical solution
to equations \eqref{E:OUTGOINGRIEMANNINVARIANTEVOLUTION}--\eqref{E:INGOINGRIEMANNINVARIANTEVOLUTION} 
exists on $\mathbf{GH}_{\textnormal{Boot};T_{\textsf{Local}}}^{{\textnormal{INT}}}$, which is globally hyperbolic
in the sense of Def.\,\ref{D:GHINTERIORBOOTSTRAP},
and on which the bootstrap assumptions of Sect.\,\ref{SSS:INTERIORBOOTSTRAPASSUMPTIONS} hold with 
$\eps := \datasize^{3/4}$. We clarify that the ``data surfaces'' are the aforementioned
portion of the outgoing characteristic $\lbrace (t,r) \ | \ u(t,r) = \frac{\tstar}{2} \rbrace$
in union with the portion of $\Sigma_{\tstar}$ on which $0 \leq r \leq \mathfrak{r}(\tstar)$,
i.e., it is a spacelike-characteristic initial value problem. While in general one must check that the data for such an initial value problem
satisfy compatibility conditions at the corner, the needed
compatibility conditions are satisfied because we have already constructed a $C^3$ solution in the exterior region,
which contains the corner and the characteristic portion; see Remark~\ref{R:GHDONOTHERHYPERSURFACES}.
The estimates of Prop.\,\ref{P:APRIORIESTIMATESININTERIORREGION} hold on 
$\mathbf{GH}_{\textnormal{Boot};T_{\textsf{Local}}}^{{\textnormal{INT}}}$,
and in particular, no bootstrap assumptions are saturated on $\mathbf{GH}_{\textnormal{Boot};T_{\textsf{Local}}}^{{\textnormal{INT}};I}$,
$t \leq 2 \Tcrease + C$, and
the $C_{t,r}^1$ norm of $(\RRiemann,\LRiemann)$ is uniformly bounded on 
$\mathbf{GH}_{\textnormal{Boot};T_{\textsf{Local}}}^{{\textnormal{INT}}}$. Because the a priori estimates of 
Prop.\,\ref{P:APRIORIESTIMATESININTERIORREGION} 
hold on any region $\mathbf{GH}_{\textnormal{Boot};t}^{{\textnormal{INT}}}$
of classical existence with $t \in [\tstar,\TCH\left(\frac{\tstar}{2} \right)]$,
and because the solution has already been constructed in the Exterior Region,
standard continuation criteria allow us to uniquely continue the solution 
to the closure of the region
$\mathbf{GH}_{\textnormal{Boot}; \TCH\left(\frac{\tstar}{2} \right)}^{{\textnormal{INT}}}$,
on which the estimates of Prop.\,\ref{P:APRIORIESTIMATESININTERIORREGION} hold.
This closure is precisely region $I$.

\medskip
\noindent \underline{\textbf{Existence in region $II$}}:
For any $T'$ satisfying $\TCH\left(\frac{\tstar}{2} \right) \leq T' \leq 4 \Tcrease$, let
$\underline{\upgamma}_{[\TCH\left(\frac{\tstar}{2} \right),T']}$ denote the portion of the $t$-parameterized integral curve of $\uLunit$
that emanates from the point $(\TCH\left(\frac{\tstar}{2} \right),\mathfrak{r}(\TCH\left(\frac{\tstar}{2} \right)))$
and along which $t \in [\TCH\left(\frac{\tstar}{2} \right),T']$. Let $\mathfrak{R}(t)$ denote the value of $r$ along this integral curve,
i.e., $\underline{\upgamma}_{[\TCH\left(\frac{\tstar}{2} \right),T']}$ is parameterized by $t \rightarrow (t,\mathfrak{R}(t))$, where $\mathfrak{R}(\TCH\left(\frac{\tstar}{2} \right))
= \mathfrak{r}(\TCH\left(\frac{\tstar}{2} \right))$; see Fig.\,\ref{F:EXISTENCEANDCONTINUATION}; 
$\underline{\upgamma}_{[\TCH\left(\frac{\tstar}{2} \right),T']}$ is a portion of the Cauchy horizon lying to the future of 
the portion we constructed in the Exterior Region.
On any region on which the bootstrap assumptions hold, 
by 
\eqref{E:NULLVECTORFIELDS},
\eqref{E:SPEEDPOINTWISEINTERIOR},
and the bootstrap assumptions,
we have that 
$\uLunit r = v^r - \Speed = - 1 + \frac{\mathcal{O}(\eps)}{1 + t + r}$, i.e.,
$\mathfrak{R}$ is a decreasing function.
Consider the following bootstrap region, whose future-boundary is the characteristic curve portion 
$\underline{\upgamma}_{[\TCH\left(\frac{\tstar}{2} \right),T']}$:
\begin{align} \label{E:TYPEBINTERIORREGIONOFEXISTENCE}
	\widetilde{\mathbf{GH}}_{\textnormal{Boot};T'}^{{\textnormal{INT}}}
	& := \lbrace (t,r) \ | \ \TCH\left(\frac{\tstar}{2} \right) \leq t < T' \rbrace \cap \lbrace 0 \leq r \leq \mathfrak{R}(t) \rbrace.
\end{align}
The results we derived above for region $I$ in particular show that the estimates of
Prop.\,\ref{P:APRIORIESTIMATESININTERIORREGION} hold
on the spacelike hypersurface portion $\Sigma_{\TCH\left(\frac{\tstar}{2} \right)} \cap \lbrace (t,r) \ | \ r \leq \mathfrak{R}(\TCH\left(\frac{\tstar}{2} \right)) \rbrace$,
which is a Cauchy hypersurface for regions of type \eqref{E:TYPEBINTERIORREGIONOFEXISTENCE}.
Hence, by standard local well-posedness for quasilinear transport systems with data given on a spacelike hypersurface portion,
there is a $T_{\textsf{Local}}' \in (\TCH\left(\frac{\tstar}{2} \right), 4 \Tcrease]$
such that a unique classical solution
to equations \eqref{E:OUTGOINGRIEMANNINVARIANTEVOLUTION}--\eqref{E:INGOINGRIEMANNINVARIANTEVOLUTION} 
on $\widetilde{\mathbf{GH}}_{\textnormal{Boot};T_{\textsf{Local}}'}^{{\textnormal{INT}}}$, which is globally hyperbolic
in the sense of Def.\,\ref{D:GHINTERIORBOOTSTRAP},
and on which the bootstrap assumptions of Sect.\,\ref{SSS:INTERIORBOOTSTRAPASSUMPTIONS} hold with 
$\eps := \datasize^{3/4}$, and on which $\mathfrak{R}(t) > 0$;
see Fig.\,\ref{F:EXISTENCEANDCONTINUATION}. 
The estimates of Prop.\,\ref{P:APRIORIESTIMATESININTERIORREGION} hold on 
$\widetilde{\mathbf{GH}}_{\textnormal{Boot};T_{\textsf{Local}}'}^{{\textnormal{INT}}}$,
and in particular, no bootstrap assumptions are saturated on 
$\widetilde{\mathbf{GH}}_{\textnormal{Boot};T_{\textsf{Local}}'}^{{\textnormal{INT}}}$,
$T_{\textsf{Local}}' \leq 2 \Tcrease + C$,
and the $C_{t,r}^1$ norm of $(\RRiemann,\LRiemann)$ is uniformly bounded on 
$\widetilde{\mathbf{GH}}_{\textnormal{Boot};T_{\textsf{Local}}'}^{{\textnormal{INT}}}$. 
Because the a priori estimates of Prop.\,\ref{P:APRIORIESTIMATESININTERIORREGION} 
hold on any region $\widetilde{\mathbf{GH}}_{\textnormal{Boot};T'}^{{\textnormal{INT}}}$
of classical existence with $T' \in [\TCH\left(\frac{\tstar}{2} \right), 4 \Tcrease]$,
standard continuation criteria allow us to uniquely continue the solution to 
$\widetilde{\mathbf{GH}}_{\textnormal{Boot};T_{\textsf{Max}}}^{{\textnormal{INT}}}$
on which the estimates of Prop.\,\ref{P:APRIORIESTIMATESININTERIORREGION} hold,
where $T_{\textsf{Max}}$ is the minimum of the two numbers $4 \Tcrease$ and
the smallest value of $t$ for which $\mathfrak{R}(t) = 0$
(i.e., the value of $t$ at which the integral curve of $\uLunit$ emanating from 
$(\TCH\left(\frac{\tstar}{2} \right),\mathfrak{R}(\TCH\left(\frac{\tstar}{2} \right)))$ intersects the time axis).
By \eqref{E:INTERIORREGIONIMPROVEDMAXIMUMBOUNDONT}, we have
$
T_{\textsf{Max}}
< 2 \Tcrease + \mathcal{O}(1)$,
which implies \eqref{E:MAINTHEOREMESTIMATEFORTMAXININTERIORREGION}
as well as the relation $\mathfrak{R}(T_{\textsf{Max}}) = 0$,
i.e, that $(T_{\textsf{Max}},0)$ is a point on the time axis.

\medskip
\noindent \underline{\textbf{Proof that $\mathbf{M}^{\textsf{Future}}$ is globally hyperbolic}}:
It suffices to show that $\mathbf{M}^{\textsf{Future}} \cap \lbrace t \geq \tstar \rbrace$ is globally hyperbolic with Cauchy hypersurface 
$\Sigma_{\tstar}$ in the sense of Point~5 of the theorem;
it is a standard result that the remaining ``small-time'' region in between $\Sigma_0$ and $\Sigma_{\tstar}$ is globally hyperbolic with Cauchy hypersurface 
$\Sigma_0$, which, when combined with the global hyperbolicity of $\mathbf{M}^{\textsf{Future}} \cap \lbrace t \geq \tstar \rbrace$, yields the desired global hyperbolicity of $\mathbf{M}^{\textsf{Future}}$ with Cauchy hypersurface $\Sigma_0$.

To prove the global hyperbolicity of $\mathbf{M}^{\textsf{Future}} \cap \lbrace t \geq \tstar \rbrace$, we start by recalling that
we already analyzed the exterior region in the proof of Prop.\,\ref{P:EXISTENCEUPTOCREASEANDSINGULARBOUNDARYANDPORTIONOFCAUCHYHORIZON}. Hence, 
it suffices to consider the Interior Region. Let $q$ be a point in the Interior Region that is disjoint from the Cauchy horizon and disjoint from the characteristic curve
$\lbrace u = \frac{\tstar}{2} \rbrace$, which separates the Interior Region from the Exterior Region. 

Consider the backwards integral curve $t \rightarrow \underline{\upgamma}_q(t)$ 
of $\uLunit$ emanating from $q$. Since $\uLunit t = 1$ and since 
\eqref{E:SPEEDOFSOUNDEXPANSION}--\eqref{E:SPEEDOFSOUNDERRORFUNCTIONVANISHESATORIGIN},
\eqref{E:NULLVECTORFIELDS}, and
Prop.\,\ref{P:APRIORIESTIMATESININTERIORREGION} imply that
$\uLunit r = v^r - \Speed = -1 + \frac{\mathcal{O}(\datasize)}{1 + t + r}$, it follows that
$t$ decreases and $r$ increases towards the past along $\underline{\upgamma}_q$. 
Because $r$ increases, $\underline{\upgamma}_q$ cannot intersect the time axis. By ODE uniqueness,
$\underline{\upgamma}_q$ cannot intersect the Cauchy horizon because both curves are integral curves of $\uLunit$ in regions where the solution is smooth.
If $\underline{\upgamma}_q$ intersects $\Sigma_{\tstar}$ before exiting the Interior Region, then that is the desired conclusion. Since the Interior Region is pre-compact, the only other possibility is that $\underline{\upgamma}_q$ intersects the outgoing characteristic $\lbrace u = \frac{\tstar}{2} \rbrace$, 
which separates the Interior Region from the exterior one and which is transversal to $\underline{\upgamma}_q$,
since its tangent vector $\Lunit$ everywhere satisfies 
(again, by \eqref{E:SPEEDOFSOUNDEXPANSION}--\eqref{E:SPEEDOFSOUNDERRORFUNCTIONVANISHESATORIGIN},
\eqref{E:NULLVECTORFIELDS}, and
Prop.\,\ref{P:APRIORIESTIMATESININTERIORREGION}) 
$\Lunit t = 1$ and $\Lunit r = 1 + \frac{\mathcal{O}(\datasize)}{1 + t + r}$. After $\underline{\upgamma}_q$
intersects $\lbrace u = \frac{\tstar}{2} \rbrace$ and penetrates the Exterior Region, we can apply the Exterior Region analysis from the proof of 
Prop.\,\ref{P:EXISTENCEUPTOCREASEANDSINGULARBOUNDARYANDPORTIONOFCAUCHYHORIZON} to conclude that $\underline{\upgamma}_q$ 
terminates at a point in $\Sigma_{\tstar}$.

Let now $t \rightarrow \upgamma_q(t)$ be the integral curve of $\Lunit$ emanating from $q$. As $t$ decreases, 
by ODE uniqueness, $\upgamma_q(t)$ cannot intersect  
$\lbrace u = \frac{\tstar}{2} \rbrace$, because both curves are integral curves of $\Lunit$ in regions where the solution is smooth.
If $\upgamma_q(t)$ intersects $\Sigma_{\tstar}$ before intersecting the time
axis, then that is the desired conclusion. $\upgamma_q(t)$ cannot intersect the future Cauchy horizon $\futureCauchyhor$ as $t$ decreases
because no points in $\futureCauchyhor$ belong to the causal past of $q$.
Since the Interior Region is pre-compact, the only other possibility is that $\upgamma_q$ intersects the remaining boundary component 
(in $(t,r)$ coordinates)
of the Interior Region, namely the time axis, which is the desired conclusion. 

\medskip
\noindent \underline{\textbf{Uniqueness of the classical solution on $\mathbf{M}^{\textsf{Future}}$}}:
This is a standard result, whose proof we briefly describe.
Let $(\widehat{\mathcal{R}}_+,\widehat{\mathcal{R}}_-)$ be another $C^2$ 
solution to equations \eqref{E:OUTGOINGRIEMANNINVARIANTEVOLUTION}--\eqref{E:INGOINGRIEMANNINVARIANTEVOLUTION}  
on $\mathbf{M}^{\textsf{Future}}$ with the same initial data on
$\Sigma_0$ as $(\RRiemann,\LRiemann)$. Let 
$(\RRiemann^{\Delta},\LRiemann^{\Delta}) := (\RRiemann,\LRiemann) - (\widehat{\mathcal{R}}_+,\widehat{\mathcal{R}}_-)$
be the difference of the two solutions.
It is straightforward to compute that on $\mathbf{M}^{\textsf{Future}}$, 
$(\RRiemann^{\Delta},\LRiemann^{\Delta})$ satisfies an evolution system of the form 
$\Lunit \RRiemann^{\Delta} = \cdots$
and $\uLunit \LRiemann^{\Delta} = \cdots$, where the inhomogeneous terms are 
\emph{linear} in $(\RRiemann^{\Delta},\LRiemann^{\Delta})$ with coefficients that depend on the up-to-first-order derivatives of
$(\RRiemann,\LRiemann)$ and $(\widehat{\mathcal{R}}_+,\widehat{\mathcal{R}}_-)$
and whose data on $\Sigma_0$ is vanishing. It is straightforward to adapt the analysis used in proving
Props.\,\ref{P:RIEMANNINVARIANTSAPRIORIEXTERIORREGIONESTIMATES} and Prop.\,\ref{P:APRIORIESTIMATESININTERIORREGION}
to conclude, via Gr\"{o}nwall's inequality, that $(\RRiemann^{\Delta},\LRiemann^{\Delta}) \equiv (0,0)$ on $\mathbf{M}^{\textsf{Future}}$,
as is desired. 

\medskip 
\noindent \underline{\textbf{Putting some results together}}:
Also recalling Prop.\,\ref{P:EXISTENCEUPTOCREASEANDSINGULARBOUNDARYANDPORTIONOFCAUCHYHORIZON}, 
we see that we have proved Points~1-5 in the theorem. As we described in the opening remarks, Point~6 follows from similar considerations.

\medskip \noindent \underline{\textbf{Proof that $\mathbf{M}^{\textsf{Max}}$ is globally hyperbolic, Part $A$}}:
Point~7 in the theorem follows easily from the global hyperbolicity of 
$\mathbf{M}^{\textsf{Future}}$ and $\mathbf{M}^{\textsf{Past}}$.

\medskip \noindent \underline{\textbf{Proof that $\mathbf{M}^{\textsf{Max}}$ is globally hyperbolic, Part $B$}}:
Point~8 follows from Point~7 and standard arguments in Lorentzian geometry. For the reader's convenience, 
we provide a proof in Prop.\,\ref{P:EQUIVALENCEOFGLOBALHYPERBOLICITYMMAX},
of Appendix~\ref{A:GLOBALHYPERBOLICITYANDMGHDS}, after providing some background material on Lorentzian geometry.

\medskip
\noindent \underline{\textbf{Maximality of $\mathbf{M}^{\textsf{Max}}$}}:
We now prove Point~9 stated in the theorem.
We argue by contradiction. Suppose that there exists a $\widetilde{\mathbf{M}}$ with the stated properties
but with $\mathbf{M}^{\textsf{Max}} \subsetneq \widetilde{\mathbf{M}}$. 
Since $\widetilde{\mathbf{M}}$ equipped with the acoustical metric
$\mathbf{g}$ (see \eqref{E:MAINTHEOREMACOUSTICALMETRIC}) is assumed to be globally hyperbolic
in the sense\footnote{This definition only directly refers to $\mathbf{g}$-timelike curves, 
while Lemma~\ref{L:CAUCHYHYPERSURFACEDODISENTIREMANIFOLD} has implications for any inextendible $\mathbf{g}$-causal curve.} of Def.\,\ref{D:CAUCHYHYPERSURFACESANDGLOBALLYHYPERBOLICREGIONS} with $\lbrace t=  0 \rbrace$ as a Cauchy hypersurface, 
Lemma~\ref{L:CAUCHYHYPERSURFACEDODISENTIREMANIFOLD} (see also Remark~\ref{R:CAUSALCURVESINTERSECTCAUCHYHYPEERSURFACES}) guarantees 
that all inextendible causal curves in $(\widetilde{\mathbf{M}},\mathbf{g})$ must intersect $\lbrace t = 0 \rbrace$.

Note that $\widetilde{\mathbf{M}}$ and the solution on it are not assumed to be spherically symmetric in 
$\widetilde{\mathbf{M}} \backslash \mathbf{M}^{\textsf{Max}}$. However, we will now prove that $\widetilde{\mathbf{M}}$ and the solution on it
\emph{are} spherically symmetric.
By \cite{aBmS2006}, $\widetilde{\mathbf{M}}$ is foliated by the level sets of a smooth temporal function $\uptau$
(i.e., all level sets of $\uptau$ are $\mathbf{g}$-spacelike and the gradient of $\uptau$ is $\mathbf{g}$-timelike and everywhere non-vanishing)
such that $\lbrace \uptau = 0 \rbrace = \lbrace t = 0 \rbrace$.
Since the initial data for \eqref{E:INTROTRANSPORTVI}--\eqref{E:INTROTRANSPORTDENSITY} are spherically symmetric by assumption, 
we can use a standard argument based on energy estimates\footnote{More precisely, one could commute equations 
\eqref{E:INTROTRANSPORTVI}--\eqref{E:INTROTRANSPORTDENSITY} 
with the standard Euclidean rotation vectorfields.
The initial $L^2$ energy of the non-radial components of the solution, including the rotationally-differentiated components, is vanishing in view of the spherical symmetry of the data. 
Moreover, since the equations are invariant under the Euclidean rotation symmetry group,
the commuted equations are linear in the components of the non-radial components of the solution, i.e., the components that one would like to show 
are identically vanishing. Hence, the desired vanishing of the non-radial components of the solution follows from standard $L^2$-type energy identities and
Gr\"{o}nwall's inequality with respect to the time function $\uptau$. \label{FN:RADIALSYMMETRYOFDATAISPRESERVEDINGHDS}} 
on the level sets of $\uptau$ 
to show that the solution is spherically symmetric on all of 
$\widetilde{\mathbf{M}}$. By extending the solution if necessary to be independent of angular coordinates relative to standard spherical
coordinates on $\mathbb{R}^{1+3}$, we can assume that the region $\widetilde{\mathbf{M}}$ is spherically symmetric. Hence, for the rest of the argument,
we can use $(t,r)$-coordinates and spherical Riemann invariants, and we can view 
$\mathbf{M}^{\textsf{Max}}$ and $\widetilde{\mathbf{M}}$ as subsets of $(t,r)$-coordinate space.

Let $q \in \widetilde{\mathbf{M}} \backslash \mathbf{M}^{\textsf{Max}}$. 
Without loss of generality, we will assume that $t$ is positive at $q$; the case in which 
$t$ is negative at $q$ can be handled using similar arguments.
Let $t \rightarrow \underline{\upgamma}_q(t)$ be the maximally extended, $t$-parameterized 
integral curve of $\uLunit$ in $ \widetilde{\mathbf{M}}$ that emanates from $q$. 
Since $\underline{\upgamma}_q(t)$ is a $\mathbf{g}$-causal curve, the aforementioned 
Lemma~\ref{L:CAUCHYHYPERSURFACEDODISENTIREMANIFOLD} guarantees that $t \rightarrow \underline{\upgamma}_q(t)$
\emph{must} intersect $\Sigma_0$ (at time $0$) as $t$ decreases. Recall that the future-boundary of $\mathbf{M}^{\textsf{Max}}$, 
i.e., $\futureCauchyhor \cup \Upsilon(\futuresinghyp)$,
is the graph of a continuous function and that $\mathbf{M}^{\textsf{Max}}$ lies below this graph; see Fig.\,\ref{F:MGHD}. 
Hence,
$\underline{\upgamma}_q$ must intersect $\futureCauchyhor \cup \Upsilon(\futuresinghyp)$ as $t$ decreases.
We will show that this is impossible. 
$\underline{\upgamma}_q(t)$ cannot intersect $\Upsilon(\futuresinghyp)$ as $t$ decreases because $|\partial_r \RRiemann|$ blows up as 
$\Upsilon(\futuresinghyp)$ is
approached within $\mathbf{M}^{\textsf{Max}}$, i.e., the solution cannot be classically continued to $\Upsilon(\futuresinghyp)$, 
and thus $\Upsilon(\futuresinghyp)$ cannot belong to $\widetilde{\mathbf{M}}$.
If $\underline{\upgamma}_q(t)$ intersects $\futureCauchyhor$ as $t$ decreases,
then since the solution $(\RRiemann,\LRiemann)$ uniquely extends to be 
$C^2$ on $\futureCauchyhor$,
by ODE uniqueness, this integral curve must be a parameterization of the future-Cauchy horizon $\futureCauchyhor$.
Now before hitting $\Sigma_0$, $\futureCauchyhor$ terminates (to the past) on
the future-crease $\Upsilon(\futurecrease)$, where $|\partial_r \RRiemann|$ blows up.
This contradicts the fact that $t \rightarrow \underline{\upgamma}_q(t)$ intersects $\Sigma_0$.
We have therefore shown that no point $q \in \widetilde{\mathbf{M}} \backslash \mathbf{M}^{\textsf{Max}}$ exists.

\medskip
\noindent \underline{\textbf{Uniqueness of $\mathbf{M}^{\textsf{Max}}$}}:
We now prove the final Point~10 stated in the theorem.
As we noted above, the MGHD manifold $\mathbf{M}^{\textsf{Max}}$ that we constructed has the following key structural property:
its future-boundary $\futureCauchyhor \cup \Upsilon(\futuresinghyp)$, viewed as a spherically symmetric subset of $\mathbb{R}^{1+3}$,
is the graph of a continuous function
(see Point $8$ of Prop.\,\ref{P:EXISTENCEUPTOCREASEANDSINGULARBOUNDARYANDPORTIONOFCAUCHYHORIZON}),
and similar remarks apply to its past-boundary.
The ``uniqueness'' conclusions stated in the final point of Theorem~\ref{T:MAINMGHDEXISTENCETHEOREM} 
therefore follow from \cite{fEhRjS2019}*{Theorem~4.82}, modulo the following remark: 
While strictly speaking, the results of \cite{fEhRjS2019}*{Theorem~4.82} were derived for scalar quasilinear wave equations on 
$\mathbb{R}^{1+3}$ 
of the form $(\mathbf{g}^{-1})^{\alpha \beta}(\partial \phi) \partial_{\alpha} \partial_{\beta} \phi = \mathcal{N}(\phi,\partial \phi)$,
with $\mathcal{N}$ a nonlinearity,
for $C^2$ solutions,
the irrotational and isentropic compressible Euler equations are equivalent to a quasilinear wave equation system,
where $\phi$ is an array comprising the density $\varrho$ and the Cartesian velocity components $v^1,v^2,v^3$;
see \cites{jLjS2020a,jS2019c}. The fact that the fluid wave equations from \cites{jLjS2020a,jS2019c} 
are system with a common wave operator, 
as opposed to the scalar wave equations treated in
\cite{fEhRjS2019}*{Theorem~4.82}, does not have any affect on the proof or conclusions of \cite{fEhRjS2019}*{Theorem~4.82}.
For the reader's convenience, we discuss \cite{fEhRjS2019}*{Theorem~4.82} in more detail in Sect.\,\ref{SSS:UNIQUEMGHDSWITHUNJUSITFIEDCONNECTEDNESS}.

This concludes our proof of Theorem~\ref{T:MAINMGHDEXISTENCETHEOREM}.

\section{Global existence}
\label{S:GLOBALEXISTENCE}
In this section, we prove Theorem~\ref{T:MAINGLOBALEXISTENCETHEOREM}, which yields global existence for an open set of small ``rarefactive'' initial data that
essentially has the opposite sign (see also Remark~\ref{R:CHANGINGTHESIGNFORGLOBALEXISTENCE}) 
of the ``compressive'' data from Theorem~\ref{T:MAINMGHDEXISTENCETHEOREM}.

\subsection{Data at time \texorpdfstring{$0$}{zero}}
\label{SS:GLOBALEXISTENCEDATAATTIME0}
As we described in the introduction, to prove global existence in the ``rarefactive regime,'' we will consider perturbations of the following profile,
whose sign is the opposite of the ``compressive profile'' in \eqref{E:BACKGROUNDPROFILEWITHSMALLAMPLITUDEFACTOR}:
\begin{align} \label{E:GLOBALEXISTENCEBACKGROUNDPROFILEWITHSMALLAMPLITUDEFACTOR}
	\globalprofileoverr(r)
	& := \datasize \frac{\arctan(r)}{r}.
\end{align}

Much like in our proof of Theorem~\ref{T:MAINMGHDEXISTENCETHEOREM}, to prove our main global existence results, we will assume that $\datasize > 0$ is sufficiently small, and
that $(\dataRRiemann,\dataLRiemann)$
is close to $(\globalprofileoverr,\globalprofileoverr)$ in the following sense:
\begin{subequations}
\begin{align}
		\| \dataRRiemann - \globalprofileoverr  \|_{C_*^2(\Sigma_0)}
		& \leq \datasize^{1+},
			 \label{E:GLOBALEXISTENCERPLUSDATAISCLOSETOBACKGROUNDATTIME0} \\
		\| \dataLRiemann - \globalprofileoverr  \|_{C_*^2(\Sigma_0)}
		& \leq \datasize^{1+},
		\label{E:GLOBAEXISTENCERMINUSDATAISCLOSETOBACKGROUNDATTIME0}
\end{align}
\end{subequations}
where: 
\begin{align} \label{E:GLOBALEXISTENCEPERTURBATIONTIME0NORM}
	\| \phi \|_{C_*^2(\Sigma_0)}
	& := 
		\sum_{K=0}^2
		\sup_{r \geq 0} 
		|(1 + r)^{2+K} \partial_r^K \phi(r)|.
\end{align}

\begin{remark}[Don't need the extra derivative for global existence]
\label{R:C2GLOBALEXISTENCE}
Note that \eqref{E:GLOBALEXISTENCEPERTURBATIONTIME0NORM} involves one fewer derivative compared to \eqref{E:PERTURBATIONTIME0NORM}.
The reason is that in Theorem~\ref{T:MAINMGHDEXISTENCETHEOREM}, the extra derivative was used only to derive the structure of the singular boundary 
and Cauchy horizon, which are not present for the global solutions that we study in Theorem~\ref{T:MAINGLOBALEXISTENCETHEOREM}.
\end{remark}

As before, to ensure that our solutions in $(t,r)$-coordinates correspond to smooth solutions on $\mathbb{R}^{1+3}$,
we assume that the Riemann invariant data satisfy the following boundary conditions:
\begin{align} \label{E:GLOBALEXISTENCERIEMANNINVARIANTMATCHINGCONDITION}
	\partial_r^{2k} 
	\RRiemann \restriction_{\lbrace r=0 \rbrace}
	& 
	= 
	\partial_r^{2k} \LRiemann \restriction_{\lbrace r=0 \rbrace},
	& k=0,1 \\
	\label{E:GLOBALEXISTENCERIEMANNINVARIANTRDERIVATIVEMATCHINGCONDITION}
	\partial_r \RRiemann \restriction_{\lbrace r=0 \rbrace}
	& 
	= 
	- \partial_r \LRiemann \restriction_{\lbrace r=0 \rbrace}.
	&&
\end{align}

\begin{remark}[Changing the sign for global existence]
\label{R:CHANGINGTHESIGNFORGLOBALEXISTENCE}
For many of the estimates in Sect.\,\ref{S:GLOBALEXISTENCE}, the only important difference compared to the rest of the paper is that signed terms involving
$\datasize$ have the opposite sign, thanks to our assumption \eqref{E:GLOBALEXISTENCEBACKGROUNDPROFILEWITHSMALLAMPLITUDEFACTOR}.
\end{remark}

\subsection{Statement of the main global existence theorem}
\label{SS:STATEMENTOFMAINGLOBALEXISTENCETHEOREM}
In this section, we state our main global existence theorem. We provide the proof in Sect.\,\ref{SS:PROOFOFMAINGLOBALEXISTENCETHEOREM}, 
after deriving the a priori estimates that lie at the heart of the result.

\begin{theorem}[Main global existence theorem]
\label{T:MAINGLOBALEXISTENCETHEOREM}
Consider the compressible Euler equations \eqref{E:OUTGOINGRIEMANNINVARIANTEVOLUTION}--\eqref{E:INGOINGRIEMANNINVARIANTEVOLUTION}
under any $C^3$ equation of state with positive sound speed, including the Chaplygin gas.
Without loss of generality, assume the positivity condition\footnote{If $\lifespanconstant < 0$, 
then we can just change the sign of $\datasize$, and the theorem still holds. \label{E:AGAINCHANGESIGNOFDATAFOROPPOSITESIGNEDNONLINEARITIES}} 
\eqref{E:LIFESPANCONSTANT} on the constant $\lifespanconstant$, except in the case of the Chaplygin gas, for which 
$\lifespanconstant = 0$.
Assume that the Riemann invariant data functions $(\dataRRiemann,\dataLRiemann)$ satisfy
\eqref{E:GLOBALEXISTENCERPLUSDATAISCLOSETOBACKGROUNDATTIME0}--\eqref{E:GLOBAEXISTENCERMINUSDATAISCLOSETOBACKGROUNDATTIME0}
(see Appendix~\ref{A:EXTENDRESULTSTOOTHERPROFILES} for a much larger class of data to which our results can be extended).
If $\datasize$ is sufficiently small, then the 
following results hold.

\begin{itemize}
	\item (\textbf{Global existence and uniqueness}).
The initial data launch a unique classical solution $(\RRiemann,\LRiemann)$ to the compressible Euler equations
\eqref{E:OUTGOINGRIEMANNINVARIANTEVOLUTION}--\eqref{E:INGOINGRIEMANNINVARIANTEVOLUTION} that
satisfies the boundary conditions
\eqref{E:RIEMANNINVARIANTMATCHINGCONDITION}--\eqref{E:RIEMANNINVARIANTRDERIVATIVEMATCHINGCONDITION}
and exists globally in space and time,
i.e., for $(t,r) \in (-\infty,\infty) \times [0,\infty)$.
\item (\textbf{Future-estimates}). The
eikonal function $u$ and $\upmu$ are classical solutions to \eqref{E:GLOBALEXISTENCEODEFOREXTERIOREIKONALFUNCTION} and
\eqref{E:MUEVOLUTION} respectively in the away-from-time-axis region
$(t,r) \in [0,\infty) \times [1,\infty)$.
Moreover, for $t \geq 0$,  $(\RRiemann,\LRiemann)$, $u$, and $\upmu$ 
obeys the estimates proved in
Props.\,\ref{P:GLOBALEXISTENCERIEMANNINVARIANTSAPRIORIEXTERIORREGIONESTIMATES}
and \ref{P:GLOBALEXISTENCESHARPESTIMATESFORMU}. 
\item (\textbf{Future-diffeomorphism property of the change of variables map}). 
Let $\Datafunctionforeikonal(t)$ be the function constructed in
Sect.\,\ref{SS:GLOBALEXISTENCEEIKONALFUNCTION}
(note that \eqref{E:TIMELILKEEIKONALDATAISAPERTURBATIONOFTMINUS1}--\eqref{E:C2ESTIMATEFORTIMELILKEEIKONALDATAPERTURBATION} imply that 
$\Datafunctionforeikonal(t)$ is a perturbation of $t-1$).
Then the change of variables map $\Upsilon(t,u) = (t,r)$
is a $C^2$ diffeomorphism from the closed set
$\lbrace (t,u) \ 0 \leq t < \infty, u \leq \Datafunctionforeikonal(t) \rbrace$
onto the closed image set $\lbrace r \geq 1 \rbrace \cap \lbrace 0 \leq t < \infty \rbrace$.
\item (\textbf{Similar results to the past}). 
Finally, for $t < 0$, the solution obeys analogous results, where the 
estimates corresponding to
Props.\,\ref{P:GLOBALEXISTENCERIEMANNINVARIANTSAPRIORIEXTERIORREGIONESTIMATES}
and \ref{P:GLOBALEXISTENCESHARPESTIMATESFORMU}
are obtained by making the following three changes: 
i) replacing $t$ with $-t = |t|$; ii) interchanging $\Lunit$ and $\uLunit$; iii) interchanging $\RRiemann$ and $\LRiemann$.

\end{itemize}

\end{theorem}

\subsection{Estimates at time \texorpdfstring{$0$}{zero}}
\label{SS:GLOBALEXISTENCEESTIMATESATTIME0}
In this section, we derive various estimates for the Riemann invariant variables at time $0$. 
We collect these estimates into the following lemma.
We will use them later on when deriving a priori estimates.

\begin{lemma}[Estimates for the Riemann invariants at time $0$]
\label{L:GLOBALEXISTENCERIEMANNINVARIANTESTIMATESATTIME0}
Assume that the Riemann invariant data functions $(\dataRRiemann,\dataLRiemann)$ satisfy
\eqref{E:GLOBALEXISTENCERPLUSDATAISCLOSETOBACKGROUNDATTIME0}--\eqref{E:GLOBAEXISTENCERMINUSDATAISCLOSETOBACKGROUNDATTIME0}.
If $\datasize$ is sufficiently small, then the solution
$(\RRiemann,\LRiemann)$ to
\eqref{E:OUTGOINGRIEMANNINVARIANTEVOLUTION}--\eqref{E:INGOINGRIEMANNINVARIANTEVOLUTION}
obeys the following estimates at time $0$.

\noindent \underline{\textbf{Crude estimates for the Riemann invariants}}.

\begin{subequations}
\begin{align}
	|\RRiemann(0,r)|
	& \lesssim \frac{\datasize}{1 + r},
		\label{E:GLOBALEXISTENCERPLUSPOINTWISEESTIMATEATTIME0} \\
		|\Lunit \RRiemann(0,r)|
	& \lesssim \frac{\datasize}{(1 + r)^3},
	\label{E:GLOBALEXISTENCELUNITRPLUSPOINTWISEESTIMATEATTIME0}
		\\
	|\uLunit \RRiemann(0,r)|,
		\,
	|\partial_t \RRiemann(0,r)|
	& \lesssim \frac{\datasize}{(1 + r)^2},
		\label{E:GLOBALEXISTENCEULUNITORPARTIALTRPLUSPOINTWISEESTIMATEATTIME0} 
\end{align}
\end{subequations}

\begin{subequations}
\begin{align}
	|\LRiemann(0,r)|
	& \lesssim \frac{\datasize}{1 + r},
		\label{E:GLOBALEXISTENCERMINUSPOINTWISEESTIMATEATTIME0} 
			\\
	|\Lunit \LRiemann(0,r)|,
		\,
	|\partial_t \LRiemann(0,r)|
	& \lesssim \frac{\datasize}{(1 + r)^2},
	\label{E:GLOBALEXISTENCELUNITORPARTIALTRMINUSPOINTWISEESTIMATEATTIME0}
		\\
	|\uLunit \LRiemann(0,r)|
	& \lesssim \frac{\datasize}{(1 + r)^3}.
		\label{E:GLOBALEXISTENCEULUNITRMINUSPOINTWISEESTIMATEATTIME0} 
\end{align}
\end{subequations}

\medskip

\noindent \underline{\textbf{Sharp estimates for special combinations involving the Riemann invariants}}.
For $0 \leq J + K \leq 1$, we have:
\begin{align} \begin{split} \label{E:GLOBALEXISTENCETIME0TRCOORDINATESALLDERIVATIVESTRANSPORTEDMODIFIEDULUNITRPLUSDATA}
	\Lunit^J
	\uLunit^K
	\left\lbrace
			r
			\uLunit
			\RRiemann(0,r)
			- 2 \RRiemann(0,r)
		\right\rbrace
		&	= 
			- 2 \datasize 
				\left\lbrace
				\left(\partial_t + \partial_r \right)^J
				\left( \partial_t - \partial_r \right)^K
				\frac{1}{1 + (t-r)^2}
				\right\rbrace \restriction_{t = 0}
				\\
		& \ \
			+
			\mathcal{O}(\datasize^{1^+}) \frac{1}{(1 + r)^{2+J+K}},
\end{split}
	\\
\begin{split}  \label{E:GLOBALEXISTENCETIME0TRCOORDINATESPARTIALTDERIVATIVETRANSPORTEDMODIFIEDULUNITRPLUSDATA}
		\left\lbrace
			r
			\uLunit
			\partial_t \RRiemann(0,r)
			- 2 \partial_t \RRiemann(0,r)
		\right\rbrace
		&	= 
			- 2 \datasize 
				\left\lbrace
				\partial_t
				\left( \frac{1}{1 + (t-r)^2} \right)
				\right\rbrace \restriction_{t= 0}
				\\
		& \ \
			+
			\mathcal{O}(\datasize^{1^+}) \frac{1}{(1 + r)^3}.
\end{split}
\end{align}

For $J=0,1$, we have:
\begin{align} 
\begin{split} \label{E:GLOBALEXISTENCETIME0TRCOORDINATESPURELDERIVATIVESTRANSPORTEDMODIFIEDLUNITRMINUSDATA}
	\Lunit^J
	\left\lbrace
			r
			\Lunit \LRiemann(0,r)
			+ 
			2 \LRiemann(0,r)
		\right\rbrace 
		& = 
				2 [\datasize + \mathcal{O}(\datasize^{1^+})]
				\left\lbrace
				\left(\partial_t + \partial_r \right)^J
				\frac{1}{1 + (t + r)^2} 
				\right\rbrace
				\restriction_{t=0}.
\end{split}
\end{align}

For $0 \leq J + K' \leq 1$, we have:
\begin{align} 
\begin{split} \label{E:GLOBALEXISTENCETIME0TRCOORDINATESATLEASTONELBARDERIVATIVETRANSPORTEDMODIFIEDLUNITRMINUSDATA}
	\Lunit^J
	\uLunit^{K'}
	\uLunit
	\left\lbrace
			r
			\Lunit \LRiemann(0,r)
			+ 
			2 \LRiemann(0,r)
		\right\rbrace 
		& = 
			\mathcal{O}(\datasize^{1^+}) \frac{1}{(1 + r)^{3+J+K'}}.
\end{split}
\end{align}

Finally, we have:
\begin{align}
\begin{split} \label{E:GLOBALEXISTENCETIME0TRCOORDINATESPARTIALDDERIVATIVETRANSPORTEDMODIFIEDLUNITRMINUSDATA}
	\left\lbrace
			r
			\Lunit \partial_t \LRiemann(0,r)
			+ 
			2 \partial_t \LRiemann(0,r)
		\right\rbrace 
		& = 
				2 
				[\datasize + \mathcal{O}(\datasize^{1^+})]
				\left\lbrace
					\partial_t
				 \left( \frac{1}{1 + (t + r)^2} \right)
				\right\rbrace
				\restriction|_{t=0}
					\\
		& \ \
		+ 
		\mathcal{O}(\datasize^{1^+}) \frac{1}{(1 + r)^3}.
\end{split}
\end{align}

\end{lemma}

\begin{proof}[Discussion of the proof]
	The lemma can be proved via the same arguments that we used in the proof of Lemma~\ref{L:RIEMANNINVARIANTESTIMATESATTIME0},
	taking into account Remark~\ref{R:CHANGINGTHESIGNFORGLOBALEXISTENCE}.
\end{proof}

\subsection{Eikonal function construction}
\label{SS:GLOBALEXISTENCEEIKONALFUNCTION}
As before, to see the ``good null structure'' in the evolution equations,
we work relative to geometric coordinates $(t,u)$, where the eikonal function $u$ once again
solves the following transport equation:
\begin{align} \label{E:GLOBALEXISTENCEODEFOREXTERIOREIKONALFUNCTION}
\Lunit u 
	& = 0.
\end{align}
However, we now prescribe data for $u$ so that it will be defined for $r \geq 1$. To achieve this, we
prescribe data on the manifold $\left(\Sigma_0 \cap \lbrace r \geq 1 \rbrace \right) \cup \lbrace r = 1 \rbrace$,
which has a corner at $(t,r) = (0,1)$. Specifically, we set:
\begin{subequations}
\begin{align}
u \restriction_{\Sigma_0 \cap \lbrace r \geq 1 \rbrace} 
	& 
	= - r, \label{E:GLOBALEXISTENCESPACELIKEEIKONALFUNCTIONDATA}
		\\
	u \restriction_{\lbrace r = 1 \rbrace }
	& = \Datafunctionforeikonal(t),
	\label{E:GLOBALEXISTENCETIMELKEEIKONALFUNCTIONDATA}
\end{align}
\end{subequations}
where in \eqref{E:TIMELIKEDATAFORGLOBALEXISTENCEEIKONALFUNCTION}, 
we construct the function $\Datafunctionforeikonal(t)$ so that 
$u$ is $C^2$ at the corner. 
To achieve this $C^2$ corner regularity, we must ensure that the data for $u$ satisfies corner compatibility conditions,
i.e., that the time derivatives of $\Datafunctionforeikonal$ up to second order at the corner match
the corresponding time derivatives determined by the eikonal equation and the spatial derivative data
in \eqref{E:GLOBALEXISTENCESPACELIKEEIKONALFUNCTIONDATA}. 
Hence, to proceed, we use the
eikonal equation \eqref{E:GLOBALEXISTENCEODEFOREXTERIOREIKONALFUNCTION},
equation \eqref{E:NULLVECTORFIELDS} for $\Lunit$,
and the spacelike data-assumption \eqref{E:GLOBALEXISTENCESPACELIKEEIKONALFUNCTIONDATA} to compute the time Taylor coefficients of
$u$ at the corner $(0,1)$:
\begin{subequations}
\begin{align}
	u(0,1) & = -1,
		\label{E:0THTAYLORCOEFFICIENTFORGLOBALEXISTENCEEIKONALFUNCTION} 
		\\
	\partial_t u(0,1)
	& = [v^r + \Speed](0,1)
		:= A_1,
		\label{E:FIRSTTAYLORCOEFFICIENTFORGLOBALEXISTENCEEIKONALFUNCTION} \\
 \partial_t^2 u(0,1)
	& = [\partial_t v^r + \partial_t \Speed](0,1)
		- [(v^r + \Speed) \partial_r (v^r + \Speed)](0,1)
	:= A_2.
	\label{E:SECONDTAYLORCOEFFICIENTFORGLOBALEXISTENCEEIKONALFUNCTION}
\end{align}
\end{subequations}
We then let $P(t)$ be the corresponding second-order Taylor polynomial:
\begin{align} \label{E:SECONDORDERTAYLORPOLYFOREIKONALCORNERCOMPATIBILITY}
	P(t) & := 
		-
		1
		+
		A_1 t
		+
		\frac{A_2}{2} t^2.
\end{align}
Using the data-assumptions 
\eqref{E:GLOBALEXISTENCERPLUSDATAISCLOSETOBACKGROUNDATTIME0}--\eqref{E:GLOBAEXISTENCERMINUSDATAISCLOSETOBACKGROUNDATTIME0},
\eqref{E:SPEEDOFSOUNDEXPANSION}--\eqref{E:SPEEDOFSOUNDERRORFUNCTIONVANISHESATORIGIN},
\eqref{E:PARTIALRINTERMSOFLUNITANDULUNIT},
and \eqref{E:PARTIALTINTERMSOFLANDLBAR},
it is straightforward to compute that when $\datasize$ is sufficiently small, we have:
\begin{subequations}
\begin{align}
	A_1 & = 1 + \mathcal{O}(\datasize),
	\label{E:POINTWISEESTIMATEFORFIRSTTAYLORCOEFFICIENTFORGLOBALEXISTENCEEIKONALFUNCTION}	\\
	A_2 & = \mathcal{O}(\datasize).
	\label{E:POINTWISEESTIMATEFORSECONDTAYLORCOEFFICIENTFORGLOBALEXISTENCEEIKONALFUNCTION}	
\end{align}
\end{subequations}
We then fix a $C^{\infty}$, non-negative cut-off function $\phi$ satisfying:
\begin{align} \label{E:GLOBALEXISTENCEEIKONALFUNCTIONCUTOFFFUNCTIONFORTIMELIKEDATA}
	\phi(t)
	& = \begin{cases}	
				1, & \mbox{if } 0 \leq t \leq 1,
					\\
				0 & \mbox{if } t \geq 2.
			\end{cases}
\end{align}
To construct $\Datafunctionforeikonal(t)$, we simply use $\phi(t)$ to glue $P(t)$ to the function $t-1$:
\begin{align} \label{E:TIMELIKEDATAFORGLOBALEXISTENCEEIKONALFUNCTION}
	\Datafunctionforeikonal(t)
	& := \phi(t) P(t) 
		+ 
		(1 - \phi(t))(t - 1).
\end{align}
By construction, 
$\Datafunctionforeikonal(0)$ is equal to RHS~\eqref{E:0THTAYLORCOEFFICIENTFORGLOBALEXISTENCEEIKONALFUNCTION},
$\partial_t \Datafunctionforeikonal(0)$ is equal to RHS~\eqref{E:FIRSTTAYLORCOEFFICIENTFORGLOBALEXISTENCEEIKONALFUNCTION},
and
$\partial_t^2 \Datafunctionforeikonal(0)$
is equal to RHS~\eqref{E:SECONDTAYLORCOEFFICIENTFORGLOBALEXISTENCEEIKONALFUNCTION},
i.e., the data satisfy the desired $C^2$ corner compatibility conditions.

We now derive some basic estimates for $\Datafunctionforeikonal$.
Using \eqref{E:SECONDORDERTAYLORPOLYFOREIKONALCORNERCOMPATIBILITY}--\eqref{E:TIMELIKEDATAFORGLOBALEXISTENCEEIKONALFUNCTION},
we compute that:
\begin{align} \label{E:TIMELILKEEIKONALDATAISAPERTURBATIONOFTMINUS1}
	\Datafunctionforeikonal(t)
	& = t - 1
		+
		g(t),
\end{align}
where the function $g$ is supported in $[1,2]$ and satisfies:
\begin{align} \label{E:C2ESTIMATEFORTIMELILKEEIKONALDATAPERTURBATION}
	\| g \|_{C^2([1,2])}
	& \lesssim \datasize.
\end{align}
In particular, on $\lbrace r = 1 \rbrace \cap \lbrace t \geq 2 \rbrace$, we have:
\begin{align} \label{E:UISTMINUS1FORTBIGGERTHAN2}
	u(t,1) & = \Datafunctionforeikonal(t) = t - 1 = t-r.
\end{align}
Note also that for any number $\Tboot > 0$, the following two subsets of $(t,r)$-space are equal:
\begin{align} \label{E:RBIGGERTHANONEISUBIGGERTHANUOFT}
	[0,\Tboot) \times [1,\infty)
	& = \lbrace (t,r) \ | \ 0 \leq t < \Tboot, \, u(t,r) \leq \Datafunctionforeikonal(t) \rbrace.
\end{align}
We will refer to the subsets in \eqref{E:RBIGGERTHANONEISUBIGGERTHANUOFT} as the ``away-from-time-axis'' region, i.e., away from $r=0$, where
the equations exhibit a coordinate degeneracy.

Finally, for future use, we use 
\eqref{E:NULLVECTORFIELDS},
\eqref{E:PARTIALRINTERMSOFLUNITANDULUNIT}--\eqref{E:PARTIALTINTERMSOFLANDLBAR}, 
\eqref{E:DEFOFINVERSEFOLIATIONDENSITY},
and
\eqref{E:GLOBALEXISTENCEODEFOREXTERIOREIKONALFUNCTION}--\eqref{E:GLOBALEXISTENCETIMELKEEIKONALFUNCTIONDATA}
to derive the following identities for the data of $\upmu$:
\begin{subequations}
\begin{align} \label{E:IDENTITYFORMUATTIMEZERO}	
	\upmu \restriction_{\Sigma_0 \cap \lbrace r \geq 1 \rbrace}
	& = \frac{1}{\Speed} \restriction_{\Sigma_0 \cap \lbrace r \geq 1 \rbrace} 
			\\
	\label{E:IDENTITYFORMUALONGREQUALS1}	
	\upmu \restriction_{\lbrace r = 1 \rbrace}
	& = \frac{1 + \frac{v^r}{\Speed}}{\partial_t \Datafunctionforeikonal(t)} \restriction_{\lbrace r = 1 \rbrace}.
\end{align}
\end{subequations}

\subsection{Bootstrap assumptions and the away-from-time-axis and near-time-axis regions}
\label{SS:GLOBALEXISTENCEBOOTSTRAPANDTWOREGIONS}
To prove global existence, we will make bootstrap assumptions on a spacetime slab for a bootstrap time $\Tboot > 0$:
\begin{align} \label{E:GLOBALLEXISTENCEBOOTSTRAPSLAB}
	\globalexistenceslab{\Tboot}
	& := \lbrace r \geq 0 \rbrace \cap \lbrace 0 \leq t < \Tboot \rbrace.
\end{align}
We divide the bootstrap assumptions and the PDE analysis into the away-from-time-axis region: 
\begin{align} \label{E:AWAYFROMTIMEAXISGLOBALLEXISTENCEBOOTSTRAPSLAB}
	\awayglobalexistenceslab{\Tboot}
	& := \lbrace r \geq 1 \rbrace \cap \lbrace 0 \leq t < \Tboot \rbrace
	=
	\lbrace (t,r) \ | \ 0 \leq t < \Tboot \rbrace, u(t,r) \in (-\infty,\Datafunctionforeikonal(t)] \rbrace,
\end{align}
where the acoustic geometry and $(t,u)$ coordinates are defined,
and the near-time-axis region:
\begin{align} \label{E:NEARTIMEAXISGLOBALLEXISTENCEBOOTSTRAPSLAB}
	\nearglobalexistenceslab{\Tboot}
	& := \lbrace 0 \leq r \leq 1 \rbrace \cap \lbrace 0 \leq t < \Tboot \rbrace,
\end{align}
where the solution exhibits a lot of time-decay, enabling us to work only in $(t,r)$ coordinates.

\subsubsection{The away-from-time-axis region} 
\label{SSS:GLOBALEXISTENCERBIGGERTHANONEBOOTSTRAP}
We assume that on $\awayglobalexistenceslab{\Tboot}$,
the following inequalities hold for some number $\eps$ satisfying:
\begin{align} \label{E:GLOBALEXISTENCEBOUNDONBOOTSTRAPPARAMETERSIZE}
	0 & < \eps \leq \datasize^{3/4}.
\end{align}

\medskip 
\noindent \underline{\textbf{Bootstrap assumptions for $\RRiemann$ and $\LRiemann$}}.

\begin{subequations}
\begin{align} \label{E:GLOBALEXISTENCEFARFROMTIMEAXISPOINTWISEBOOTSTRAPRPLUS}
	|\RRiemann|
	& 
	\leq \frac{\eps}{1 + t + |u|},
		\\
|\partial_t \RRiemann|
		& 
	\leq \frac{\eps}{(1 + t + |u|)(1 + |u|)},
		\label{E:GLOBALEXISTENCEFARFROMTIMEAXISPOINTWISEBOOTSTRAPPARTIALTRPLUS}
			\\
	|\muuLunit \RRiemann|
	& 
	\leq \frac{\eps}{(1 + t + |u|)(1 + |u|)},
		\label{E:GLOBALEXISTENCEFARFROMTIMEAXISPOINTWISEBOOTSTRAPMULBARRPLUS}
		\\
|r \muuLunit \RRiemann|
	& 
	\leq
	\frac{\eps}{1 + t + |u|} 
	+
	\frac{\eps}{1 + u^2},
		\label{E:GLOBALEXISTENCEFARFROMTIMEAXISPOINTWISEBOOTSTRAPRTIMESMUULUNITRPLUS}
		\\
|r \uLunit \RRiemann|
	& 
	\leq
	C \frac{\eps}{1 + t + |u|} 
	+
	C \frac{\eps}{1 + u^2},
		\label{E:GLOBALEXISTENCEFARFROMTIMEAXISPOINTWISEBOOTSTRAPRTIMESULUNITRPLUS}
			\\
|\Lunit \RRiemann|
	& 
	\leq 
		\frac{\eps}{(1 + t + |u|)^2},
		\label{E:GLOBALEXISTENCEFARFROMTIMEAXISPOINTWISEBOOTSTRAPLRPLUS}
			\\
|r \Lunit \RRiemann|
	& 
	\leq 
		\frac{\eps}{1 + t + |u|},
		\label{E:GLOBALEXISTENCEFARFROMTIMEAXISPOINTWISEBOOTSTRAPRTIMESLRPLUS}
			\\
|r \muuLunit (\upmu \partial_t \RRiemann)|
		& 
	\leq  
		\eps 
		\frac{1}{(1 + t + |u|)(1 + |u|)}
		+
		\eps \frac{\ln_+(|u|)}{(1 + |u|)^3},
		\label{E:GLOBALEXISTENCEFARFROMTIMEAXISPOINTWISEBOOTSTRAPRTIMESMULBARMUPARTIALTRPLUS}
			\\
|r \Lunit \muuLunit \RRiemann|
	& 
	\leq 
	\frac{\eps}{(1 + t + |u|)(1 + |u|)},
	\label{E:GLOBALEXISTENCEFARFROMTIMEAXISPOINTWISEBOOTSTRAPTIMESLUNITMUULUNITRPLUS}
		\\
|r \muuLunit \muuLunit \RRiemann|
	& 
	\leq 
		\eps 
		\frac{1}{(1 + t + |u|)(1 + |u|)}
		+
		\eps \frac{\ln_+(|u|)}{(1 + |u|)^3},
	\label{E:GLOBALEXISTENCEFARFROMTIMEAXISPOINTWISEBOOTSTRAPTIMESMUULUNITMUULUNITRPLUS}
	\\
|r \Lunit \Lunit \RRiemann|
	& 
	\leq 
	\frac{\eps}{(1 + t + |u|)^2},
	\label{E:GLOBALEXISTENCEFARFROMTIMEAXISPOINTWISEBOOTSTRAPTIMESLUNITLUNITRPLUS}
\end{align}
\end{subequations}

\begin{subequations}
\begin{align} \label{E:GLOBALEXISTENCEFARFROMTIMEAXISPOINTWISEBOOTSTRAPRMINUS}
	|\LRiemann|
	& 
	\leq \frac{\eps}{1 + t + |u|},
		\\
|\partial_t \LRiemann|
		& 
	\leq \frac{\eps}{(1 + t + |u|)^2},
		\label{E:GLOBALEXISTENCEFARFROMTIMEAXISPOINTWISEBOOTSTRAPPARTIALTRMINUS}
			\\
|\uLunit \LRiemann|
		& 
	\leq \frac{\eps}{(1 + t + |u|)^2},
		\label{E:GLOBALEXISTENCEFARFROMTIMEAXISPOINTWISEBOOTSTRAPLBARRMINUS}
		\\
	|r \uLunit \LRiemann|
	& 
	\leq \frac{\eps}{1 + t + |u|},
		\label{E:GLOBALEXISTENCEFARFROMTIMEAXISPOINTWISEBOOTSTRAPRTIMESLBARRMINUS}
		\\
|\Lunit \LRiemann|
	& 
	\leq 
		\eps \frac{1}{(1 + t + |u|)^2},
		\label{E:GLOBALEXISTENCEFARFROMTIMEAXISPOINTWISEBOOTSTRAPLRMINUS}
			\\
|r \Lunit \LRiemann|
	& 
	\leq 
		\eps \frac{1}{1 + t + |u|},
		\label{E:GLOBALEXISTENCEFARFROMTIMEAXISPOINTWISEBOOTSTRAPRTIMESLRMINUS}
			\\
	|r \Lunit \Lunit \LRiemann|
	& 
	\leq 
	\eps \frac{1}{(1 + t + |u|)^2},
	\label{E:GLOBALEXISTENCEFARFROMTIMEAXISPOINTWISEBOOTSTRAPRTIMESLUNITLUNITRMINUS}
		\\
	|r \Lunit \partial_t \LRiemann|
		& 
	\leq 
	\frac{\eps}{(1 + t + |u|)^2},
		\label{E:GLOBALEXISTENCEFARFROMTIMEAXISPOINTWISEBOOTSTRAPRTIMESLUNITPARTIALTRMINUS}
			\\
	|r \muuLunit \Lunit \LRiemann|
	& 
	\leq 
	\frac{\eps}{(1 + t + |u|)^2}
	+
	\frac{\eps \ln_+(t-u)}{(1 + t + |u|)^2(1 + |u|)^2},
	\label{E:GLOBALEXISTENCEFARFROMTIMEAXISPOINTWISEBOOTSTRAPRTIMESMUULUNITLUNITRMINUS}
		\\
	|r \muuLunit \uLunit \LRiemann|
		& 
	\leq 
	\frac{\eps}{(1 + t + |u|)(1 + |u|)}.
		\label{E:GLOBALEXISTENCEFARFROMTIMEAXISPOINTWISEBOOTSTRAPRTIMESMUULUNITULUNITRMINUS}
\end{align}
\end{subequations}

\medskip

\medskip 
\noindent \underline{\textbf{Bootstrap assumptions for $\upmu$}}.

\begin{subequations}
\begin{align}
	 \label{E:GLOBALEXISTENCEMUSIMPLERBOUNDSBOOTSTRAP}
	1
	- 
	\eps
	& 
	\leq
	\upmu 
	\leq 1 
			+
			\eps
			+ 
			\frac{\eps}{(1 + u^2)} \ln_+(t - u),
				\\
	|\muuLunit \upmu|
	& \leq  
			\eps \frac{\ln_+(u_+)}{1 + |u|^2}
			+
			\eps
			\frac{\ln_+ \left(t-u \right) \ln_+(|u|)}{(1 + |u|)^3}.
	\label{E:GLOBALEXISTENCEMUXMUSIMPLERUPPERBOUNDBOOTSTRAP}		
\end{align}
\end{subequations}

\subsubsection{The near-time-axis region} 
\label{SSS:GLOBALEXISTENCELESSTHANONEBOOTSTRAP}
We assume that on $\nearglobalexistenceslab{\Tboot}$,
the following inequalities hold for some number $\eps$ satisfying \eqref{E:GLOBALEXISTENCEBOUNDONBOOTSTRAPPARAMETERSIZE}.

\medskip

\noindent \underline{\textbf{Bootstrap assumptions for the Riemann invariants}}.

\begin{subequations}
\begin{align} \label{E:GLOBALEXISTENCENEARTIMEAXISPOINTWISEBOOTSTRAPRPLUS}
	|\RRiemann|
	& 
	\leq \frac{\eps}{(1 + t + r)},
		\\
|\partial_t \RRiemann|
		& 
	\leq C \frac{\eps}{(1 + t + r)(1 + |t-r|)},
		\label{E:GLOBALEXISTENCENEARTIMEAXISPOINTWISEBOOTSTRAPPARTIALTRPLUS}
			\\
	|\uLunit \RRiemann|
	& 
	\leq \frac{\eps}{(1 + t + r)(1 + |t-r|)},
		\label{E:GLOBALEXISTENCENEARTIMEAXISPOINTWISEBOOTSTRAPLBARRPLUS}
		\\
|r \uLunit \RRiemann|
	& 
	\leq
	\frac{\eps}{1 + t + r}
	+
	\frac{\eps}{(1 + |t - r|)^2},
		\label{E:GLOBALEXISTENCENEARTIMEAXISPOINTWISEBOOTSTRAPRTIMESLBARRPLUS}
			\\
|\Lunit \RRiemann|
	& 
	\leq 
		\frac{\eps}{(1 + t + r)^2},
		\label{E:GLOBALEXISTENCENEARTIMEAXISPOINTWISEBOOTSTRAPLRPLUS}
			\\
|r \Lunit \RRiemann|
	& 
	\leq 
		\frac{\eps}{1 + t + r},
		\label{E:GLOBALEXISTENCENEARTIMEAXISPOINTWISEBOOTSTRAPRTIMESLRPLUS}
			\\
|r \uLunit \partial_t \RRiemann|
		& 
	\leq  
		\frac{\eps}{(1 + t + r)(1 + |t-r|)},
		\label{E:GLOBALEXISTENCENEARTIMEAXISPOINTWISEBOOTSTRAPRTIMESLBARPARTIALTRPLUS}
			\\
|r \Lunit \uLunit \RRiemann|
	& 
	\leq 
	\frac{\eps}{(1 + t + r)(1 + |t-r|)},
	\label{E:GLOBALEXISTENCENEARTIMEAXISPOINTWISEBOOTSTRAPTIMESLUNITLBARRPLUS}
		\\
|r \uLunit \uLunit \RRiemann|
	& 
	\leq 
	\frac{\eps}{(1 + t + r)(1 + |t-r|)},
	\label{E:GLOBALEXISTENCENEARTIMEAXISPOINTWISEBOOTSTRAPTIMESLBARLBARRPLUS}
	\\
|r \Lunit \Lunit \RRiemann|
	& 
	\leq 
	\frac{\eps}{(1 + t + r)^2},
	\label{E:GLOBALEXISTENCENEARTIMEAXISPOINTWISEBOOTSTRAPTIMESLUNITLUNITRPLUS}
\end{align}
\end{subequations}

\begin{subequations}
\begin{align} \label{E:GLOBALEXISTENCENEARTIMEAXISPOINTWISEBOOTSTRAPRMINUS}
	|\LRiemann|
	& 
	\leq \frac{\eps}{(1 + t + r)},
		\\
|\partial_t \LRiemann|
		& 
	\leq
		\frac{\eps}{(1 + t + r)^3}	
		+
		\frac{\eps^{1^+} \ln_+(t+r)}{(1 + t + r)^3},
		\label{E:GLOBALEXISTENCENEARTIMEAXISPOINTWISEBOOTSTRAPPARTIALTRMINUS}
			\\
|\uLunit \LRiemann|
		& 
	\leq 
		\frac{\eps}{(1 + t + r)^3}	
		+
		\frac{\eps^{1^+} \ln_+(t+r)}{(1 + t + r)^3},
		\label{E:GLOBALEXISTENCENEARTIMEAXISPOINTWISEBOOTSTRAPLBARRMINUS}
		\\
	|r \uLunit \LRiemann|
	& 
	\leq 
		\frac{\eps r}{(1 + t + r)^3}	
		+
		\frac{\eps^{1^+} \ln_+(t+r) r}{(1 + t + r)^3},
		\label{E:GLOBALEXISTENCENEARTIMEAXISPOINTWISEBOOTSTRAPRTIMESLBARRMINUS}
		\\
|\Lunit \LRiemann|
	& 
	\leq 
		\frac{\eps}{(1 + t + r)^3}	
		+
		\frac{\eps^{1^+} \ln_+(t+r)}{(1 + t + r)^3},
		\label{E:GLOBALEXISTENCENEARTIMEAXISPOINTWISEBOOTSTRAPLRMINUS}
			\\
|r \Lunit \LRiemann|
	& 
	\leq 
		\frac{\eps \ln(t+r)}{(1 + t + r)^3}	
		+
		\frac{\eps^{1^+} [\ln(t+r)]^2}{(1 + t + r)^3},
		\label{E:GLOBALEXISTENCENEARTIMEAXISPOINTWISEBOOTSTRAPRTIMESLRMINUS}
			\\
	|r \Lunit \Lunit \LRiemann|
	& 
	\leq 
		\frac{\eps}{(1 + t + r)^3}	
		+
		\frac{\eps^{1^+}\ln_+(t+r)}{(1 + t + r)^3},
	\label{E:GLOBALEXISTENCENEARTIMEAXISPOINTWISEBOOTSTRAPRTIMESLUNITLUNITRMINUS}
		\\
	|r \Lunit \partial_t \LRiemann|
		& 
	\leq 
	\frac{\eps}{(1 + t + r)^3}	
		+
		\frac{\eps^{1^+} \ln_+(t+r)}{(1 + t + r)^3},
		\label{E:GLOBALEXISTENCENEARTIMEAXISPOINTWISEBOOTSTRAPRTIMESLUNITPARTIALTRMINUS}
			\\
	|r \uLunit \Lunit \LRiemann|
	& 
	\leq 
	\frac{\eps}{(1 + t + r)^2},
	\label{E:GLOBALEXISTENCENEARTIMEAXISPOINTWISEBOOTSTRAPRTIMESLBARLUNITRMINUS}
		\\
	|r \uLunit \uLunit \LRiemann|
		& 
	\leq 
	\frac{\eps}{(1 + t + r)(1 + |t-r|)}.
		\label{E:GLOBALEXISTENCENEARTIMEAXISPOINTWISEBOOTSTRAPRTIMESLBARLBARRMINUS}
\end{align}
\end{subequations}

\subsection{Preliminary estimates}
\label{SS:GLOBALEXISTENCEPRELIMINARYESTIMATES}
To prove our main global existence results,
we will rely on some general estimates for the coordinate functions and integral inequalities,
which we derive in the next lemma with the help of the bootstrap assumptions.

\begin{lemma}[Coordinate function estimates and integral inequalities]
	\label{L:COORDINATEFUNCTIONESTIMATESANDINTEGRALINEQUALITIES}
 Let $\ln_+$ be the function defined by:
	\begin{align} \label{E:GLOBALEXISTENCELOGPLUS}
		\ln_+(z) 
		& := \ln(e + z).	
	\end{align}
Under the data assumptions of Sect.\,\ref{SS:GLOBALEXISTENCEDATAATTIME0} and the bootstrap assumptions of
Sect.\,\ref{SS:GLOBALEXISTENCEBOOTSTRAPANDTWOREGIONS}, 
if $\datasize > 0$ is sufficiently small, then
the following estimates hold.
	
\noindent \underline{\textbf{Estimates for the eikonal function in} $\awayglobalexistenceslab{\Tboot}$}.
The eikonal function $u$ satisfies the following estimates in the away-from-time-axis region
$\awayglobalexistenceslab{\Tboot}$ defined in \eqref{E:AWAYFROMTIMEAXISGLOBALLEXISTENCEBOOTSTRAPSLAB}:
\begin{align} 
	\begin{split} \label{E:GLOBALEXISTENCESHARPPOINWISECOMPARISONBETWEENUANDTMINUSR}
	u 
	& = 
	\begin{cases}
	t - r 
	+
	\mathcal{O}(\eps)
	\ln \left( \frac{1 + t + |u|}{1 + |u|} \right),
	&
	u \leq - 1,
		\\
	t - r 
	+
	\mathcal{O}(\eps)
	\ln_+ \left( \frac{1 + t + |u|}{2 + u + |u|} \right),
	& u \geq -1
\end{cases} 
= t - r + \mathcal{O}(\eps)
	\ln_+ \left( \frac{1 + t + |u|}{1 + |u| + (1 + u)_+} \right).
\end{split}
\end{align}

\medskip

\noindent \underline{\textbf{Comparison estimates for various coordinate functions in} $\awayglobalexistenceslab{\Tboot}$}.
The following estimates hold
 in the away-from-time-axis region
$\awayglobalexistenceslab{\Tboot}$ defined in \eqref{E:AWAYFROMTIMEAXISGLOBALLEXISTENCEBOOTSTRAPSLAB}:
\begin{subequations}
	\begin{align}  \label{E:GLOBALEXISTENCEUSMALLERTHANTMINUS1}
		u & \leq 
		\begin{cases}
			t - 1 + \mathcal{O}(\datasize), & \mbox{if } u \in (-1,1), 
				\\
			t - 1, & \mbox{if } u \in (-\infty,-1] \cup [1,\infty), 
		\end{cases}
			\\
		1 
		& 
		\leq 
		t - u
		+
		\mathcal{O}(\datasize)
		\approx 1 + t - u
		\leq
		1 + t + |u|
		\approx 2t - u
		\approx 1 + t + r,
		\label{E:GLOBALEXISTENCECOMPARISONBETWEENTOVERUTMINUSUANDTPLUSMODU}
			\\
		r & = [1 + \mathcal{O}(\eps)](t-u).
			\label{E:GLOBALEXISTENCERISANALMOSTUNITYMULIPLEOFTMINUSU}
	\end{align}
\end{subequations}

\medskip

\noindent \underline{\textbf{Cartesian coordinate function estimates along integral curves of $\uLunit$ in the entire bootstrap slab}}.
Let $(t_0,r_0) \in \globalexistenceslab{\Tboot}$ and $(t_1,r_1) \in \globalexistenceslab{\Tboot}$,
and assume that 
$(t_1,r_1)$ lies on the future-directed integral curve of $\uLunit$ emanating from
$(t_0,r_0)$. 
Then the following estimates hold:
\begin{subequations}
\begin{align} \label{E:GLOBALEXISTENCEINTERIORREGIONPOINTWISEESTIMATECHANGEINTPLUSRALONGINTEGRALCURVESOFLBAR}
	t_1 + r_1 
	& = t_0 + r_0
	+ \mathcal{O}(\eps) \frac{t_1 - t_0}{(t_1 + r_1)},
		\\
	t_1 - r_1
	& = t_0 - r_0
		+ 
		2(t_1 - t_0)
		+
		\mathcal{O}(\eps) \frac{t_1 - t_0}{(t_1 + r_1)},
		\label{E:GLOBALEXISTENCEINTERIORREGIONPOINTWISEESTIMATECHANGEINTMINUSRALONGINTEGRALCURVESOFLBAR}
			\\
\frac{1}{t_1 + r_1}
	& = \frac{1}{t_0 + r_0}
	+ 
	\mathcal{O}(\eps) \frac{t_1 - t_0}{(t_1 + r_1)^2(t_0 + r_0)},
	\label{E:GLOBALEXISTENCERECIPROCALTPLUSRNEARLYCONSTANTALONGINTEGRALCURVESOFULUNIT}
		\\
	\frac{1}{2(1 + t_0 + r_0)}
	& \leq	
	\frac{1}{1 + t_1 + r_1}
	\leq \frac{2}{1 + t_0 + r_0}.
	\label{E:GLOBALEXISTENCERECIPROCALONEPLUSTPLUSRNEARLYCONSTANTALONGINTEGRALCURVESOFULUNIT}
\end{align}
\end{subequations}

\medskip
	\noindent \underline{\textbf{Geometric coordinate function estimates along integral curves of $\uLunit$ in} $\awayglobalexistenceslab{\Tboot}$}.
Let $(t_0,u_0), (t_1,u_1)$ be a pair of points in $\awayglobalexistenceslab{\Tboot}$ expressed relative to the geometric coordinates, 
and assume that $(t_1,u_1)$ lies on the future-directed integral curve of $\uLunit$ emanating from
$(t_0,u_0)$. Then the following estimates hold:
\begin{subequations}
	\begin{align} 
		\begin{split} \label{E:GLOBALEXISTENCETWOTMINUSUALONGINTEGRALCURVESOFULUNIT}
		2t_1 - u_1
		& = 2 t_0 - u_0 
		+
		\mathcal{O}(\eps) \ln_+ \left( \frac{1 + t_0 + |u_0|}{1 + |u_0| + (1 + u_0)_+} \right)
		+
		\mathcal{O}(\eps) \ln_+ \left( \frac{1 + t_1 + |u_1|}{1 + |u_1| + (1 + u_1)_+} \right)
			\\
		& \ \
		+
		\mathcal{O}(\eps) \frac{t_1 - t_0}{(1 + t_1 + |u_1|)},
	\end{split}
		\\
		\frac{1}{2 t_0 - u_0}
		& = \frac{1}{2 t_1 - u_1}
			+ 
			\mathcal{O}(\eps) \frac{\ln_+(t_1)}{(2t_1 - u_1)^2}.	
		\label{E:ONEOVERGLOBALEXISTENCETWOTMINUSUALONGINTEGRALCURVESOFULUNIT}
	\end{align}
	\end{subequations}
	
	
	
	
	\medskip

\noindent \underline{\textbf{Estimates for integrals along integral curves of $\uLunit$ in} $\nearglobalexistenceslab{\Tboot}$}.
Let $t \rightarrow (t,\underline{\mathfrak{r}}(t))$ be the integral curve of $\uLunit$
joining the point $(t_1,r_1)$ to the point $(t_2,r_2)$ lying to its future.
Assume that $(t_1,r_1)$ and $(t_2,r_2)$ both lie in the near-time axis region 
$\nearglobalexistenceslab{\Tboot}$ defined in \eqref{E:NEARTIMEAXISGLOBALLEXISTENCEBOOTSTRAPSLAB}.
Let $A \geq 0$ and $B \geq 0$ be constants.
Then the following estimate holds, where the implicit constants
depend on $A$ and $B$:
\begin{align} \label{E:GLOBALEXISTENCENEARTIMEAXISINTEGRALESTIMATEONEOVERONEPLUSTPLUSRTOBTIMESONEOVER1PLUSTMINUSRTOA}
	\int_{t_1}^{t_2} 
		\frac{1}{[1 + t + \underline{\mathfrak{r}}(t)]^B
		(1 + |t - \underline{\mathfrak{r}}(t)|)^A} 
	\, \mathrm{d} t
	& \lesssim 
		\frac{1}{(1 + t_2 + r_2)^{A+B}}.
\end{align}

\medskip
	\noindent \underline{\textbf{Estimates of integrals involving the geometric coordinates in} $\awayglobalexistenceslab{\Tboot}$}.	
	Let $u \rightarrow (\mathfrak{t}_{t_0,u_0}(u),u)$ be the $u$-parameterized integral curve of 
	$\muuLunit = \upmu \frac{\partial}{\partial t} + 2 \frac{\partial}{\partial u}$ emanating from the point
	$(t_0,u_0) \in \Sigma_{t_0}$, i.e., $\mathfrak{t}_{t_0,u_0}$ is the solution to the following initial value problem:
	\begin{align} \label{E:GLOBALEXISTENCEIVPFORINTEGRALCURVESOFMUULUNITPARAMETERIZEDBYU}
		\frac{d}{d u}
		\mathfrak{t}_{t_0,u_0}(u)
		& = \frac{1}{2} \upmu(\mathfrak{t}(u),u),
		&&
		\mathfrak{t}_{t_0,u_0}(u_0) 
		= t_0. 
	\end{align}
	Let $u_1 \geq u_0$, set $t_1 := \mathfrak{t}_{t_0,u_0}(u_1)$, let $r_0 = r(t_0,u_0)$ be the radial value 
	corresponding to the point with geometric coordinates $(t_0,u_0)$, and let
	let $r_1 := r(t_1,u_1)$ be the corresponding
	radial value at the point with geometric coordinates $(t_1,u_1)$.
	Assume that $(t_0,u_0)$ and $(t_1,u_1)$
	lie in the 
	away-from-time-axis region
	$\awayglobalexistenceslab{\Tboot}$ defined in \eqref{E:AWAYFROMTIMEAXISGLOBALLEXISTENCEBOOTSTRAPSLAB}.
	Let $A \geq 0$ and $\updelta > 0$ be constants.
	Then there exists a $C > 0$, depending on $A$ and $\updelta$, such that the following estimates hold:
	\begin{subequations}
	\begin{align} \label{E:GLOBALEXISTENCEINTEGRALESTIMATEALONGMUULUNITONEOVERTPLUSUTOPOWERATIMESONEOVER1PLUSU}
	\int_{u_0}^{u_1}
		\frac{1}{(1 + \mathfrak{t}_{u_0}(u) + |u|)^A (1 + |u|)}
	\, \mathrm{d} u
	& \leq 
			C
			\frac{\ln_+(|u_1| + |u_0|)}{(1 + t_1 + |u_1|)^A},
				\\
	\int_{u_0}^{u_1}
		\frac{1}{(1 + \mathfrak{t}_{u_0}(u) + |u|)^A (1 + |u|)^{1 + \updelta}}
	\, \mathrm{d} u
	& \leq 
			C
			\frac{1}{(1 + t_1 + |u_1|)^A}.
			 \label{E:GLOBALEXISTENCEINTEGRALESTIMATEALONGMUULUNITONEOVERTPLUSUTOPOWERATIMSINTEGRABLEINU}
	\end{align}
	\end{subequations}

\medskip

\noindent \underline{\textbf{Behavior of coordinate functions along integral curves of $\Lunit$ in} $\globalexistenceslab{\Tboot}$}.
Let $(t_0,r_0) \in \globalexistenceslab{\Tboot}$ and $(t_1,r_1) \in \globalexistenceslab{\Tboot}$,
and assume that $(t_1,r_1)$ lies on the future-directed integral curve of $\Lunit$ emanating from
$(t_0,r_0)$. Then the following estimates hold:

\begin{subequations}
\begin{align} 
\begin{split} \label{E:GLOBALEXISTENCEBOUNDFORTPLUSRALONGINTEGRALCURVESOFL}
	t_1 + r_1 
	& = 
	t_0 + r_0
	+ 
	2(t_1 - t_0)
	+
	\mathcal{O}(\eps)
	\ln \left( \frac{1 + t_1}{1 + t_0} \right),
	\end{split}
	\\
	\begin{split} \label{E:GLOBALEXISTENCEBOUNDFORTMINUSRALONGINTEGRALCURVESOFL}
	t_1 - r_1
	& = t_0 - r_0
		+
		\mathcal{O}(\eps)
		\ln \left( \frac{1 + t_1}{1 + t_0} \right).
	\end{split}
\end{align}
\end{subequations}

\medskip

\noindent \underline{\textbf{Estimates for integrals along integral curves of $\Lunit$ in} $\globalexistenceslab{\Tboot}$ 
\textbf{with a smallness assumption}}.
Let $(t_0,r_0) \in \globalexistenceslab{\Tboot}$ and $(t_1,r_1) \in \globalexistenceslab{\Tboot}$,
and let $t \rightarrow (t,\mathfrak{r}(t))$ be the integral curve of $\Lunit$ joining the point 
$(t_0,r_0)$ to the point
$(t_1,r_1)$ lying to its future. 
Assume that $r_1 \leq 1 + t_0$ or $t_1 \leq 1 + r_0$.
Then for any constants $A \geq 0$ and $B \geq 0$, the following
estimates hold, where the implicit constants depend on $A$ and $B$:
\begin{subequations}
\begin{align} \label{E:GLOBALEXISTENCENEARTIMEAXISINTEGRALESTIMATEROVER1PLUSTPLUSRPOWERAALONGINTEGRALCURVESOFL}
	\int_{t_0}^{t_1} \frac{\mathfrak{r}^A(t)}{[1 + t + \mathfrak{r}(t)]^B} \, \mathrm{d} t
	& \lesssim 
	\frac{r_1^{1+A}}{(1 + t_1 + r_1)^B},
		\\
	\int_{t_0}^{t_1} \frac{\mathfrak{r}^A(t) \ln_+(t + \mathfrak{r}(t))}{[1 + t + \mathfrak{r}(t)]^B} \, \mathrm{d} t
	& \lesssim 
	\frac{r_1^{1+A} \ln_+(t_1 + r_1)}{(1 + t_1 + r_1)^B},
	\label{E:GLOBALEXISTENCENEARTIMEAXISINTEGRALESITMATELOGTPLUSRTIMESRTOAOVERONEPLUSTPLUSRTOBALONGINTEGRALCURVESOFL}
		\\
\int_{t_0}^{t_1} \frac{\mathfrak{r}^A(t)}{[1 + t + \mathfrak{r}(t)]^{1+B}[1 + \mathfrak{r}(t)]} \, \mathrm{d} t
	& \lesssim 
	\frac{r_1^{1+A}}{(1 + t_1 + r_1)^{1+B}(1+r_0)}.
	\label{E:GLOBALEXISTENCENEARTIMEAXISINTEGRALESTIMATEROVER1PLUSTPLUSRPOWERATIMESONEPLUSTPLUSRALONGINTEGRALCURVESOFL}
\end{align}
\end{subequations}

\noindent \underline{\textbf{Estimates for integrals along integral curves of $\Lunit$ in} $\globalexistenceslab{\Tboot}$ \textbf{without a smallness assumption}}.
Let $(t_1,r_1) \in \globalexistenceslab{\Tboot}$ and $(t_2,r_2) \in \globalexistenceslab{\Tboot}$,
and let $t \rightarrow (t,\mathfrak{r}(t))$ be the integral curve of $\Lunit$ joining the point 
$(t_1,r_1)$ to the point
$(t_2,r_2)$ lying to its future. 
Let $\updelta > 0$ be a constant.
Then the following estimates hold, where the implicit constants depend on $\updelta$:
\begin{subequations}
\begin{align} \label{E:GLOBALEXISTENCEFARFROMTIMEAXISINTEGRALESTIMATEROVER1PLUSTPLUSRSQUAREDALONGINTEGRALCURVESOFL}
	\int_{t_1}^{t_2} \frac{\mathfrak{r}(t)}{(1 + t + \mathfrak{r}(t))^2} \, \mathrm{d} t
	& \lesssim 
		\ln \left( \frac{1 + r_2}{1 + r_1} \right),
		\\
	\int_{t_1}^{t_2} \frac{1}{1 + t + \mathfrak{r}(t)} \, \mathrm{d} t
	& \lesssim 
			\ln\left( \frac{1 + r_2}{1 + r_1} \right),
	\label{E:GLOBALEXISTENCEFARFROMTIMEAXISONEOVERONEPLUSTPLUSRALONGINTEGRALCURVESOFL}
		\\
	\int_{t_1}^{t_2} \frac{\ln(1+ t + \mathfrak{r}(t))}{1 + t + \mathfrak{r}(t)} \, \mathrm{d} t
	& \lesssim 
			[\ln(1 + t_2 + r_2)]^2 - [\ln(1 + t_1 + r_1)]^2,
			\label{E:GLOBALEXISTENCEFARFROMTIMEAXISLOGTPLUSROVERONEPLUSTPLUSRALONGINTEGRALCURVESOFL}
				\\
	\int_{t_1}^{t_2} \frac{1}{(1 + t + \mathfrak{r}(t))^{1 + \updelta}} \, \mathrm{d} t
	& \lesssim 
			\frac{1}{(1 + t_1 + r_1)^{\updelta}},
	\label{E:GLOBALEXISTENCEFARFROMTIMEAXISONEOVERONEPLUSTPLUSRTOPOWERONEPLUSDELTAALONGINTEGRALCURVESOFL}
		\\
	\int_{t_1}^{t_2} \frac{\ln(1+ t + \mathfrak{r}(t))}{(1 + t + \mathfrak{r}(t))^{1 + \updelta}} \, \mathrm{d} t
	& \lesssim 
			\frac{\ln(1+ t_1 + r_1)}{(1 + t_1 + r_1)^{\updelta}}.
	\label{E:GLOBALEXISTENCEFARFROMTIMEAXISLOGTPLUSROVERONEPLUSTPLUSRTOPOWERONEPLUSDELTAALONGINTEGRALCURVESOFL}
\end{align}
\end{subequations}


\end{lemma}

\begin{proof}

\noindent \textbf{Proof of \eqref{E:GLOBALEXISTENCESHARPPOINWISECOMPARISONBETWEENUANDTMINUSR} and \eqref{E:GLOBALEXISTENCEUSMALLERTHANTMINUS1}}:
To prove \eqref{E:GLOBALEXISTENCESHARPPOINWISECOMPARISONBETWEENUANDTMINUSR},
	we first repeat the proof of \eqref{E:POINTWISEESTIMATEFORLDERIVATIVEOFUMINUSTMINUSR} to reach the same conclusion:
	\begin{align} \label{E:GLOBALEXISTENCEPOINTWISEESTIMATEFORLDERIVATIVEOFUMINUSTMINUSR}
	\Lunit [u - (t-r)] 
	& = \frac{\mathcal{O}(\eps)}{(1 + t + |u|)}.
\end{align}
Recalling that $\Lunit = \frac{\partial}{\partial t}$ in geometric coordinates,
we integrate 
\eqref{E:GLOBALEXISTENCEPOINTWISEESTIMATEFORLDERIVATIVEOFUMINUSTMINUSR}
and use the data assumptions
\eqref{E:GLOBALEXISTENCESPACELIKEEIKONALFUNCTIONDATA}--\eqref{E:GLOBALEXISTENCETIMELKEEIKONALFUNCTIONDATA}
as well as \eqref{E:TIMELILKEEIKONALDATAISAPERTURBATIONOFTMINUS1}--\eqref{E:C2ESTIMATEFORTIMELILKEEIKONALDATAPERTURBATION},
thereby concluding \eqref{E:GLOBALEXISTENCESHARPPOINWISECOMPARISONBETWEENUANDTMINUSR}.

Similarly, \eqref{E:GLOBALEXISTENCEUSMALLERTHANTMINUS1} follows from integrating the identity $\Lunit [u - (t-1)] = -1$ and using the data assumptions for $u$.

\medskip

\noindent \textbf{Proof of \eqref{E:GLOBALEXISTENCECOMPARISONBETWEENTOVERUTMINUSUANDTPLUSMODU}}:
The inequalities in \eqref{E:GLOBALEXISTENCECOMPARISONBETWEENTOVERUTMINUSUANDTPLUSMODU} follow in a straightforward fashion
from \eqref{E:GLOBALEXISTENCESHARPPOINWISECOMPARISONBETWEENUANDTMINUSR} and \eqref{E:GLOBALEXISTENCEUSMALLERTHANTMINUS1}.

\medskip

\noindent \textbf{Proof of \eqref{E:GLOBALEXISTENCERISANALMOSTUNITYMULIPLEOFTMINUSU}}:
To prove \eqref{E:GLOBALEXISTENCERISANALMOSTUNITYMULIPLEOFTMINUSU}, we let $\widetilde{\Lunit} := \frac{1}{v^r + \Speed} \Lunit$
be a rescaled version of $\Lunit$. By \eqref{E:NULLVECTORFIELDS}, we have $\widetilde{\Lunit} r = 1$.
Multiplying \eqref{E:GLOBALEXISTENCEPOINTWISEESTIMATEFORLDERIVATIVEOFUMINUSTMINUSR} by $\frac{1}{v^r + \Speed}$ and using
\eqref{E:SPEEDOFSOUNDEXPANSION}--\eqref{E:SPEEDOFSOUNDERRORFUNCTIONVANISHESATORIGIN},
\eqref{E:GLOBALEXISTENCECOMPARISONBETWEENTOVERUTMINUSUANDTPLUSMODU}, 
and the bootstrap assumptions, we see that:
\begin{align} \label{E:GLOBALEXISTENCEPOINTWISEESTIMATEFORRESCALEDLDERIVATIVEOFUMINUSTMINUSR}
	\widetilde{\Lunit} [u - (t-r)] 
	& = \frac{\mathcal{O}(\eps)}{(1 + r)}.
\end{align}
Integrating \eqref{E:GLOBALEXISTENCEPOINTWISEESTIMATEFORRESCALEDLDERIVATIVEOFUMINUSTMINUSR},
using that $\widetilde{\Lunit} r = 1$,
and using the data assumptions \eqref{E:GLOBALEXISTENCESPACELIKEEIKONALFUNCTIONDATA}--\eqref{E:GLOBALEXISTENCETIMELKEEIKONALFUNCTIONDATA},
we find that 
$
u - (t-r) = \mathcal{O}(\datasize) + \mathcal{O}(\eps) \ln_+(1 + r)
$,
i.e., 
$
r = t - u + \mathcal{O}(\datasize) + \mathcal{O}(\eps) \ln_+(1 + r)
$.
Dividing this equation by $r$ and using that $r \geq 1$ in the away-from-time-axis region, we conclude that:
$
\frac{t-u}{r}  = 1 + \mathcal{O}(\datasize) + \mathcal{O}(\eps) \frac{\ln_+(1 + r)}{r} = 1 + \mathcal{O}(\eps),
$
which yields \eqref{E:GLOBALEXISTENCERISANALMOSTUNITYMULIPLEOFTMINUSU}.

\medskip

\noindent \textbf{Proof of \eqref{E:GLOBALEXISTENCEINTERIORREGIONPOINTWISEESTIMATECHANGEINTPLUSRALONGINTEGRALCURVESOFLBAR}--\eqref{E:GLOBALEXISTENCERECIPROCALONEPLUSTPLUSRNEARLYCONSTANTALONGINTEGRALCURVESOFULUNIT}}:
Thanks to the availability of \eqref{E:GLOBALEXISTENCECOMPARISONBETWEENTOVERUTMINUSUANDTPLUSMODU},
\eqref{E:GLOBALEXISTENCEINTERIORREGIONPOINTWISEESTIMATECHANGEINTPLUSRALONGINTEGRALCURVESOFLBAR}--\eqref{E:GLOBALEXISTENCERECIPROCALONEPLUSTPLUSRNEARLYCONSTANTALONGINTEGRALCURVESOFULUNIT}
can be proved using the same arguments we used to prove
\eqref{E:INTERIORREGIONPOINTWISEESTIMATECHANGEINTPLUSRALONGINTEGRALCURVESOFLBAR}--\eqref{E:RECIPROCALONEPLUSTPLUSRNEARLYCONSTANTALONGINTEGRALCURVESOFULUNIT}.

\medskip
\noindent \textbf{Proof of \eqref{E:GLOBALEXISTENCETWOTMINUSUALONGINTEGRALCURVESOFULUNIT} and \eqref{E:ONEOVERGLOBALEXISTENCETWOTMINUSUALONGINTEGRALCURVESOFULUNIT}}:
To prove \eqref{E:GLOBALEXISTENCETWOTMINUSUALONGINTEGRALCURVESOFULUNIT},
we use \eqref{E:GLOBALEXISTENCESHARPPOINWISECOMPARISONBETWEENUANDTMINUSR}
twice to substitute in \eqref{E:GLOBALEXISTENCEINTERIORREGIONPOINTWISEESTIMATECHANGEINTPLUSRALONGINTEGRALCURVESOFLBAR}, once
for $t_1$ and once for $r_0$, and we also use \eqref{E:GLOBALEXISTENCECOMPARISONBETWEENTOVERUTMINUSUANDTPLUSMODU}.

\eqref{E:ONEOVERGLOBALEXISTENCETWOTMINUSUALONGINTEGRALCURVESOFULUNIT} follows from
\eqref{E:GLOBALEXISTENCETWOTMINUSUALONGINTEGRALCURVESOFULUNIT}
and the following Taylor expansion estimate,  
valid for real numbers satisfying 
$\upbeta = \upalpha + \upkappa$ with $\upbeta > 0$ and $\frac{|\upkappa|}{\upbeta} \leq \frac{1}{2}$:
$\frac{1}{\upalpha} = \frac{1}{\upbeta} + \mathcal{O}(\upkappa) \frac{1}{\upbeta^2}$.
We clarify that to apply this Taylor expansion estimate to obtain \eqref{E:ONEOVERGLOBALEXISTENCETWOTMINUSUALONGINTEGRALCURVESOFULUNIT},
the role of $\upalpha$ is played by $2 t_0 - u_0$, 
the role of $\upbeta$ is played by $2 t_1 - u_1$, 
and the role of the equation $\upbeta = \upalpha + \upkappa$ 
is played by \eqref{E:GLOBALEXISTENCETWOTMINUSUALONGINTEGRALCURVESOFULUNIT}.

\medskip
\noindent \textbf{Proof of \eqref{E:GLOBALEXISTENCEINTEGRALESTIMATEALONGMUULUNITONEOVERTPLUSUTOPOWERATIMSINTEGRABLEINU}}:
To prove \eqref{E:GLOBALEXISTENCEINTEGRALESTIMATEALONGMUULUNITONEOVERTPLUSUTOPOWERATIMSINTEGRABLEINU}, we first use 
\eqref{E:GLOBALEXISTENCECOMPARISONBETWEENTOVERUTMINUSUANDTPLUSMODU} and
\eqref{E:ONEOVERGLOBALEXISTENCETWOTMINUSUALONGINTEGRALCURVESOFULUNIT}
to bound the integrand factor 
$
\frac{1}{(1 + \mathfrak{t}_{u_0}(u) + |u|)^A}
$
on
LHS~\eqref{E:GLOBALEXISTENCEINTEGRALESTIMATEALONGMUULUNITONEOVERTPLUSUTOPOWERATIMSINTEGRABLEINU}
by 
$\lesssim
\frac{1}{(1 + 2 \mathfrak{t}_{u_0}(u) - u)^A}
\lesssim 
\frac{1}{(1 + 2 t_1 - u_1)^A}
\lesssim 
\frac{1}{(1 + t_1 + |u_1|)^A}
$.
From this bound, we see that the integrand
on LHS~\eqref{E:GLOBALEXISTENCEINTEGRALESTIMATEALONGMUULUNITONEOVERTPLUSUTOPOWERATIMSINTEGRABLEINU} is bounded by:
\begin{align} \label{E:PROOFSTEPGLOBALEXISTENCEINTEGRALESTIMATEALONGMUULUNITONEOVERTPLUSUTOPOWERATIMSINTEGRABLEINU}
\lesssim
	\frac{1}{(1 + t_1 + |u_1|)^A}
	\int_{u_0}^{u_1}
		\frac{1}{(1 + |u|)^{1 + \updelta}}
	\, \mathrm{d} u
	& \lesssim
			\frac{1}{(1 + t_1 + |u_1|)^A},
\end{align}
which yields \eqref{E:GLOBALEXISTENCEINTEGRALESTIMATEALONGMUULUNITONEOVERTPLUSUTOPOWERATIMSINTEGRABLEINU}.

\medskip
\noindent \textbf{Proof of \eqref{E:GLOBALEXISTENCEINTEGRALESTIMATEALONGMUULUNITONEOVERTPLUSUTOPOWERATIMESONEOVER1PLUSU}}:
The estimate \eqref{E:GLOBALEXISTENCEINTEGRALESTIMATEALONGMUULUNITONEOVERTPLUSUTOPOWERATIMESONEOVER1PLUSU} can be proved using arguments similar to the ones we used to prove \eqref{E:GLOBALEXISTENCEINTEGRALESTIMATEALONGMUULUNITONEOVERTPLUSUTOPOWERATIMSINTEGRABLEINU}, and we omit the details.

\medskip
\noindent \textbf{Proof of \eqref{E:GLOBALEXISTENCEBOUNDFORTPLUSRALONGINTEGRALCURVESOFL}--\eqref{E:GLOBALEXISTENCEBOUNDFORTMINUSRALONGINTEGRALCURVESOFL}}:
\eqref{E:GLOBALEXISTENCEBOUNDFORTPLUSRALONGINTEGRALCURVESOFL}--\eqref{E:GLOBALEXISTENCEBOUNDFORTMINUSRALONGINTEGRALCURVESOFL}
can be proved using the same arguments we used to prove
the first equalities in
\eqref{E:BOUNDFORTPLUSRALONGINTEGRALCURVESOFL}--\eqref{E:BOUNDFORTMINUSRALONGINTEGRALCURVESOFL}.

\medskip
\noindent \textbf{Proof of \eqref{E:GLOBALEXISTENCENEARTIMEAXISINTEGRALESTIMATEROVER1PLUSTPLUSRPOWERAALONGINTEGRALCURVESOFL}--\eqref{E:GLOBALEXISTENCENEARTIMEAXISINTEGRALESTIMATEROVER1PLUSTPLUSRPOWERATIMESONEPLUSTPLUSRALONGINTEGRALCURVESOFL}}:
We will prove \eqref{E:GLOBALEXISTENCENEARTIMEAXISINTEGRALESTIMATEROVER1PLUSTPLUSRPOWERAALONGINTEGRALCURVESOFL}
only in the case $r_1 \leq 1 + t_0$ since the case
$t_1 \leq 1 + r_0$ can be handled using similar arguments.
To proceed, we use	
\eqref{E:GLOBALEXISTENCEBOUNDFORTMINUSRALONGINTEGRALCURVESOFL}
and our assumption that $r_1 \leq 1 + t_0$ 
to deduce that for all $t \in [t_0,t_1]$, we have
$1 + t + \mathfrak{r}(t) \approx 1 + t_0 \approx 1 + t_1 + r_1$.
This allows us to bound LHS~\eqref{E:GLOBALEXISTENCENEARTIMEAXISINTEGRALESTIMATEROVER1PLUSTPLUSRPOWERAALONGINTEGRALCURVESOFL} by:
\begin{align} \label{E:FIRSTPROOFSTEPGLOBALEXISTENCENEARTIMEAXISINTEGRALESTIMATEROVER1PLUSTPLUSRPOWERAALONGINTEGRALCURVESOFL}
	&
	\lesssim
	\frac{1}{(1 + t_1 + r_1)^B} 
	\int_{t_0}^{t_1} \mathfrak{r}^A(t) 
	\, \mathrm{d} t.
\end{align}
Next, we note that in the present context, the estimate \eqref{E:LUNITRISONEPLUSSMALLERROR} 
holds for the same reasons as before.
We can therefore change variables to $r = \mathfrak{r}(t)$ 
in the integral in \eqref{E:FIRSTPROOFSTEPGLOBALEXISTENCENEARTIMEAXISINTEGRALESTIMATEROVER1PLUSTPLUSRPOWERAALONGINTEGRALCURVESOFL},
using \eqref{E:LUNITRISONEPLUSSMALLERROR} to deduce that $\mathrm{d} r = (1 + \mathcal{O}(\eps)) \mathrm{d} t$,
thereby bounding RHS~\eqref{E:FIRSTPROOFSTEPGLOBALEXISTENCENEARTIMEAXISINTEGRALESTIMATEROVER1PLUSTPLUSRPOWERAALONGINTEGRALCURVESOFL} by: 
\begin{align} \label{E:SECONDPROOFSTEPNEARTIMEAXISINTEGRALESTIMATEROVER1PLUSTPLUSRPOWERBALONGINTEGRALCURVESOFL}
	&
	\lesssim
	\frac{1}{(1 + t_1 + r_1)^B} 
	\int_{r_0}^{r_1} r^A
	\, \mathrm{d} r
	  \lesssim 
		\frac{r_1^{A+1}}{(1 + t_1 + r_1)^B}.
\end{align}
We have therefore proved \eqref{E:GLOBALEXISTENCENEARTIMEAXISINTEGRALESTIMATEROVER1PLUSTPLUSRPOWERAALONGINTEGRALCURVESOFL}.

\eqref{E:GLOBALEXISTENCENEARTIMEAXISINTEGRALESITMATELOGTPLUSRTIMESRTOAOVERONEPLUSTPLUSRTOBALONGINTEGRALCURVESOFL} 
follows from \eqref{E:GLOBALEXISTENCENEARTIMEAXISINTEGRALESTIMATEROVER1PLUSTPLUSRPOWERAALONGINTEGRALCURVESOFL} and the simple fact that
by \eqref{E:GLOBALEXISTENCEBOUNDFORTPLUSRALONGINTEGRALCURVESOFL}, 
the integrand factor $\ln_+(t + \mathfrak{r}(t))$ on the LHS increases along the future-directed integral curves of $\Lunit$.

\eqref{E:GLOBALEXISTENCENEARTIMEAXISINTEGRALESTIMATEROVER1PLUSTPLUSRPOWERATIMESONEPLUSTPLUSRALONGINTEGRALCURVESOFL}
follows from \eqref{E:GLOBALEXISTENCENEARTIMEAXISINTEGRALESTIMATEROVER1PLUSTPLUSRPOWERAALONGINTEGRALCURVESOFL} and the simple fact that
by \eqref{E:GLOBALEXISTENCEBOUNDFORTPLUSRALONGINTEGRALCURVESOFL}, 
the integrand factor $\frac{1}{1 + \mathfrak{r}(t)}$ on the LHS decreases along the future-directed integral curves of $\Lunit$.

\medskip
\noindent \textbf{Proof of \eqref{E:GLOBALEXISTENCEFARFROMTIMEAXISINTEGRALESTIMATEROVER1PLUSTPLUSRSQUAREDALONGINTEGRALCURVESOFL}--\eqref{E:GLOBALEXISTENCEFARFROMTIMEAXISLOGTPLUSROVERONEPLUSTPLUSRALONGINTEGRALCURVESOFL}}:
To prove \eqref{E:GLOBALEXISTENCEFARFROMTIMEAXISONEOVERONEPLUSTPLUSRALONGINTEGRALCURVESOFL},
we use the change of variables $\mathrm{d} r = (1 + \mathcal{O}(\eps)) \mathrm{d} t$ from above to bound
LHS~\eqref{E:GLOBALEXISTENCEFARFROMTIMEAXISONEOVERONEPLUSTPLUSRALONGINTEGRALCURVESOFL} by:
\begin{align} \label{E:FIRSTPROOFSTEPGLOBALEXISTENCEFARFROMTIMEAXISONEOVERONEPLUSTPLUSRALONGINTEGRALCURVESOFL}
	&
	\lesssim
	\int_{r_0}^{r_1}
	\frac{1}{1 + r} 
	\, \mathrm{d} r
	= \ln\left( \frac{1 + r_1}{1+r_0} \right),
\end{align}
as is desired.

The estimates
\eqref{E:GLOBALEXISTENCEFARFROMTIMEAXISINTEGRALESTIMATEROVER1PLUSTPLUSRSQUAREDALONGINTEGRALCURVESOFL} and
\eqref{E:GLOBALEXISTENCEFARFROMTIMEAXISLOGTPLUSROVERONEPLUSTPLUSRALONGINTEGRALCURVESOFL}
can be proved via similar arguments, and we omit the details.

\medskip

\textbf{Proof of \eqref{E:GLOBALEXISTENCEFARFROMTIMEAXISONEOVERONEPLUSTPLUSRTOPOWERONEPLUSDELTAALONGINTEGRALCURVESOFL}--\eqref{E:GLOBALEXISTENCEFARFROMTIMEAXISLOGTPLUSROVERONEPLUSTPLUSRTOPOWERONEPLUSDELTAALONGINTEGRALCURVESOFL}}:
To prove \eqref{E:GLOBALEXISTENCEFARFROMTIMEAXISONEOVERONEPLUSTPLUSRTOPOWERONEPLUSDELTAALONGINTEGRALCURVESOFL},
we make the change of variables $z := t + \mathfrak{r}(t)$. Then using \eqref{E:NULLVECTORFIELDS}, 
\eqref{E:SPEEDOFSOUNDEXPANSION}--\eqref{E:SPEEDOFSOUNDERRORFUNCTIONVANISHESATORIGIN}, and the bootstrap assumptions,
we see that $\Lunit t = 1$ and
$\Lunit r = v^r + \Speed = 1 + \mathcal{O}(\eps)$ and hence
$\mathrm{d} z = 2 (1 + \mathcal{O}(\eps)) \mathrm{d} t$. Setting $z_1 := t_1 + r_1$ and $z_2 := t_2 + r_2$, we change variables in the integral to bound it by:
\begin{align} \label{E:FIRSTPROOFSTEPGLOBALEXISTENCEFARFROMTIMEAXISONEOVERONEPLUSTPLUSRTOPOWERONEPLUSDELTAALONGINTEGRALCURVESOFL}
	& \lesssim
		\int_{z_1}^{z_2} \frac{1}{(1 + z)^{1 + \updelta}} \, \mathrm{d} z
		\lesssim \frac{1}{(1 + z_1)^{\updelta}}
		= \frac{1}{(1 + t_1 + r_1)^{\updelta}},
\end{align}
as is desired.

The estimate \eqref{E:GLOBALEXISTENCEFARFROMTIMEAXISLOGTPLUSROVERONEPLUSTPLUSRTOPOWERONEPLUSDELTAALONGINTEGRALCURVESOFL}
can be proved via a similar argument, and we omit the details.
\end{proof}

\subsection{The main a priori estimates for the Riemann invariants}
\label{SS:GLOBALEXISTENCEMAINAPRIORIFORRIEMANNINVARIANTS}
In the next proposition, we derive our main a priori estimates for the Riemann invariants.
The proof relies on a combination of methods used in the proofs of 
Props.\,\ref{P:APRIORIESTIMATESININTERIORREGION} and \ref{P:APRIORIESTIMATESININTERIORREGION}, 
but modified to account for the fact that $\upmu$ grows logarithmically\footnote{In the case of the Chaplygin gas equation of state, 
the constant $\lifespanconstant$ from \eqref{E:LIFESPANCONSTANT} vanishes, and as a consequence,
$\upmu$ stays close to $1$ in the entire subset of spacetime where it is defined; see \eqref{E:PROOFSTEPGLOBALEXISTENCEMUPOINTWISESTIMATE}.} 
in time and to account for the fact that we are starting the analysis from time $0$, 
rather than time $\tstar$.

\begin{proposition}[A priori estimates and strict improvement of the bootstrap assumptions for the Riemann invariants]
\label{P:GLOBALEXISTENCERIEMANNINVARIANTSAPRIORIEXTERIORREGIONESTIMATES}
Under the data assumptions of Section~\ref{SS:GLOBALEXISTENCEDATAATTIME0} and the bootstrap assumptions of
Sect.\,\ref{SS:GLOBALEXISTENCEBOOTSTRAPANDTWOREGIONS}, 
if $\datasize$ is sufficiently small, then the following estimates hold.

\medskip

\noindent \underline{\textbf{Estimates in the away-from-time-axis region $\lbrace r \geq 1 \rbrace$}}.
The following estimates hold with respect to the $(t,u)$ coordinates in the region $\lbrace r \geq 1 \rbrace \cap \lbrace 0 \leq t < \Tboot \rbrace$:

\medskip 
\begin{subequations}
\begin{align} \label{E:GLOBALEXISTENCEFARFROMTIMEAXISPOINTWISEAPRIORIRPLUS}
	|\RRiemann|
	& 
	\leq C \frac{\datasize}{(1 + t + |u|)},
		\\
|\partial_t \RRiemann|
		& 
	\leq C \frac{\datasize}{(1 + t + |u|)(1 + |u|)},
		\label{E:GLOBALEXISTENCEFARFROMTIMEAXISPOINTWISEAPRIORIPARTIALTRPLUS}
			\\
	|\muuLunit \RRiemann|
	& 
	\leq C \frac{\datasize}{(1 + t + |u|)(1 + |u|)},
		\label{E:GLOBALEXISTENCEFARFROMTIMEAXISPOINTWISEAPRIORIMULBARRPLUS}
		\\
|r \muuLunit \RRiemann|
	& 
	\leq
	C \frac{\datasize}{1 + t + |u|} 
	+
	C \frac{\datasize}{1 + u^2},
		\label{E:GLOBALEXISTENCEFARFROMTIMEAXISPOINTWISEAPRIORIRTIMESMUULUNITRPLUS}
		\\
|r \uLunit \RRiemann|
	& 
	\leq
	C \frac{\datasize}{1 + t + |u|} 
	+
	C \frac{\datasize}{1 + u^2},
		\label{E:GLOBALEXISTENCEFARFROMTIMEAXISPOINTWISEAPRIORIRTIMESULUNITRPLUS}
			\\
|\Lunit \RRiemann|
	& 
	\leq 
		C \frac{\datasize}{(1 + t + |u|)^2},
		\label{E:GLOBALEXISTENCEFARFROMTIMEAXISPOINTWISEAPRIORILRPLUS}
			\\
|r \Lunit \RRiemann|
	& 
	\leq 
		C \frac{\datasize}{1 + t + |u|},
		\label{E:GLOBALEXISTENCEFARFROMTIMEAXISPOINTWISEAPRIORIRTIMESLRPLUS}
			\\
|r \muuLunit (\upmu \partial_t \RRiemann)|
		& 
	\leq  
		\datasize 
		\frac{1}{(1 + t + |u|)(1 + |u|)}
		+
		\datasize \frac{1}{(1 + |u|)^3}
		+
		\datasize^{1^+} \frac{\ln_+(|u|)}{(1 + |u|)^3},
		\label{E:GLOBALEXISTENCEFARFROMTIMEAXISPOINTWISEAPRIORIRTIMESMULBARMUPARTIALTRPLUS}
			\\
|r \Lunit \muuLunit \RRiemann|
	& 
	\leq 
	C \frac{\datasize}{(1 + t + |u|)(1 + |u|)},
	\label{E:GLOBALEXISTENCEFARFROMTIMEAXISPOINTWISEAPRIORITIMESLUNITMUULUNITRPLUS}
		\\
|r \muuLunit \muuLunit \RRiemann|
	& 
	\leq 
		C
		\datasize 
		\frac{1}{(1 + t + |u|)(1 + |u|)}
		+
		C
		\datasize \frac{1}{(1 + |u|)^3}
		+
		C
		\datasize^{1^+} \frac{\ln_+(|u|)}{(1 + |u|)^3},
	\label{E:GLOBALEXISTENCEFARFROMTIMEAXISPOINTWISEAPRIORITIMESMUULUNITMUULUNITRPLUS}
	\\
|r \Lunit \Lunit \RRiemann|
	& 
	\leq 
	C \frac{\datasize}{(1 + t + |u|)^2},
	\label{E:GLOBALEXISTENCEFARFROMTIMEAXISPOINTWISEAPRIORITIMESLUNITLUNITRPLUS}
\end{align}
\end{subequations}

\begin{subequations}
\begin{align} \label{E:GLOBALEXISTENCEFARFROMTIMEAXISPOINTWISEAPRIORIRMINUS}
	|\LRiemann|
	& 
	\leq C \frac{\datasize}{(1 + t + |u|)},
		\\
|\partial_t \LRiemann|
		& 
	\leq C \frac{\datasize}{(1 + t + |u|)^2},
		\label{E:GLOBALEXISTENCEFARFROMTIMEAXISPOINTWISEAPRIORIPARTIALTRMINUS}
			\\
|\uLunit \LRiemann|
		& 
	\leq C \frac{\datasize}{(1 + t + |u|)^2},
		\label{E:GLOBALEXISTENCEFARFROMTIMEAXISPOINTWISEAPRIORILBARRMINUS}
		\\
	|r \uLunit \LRiemann|
	& 
	\leq C \frac{\datasize}{1 + t + |u|},
		\label{E:GLOBALEXISTENCEFARFROMTIMEAXISPOINTWISEAPRIORIRTIMESLBARRMINUS}
		\\
|\Lunit \LRiemann|
	& 
	\leq 
		C \datasize \frac{1}{(1 + t + |u|)^2},
		\label{E:GLOBALEXISTENCEFARFROMTIMEAXISPOINTWISEAPRIORILRMINUS}
			\\
|r \Lunit \LRiemann|
	& 
	\leq 
		C \datasize \frac{1}{1 + t + |u|},
		\label{E:GLOBALEXISTENCEFARFROMTIMEAXISPOINTWISEAPRIORIRTIMESLRMINUS}
			\\
	|r \Lunit \Lunit \LRiemann|
	& 
	\leq 
	C \datasize \frac{1}{(1 + t + |u|)^2},
	\label{E:GLOBALEXISTENCEFARFROMTIMEAXISPOINTWISEAPRIORIRTIMESLUNITLUNITRMINUS}
		\\
	|r \Lunit \partial_t \LRiemann|
		& 
	\leq 
	C \frac{\datasize}{(1 + t + |u|)^2},
		\label{E:GLOBALEXISTENCEFARFROMTIMEAXISPOINTWISEAPRIORIRTIMESLUNITPARTIALTRMINUS}
			\\
	|r \muuLunit \Lunit \LRiemann|
	& 
	\leq 
	C \frac{\datasize}{(1 + t + |u|)^2}
	+
	C \frac{\datasize \ln_+(t-u)}{(1 + t + |u|)^2(1 + |u|)^2},
	\label{E:GLOBALEXISTENCEFARFROMTIMEAXISPOINTWISEAPRIORIRTIMESMUULUNITLUNITRMINUS}
		\\
	|r \muuLunit \uLunit \LRiemann|
		& 
	\leq 
	C \frac{\datasize}{(1 + t + |u|)(1 + |u|)}.
		\label{E:GLOBALEXISTENCEFARFROMTIMEAXISPOINTWISEAPRIORIRTIMESMUULUNITULUNITRMINUS}
\end{align}
\end{subequations}

\medskip

\noindent \underline{\textbf{Estimates in the near-time-axis region $\lbrace 0 \leq r < 1 \rbrace$}}.
The following estimates hold with respect to the $(t,r)$ coordinates on the region $\lbrace 0 \leq r < 1 \rbrace \cap \lbrace 0 \leq t < \Tboot \rbrace$:

\begin{subequations}
\begin{align} \label{E:GLOBALEXISTENCENEARTIMEAXISPOINTWISEAPRIORIRPLUS}
	|\RRiemann|
	& 
	\leq C \frac{\datasize}{1 + t + r},
		\\
|\partial_t \RRiemann|
		& 
	\leq C \frac{\datasize}{(1 + t + r)(1 + |t-r|)},
		\label{E:GLOBALEXISTENCENEARTIMEAXISPOINTWISEAPRIORIPARTIALTRPLUS}
			\\
	|\uLunit \RRiemann|
	& 
	\leq C \frac{\datasize}{(1 + t + r)(1 + |t-r|)},
		\label{E:GLOBALEXISTENCENEARTIMEAXISPOINTWISEAPRIORILBARRPLUS}
		\\
|r \uLunit \RRiemann|
	& 
	\leq
	C
	\frac{\datasize}{1 + t + r}
	+
	C \frac{\datasize}{(1 + |t - r|)^2},
		\label{E:GLOBALEXISTENCENEARTIMEAXISPOINTWISEAPRIORIRTIMESLBARRPLUS}
			\\
|\Lunit \RRiemann|
	& 
	\leq 
		C \frac{\datasize}{(1 + t + r)^2},
		\label{E:GLOBALEXISTENCENEARTIMEAXISPOINTWISEAPRIORILRPLUS}
			\\
|r \Lunit \RRiemann|
	& 
	\leq 
		C \frac{\datasize}{1 + t + r},
		\label{E:GLOBALEXISTENCENEARTIMEAXISPOINTWISEAPRIORIRTIMESLRPLUS}
			\\
|r \uLunit \partial_t \RRiemann|
		& 
	\leq  
		C \frac{\datasize}{(1 + t + r)(1 + |t-r|)},
		\label{E:GLOBALEXISTENCENEARTIMEAXISPOINTWISEAPRIORIRTIMESLBARPARTIALTRPLUS}
			\\
|r \Lunit \uLunit \RRiemann|
	& 
	\leq 
	C \frac{\datasize}{(1 + t + r)(1 + |t-r|)},
	\label{E:GLOBALEXISTENCENEARTIMEAXISPOINTWISEAPRIORITIMESLUNITLBARRPLUS}
		\\
|r \uLunit \uLunit \RRiemann|
	& 
	\leq 
	C \frac{\datasize}{(1 + t + r)(1 + |t-r|)},
	\label{E:GLOBALEXISTENCENEARTIMEAXISPOINTWISEAPRIORITIMESLBARLBARRPLUS}
	\\
|r \Lunit \Lunit \RRiemann|
	& 
	\leq 
	C \frac{\datasize}{(1 + t + r)^2},
	\label{E:GLOBALEXISTENCENEARTIMEAXISPOINTWISEAPRIORITIMESLUNITLUNITRPLUS}
\end{align}
\end{subequations}

\begin{subequations}
\begin{align} \label{E:GLOBALEXISTENCENEARTIMEAXISPOINTWISEAPRIORIRMINUS}
	|\LRiemann|
	& 
	\leq C \frac{\datasize}{1 + t + r},
		\\
|\partial_t \LRiemann|
		& 
	\leq 	
		C \frac{\datasize \ln(t+r)}{(1 + t + r)^3}	
		+
		C \frac{\datasize^{1^+} [\ln(t+r)]^2}{(1 + t + r)^3},
		\label{E:GLOBALEXISTENCENEARTIMEAXISPOINTWISEAPRIORIPARTIALTRMINUS}
			\\
|\uLunit \LRiemann|
		& 
	\leq C \frac{\datasize}{(1 + t + r)^2},
		\label{E:GLOBALEXISTENCENEARTIMEAXISPOINTWISEAPRIORILBARRMINUS}
		\\
	|r \uLunit \LRiemann|
	& 
	\leq C \frac{\datasize}{1 + t + r},
		\label{E:GLOBALEXISTENCENEARTIMEAXISPOINTWISEAPRIORIRTIMESLBARRMINUS}
		\\
|\Lunit \LRiemann|
	& 
	\leq 
		C
		\datasize
		\frac{1}{(1 + t + r)^3}	
		+
		C
		\datasize^{1^+} 
		\frac{\ln_+(t+r)}{(1 + t + r)^3},
		\label{E:GLOBALEXISTENCENEARTIMEAXISPOINTWISEAPRIORILRMINUS}
			\\
|r \Lunit \LRiemann|
	& 
	\leq 
		C \frac{\datasize \ln(t+r)}{(1 + t + r)^3}	
		+
		C \frac{\datasize^{1^+} [\ln(t+r)]^2}{(1 + t + r)^3},
		\label{E:GLOBALEXISTENCENEARTIMEAXISPOINTWISEAPRIORIRTIMESLRMINUS}
			\\
	|r \Lunit \Lunit \LRiemann|
	& 
	\leq 
		C
		\datasize
		\frac{1}{(1 + t + r)^3}	
		+
		C		
		\datasize^{1^+} 
		\frac{\ln_+(t+r)}{(1 + t + r)^3},
	\label{E:GLOBALEXISTENCENEARTIMEAXISPOINTWISEAPRIORIRTIMESLUNITLUNITRMINUS}
		\\
	|r \Lunit \partial_t \LRiemann|
		& 
	\leq 
		C \frac{\datasize \ln(t+r)}{(1 + t + r)^3}	
		+
		C \frac{\datasize^{1^+} [\ln(t+r)]^2}{(1 + t + r)^3},
		\label{E:GLOBALEXISTENCENEARTIMEAXISPOINTWISEAPRIORIRTIMESLUNITPARTIALTRMINUS}
			\\
	|r \uLunit \Lunit \LRiemann|
	& 
	\leq 
	C \frac{\datasize}{(1 + t + r)^2},
	\label{E:GLOBALEXISTENCENEARTIMEAXISPOINTWISEAPRIORIRTIMESLBARLUNITRMINUS}
		\\
	|r \uLunit \uLunit \LRiemann|
		& 
	\leq 
	C \frac{\datasize}{(1 + t + r)(1 + |t-r|)}.
		\label{E:GLOBALEXISTENCENEARTIMEAXISPOINTWISEAPRIORIRTIMESLBARLBARRMINUS}
\end{align}
\end{subequations}


\begin{remark}\label{Remark:improveboundsR+-:global}
 As in Remarks~\ref{Remark:improveboundsR+-} and \ref{Remark:improveboundsR+-:interior},
some decay rates in Prop.\,\ref{P:GLOBALEXISTENCERIEMANNINVARIANTSAPRIORIEXTERIORREGIONESTIMATES} 
could be improved, but we don't need such improvements to close our bootstrap argument.
\end{remark}

\end{proposition}

\begin{proof}
Throughout the proof, we silently use the estimates 
\eqref{E:GLOBALEXISTENCESHARPPOINWISECOMPARISONBETWEENUANDTMINUSR}
and 
\eqref{E:GLOBALEXISTENCEUSMALLERTHANTMINUS1}--\eqref{E:GLOBALEXISTENCERISANALMOSTUNITYMULIPLEOFTMINUSU}.
We will also silently use \eqref{E:SPEEDOFSOUNDEXPANSION}--\eqref{E:SPEEDOFSOUNDERRORFUNCTIONVANISHESATORIGIN},
which allow us to replace $\Speed$ with $1$ in various estimates, up to harmless error terms

\medskip
\noindent \textbf{Proof of \eqref{E:GLOBALEXISTENCEFARFROMTIMEAXISPOINTWISEAPRIORIRMINUS}, 
\eqref{E:GLOBALEXISTENCEFARFROMTIMEAXISPOINTWISEAPRIORIRTIMESLRMINUS},
\eqref{E:GLOBALEXISTENCENEARTIMEAXISPOINTWISEAPRIORIRMINUS},
and \eqref{E:GLOBALEXISTENCENEARTIMEAXISPOINTWISEAPRIORIRTIMESLRMINUS}}:
We first use the wave equation \eqref{E:REVAMPEDRMINUSWAVEEQUATION},
and the transport equation \eqref{E:MUEVOLUTION},
carefully noting that the products  
$2 \LRiemann \muuLunit \Speed$ cancel from each side of the equation,
and 
use \eqref{E:SPEEDOFSOUNDEXPANSION}--\eqref{E:SPEEDOFSOUNDERRORFUNCTIONVANISHESATORIGIN} and
the bootstrap assumptions, thereby deducing that
	the following preliminary estimate holds in the away-from-time-axis region 
	$\lbrace (t,u ) \ | \ 0 \leq t < \Tboot, u \in (-\infty,\Datafunctionforeikonal(t)] \rbrace$,
	where $\Datafunctionforeikonal$ is as in \eqref{E:TIMELILKEEIKONALDATAISAPERTURBATIONOFTMINUS1}:
	\begin{align} \label{E:GLOBALEXISTENCEWAVEEQUATIONRMINUSINHOMOGENEOUSTERMBOUND}
	\left|
	\muuLunit
	\left\lbrace
		r \Lunit \LRiemann
		+
		2 \LRiemann
	\right\rbrace
	\right|
	& \lesssim \frac{\eps^2}{(1 + t + |u|)^3}
		+
		\frac{\eps^2 }{(1 + t + |u|)^2(1 + |u|^2)}.
	\end{align}
Similarly, using the wave equation \eqref{E:LEFTRIEMANNEQUATIONINTERIORREGION}
and the bootstrap assumptions,
we further deduce that 
the following bound holds in the near-time-axis region $(t,r) \in [0,\Tboot) \times [0,1]$:
\begin{align} \label{E:NEARTIMEAXISGLOBALEXISTENCEWAVEEQUATIONRMINUSINHOMOGENEOUSTERMBOUND}
	\left|
	\uLunit
	\left\lbrace
		r \Lunit \LRiemann
		+
		2 \LRiemann
	\right\rbrace
	\right|
	& \lesssim \frac{\eps^2}{(1 + t + r)^3}.
	\end{align}

	Next, we integrate \eqref{E:GLOBALEXISTENCEWAVEEQUATIONRMINUSINHOMOGENEOUSTERMBOUND}
	along the integral curve of $\muuLunit$ that connects $(t,u)$ to 
	the data point
	$(0,u_0)$
	and repeat the proof of \eqref{E:EXTERIORREGIONPOINTWISEESTIMATEFORALLLUNITDERIVATIVESMODIFIEDLRMINUSDERIVATIVE},
	taking into account the integral estimate \eqref{E:GLOBALEXISTENCEINTEGRALESTIMATEALONGMUULUNITONEOVERTPLUSUTOPOWERATIMSINTEGRABLEINU},
	\eqref{E:GLOBALEXISTENCETWOTMINUSUALONGINTEGRALCURVESOFULUNIT},
	and the data estimate \eqref{E:GLOBALEXISTENCETIME0TRCOORDINATESPURELDERIVATIVESTRANSPORTEDMODIFIEDLUNITRMINUSDATA} (with $J=0$),
	thereby deducing that the following estimate holds
	in the away-from-time-axis region 
	$\lbrace (t,u) \ | \ 0 \leq t < \Tboot, u \in (-\infty,\Datafunctionforeikonal(t)] \rbrace$:
	\begin{align} \label{E:GLOBALEXISTENCEEXTERIORREGIONPOINTWISEESTIMATEFORMODIFIEDLRMINUSDERIVATIVE}
	\begin{split}
		[r \Lunit \LRiemann
		+
		2 \LRiemann]
		& 
		=
		2 \datasize
		\frac{1}{1 + (u_0)^2} 
		+
		\mathcal{O}\left(\frac{\datasize^{1^+}}{(1 + t + |u|)^2} \right)
			\\
	& =
		2 
		\left\lbrace
			\datasize 
			+
		\mathcal{O}(\datasize^{1^+})
		\right\rbrace
		\frac{1}{1 + (2t - u)^2}.
	\end{split}
	\end{align}
Moreover, we can integrate \eqref{E:NEARTIMEAXISGLOBALEXISTENCEWAVEEQUATIONRMINUSINHOMOGENEOUSTERMBOUND} 
and use \eqref{E:GLOBALEXISTENCEINTERIORREGIONPOINTWISEESTIMATECHANGEINTPLUSRALONGINTEGRALCURVESOFLBAR}
and
the integral estimate \eqref{E:GLOBALEXISTENCENEARTIMEAXISINTEGRALESTIMATEONEOVERONEPLUSTPLUSRTOBTIMESONEOVER1PLUSTMINUSRTOA},
as well as the
initial conditions on $\lbrace r=1 \rbrace$ furnished by \eqref{E:GLOBALEXISTENCEEXTERIORREGIONPOINTWISEESTIMATEFORMODIFIEDLRMINUSDERIVATIVE},
\eqref{E:GLOBALEXISTENCETIMELKEEIKONALFUNCTIONDATA}, 
and \eqref{E:TIMELILKEEIKONALDATAISAPERTURBATIONOFTMINUS1}--\eqref{E:C2ESTIMATEFORTIMELILKEEIKONALDATAPERTURBATION},
thereby deducing that the following bound holds in the near-time-axis region $(t,r) \in [0,\Tboot) \times [0,1]$:
\begin{align} \label{E:NEARTIMEAXISREGIONGLOBALEXISTENCEEXTERIORREGIONPOINTWISEESTIMATEFORMODIFIEDLRMINUSDERIVATIVE}
	\begin{split}
		[r \Lunit \LRiemann
		+
		2 \LRiemann]
		& 
		=
		2 
		\frac{\left\lbrace
			\datasize 
			+
		\mathcal{O}(\datasize^{1^+})
		\right\rbrace}{1 + (t + r)^2}.
	\end{split}
	\end{align}
The estimate \eqref{E:NEARTIMEAXISREGIONGLOBALEXISTENCEEXTERIORREGIONPOINTWISEESTIMATEFORMODIFIEDLRMINUSDERIVATIVE} in particular implies
that:
\begin{align} \label{E:GLOBALEXISTENCELRIEMANNATREQUAS0BOUND}
	\LRiemann \restriction_{\lbrace r = 0 \rbrace}
	& = \frac{\datasize  + \mathcal{O}(\datasize^{1^+})}{1 + t^2}.
\end{align}	

We next multiply \eqref{E:NEARTIMEAXISREGIONGLOBALEXISTENCEEXTERIORREGIONPOINTWISEESTIMATEFORMODIFIEDLRMINUSDERIVATIVE} by $r$
and argue as in the proof of \eqref{E:INTERIORREGIONPPOINWISEESTIMATEFORRTIMESLUNITRSQUAREDRMINUS}
to deduce the following pointwise estimate, valid
in the near-time-axis region $(t,r) \in [0,\Tboot) \times [0,1]$:
\begin{align} \label{E:GLOBALEXISTENCEPOINWISEESTIMATENEARTIMEAXISFORRTIMESLUNITRSQUAREDRMINUS}
	\Lunit(r^2 \LRiemann)
	& =
		\left\lbrace
			2 \datasize + \mathcal{O}(\datasize^{1^+})
		\right\rbrace
		\frac{r}{1 + (t + r)^2}.
\end{align}
Similarly, we multiply \eqref{E:GLOBALEXISTENCEEXTERIORREGIONPOINTWISEESTIMATEFORMODIFIEDLRMINUSDERIVATIVE} by $r$
and use \eqref{E:GLOBALEXISTENCESHARPPOINWISECOMPARISONBETWEENUANDTMINUSR} to deduce the following pointwise estimate in the away-from-time-axis region 
	$\lbrace (t,u ) \ | \ 0 \leq t < \Tboot, u \in (-\infty,\Datafunctionforeikonal(t)] \rbrace$:
	\begin{align} \label{E:GLOBALEXISTENCEEXTERIORREGIONPOINTWISEESTIMATEFORLUNITDERIVATIVERSQUAREDTIMERMINUS}
	\begin{split}
		\Lunit(r^2 \LRiemann)
		& =
		2 
		\left\lbrace
			\datasize 
			+
		\mathcal{O}(\datasize^{1^+})
		\right\rbrace
	\frac{r}{1 + (2t - u)^2}
	=
	\left\lbrace
			-2 \datasize + \mathcal{O}(\datasize^{1^+})
		\right\rbrace
		\frac{r}{1 + (t + r)^2}.
	\end{split}
	\end{align}
\eqref{E:GLOBALEXISTENCEPOINWISEESTIMATENEARTIMEAXISFORRTIMESLUNITRSQUAREDRMINUS} and 
\eqref{E:GLOBALEXISTENCEEXTERIORREGIONPOINTWISEESTIMATEFORLUNITDERIVATIVERSQUAREDTIMERMINUS}
together imply that the following pointwise estimate holds in the entire bootstrap region
$(t,r) \in [0,\Tboot) \times [0,\infty)$:
\begin{align} \label{E:GLOBALEXISTENCEPOINWISEESTIMATEWHOLEBOOTSTRAPREGIONFORRTIMESLUNITRSQUAREDRMINUS}
	\Lunit(r^2 \LRiemann)
	& =
		\left\lbrace
			2 \datasize + \mathcal{O}(\datasize^{1^+})
		\right\rbrace
		\frac{r}{1 + (t + r)^2}.
\end{align}

We are now ready to prove \eqref{E:GLOBALEXISTENCEFARFROMTIMEAXISPOINTWISEAPRIORIRMINUS} and \eqref{E:GLOBALEXISTENCENEARTIMEAXISPOINTWISEAPRIORIRMINUS}.
The proof splits into two cases.

\noindent \textit{Case $1$: Estimates along integral curves of $\Lunit$ that emanate from the time axis}.
Let $(t_0,0)$ be a point on the time axis with $t_0 \geq 0$, 
and let $(t_1,r_1)$ be a point lying on the integral curve $t \rightarrow \upgamma_{t_0}(t)$ of $\Lunit$ emanating from
$(t_0,0)$ with $t_1 \geq t_0$. 

\noindent \textit{Sub-case $1a$: $r_1 \leq 1 + t_0$}. 
From \eqref{E:GLOBALEXISTENCEPOINWISEESTIMATEWHOLEBOOTSTRAPREGIONFORRTIMESLUNITRSQUAREDRMINUS},
the fundamental theorem of calculus, 
and the integral estimate \eqref{E:GLOBALEXISTENCENEARTIMEAXISINTEGRALESTIMATEROVER1PLUSTPLUSRPOWERAALONGINTEGRALCURVESOFL},
we deduce that:
\begin{align} \label{E:GLOBALEXISTENCENEARTIMEAXISPOINTWISEESTIMATESFORRSQUAREDTIMESRMINUS}
	|r_1^2 \LRiemann(t_1,r_1)|
	& \lesssim 
		\datasize \frac{r_1^2}{(1 + t_1 + r_1)^2}.
\end{align}
Dividing \eqref{E:GLOBALEXISTENCENEARTIMEAXISPOINTWISEESTIMATESFORRSQUAREDTIMESRMINUS} by $r_1^2$,
we conclude that:
\begin{align} \label{E:GLOBALEXISENCENEARTIMEAXISPOINTWISEESTIMATESFORRMINUS}
	|\LRiemann(t_1,r_1)|
	& \lesssim 
		 \frac{\datasize}{(1 + t_1 + r_1)^2},
\end{align}
which in particular implies a stronger estimate than 
\eqref{E:GLOBALEXISTENCEFARFROMTIMEAXISPOINTWISEAPRIORIRMINUS} and \eqref{E:GLOBALEXISTENCENEARTIMEAXISPOINTWISEAPRIORIRMINUS}.

\noindent \textit{Sub-case $1b$: $r_2 \geq 1 + t_0$}. 
We now consider the remaining sub-case, in which the point $(t_2,r_2)$ on $\upgamma_{t_0}$ 
satisfies $r_2 \geq 1 + t_0$. We fix $(t_1,r_1)$ to be the point on $\upgamma_{t_0}$ in which $r_1 = 1 + t_0$; the point $(t_1,r_1)$ was handled in Sub-case $1a$.
Using \eqref{E:GLOBALEXISTENCEBOUNDFORTMINUSRALONGINTEGRALCURVESOFL} 
with $(t_2,r_2)$ in the role of $(t_1,r_1)$ and $0$ in the role of $r_0$, we deduce that
$.9 t_2 \leq t_0 + r_2 \leq 1.1 t_2$. Also using our assumption $r_2 \geq 1 + t_0$ and the simple
bound $t_2 > t_0$ (which follows from \eqref{E:GLOBALEXISTENCEBOUNDFORTPLUSRALONGINTEGRALCURVESOFL}
with $(t_2,r_2)$ in the role of $(t_1,r_1)$ and $0$ in the role of $r_0$),
we deduce that:
\begin{align} \label{E:GLOBALEXISTENCETANDRCOMPARABLEININTERIORREGIONAWAYFROMTIMEAXIS}
	r_2 
	& \approx
	1 + t_0 + r_2
	\approx
	1 + t_2 + r_2.
\end{align} 
To proceed, we integrate \eqref{E:GLOBALEXISTENCEPOINWISEESTIMATEWHOLEBOOTSTRAPREGIONFORRTIMESLUNITRSQUAREDRMINUS}
along $\upgamma_{t_0}$ from time $t_1$ to time $t_2$
and use the already proven bound \eqref{E:GLOBALEXISTENCENEARTIMEAXISPOINTWISEESTIMATESFORRSQUAREDTIMESRMINUS} at $(t_1,r_1)$,
the fundamental theorem of calculus, 
and the integral estimate \eqref{E:GLOBALEXISTENCEFARFROMTIMEAXISINTEGRALESTIMATEROVER1PLUSTPLUSRSQUAREDALONGINTEGRALCURVESOFL},
thereby deducing that:
\begin{align} \label{E:GLOBALEXISTENCEFARFROMTIMEAXISPOINTWISEESTIMATESFORRSQUAREDTIMESRMINUS}
	\begin{split}
	|r_2^2 \LRiemann(t_2,r_2)|
	&
	\leq
	|r_1^2 \LRiemann(t_1,r_1)|
	+
		C
		\datasize
		\ln \left( \frac{1 + r_2}{1 + r_1} \right)
		\\
	& \lesssim 
		\datasize
		\ln(1 + r_2).
\end{split}
\end{align}
Dividing \eqref{E:GLOBALEXISTENCEFARFROMTIMEAXISPOINTWISEESTIMATESFORRSQUAREDTIMESRMINUS} by $r_2^2$
and using \eqref{E:GLOBALEXISTENCETANDRCOMPARABLEININTERIORREGIONAWAYFROMTIMEAXIS},
we conclude that:
\begin{align} \label{E:GLOBALEXISTENCEFARFROMTIMEAXISPOINTWISEESTIMATESFORRMINUS}
	\begin{split}
	|\LRiemann(t_2,r_2)|
	& \lesssim 
		\datasize
		\frac{\ln(1 + r_2)}{(1 + t_2 + r_2)^2},
\end{split}
\end{align}
which in particular implies a stronger estimate than \eqref{E:GLOBALEXISTENCEFARFROMTIMEAXISPOINTWISEAPRIORIRMINUS} 
and \eqref{E:GLOBALEXISTENCENEARTIMEAXISPOINTWISEAPRIORIRMINUS}.
We have therefore proved \eqref{E:GLOBALEXISTENCEFARFROMTIMEAXISPOINTWISEAPRIORIRMINUS} and \eqref{E:GLOBALEXISTENCENEARTIMEAXISPOINTWISEAPRIORIRMINUS} in Case $1$.

\medskip

\noindent \textit{Case 2: Estimates along integral curves of $\Lunit$ that emanate from $\Sigma_0$}.
Let $(0,r_0) \in \Sigma_0$,
and let $(t_1,r_1)$ be a point lying on the integral curve $t \rightarrow \upgamma_{r_0}(t)$ of $\Lunit$ emanating from
$(0,r_0)$ with $t_1 \geq 0$. 

\noindent \textit{Sub-case 2a: $t_1 \leq 1 + r_0$}. 
We argue as in the proof of \eqref{E:GLOBALEXISTENCENEARTIMEAXISPOINTWISEESTIMATESFORRSQUAREDTIMESRMINUS}, 
but this time we incur a non-vanishing data term $r_0^2 \LRiemann(0,r_0)$.
This yields the following bound:
\begin{align} \label{E:GLOBALEXISTENCENEARSIGMA0POINTWISEESTIMATESFORRSQUAREDTIMESRMINUS}
	|r_1^2 \LRiemann(t_1,r_1)|
	& \lesssim 
		|r_0^2 \LRiemann(0,r_0)|
		+
		\datasize 
		\frac{r_1^2}{(1 + t_1 + r_1)^2}.
\end{align}
Dividing \eqref{E:GLOBALEXISTENCENEARSIGMA0POINTWISEESTIMATESFORRSQUAREDTIMESRMINUS} by $r_1^2$,
we find that:
\begin{align} \label{E:FIRSTVERSIONNEARTIMEAXISNEARSIGMA0POINTWISEESTIMATESFORRMINUS}
	|\LRiemann(t_1,r_1)|
	& \lesssim 
		\frac{r_0^2}{r_1^2} |\LRiemann(0,r_0)| 
		+
		\datasize 
		\frac{1}{(1 + t_1 + r_1)^2}.
\end{align}
Inequalities \eqref{E:GLOBALEXISTENCEBOUNDFORTPLUSRALONGINTEGRALCURVESOFL}--\eqref{E:GLOBALEXISTENCEBOUNDFORTMINUSRALONGINTEGRALCURVESOFL} 
and our assumption that $t_1 \leq 1 + r_0$ collectively imply that:
\begin{subequations}
\begin{align} \label{E:GLOBALEXISTENCENEARSIGMA0R1GREATERTHANR0}
	r_0 & \leq r_1, 
		\\
	1 + r_0 
	& \approx
	1 + r_0 + t_1
	\approx
	1 + r_1 + t_1.
	\label{E:GLOBALEXISTENCE1PLUSR0APPROXIMATELYEQUALTO1PLUSR1PLUST1}
	\end{align}
\end{subequations}
Combining 
\eqref{E:GLOBALEXISTENCENEARSIGMA0R1GREATERTHANR0}--\eqref{E:GLOBALEXISTENCE1PLUSR0APPROXIMATELYEQUALTO1PLUSR1PLUST1}
 with \eqref{E:FIRSTVERSIONNEARTIMEAXISNEARSIGMA0POINTWISEESTIMATESFORRMINUS} and the data bound 
\eqref{E:GLOBALEXISTENCERMINUSPOINTWISEESTIMATEATTIME0},
we conclude that:
\begin{align} 
	\begin{split} \label{E:GLOBALEXISTENCENEARSIGMA0POINTWISEESTIMATESFORRMINUS}
	|\LRiemann(t_1,r_1)|
	& \lesssim
		\datasize \frac{1}{1 + t_1 + r_1},
\end{split}
\end{align}
which yields \eqref{E:GLOBALEXISTENCEFARFROMTIMEAXISPOINTWISEAPRIORIRMINUS} and \eqref{E:GLOBALEXISTENCENEARTIMEAXISPOINTWISEAPRIORIRMINUS} in
this sub-case.

\noindent \textit{Sub-case 2b: $t_2 \geq 1 + r_0$}. 
We now consider the final sub-case in which the point $(t_2,r_2)$ on $\upgamma_{r_0}$ 
satisfies $t_2 \geq 1 + r_0$. We fix $(t_1,r_1)$ to be the point on $\upgamma_{r_0}$ in which $t_1 = 1 + r_0$; 
the point $(t_1,r_1)$ was handled in Sub-case 2a.
The proof of the first inequality in \eqref{E:GLOBALEXISTENCEFARFROMTIMEAXISPOINTWISEESTIMATESFORRSQUAREDTIMESRMINUS} goes through verbatim,
which yields:
\begin{align} \label{E:GLOBALEXISTENCEFARFROMSIGMA0POINTWISEESTIMATESFORRSQUAREDTIMESRMINUS}
	|r_2^2 \LRiemann(t_2,r_2)|
	&
	\leq
	|r_1^2 \LRiemann(t_1,r_1)|
	+
		C
		\datasize
		\ln \left( \frac{1 + r_2}{1 + r_1} \right).
\end{align}
Multiplying \eqref{E:GLOBALEXISTENCENEARSIGMA0POINTWISEESTIMATESFORRMINUS} by $r_1^2$ and using the resulting bound  
to control the first term on RHS~\eqref{E:GLOBALEXISTENCEFARFROMSIGMA0POINTWISEESTIMATESFORRSQUAREDTIMESRMINUS},
we further deduce that:
\begin{align} \label{E:SECONDVERSIONGLOBALEXISTENCEFARFROMSIGMA0POINTWISEESTIMATESFORRSQUAREDTIMESRMINUS}
	|r_2^2 \LRiemann(t_2,r_2)|
	&
	\lesssim
		\datasize r_1
		+
		\datasize
		\ln \left(\frac{1 + r_2}{1 + r_1} \right).
\end{align}
Using \eqref{E:GLOBALEXISTENCEBOUNDFORTMINUSRALONGINTEGRALCURVESOFL} with $(t_2,r_2)$ in the role of $(t_1,r_1)$ and $0$ in the role of $t_0$, 
and our assumption that $t_2 \geq 1 + r_0$, we deduce that:
\begin{subequations}
\begin{align}  \label{E:GLOBALEXISTENCETANDRCOMPARABLEAWAYFROMSIGMA0}
	r_2
	& \approx
	r_0 + t_2
	\approx
	1 + r_0 + t_2
	\approx
	1 + t_2 + r_2.
\end{align}
Moreover, using \eqref{E:GLOBALEXISTENCEBOUNDFORTMINUSRALONGINTEGRALCURVESOFL},
we deduce the following simple estimate:
\begin{align} 	\label{E:GLOBALEXISTENCER2LARGERTHANR1ALONGINTEGRALCURVEOFL}
	r_1 & \leq r_2.
\end{align}	
\end{subequations}
Using \eqref{E:SECONDVERSIONGLOBALEXISTENCEFARFROMSIGMA0POINTWISEESTIMATESFORRSQUAREDTIMESRMINUS},
\eqref{E:GLOBALEXISTENCETANDRCOMPARABLEAWAYFROMSIGMA0},
and \eqref{E:GLOBALEXISTENCER2LARGERTHANR1ALONGINTEGRALCURVEOFL},
we divide
\eqref{E:SECONDVERSIONGLOBALEXISTENCEFARFROMSIGMA0POINTWISEESTIMATESFORRSQUAREDTIMESRMINUS} by $r_2^2$ to conclude that:
\begin{align} 
\begin{split} \label{E:GLOBALEXISTENCEFARFROMSIGMA0POINTWISEESTIMATESFORRMINUS}
	|\LRiemann(t_2,r_2)|
	&
	\lesssim
		\datasize \frac{r_1}{r_2^2}
		+
		\datasize
		\frac{\ln \left( \frac{1 + r_2}{1 + r_1} \right)}{r_2^2}
			\\
	& \lesssim
	\frac{\datasize}{r_2}
	\lesssim
	\frac{\datasize}{1 + t_2 + r_2}.
\end{split}
\end{align}
We have therefore proved \eqref{E:GLOBALEXISTENCEFARFROMTIMEAXISPOINTWISEAPRIORIRMINUS} and \eqref{E:GLOBALEXISTENCENEARTIMEAXISPOINTWISEAPRIORIRMINUS}  in Case $2$, which finishes the proof of these two estimates.

\eqref{E:GLOBALEXISTENCEFARFROMTIMEAXISPOINTWISEAPRIORIRTIMESLRMINUS} and \eqref{E:GLOBALEXISTENCENEARTIMEAXISPOINTWISEAPRIORIRTIMESLRMINUS}
then follow from 
\eqref{E:GLOBALEXISTENCEFARFROMTIMEAXISPOINTWISEAPRIORIRMINUS}, 
\eqref{E:GLOBALEXISTENCENEARTIMEAXISPOINTWISEAPRIORIRMINUS},
\eqref{E:GLOBALEXISTENCEEXTERIORREGIONPOINTWISEESTIMATEFORMODIFIEDLRMINUSDERIVATIVE},
and
\eqref{E:NEARTIMEAXISREGIONGLOBALEXISTENCEEXTERIORREGIONPOINTWISEESTIMATEFORMODIFIEDLRMINUSDERIVATIVE}.

\medskip

\noindent \textbf{Proof of \eqref{E:GLOBALEXISTENCEFARFROMTIMEAXISPOINTWISEAPRIORILRMINUS}, \eqref{E:GLOBALEXISTENCENEARTIMEAXISPOINTWISEAPRIORILRMINUS}, \eqref{E:GLOBALEXISTENCEFARFROMTIMEAXISPOINTWISEAPRIORIRTIMESLUNITLUNITRMINUS},   
and \eqref{E:GLOBALEXISTENCEFARFROMTIMEAXISPOINTWISEAPRIORIRTIMESLUNITLUNITRMINUS}}:
With the help of Lemma~\ref{L:VECTORFIELDSINTERMSOFGEOMETRICCOORDINATES}, 
\eqref{E:NULLVECTORFIELDS},
and \eqref{E:MUEVOLUTION}, 
we commute the wave equation \eqref{E:REVAMPEDRMINUSWAVEEQUATION} with $\Lunit$
and use the bootstrap assumptions,
carefully noting that the products $2 \LRiemann \Lunit \muuLunit \Speed$
(which are ``dangerous'' in that they involve a second derivative that lacks an $r$-weight)
cancel from both sides of the resulting equation,
thereby deducing that
	the following preliminary estimate holds in the away-from-time-axis region 
	$\lbrace (t,u ) \ | \ 0 \leq t < \Tboot, u \in (-\infty,\Datafunctionforeikonal(t)] \rbrace$:
	\begin{align} \label{E:GLOBALEXISTENCEWAVEEQUATIONRMINUSONELCOMMUTATIONINHOMOGENEOUSTERMBOUND}
	\left|
	\muuLunit
	\Lunit
	\left\lbrace
		r \Lunit \LRiemann
		+
		2 \Speed \LRiemann
	\right\rbrace
	\right|
	& \lesssim \frac{\eps^2}{(1 + t + |u|)^3(1 + |u|)}.
	\end{align}

Similarly, we commute the wave equation \eqref{E:LEFTRIEMANNEQUATIONINTERIORREGION} with $\Lunit$ and use
the bootstrap assumptions to further deduce that 
the following pointwise estimate holds in the near-time-axis region $(t,r) \in [0,\Tboot) \times [0,1]$:
\begin{align} \label{E:NEARTIMEAXISGLOBALEXISTENCEWAVEEQUATIONONELCOMMUTATIONRMINUSINHOMOGENEOUSTERMBOUND}
	\left|
	\uLunit
	\Lunit
	\left\lbrace
		r \Lunit \LRiemann
		+
		2 \LRiemann
	\right\rbrace
	\right|
	& \lesssim \frac{\eps^2}{(1 + t + r)^4}.
	\end{align}

Next, we integrate \eqref{E:GLOBALEXISTENCEWAVEEQUATIONRMINUSONELCOMMUTATIONINHOMOGENEOUSTERMBOUND} 
along the integral curve of $\muuLunit$ that connects $(t,u)$ to 
	the data point
	$(0,u_0)$
	and use the integral estimate \eqref{E:GLOBALEXISTENCEINTEGRALESTIMATEALONGMUULUNITONEOVERTPLUSUTOPOWERATIMESONEOVER1PLUSU},
	the data estimate \eqref{E:GLOBALEXISTENCETIME0TRCOORDINATESPURELDERIVATIVESTRANSPORTEDMODIFIEDLUNITRMINUSDATA} (with $J=1$),
	\eqref{E:GLOBALEXISTENCESPACELIKEEIKONALFUNCTIONDATA},
	and \eqref{E:GLOBALEXISTENCETWOTMINUSUALONGINTEGRALCURVESOFULUNIT}
	to deduce that the following pointwise estimate holds in the away-from-time-axis region 
	$\lbrace (t,u ) \ | \ 0 \leq t < \Tboot, u \in (-\infty,\Datafunctionforeikonal(t)] \rbrace$:
	\begin{align} \label{E:GLOBALEXISTENCEINTEGRATEDWAVEEQUATIONRMINUSONELCOMMUTATIONINHOMOGENEOUSTERMBOUND} 
	\begin{split}
		\Lunit
		[r \Lunit \LRiemann
		+
		2 \LRiemann]
		& 
		=
		-
		4 [\datasize + \mathcal{O}(\datasize^{1^+})]
		\frac{(2t-u_0)}{[1 + (2t - u_0)^2]^2} \restriction_{t=0}
		+
		\mathcal{O}(\datasize^{1^+}) \frac{\ln_+(|u_0| + |u|)}{(1 + t + |u|)^3}
			\\
	& =
		-
		4 \datasize 
		\frac{(2t-u)}{[1 + (2t - u)^2]^2}
		+
		\mathcal{O}(\datasize^{1^+}) \frac{\ln_+(t + |u|)}{(1 + t + |u|)^3}.
	\end{split}
	\end{align}
	Commuting in the outer factor of $\Lunit$ on LHS~\eqref{E:GLOBALEXISTENCEINTEGRATEDWAVEEQUATIONRMINUSONELCOMMUTATIONINHOMOGENEOUSTERMBOUND} 
	and using the identity 
	$\Lunit r = v^r + \Speed$ (see \eqref{E:NULLVECTORFIELDS}),
	\eqref{E:SPEEDOFSOUNDEXPANSION}--\eqref{E:SPEEDOFSOUNDERRORFUNCTIONVANISHESATORIGIN},
	and the bootstrap assumptions,
	we further deduce that
	the following pointwise estimate holds 
	in the away-from-time-axis region 
	$\lbrace (t,u ) \ | \ 0 \leq t < \Tboot, u \in (-\infty,\Datafunctionforeikonal(t)] \rbrace$:
	\begin{align} \label{E:GLOBALEXISTENCEINTEGRATEDWAVEEQUATIONRMINUSONELCOMMUTATIONALLTHEWAYINPOINTWISE} 
	\begin{split}
		[r \Lunit \Lunit \LRiemann
		+
		3 \Lunit \LRiemann]
		& =
		-
		4 \datasize 
		\frac{(2t-u)}{[1 + (2t - u)^2]^2}
		+
		\mathcal{O}(\datasize^{1^+}) \frac{\ln_+(t + |u|)}{(1 + t + |u|)^3}.
	\end{split}
	\end{align}
	Multiplying \eqref{E:GLOBALEXISTENCEINTEGRATEDWAVEEQUATIONRMINUSONELCOMMUTATIONALLTHEWAYINPOINTWISE} by $r^2$ 
	and again using the identity $\Lunit r = v^r + \Speed$, 
	\eqref{E:SPEEDOFSOUNDEXPANSION}--\eqref{E:SPEEDOFSOUNDERRORFUNCTIONVANISHESATORIGIN},
	and the bootstrap assumptions,
	we find that
	the following pointwise estimate holds 
	in the away-from-time-axis region 
	$\lbrace (t,u ) \ | \ 0 \leq t < \Tboot, u \in (-\infty,\Datafunctionforeikonal(t)] \rbrace$:
	\begin{align} \label{E:GLOBALEXISTENCEPOINTWISELRCUBEDRMINUS} 
	\begin{split}
		\Lunit (r^3 \Lunit \LRiemann)
		& =
		-
		4 \datasize 
		\frac{(2t-u)r^2}{[1 + (2t - u)^2]^2}
		+
		\mathcal{O}(\datasize^{1^+}) \frac{\ln_+(t + |u|)r^2}{(1 + t + |u|)^3}.
	\end{split}
	\end{align}

	Similarly, 
	we can integrate \eqref{E:NEARTIMEAXISGLOBALEXISTENCEWAVEEQUATIONONELCOMMUTATIONRMINUSINHOMOGENEOUSTERMBOUND}
	along integral curves of $\uLunit$
and use the integral estimate \eqref{E:GLOBALEXISTENCENEARTIMEAXISINTEGRALESTIMATEONEOVERONEPLUSTPLUSRTOBTIMESONEOVER1PLUSTMINUSRTOA},
\eqref{E:GLOBALEXISTENCEINTERIORREGIONPOINTWISEESTIMATECHANGEINTPLUSRALONGINTEGRALCURVESOFLBAR}, and
the initial conditions on $\lbrace r=1 \rbrace$ furnished by 
\eqref{E:GLOBALEXISTENCEINTEGRATEDWAVEEQUATIONRMINUSONELCOMMUTATIONINHOMOGENEOUSTERMBOUND},
thereby deducing that
the following bound holds in the near-time-axis region $(t,r) \in [0,\Tboot) \times [0,1]$:
\begin{align} \label{E:NEARTIMEAXISGLOBALEXISTENCEINTEGRATEDWAVEEQUATIONRMINUSONELCOMMUTATIONINHOMOGENEOUSTERMBOUND}
	\begin{split}
		\Lunit
		[r \Lunit \LRiemann
		+
		2 \LRiemann]
		& 
		=
		-
		4 
		\datasize 
		\frac{t+r}{[1 + (t + r)^2]^2}
		+
		\mathcal{O}(\datasize^{1^+}) \frac{\ln_+(t + r)}{(1 + t + r)^3}.
	\end{split}
	\end{align}
		Commuting in the outer factor of $\Lunit$ on LHS~\eqref{E:NEARTIMEAXISGLOBALEXISTENCEINTEGRATEDWAVEEQUATIONRMINUSONELCOMMUTATIONINHOMOGENEOUSTERMBOUND}
	and again using the identity $\Lunit r = v^r + \Speed$,
	\eqref{E:SPEEDOFSOUNDEXPANSION}--\eqref{E:SPEEDOFSOUNDERRORFUNCTIONVANISHESATORIGIN},
	and the bootstrap assumptions,
we further deduce that the following bound holds in the near-time-axis region $(t,r) \in [0,\Tboot) \times [0,1]$:
\begin{align} \label{E:NEARTIMEAXISGLOBALEXISTENCEINTEGRATEDWAVEEQUATIONRMINUSONELCOMMUTATIONALLTHEWAYINPOINTWISE}
	\begin{split}
		[r \Lunit \Lunit \LRiemann
		+
		3 \Lunit \LRiemann]
		& 
		=
		-
		4 
		\datasize 
		\frac{t+r}{[1 + (t + r)^2]^2}
		+
		\mathcal{O}(\datasize^{1^+}) \frac{\ln_+(t + r)}{(1 + t + r)^3}.
	\end{split}
	\end{align}
	In particular, \eqref{E:NEARTIMEAXISGLOBALEXISTENCEINTEGRATEDWAVEEQUATIONRMINUSONELCOMMUTATIONALLTHEWAYINPOINTWISE} implies that:
	\begin{align} \label{E:GLOBALEXISTENCELUNITLRIEMANNATREQUAS0BOUND}
	\Lunit \LRiemann \restriction_{\lbrace r = 0 \rbrace}
	& = -
		\frac{4}{3} 
		\datasize 
		\frac{t}{[1 + t^2]^2}
		+
		\mathcal{O}(\datasize^{1^+}) \frac{\ln_+(t)}{(1 + t)^3}.
\end{align}	
	
	Multiplying \eqref{E:NEARTIMEAXISGLOBALEXISTENCEINTEGRATEDWAVEEQUATIONRMINUSONELCOMMUTATIONALLTHEWAYINPOINTWISE} by $r^2$ 
	and again using the identity $\Lunit r = v^r + \Speed$, 
	\eqref{E:SPEEDOFSOUNDEXPANSION}--\eqref{E:SPEEDOFSOUNDERRORFUNCTIONVANISHESATORIGIN},
	and the bootstrap assumptions,
	we find that the following bound holds for $(t,r) \in [0,\Tboot) \times [0,1]$:
\begin{align} \label{E:NEARTIMEAXISGLOBALEXISTENCEPOINTWISELRCUBEDRMINUS}
	\begin{split}
		\Lunit (r^3 \Lunit \LRiemann)
		& 
		=
		-
		4 
		\datasize 
		\frac{(t+r)r^2}{[1 + (t + r)^2]^2}
		+
		\mathcal{O}(\datasize^{1^+}) \frac{\ln_+(t + r)r^2}{(1 + t + r)^3}.
	\end{split}
	\end{align}
	Combining
	\eqref{E:GLOBALEXISTENCEPOINTWISELRCUBEDRMINUS} and
	\eqref{E:NEARTIMEAXISGLOBALEXISTENCEPOINTWISELRCUBEDRMINUS}
	and using \eqref{E:GLOBALEXISTENCESHARPPOINWISECOMPARISONBETWEENUANDTMINUSR},
	we find that the following
	pointwise estimate holds in the entire bootstrap region $(t,r) \in [0,\Tboot) \times [0,\infty)$:
	\begin{align} \label{E:GLOBALEXISTENCEENTIREBOOTSTRAPREGIONPOINTWISELRCUBEDRMINUS}
	\begin{split}
		\Lunit (r^3 \Lunit \LRiemann)
		& 
		=
		-
		4 
		\datasize 
		\frac{(t+r)r^2}{[1 + (t + r)^2]^2}
		+
		\mathcal{O}(\datasize^{1^+}) \frac{\ln_+(t + r)r^2}{(1 + t + r)^3}.
	\end{split}
	\end{align}
	Starting from \eqref{E:GLOBALEXISTENCEENTIREBOOTSTRAPREGIONPOINTWISELRCUBEDRMINUS},
	we can integrate and use arguments similar to the ones we used to prove 
	\eqref{E:GLOBALEXISTENCEFARFROMTIMEAXISPOINTWISEAPRIORIRMINUS} and \eqref{E:GLOBALEXISTENCENEARTIMEAXISPOINTWISEAPRIORIRMINUS},
	based on splitting the argument into Case 1 and Case 2 and using
	the data estimates \eqref{E:GLOBALEXISTENCELUNITORPARTIALTRMINUSPOINTWISEESTIMATEATTIME0} 
	and \eqref{E:GLOBALEXISTENCELUNITLRIEMANNATREQUAS0BOUND}
	and the integral estimates
	\eqref{E:GLOBALEXISTENCENEARTIMEAXISINTEGRALESTIMATEROVER1PLUSTPLUSRPOWERAALONGINTEGRALCURVESOFL}--\eqref{E:GLOBALEXISTENCENEARTIMEAXISINTEGRALESITMATELOGTPLUSRTIMESRTOAOVERONEPLUSTPLUSRTOBALONGINTEGRALCURVESOFL}
	and \eqref{E:GLOBALEXISTENCEFARFROMTIMEAXISONEOVERONEPLUSTPLUSRALONGINTEGRALCURVESOFL}--\eqref{E:GLOBALEXISTENCEFARFROMTIMEAXISLOGTPLUSROVERONEPLUSTPLUSRALONGINTEGRALCURVESOFL},
	thereby concluding that
	the following pointwise estimate holds 
	in the near-time-axis region $(t,r) \in [0,\Tboot) \times [0,1]$:
	\begin{align} \label{E:GLOBALEXISTENCENEARTIMEAXISPOINTWISERCUBEDRMINUS} 
	\begin{split}
		|\Lunit \LRiemann|
		& 
		\lesssim
		\datasize
		\frac{1}{(1 + t + r)^3}	
		+
		\datasize^{1^+} 
		\frac{\ln_+(t+r)}{(1 + t + r)^3},
	\end{split}
	\end{align}
	and that the following pointwise holds  
	in the away-from-time-axis region 
	$\lbrace (t,u ) \ | \ 0 \leq t < \Tboot, u \in (-\infty,\Datafunctionforeikonal(t)] \rbrace$:
	\begin{align} \label{E:GLOBALEXISTENCEFARFROMTIMEAXISPOINTWISERCUBEDRMINUS} 
	\begin{split}
		|\Lunit \LRiemann|
		& 
		\lesssim 
		\frac{\datasize}{(1 + t + r)^2}.
	\end{split}
	\end{align}	
	We have therefore proved
	\eqref{E:GLOBALEXISTENCEFARFROMTIMEAXISPOINTWISEAPRIORILRMINUS}
	and
	\eqref{E:GLOBALEXISTENCENEARTIMEAXISPOINTWISEAPRIORILRMINUS}.

	 \eqref{E:GLOBALEXISTENCEFARFROMTIMEAXISPOINTWISEAPRIORIRTIMESLUNITLUNITRMINUS}
	and \eqref{E:GLOBALEXISTENCENEARTIMEAXISPOINTWISEAPRIORIRTIMESLUNITLUNITRMINUS}
	then follow from
	\eqref{E:GLOBALEXISTENCEFARFROMTIMEAXISPOINTWISEAPRIORILRMINUS},
	\eqref{E:GLOBALEXISTENCENEARTIMEAXISPOINTWISEAPRIORILRMINUS},
	\eqref{E:GLOBALEXISTENCEINTEGRATEDWAVEEQUATIONRMINUSONELCOMMUTATIONALLTHEWAYINPOINTWISE},
	and \eqref{E:NEARTIMEAXISGLOBALEXISTENCEINTEGRATEDWAVEEQUATIONRMINUSONELCOMMUTATIONALLTHEWAYINPOINTWISE}.

\medskip

\noindent \textbf{Proof of \eqref{E:GLOBALEXISTENCEFARFROMTIMEAXISPOINTWISEAPRIORIPARTIALTRMINUS}, \eqref{E:GLOBALEXISTENCENEARTIMEAXISPOINTWISEAPRIORIPARTIALTRMINUS}, \eqref{E:GLOBALEXISTENCEFARFROMTIMEAXISPOINTWISEAPRIORIRTIMESLUNITPARTIALTRMINUS}, and \eqref{E:GLOBALEXISTENCENEARTIMEAXISPOINTWISEAPRIORIRTIMESLUNITPARTIALTRMINUS}}:
	First, we use \eqref{E:PARTIALTINTERMSOFLANDLBAR},
	\eqref{E:GLOBALEXISTENCEWAVEEQUATIONRMINUSINHOMOGENEOUSTERMBOUND},
	\eqref{E:GLOBALEXISTENCEINTEGRATEDWAVEEQUATIONRMINUSONELCOMMUTATIONINHOMOGENEOUSTERMBOUND},
	and the bootstrap assumptions
	to deduce that the following pointwise estimate holds in the away-from-time-axis region 
	$\lbrace (t,u ) \ | \ 0 \leq t < \Tboot, u \in (-\infty,\Datafunctionforeikonal(t)] \rbrace$:
		\begin{align} \label{E:GLOBALEXISTENCEEXTERIORREGIONPOINTWISEESTIMATEEXACTLYONEFLATPARTIALTDERIVATIVEMODIFIEDLRMINUSDERIVATIVE}
	\begin{split}
		\partial_t
		[r \Lunit \LRiemann
		+
		2 \LRiemann]
		& =
		-
		\datasize 
		\frac{2(2t - u)}{[1 + (2t - u)^2]^2}
			\\
		& \ \
		+
		\mathcal{O}(\datasize^{1^+}) \frac{1}{(1 + t + |u|)^2(1+|u|)^2}
		+
		\mathcal{O}(\datasize^{1^+}) \frac{\ln_+(t + |u|)}{(1 + t + |u|)^3}.
	\end{split}
	\end{align}	
	Commuting the operator $\partial_t$ on 
	LHS~\eqref{E:GLOBALEXISTENCEEXTERIORREGIONPOINTWISEESTIMATEEXACTLYONEFLATPARTIALTDERIVATIVEMODIFIEDLRMINUSDERIVATIVE},
	we further deduce, with the help of \eqref{E:PARTIALTINTERMSOFLANDLBAR},
	\eqref{E:COMMUTATOROFLANDULUNIT},
	and
	the bootstrap assumptions, that the following pointwise estimate holds in the away-from-time-axis region 
	$\lbrace (t,u ) \ | \ 0 \leq t < \Tboot, u \in (-\infty,\Datafunctionforeikonal(t)] \rbrace$:
	\begin{align} \label{E:GLOBALEXISTENCEEXTERIORREGIONPOINTWISEESTIMATEEXACTLYONEFLATPARTIALTDERIVATIVECOMMUTEDINSIDEMODIFIEDLRMINUSDERIVATIVE}
	\begin{split}
		[r \Lunit \partial_t \LRiemann
		+
		2 \partial_t \LRiemann]
		& =
		-
		\datasize 
		\frac{2(2t - u)}{[1 + (2t - u)^2]^2}
			\\
		& \ \
		+
		\mathcal{O}(\datasize^{1^+}) \frac{1}{(1 + t + |u|)^2(1+|u|)^2}
		+
		\mathcal{O}(\datasize^{1^+}) \frac{\ln_+(t + |u|)}{(1 + t + |u|)^3}.
	\end{split}
	\end{align}	
	Multiplying \eqref{E:GLOBALEXISTENCEEXTERIORREGIONPOINTWISEESTIMATEEXACTLYONEFLATPARTIALTDERIVATIVECOMMUTEDINSIDEMODIFIEDLRMINUSDERIVATIVE} by $r$
	and using the identity $\Lunit r = v^r + \Speed$, 
	\eqref{E:SPEEDOFSOUNDEXPANSION}--\eqref{E:SPEEDOFSOUNDERRORFUNCTIONVANISHESATORIGIN},
	and the bootstrap assumptions,
	we further deduce that the following pointwise estimate holds in the away-from-time-axis region 
	$\lbrace (t,u ) \ | \ 0 \leq t < \Tboot, u \in (-\infty,\Datafunctionforeikonal(t)] \rbrace$:
	\begin{align} \label{E:GLOBALEXISTENCEEXTERIORREGIONPOINTWISEESTIMATELRSQUAREDPARTIALTRMINUS}
	\begin{split}
		\Lunit (r^2 \partial_t \LRiemann)
		& =
		-
		\datasize 
		\frac{2(2t - u)(t-u)}{[1 + (2t - u)^2]^2}
			\\
		& \ \
		+
		\mathcal{O}(\datasize^{1^+}) \frac{(t-u)}{(1 + t + |u|)^2(1+|u|)^2}
		+
		\mathcal{O}(\datasize^{1^+}) \frac{\ln_+(t + |u|)(t-u)}{(1 + t + |u|)^3}.
	\end{split}
	\end{align}

	Similarly, we use \eqref{E:PARTIALTINTERMSOFLANDLBAR}, 
	\eqref{E:NEARTIMEAXISGLOBALEXISTENCEWAVEEQUATIONRMINUSINHOMOGENEOUSTERMBOUND}, 	
	\eqref{E:GLOBALEXISTENCEEXTERIORREGIONPOINTWISEESTIMATEFORMODIFIEDLRMINUSDERIVATIVE},
	and the bootstrap assumptions to deduce that the following pointwise estimate holds in the near-time-axis region 
	$(t,r) \in [0,\Tboot) \times [0,1]$:
		\begin{align} \label{E:NEARTIMEAXISREGIONGLOBALEXISTENCEEXTERIORREGIONPOINTWISEESTIMATEEXACTLYONEFLATPARTIALTDERIVATIVEMODIFIEDLRMINUSDERIVATIVE}
	\begin{split}
		\partial_t
		[r \Lunit \LRiemann
		+
		2 \LRiemann]
		& =
		-
		2 
		\datasize 
		\frac{t+r}{[1 + (t + r)^2]^2}
		+
		\mathcal{O}(\datasize^{1^+}) \frac{\ln_+(t + r)}{(1 + t + r)^3}.
	\end{split}
	\end{align}	
Commuting the operator $\partial_t$ 
on LHS~\eqref{E:NEARTIMEAXISREGIONGLOBALEXISTENCEEXTERIORREGIONPOINTWISEESTIMATEEXACTLYONEFLATPARTIALTDERIVATIVEMODIFIEDLRMINUSDERIVATIVE},
	we further deduce, with the help of \eqref{E:PARTIALTINTERMSOFLANDLBAR},
	\eqref{E:COMMUTATOROFLANDULUNIT},
	and the bootstrap assumptions, that the following pointwise estimate holds
	in the near-time-axis region 
	$(t,r) \in [0,\Tboot) \times [0,1]$:
	\begin{align} \label{E:NEARTIMEAXISGLOBALEXISTENCEEXTERIORREGIONPOINTWISEESTIMATEEXACTLYONEFLATPARTIALTDERIVATIVECOMMUTEDINSIDEMODIFIEDLRMINUSDERIVATIVE}
	\begin{split}
		[r \Lunit \partial_t \LRiemann
		+
		2 \partial_t \LRiemann]
		& =
		-
		2 
		\datasize 
		\frac{t+r}{[1 + (t + r)^2]^2}
		+
		\mathcal{O}(\datasize^{1^+}) \frac{\ln_+(t + r)}{(1 + t + r)^3}.
	\end{split}
	\end{align}	
	The estimate \eqref{E:NEARTIMEAXISGLOBALEXISTENCEEXTERIORREGIONPOINTWISEESTIMATEEXACTLYONEFLATPARTIALTDERIVATIVECOMMUTEDINSIDEMODIFIEDLRMINUSDERIVATIVE} 
	in particular implies that:
	\begin{align} \label{E:GLOBALEXISTENCEPARTIALTLRIEMANNATREQUAS0BOUND}
		\partial_t \LRiemann \restriction_{\lbrace r = 0 \rbrace}
	& = 
		-
		\datasize 
		\frac{t}{[1 + t^2]^2}
		+
		\mathcal{O}(\datasize^{1^+}) \frac{\ln_+(t)}{(1 + t)^3}.
\end{align}	
	
	Multiplying \eqref{E:NEARTIMEAXISGLOBALEXISTENCEEXTERIORREGIONPOINTWISEESTIMATEEXACTLYONEFLATPARTIALTDERIVATIVECOMMUTEDINSIDEMODIFIEDLRMINUSDERIVATIVE} 
	by $r$ and again using the identity $\Lunit r = v^r + \Speed$, 
	\eqref{E:SPEEDOFSOUNDEXPANSION}--\eqref{E:SPEEDOFSOUNDERRORFUNCTIONVANISHESATORIGIN},
	and the bootstrap assumptions,
	we further deduce that the following pointwise estimate holds 
	in the near-time-axis region 
	$(t,r) \in [0,\Tboot) \times [0,1]$:
	\begin{align} \label{E:NEARTIMEAXISGLOBALEXISTENCEEXTERIORREGIONPOINTWISEESTIMATELRSQUAREDPARTIALTRMINUS}
	\begin{split}
		\Lunit (r^2 \partial_t \LRiemann)
		& =
		-
		2 
		\datasize 
		\frac{(t+r) r}{[1 + (t + r)^2]^2}
		+
		\mathcal{O}(\datasize^{1^+}) \frac{\ln_+(t + r) r}{(1 + t + r)^3}.
	\end{split}
	\end{align}	
	
\eqref{E:GLOBALEXISTENCEEXTERIORREGIONPOINTWISEESTIMATELRSQUAREDPARTIALTRMINUS} and 
\eqref{E:NEARTIMEAXISGLOBALEXISTENCEEXTERIORREGIONPOINTWISEESTIMATELRSQUAREDPARTIALTRMINUS}
together imply that the following pointwise estimate holds in the entire bootstrap region
$(t,r) \in [0,\Tboot) \times [0,\infty)$:
\begin{align} \label{E:GLOBALEXISTENCEPOINWISEESTIMATEWHOLEBOOTSTRAPREGIONFORLUNITRSQUAREDPARTIALTRMINUS}
	\Lunit(r^2 \partial_t \LRiemann)
	& =
		-
		2 
		\datasize 
		\frac{(t+r) r}{[1 + (t + r)^2]^2}
		+
		\mathcal{O}(\datasize^{1^+}) \frac{\ln_+(t + r) r}{(1 + t + r)^3}.
\end{align}

Starting from \eqref{E:GLOBALEXISTENCEPOINWISEESTIMATEWHOLEBOOTSTRAPREGIONFORLUNITRSQUAREDPARTIALTRMINUS},
	we can integrate and use arguments similar to the ones we used to prove 
	\eqref{E:GLOBALEXISTENCEFARFROMTIMEAXISPOINTWISEAPRIORIRMINUS} and \eqref{E:GLOBALEXISTENCENEARTIMEAXISPOINTWISEAPRIORIRMINUS},
	based on splitting the argument into Case 1 and Case 2 and using
	the data estimates 
	\eqref{E:GLOBALEXISTENCEULUNITORPARTIALTRPLUSPOINTWISEESTIMATEATTIME0} and
	\eqref{E:GLOBALEXISTENCEPARTIALTLRIEMANNATREQUAS0BOUND}
	and the integral estimates
	\eqref{E:GLOBALEXISTENCENEARTIMEAXISINTEGRALESTIMATEROVER1PLUSTPLUSRPOWERAALONGINTEGRALCURVESOFL}--\eqref{E:GLOBALEXISTENCENEARTIMEAXISINTEGRALESITMATELOGTPLUSRTIMESRTOAOVERONEPLUSTPLUSRTOBALONGINTEGRALCURVESOFL}
	and \eqref{E:GLOBALEXISTENCEFARFROMTIMEAXISONEOVERONEPLUSTPLUSRTOPOWERONEPLUSDELTAALONGINTEGRALCURVESOFL}--\eqref{E:GLOBALEXISTENCEFARFROMTIMEAXISLOGTPLUSROVERONEPLUSTPLUSRTOPOWERONEPLUSDELTAALONGINTEGRALCURVESOFL},
	thereby concluding that
	the following pointwise estimate holds 
	in the near-time-axis region $(t,r) \in [0,\Tboot) \times [0,1]$:
	\begin{align} \label{E:GLOBALEXISTENCENEARTIMEAXISPOINTWISEPARTIALTRMINUS} 
	\begin{split}
		|\partial_t \LRiemann|
		& 
		\lesssim
		\datasize
		\frac{1}{(1 + t + r)^3}	
		+
		\datasize^{1^+} 
		\frac{\ln_+(t+r)}{(1 + t + r)^3},
	\end{split}
	\end{align}
	and that the following pointwise holds  
	in the away-from-time-axis region 
	$\lbrace (t,u ) \ | \ 0 \leq t < \Tboot, u \in (-\infty,\Datafunctionforeikonal(t)] \rbrace$:
	\begin{align} \label{E:GLOBALEXISTENCEFARFROMTIMEAXISPOINTWISEPARTIALTRMINUS}  
	\begin{split}
		|\partial_t \LRiemann|
		& 
		\lesssim 
		\frac{\datasize}{(1 + t + r)^2}.
	\end{split}
	\end{align}	
	We have therefore proved
	\eqref{E:GLOBALEXISTENCEFARFROMTIMEAXISPOINTWISEAPRIORIPARTIALTRMINUS} and \eqref{E:GLOBALEXISTENCENEARTIMEAXISPOINTWISEAPRIORIPARTIALTRMINUS}.
	
	\eqref{E:GLOBALEXISTENCEFARFROMTIMEAXISPOINTWISEAPRIORIRTIMESLUNITPARTIALTRMINUS} and \eqref{E:GLOBALEXISTENCENEARTIMEAXISPOINTWISEAPRIORIRTIMESLUNITPARTIALTRMINUS}
	then follow from
	\eqref{E:GLOBALEXISTENCEFARFROMTIMEAXISPOINTWISEAPRIORIPARTIALTRMINUS},
	\eqref{E:GLOBALEXISTENCENEARTIMEAXISPOINTWISEAPRIORIPARTIALTRMINUS},
	\eqref{E:GLOBALEXISTENCEEXTERIORREGIONPOINTWISEESTIMATEEXACTLYONEFLATPARTIALTDERIVATIVECOMMUTEDINSIDEMODIFIEDLRMINUSDERIVATIVE},
	and \eqref{E:NEARTIMEAXISGLOBALEXISTENCEEXTERIORREGIONPOINTWISEESTIMATEEXACTLYONEFLATPARTIALTDERIVATIVECOMMUTEDINSIDEMODIFIEDLRMINUSDERIVATIVE}.

\medskip
\noindent \underline{\textbf{Proof of 
\eqref{E:GLOBALEXISTENCEFARFROMTIMEAXISPOINTWISEAPRIORILRPLUS},
\eqref{E:GLOBALEXISTENCEFARFROMTIMEAXISPOINTWISEAPRIORILBARRMINUS}, 
\eqref{E:GLOBALEXISTENCENEARTIMEAXISPOINTWISEAPRIORILRPLUS},
and \eqref{E:GLOBALEXISTENCENEARTIMEAXISPOINTWISEAPRIORILBARRMINUS}}}.
Using \eqref{E:PARTIALTINTERMSOFLANDLBAR}, 
\eqref{E:SPEEDOFSOUNDEXPANSION}--\eqref{E:SPEEDOFSOUNDERRORFUNCTIONVANISHESATORIGIN},
and the bootstrap assumptions,
we deduce that $|\uLunit \LRiemann| \lesssim |\partial_t \LRiemann| + |\Lunit \LRiemann|$.
From this estimate,
\eqref{E:GLOBALEXISTENCEFARFROMTIMEAXISPOINTWISEAPRIORIPARTIALTRMINUS},
and \eqref{E:GLOBALEXISTENCEFARFROMTIMEAXISPOINTWISEAPRIORILRMINUS},
we conclude \eqref{E:GLOBALEXISTENCEFARFROMTIMEAXISPOINTWISEAPRIORILBARRMINUS}.
\eqref{E:GLOBALEXISTENCENEARTIMEAXISPOINTWISEAPRIORILBARRMINUS} follows from a similar argument based on
\eqref{E:GLOBALEXISTENCENEARTIMEAXISPOINTWISEAPRIORIPARTIALTRMINUS}
and
\eqref{E:GLOBALEXISTENCENEARTIMEAXISPOINTWISEAPRIORILRMINUS}.

\eqref{E:GLOBALEXISTENCEFARFROMTIMEAXISPOINTWISEAPRIORILRPLUS} follows from
\eqref{E:GLOBALEXISTENCEFARFROMTIMEAXISPOINTWISEAPRIORILBARRMINUS} and the identity
$\Lunit \RRiemann = \uLunit \LRiemann$,
which follows from \eqref{E:OUTGOINGRIEMANNINVARIANTEVOLUTION}--\eqref{E:INGOINGRIEMANNINVARIANTEVOLUTION}.
\eqref{E:GLOBALEXISTENCENEARTIMEAXISPOINTWISEAPRIORILRPLUS} follows from a similar argument based on 
\eqref{E:GLOBALEXISTENCENEARTIMEAXISPOINTWISEAPRIORILBARRMINUS}.

\medskip
\noindent \underline{\textbf{Proof of \eqref{E:GLOBALEXISTENCEFARFROMTIMEAXISPOINTWISEAPRIORIRPLUS} and 
\eqref{E:GLOBALEXISTENCENEARTIMEAXISPOINTWISEAPRIORIRPLUS}}}.
Equation \eqref{E:OUTGOINGRIEMANNINVARIANTEVOLUTION} and the bootstrap assumptions imply that
$|\RRiemann| \lesssim |r \Lunit \RRiemann| + |\LRiemann|$.
From this estimate,
\eqref{E:GLOBALEXISTENCERISANALMOSTUNITYMULIPLEOFTMINUSU},
\eqref{E:GLOBALEXISTENCEFARFROMTIMEAXISPOINTWISEAPRIORILRPLUS},
and \eqref{E:GLOBALEXISTENCEFARFROMTIMEAXISPOINTWISEAPRIORIRMINUS},
we conclude the desired bound \eqref{E:GLOBALEXISTENCEFARFROMTIMEAXISPOINTWISEAPRIORIRPLUS}.
\eqref{E:GLOBALEXISTENCENEARTIMEAXISPOINTWISEAPRIORIRPLUS} follows from a similar
argument based on	
\eqref{E:GLOBALEXISTENCENEARTIMEAXISPOINTWISEAPRIORILRPLUS}
and
\eqref{E:GLOBALEXISTENCENEARTIMEAXISPOINTWISEAPRIORIRMINUS}.

\medskip
\noindent \underline{\textbf{Proof of 
\eqref{E:GLOBALEXISTENCEFARFROMTIMEAXISPOINTWISEAPRIORIRTIMESMUULUNITRPLUS},
\eqref{E:GLOBALEXISTENCEFARFROMTIMEAXISPOINTWISEAPRIORIRTIMESULUNITRPLUS}, 
and 
\eqref{E:GLOBALEXISTENCENEARTIMEAXISPOINTWISEAPRIORIRTIMESLBARRPLUS}}}.
	We consider the wave equation \eqref{E:REVAMPEDRPLUSWAVEEQUATION}
	in the away-from-time-axis region $\lbrace (t,u) \ | \ 0 \leq t < \Tboot, u \in (-\infty,\Datafunctionforeikonal(t)] \rbrace$.
	The boxed term on the right cancels the second term on the left, 
	up to the term $- 2 \upmu \Speed \Lunit \RRiemann$, which appears on the left.
	The bootstrap assumptions (see especially \eqref{E:GLOBALEXISTENCEMUSIMPLERBOUNDSBOOTSTRAP}) 
	and \eqref{E:SPEEDOFSOUNDEXPANSION}--\eqref{E:SPEEDOFSOUNDERRORFUNCTIONVANISHESATORIGIN}
	imply that $- 2 \upmu \Speed \Lunit \RRiemann = - 2 \upmu \Lunit \RRiemann$
	plus an error term that is $\lesssim \eps^2 \frac{1}{(1 + t + |u|)^3} 
	+ 
	\eps^2 \frac{\ln(t-u)}{(1 + t + |u|)^3(1 + |u|)^2}$.
	Moreover, $- 2 \upmu \Lunit \RRiemann = \Lunit(- 2 \upmu \RRiemann)$
	plus the error term
	$2 (\Lunit \upmu) \cdot \RRiemann$,
	which by the transport equation \eqref{E:MUEVOLUTION}, the
	bootstrap assumptions (see especially \eqref{E:GLOBALEXISTENCEFARFROMTIMEAXISPOINTWISEBOOTSTRAPRTIMESMUULUNITRPLUS}),
	and \eqref{E:GLOBALEXISTENCECOMPARISONBETWEENTOVERUTMINUSUANDTPLUSMODU}--\eqref{E:GLOBALEXISTENCERISANALMOSTUNITYMULIPLEOFTMINUSU},
	is bounded in magnitude by 
	$\lesssim \frac{\eps^2}{(1 + t + |u|)^2(1 + t-u)} 
	+
	\frac{\eps^2}{(1 + t + |u|)(t-u)(1 + u^2)}$.
	Moreover, the bootstrap assumptions 
	and the transport equation \eqref{E:MUEVOLUTION}
	imply that the remaining (non-boxed)
	terms on RHS~\eqref{E:REVAMPEDRPLUSWAVEEQUATION}
	are bounded in magnitude by 
	$\lesssim \eps^2 \frac{1}{(1 + t + |u|)^3}
		+
	\eps^2 \frac{1}{(1 + t + |u|)^2(1 + |u|)^2}$.
	In total, 
	we find that the following pointwise estimate holds in the away-from-time-axis region
	$\lbrace (t,u) \ | \ 0 \leq t < \Tboot, u \in (-\infty,\Datafunctionforeikonal(t)] \rbrace$:
  \begin{align} \label{E:GLOBALEXISTENCEWAVEEQUATIONRPLUSINHOMOGENEOUSTERMBOUND}
		\left|
		\Lunit 
		\left\lbrace
			r
			\muuLunit \RRiemann
			- 
			2 \upmu \RRiemann
		\right\rbrace
		\right|
		& \lesssim 
		\eps^2 \frac{1}{(1 + t + |u|)^3}
		+
		\eps^2 \frac{1}{(1 + t + |u|)(1 + t-u)(1 + u^2)}.
	\end{align}

Similarly, using the wave equation \eqref{E:RIGHTRIEMANNEQUATIONINTERIORREGION}
and the bootstrap assumptions,
we further deduce (through a simpler argument) that 
the following bound holds in the near-time-axis region $(t,r) \in [0,\Tboot) \times [0,1]$:
\begin{align} \label{E:NEARTIMEAXISGLOBALEXISTENCEWAVEEQUATIONRPLUSINHOMOGENEOUSTERMBOUND}
	\left|
		\Lunit 
		\left\lbrace
			r
			\uLunit \RRiemann
			- 
			2 \RRiemann
		\right\rbrace
		\right|
		& \lesssim 
		\frac{\eps^2}{(1 + t)^3}.
	\end{align}

We will now show that the following pointwise estimate holds in the near-time-axis region $(t,r) \in [0,\Tboot) \times [0,1]$:
\begin{align} \label{E:GLOBALEXISTENCENEARTIMEAXISPOINTWISESPECIALCOMBINATIONRMUULUNITRPLUS}
		r
		\uLunit \RRiemann
		- 
		2 \RRiemann
		&
		=
		 \frac{\left\lbrace
			-2 \datasize 
			+
		\mathcal{O}(\datasize^{1^+})
		\right\rbrace}{1 + (t - r)^2}.
	\end{align}
The proof of \eqref{E:GLOBALEXISTENCENEARTIMEAXISPOINTWISESPECIALCOMBINATIONRMUULUNITRPLUS} splits into two cases.

\noindent \textit{Case $1$: Estimates along integral curves of $\Lunit$ that emanate from the time axis}.
Let $(t_0,0)$ be a point on the time axis with $t_0 \geq 0$, 
and let $(t_1,r_1)$ be a point lying on the integral curve $t \rightarrow \upgamma_{t_0}(t)$ of $\Lunit$ emanating from
$(t_0,0)$ with $t_1 \geq t_0$ with $0 \leq r_1 \leq 1$. 
First, we note that the estimate \eqref{E:GLOBALEXISTENCELRIEMANNATREQUAS0BOUND} and the boundary condition \eqref{E:RIEMANNINVARIANTMATCHINGCONDITION}
imply that:
\begin{align} \label{E:GLOBALEXISTENCERRIEMANNESTIMATEONTIMEAXIS}
	\RRiemann \restriction_{\lbrace r = 0 \rbrace}
	& = \frac{\datasize  + \mathcal{O}(\datasize^{1^+})}{1 + t^2}.
\end{align}	
Next, recalling that $\Lunit t = 1$, we can integrate \eqref{E:NEARTIMEAXISGLOBALEXISTENCEWAVEEQUATIONRPLUSINHOMOGENEOUSTERMBOUND} 
along $\upgamma_{t_0}$ from time $t_0$ to time $t_1$
and use the fundamental theorem of calculus, 
the data bound \eqref{E:GLOBALEXISTENCERRIEMANNESTIMATEONTIMEAXIS},
our assumption that $r_1 \leq 1$, and \eqref{E:GLOBALEXISTENCEBOUNDFORTMINUSRALONGINTEGRALCURVESOFL}
to deduce that:
\begin{align} 
\begin{split} \label{E:GLOBALEXISTENCENEARTIMEAXISEMANATINGFROMTIMEAXISPOINTWISESPECIALCOMBINATIONRMUULUNITRPLUS}
		r_1
		\uLunit \RRiemann(t_1,r_1)
		- 
		2 \RRiemann(t_1,r_1)
		&
		=
		- 
		2 \RRiemann(t_0,0)
		+
		\mathcal{O}(\datasize^{1^+}) \frac{1}{1 + t_0^2}
			\\
		& 
		= 
		\frac{\left\lbrace
			-2 \datasize 
			+
		\mathcal{O}(\datasize^{1^+})
		\right\rbrace}{1 + t_0^2}
			\\
	& = \frac{\left\lbrace
			-2 \datasize 
			+
		\mathcal{O}(\datasize^{1^+})
		\right\rbrace}{1 + (t_1 - r_1)^2},
\end{split}
\end{align}
which yields \eqref{E:GLOBALEXISTENCENEARTIMEAXISPOINTWISESPECIALCOMBINATIONRMUULUNITRPLUS} in Case $1$.

\medskip

\noindent \textit{Case 2: Estimates along integral curves of $\Lunit$ that emanate from $\Sigma_0$}.
Let $(0,r_0) \in \Sigma_0$,
and let $(t_1,r_1)$ be a point lying on the integral curve $t \rightarrow \upgamma_{r_0}(t)$ of $\Lunit$ emanating from
$(0,r_0)$ with $0 \leq r_1 \leq 1$. 
We again integrate \eqref{E:NEARTIMEAXISGLOBALEXISTENCEWAVEEQUATIONRPLUSINHOMOGENEOUSTERMBOUND}, using the fundamental theorem of calculus, 
the data bound \eqref{E:GLOBALEXISTENCETIME0TRCOORDINATESALLDERIVATIVESTRANSPORTEDMODIFIEDULUNITRPLUSDATA} (with $J=K=0$),
our assumption that $r_1 \leq 1$, and \eqref{E:GLOBALEXISTENCEBOUNDFORTMINUSRALONGINTEGRALCURVESOFL} to deduce that:
\begin{align} 
\begin{split} \label{E:GLOBALEXISTENCENEARTIMEAXISEMANATINGFROMSIGMA0POINTWISESPECIALCOMBINATIONRMUULUNITRPLUS}
		r_1
		\uLunit \RRiemann(t_1,r_1)
		- 
		2 \RRiemann(t_1,r_1)
		&
		=
		r_1
		\uLunit \RRiemann(0,r_0)
		- 
		2 \RRiemann(0,r_0)
		+
		\mathcal{O}(\datasize^{1^+})
			\\
		& 
		= 
		\frac{\left\lbrace
			-2 \datasize 
			+
		\mathcal{O}(\datasize^{1^+})
		\right\rbrace}{1 + r_0^2}
			\\
	& = \frac{\left\lbrace
			-2 \datasize 
			+
		\mathcal{O}(\datasize^{1^+})
		\right\rbrace}{1 + (t_1 - r_1)^2},
\end{split}
\end{align}
which yields \eqref{E:GLOBALEXISTENCENEARTIMEAXISPOINTWISESPECIALCOMBINATIONRMUULUNITRPLUS} in Case $2$, thereby finishing its proof.

The desired bound \eqref{E:GLOBALEXISTENCENEARTIMEAXISPOINTWISEAPRIORIRTIMESLBARRPLUS} 
now follows from 
\eqref{E:GLOBALEXISTENCENEARTIMEAXISPOINTWISEAPRIORIRPLUS}
and
\eqref{E:GLOBALEXISTENCENEARTIMEAXISPOINTWISESPECIALCOMBINATIONRMUULUNITRPLUS}.

Next, recalling that $\Lunit = \frac{\partial}{\partial t}$ in geometric coordinates, we integrate
\eqref{E:GLOBALEXISTENCEWAVEEQUATIONRPLUSINHOMOGENEOUSTERMBOUND} in time and use the bootstrap assumptions,
the data bound \eqref{E:GLOBALEXISTENCETIME0TRCOORDINATESALLDERIVATIVESTRANSPORTEDMODIFIEDULUNITRPLUSDATA},
the initial conditions on $\lbrace r=1 \rbrace$ furnished by 
\eqref{E:GLOBALEXISTENCENEARTIMEAXISEMANATINGFROMSIGMA0POINTWISESPECIALCOMBINATIONRMUULUNITRPLUS},
\eqref{E:GLOBALEXISTENCESPACELIKEEIKONALFUNCTIONDATA}--\eqref{E:GLOBALEXISTENCETIMELKEEIKONALFUNCTIONDATA},
and \eqref{E:TIMELILKEEIKONALDATAISAPERTURBATIONOFTMINUS1}--\eqref{E:C2ESTIMATEFORTIMELILKEEIKONALDATAPERTURBATION}
to deduce that the following pointwise estimate holds in the away-from-time-axis region
	$\lbrace (t,u) \ | \ 0 \leq t < \Tboot, u \in (-\infty,\Datafunctionforeikonal(t)] \rbrace$:
  \begin{align} \label{E:GLOBALEXISTENCEPOINTWISESPECIALCOMBINATIONRMUULUNITRPLUS}
			r
			\muuLunit \RRiemann
			- 
			2 \upmu \RRiemann
		& 
		=
	 \frac{\left\lbrace
			- 2 \datasize 
			+
		\mathcal{O}(\datasize^{1^+})
		\right\rbrace}{1 + u^2}.
	\end{align}
	From \eqref{E:GLOBALEXISTENCEPOINTWISESPECIALCOMBINATIONRMUULUNITRPLUS}, \eqref{E:GLOBALEXISTENCEFARFROMTIMEAXISPOINTWISEAPRIORIRPLUS},
	and the bootstrap assumptions, we conclude that:
	\begin{align}
	\begin{split} \label{E:PROOFSTEPSGLOBALEXISTENCEFARFROMTIMEAXISPOINTWISEAPRIORIRTIMESMUULUNITRPLUS}
	|r \muuLunit \RRiemann|
	&
	\lesssim 
	\frac{\datasize}{1 + u^2}
	+
	\upmu |\RRiemann|
		\\
	& \lesssim 
	\frac{\datasize}{1 + u^2} 
	+ 
	\frac{\datasize}{1 + t + |u|} 
	+ 
	\frac{\datasize^{1^+} \ln_+(t-u)}{(1 + t + |u|)(1 + u^2)}
		\\
& \lesssim 
	\frac{\datasize}{1 + u^2} 
	+ 
	\frac{\datasize}{1 + t + |u|},
	\end{split}
	\end{align}
	which yields the desired bound \eqref{E:GLOBALEXISTENCEFARFROMTIMEAXISPOINTWISEAPRIORIRTIMESMUULUNITRPLUS}.
	\eqref{E:GLOBALEXISTENCEFARFROMTIMEAXISPOINTWISEAPRIORIRTIMESULUNITRPLUS} follows from
	a similar argument based on \eqref{E:GLOBALEXISTENCEPOINTWISESPECIALCOMBINATIONRMUULUNITRPLUS}.

\medskip

\noindent \underline{\textbf{Proof of 
\eqref{E:GLOBALEXISTENCEFARFROMTIMEAXISPOINTWISEAPRIORIRTIMESLRPLUS},
\eqref{E:GLOBALEXISTENCEFARFROMTIMEAXISPOINTWISEAPRIORIRTIMESLBARRMINUS},
\eqref{E:GLOBALEXISTENCENEARTIMEAXISPOINTWISEAPRIORIRTIMESLRPLUS}, and
\eqref{E:GLOBALEXISTENCENEARTIMEAXISPOINTWISEAPRIORIRTIMESLBARRMINUS}
}}.
These estimates follows from multiplying equations
\eqref{E:OUTGOINGRIEMANNINVARIANTEVOLUTION}--\eqref{E:INGOINGRIEMANNINVARIANTEVOLUTION} 
by $r$ and using the bound $\Speed \lesssim 1$ (which follows easily from the bootstrap assumptions),
\eqref{E:GLOBALEXISTENCEFARFROMTIMEAXISPOINTWISEAPRIORIRPLUS},
\eqref{E:GLOBALEXISTENCEFARFROMTIMEAXISPOINTWISEAPRIORIRMINUS},
\eqref{E:GLOBALEXISTENCENEARTIMEAXISPOINTWISEAPRIORIRPLUS},
and \eqref{E:GLOBALEXISTENCENEARTIMEAXISPOINTWISEAPRIORIRMINUS}.

\medskip
\noindent \underline{\textbf{Proof of \eqref{E:GLOBALEXISTENCEFARFROMTIMEAXISPOINTWISEAPRIORIMULBARRPLUS} and \eqref{E:GLOBALEXISTENCENEARTIMEAXISPOINTWISEAPRIORILBARRMINUS}}}.
	From \eqref{E:GLOBALEXISTENCEWAVEEQUATIONRPLUSINHOMOGENEOUSTERMBOUND}, 
	\eqref{E:MUEVOLUTION},
	the already proven bound \eqref{E:GLOBALEXISTENCEFARFROMTIMEAXISPOINTWISEAPRIORILRPLUS},
	and the bootstrap assumptions, we
	deduce that the following pointwise estimate holds
	in the away-from-time-axis region $\lbrace (t,u) \ | \ 0 \leq t < \Tboot, u \in (-\infty,\Datafunctionforeikonal(t)] \rbrace$:
  \begin{align} 
		\begin{split} \label{E:GLOBALEXISTENCEPOINTWISELRMUULUNITRPLUS}
		\left|
		\Lunit 
		\left(
			r
			\muuLunit \RRiemann
		\right)
		\right|
		& \lesssim
		|\Lunit \upmu||\RRiemann|
		+
		\upmu |\Lunit \RRiemann|
		+
		\datasize^{1^+} \frac{1}{(1 + t + |u|)^3}
		+
		\datasize^{1^+} \frac{1}{(1 + t + |u|)(1 + t - u)(1 + |u|)^2}
			\\
	& \lesssim
		\datasize \frac{1}{(1 + t + |u|)^2}
		+
		\datasize^{1^+} \frac{\ln_+(t-u)}{(1 + t + |u|)^2(1 + |u|)^2}
		+
		\datasize^{1^+} \frac{1}{(1 + t + |u|)(1 + t - u)(1 + |u|)^2}.
	\end{split}
	\end{align}

	Similarly, using \eqref{E:NEARTIMEAXISGLOBALEXISTENCEWAVEEQUATIONRPLUSINHOMOGENEOUSTERMBOUND}
and the already proven bound \eqref{E:GLOBALEXISTENCENEARTIMEAXISPOINTWISEAPRIORILRPLUS},
we deduce that the following bound holds in the near-time-axis region $(t,r) \in [0,\Tboot) \times [0,1]$:
\begin{align} \label{E:NEARTIMEAXISGLOBALEXISTENCEPOINTWISELRULUNITRPLUS}
		\left|
		\Lunit 
		\left(
			r
			\uLunit \RRiemann
		\right)
		\right|
		& \lesssim 
		\frac{\datasize}{(1 + t + r)^2}.
	\end{align}
	Starting from \eqref{E:NEARTIMEAXISGLOBALEXISTENCEPOINTWISELRULUNITRPLUS},
	we can integrate and use arguments similar to the ones we used to prove 
	\eqref{E:GLOBALEXISTENCENEARTIMEAXISPOINTWISEAPRIORIRMINUS},
	based on splitting the argument into Case 1 and Case 2 and using
	the data estimate \eqref{E:GLOBALEXISTENCEULUNITORPARTIALTRPLUSPOINTWISEESTIMATEATTIME0}, 
	thereby concluding that
	the following pointwise estimate holds 
	in the near-time-axis region $(t,r) \in [0,\Tboot) \times [0,1]$:
	\begin{align} \label{E:GLOBALEXISTENCENEARTIMEAXISPOINTWISEULUNITRPLUS} 
	\begin{split}
		|\uLunit \RRiemann|
		& 
		\lesssim
		\frac{\datasize}{(1 + t + r)^2},
	\end{split}
	\end{align}
	which yields \eqref{E:GLOBALEXISTENCENEARTIMEAXISPOINTWISEAPRIORILBARRMINUS}.
	
	The proof of \eqref{E:GLOBALEXISTENCEFARFROMTIMEAXISPOINTWISEAPRIORIMULBARRPLUS} splits into two cases.
	
\noindent \textit{Case $i$: Estimates along integral curves of $\Lunit$ that emanate from $\lbrace r = 1 \rbrace$}.
Relative to the geometric coordinates, let $(t_0,u)$ be a point on $\lbrace r = 1 \rbrace$, 
and let $(t_1,u)$ be a point lying on the $t-$parameterized integral curve $t \rightarrow \upgamma_{t_0}(t)$ of $\Lunit$ emanating from
$(t_0,u)$ with $t_1 \geq t_0$.
From \eqref{E:TIMELILKEEIKONALDATAISAPERTURBATIONOFTMINUS1}--\eqref{E:C2ESTIMATEFORTIMELILKEEIKONALDATAPERTURBATION},
we see that:
\begin{align} \label{E:GLOBALEXISTENCEUTRELATIONSHIPESTIMATEONREQUALS1}
	u & = t_0 - 1 + \mathcal{O}(\datasize).
\end{align}
For use below, let $r_1 := r(t_1,u)$ denote the radial coordinate evaluated at $(t_1,u)$. 

\noindent \textit{Sub-case $ia$: $r_1 \leq 1 + t_0$}. 
Recalling that $\Lunit = \frac{\partial}{\partial t}$ in geometric coordinates, we integrate \eqref{E:GLOBALEXISTENCEPOINTWISELRMUULUNITRPLUS}
along $\rightarrow \upgamma_{t_0}$ 
from time $t_0$ to time $t_1$. 
\eqref{E:GLOBALEXISTENCECOMPARISONBETWEENTOVERUTMINUSUANDTPLUSMODU} implies that the first term
$\frac{1}{(1 + t + |u|)^2}$  
on RHS~\eqref{E:GLOBALEXISTENCEPOINTWISELRMUULUNITRPLUS} is 
$\lesssim \datasize \frac{1}{(1 + t + r)^2}$. Also using \eqref{E:GLOBALEXISTENCENEARTIMEAXISINTEGRALESTIMATEROVER1PLUSTPLUSRPOWERAALONGINTEGRALCURVESOFL},
we see that
the corresponding integral is $\lesssim \datasize \frac{r_1}{(1 + t_1 + r_1)^2}
\lesssim  \datasize \frac{r_1}{(1 + t_1 + |u|)^2}$.
Similarly, the second term
$\datasize^{1^+} \frac{\ln_+(t-u)}{(1 + t + |u|)^2(1 + |u|)^2}$ 
on RHS~\eqref{E:GLOBALEXISTENCEPOINTWISELRMUULUNITRPLUS} is 
$\lesssim \datasize^{1^+} \frac{\ln_+(r)}{(1 + t + r)^2(1 + |u|)^2}$
and thus, by \eqref{E:GLOBALEXISTENCENEARTIMEAXISINTEGRALESITMATELOGTPLUSRTIMESRTOAOVERONEPLUSTPLUSRTOBALONGINTEGRALCURVESOFL},
the corresponding integral is $\lesssim \datasize^{1^+} \frac{\ln_+(t_1 + r_1) r_1}{(1 + t_1 + r_1)^2(1 + |u|)^2}
\lesssim \datasize^{1^+} \frac{\ln_+(t_1 + |u|) r_1}{(1 + t_1 + |u|)^2(1 + |u|)^2}$.
Similarly, the third term
$\datasize^{1^+} \frac{1}{(1 + t + |u|)(1 + t - u)(1 + |u|)^2}$ 
on RHS~\eqref{E:GLOBALEXISTENCEPOINTWISELRMUULUNITRPLUS} is 
$\lesssim \datasize^{1^+} \frac{1}{(1 + t + r)(1+r)(1 + |u|)^2}$
and thus, by \eqref{E:GLOBALEXISTENCENEARTIMEAXISINTEGRALESTIMATEROVER1PLUSTPLUSRPOWERATIMESONEPLUSTPLUSRALONGINTEGRALCURVESOFL}
the corresponding integral is $\lesssim \datasize^{1^+} \frac{r_1}{(1 + t_1 + r_1)(1 + |u|)^2}$.
Moreover, the initial conditions on $\lbrace r=1 \rbrace$ furnished by \eqref{E:GLOBALEXISTENCENEARTIMEAXISPOINTWISEULUNITRPLUS},
\eqref{E:GLOBALEXISTENCEUTRELATIONSHIPESTIMATEONREQUALS1},
the estimates 
\eqref{E:GLOBALEXISTENCECOMPARISONBETWEENTOVERUTMINUSUANDTPLUSMODU} and \eqref{E:GLOBALEXISTENCEBOUNDFORTMINUSRALONGINTEGRALCURVESOFL},
the bootstrap assumptions,
and our assumption that
$r_1 \leq 1 + t_0$ collectively imply that:
$
|\muuLunit \RRiemann|(t_0,u)
\lesssim
\datasize \frac{1}{(1 + t_0 + |u|)^2}
\lesssim
\datasize \frac{1}{(1 + t_1 + r_1)^2}
\lesssim
\datasize \frac{1}{(1 + t_1 + |u|)^2}
$.
Combining these estimates and applying the fundamental theorem of calculus
to \eqref{E:GLOBALEXISTENCEPOINTWISELRMUULUNITRPLUS},
we deduce the following pointwise estimate, where $r_0 = 1$:
\begin{align} 
	\begin{split} \label{E:FIRSTVERSIONNEARTIMEAXISPOINTWISEESTIMATESFORRTIMESMUULUNITRPLUS}
	r_1 |\muuLunit \RRiemann|(t_1,u)
	& \lesssim
		r_0|\muuLunit \RRiemann|(t_0,u)
		+
		\datasize \frac{r_1}{(1 + t_1 + |u|)^2}
		+
		\datasize^{1^+} \frac{r_1}{(1 + t_1 + |u|)(1 + |u|)^2}
			\\
	& \lesssim 	
		\datasize \frac{1}{(1 + t_1 + |u|)^2}
		+
		\datasize \frac{r_1}{(1 + t_1 + |u|)^2}
		+
		\datasize^{1^+} \frac{r_1}{(1 + t_1 + |u|)(1 + |u|)^2}.
	\end{split}
\end{align}
Dividing \eqref{E:FIRSTVERSIONNEARTIMEAXISPOINTWISEESTIMATESFORRTIMESMUULUNITRPLUS} by $r_1$ 
and using the basic fact that $r_1 \geq 1$, we conclude that:
\begin{align} \label{E:NEARTIMEAXISPOINTWISEESTIMATESFORRTIMESMUULUNITRPLUS}
	|\muuLunit \RRiemann|(t_1,u)
	& \lesssim 
		\datasize \frac{1}{(1 + t_1 + |u|)^2}
		+
		\datasize^{1^+} \frac{1}{(1 + t_1 + |u|)(1 + |u|)^2},
\end{align}
which is stronger than the bound stated in \eqref{E:GLOBALEXISTENCEFARFROMTIMEAXISPOINTWISEAPRIORIMULBARRPLUS}.

\noindent \textit{Sub-case $ib$: $r_1 \geq 1 + t_0$}. 	
In this case,
the estimates \eqref{E:GLOBALEXISTENCEBOUNDFORTMINUSRALONGINTEGRALCURVESOFL} 
and
\eqref{E:GLOBALEXISTENCECOMPARISONBETWEENTOVERUTMINUSUANDTPLUSMODU}
imply that
$r_1 \approx 1 + t_0 + r_1 \approx 1 + t_1 + r_1 \approx 1 + t_1 + |u|$.
From this estimate and
the already proven bound \eqref{E:GLOBALEXISTENCEFARFROMTIMEAXISPOINTWISEAPRIORIRTIMESMUULUNITRPLUS} for 
$r \muuLunit \RRiemann$,
we deduce that: $|\muuLunit \RRiemann|(t_1,u) \lesssim \frac{1}{r_1} \left\lbrace \frac{\datasize}{1 + t_1 + |u|} + \frac{\datasize}{1 + u^2} \right\rbrace \lesssim 
\frac{\datasize}{(1 + t_1 + |u|)^2}
+
 \frac{\datasize}{(1 + t_1 + |u|)(1 + u^2)}
$,
which is stronger than the bound stated in \eqref{E:GLOBALEXISTENCEFARFROMTIMEAXISPOINTWISEAPRIORIMULBARRPLUS}.

\medskip

\noindent \textit{Case $ii$: Estimates along integral curves of $\Lunit$ that emanate from $\Sigma_0 \cap \lbrace r \geq 1 \rbrace$}.
Let $(0,r_0) \in \Sigma_0 \cap \lbrace r \geq 1 \rbrace$,
and let $(t_1,r_1)$ be a point lying on the integral curve $t \rightarrow \upgamma_{r_0}(t)$ of $\Lunit$ emanating from
$(0,r_0)$ with $t_1 \geq 0$. 

\noindent \textit{Sub-case $iia$: $t_1 \leq 1 + r_0$}. 
In this case, the desired bound \eqref{E:GLOBALEXISTENCEFARFROMTIMEAXISPOINTWISEAPRIORIMULBARRPLUS} follows from combining arguments
nearly identical to the ones we used to prove the first line of
\eqref{E:FIRSTVERSIONNEARTIMEAXISPOINTWISEESTIMATESFORRTIMESMUULUNITRPLUS}
with the initial data bound \eqref{E:GLOBALEXISTENCEULUNITORPARTIALTRPLUSPOINTWISEESTIMATEATTIME0},
which, in conjunction with the relation $r_0 = r \restriction_{t=0} = |u|$, implies that: 
$|\muuLunit \RRiemann|(0,r_0) \lesssim \datasize\frac{1}{1 + r_0^2} 
\lesssim 
 \datasize \frac{1}{(1 + r_0 + t_1)^2} =  \datasize \frac{1}{(1 + t_1 + |u|)^2} 
$.

\noindent \textit{Sub-case $iib$: $t_2 \geq 1 + r_0$}. 
In this case,
the estimate \eqref{E:GLOBALEXISTENCEBOUNDFORTMINUSRALONGINTEGRALCURVESOFL} 
and
\eqref{E:GLOBALEXISTENCECOMPARISONBETWEENTOVERUTMINUSUANDTPLUSMODU}
with $(t_2,r_2)$ in the role of $(t_1,r_1)$
and $(0,r_0)$ in the role of $(t_0,r_0)$
imply that:
$
r_2
\approx 
t_2 + r_0
\approx 1 + r_0 + t_2 \approx 1 + r_2 + t_2 \approx 1 + t_2 + |u|$.
Then the same arguments we used in Sub-case $ib$ yield the desired bound,
thereby finishing the proof of \eqref{E:GLOBALEXISTENCEFARFROMTIMEAXISPOINTWISEAPRIORIMULBARRPLUS}.

\medskip

\noindent \underline{\textbf{Proof of \eqref{E:GLOBALEXISTENCEFARFROMTIMEAXISPOINTWISEAPRIORIRTIMESMULBARMUPARTIALTRPLUS} and 
\eqref{E:GLOBALEXISTENCENEARTIMEAXISPOINTWISEAPRIORIRTIMESLBARPARTIALTRPLUS}}}.
We first note that equation \eqref{E:PARTIALTINTERMSOFLANDLBAR} implies that:
\begin{align}	 \label{E:MUPARTIALTINTERMSOFLANDLBAR}
	\upmu \partial_t
& = 
			\frac{1}{2}
			\upmu
			\left\lbrace
				1 
				-
				\frac{v^r}{\Speed}
			\right\rbrace
			\Lunit
			+
			\frac{1}{2}
			\left\lbrace
				1 
				+  
				\frac{v^r}{\Speed}
			\right\rbrace 
			\muuLunit.
\end{align}
	We then commute equation \eqref{E:REVAMPEDRPLUSWAVEEQUATION} with $\upmu \partial_t$, 
carefully noting that the products $2 \upmu \partial_t \left\lbrace 2 \RRiemann  \Lunit (\upmu \Speed) \right\rbrace$
(which are ``dangerous'' in that they involve a second derivative that lacks an $r$-weight)
cancel from both sides of the resulting equation.
Hence, using
\eqref{E:SPEEDOFSOUNDEXPANSION}--\eqref{E:SPEEDOFSOUNDERRORFUNCTIONVANISHESATORIGIN},
\eqref{E:MUPARTIALTINTERMSOFLANDLBAR},
	the commutation formula \eqref{E:COMMUTATOROFLANDMUULUNIT},
	the evolution equation \eqref{E:MUEVOLUTION},
	the bootstrap assumptions, and \eqref{E:GLOBALEXISTENCEWAVEEQUATIONRPLUSINHOMOGENEOUSTERMBOUND},
	we deduce that the following pointwise estimate holds
	in the away-from-time-axis region 
	$\lbrace (t,u ) \ | \ 0 \leq t < \Tboot, u \in (-\infty,\Datafunctionforeikonal(t)] \rbrace$:
  \begin{align} \label{E:GLOBALEXISTENCEWAVEEQUATIONRPLUSMUPARTIALTCOMMUTEDINHOMOGENEOUSTERMBOUND}
		\left|
		\Lunit 
		\left(
		\left\lbrace
			r
			\upmu \partial_t
			(\muuLunit \RRiemann)
			- 
			2 \upmu \Speed (\upmu \partial_t \RRiemann)
		\right\rbrace
		\right)
		\right|
		& \lesssim 
		\datasize^{1^+} 
		\frac{1}{(1 + t + |u|)^2 (1 + |u|)^2}.
	\end{align}
	Commuting the operator $\upmu \partial_t$ 
	on LHS~\eqref{E:GLOBALEXISTENCEWAVEEQUATIONRPLUSMUPARTIALTCOMMUTEDINHOMOGENEOUSTERMBOUND},
	we further deduce, with the help of \eqref{E:MUPARTIALTINTERMSOFLANDLBAR} and the bootstrap assumptions, 
	that the following pointwise estimate holds
	in the away-from-time-axis region $\lbrace (t,u ) \ | \ 0 \leq t < \Tboot, u \in (-\infty,\Datafunctionforeikonal(t)] \rbrace$:
  \begin{align} \label{E:GLOBALEXISTENCEWAVEEQUATIONRPLUSALLTHEWAYINMUPARTIALTCOMMUTEDINHOMOGENEOUSTERMBOUND}
		\left|
		\Lunit 
		\left(
		\left\lbrace
			r
			\muuLunit
			(\upmu \partial_t \RRiemann)
			- 
			2 \upmu \Speed (\upmu \partial_t \RRiemann)
		\right\rbrace
		\right)
		\right|
		& \lesssim 
		\datasize^{1^+} 
		\frac{1}{(1 + t + |u|)^2 (1 + |u|)^2}.
	\end{align}

Similarly, we commute equation \eqref{E:RIGHTRIEMANNEQUATIONINTERIORREGION} with $\partial_t$
and use the bootstrap assumptions to deduce that 
the following bound holds in the near-time-axis region  $(t,r) \in [0,\Tboot) \times [0,1]$:
\begin{align} \label{E:FOLLOWUPGLOBALEXISTENCEWAVEEQUATIONRPLUSPARTIALTCOMMUTEDINHOMOGENEOUSTERMBOUND}
	\left|
		\Lunit 
		\left\lbrace
			r
			\partial_t
			\uLunit \RRiemann
			- 
			2 \Speed \partial_t \RRiemann
		\right\rbrace
		\right|
		& \lesssim 
		\frac{\mathcal{O}(\datasize^{1^+})}{(1 + t)^4}.
	\end{align}

Next, we note that the estimate \eqref{E:GLOBALEXISTENCEPARTIALTLRIEMANNATREQUAS0BOUND},
\eqref{E:SPEEDOFSOUNDEXPANSION}--\eqref{E:SPEEDOFSOUNDERRORFUNCTIONVANISHESATORIGIN},
the bootstrap assumptions, and the boundary condition
\eqref{E:RIEMANNINVARIANTMATCHINGCONDITION} imply that:
	\begin{align} \label{E:GLOBALEXISTENCEPARTIALTRRIEMANNATREQUAS0BOUND}
		\partial_t \RRiemann \restriction_{\lbrace r = 0 \rbrace}
	& = 
		-
		\datasize 
		\frac{t}{[1 + t^2]^2}
		+
		\mathcal{O}(\datasize^{1^+}) \frac{\ln_+(t)}{(1 + t)^3}.
\end{align}	
	
Starting from \eqref{E:FOLLOWUPGLOBALEXISTENCEWAVEEQUATIONRPLUSPARTIALTCOMMUTEDINHOMOGENEOUSTERMBOUND}, we
can integrate and argue as in the proof of \eqref{E:GLOBALEXISTENCENEARTIMEAXISPOINTWISESPECIALCOMBINATIONRMUULUNITRPLUS},
based on splitting the argument into Case 1 and Case 2 and using
the data estimates 
\eqref{E:GLOBALEXISTENCETIME0TRCOORDINATESPARTIALTDERIVATIVETRANSPORTEDMODIFIEDULUNITRPLUSDATA}
and
\eqref{E:GLOBALEXISTENCEPARTIALTRRIEMANNATREQUAS0BOUND}, 
and also using \eqref{E:SPEEDOFSOUNDEXPANSION}--\eqref{E:SPEEDOFSOUNDERRORFUNCTIONVANISHESATORIGIN},
thereby deducing that
the following pointwise estimate holds in the near-time-axis region $(t,r) \in [0,\Tboot) \times [0,1]$:
\begin{align} \label{E:ALMOSTDONEGLOBALEXISTENCENEARTIMEAXISPOINTWISESPECIALCOMBINATIONPARTIALTRPLUS}
	\left|
			r
			\uLunit \partial_t \RRiemann
			- 
			2 \partial_t \RRiemann
		\right|
		& \lesssim 
			\datasize \frac{1}{(1 + t)^3}
			+
			\datasize^{1^+} \frac{\ln_+(t)}{(1 + t)^3}.
	\end{align}
From \eqref{E:ALMOSTDONEGLOBALEXISTENCENEARTIMEAXISPOINTWISESPECIALCOMBINATIONPARTIALTRPLUS}
and
\eqref{E:GLOBALEXISTENCENEARTIMEAXISPOINTWISEAPRIORIPARTIALTRPLUS}, we conclude that
the following pointwise estimate holds for $(t,r) \in [0,\Tboot) \times [0,1]$:
\begin{align} 
\begin{split} \label{E:GLOBALEXISTENCENEARTIMEAXISPOINTWISESPECIALCOMBINATIONPARTIALTRPLUS}
	\left|
			r
			\uLunit \partial_t \RRiemann
		\right|
		& \lesssim 
			2 |\partial_t \RRiemann|
			+
			\datasize \frac{1}{(1 + t)^3}
			+
			\datasize^{1^+} \frac{\ln_+(t)}{(1 + t)^3}
				\\
		& \lesssim 
			\datasize
			\frac{1}{(1 + t + r)(1 + |t-r|)},
	\end{split}
	\end{align}
	which yields the desired bound \eqref{E:GLOBALEXISTENCENEARTIMEAXISPOINTWISEAPRIORIRTIMESLBARPARTIALTRPLUS}.

	Next, recalling that $\Lunit = \frac{\partial}{\partial t}$ in geometric coordinates, 
	we integrate \eqref{E:GLOBALEXISTENCEWAVEEQUATIONRPLUSALLTHEWAYINMUPARTIALTCOMMUTEDINHOMOGENEOUSTERMBOUND} in time
	and use the bootstrap assumptions, the data on $\lbrace r = 1 \rbrace$ induced by
	\eqref{E:ALMOSTDONEGLOBALEXISTENCENEARTIMEAXISPOINTWISESPECIALCOMBINATIONPARTIALTRPLUS}, and the data 
	estimate \eqref{E:GLOBALEXISTENCETIME0TRCOORDINATESPARTIALTDERIVATIVETRANSPORTEDMODIFIEDULUNITRPLUSDATA},
	thereby deducing that the following pointwise estimate holds
	in the away-from-time-axis region 
	$\lbrace (t,u ) \ | \ \tstar \leq t < \Tboot, u \in (-\infty,\Datafunctionforeikonal(t)] \rbrace$:
	\begin{align} \label{E:ALMOSTDONEGLOBALEXISTENCEFARFROMTIMEAXISPOINTWISESPECIALCOMBINATIONMUPARTIALTRPLUS}
	\left|
			r
			\muuLunit
			(\upmu \partial_t \RRiemann)
			- 
			2 \upmu \Speed (\upmu \partial_t \RRiemann)
		\right|
		& \lesssim 
			\datasize \frac{1}{(1 + |u|)^3}
			+
			\datasize^{1^+} \frac{\ln_+(|u|)}{(1 + |u|)^3}.
	\end{align}
	From \eqref{E:ALMOSTDONEGLOBALEXISTENCEFARFROMTIMEAXISPOINTWISESPECIALCOMBINATIONMUPARTIALTRPLUS},
	\eqref{E:SPEEDOFSOUNDEXPANSION}--\eqref{E:SPEEDOFSOUNDERRORFUNCTIONVANISHESATORIGIN},
	\eqref{E:MUPARTIALTINTERMSOFLANDLBAR}, the bootstrap assumptions, 
	and the already proven estimates \eqref{E:GLOBALEXISTENCEFARFROMTIMEAXISPOINTWISEAPRIORIMULBARRPLUS} and \eqref{E:GLOBALEXISTENCEFARFROMTIMEAXISPOINTWISEAPRIORILRPLUS},
	 we conclude that
the following pointwise estimate holds in the away-from-time-axis region 
	$\lbrace (t,u ) \ | \ \tstar \leq t < \Tboot, u \in (-\infty,\Datafunctionforeikonal(t)] \rbrace$:
	\begin{align} 
	\begin{split} \label{E:GLOBALEXISTENCEFARFROMTIMEAXISPOINTWISESPECIALCOMBINATIONMUPARTIALTRPLUS}
	\left|
			r
			\muuLunit
			(\upmu \partial_t \RRiemann)
		\right|
		& \lesssim 
			\upmu |\upmu \Lunit \RRiemann|
			+
			\upmu |\muuLunit \RRiemann|
			+
			\datasize \frac{1}{(1 + |u|)^3}
			+
			\datasize^{1^+} \frac{\ln_+(|u|)}{(1 + |u|)^3}
				\\
		& \lesssim 
		\datasize 
		\frac{1}{(1 + t + |u|)(1 + |u|)}
		+
		\datasize \frac{1}{(1 + |u|)^3}
		+
		\datasize^{1^+} \frac{\ln_+(|u|)}{(1 + |u|)^3},
	\end{split}
	\end{align}
	which yields the desired bound 
	\eqref{E:GLOBALEXISTENCEFARFROMTIMEAXISPOINTWISEAPRIORIRTIMESMULBARMUPARTIALTRPLUS}.
	
	\medskip
\noindent \underline{\textbf{Proof of \eqref{E:GLOBALEXISTENCEFARFROMTIMEAXISPOINTWISEAPRIORITIMESLUNITMUULUNITRPLUS}, \eqref{E:GLOBALEXISTENCEFARFROMTIMEAXISPOINTWISEAPRIORIRTIMESMUULUNITLUNITRMINUS}, \eqref{E:GLOBALEXISTENCENEARTIMEAXISPOINTWISEAPRIORITIMESLUNITLBARRPLUS}, and \eqref{E:GLOBALEXISTENCENEARTIMEAXISPOINTWISEAPRIORIRTIMESLBARLUNITRMINUS}}}.

To prove \eqref{E:GLOBALEXISTENCEFARFROMTIMEAXISPOINTWISEAPRIORITIMESLUNITMUULUNITRPLUS}, 
we first use \eqref{E:GLOBALEXISTENCEWAVEEQUATIONRPLUSINHOMOGENEOUSTERMBOUND},
the identity $\Lunit r = v^r + \Speed$ (see \eqref{E:NULLVECTORFIELDS}),
and the bootstrap assumptions to deduce that
$|r \Lunit \muuLunit \RRiemann| 
\lesssim 	
|\muuLunit \RRiemann|
+
\upmu |\Lunit \RRiemann|
+
|\Lunit \upmu| |\RRiemann|
+
\datasize^{1^+} \frac{1}{(1 + t + |u|)^3}
		+
		\datasize^{1^+} \frac{1}{(1 + t + |u|)(1 + t - u)(1 + |u|)^2}$.
From this bound, 
the transport equation \eqref{E:MUEVOLUTION},
the bootstrap assumptions, and the already proven bounds \eqref{E:GLOBALEXISTENCEFARFROMTIMEAXISPOINTWISEAPRIORIMULBARRPLUS} and
\eqref{E:GLOBALEXISTENCEFARFROMTIMEAXISPOINTWISEAPRIORILRPLUS},
we conclude \eqref{E:GLOBALEXISTENCEFARFROMTIMEAXISPOINTWISEAPRIORITIMESLUNITMUULUNITRPLUS}.

\eqref{E:GLOBALEXISTENCENEARTIMEAXISPOINTWISEAPRIORITIMESLUNITLBARRPLUS} follows from a similar argument based on
\eqref{E:NEARTIMEAXISGLOBALEXISTENCEWAVEEQUATIONRPLUSINHOMOGENEOUSTERMBOUND},
\eqref{E:GLOBALEXISTENCENEARTIMEAXISPOINTWISEAPRIORILBARRPLUS},
and \eqref{E:GLOBALEXISTENCENEARTIMEAXISPOINTWISEAPRIORILRPLUS}.
	
\eqref{E:GLOBALEXISTENCEFARFROMTIMEAXISPOINTWISEAPRIORIRTIMESMUULUNITLUNITRMINUS}
follows from a similar argument based on
\eqref{E:GLOBALEXISTENCEWAVEEQUATIONRMINUSINHOMOGENEOUSTERMBOUND},
\eqref{E:GLOBALEXISTENCEFARFROMTIMEAXISPOINTWISEAPRIORILBARRMINUS},
and \eqref{E:GLOBALEXISTENCEFARFROMTIMEAXISPOINTWISEAPRIORILRMINUS}.

\eqref{E:GLOBALEXISTENCENEARTIMEAXISPOINTWISEAPRIORIRTIMESLBARLUNITRMINUS}
follows from a similar argument based on \eqref{E:NEARTIMEAXISGLOBALEXISTENCEWAVEEQUATIONRMINUSINHOMOGENEOUSTERMBOUND},
\eqref{E:GLOBALEXISTENCENEARTIMEAXISPOINTWISEAPRIORILBARRMINUS}, and \eqref{E:GLOBALEXISTENCENEARTIMEAXISPOINTWISEAPRIORILRMINUS}.

	\medskip
\noindent \underline{\textbf{Proof of \eqref{E:GLOBALEXISTENCEFARFROMTIMEAXISPOINTWISEAPRIORITIMESLUNITLUNITRPLUS}, \eqref{E:GLOBALEXISTENCEFARFROMTIMEAXISPOINTWISEAPRIORIRTIMESMUULUNITULUNITRMINUS}, \eqref{E:GLOBALEXISTENCENEARTIMEAXISPOINTWISEAPRIORITIMESLUNITLUNITRPLUS}, and \eqref{E:GLOBALEXISTENCENEARTIMEAXISPOINTWISEAPRIORIRTIMESLBARLBARRMINUS}}}.

\eqref{E:GLOBALEXISTENCEFARFROMTIMEAXISPOINTWISEAPRIORITIMESLUNITLUNITRPLUS} 
and \eqref{E:GLOBALEXISTENCENEARTIMEAXISPOINTWISEAPRIORITIMESLUNITLUNITRPLUS}
follow from multiplying equation
\eqref{E:OUTGOINGRIEMANNINVARIANTEVOLUTION}
by $r$, taking an $\Lunit$ derivative, and using 
the identity $\Lunit r = v^r + \Speed$ (see \eqref{E:NULLVECTORFIELDS}),
\eqref{E:SPEEDOFSOUNDEXPANSION}--\eqref{E:SPEEDOFSOUNDERRORFUNCTIONVANISHESATORIGIN},
the bootstrap assumptions,
and the already proven bounds
\eqref{E:GLOBALEXISTENCEFARFROMTIMEAXISPOINTWISEAPRIORILRPLUS},
\eqref{E:GLOBALEXISTENCEFARFROMTIMEAXISPOINTWISEAPRIORILRMINUS},
\eqref{E:GLOBALEXISTENCENEARTIMEAXISPOINTWISEAPRIORILRPLUS},
and
\eqref{E:GLOBALEXISTENCENEARTIMEAXISPOINTWISEAPRIORILRMINUS}.

\eqref{E:GLOBALEXISTENCEFARFROMTIMEAXISPOINTWISEAPRIORIRTIMESMUULUNITULUNITRMINUS}	
and
\eqref{E:GLOBALEXISTENCENEARTIMEAXISPOINTWISEAPRIORIRTIMESLBARLBARRMINUS}
follow from similar arguments based on equation
\eqref{E:INGOINGRIEMANNINVARIANTEVOLUTION}
and the already proven estimates
\eqref{E:GLOBALEXISTENCEFARFROMTIMEAXISPOINTWISEAPRIORIMULBARRPLUS},
\eqref{E:GLOBALEXISTENCEFARFROMTIMEAXISPOINTWISEAPRIORILBARRMINUS},
\eqref{E:GLOBALEXISTENCENEARTIMEAXISPOINTWISEAPRIORILBARRPLUS},
and
\eqref{E:GLOBALEXISTENCENEARTIMEAXISPOINTWISEAPRIORILBARRMINUS}.



\medskip
\noindent \underline{\textbf{Proof of \eqref{E:GLOBALEXISTENCEFARFROMTIMEAXISPOINTWISEAPRIORITIMESMUULUNITMUULUNITRPLUS}
and \eqref{E:GLOBALEXISTENCEFARFROMTIMEAXISPOINTWISEAPRIORIRTIMESLUNITLUNITRMINUS}}}.
These are the only estimates in the proposition that have not yet been proved.
They easily from the identity \eqref{E:PARTIALTINTERMSOFLANDLBAR}, the bootstrap assumptions, 
and the other estimates in the proposition.

\end{proof}

\subsection{Estimates for the inverse foliation density}
\label{SS:GLOBALEXISTENCEINVERSEFOLIATIONDENSITYESTIMATES}
In this section, we prove a proposition that yields sharp pointwise estimates for $\upmu$. These estimates in particular provide 
an improvement of the 
bootstrap assumptions
\eqref{E:GLOBALEXISTENCEMUSIMPLERBOUNDSBOOTSTRAP}--\eqref{E:GLOBALEXISTENCEMUXMUSIMPLERUPPERBOUNDBOOTSTRAP}.

\begin{proposition}[Sharp estimates for the eikonal function and the inverse foliation density]
\label{P:GLOBALEXISTENCESHARPESTIMATESFORMU}
Recall that $\lifespanconstant$ is the constant
from \eqref{E:NULLCONDITIONFAILURECONSTANT}, and that we are assuming that it satisfies\footnote{In the case of the
Chaplygin gas equation of state, we have $\lifespanconstant = 0$. The estimates are easier in this case, because several terms identically vanish, 
thereby leading to the absence of various terms that would otherwise have grown logarithmically in time. \label{FN:CHAPLYGINGASEASIER}}  
$\lifespanconstant \geq 0$.
Under the data assumptions of Section~\ref{SS:GLOBALEXISTENCEDATAATTIME0} and the bootstrap assumptions of
Sect.\,\ref{SS:GLOBALEXISTENCEBOOTSTRAPANDTWOREGIONS}, 
if $\datasize$ is sufficiently small, then
the following estimates hold on the away-from-time-axis region
	$\lbrace (t,u) \ | \ 0 \leq t < \Tboot, u \in (-\infty,\Datafunctionforeikonal(t)] \rbrace$:

\begin{subequations}
	\begin{align}  \label{E:GLOBALEXISTENCEMUPOINTWISESTIMATE} 
	\begin{split}
	\upmu
	& = 
			\begin{cases}
			1 
			+
			\mathcal{O}(\datasize)
			+
			\frac{\lifespanconstant \datasize + \mathcal{O}(\datasize^{1^+})}{(1 + u^2)} \ln \left( \frac{t-u}{-u} \right),
			& u \leq - 1,
				\\
			1 
			+
			\mathcal{O}(\datasize)
			+
			\frac{\lifespanconstant \datasize + \mathcal{O}(\datasize^{1^+})}{(1 + u^2)} \ln (t-u),
			& u \geq -1,
			\end{cases}
	\end{split}
		\\
	|\Lunit \upmu|  \label{E:GLOBALEXISTENCELUNITMUPOINTWISESTIMATE} 
	& \leq C \datasize \frac{1}{(1 + t + |u|)(1 + |u|)},
		\\
		\label{E:GLOBALEXISTENCEMUULUNITMUPOINTWISESTIMATE}
	\begin{split}
	|\muuLunit \upmu|
	& \lesssim 
			\begin{cases}
			\datasize \frac{1}{1 + |u|^2}
			+
			\datasize \frac{\ln \left( \frac{t-u}{-u} \right)}{(1 + |u|)^3}
		+
		\datasize^{1^+} \frac{\ln \left( \frac{t-u}{-u} \right) \ln_+(|u|)}{(1 + |u|)^3},
			& u \leq - 1,
				\\
			\datasize \frac{\ln_+(u_+)}{1 + |u|^2}
			+
			\datasize
			\frac{\ln \left(t-u \right)}{(1 + |u|)^3} 
			+
			\datasize^{1^+}
			\frac{\ln \left(t-u \right) \ln_+(|u|)}{(1 + |u|)^3},
			& u \geq -1.
			\end{cases}
	\end{split}
\end{align}
\end{subequations}

In particular, we have:
\begin{subequations}
\begin{align}
	 \label{E:GLOBALEXISTENCEMUSIMPLERBOUNDS}
	1
	- 
	C \datasize
	& 
	\leq
	\upmu 
	\leq 1 
			+
			C \datasize
			+ 
			C \frac{\datasize}{(1 + u^2)} \ln_+(t - u),
				\\
	|\muuLunit \upmu|
	& \leq  
			C
			\datasize \frac{\ln_+(u_+)}{1 + |u|^2}
			+
			C \datasize
			\frac{\ln_+ \left(t-u \right) \ln_+(|u|)}{(1 + |u|)^3}.
	\label{E:GLOBALEXISTENCEMUXMUSIMPLERUPPERBOUND}		
\end{align}
\end{subequations}

\end{proposition}

\begin{remark}[Strict improvement of the bootstrap assumptions]
	\label{R:GLOBALEXISTENCESTRICTIMPROVEMENTOFBOOTSTRAP}
	Recall that we assumed only that $0 < \eps \leq \datasize^{3/4}$
	in the bootstrap assumptions that we made
	in Section~\ref{S:EXTERIORGLOBALLYHYPERBOLICBOOTSTRAPASSUMPTIONS}; see \eqref{E:GLOBALEXISTENCEBOUNDONBOOTSTRAPPARAMETERSIZE}.
	It is therefore straightforward to see that 
	the estimates of Props.\,\ref{P:GLOBALEXISTENCERIEMANNINVARIANTSAPRIORIEXTERIORREGIONESTIMATES}
	and \ref{P:GLOBALEXISTENCESHARPESTIMATESFORMU}
	collectively imply strict improvements of all 
	the bootstrap assumptions of Sect.\,\ref{SS:GLOBALEXISTENCEBOOTSTRAPANDTWOREGIONS}.
\end{remark}

\begin{proof}[Proof of Prop.\,\ref{P:GLOBALEXISTENCESHARPESTIMATESFORMU}]
 \ \\
	\medskip
	\noindent \underline{\textbf{Data estimates}}.
	We will first show that the following ``data estimates'' hold for $\upmu$:
	\begin{subequations}
\begin{align}
	\upmu \restriction_{\Sigma_0 \cap \lbrace r \geq 1 \rbrace}
	& = 1 + \mathcal{O}(\datasize) \frac{1}{1 + |u|},
		\label{E:MUESTIMATEATTIME0} 
			\\
		\upmu \restriction_{\lbrace r = 1 \rbrace}
	& = 1 + \mathcal{O}(\datasize)\frac{1}{(1+t)^2}
		= 1 + \mathcal{O}(\datasize)\frac{1}{1+|u|^2},
		\label{E:MUESTIMATEATREQUALS1} 
\end{align}
\end{subequations}

\begin{subequations}
\begin{align}
	\muuLunit \upmu \restriction_{\Sigma_0 \cap \lbrace r \geq 1 \rbrace}
	& = \mathcal{O}(\datasize) \frac{1}{1 + |u|^2},
		\label{E:MUULUNITMUESTIMATEATTIME0} 
			\\
	\muuLunit \upmu \restriction_{\lbrace r = 1 \rbrace}
	& = \mathcal{O}(\datasize) \frac{1}{(1+t)^2}
		= \mathcal{O}(\datasize) \frac{1}{1+|u|^2}.
		\label{E:MUULUNITMUESTIMATEATREQUALS1} 
\end{align}
\end{subequations}
\eqref{E:MUESTIMATEATTIME0}--\eqref{E:MUESTIMATEATREQUALS1} 
follow from the identities \eqref{E:IDENTITYFORMUATTIMEZERO}--\eqref{E:IDENTITYFORMUALONGREQUALS1},
\eqref{E:SPEEDOFSOUNDEXPANSION}--\eqref{E:SPEEDOFSOUNDERRORFUNCTIONVANISHESATORIGIN},
\eqref{E:TIMELILKEEIKONALDATAISAPERTURBATIONOFTMINUS1}--\eqref{E:C2ESTIMATEFORTIMELILKEEIKONALDATAPERTURBATION},
and the estimates \eqref{E:GLOBALEXISTENCEFARFROMTIMEAXISPOINTWISEAPRIORIRPLUS} and \eqref{E:GLOBALEXISTENCEFARFROMTIMEAXISPOINTWISEAPRIORIRMINUS}.

To prove \eqref{E:MUULUNITMUESTIMATEATTIME0}, 
we first use 
\eqref{E:PARTIALRINTERMSOFLUNITANDULUNIT}
and
\eqref{E:IDENTITYFORMUATTIMEZERO}
to deduce the following identity: 
\begin{align} \label{E:IDENTITYFORMUULUNITMUATTIME0}
\begin{split}
	\muuLunit \upmu \restriction_{\Sigma_0 \cap \lbrace r \geq 1 \rbrace}
	& = \upmu \Lunit \upmu \restriction_{\Sigma_0 \cap \lbrace r \geq 1 \rbrace}
		-
		2 \upmu \Speed \partial_r \left(\frac{1}{\Speed} \right) \restriction_{\Sigma_0 \cap \lbrace r \geq 1 \rbrace}.
\end{split}
\end{align}
 From \eqref{E:IDENTITYFORMUULUNITMUATTIME0},
\eqref{E:SPEEDOFSOUNDEXPANSION}--\eqref{E:SPEEDOFSOUNDERRORFUNCTIONVANISHESATORIGIN},
\eqref{E:PARTIALRINTERMSOFLUNITANDULUNIT},
\eqref{E:MUEVOLUTION},
\eqref{E:IDENTITYFORMUATTIMEZERO},
the estimates of Lemma~\ref{L:GLOBALEXISTENCERIEMANNINVARIANTESTIMATESATTIME0},
and the estimate \eqref{E:MUESTIMATEATTIME0},
we conclude the desired bound \eqref{E:MUULUNITMUESTIMATEATTIME0}.

To prove \eqref{E:MUULUNITMUESTIMATEATREQUALS1},
we first use \eqref{E:NULLVECTORFIELDS} to deduce the following identity: 
$\uLunit =\left\lbrace 1-  \frac{v^r - \Speed}{v^r + \Speed} \right\rbrace \partial_t
+
\frac{v^r - \Speed}{v^r + \Speed} \Lunit$.
From this identity and \eqref{E:IDENTITYFORMUALONGREQUALS1},
we deduce the following identity:
\begin{align} \label{E:IDENTITYFORMUULUNITONREQUAS1}
	\muuLunit \upmu \restriction_{\lbrace r = 1 \rbrace}
	& = 
			\upmu
			\left\lbrace 1-  \frac{v^r - \Speed}{v^r + \Speed} \right\rbrace \partial_t
			\left(
			\frac{1 + \frac{v^r}{\Speed}}{\partial_t \Datafunctionforeikonal(t)} \right)
			\restriction_{\lbrace r = 1 \rbrace}
			+
			\upmu
			\frac{v^r - \Speed}{v^r + \Speed} \Lunit \upmu \restriction_{\lbrace r = 1 \rbrace}.
\end{align}

From \eqref{E:TIMELILKEEIKONALDATAISAPERTURBATIONOFTMINUS1} 
and
\eqref{E:C2ESTIMATEFORTIMELILKEEIKONALDATAPERTURBATION}, we see that:
\begin{align} \label{E:FIRSTTIMEDERIVATIVEOFTIMELIKEEIKONALFUNCTIONDATAISCLOSETOUNITY}
	\partial_t \Datafunctionforeikonal(t)
	& = 1 + \mathcal{O}(\datasize).
\end{align}

Moreover, \eqref{E:C2ESTIMATEFORTIMELILKEEIKONALDATAPERTURBATION} and \eqref{E:UISTMINUS1FORTBIGGERTHAN2} imply that 
$\partial_t^2 \Datafunctionforeikonal(t) = 0 $ for $t \geq 2$ (i.e., when $u \geq 1$ along $\lbrace r = 1 \rbrace$)
and that $|\partial_t^2 \Datafunctionforeikonal(t)| \lesssim \datasize $ for $0 \leq t < 2$ (i.e., when $-1 \leq u < 1$ along $\lbrace r = 1 \rbrace$).
Together, these bounds imply that:
\begin{align}  \label{E:SIMPLEBOUNDFORSECONDTIMEDERIVATIVEOFTIMELIKEEIKONALFUNCTIONDATA}
	|\partial_t^2 \Datafunctionforeikonal| \restriction_{\lbrace r = 1 \rbrace \cap \lbrace t \geq 0 \rbrace}
	& \lesssim \frac{\datasize}{(1 + |u|)^2}.
\end{align}

Next, we note that the estimates of Prop.\,\ref{P:GLOBALEXISTENCERIEMANNINVARIANTSAPRIORIEXTERIORREGIONESTIMATES},
the identities \eqref{E:PARTIALRINTERMSOFLUNITANDULUNIT} and \eqref{E:PARTIALTINTERMSOFLANDLBAR},
and \eqref{E:C2ESTIMATEFORTIMELILKEEIKONALDATAPERTURBATION} imply that
imply that on $\lbrace r = 1 \rbrace$,
the first derivatives of $(\RRiemann,\LRiemann)$ with respect any of $\partial_t, \partial_r, \Lunit, \uLunit$ are
are bounded in magnitude by
$\lesssim \frac{\datasize}{(1 + |u|)^2}$.
From this bound, the estimates \eqref{E:FIRSTTIMEDERIVATIVEOFTIMELIKEEIKONALFUNCTIONDATAISCLOSETOUNITY}--\eqref{E:SIMPLEBOUNDFORSECONDTIMEDERIVATIVEOFTIMELIKEEIKONALFUNCTIONDATA},
the estimates of Prop.\,\ref{P:GLOBALEXISTENCERIEMANNINVARIANTSAPRIORIEXTERIORREGIONESTIMATES},
the estimate \eqref{E:MUESTIMATEATREQUALS1},
\eqref{E:SPEEDOFSOUNDEXPANSION}--\eqref{E:SPEEDOFSOUNDERRORFUNCTIONVANISHESATORIGIN},
and equation \eqref{E:MUEVOLUTION},
we conclude that RHS~\eqref{E:IDENTITYFORMUULUNITONREQUAS1} is bounded in magnitude by
$\lesssim \frac{\datasize}{(1 + |u|)^2}$, which yields the desired bound \eqref{E:MUULUNITMUESTIMATEATREQUALS1}.

	\medskip
	\noindent \underline{\textbf{Proof of \eqref{E:GLOBALEXISTENCEMUPOINTWISESTIMATE}}}.
	To prove \eqref{E:GLOBALEXISTENCEMUPOINTWISESTIMATE}, we first use
	equation \eqref{E:MUEVOLUTION},
\eqref{E:SPEEDOFSOUNDEXPANSION}--\eqref{E:SPEEDOFSOUNDERRORFUNCTIONVANISHESATORIGIN},
	\eqref{E:NULLCONDITIONFAILURECONSTANT},
	the bootstrap assumptions,
	and the estimates of Prop.\,\ref{P:GLOBALEXISTENCERIEMANNINVARIANTSAPRIORIEXTERIORREGIONESTIMATES} to deduce
	that the following pointwise estimate holds in the away-from-time-axis region 
	$\lbrace (t,u) \ | \ 0 \leq t < \Tboot, u \in (-\infty,\Datafunctionforeikonal(t)] \rbrace$:
	\begin{align} 
	\begin{split} \label{E:GLOBALEXISTENCELUNITMUFIRSTPOINTWISEESTIMATE}
		\Lunit \upmu
		& = 
			-
			\frac{1}{2 \Speed}
			\left\lbrace
				1
				+
				\Speed' 
				\InverseRiemannfunction'
			\right\rbrace
			\muuLunit 
			\RRiemann
			+
			\mathcal{O}(\datasize) \frac{1}{(1 + t + |u|)^2}
			+
			\mathcal{O}(\datasize^{1^+}) \frac{\ln_+(t-u)}{(1 + t + |u|)^2 (1 + |u|)^2}
				\\
		& = - \frac{1}{2}
					\lifespanconstant
					\muuLunit \RRiemann
					+
			\mathcal{O}(\datasize) \frac{1}{(1 + t + |u|)^2}
			+
			\mathcal{O}(\datasize^{1^+}) \frac{\ln_+(t-u)}{(1 + t + |u|)^2 (1 + |u|)^2}.
		\end{split}
	\end{align}
	Also using 
	\eqref{E:GLOBALEXISTENCERISANALMOSTUNITYMULIPLEOFTMINUSU},
	\eqref{E:GLOBALEXISTENCEPOINTWISESPECIALCOMBINATIONRMUULUNITRPLUS},
	the bootstrap assumptions, 
	and the estimates of Prop.\,\ref{P:GLOBALEXISTENCERIEMANNINVARIANTSAPRIORIEXTERIORREGIONESTIMATES}
	to estimate the term $\muuLunit \RRiemann$ on RHS~\eqref{E:GLOBALEXISTENCELUNITMUFIRSTPOINTWISEESTIMATE},
	we further deduce
	that the following pointwise estimate holds in the away-from-time-axis region 
	$\lbrace (t,u) \ | \ 0 \leq t < \Tboot, u \in (-\infty,\Datafunctionforeikonal(t)] \rbrace$:
	\begin{align}	
		\begin{split} \label{E:GLOBALEXISTENCELUNITMUSECONDPOINTWISEESTIMATE}
		\Lunit \upmu
		& = 
			\lifespanconstant
			\frac{\datasize + \mathcal{O}(\datasize^{1^+})}{(t-u)(1 + u^2)}
			+
			\mathcal{O}(\datasize) \frac{1}{(t-u)(1 + t + |u|)}
			+
			\mathcal{O}(\datasize) \frac{\ln_+(t-u)}{(t-u)(1 + t + |u|)(1 + |u|)^2}.
	\end{split}
	\end{align}
	Recalling that $\Lunit = \frac{\partial}{\partial t}$ in geometric coordinates, we integrate 
	\eqref{E:GLOBALEXISTENCELUNITMUSECONDPOINTWISEESTIMATE} with respect to time and use the 
	data estimates
	\eqref{E:MUESTIMATEATTIME0}--\eqref{E:MUESTIMATEATREQUALS1} for $\upmu$ to conclude the following pointwise estimate:
	\begin{align}  \label{E:PROOFSTEPGLOBALEXISTENCEMUPOINTWISESTIMATE} 
	\begin{split}
	\upmu
	& = 
			\begin{cases}
			1 
			+
			\mathcal{O}(\datasize)
			+
			\frac{\lifespanconstant \datasize + \mathcal{O}(\datasize^{1^+})}{(1 + u^2)} \ln \left( \frac{t-u}{-u} \right),
			& u \leq - 1,
				\\
			1 
			+
			\mathcal{O}(\datasize)
			+
			\frac{\lifespanconstant \datasize + \mathcal{O}(\datasize^{1^+})}{(1 + u^2)} \ln (t-u),
			& u \geq -1,
			\end{cases}
	\end{split}
	\end{align}
	which yields \eqref{E:GLOBALEXISTENCEMUPOINTWISESTIMATE}.
	
	\medskip
	\noindent \underline{\textbf{Proof of \eqref{E:GLOBALEXISTENCELUNITMUPOINTWISESTIMATE}}}.
	\eqref{E:GLOBALEXISTENCELUNITMUPOINTWISESTIMATE} can be proved via arguments that are similar to but simpler than the ones we used to prove 
	\eqref{E:GLOBALEXISTENCELUNITMUSECONDPOINTWISEESTIMATE}, where we rely on \eqref{E:GLOBALEXISTENCEFARFROMTIMEAXISPOINTWISEAPRIORIMULBARRPLUS} 
	instead of \eqref{E:GLOBALEXISTENCEPOINTWISESPECIALCOMBINATIONRMUULUNITRPLUS}
	to estimate the term $\muuLunit \RRiemann$ on RHS~\eqref{E:GLOBALEXISTENCELUNITMUFIRSTPOINTWISEESTIMATE}.

	\medskip
	\noindent \underline{\textbf{Proof of \eqref{E:GLOBALEXISTENCEMUULUNITMUPOINTWISESTIMATE}}}.
	To prove \eqref{E:GLOBALEXISTENCEMUULUNITMUPOINTWISESTIMATE}, we commute
	equation \eqref{E:MUEVOLUTION} with $\muuLunit$ and use 
	the commutation formula \eqref{E:COMMUTATOROFLANDMUULUNIT},
	the evolution equation \eqref{E:MUEVOLUTION},
	the bootstrap assumptions,
	\eqref{E:GLOBALEXISTENCERISANALMOSTUNITYMULIPLEOFTMINUSU},
	and the estimates of Prop.\,\ref{P:GLOBALEXISTENCERIEMANNINVARIANTSAPRIORIEXTERIORREGIONESTIMATES}.
	From the point of view of decay, the main term that arises in the analysis is the $\muuLunit$ derivative of the first product on
	RHS~\eqref{E:MUEVOLUTION}, which leads to a term of type
	$\muuLunit \muuLunit \RRiemann$ that we bound in magnitude by
	$
	\lesssim 
	\datasize 
		\frac{1}{(1 + t + |u|)(t - u)(1 + |u|)}
		+
		\datasize \frac{1}{(t - u) (1 + |u|)^3}
		+
		\datasize^{1^+} \frac{\ln_+(|u|)}{(t - u) (1 + |u|)^3}
	$
	with the help of \eqref{E:GLOBALEXISTENCEFARFROMTIMEAXISPOINTWISEAPRIORITIMESMUULUNITMUULUNITRPLUS}
	and \eqref{E:GLOBALEXISTENCERISANALMOSTUNITYMULIPLEOFTMINUSU}. The remaining terms that arise decay even faster.
	In total, we deduce that the following pointwise estimate holds 
	in the away-from-time-axis region 
	$\lbrace (t,u) \ | \ 0 \leq t < \Tboot, u \in (-\infty,\Datafunctionforeikonal(t)] \rbrace$:
	\begin{align}	
		\begin{split} \label{E:GLOBALEXISTENCELUNITMUULUNITMUFIRSTPOINTWISEESTIMATE}
		|\Lunit \muuLunit \upmu|
		& \lesssim 
		\datasize 
		\frac{1}{(1 + t + |u|)(t - u)(1 + |u|)}
			\\
		&  \ \
		+
		\datasize \frac{1}{(t - u) (1 + |u|)^3}
		+
		\datasize^{1^+} \frac{\ln_+(|u|)}{(t - u) (1 + |u|)^3}.
	\end{split}
	\end{align}
	Again recalling that $\Lunit = \frac{\partial}{\partial t}$ in geometric coordinates,
	we integrate \eqref{E:GLOBALEXISTENCELUNITMUULUNITMUFIRSTPOINTWISEESTIMATE} with respect to time and use the 
	data estimates
	\eqref{E:MUULUNITMUESTIMATEATTIME0}--\eqref{E:MUULUNITMUESTIMATEATREQUALS1}
	for $\muuLunit \upmu$ to deduce that the following pointwise estimate holds 
	in the away-from-time-axis region 
	$\lbrace (t,u) \ | \ 0 \leq t < \Tboot, u \in (-\infty,\Datafunctionforeikonal(t)] \rbrace$:
	\begin{align}  \label{E:PROOFSTEPGLOBALEXISTENCEMUULUNITMUPOINTWISESTIMATE}
	\begin{split}
	|\muuLunit \upmu|
	& \lesssim 
			\begin{cases}
			\datasize \frac{1}{1 + |u|^2}
			+
			\datasize \frac{\ln \left( \frac{t-u}{-u} \right)}{(1 + |u|)^3}
		+
		\datasize^{1^+} \frac{\ln \left( \frac{t-u}{-u} \right) \ln_+(|u|)}{(1 + |u|)^3},
			& u \leq - 1,
				\\
			\datasize \frac{\ln_+(|u|)}{1 + |u|^2}
			+
			\datasize
			\frac{\ln \left(t-u \right)}{(1 + |u|)^3} 
			+
			\datasize^{1^+}
			\frac{\ln \left(t-u \right) \ln_+(|u|)}{(1 + |u|)^3},
			& u \geq -1,
			\end{cases}
	\end{split}
	\end{align}
	which yields \eqref{E:GLOBALEXISTENCEMUULUNITMUPOINTWISESTIMATE}.

	\medskip
	\noindent \underline{\textbf{Proof of \eqref{E:GLOBALEXISTENCEMUSIMPLERBOUNDS}--\eqref{E:GLOBALEXISTENCEMUXMUSIMPLERUPPERBOUND}}}.
	The estimates \eqref{E:GLOBALEXISTENCEMUSIMPLERBOUNDS}--\eqref{E:GLOBALEXISTENCEMUXMUSIMPLERUPPERBOUND}
	are straightforward consequences of \eqref{E:GLOBALEXISTENCEMUPOINTWISESTIMATE}
	and 
	\eqref{E:GLOBALEXISTENCEMUULUNITMUPOINTWISESTIMATE}
	respectively.

\end{proof}

\subsection{Diffeomorphism properties of the change of variables map}
\label{SS:GLOBALEXISTENCECHANGEOFVARIABLESMAP}
In this section, we prove a simple lemma showing that the change of variables map from geometric to Cartesian coordinates is a diffeomorphism in the 
entire bootstrap region.

\begin{lemma}[Diffeomorphism properties of the change of variables map]
\label{L:GLOBALEXISTENCECHANGEOFVARIABLESMAP}
Assume that the bootstrap assumptions of Sect.\,\ref{SS:GLOBALEXISTENCEBOOTSTRAPANDTWOREGIONS}
hold on the open-at-the-top region 
$\awayglobalexistenceslab{\Tboot} = \lbrace r \geq 1 \rbrace \cap \lbrace 0 \leq t < \Tboot \rbrace$
defined in \eqref{E:AWAYFROMTIMEAXISGLOBALLEXISTENCEBOOTSTRAPSLAB}.
Then if $\datasize$ is sufficiently small and $\Tboot < \infty$, 
the change of variables map
$\Upsilon(t,u) := (t,r)$ extends to a global $C^2$ diffeomorphism from the closed set
$\lbrace (t,u) \ 0 \leq t \leq \Tboot, u \leq \Datafunctionforeikonal(t) \rbrace$
(where $\Datafunctionforeikonal(t)$ is defined in Sect.\,\ref{SS:GLOBALEXISTENCEEIKONALFUNCTION})
onto the closed set $
\overline{\awayglobalexistenceslab{\Tboot}}
=
\lbrace r \geq 1 \rbrace \cap \lbrace 0 \leq t \leq \Tboot \rbrace$.

Similarly, if $\Tboot = \infty$, then $\Upsilon$
is a $C^2$ diffeomorphism from the closed set
$\lbrace (t,u) \ 0 \leq t < \infty, u \leq \Datafunctionforeikonal(t) \rbrace$
onto the closed image set $\lbrace r \geq 1 \rbrace \cap \lbrace 0 \leq t < \infty \rbrace$.

\end{lemma}

\begin{proof}
	We will give the proof only in the case $\Tboot < \infty$ since the case $\Tboot = \infty$ can be handled with similar but simpler arguments.

	First, $\Upsilon$ extends to a $C^2$ function on the domain closure $\lbrace (t,u) \ 0 \leq t \leq \Tboot, u \leq \Datafunctionforeikonal(t) \rbrace$
	because equation \eqref{E:CHOVJACOBIAN},
	the bootstrap assumptions, and the transport equation \eqref{E:MUEVOLUTION}
	imply that $\Lunit \mathrm{d} \Upsilon = \frac{\partial}{\partial t} \mathrm{d} \Upsilon$ is uniformly bounded 
	in $C_{t,u}^1$
	on the open-at-the-top region $\lbrace (t,u) \ 0 \leq t < \Tboot, u \leq \Datafunctionforeikonal(t) \rbrace$.
	Moreover, by \eqref{E:CHOVJACOBIAN}, we have the following formula, which holds on the closed set 
	$\lbrace (t,u) \ 0 \leq t \leq \Tboot, u \leq \Datafunctionforeikonal(t) \rbrace$:
		\begin{align} \label{E:UDERIVATIVEOFRADIALCOORDINATE}
			\frac{\partial}{\partial u} r 
			& = - \Speed \upmu. 
		\end{align}	
		Equations \eqref{E:SPEEDOFSOUNDEXPANSION}--\eqref{E:SPEEDOFSOUNDERRORFUNCTIONVANISHESATORIGIN}
		and the bootstrap assumptions (see in particular \eqref{E:GLOBALEXISTENCEMUSIMPLERBOUNDSBOOTSTRAP})
		imply that RHS~\eqref{E:UDERIVATIVEOFRADIALCOORDINATE} is bounded from above and below by strictly negative constants.
		It follows that for each fixed $t \in [0,\Tboot]$, the map $u \rightarrow r(t,u)$ is strictly decreasing
		and thus the map $\Upsilon$ is injective on $\lbrace (t,u) \ 0 \leq t \leq \Tboot, u \leq \Datafunctionforeikonal(t) \rbrace$.
		In view of \eqref{E:CHOVJACOBIANDTERMINANT}, we conclude that the extended $\Upsilon$ is a $C^2$ diffeomorphism on
		the domain closure $\lbrace (t,u) \ 0 \leq t \leq \Tboot, u \leq \Datafunctionforeikonal(t) \rbrace$, as is desired.
\end{proof}

\subsection{Proof of the main global existence theorem}
\label{SS:PROOFOFMAINGLOBALEXISTENCETHEOREM}
We are now ready to prove Theorem~\ref{T:MAINGLOBALEXISTENCETHEOREM}.
Recall Remark~\ref{R:DISCRETESYMMETRY}, which implies that it suffices to prove global existence and uniqueness for future-times $t \geq 0$. 

By standard local well-posedness and continuation criteria for quasilinear transport equations, 
if $\datasize > 0$ is sufficiently small, then
there exists a maximum time $T_{\textsf{Maximal}} > 0$ such that 
$(\RRiemann,\LRiemann)$ are classical solutions to 
\eqref{E:OUTGOINGRIEMANNINVARIANTEVOLUTION}--\eqref{E:INGOINGRIEMANNINVARIANTEVOLUTION}
for $(t,r) \in [0,T_{\textsf{Maximal}}) \times [0,\infty)$,
such that the eikonal function $u$ and $\upmu$ are classical solutions to \eqref{E:GLOBALEXISTENCEODEFOREXTERIOREIKONALFUNCTION} and
\eqref{E:MUEVOLUTION} respectively
for $(t,r) \in [0,T_{\textsf{Maximal}}) \times [1,\infty)$,
and such that the bootstrap assumptions of Sect.\,\ref{SS:GLOBALEXISTENCEBOOTSTRAPANDTWOREGIONS}
hold with $\Tboot := T_{\textsf{Maximal}}$ and $\eps := \datasize^{3/4}$.
Since the estimates of Props.\,\ref{P:GLOBALEXISTENCERIEMANNINVARIANTSAPRIORIEXTERIORREGIONESTIMATES}
and \ref{P:GLOBALEXISTENCESHARPESTIMATESFORMU} hold with $\Tboot := T_{\textsf{Maximal}}$, if
$T_{\textsf{Maximal}} = \infty$, then also considering Lemma~\ref{L:GLOBALEXISTENCECHANGEOFVARIABLESMAP}, 
we conclude Theorem~\ref{T:MAINGLOBALEXISTENCETHEOREM}.
Hence, we assume for the sake of deriving a contradiction that $T_{\textsf{Maximal}} < \infty$.
Since the estimates of Props.\,\ref{P:GLOBALEXISTENCERIEMANNINVARIANTSAPRIORIEXTERIORREGIONESTIMATES}
and \ref{P:GLOBALEXISTENCESHARPESTIMATESFORMU} hold with $\Tboot := T_{\textsf{Maximal}}$, none of the bootstrap assumptions are saturated
on the bootstrap domain $(t,r) \in [0,T_{\textsf{Maximal}}) \times [0,\infty)$; see Remark~\ref{R:GLOBALEXISTENCESTRICTIMPROVEMENTOFBOOTSTRAP}.
Moreover, Lemma~\ref{L:GLOBALEXISTENCECHANGEOFVARIABLESMAP} guarantees that the change of variables map $\Upsilon(t,u) := (t,r)$ extends
to a $C^2$ diffeomorphism onto the closed set 
$(t,r) \in [0,T_{\textsf{Maximal}}] \times [1,\infty)$, and the formula
\eqref{E:CHOVJACOBIAN} implies that the $C^2$ norm of $\Upsilon$ is uniformly bounded on the corresponding domain
$\lbrace (t,u) \ 0 \leq t \leq T_{\textsf{Maximal}}, u \leq \Datafunctionforeikonal(t) \rbrace$.  
In total, we see that the $C_{t,r}^2$ norm of $(\RRiemann,\LRiemann)$ is uniformly bounded on
$(t,r) \in [0,T_{\textsf{Maximal}}] \times [1,\infty)$, and also obeys the uniform ($r$-weighted-at-the-top-order) $C_{t,r}^2$ estimates
provided by Prop.\,\ref{P:GLOBALEXISTENCERIEMANNINVARIANTSAPRIORIEXTERIORREGIONESTIMATES} in the complementary near-time-axis region
$(t,r) \in [0,T_{\textsf{Maximal}}] \times [0,1]$. Hence, standard continuation criteria and continuity allow us to extend 
$(\RRiemann,\LRiemann)$, $u$, and $\upmu$ as classical solutions to a larger slab 
$(t,r) \in [0,T_{\textsf{Maximal}} + \Delta) \times [0,\infty)$ (for some $\Delta > 0$) on which the bootstrap assumptions hold,
thereby contradicting the definition of $T_{\textsf{Maximal}}$.
We have therefore shown that $T_{\textsf{Maximal}} = \infty$ and 
proved the future-global existence results and estimates stated in Theorem~\ref{T:MAINGLOBALEXISTENCETHEOREM}.

The uniqueness of the solution follows from subtracting two solutions to 
\eqref{E:OUTGOINGRIEMANNINVARIANTEVOLUTION}--\eqref{E:INGOINGRIEMANNINVARIANTEVOLUTION}
that share the same data and applying a standard argument, 
based on Gr\"{o}nwall's inequality and estimates similar to the ones we used to prove 
Props.\,\ref{P:GLOBALEXISTENCERIEMANNINVARIANTSAPRIORIEXTERIORREGIONESTIMATES}
and \ref{P:GLOBALEXISTENCESHARPESTIMATESFORMU},
to conclude that the difference of the solutions vanishes for all time; we omit the details.

\appendix

\section{Extending the results to other profiles}
\label{A:EXTENDRESULTSTOOTHERPROFILES}
In the gradient-blowup results of Theorem~\ref{T:MAINMGHDEXISTENCETHEOREM}, we assumed that the Riemann invariant data were perturbations of
$- \datasize \frac{\arctan(r)}{r}$ for $\datasize > 0$ sufficiently small, 
while in the global existence results of Theorem~\ref{T:MAINGLOBALEXISTENCETHEOREM}, we assumed that 
the Riemann invariant data were perturbations of 
$\datasize \frac{\arctan(r)}{r}$; see 
and
\eqref{E:RPLUSDATAISCLOSETOBACKGROUNDATTIME0}--\eqref{E:RMINUSDATAISCLOSETOBACKGROUNDATTIME0}
and
\eqref{E:GLOBALEXISTENCERPLUSDATAISCLOSETOBACKGROUNDATTIME0}--\eqref{E:GLOBAEXISTENCERMINUSDATAISCLOSETOBACKGROUNDATTIME0}. 
What was fundamentally important for the PDE analysis was not the fine details of the $\arctan$ function, 
but rather only a few key features of it. Hence, in this appendix, we explain how to extend the results of those two theorems to 
other ``profiles.'' 

Our starting point is the following definition of an Admissible Profile Seed function; the function $r \rightarrow \arctan(r)$ 
that we handled in detail in the bulk of the paper is a special case of such a function.

\begin{definition}[Admissible Profile Seed]
\label{D:ADMISSIBLEPROFILESEED}
A function $\similarprofile: \mathbb{R} \rightarrow \mathbb{R}$
is said to be an Admissible Profile Seed if it satisfies the following properties:
\begin{enumerate}
	\item (Regularity of the density at the center of symmetry)
		$\similarprofile$ is odd.
	\item (Sign for positive $r$ values) $\similarprofile(r) > 0$ when $r > 0$.
	\item (Global regularity of the profile\footnote{While our proof of Theorem~\ref{T:MAINMGHDEXISTENCETHEOREM} assumed $C^3$ Riemann invariant data,
	 Theorem~\ref{T:MAINGLOBALEXISTENCETHEOREM} only requires the Riemann invariant data to be $C^2$.}) $\similarprofile \in C^3(\mathbb{R})$.
	\item (Quantitative Monotonicity of the first derivative) There exists a $C \geq 1$ such that
		$\frac{1}{C (1 + r)^2} < \similarprofile'(r) \leq \frac{C}{(1 + r)^2}$ for $r \in [0,\infty)$.
        
	\item (Uniqueness of the crease) On the half-line $[0,\infty)$,
		$\similarprofile'$ has a unique, non-degenerate maximum at a (non-negative) point $r_0$.
	\item (Upper bounds for the second and third derivatives). There exists a $C > 0$ such that for $r \in [0,\infty)$, we have
		$|\similarprofile''(r)| \leq \frac{C}{(1 + r)^3}$ 
		and
		$|\similarprofile'''(r)| \leq \frac{C}{(1 + r)^4}$.
\end{enumerate}
\end{definition}

The following properties are simple consequences of Def.\,\ref{D:ADMISSIBLEPROFILESEED}:

\begin{itemize}
	\item (Non-trivial limit at infinity) 
		$\lim_{r \to \infty} \similarprofile(r)$ exists and is positive and finite. Hence, upon multiplying $\similarprofile$ by a positive constant, 
		we can assume that $\lim_{r \to \infty} \similarprofile(r) = 1$.
		Property $4$ above then implies that $|\similarprofile(r) - 1| \lesssim \frac{1}{1+r}$, in particular, convergence to the limit at the rate 
		$\frac{1}{1 + r}$. 
	\item (Properties relevant for the structure of the singular boundary near the crease) 
		$\similarprofile'(r_0) > 0$, $\similarprofile''(r_0) = 0$, and $\similarprofile'''(r_0) < 0$.
\end{itemize}	

\begin{remark}[Other tail rates]
\label{R:OTHERTAILRATES}
It is conceivable that Theorems~\ref{T:MAINMGHDEXISTENCETHEOREM} and \ref{T:MAINGLOBALEXISTENCETHEOREM} might remain valid if one makes the following change to Properties 4 and 6: for a fixed constant $\lambda\in(1,2]$, assume that there exists a $C\geq 1$ such that $\frac{1}{C (1 + r)^\lambda} < \similarprofile'(r) \leq \frac{C}{(1 + r)^\lambda}$, that $|\similarprofile''(r)|\leq \frac{C}{(1+r)^{\lambda+1}}$, and that $|\similarprofile'''(r)|\leq \frac{C}{(1+r)^{\lambda+2}}$. However, any attempt to prove the theorems with these altered hypotheses would require substantial modifications to the arguments, and here we do not make any attempt to rigorously pursue this line of investigation. 
On the other hand, our proof approach will not work if we assume that $\similarprofile'(r)$ decays faster than $\frac{1}{(1+r)^2}$ as $r\to\infty$. For example, the remainder term in \eqref{E:EXTERIORREGIONPOINTWISEESTIMATEFORMODIFIEDLBARRPLUSDERIVATIVE} cannot be improved,
even if we assume that $\similarprofile'$ has a faster decay; 
this remainder term comes from integrating the second term on 
RHS~\eqref{E:EXTERIORREGIONPOINTWISEESTIMATEFORLDERIVATIMEOFMODIFIEDLBARRPLUSDERIVATIVE} with respect to time. 
That is, we have:
\begin{align} \label{E:EXAMPLEESTIMATE}
    r\muuLunit\RRiemann-2\upmu\RRiemann&=2\datasize\similarprofile'(u)+\frac{\mcl{O}(\datasize^{1^+})}{1+u^2}.
\end{align}
For our proof of Theorem~\ref{T:MAINMGHDEXISTENCETHEOREM} to go through, RHS~\eqref{E:EXAMPLEESTIMATE} must be non-negative.
Thus, we must have $\similarprofile'(r)\geq \frac{1}{C(1+r)^2}$ for some $C>0$.
\end{remark}

\subsection{Extending the results in the shock-forming case}
\label{ASS:MGHDEXTENSION}
In this section, we describe how to extend Theorem~\ref{T:MAINMGHDEXISTENCETHEOREM} to a larger set of initial data tied to an Admissible Seed Profile
$\similarprofile$.

To motivate the discussion, much like in Sect.\,\ref{SSS:MOTIVATIONVIALINEARIZEDSYSTEM},
we consider the following initial data for the linear system \eqref{E:LINEARIZEDOUTGOINGRIEMANNINVARIANTEVOLUTION}--\eqref{E:LLINEARIZEDINGOINGRIEMANNINVARIANTEVOLUTION}:
\begin{align} \label{E:OTHERPROFILEBACKGROUNDDATAFORLINEARSYSYSTEM}
	\flatRRiemann(0,r)
	& = \flatLRiemann(0,r)
	:= - \frac{\datasize \similarprofile(r)}{r}
	:= \similarprofileoverr(r),
\end{align}
where $\datasize > 0$. We refer to $\similarprofileoverr$ as the ``profile.''
As further motivation, we will derive a closed formula for the solution to the linear system. The formula involves the following function:
\begin{align} \label{E:MODIFIEDANTIDERIVATIVEOFOTHERPROFILE}
	\antiderivativesimilarprofile(z)
	& := \int_0^z s \similarprofile'(s) \, \mathrm{d}s.
\end{align}
Note that $\antiderivativesimilarprofile$ is even and satisfies:
\begin{subequations}
\begin{align} \label{E:VANISHESATORIGINMODIFIEDANTIDERIVATIVEOFOTHERPROFILE}
	\antiderivativesimilarprofile(0)
	& = 0,	
		\\
	\antiderivativesimilarprofile'(0)
	& = 0.
	\label{E:VANISHESATORIGINDERIVATIVEOFMODIFIEDANTIDERIVATIVEOFOTHERPROFILE}
\end{align}
\end{subequations}

It is straightforward to check that the 
solution to the linear system 
\eqref{E:LINEARIZEDOUTGOINGRIEMANNINVARIANTEVOLUTION}--\eqref{E:LLINEARIZEDINGOINGRIEMANNINVARIANTEVOLUTION}
launched by the data \eqref{E:OTHERPROFILEBACKGROUNDDATAFORLINEARSYSYSTEM} is:
\begin{subequations}
\begin{align} \label{E:OTHERPROFILELINEARRPLUSSOLUTIONFORMULA}
	\flatRRiemann 
	& = \frac{1}{r^2}
					\frac{\datasize}{2} 
					\left\lbrace
						(t + r) \similarprofile(t-r)
						-
						(t + r) \similarprofile(t+r)
						+
						\antiderivativesimilarprofile(t+r) 
						- 
						\antiderivativesimilarprofile(t-r)
					\right\rbrace
			\\
	\flatLRiemann  \label{E:OTHERPROFILELINEARRMINUSSOLUTIONFORMULA}
	& = 
			\frac{\datasize}{2}
				\frac{1}{r^2}
				\left\lbrace
					(t - r) \similarprofile ((t + r))
					- 
					(t - r) \similarprofile ((t - r))
					+
					\antiderivativesimilarprofile(t-r) 
					- 
					\antiderivativesimilarprofile(t+r)
				\right\rbrace,
\end{align}				
\end{subequations}
and that the following identities hold:
\begin{subequations}
\begin{align} \label{E:OTHERPROFILELINEARIZEDSOLUTIONSPECIALRPLUSCOMBINATIONFORMULA} 
	r
	\flatuLunit
	\flatRRiemann
	- 
	2 \flatRRiemann
	& =  
			2 \datasize \similarprofile'(t-r),
				\\
	r
	\flatLunit
	\flatLRiemann
	+ 
	2 \flatLRiemann
	& =  
		-2 \datasize \similarprofile'(t+r).
		\label{E:OTHERPROFILELINEARIZEDSOLUTIONSPECIALRMINUSCOMBINATIONFORMULA} 
\end{align}
\end{subequations}

We now briefly explain why Theorem~\ref{T:MAINMGHDEXISTENCETHEOREM} extends 
to the case in which the profile \eqref{E:BACKGROUNDPROFILEWITHSMALLAMPLITUDEFACTOR}
is replaced with the profile on RHS~\eqref{E:OTHERPROFILEBACKGROUNDDATAFORLINEARSYSYSTEM},
where $\datasize > 0$ is sufficiently small.
As before, we can perturb this profile within the class of data that 
do not change the decay-in-$r$-rates, as in \eqref{E:RPLUSDATAISCLOSETOBACKGROUNDATTIME0}--\eqref{E:RMINUSDATAISCLOSETOBACKGROUNDATTIME0};
we will ignore these data perturbations because they have no substantial effect on the analysis.

The key point is that in geometric $(t,u)$ coordinates, the solution to the nonlinear problem behaves very similarly to the corresponding solution to the
linear system, where the role of $t-r$ is played by the eikonal function $u$; this is essentially what we showed in
Prop.\,\ref{P:GLOBALEXISTENCERIEMANNINVARIANTSAPRIORIEXTERIORREGIONESTIMATES} in the case that 	$\similarprofile(r) := \arctan(r)$.
In particular, because the decay-in-$r$ properties of the profile $\similarprofileoverr(r)$ are the same as those for
$\profileoverr(r) := \datasize \frac{\arctan(r)}{r}$, for the nonlinear gradient-blowup problem, the decay rates for the nonlinear solution relative to the
geometric coordinates $(t,u)$ will be the same as in Prop.\,\ref{P:GLOBALEXISTENCERIEMANNINVARIANTSAPRIORIEXTERIORREGIONESTIMATES}, 
which allows us to control the nonlinear terms perturbatively. 

Next, we recall
that for the nonlinear problem, 
the gradient-blowup is driven by the term on RHS~\eqref{E:OTHERPROFILELINEARIZEDSOLUTIONSPECIALRPLUSCOMBINATIONFORMULA},
where the role of the argument $t-r$ is played by the eikonal function $u$.
Since $\datasize > 0$ by assumption, the $u$-value at which the gradient-blowup first occurs in the nonlinear problem 
is therefore approximately the value of $u$
at which $\similarprofile'(u)$ is as large as possible. 
In view of Def.\,\ref{D:ADMISSIBLEPROFILESEED}, we see that this $u$-value is unique and is approximately equal to $-r_0$;
in Theorem~\ref{T:MAINMGHDEXISTENCETHEOREM}, we had $\similarprofile(r) := \arctan(r)$ and $r_0=0$; see \eqref{E:UCREASEISSMALL}.
Moreover, in the general case in which the data are \eqref{E:OTHERPROFILEBACKGROUNDDATAFORLINEARSYSYSTEM},
the singular boundary once again will extend all the way out to spatial infinity because
by assumption, $\similarprofile'(z)$ is positive for all $z > 0$, and the gradient-blowup along each characteristic is driven by the positivity of 
RHS~\eqref{E:OTHERPROFILELINEARIZEDSOLUTIONSPECIALRPLUSCOMBINATIONFORMULA}.
To obtain the full structure of the singular boundary and crease, the only other important ingredient is the non-degeneracy of the maximum of
$\similarprofile'$, which in the context of Def.\,\ref{D:ADMISSIBLEPROFILESEED} is the property that $\similarprofile'''(r_0) < 0$.
This assumption leads to the positivity of the analog of the factor
$	\frac{1 - 3u^2}{(1 + u^2)^3}$ on RHS~\eqref{E:TWOXBREVEMUPOINTWISEESTIMATEEXTERIOR} near the crease,
which was crucial for the analysis near the crease, carried out in 
Prop.\,\ref{P:EXISTENCEUPTOCREASEANDSINGULARBOUNDARYANDPORTIONOFCAUCHYHORIZON}
(recall that in Theorem~\ref{T:MAINMGHDEXISTENCETHEOREM}, the $u$-value of the crease was approximately $0$, and this factor is therefore positive there).


\subsection{Extending the results in the global existence case}
\label{ASS:GLOBALEXISTENCEEXTENSION}
In this section, we describe how to extend Theorem~\ref{T:MAINGLOBALEXISTENCETHEOREM} to a larger set of initial data tied to an Admissible Seed Profile.

Let $\similarprofile$ be any admissible seed profile from Def.\ref{D:ADMISSIBLEPROFILESEED}.
Much like in Sect.\,\ref{SSS:MOTIVATIONVIALINEARIZEDSYSTEM},
we consider the following initial data for the linear system \eqref{E:LINEARIZEDOUTGOINGRIEMANNINVARIANTEVOLUTION}--\eqref{E:LLINEARIZEDINGOINGRIEMANNINVARIANTEVOLUTION}
that has the opposite sign compared to \eqref{E:OTHERPROFILEBACKGROUNDDATAFORLINEARSYSYSTEM}:
\begin{align} \label{E:GLOBALEXISTENCEOTHERPROFILEBACKGROUNDDATAFORLINEARSYSYSTEM}
	\flatRRiemann(0,r)
	& = \flatLRiemann(0,r)
	:= \frac{\datasize \similarprofile(r)}{r}
	:= \globalsimilarprofileoverr(r).
\end{align}
We again refer to $\globalsimilarprofileoverr$ as the ``profile.''

For the same reasons as above, the solution to the linear system \eqref{E:LINEARIZEDOUTGOINGRIEMANNINVARIANTEVOLUTION}--\eqref{E:LLINEARIZEDINGOINGRIEMANNINVARIANTEVOLUTION}
launched by the data \eqref{E:GLOBALEXISTENCEOTHERPROFILEBACKGROUNDDATAFORLINEARSYSYSTEM} satisfies
the identities \eqref{E:OTHERPROFILELINEARIZEDSOLUTIONSPECIALRPLUSCOMBINATIONFORMULA}--\eqref{E:OTHERPROFILELINEARIZEDSOLUTIONSPECIALRMINUSCOMBINATIONFORMULA},
but with the signs on the RHSs changed:
\begin{subequations}
\begin{align} \label{E:GLOBALEXISTENCEOTHERPROFILELINEARIZEDSOLUTIONSPECIALRPLUSCOMBINATIONFORMULA} 
	r
	\flatuLunit
	\flatRRiemann
	- 
	2 \flatRRiemann
	& =  
			- 2 \datasize \similarprofile'(t-r),
				\\
	r
	\flatLunit
	\flatLRiemann
	+ 
	2 \flatLRiemann
	& =  
		2 \datasize \similarprofile'(t+r).
		\label{E:GLOBALEXISTENCEOTHERPROFILELINEARIZEDSOLUTIONSPECIALRMINUSCOMBINATIONFORMULA} 
\end{align}
\end{subequations}

We now briefly explain how to extend Theorem~\ref{T:MAINGLOBALEXISTENCETHEOREM}
to the case in which the profile \eqref{E:GLOBALEXISTENCEBACKGROUNDPROFILEWITHSMALLAMPLITUDEFACTOR}
is replaced with the profile on RHS~\eqref{E:GLOBALEXISTENCEOTHERPROFILEBACKGROUNDDATAFORLINEARSYSYSTEM}.
As before, we can perturb this profile within the class of data that  
do not change the decay-in-$r$-rates, as in \eqref{E:GLOBALEXISTENCERPLUSDATAISCLOSETOBACKGROUNDATTIME0}--\eqref{E:GLOBAEXISTENCERMINUSDATAISCLOSETOBACKGROUNDATTIME0};
we will ignore these data perturbations because they have no substantial effect on the analysis.
The main observation is that as in Theorem~\ref{T:MAINGLOBALEXISTENCETHEOREM}, 
since $\datasize > 0$ by assumption, the ``dangerous'' term $- 2 \datasize \similarprofile'$
on RHS~\eqref{E:GLOBALEXISTENCEOTHERPROFILELINEARIZEDSOLUTIONSPECIALRPLUSCOMBINATIONFORMULA} 
has a good sign, i.e., in view of Def.\,\ref{D:ADMISSIBLEPROFILESEED}, we see that it is
everywhere negative. As in Theorem~\ref{T:MAINGLOBALEXISTENCETHEOREM}, this leads to global-in-space
rarefaction, i.e., along each characteristic, $\upmu$ can grow logarithmically in time, but not shrink, as in \eqref{E:GLOBALEXISTENCEMUPOINTWISESTIMATE}.
Hence, when $\datasize > 0$ is sufficiently small, 
the characteristics will never intersect, and the solution exists globally
and exhibits the same decay rates relative to the $(t,u)$-coordinates as in Theorem~\ref{T:MAINGLOBALEXISTENCETHEOREM}.

\section{Global hyperbolicity and MGHDs}
\label{A:GLOBALHYPERBOLICITYANDMGHDS}
In this appendix, we provide an introduction to global hyperbolicity and maximal globally hyperbolic developments. 
The material that is most relevant for the bulk of the paper is Prop.\,\ref{P:EQUIVALENCEOFGLOBALHYPERBOLICITYMMAX},
in which we prove that the ``spherically symmetric'' 
notion of global hyperbolicity of the maximal development that we proved in Point~7 of Theorem~\ref{T:MAINMGHDEXISTENCETHEOREM}
is equivalent to the standard notion of global hyperbolicity from Lorentzian geometry,
i.e., we prove Point~8. We also provide a standard result in Lorentzian geometry, Lemma~\ref{L:CAUCHYHYPERSURFACEDODISENTIREMANIFOLD},
which we used in proving Point~9 of Theorem~\ref{T:MAINMGHDEXISTENCETHEOREM}.
The remaining parts of this appendix provide background material  
and results that illustrate why MGHD uniqueness cannot be taken for granted.
Spherical symmetry is not important for most of this appendix, so unless otherwise stated, symmetry is not assumed. 

\subsection{Lorentzian metrics tied to wave equations and the compressible Euler equations}
\label{SS:LORENTZIANMETRICS}
To keep the exposition short, we mostly focus the discussion on quasilinear wave systems on $\mathbb{R}^{1+3}$ of the following form:\footnote{In \eqref{E:WAVESYSTEM},
we are using Einstein's summation convention, in which lowercase Greek indices are summed from $0$ to $3$. \label{FN:GREEKSPACETIMEEINSTEINSUMMATION}}
\begin{align} \label{E:WAVESYSTEM}
	(\mathbf{g}^{-1})^{\alpha \beta}[\Phi,\partial \Phi] \partial_{\alpha} \partial_{\beta} \Phi
	& = \mathcal{N}[\Phi,\partial \Phi],
\end{align}
where $\Phi$ is the array of unknowns,
$\mathcal{N}$ is a smooth nonlinearity,
and $\mathbf{g}$ is a solution-dependent Lorentzian metric, that is, 
a bilinear form of signature $(-,+,+,+)$.
Below we explain why the $3D$ isentropic and irrotational compressible Euler equations are 
equivalent to a special case of \eqref{E:WAVESYSTEM}.

It is well-known that wave equations exhibit finite speed of propagation. 
By this, we mean that if $\Omega \subset \R^3$ is an open subset on which smooth data are given and  
$\mathscr{D}(\Omega) \subset \R^{1+3}$ is a corresponding domain of dependence (see Def.\,\ref{D:DOD} for a precise definition) 
on which a smooth solution exists, then any other smooth solution that exists on the same set 
$\mathscr{D}(\Omega)$ and that has the same data on $\Omega$ must in fact agree with the first solution on all of  
$\mathscr{D}(\Omega)$. In the case of linear or semilinear wave equations (i.e., when $(\mathbf{g}^{-1})^{\alpha\beta}$ is a constant-coefficient Lorentzian metric in \eqref{E:WAVESYSTEM}), $\mathscr{D}(\Omega)$ does not depend on the solution, but rather  
is completely determined by the set $\Omega$ and the ``known'' coefficients $(\mathbf{g}^{-1})^{\alpha\beta}$,\footnote{For example, if $\mathbf{g}^{-1} = \mbox{diag}(-1,1,1,1)$ is the Minkowski metric, 
as is the case for the standard linear wave operator, 
then $\mathscr{D}(\Omega)$ 
is the set of points $p$ in $\mathbb{R}^{1+3}$ such that all straight lines contained in 
the solid (forward + backward) light cone with tip at $p$ intersect $\Omega$.}
and classical solution non-uniqueness of the type we discuss in Sect.\,\ref{SSS:EPREXAMPLE} never occurs.
However, if the equation is \emph{quasilinear} (i.e., $(\mathbf{g}^{-1})^{\alpha\beta}$ depends on the solution), 
then the \emph{solution} and the corresponding dynamic Lorentzian geometry implicitly determine the domain of dependence. 
In this appendix, especially in Sect.\,\ref{SSS:EPREXAMPLE},
we will exhibit how $\mathscr{D}(\Omega)$'s dependence on the solution introduces surprising subtleties concerning uniqueness and existence of classical solutions to 
\eqref{E:WAVESYSTEM}.

Next, we highlight that in \cites{jLjS2020a,jS2019c}, it was shown that the $3D$ irrotational and isentropic compressible Euler equations are equivalent to a system
of type \eqref{E:WAVESYSTEM},\footnote{The compressible Euler equations \eqref{E:INTROTRANSPORTVI}--\eqref{E:INTROTRANSPORTDENSITY}
cannot generally be recast in the  
form \eqref{E:WAVESYSTEM}, due to the presence of vorticity. Nevertheless, for all classical solutions,
the acoustical metric $\mathbf{g}$ defined in \eqref{E:ACOUSTICALMETRIC} determines the sound wave characteristics in the system, which are the fastest characteristic; vorticity is transported in the direction of the material derivative vectorfield 
\eqref{E:MATERIALDERIVATIVEVECOTRFIELD}, which can be shown to always be timelike with respect to $\mathbf{g}$, i.e., vorticity travels more slowly than sound waves. For these reasons, the notion of global hyperbolicity with respect to the acoustical metric still applies to all classical solutions 
to equations \eqref{E:INTROTRANSPORTVI}--\eqref{E:INTROTRANSPORTDENSITY}, even those with non-vanishing vorticity. \label{FN:GLOBALHYPERBOLICITYMAKESENSEFORALLEULERSOLUTIONS}} where the solution-array $\Phi$ comprises the density\footnote{More precisely, in \cites{jLjS2020a,jS2019c}, the authors derived a wave equation for the log of the density. However, it is straightforward (based on the chain rule) to use that equation to derive a wave equation
for the density itself.} 
$\varrho$ and the Cartesian velocity components $v^1,v^2,v^3$,
and $\mathbf{g} = \mathbf{g}[\Phi]$ is the \emph{acoustical metric}:
\begin{align} \label{E:ACOUSTICALMETRIC}
	\mathbf{g} 
		& := 
		- \mathrm{d} t \otimes \mathrm{d}t
			+ 
			\Speed^{-2} 
			\sum_{a=1}^3(\mathrm{d}x^a - v^a \mathrm{d}t) \otimes (\mathrm{d}x^a - v^a \mathrm{d}t).
\end{align}

The acoustical metric in \eqref{E:ACOUSTICALMETRIC} naturally arises in the study of  
the characteristics of the compressible Euler equations \eqref{E:INTROTRANSPORTVI}--\eqref{E:INTROTRANSPORTDENSITY} without symmetry or irrotationality assumptions. More precisely, straightforward computations (see e.g.\ \cite{dC2007}*{Chapter~1} for the details in the case of the relativistic Euler equations) 
yield that at each point in spacetime, the characteristic subset of the tangent space is
the following subset of vectors, where $\mathbf{g}$ is the acoustical metric and $\Transport$ is the material derivative
vectorfield from \eqref{E:MATERIALDERIVATIVEVECOTRFIELD}:
\begin{align} \label{E:COMPRESSIBLEEUELRCHARSUBSET}
		\lbrace Y \ | \ \mathbf{g}_{\alpha \beta} Y^{\alpha} Y^{\beta}  = 0 \rbrace
		\cup
		\lbrace \mbox{span}(\Transport) \rbrace.
\end{align}
The subset $\lbrace Y \ | \ \mathbf{g}_{\alpha \beta} Y^{\alpha} Y^{\beta}  = 0 \rbrace$ in \eqref{E:COMPRESSIBLEEUELRCHARSUBSET} is the \emph{sound cone} in the tangent space, i.e., the set of vectors that are null with respect to $\mathbf{g}$.
It is straightforward to check that $\Transport$ 
satisfies $\mathbf{g}(\Transport,\Transport)=-1$ and hence the subset $\lbrace \mbox{span}(\Transport) \rbrace$ in \eqref{E:COMPRESSIBLEEUELRCHARSUBSET}
lies strictly in the interior of the sound cone, a geometric representation of the fact that vorticity (which is transported along the integral curves of $\Transport$)
travels more slowly than sound waves. For these reasons, even outside the class of spherically symmetric and irrotational solutions,
the intrinsic geometry of compressible Euler flow is governed by the acoustical metric. Hence, all of the metric-dependent definitions that we introduce below for the wave system~\eqref{E:WAVESYSTEM} extend to the compressible Euler equations \eqref{E:INTROTRANSPORTVI}--\eqref{E:INTROTRANSPORTDENSITY} without symmetry or irrotationality assumptions. 
However, to keep the exposition short, we will often focus the discussion only on \eqref{E:WAVESYSTEM}.

For future use, we note that straightforward computations imply that in spherical symmetry, 
we have the following identity for the inverse acoustical metric:
\begin{align} \label{E:INVERSEACOUSTICALMETRICINSPHERICALSYMMETRY}
	\mathbf{g}^{-1}
	& = - \frac{1}{2} \Lunit \otimes \uLunit
			- 
			\frac{1}{2} \uLunit \otimes \Lunit
			+
			\Speed^2 r^{-2} e_{\mathbb{S}^2}^{-1},
\end{align}
where $\Lunit$ and $\uLunit$ are the null vectorfields from \eqref{E:NULLVECTORFIELDS},
$e_{\mathbb{S}^2}$ is the standard round metric on the Euclidean-unit sphere $\mathbb{S}^2 \subset \mathbb{R}^3$,
and $e_{\mathbb{S}^2}^{-1}$ is the corresponding inverse metric.

\subsection{Basic definitions from Lorentzian geometry}
\label{SS:GHDSBASICDEFINITIONS}
We now provide some basic definitions needed to rigorously define MGHDs. Throughout this section, we assume that
$\mathbf{g}$ is a $C^1$ Lorentzian metric defined on an open set $\mathbf{M} \subset \mathbb{R}^{1+3}$.

\begin{definition}[Spacelike, timelike, null, and extendibility] \label{D:SPACELIKETIMELIKEANDNULL} \hfill
	\begin{itemize} 
		\item \textbf{Vectors:} For any $p \in \mathbf{M}$, we say a tangent spacetime vector $v \in T_p \mathbf{M}$ is:
			\begin{itemize}
				\item spacelike (with respect to $\mathbf{g}$) if $\mathbf{g}(v,v) := \mathbf{g}_{\alpha\beta}(p) v^\alpha v^\beta > 0$.
				\item timelike if $ \mathbf{g}_{\alpha\beta}(p) v^\alpha v^\beta < 0$.
				\item null if $ \mathbf{g}_{\alpha\beta}(p) v^\alpha v^\beta = 0$.
			\end{itemize}
		\item \textbf{Hypersurfaces:} Let $\Sigma \subset \mathbf{M}$ be a hypersurface. Then we say $\Sigma$ is:
			\begin{itemize}
				 \item spacelike if for every $p \in \Sigma$, the normal $\mathbf{N}_p$ (with respect to $\mathbf{g}$)  is a timelike vector. 
				 \item  timelike if for every $p \in \Sigma$, the normal $\mathbf{N}_p$ is a spacelike vector.
				\item null if for every $p \in \Sigma$, the normal $\mathbf{N}_p$ is a null vector.
			\end{itemize}
		\item \textbf{Curves:} Let $\upgamma: I \to \mathbf{M}$ be a piecewise $C^1$ curve,\footnote{By definition, such curves are continuous on $I$, and
			continuously differentiable up to the (possibly empty) boundary of $I$ at all but a finite number of points, and for all $\uplambda_0 \in \mbox{interior}(I)$, $\lim_{\uplambda \downarrow \uplambda_0} \dot{\upgamma}(\uplambda)$
			and $\lim_{\uplambda \uparrow \uplambda_0} \dot{\upgamma}(\uplambda)$ exist. \label{FT:PIECEWISEC1}} 
		where 
		$I \subseteq \R$ is an interval of the form $[a,b], (a,b), [a,b)$, or $(a,b]$ where $-\infty \le a < b \le \infty$. Then:
		\begin{itemize}
			\item We say that $\upgamma$ is timelike (null) if all its tangent vectors $\dot{\upgamma}$
			(on both sides at points where $\dot{\upgamma}$ exhibits a jump discontinuity)
			are timelike (null).
			We say $\upgamma$ is causal if all its tangent vectors (including on both sides at points where $\dot{\upgamma}$ exhibits a jump discontinuity)
			are timelike or null.
			\item We say that $\upgamma$ is \emph{extendible} in $\mathbf{M}$ if at least one of the limits 
			$\lim_{s \to a} \gamma(s)$ or $\lim_{s \to b} \gamma(s)$ exist \emph{in $\mathbf{M}$}. 
			We stress that this notion of extendibility depends on the \emph{domain} $\mathbf{M}$. 
			\item If $\upgamma$ is not extendible, we say that it is inextendible.
		\end{itemize}
	\end{itemize}
\end{definition}

For the next definition, we note that the set of timelike vectors at $p$, i.e., 
$\mathcal{T}(p):= \{ v \in T_p\mathbf{M} \, | \, \mathbf{g}(v,v) < 0\}$, forms a disconnected double-cone. 

\begin{definition}[Time orientation, future, and past] A \emph{time orientation} of the Lorentzian manifold $(\mathbf{M},\mathbf{g})$ is a continuous choice of a positive component $\mathcal{T}(p)$ for each $p \in \mathbf{M}$. We call causal vectors that lie in the closure of the positive component \emph{future-directed}, and those that lie in the closure of the negative component \emph{past-directed}. For each $p \in \mathbf{M}$, we define the timelike (causal) \emph{future} $\mathscr{I}^+(p)$ ($\mathscr{J}^+(p)$) to be the set of points which can be reached from $p$ by a future-directed timelike (causal) curve, i.e., a piecewise $C^1$, timelike (causal) curve whose tangent vectors are always future-directed. 
The timelike (causal) \emph{past} $\mathscr{I}^-(p)$ ($\mathscr{J}^-(p)$) is defined analogously. 
\end{definition}

\begin{definition}[Domain of dependence]
\label{D:DOD}
Let $\mathbf{M} \subset \mathbb{R}^{1+3}$ be an open set, let $\mathbf{g}$ be a $C^1$ Lorentzian metric on $\mathbf{M}$,
and let $\Omega \subset \mathbf{M}$.
We define $\mathscr{D}(\Omega)$, the \emph{domain of dependence} of $\Omega$ (in $\mathbf{M}$) as follows:
\begin{align} \label{E:DODDEFINITION}
	\mathscr{D}(\Omega)
	& := \lbrace q \in \mathbf{M} \ | \ \mbox{every inextendible, piecewise $C^1$, causal curve through $q$ intersects $\Omega$} \rbrace.
\end{align}

\end{definition}

With the ``causality notions'' from above in hand, we can now define global hyperbolicity and Cauchy hypersurfaces.

\begin{definition}[Cauchy hypersurfaces and Globally hyperbolic regions] \label{D:CAUCHYHYPERSURFACESANDGLOBALLYHYPERBOLICREGIONS}
	Let $\mathbf{M} \subset \mathbb{R}^{1+3}$ be an open set, and let $\mathbf{g}$ be a $C^1$ Lorentzian metric on $\mathbf{M}$.
		\begin{itemize}
			\item  We say that a hypersurface $\Sigma \subset \mathbf{M}$ is a \emph{Cauchy hypersurface} 
			(with respect to $\mathbf{g}$) for $\mathbf{M}$ if every timelike, piecewise $C^1$, time-oriented,\footnote{A time-oriented piecewise $C^1$ curve $\upgamma : I \to \mathbf{M}$ is one where, for all $\uplambda \in \textnormal{interior}(I)$, the vectors $ \lim_{\uplambda \downarrow \uplambda_0} \dot{\upgamma}(\uplambda)$ and $\lim_{\uplambda \uparrow \uplambda_0}\dot{\upgamma}(\uplambda)$ all lie in the same connected component of $\mathscr{I}(\upgamma(\uplambda_0))$.} inextendible curve in $\mathbf{M}$ intersects $\Sigma$ exactly once.
			\item We say that $(\mathbf{M},\mathbf{g})$ is \emph{globally hyperbolic} whenever $\mathbf{M}$ contains a Cauchy hypersurface.
		\end{itemize}
\end{definition}

The following is a basic result from Lorentzian geometry; see \cite{bO1983}*{Chapter 14, Lemma~29}
for a proof.

\begin{lemma}[$\mathscr{D}(\Sigma) = \mathbf{M}$ for Cauchy hypersurfaces]
\label{L:CAUCHYHYPERSURFACEDODISENTIREMANIFOLD}
Let $\mathbf{M} \subset \mathbb{R}^{1+3}$ be an open set, and let $\mathbf{g}$ be a $C^1$ Lorentzian metric on $\mathbf{M}$.
Assume that $\Sigma \subset \mathbf{M}$ is a Cauchy hypersurface in the sense
of Def.\,\ref{D:CAUCHYHYPERSURFACESANDGLOBALLYHYPERBOLICREGIONS}. Then $\mathscr{D}(\Sigma) = \mathbf{M}$.
\end{lemma}

\begin{remark}[Inextendible causal curves intersect Cauchy hypersurfaces]
	\label{R:CAUSALCURVESINTERSECTCAUCHYHYPEERSURFACES}
	In our proof of Theorem~\ref{T:MAINMGHDEXISTENCETHEOREM}, we needed Lemma~\ref{L:CAUCHYHYPERSURFACEDODISENTIREMANIFOLD} only because of the
	following simple implication, which follows from
	Def.\,\ref{D:DOD} and the lemma: if $\Sigma$ is a Cauchy hypersurface,
	then any inextendible causal curve in $\mathbf{M}$ must intersect $\Sigma$.
	\end{remark}

\begin{definition}[Globally hyperbolic developments of data on spacelike hypersurfaces] \label{D:GHD}
For equation\footnote{A GHD for the compressible Euler equations \eqref{E:INTROTRANSPORTVI}--\eqref{E:INTROTRANSPORTDENSITY}
is defined in an analogous fashion, where the metric is the acoustical metric \eqref{E:ACOUSTICALMETRIC} and, since the system is first-order, 
the data are $(\varrho,v^1,v^2,v^3)\restriction_{\Sigma}$. The remaining definitions that we provide in this appendix 
for equation \eqref{E:WAVESYSTEM} also extend to
equations \eqref{E:INTROTRANSPORTVI}--\eqref{E:INTROTRANSPORTDENSITY} in an obvious fashion. \label{FN:DEFINITIONSEXTENDTOEULER}} \eqref{E:WAVESYSTEM},
a globally hyperbolic development (GHD) 
of initial data $(\mathring{\Phi}_1,\mathring{\Phi}_2)$ posed on a spacelike hypersurface $\Sigma$ is a pair $(\mathbf{M},\Phi)$, 
where $\mathbf{M}$ is an open neighborhood of $\Sigma$, $\Phi$ is a classical solution\footnote{Here, classical means ``$C^2$.''
One can define GHDs for solutions other than classical ones, e.g., for solutions belonging to function spaces in which local well-posedness holds.} 
to \eqref{E:WAVESYSTEM} defined on $\mathbf{M}$ such that 
$(\mathbf{M},\mathbf{g}[\Phi,\partial\Phi])$ is globally hyperbolic with $\Sigma$ as a Cauchy hypersurface, 
and $(\Phi\restriction_{\Sigma},\mathbf{N}\Phi\restriction_{\Sigma}) = (\mathring{\Phi}_1,\mathring{\Phi}_2)$, 
where $\mathbf{N}$ is the $\mathbf{g}$-unit normal to $\Sigma$.
\end{definition}

The following theorem shows that GHDs always exist for smooth data on spacelike Cauchy hypersurfaces. 
The proof relies on standard local well-posedness for the Cauchy problem for the wave equation \eqref{E:WAVESYSTEM}.

\begin{theorem}[GHDs exist for data on Cauchy hypersurfaces, \cite{fEhRjS2019}*{Theorem 4.75}]
\label{T:EXISTENCEOFGHDS}
	Given smooth initial data $(\mathring{\Phi}_1,\mathring{\Phi}_2)$ for the wave equation \eqref{E:WAVESYSTEM} on a spacelike hypersurface $\Sigma \subset \R^{1+3}$ as in Def.\,\ref{D:GHD}, there exists a GHD $(\mathbf{M},\Phi)$.
\end{theorem}

\begin{remark}[Globally hyperbolic developments for data on other families of hypersurfaces] \label{R:GHDONOTHERHYPERSURFACES}
In Def.\,\ref{D:GHD} for the sake of brevity, we chose to provide the definition of GHD only for \emph{spacelike} Cauchy hypersurfaces,
since our main results concern GHDs of data given on the spacelike hypersurface $\lbrace t = 0 \rbrace$.
However, this
definition could be extended to other kinds of causal Cauchy hypersurfaces. 
For example, a common physically relevant choice is a spacelike-characteristic Cauchy hypersurface, i.e., a piecewise-smooth surface consisting of a spacelike portion and a null portion that intersect transversally along a co-dimension $2$ submanifold. This is the type of Cauchy hypersurface that we studied in the region $I$ in the proof of 
Theorem~\ref{T:MAINMGHDEXISTENCETHEOREM}.  
Another well-known choice of a Cauchy hypersurface consists of a pair of two transversal null hypersurfaces that intersect transversally in a co-dimension $2$ topological sphere. 
In the latter case, the initial value problem for \eqref{E:WAVESYSTEM} is called the \emph{characteristic initial value problem}.

We point out that for non-spacelike Cauchy hypersurfaces, the initial data typically have to be posed in a different way compared to the spacelike case.\footnote{For example, the normal vectors which are null cannot be made to be of unit length, so $(\Phi\restriction_{\Sigma}, \mathbf{N}\Phi\restriction_{\Sigma}) = (\mathring\Phi_1,\mathring\Phi_2)$ does not make sense for characteristic initial data.} Although a full treatment of these issues is beyond the scope of this article, we note that data posed on spacelike-characteristic and characteristic-characteristic Cauchy hypersurfaces typically must satisfy constraint equations as well as compatibility conditions on the co-dimension 2 submanifold on which they intersect transversally.\footnote{For example, let $C_0 := \{t = x^1\}$, $\underline{C}_0 := \{ t = -x^1\}$, and consider $f \in C^2(C_0)$ and $g \in C^2(\underline{C}_0)$. Then it is easy to check that the linear wave equation $-\partial_t^2 \Phi + \partial_{x^1}^2 \Phi = 0$ with $\Phi|_{C_0}=f$ and $\Phi|_{\underline{C}_0} = g$ can be solved classically only if the compatibility condition $f(0,0) = g(0,0)$ holds.} These compatibility requirements would have to be included in the corresponding GHD definitions. However, there is a specific case in which these additional constraints and compatibility conditions are automatically satisfied: \emph{whenever it is already known that the solution exists in an open neighborhood of the Cauchy hypersurface in question}. In this case, the constraint equations and compatibility conditions can be immediately inferred by simply tracing the known solution on the Cauchy hypersurface. We have highlighted these issues because they appeared (in a very mild form in the context of spherical symmetry) in our proof of Theorem~\ref{T:MAINMGHDEXISTENCETHEOREM}.
\end{remark}

The following definition, taken from \cite{fEhRjS2019}*{Section~4.7}, captures a form of uniqueness for a GHD.

\begin{definition}[Unique GHDs]
\label{D:UNIQUEGHD}
Given initial data $(\mathring{\Phi}_1,\mathring{\Phi}_2)$ for equation \eqref{E:WAVESYSTEM}, 
we call a GHD $(\mathbf{M},\Phi)$ a unique globally hyperbolic development (UGHD) iff for
all other GHDs $(\widetilde{\mathbf{M}},\widetilde{\Phi})$ of the same data, we have
we have $\Phi = \widetilde{\Phi}$ on $\mathbf{M} \cap \widetilde{\mathbf{M}}$.
\end{definition}

Intuitively, failure of MGHD uniqueness might happen when $\mathbf{M}$ is too ``big.'' 
However, as we will explain in Sect.\,\ref{SSS:NOGHD},
in principle, there might be smooth data such that not even a single UGHD, however ``small,'' exists!

\begin{definition}[GHD Extendibility]
\label{D:EXTENDIBILITY}
For equation~\eqref{E:WAVESYSTEM},
a GHD $(\mathbf{M},\Phi)$ 
of initial data $(\mathring{\Phi}_1,\mathring{\Phi}_2)$ posed on a spacelike hypersurface $\Sigma$ 
is said to be \emph{extendible} if there exists 
another GHD $(\widetilde{\mathbf{M}},\widetilde{\Phi})$
of the (same) data on $\Sigma$ such that the following hold:
\begin{itemize}
	\item $\mathbf{M} \subsetneq \widetilde{\mathbf{M}}$
	\item $\widetilde{\Phi} \restriction_{\mathbf{M}} = \Phi$
\end{itemize}	
\end{definition}

\begin{definition}[Maximal globally hyperbolic development]  \label{D:MGHD}
For equation~\eqref{E:WAVESYSTEM},
a GHD $(\mathbf{M}^{\textsf{Max}},\Phi^{\textsf{Max}})$ of data on $\Sigma$ 
is said to be a \textbf{maximal globally hyperbolic development} (MGHD) 
if it fails to be extendible in the sense of Def.\,\ref{D:EXTENDIBILITY}.
\end{definition}

The following theorem shows that MGHDs always exist for smooth data on spacelike Cauchy hypersurfaces. 
The proof is based on applying Zorn's lemma to chains of GHDs.

\begin{theorem}[At least one MGHD exists for data on Cauchy hypersurfaces]
\label{T:EXISTENCEOFMGHDS}
	If $(\mathbf{M},\Phi)$ is a GHD of smooth initial data on a Cauchy hypersurface $\Sigma$ as in Def.\,\ref{D:GHD} (and as afforded by Theorem~\ref{T:EXISTENCEOFGHDS}), 
	then there exists at least one MGHD
	$(\mathbf{M}^{\textsf{Max}},\Phi^{\textsf{Max}})$ in the sense of Def.\,\ref{D:MGHD} that 
	is an extension of $(\mathbf{M},\Phi)$ in the sense of Def.\,\ref{D:EXTENDIBILITY}.
	\end{theorem}

\begin{proof}
	Let $\mathfrak{G}$ denote the set of globally hyperbolic developments of the initial data on $\Sigma$
	that are extensions in the sense of Def.\,\ref{D:EXTENDIBILITY} of the fixed GHD $(\mathbf{M},\Phi)$.
 $\mathfrak{G}$ is non-empty since $(\mathbf{M},\Phi)$ is an element.
	Partially order the elements of $\mathfrak{G}$ by extension, i.e., 
	the partial ordering is defined by $(\mathbf{M}_1,\Phi_1) \leq (\mathbf{M}_2,\Phi_2)$
	$\iff$ $\mathbf{M}_1 \subset \mathbf{M}_2$ and $\Phi_2 \restriction_{\mathbf{M}_1} = \Phi_1$.
	Let $\mathfrak{C}$ be a chain, i.e., a totally ordered subset of $\mathfrak{G}$.
	Let $\mathbf{M}^{\textsf{Upper}}$ denote the union of all the solution manifolds of the elements of $\mathfrak{C}$,
	and define the function $\Phi^{\textsf{Upper}}$ by
	$\Phi^{\textsf{Upper}}(p) = \widetilde{\Phi}(p)$ whenever $p \in \widetilde{\mathbf{M}}$
	and $(\widetilde{\mathbf{M}},\widetilde{\Phi}) \in \mathfrak{C}$. Since $\mathfrak{C}$ is totally ordered by assumption, it is straightforward to
	see that $\Phi^{\textsf{Upper}}$ is well-defined.
	
	We will now prove that $(\mathbf{M}^{\textsf{Upper}},\mathbf{g}[\Phi^{\textsf{Upper}},\partial \Phi^{\textsf{Upper}}])$
	is a globally hyperbolic manifold containing $\Sigma$ as a Cauchy hypersurface. To this end, we let 
	$\upgamma:I \to \mathbf{M}^{\textnormal{Upper}}$ be an 
	inextendible $\mathbf{g}[\Phi^{\textsf{Upper}},\partial \Phi^{\textsf{Upper}}]$-timelike, time-oriented, piecewise $C^1$ 
	curve (where $I= (a,b)$ is an open interval with $-\infty \le a < b \le \infty$) in $\mathbf{M}^{\textsf{Upper}}$ that contains
	$p$ as an image point. Then, there exists an $\widetilde{\mathbf{M}}$ such that $(\widetilde{\mathbf{M}},\widetilde{\Phi}) \in \mathfrak{C}$ and
	$p \in \widetilde{\mathbf{M}}$. Let $\widetilde{I} \subset I$ be the largest sub-interval of $I$ 
	such that $p \in \upgamma(\widetilde{I}) \subset \widetilde{\mathbf{M}}$,
	and set $\widetilde{\upgamma} := \upgamma \restriction_{\widetilde{I}}$. Note that $\widetilde{I}$ must be an open interval of the form 
	$(\widetilde{a},\widetilde{b}) \subset (a,b)$
	because $\widetilde{\mathbf{M}}$ is open and $\upgamma$ is continuous.
	We claim that $\widetilde{\upgamma}$ is inextendible in $\widetilde{\mathbf{M}}$. 
	Suppose not. Then one of the limits $\lim_{s \to \widetilde{a}} \widetilde{\upgamma}(a)$ or $\lim_{s \to \widetilde{b}} \widetilde{\upgamma}(s)$ must exist in 
	$\widetilde{\mathbf{M}}$, and therefore also in $\mathbf{M}^{\textnormal{Upper}}$. Without loss of generality, we assume that the former limit exists.
	Since $\widetilde{\mathbf{M}}$ is open, there exists an open ball $U$ containing $\widetilde{\upgamma}(\widetilde{a})$ such that $U \subset \widetilde{\mathbf{M}}$. 
	By the continuity and inextendibility of $\upgamma$, for curve parameter values $s$ in an open neighborhood of $\widetilde{a}$,	
	$\upgamma(s)$ is defined and satisfies $\upgamma(s) \in U$.
	That is, there is an interval $\widetilde{I}'$ such that $\widetilde{I} \subsetneq \widetilde{I}' \subset I$ and 
	such that $\upgamma: \widetilde{I}' \to \widetilde{\mathbf{M}}$. Since $\widetilde{I}$ is a strict subset of $\widetilde{I}'$, 
	this contradicts the 
	maximality of $\widetilde{I}$. We have therefore shown by contradiction that $\widetilde{\upgamma}$ 
	is inextendible in $\widetilde{\mathbf{M}}$. 
	Moreover, since $\Phi^{\textnormal{Upper}}\restriction_{\widetilde{\mathbf{M}}} = \widetilde{\Phi}$, it follows that $\widetilde{\upgamma}$ is
	$\mathbf{g}[\widetilde{\Phi},\partial \widetilde{\Phi}]$-timelike. 
	Thus, since $(\widetilde{\mathbf{M}},\widetilde{\Phi})$ is a GHD of the data on $\Sigma$, $\widetilde{\upgamma}$ must 
	intersect $\Sigma$ exactly once. 
	Clearly this implies that the original curve $\upgamma$ intersects $\Sigma$.
	We will now show that $\upgamma$ cannot intersect $\Sigma$ two or more times. Suppose for the sake of contradiction that
	there exist (finite) numbers $A < B$ such that  
	$[A,B] \subset I$ and such that $\upgamma(A), \upgamma(B) \in \Sigma$.
	$\upgamma([A,B])$ is compact subset of $\mathbf{M}^{\textsf{Upper}}$
	and therefore must be contained in some $\widetilde{\mathbf{M}}$ with
	$(\widetilde{\mathbf{M}},\widetilde{\Phi}) \in \mathfrak{C}$. This is impossible since $\upgamma \restriction_{[A,B]}$ is time-oriented and   
	$\mathbf{g}[\widetilde{\Phi},\partial \widetilde{\Phi}]$-timelike in $\widetilde{\mathbf{M}}$, 
	and therefore can intersect $\Sigma$ at most once.
	In total, we have shown that $\Sigma$ is a Cauchy hypersurface for 
	$\mathbf{M}^{\textsf{Upper}}$
	(with respect to $\mathbf{g}[\Phi^{\textsf{Upper}},\partial \Phi^{\textsf{Upper}}]$)
	in the sense of Definition~\ref{D:CAUCHYHYPERSURFACESANDGLOBALLYHYPERBOLICREGIONS}.
	We have shown that $(\mathbf{M}^{\textsf{Upper}},\Phi^{\textsf{Upper}})$ is a GHD of the data that extends 
	$(\mathbf{M},\Phi)$, which by definition means that $(\mathbf{M}^{\textsf{Upper}},\Phi^{\textsf{Upper}}) \in \mathfrak{G}$.
	By construction, it follows that $(\mathbf{M}^{\textsf{Upper}},\Phi^{\textsf{Upper}})$ is an upper bound for $\mathfrak{C}$,
	that is, for any $(\widetilde{\mathbf{M}},\widetilde{\Phi}) \in \mathfrak{C}$, we have
	$(\widetilde{\mathbf{M}},\widetilde{\Phi}) \leq (\mathbf{M}^{\textsf{Upper}},\Phi^{\textsf{Upper}})$.
	In particular, every chain $\mathfrak{C} \subset \mathfrak{G}$ has an upper bound in $\mathfrak{G}$. 
	Hence, Zorn's lemma yields a maximal element
	$(\mathbf{M}^{\textsf{Max}},\Phi^{\textsf{Max}}) \in \mathfrak{G}$,
	i.e., a GHD of the data on $\Sigma$
	such that any other $(\widehat{\mathbf{M}},\widehat{\Phi}) \in \mathfrak{G}$ 
	that is comparable to $(\mathbf{M}^{\textsf{Max}},\Phi^{\textsf{Max}})$
	satisfies $(\widehat{\mathbf{M}},\widehat{\Phi}) \leq (\mathbf{M}^{\textsf{Max}},\Phi^{\textsf{Max}})$.
	It follows that $(\mathbf{M}^{\textsf{Max}},\Phi^{\textsf{Max}})$ is a MGHD in the
	sense of Def.\,\ref{D:MGHD}.
	
	\end{proof}

\subsection{\texorpdfstring{$\mathbf{M}^{\textsf{Max}}$ is globally hyperbolic}{Global hyperbolicity} in the sense of Sect.\,\ref{SS:GHDSBASICDEFINITIONS}}
\label{SS:EQUIVALNCEOFGHDDEFINITIONS}
In Point~7 of Theorem~\ref{T:MAINMGHDEXISTENCETHEOREM}, in spherical symmetry, 
we gave a definition of global hyperbolicity 
in terms of integral curves of $\Lunit$ and $\uLunit$ and proved that $\mathbf{M}^{\textsf{Max}}$ 
is globally hyperbolic in that sense. 
We deferred the proof of Point~8 of Theorem~\ref{T:MAINMGHDEXISTENCETHEOREM}, which is an equivalent
notion of global hyperbolicity, until the present section. 
In the next proposition, we prove Point~8, which was needed for
the proof of the final point of Theorem~\ref{T:MAINMGHDEXISTENCETHEOREM}.

\begin{proposition}[Equivalent notions of global hyperbolicity]
	\label{P:EQUIVALENCEOFGLOBALHYPERBOLICITYMMAX}
	Let $\mathbf{M}^{\textsf{Max}}$ be the maximal development region from Theorem~\ref{T:MAINMGHDEXISTENCETHEOREM},
	viewed as a subset of $(t,r)$-coordinate space, and let
	$S: \R\times[0,\infty)\times \mathbb{S}^2 \to \R^{1+3}$ be the map from spherical coordinates to Cartesian coordinates, i.e., 
\begin{align} \label{E:MAPFROMSPHERICALCOORDINATESTOCARTESIANCOORDINATES}
	S(t,r,\omega) = (t,r \omega).
\end{align}
	Then $S(\mathbf{M}^{\textsf{Max}} \times \mathbb{S}^2) \subset \mathbb{R}^{1+3}$ 
	is globally hyperbolic with respect to the acoustical metric $\mathbf{g}$ (see \eqref{E:ACOUSTICALMETRIC})
	in the sense of Def.\,\ref{D:CAUCHYHYPERSURFACESANDGLOBALLYHYPERBOLICREGIONS},
	where $S(\Sigma_0 \times \mathbb{S}^2) = \lbrace t = 0 \rbrace \subset \mathbb{R}^{1+3}$ is a Cauchy hypersurface
	in the sense of Def.\,\ref{D:CAUCHYHYPERSURFACESANDGLOBALLYHYPERBOLICREGIONS}.
\end{proposition}

\begin{proof}
\ \\
\noindent \underline{Step 1: Setup for an arbitrary inextendible, timelike, piecewise $C^1$ curve $\upgamma$}.
Let $\upgamma$ be any time-oriented, inextendible, piecewise $C^1$, $\mathbf{g}$-timelike curve in 
$S(\mathbf{M}^{\textsf{Max}} \times \mathbb{S}^2) \subset \mathbb{R}^{1+3}$. 
Note that $S(\mathbf{M}^{\textsf{Max}} \times \mathbb{S}^2)$
is a spherically symmetric open subset of $\mathbb{R}^{1+3}$, 
but the curve $\upgamma$ can have non-trivial angular dependence.
It is straightforward to check that the material derivative vectorfield $\Transport$ defined in \eqref{E:MATERIALDERIVATIVEVECOTRFIELD}
is $\mathbf{g}$-normal to constant-time surfaces and satisfies $\mathbf{g}(\Transport,\Transport)=-1$
and $\Transport t = 1$. That is, $\Transport$ is $\mathbf{g}$-timelike, and $t$ is a time function (its gradient vectorfield is everywhere non-vanishing and $\mathbf{g}$-timelike).
Hence, since $\upgamma$ is a $\mathbf{g}$-timelike curve by assumption, without loss of generality, we can assume that $\upgamma$ is parameterized by 
$t$, i.e., the curve is $t \rightarrow \upgamma(t)$, i.e., the time component $\upgamma^0$ of $\upgamma$ satisfies $\upgamma^0(t) = t$. 
Let $I$ denote the maximal interval of time on which $\upgamma$ is defined and remains in $S(\mathbf{M}^{\textsf{Max}} \times \mathbb{S}^2)$. 
In particular, $\upgamma :I \rightarrow S(\mathbf{M}^{\textsf{Max}} \times \mathbb{S}^2)$.
Without loss of generality, 
we can assume that a time $t_1 > 0$ belongs to the domain $I$ of $\upgamma$, 
and we set $q : = \upgamma(t_1)$ 
and $I^+ := I \cap [0,t_1]$; the case in which there is a $t_1 < 0$ in the domain $I$ can be handled with similar arguments.
Moreover, $\upgamma$ is continuous on $I$ and, away from a finite number of time values $t_* \in I$, is
differentiable with timelike tangent vector $\dot{\upgamma}$. Finally, for such $t_*$ belonging to the interior of $I$,
the ``upper'' tangent vector $\lim_{t \downarrow t_*} \dot{\upgamma}(t)$ and ``lower'' tangent vector $\lim_{t \uparrow t_*} \dot{\upgamma}(t)$ are both assumed to exist and be $\mathbf{g}$-timelike. 
While the upper and lower tangent vectors might disagree, this has no effect on the subsequent analysis.
We will therefore abuse notation and denote either tangent vector by $\dot{\upgamma}$; in
practice, this will not cause any confusion.

\medskip
\noindent \underline{Step 2: Continuous extension of $\upgamma$ to a closed domain}.
To prove the proposition, we must show that as $t$ decreases from $t_1$, $\upgamma(t)$ intersects $\lbrace t= 0 \rbrace$
before escaping $S\left((\mathbf{M}^{\textsf{Max}} \cap \lbrace t \geq 0 \rbrace) \times \mathbb{S}^2 \right)$,
i.e., that $I^+ = [0,t_1]$. 
To this end, we first note that by assumption, we have 
$\dot{\upgamma}^0(t) \equiv 1$ and
$\mathbf{g}(\dot{\upgamma},\dot{\upgamma}) < 0$. 
From this fact, \eqref{E:ACOUSTICALMETRIC}, and the 
estimates for the fluid yielded by Theorem~\ref{T:MAINMGHDEXISTENCETHEOREM},
it follows that for $i=1,2,3$, we have: $\sup_{t \in I^+} |\dot{\upgamma}^i(t)| \lesssim 1$.
Let $t_0$ denote the left endpoint of $I^+$.
From the above bounds, it follows that $\lim_{t \downarrow t_0} \upgamma(t) := \upgamma(t_0)$ exists,
i.e., $I^+ = [t_0,t_1]$, and $\upgamma$ extends continuously to its left endpoint and is piecewise $C^1$ and Lipschitz on $[t_0,t_1]$.
$\upgamma$ is inextendible by assumption, so $\upgamma(t_0)$ must belong to the boundary of
$S\left((\mathbf{M}^{\textsf{Max}} \cap \lbrace t \geq 0 \rbrace) \times \mathbb{S}^2 \right)$.

\medskip

\noindent \underline{Step 3: The projected curve $\hat{\upgamma}$}.
Theorem~\ref{T:MAINMGHDEXISTENCETHEOREM} shows that $\mathbf{M}^{\textsf{Future}} = \mathbf{M}^{\textsf{Max}} \cap \lbrace t \geq 0 \rbrace$.
Here, we are viewing $\mathbf{M}^{\textsf{Future}}$ as a subset of $(t,r)$-coordinate space.
It is easier to work with only in the ``$(t,r)$-components'' of $\upgamma$, which are contained in $\mathbf{M}^{\textsf{Future}}$.
To this end, we define the coordinate projection operator $\Pi(t,r,\omega) := (t,r)$,
and we define: 
\begin{align} \label{E:TIMELIKECURVEPROJECTEDINTOTRCOORDINATESPACE}
	\hat{\upgamma}(t) 
	& := \Pi S^{-1} \circ \upgamma
\end{align}
to be the projection of the curve $\upgamma(t)$ with
domain $I^+ = [t_0,t_1$]. 
Note that $\hat{\upgamma}(t)$ is obtained by expressing the point $\upgamma(t)$ in spherical coordinates and then discarding the angular components,
and that $\hat{\upgamma} : [t_0,t_1] \rightarrow \mathbf{M}^{\textsf{Future}}$.
We also define $\hat{q} : = \hat{\upgamma}(t_1)$, i.e., $\hat{q}$ is the projection of the point $q$ from above.

\medskip

\noindent \underline{Step 4: $\hat{\upgamma}(t_0) \in \partial \mathbf{M}^{\textsf{Future}}$}.
Since we showed in Step 2 that $\upgamma(t_0)$ belongs to the boundary of
$S\left(\mathbf{M}^{\textsf{Future}} \times \mathbb{S}^2 \right)$,
it follows that
$\hat{\upgamma}(t_0)$ must belong to $\partial \mathbf{M}^{\textsf{Future}}$.
By Theorem~\ref{T:MAINMGHDEXISTENCETHEOREM}, 
$\partial \mathbf{M}^{\textsf{Future}}$
comprises four pieces (see also Fig.\,\ref{F:MGHD}):
$\Sigma_0$, 
the future-Cauchy horizon, 
the future-singular boundary,
and a portion of the axis of symmetry $\lbrace r = 0 \rbrace$.
To complete the proof of the proposition, we will show that $\hat{\upgamma}(t_0) \in \Sigma_0$,
i.e., that $t_0 = 0$, and thus the curve $t \rightarrow \hat{\upgamma}(t)$ intersects $\Sigma_0$
(which implies that the original curve $t \rightarrow \upgamma(t)$ intersects $\lbrace t = 0 \rbrace$).

\medskip

\noindent \underline{Step 5: $\hat{\upgamma}$ is $\mathbf{g}$-timelike}.
By assumption, we have $\mathbf{g}(\dot{\upgamma},\dot{\upgamma}) < 0$.
By \eqref{E:INVERSEACOUSTICALMETRICINSPHERICALSYMMETRY}, in spherical coordinates,
the angular components of $\dot{\upgamma}$ make a positive contribution to the expression
$\mathbf{g}(\dot{\upgamma},\dot{\upgamma})$ and thus 
$\mathbf{g}(\dot{\hat{\upgamma}},\dot{\hat{\upgamma}}) \leq \mathbf{g}(\dot{\upgamma},\dot{\upgamma})  < 0$,
i.e., $\hat{\upgamma}$ is $\mathbf{g}$-timelike.

\medskip
\noindent \underline{Step 6: The solid past sound cone of $\hat{q}$, described in $(t,r)$ coordinates}.
Let $t \rightarrow \mathcal{J}(t)$ be the maximally extended integral curve of $\Lunit$ in $\mathbf{M}^{\textsf{Future}}$
emanating from $\hat{q}$, and let $t \rightarrow \underline{\mathcal{J}}(t)$ be the maximally extended integral curve of $\uLunit$ 
in  $\mathbf{M}^{\textsf{Future}}$ emanating from $\hat{q}$.
By Point~5 of Theorem~\ref{T:MAINMGHDEXISTENCETHEOREM}, as $t$ decreases below $t_1$,
$\underline{\mathcal{J}}(t)$ remains in $\mathbf{M}^{\textsf{Future}}$ until it intersects $\Sigma_0$ (at time $0$),
and $\mathcal{J}(t)$ remains in $\mathbf{M}^{\textsf{Future}}$ until it intersects $\Sigma_0$ (at time $0$)
or it intersects the axis $\lbrace r = 0 \rbrace$ at some time $T$ satisfying $0 < T < t_1$.

We will now construct the solid backwards sound cone $\underline{\mathfrak{C}}(\hat{q}) \subset \mathbf{M}^{\textsf{Future}}$ 
that emanates from $\hat{q}$. There are two cases.
Case $i$ is that the curve $t \rightarrow \mathcal{J}(t)$ from above never intersects $\lbrace r = 0 \rbrace$, 
other than possibly at the origin $(t=0,r=0)$; see Fig.\,\ref{F:SOUNDCONECASE1}. 
In this case, it must intersect $\Sigma_0$,
and we define $\underline{\mathfrak{C}}(\hat{q})$ to be solid region bounded on the
left by the $\Lunit$-characteristic curve 
$t \rightarrow \mathcal{J}(t)$, on the right by the $\uLunit$-characteristic curve $t \rightarrow \underline{\mathcal{J}}(t)$,
and on the bottom by the portion of $\Sigma_0$ in between $\mathcal{J}(0)$ and $\underline{\mathcal{J}}(0)$.
Case $ii$ is that curve $t \rightarrow \mathcal{J}(t)$ \emph{does} intersect $\lbrace r = 0 \rbrace$ at some positive time $T$;
see Fig.\,\ref{F:SOUNDCONECASE2}.
In this case, we define $\underline{\mathfrak{C}}(\hat{q})$ as in Case $i$, except the left boundary
is now the portion of the $\Lunit$-characteristic curve $t \rightarrow \mathcal{J}$ joining $\hat{q}$ and the axis $\lbrace r = 0 \rbrace$
in union with the axis portion $\lbrace r = 0 \rbrace \cap \lbrace 0 \leq t \leq T \rbrace$.
In both cases, $\underline{\mathfrak{C}}(\hat{q})$ is the solid backwards sound cone 
in $\mathbf{M}^{\textsf{Future}}$ emanating from $\hat{q}$.
$\underline{\mathfrak{C}}(\hat{q})$ is entirely contained in
$\mathbf{M}^{\textsf{Future}}$ because it is foliated by portions of past-directed integral curves of $\uLunit$ that emanate from the 
left boundary curve $t \rightarrow \mathcal{J}(t)$, and Point~5 of
Theorem~\ref{T:MAINMGHDEXISTENCETHEOREM} implies that all these integral curve portions are contained in 
$\mathbf{M}^{\textsf{Future}}$.

		\begin{figure} 
		\centering	
			\begin{overpic}[scale=.7, grid = false, tics=3, trim=-.5cm -1cm -1cm -.5cm, clip]{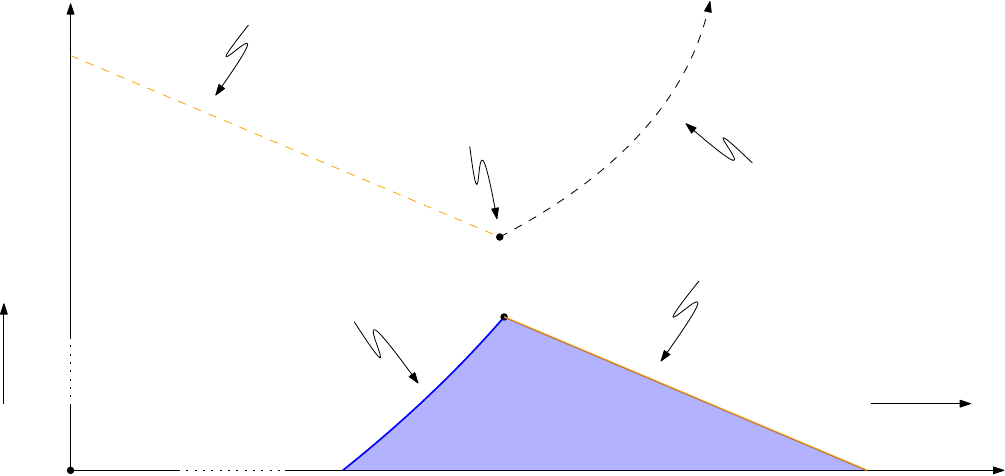} 
            \put (1,16) {\small{$t$}}
            \put (33,19) {\small{$\mathcal{J}$}}
            \put (67,23) {\small{$\underline{\mathcal{J}}$}}
            \put (72,32) {\small{$\futuresinghyp$}}
            \put (43,36) {\small{$\futurecrease$}}
            \put (26,47) {\small{$\futureCauchyhor$}}
            \put (87,13) {\small{$r$}}    
            \put (48.5,21.5) {\small{$\hat{q}$}}
			\end{overpic}
			\caption{The future portion of the MGHD $\mathbf{M}^{\textsf{Future}}$ is the region above $\{t=0\}$ and underneath, but \emph{not} including, $\futuresinghyp\cup\futurecrease\cup\futureCauchyhor$. The backwards sound cone emanating from $\hat{q} \in \mathbf{M}^{\textsf{Future}}$ in Case $i$ is the region shaded in blue.}
			 \label{F:SOUNDCONECASE1}
		\end{figure}

    	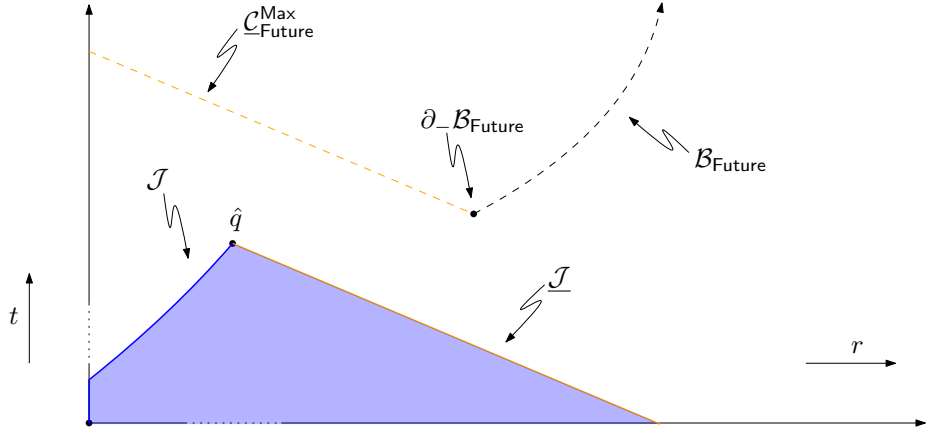
\begin{figure} 
			\centering
			\begin{overpic}[scale=.7, grid = false, tics=3, trim=-.5cm -1cm -1cm -.5cm, clip]{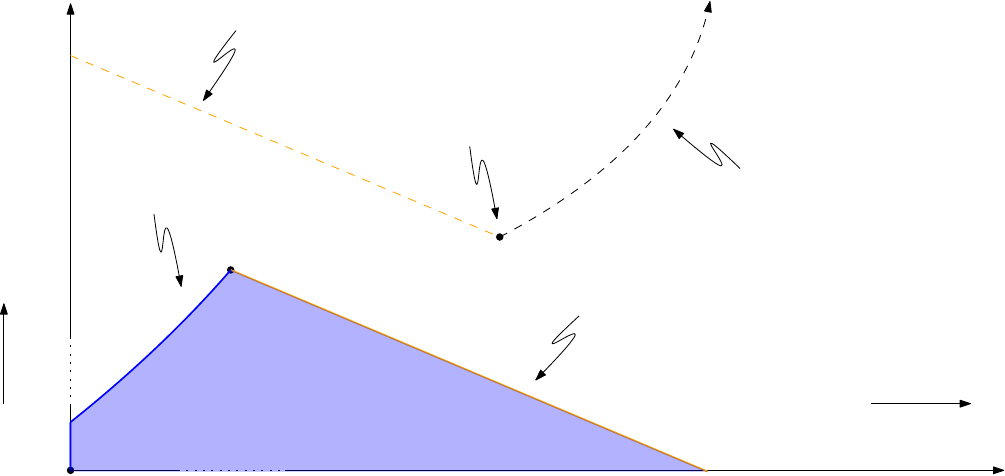} 
            \put (23.5,26) {\small{$\hat{q}$}}
            \put (1,16) {\small{$t$}}
            \put (15,30) {\small{$\mathcal{J}$}}
            \put (56,20) {\small{$\underline{\mathcal{J}}$}}
            \put (71,32) {\small{$\futuresinghyp$}}
            \put (43,36) {\small{$\futurecrease$}}
            \put (25,46) {\small{$\futureCauchyhor$}}
            \put (87,13) {\small{$r$}}
			\end{overpic}
			\caption{The backwards sound cone emanating from $\hat{q} \in \mathbf{M}^{\textsf{Future}}$ in Case $ii$ is the region shaded in blue. The solid blue curve depicting $\mathcal{J}$ is piecewise smooth, causal, and contains a portion of a past-directed integral curve of $\Lunit$ and a portion of $\{r=0\}$}.
			 \label{F:SOUNDCONECASE2}
		\end{figure}

\medskip
\noindent \underline{Step 7: $\hat{\upgamma}([t_0,t_1]) \subset \underline{\mathfrak{C}}(\hat{q})$}.
In the next paragraph, we will prove that $\hat{\upgamma}([t_0,t_1]) \subset \underline{\mathfrak{C}}(\hat{q})$.
Assuming this for a moment, we make the following two key observations:  I) in $(t,r)$ coordinate space,
$\underline{\mathfrak{C}}(\hat{q}) \cap \partial \mathbf{M}^{\textsf{Future}}$
is contained in the disjoint union of $\Sigma_0$ and a portion of the ``positive'' axis of symmetry $\lbrace r = 0 \rbrace \cap \lbrace t > 0 \rbrace$;
and II) $\hat{\upgamma}(t_0) \in \lbrace r = 0 \rbrace \cap \lbrace t > 0 \rbrace$ is impossible,
because that would imply that the original curve $\upgamma$ is extendible as a piecewise $C^1$, $\mathbf{g}$-timelike curve.
Hence, the conclusions of Step 4 about the structure of $\partial \mathbf{M}^{\textsf{Future}}$ imply 
that $\upgamma(t_0) \in \Sigma_0$, which completes the proof of the proposition.

It remains for us to show that $\hat{\upgamma}([t_0,t_1]) \subset \underline{\mathfrak{C}}(\hat{q})$. 
We recall that $\hat{\upgamma}(t_1) = \hat{q}$, 
and we note that at $\hat{q}$, the past-directed tangent vector to $\hat{\upgamma}$, being timelike, 
points into the interior of $\underline{\mathfrak{C}}(\hat{q})$ (this always happens in Case $i$ from Step 6) 
or is parallel to the axis of symmetry $\lbrace r = 0 \rbrace$ 
(this happens in Case $ii$ from Step 6 whenever $\hat{q}$ belongs to the axis of symmetry and 
$\hat{\upgamma}$ is tangent to the axis of at $\hat{q}$).
In particular, for times $t$ less than but close to $t_1$,
$\hat{\upgamma}(t) \in \underline{\mathfrak{C}}(\hat{q})$.
To complete the proof, it suffices to show that for $t$ belonging to the half-open interval $[t_0,t_1)$, $\hat{\upgamma}(t)$
never intersects either of the two characteristic boundary portions of $\underline{\mathfrak{C}}(\hat{q})$.
To this end, we assume for the sake of contradiction that as $t$ decreases strictly below $t_1$, 
$\hat{\upgamma}(t)$ intersects one of the two characteristic boundary portions 
of $\underline{\mathfrak{C}}(\hat{q})$ at some largest time value $a$ satisfying $0 < a < t_1$. 
Let $\hat{p} := \hat{\upgamma}(a)$.
Then $-\dot{\hat{\upgamma}}(a)$ is the past-directed tangent vector to $\hat{\upgamma}$ at $\hat{p}$, 
and being $\mathbf{g}$-timelike, 
$-\dot{\hat{\upgamma}}(a)$ must be transversal to the characteristic boundary.
Since $-\dot{\hat{\upgamma}}(a)$ is also past-directed, it must point \emph{inwards} into $\underline{\mathfrak{C}}(\hat{q})$,
i.e., $\dot{\hat{\upgamma}}(a)$ points \emph{outwards} of $\underline{\mathfrak{C}}(\hat{q})$.
This implies that for times $t$ slightly larger than $a$, $\hat{\upgamma}(t)$ must lie in the complement of
$\underline{\mathfrak{C}}(\hat{q})$, a contradiction. 

\end{proof}

\subsection{Subtleties of existence, uniqueness, and non-uniqueness for GHDs and MGHDs of the data}
\label{SS:GHDEXAMPLESOFEXISTENCEANDNONEXISTENCE}
There is currently no established general theorem for the proving the uniqueness of MGHDs for 
solutions to \eqref{E:WAVESYSTEM} or \eqref{E:INTROTRANSPORTVI}--\eqref{E:INTROTRANSPORTDENSITY}. 
This is in stark contrast to Einstein's equations of mathematical General Relativity, where in the historic works \cites{CB1952,cBgR1969} by Choquet-Bruhat and Choquet-Bruhat--Geroch, 
the Einstein-vacuum equations\footnote{Choquet-Bruhat--Geroch's results in fact hold for many matter models. 
\label{FN:EINSTEINMGHDFORMATTER}}
were formulated as an initial value problem starting from data on a $3D$ manifold $\Sigma$, 
and existence \underline{and} uniqueness of an MGHD was proved. 
Initial data  for Einstein's equations is a triple $(\Sigma,g,k)$, where $\Sigma$ is a connected $3$-dimensional manifold, $g$ is a Riemannian metric on $\Sigma$, and $k$ is a symmetric type $\binom{0}{2}$-tensor on $\Sigma$ satisfying certain constraints.\footnote{The tensors $g$ and $k$ must satisfy certain constraint PDEs obtained from the Gauss--Codazzi equations, Einstein’s equations \eqref{E:EINSTEIN}, and the assumptions that $\iota^* \mathbf{g} = g$ and $\iota^* \mathbf{k} = k$.} The PDE unknowns for the Einstein-vacuum equations are the following two objects: 
\renewcommand{\theenumi}{\Roman{enumi}}
\begin{enumerate}
	\item The Lorentzian manifold $\mathbf{M}$, on which the solution is defined\footnote{A posteriori, the MGHD manifold (in fact, any globally hyperbolic Lorentzian manifold)
	can be shown to be homeomorphic to $I \times \Sigma$, where $I \subset \mathbb{R}$ is an interval.}, together with an embedding $\iota : \Sigma \subset \mathbf{M}$.
	\item The components of a Lorentzian metric $\mathbf{g}$ on $\mathbf{M}$ that  
	satisfies the Einstein-vacuum equations:
	\begin{align}
		\Ric(\mathbf{g}) 
		& =0, \label{E:EINSTEIN}
	\end{align}
where $\Ric(\mathbf{g})$ is the Ricci curvature of $\mathbf{g}$.
Moreover, with $\mathbf{k}$ denoting the second fundamental form of $\iota(\Sigma)$ with respect to $\mathbf{g}$, 
and with $\iota^*$ denoting pullback by $\iota$,
we must have $\iota^* \mathbf{g} = g$ and $\iota^* \mathbf{k} = k$,
i.e., the solution must ``take on the given initial data.''
\end{enumerate}

We now highlight that Point $\textnormal{(I)}$ is very different compared to equation \eqref{E:WAVESYSTEM}, which is a PDE posed 
on the \emph{fixed background} spacetime $\R^{1+3}$. That is, for Einstein's equations, \emph{there is no given background or ``ambient'' spacetime}; the spacetime manifold must be constructed evolutionarily alongside the gravitational metric $\mathbf{g}$ which satisfies \eqref{E:EINSTEIN}. 
It is precisely this distinction that prevents one from directly extending Choquet-Bruhat--Geroch's MGHD existence and uniqueness result 
for the Einstein-vacuum equations to equation \eqref{E:WAVESYSTEM}. For similar reasons, their result is not guaranteed to extend to the compressible Euler equations \eqref{E:INTROTRANSPORTVI}--\eqref{E:INTROTRANSPORTDENSITY}. 
The rest of this appendix is dedicated to further explaining these issues and to showing that it is \underline{impossible} 
to prove a MGHD uniqueness result for general wave equations of type \eqref{E:WAVESYSTEM}. As we will explain, the basic obstacle turns out to be that for
\eqref{E:WAVESYSTEM}, two distinct GHDs of the same data can have a \emph{disconnected} intersection; 
in Sect.\,\ref{SSS:EINSTEINSEQUATIONS}, we explain why this difficulty never arises in an MGHD construction for Einstein's equations. Much of our discussion here is based on \cites{jSb2016,fEhRjS2019}.\footnote{Sbierski's paper \cite{jSb2016} proved the results of \cite{cBgR1969}, but without using Zorn's lemma. }

\subsubsection{Local vs.\ global uniqueness for GHDs for \eqref{E:WAVESYSTEM}} 
\label{SSS:GLOBALUNIQUENESSOFGHDS}
The standard local existence and uniqueness theorems for \eqref{E:WAVESYSTEM}, which are based on energy estimates, 
imply the following ``local'' uniqueness result.

\begin{proposition}[Local uniqueness, \cite{fEhRjS2019}*{Prop. 4.63}] \label{P:LOCALUNIQUENESS}
Given two GHDs $(\mathbf{M}_1,\Phi_1)$ and $(\mathbf{M}_2,\Phi_2)$ of the same initial data of \eqref{E:WAVESYSTEM}, 
there is a \emph{common globally hyperbolic development} (CGHD) $(\mathbf{V},\Phi)$ of the same data such that $\mathbf{V} \subset \mathbf{M}_1\cap \mathbf{M}_2$ with $\Phi = \Phi_1|_{\mathbf{V}}  = \Phi_2|_{\mathbf{V}}$. 
\end{proposition}

Note that Prop.\,\ref{P:LOCALUNIQUENESS} yields a rather \emph{weak} form of uniqueness, for it allows, in principle, 
the existence of a third GHD $(\mathbf{M}_3,\Phi_3)$ of the same data such that there exists a $p \in \mathbf{M}_3 \cap \mathbf{V}$ with 
$\Phi_3(p) \neq \Phi_1(p) = \Phi_2(p)$. This is clearly unsatisfying from the point of view of classical determinism. 

The following theorem is one of the main results of \cite{fEhRjS2019}, which, under a \underline{very strong assumption} that turns out to not always be true,\footnote{The paper \cite{fEhRjS2019} identified a sub-class of equations of type 
\eqref{E:WAVESYSTEM}
for which the assumption always is true. The sub-class
was called ``super-luminal.'' \label{FN:SUPERLUMINAL}} 
upgrades Prop.\,\ref{P:LOCALUNIQUENESS} to ``global'' uniqueness. As we explain in Sect.\,\ref{SSS:UNIQUEMGHDSWITHUNJUSITFIEDCONNECTEDNESS}, 
if one could find a way to guarantee that the assumption always holds,
then this notion of global uniqueness would be sufficient for \emph{both} the existence and uniqueness of the MGHD for \eqref{E:WAVESYSTEM}. 

\begin{theorem}[Global uniqueness under a connectedness assumption, \cite{fEhRjS2019}*{Theorem 4.64}]  \label{T:GLOBALUNIQUENESS}
Let $(\mathbf{M}_1,\Phi_1)$ and $(\mathbf{M}_2,\Phi_2)$  be two GHDs of the same initial data of \eqref{E:WAVESYSTEM}. 
If $\mathbf{M}_1\cap \mathbf{M}_2$ is \textbf{connected}, then $\Phi_1 = \Phi_2$ on $\mathbf{M}_1 \cap \mathbf{M}_2$. 
\end{theorem}

\begin{proof}[Sketch of Proof] 
Complete details are given in \cite{fEhRjS2019}*{Theorem 4.64}. Here, we just sketch the proof to illustrate the main ideas. 
    \begin{itemize} 
        \item[Step 1:] Consider the set $\{(\mathbf{V}_\alpha,\Phi_\alpha) \, | \, \alpha \in A\}$ of all CGHDs of the same data, where $A$ is a non-empty index set by Prop.\,\ref{P:LOCALUNIQUENESS}. Let $\mathbf{V}_{\infty} := \bigcup_{\alpha \in A} \mathbf{V}_\alpha$ and $\Phi_{\infty}(p) = \Phi_\alpha(p)$ for $p \in \mathbf{V}_\alpha$. It is clear that $\Phi_{\infty}$ is well defined and  that $(\mathbf{V}_{\infty},\Phi_{\infty})$ is a CGHD. Moreover, it is clear that $\mathbf{V}_{\infty}$ is globally hyperbolic. The goal is to prove that $\mathbf{V}_{\infty} = \mathbf{M}_1\cap \mathbf{M}_2$. 
        \item[Step 2:] Suppose not and that $\mathbf{V}_{\infty} \subsetneq \mathbf{M}_1\cap \mathbf{M}_2$. Then there exists a boundary point 
				$p \in \partial \mathbf{V}_{\infty} \cap \mathbf{M}_1 \cap \mathbf{M}_2$. 
				\emph{It is precisely the existence of this $p$ where the assumption that $\mathbf{M}_1\cap \mathbf{M}_2$ is connected is used}. 
				It can then be shown (see \cite{fEhRjS2019}*{Theorem 4.64}) that $p$ is a ``spacelike'' point in the sense 
				that there exists a spacelike hypersurface $\mathbf{\Sigma}$ containing $p$ such that
        $\mathbf{\Sigma} \subset \overline{\mathbf{V}_{\infty}} \cap \mathbf{M}_1\cap \mathbf{M}_2$.  Since $\Phi_1$ and $ \Phi_2$ and their derivatives agree on $\mathbf{\Sigma}$, we can uniquely solve the equations (via a standard method based on energy estimates and iterations) 
				in a small neighborhood of $\mathbf{\Sigma}$, obtaining a larger CGHD, thereby contradicting the maximality of $\mathbf{V}_{\infty}$.
    \end{itemize}
\end{proof}

\subsubsection{MGHD uniqueness for \eqref{E:WAVESYSTEM} if the hypotheses of Theorem~\ref{T:GLOBALUNIQUENESS} are known to always hold} 
\label{SSS:UNIQUEMGHDSWITHUNJUSITFIEDCONNECTEDNESS}
To illustrate a point, let us assume that the connectedness hypothesis of Theorem~\ref{T:GLOBALUNIQUENESS} holds for all smooth data and all MGHD pairs
$(\mathbf{M}_1,\Phi_1)$ and $(\mathbf{M}_2,\Phi_2)$. This is a \underline{very strong} assumption, not easy to verify or falsify, 
that fails for some equations; see Sect.\,\ref{SSS:EPREXAMPLE}.
In any case, under this assumption, we will explain why it is fairly easy to both construct and prove the existence and uniqueness of the MGHD for data, 
$\mathbf{M}^{\textsf{Max}} := \bigcup_{\iota \in A} \{\mathbf{M}_\iota \, | \, \iota \in A\}$, where the union is taken over all the spacetime regions of all GHDs $(\mathbf{M}_\iota,\Phi_\iota)$ of the data. Here, $\iota \in A$ is a non-empty index set because at least one GHD always exists; see Theorem~\ref{T:EXISTENCEOFGHDS}. Under our very strong assumption, we can use Theorem~\ref{T:GLOBALUNIQUENESS} to construct a well-defined function $\Phi^{\textsf{Max}}$ on $\mathbf{M}^{\textsf{Max}}$,
as was highlighted in \cite{fEhRjS2019}*{Remark 4.78}. Indeed, given initial data to \eqref{E:WAVESYSTEM}, one can define 
$\Phi^{\textsf{Max}}$ by $\Phi^{\textsf{Max}}(p) := \Phi_\iota(p)$ for $\iota \in A$, which solves
\eqref{E:WAVESYSTEM} because each $\Phi_\iota$ does. The proof of \cite{fEhRjS2019}*{Theorem 4.77} 
then implies that $(\mathbf{M}^{\textsf{Max}},\Phi^{\textsf{Max}})$ is indeed a globally hyperbolic development of the data. 
By construction, $(\mathbf{M}^{\textsf{Max}},\Phi^{\textsf{Max}})$ is maximal and unique. 

The following theorem is the main result from \cite{fEhRjS2019} that we use in this paper to prove uniqueness of the MGHD $\mathbf{M}^{\textnormal{Max}}$ in Theorem~\ref{T:MAINMGHDEXISTENCETHEOREM}. Roughly, an MGHD is unique provided that it lies on one side of its boundary.

\begin{theorem}[Uniqueness condition for MGHDs, \cite{fEhRjS2019}*{Theorem 4.82}] \label{T:UNIQUEONESIDEDNESS}
Let $(\mathbf{M}^{\textnormal{Max}},\Phi^{\textnormal{Max}})$ be an MGHD of initial data posed on a connected spacelike hypersurface $\Sigma$ 
for \eqref{E:WAVESYSTEM} (such an MGHD exists by Theorem~\ref{T:EXISTENCEOFMGHDS}). Assume that for every $p \in \partial \mathbf{M}^{\textsf{Max}}$, there exists a neighborhood $U \subset \R^{1+3}$ of $p$,
a coordinate chart $f : U \to (-\epsilon,\epsilon)^{3}$, $\epsilon > 0$, and a continuous function $g : (-\epsilon,\epsilon)^3 \to (-\epsilon,\epsilon)$ such that
the following ``one-sidedness'' assumptions hold:
    \begin{itemize}
        \item $\partial \mathbf{M}^{\textsf{Max}} \cap U = f^{-1}(\textnormal{graph $g$})$.
        \item All points \emph{below} the graph of $g$ are mapped into $\mathbf{M}^{\textsf{Max}}$.
        \item All points \emph{above} the graph of $g$ are mapped to $\R^{1+3} \setminus \mathbf{M}^{\textsf{Max}}$.
    \end{itemize} 
Then $\mathbf{M}^{\textnormal{Max}}$ is the \emph{unique} MGHD of the data. That is, any other GHD $(\widetilde{\mathbf{M}},\widetilde{\Phi})$ satisfies $\widetilde{\mathbf{M}} \subset \mathbf{M}^{\textnormal{Max}}$ and $\Phi^{\textnormal{Max}}\restriction_{\widetilde{\mathbf{M}}} = \widetilde{\Phi}$.
\end{theorem}

\begin{proof}[Discussion of the ideas in the proof] 
While we will not provide the proof of \cite{fEhRjS2019}*{Theorem 4.82} here, 
we will provide some insight into the role that the ``one-sidedness'' assumptions play in the proof. We will argue by contradiction. That is, 
we suppose that $(\mathbf{M}^{\textnormal{Max}}_1,\Phi^{\textnormal{Max}}_1)$ and $(\mathbf{M}^{\textnormal{Max}}_2,\Phi^{\textnormal{Max}}_2)$ are two 
\underline{distinct} MGHDs of the same data. In particular, $\mathbf{M}^{\textnormal{Max}}_2 \not\subseteq \mathbf{M}^{\textnormal{Max}}_1$.
Let $\mathbf{V}$ be the connected component of $\mathbf{M}^{\textnormal{Max}}_1 \cap \mathbf{M}_2^{\textnormal{Max}}$ 
that contains the Cauchy hypersurface $\Sigma$.
Note that $\partial \mathbf{V} \cap \mathbf{M}^{\textnormal{Max}}_2$ is contained in the complement of $\mathbf{M}^{\textnormal{Max}}_1$,
since points in $\mathbf{M}^{\textnormal{Max}}_1 \cap \mathbf{M}^{\textnormal{Max}}_2$ cannot belong to $\partial \mathbf{V}$.
Moreover, $\partial \mathbf{V} \cap \mathbf{M}^{\textnormal{Max}}_2$ must be non-empty, for otherwise,
$\mathbf{V} \cup \mbox{interior}(\mathbf{V}^c)$ would be a disjoint union of non-empty open sets that contains $\mathbf{M}^{\textnormal{Max}}_2$,
contradicting the connectedness\footnote{ 
$\mathbf{M}^{\textnormal{Max}}_1 \cup \mathbf{M}^{\textnormal{Max}}_2$ is 
connected as long as $\Sigma$ is, for then $\mathbf{M}^{\textnormal{Max}}_1$ and $\mathbf{M}^{\textnormal{Max}}_1$ are connected; 
see the discussion in the Step~3 of the proof discussion of Theorem~\ref{T:GLOBALUNIQEUNESSEINSTEIN}. \label{FN:CONNECTEDNESSOFMGHD}}
of  $\mathbf{M}^{\textnormal{Max}}_2$. Since $\partial \mathbf{V} \subset \partial \mathbf{M}^{\textnormal{Max}}_1 \cup  \partial \mathbf{M}^{\textnormal{Max}}_2$,
it follows that $\emptyset \neq \partial \mathbf{V} \cap \mathbf{M}^{\textnormal{Max}}_2 \subset \partial \mathbf{M}^{\textnormal{Max}}_1$.
Therefore, $(\mathbf{M}^{\textnormal{Max}}_1, \Phi^{\textnormal{Max}}_1)$ extends smoothly from $\mathbf{V}$ to a part of 
$\partial \mathbf{M}^{\textnormal{Max}}_1$. The main idea of the rest of the proof of Theorem~\ref{T:UNIQUEONESIDEDNESS} is to show that
the one-sidedness assumptions in the hypotheses allow one to extend 
$(\mathbf{M}^{\textnormal{Max}}_1, \Phi^{\textnormal{Max}}_1)$ across this portion of $\partial \mathbf{M}^{\textnormal{Max}}_1$ as a GHD,
thereby contradicting the maximality of $(\mathbf{M}^{\textnormal{Max}}_1, \Phi^{\textnormal{Max}}_1)$.
The technical part of the proof involves showing that there is a spacelike hypersurface portion $\Sigma' \subset \overline{\mathbf{V}} \cap \mathbf{M}^{\textnormal{Max}}_2$
that contains at least one of the points in $\partial \mathbf{V} \cap \mathbf{M}^{\textnormal{Max}}_2$ 
(which is a non-empty subset of $\partial \mathbf{M}^{\textnormal{Max}}_1$).
Given this, it is not difficult to extend $(\mathbf{M}^{\textnormal{Max}}_1, \Phi^{\textnormal{Max}}_1)$
to a solution on the larger region $\mathbf{M}^{\textnormal{Max}}_1 \cup \mathscr{D}(\Sigma')$
(where $\mathscr{D}(\Sigma')$ is the domain of dependence of $\Sigma'$ in $\mathbf{V}$), which
can be shown to be globally hyperbolic. 

We stress that without the one-sidedness assumptions, the extension procedure might fail because there might not be any ``room'' left on ``the other side'' of $\partial \mathbf{M}^{\textnormal{Max}}_1$ to construct any extension; see in particular Fig.\,\ref{F:NONUNIQUEMGHDJAN}, 
a context in which the assumptions fail, and the extension procedure is impossible.
\end{proof}

\subsubsection{Einstein's equations}
\label{SSS:EINSTEINSEQUATIONS}
In this section, we give an overview of why for Einstein's equations,
a global uniqueness result -- similar to the one from Theorem~\ref{T:GLOBALUNIQUENESS} for equation \eqref{E:WAVESYSTEM} 
but distinct in ways that we will explain -- always holds; see \cite{jSb2016} for the rigorous details. 
We will denote solutions to Einstein's equations by $(\mathbf{M},\mathbf{g})$ to emphasize that we are looking for GHDs of initial data that satisfy points
$(\textnormal{I})$--$(\textnormal{II})$ from Sect.\,\ref{SS:GHDEXAMPLESOFEXISTENCEANDNONEXISTENCE}. In particular, $\mathbf{M}$ itself is part of the solution.
The main subtlety is that, given two solutions $(\mathbf{M}_1,\mathbf{g}_1)$ and $(\mathbf{M}_2,\mathbf{g}_2)$, 
the intersection $\mathbf{M}_1\cap \mathbf{M}_2$ is a priori not even defined. 
This prevents one from being able to directly apply Theorem~\ref{T:GLOBALUNIQUENESS} to Einstein's equations. 

Despite this difficulty, 
there \emph{does} exist a satisfying way of formulating global uniqueness 
for Einstein's equations that allows for analog of Theorem~\ref{T:GLOBALUNIQUENESS} to be proved
\emph{without the restrictive hypothesis on the connectedness of $\mathbf{M}_1\cap \mathbf{M}_2$}. 
We state this below as Theorem~\ref{T:GLOBALUNIQEUNESSEINSTEIN}.

Before we can discuss Theorem~\ref{T:GLOBALUNIQEUNESSEINSTEIN} for Einstein's equations, we must first discuss 
what it means to ``extend'' GHDs for Einstein's equations.
The main idea is to consider isometric extensions of solutions. The key point is that unlike equations of type \eqref{E:WAVESYSTEM}, 
Einstein's equations are invariant under the diffeomorphism group. That is, if 
$(\mathbf{M},\mathbf{g})$ is a solution to the Cauchy problem $(\textnormal{I})$--$(\textnormal{II})$ 
from Sect.\,\ref{SS:GHDEXAMPLESOFEXISTENCEANDNONEXISTENCE}
and $\psi: \mathbf{M} \to \mathbf{N}$ is any smooth diffeomorphism, then $(\mathbf{N},\psi_* \mathbf{g})$ is also a solution,
where $\psi_* \mathbf{g} = (\psi^{-1})^* \mathbf{g}$ denotes the pushforward of $\mathbf{g}$ by $\psi$.
In this way, one can view solving Einstein's equations as finding an isometry class of Lorentzian manifolds satisfying 
$(\textnormal{I})$--$(\textnormal{II})$.

\begin{definition}[Extensions of GHDs of data for Einstein's equations]
	\label{D:EXTENSIONOFEINSTEINGHDS}
Given two GHDs $(\mathbf{M},\mathbf{g})$ and $(\mathbf{M}',\mathbf{g}')$ of the same initial data for Einstein's equations, 
we say $(\mathbf{M}',\mathbf{g}')$ is an \emph{extension} of $(\mathbf{M},\mathbf{g})$ 
if there exists a time-orientation-preserving isometric embedding\footnote{Isometric embeddings are such that the pullback of $\psi^* \mathbf{g}' = \mathbf{g}$, where $\psi^*$ denotes pullback by $\psi$. \label{FN:ISOMETRICEMBEDDING}} 
$\psi: \mathbf{M} \to \mathbf{M'}$ 
that preserves the initial data. 
\end{definition}

In her celebrated paper \cite{CB1952}, Choquet-Bruhat showed that given initial data for Einstein's equations, 
there exists a GHD $(\mathbf{M},\mathbf{g})$. The same paper also yields a local uniqueness result for Einstein's equations,
i.e., an analog of Prop.\,\ref{P:LOCALUNIQUENESS}, which we present as the following proposition
in the language of ``common globally hyperbolic developments'' (CGHDs). Informally, given a pair of GHDs of the same data, 
a CGHD is a third GHD that is ``contained'' in both of them.

\begin{proposition}[Local uniqueness for  Einstein's equations] \label{P:LOCALUNIQUENESSEINSTEIN}
Given two GHDs $(\mathbf{M}_1,\mathbf{g}_1)$, $(\mathbf{M}_2,\mathbf{g}_2)$ of the same data, there exists a \emph{common globally hyperbolic development} (CGHD) $(\mathbf{M},\mathbf{g})$, that is,\footnote{This definition of a CGHD breaks the diffeomorphism invariance because we are requiring that
$\mathbf{M}$ is a realized as a subset of $\mathbf{M}_1$. However, because solutions are unique only up to diffeomorphism invariance anyway, this does not cause any problems in the theory; see \cite{jSb2016}*{Remark 2.5} for a detailed discussion.} 
$\mathbf{M} \subset \mathbf{M}_1$, and
$(\mathbf{M}_1,\mathbf{g}_1)$ and $(\mathbf{M}_2,\mathbf{g}_2)$ are \emph{both} extensions
of $(\mathbf{M},\mathbf{g})$ in the sense of Def.\,\ref{D:EXTENSIONOFEINSTEINGHDS}.
\end{proposition}

\begin{remark}[Uniqueness up to isometry]
	\label{R:UNIQUENESSUPTOISOMETRY}
	Because of the role played by isometric embeddings in the conclusion of Prop.\,\ref{P:LOCALUNIQUENESSEINSTEIN},
	the proposition can be thought of as ``uniqueness up to isometry.''
\end{remark}

We now provide Theorem~\ref{T:GLOBALUNIQEUNESSEINSTEIN}, which is the main ``global uniqueness'' result for Einstein's equations.
The theorem is proved in detail in \cite{jSb2016}*{Theorem 2.7}], whereas here, we only present some of the main ideas;
we are mainly interested in the differences between Einstein's equations and equation \eqref{E:WAVESYSTEM}.
We stress again that the fundamental difference between Theorems~\ref{T:GLOBALUNIQEUNESSEINSTEIN} and \ref{T:GLOBALUNIQUENESS} is 
that the connectedness of $\mathbf{M}_1\cap\mathbf{M}_2$ 
from Theorem~\ref{T:GLOBALUNIQUENESS} \emph{does not have a natural analog} in the context of Einstein's equations.

\begin{theorem}[Global uniqueness for Einstein's equations, \cite{jSb2016}*{Theorem 2.7}] \label{T:GLOBALUNIQEUNESSEINSTEIN} 
Given two GHDs $(\mathbf{M}_1,\mathbf{g}_1)$ and $(\mathbf{M}_2,\mathbf{g}_2)$ of the same initial data for Einstein's equations, 
there exists a GHD $(\mathbf{M},\mathbf{g})$ which is an extension, in the sense of Def.\,\ref{D:EXTENSIONOFEINSTEINGHDS}, 
of both $(\mathbf{M}_1,\mathbf{g}_1)$ and $(\mathbf{M}_2,\mathbf{g}_2)$.
\end{theorem}

\begin{proof}[Sketch of the proof] \hfill
	\begin{itemize}
		\item[Step 1:] Consider the set $\{ (\mathbf{M}_\alpha,\mathbf{g}_\alpha) \, | \, \alpha \in A\}$ of all CGHDs of $(\mathbf{M}_1,\mathbf{g}_1)$ and $(\mathbf{M}_2,\mathbf{g}_2)$, where $A$ is a non-empty index set by Prop.\,\ref{P:LOCALUNIQUENESSEINSTEIN}. 
		Let $\psi_\alpha : \mathbf{M}_\alpha \to \mathbf{M}_1$ be the canonical embedding (the identity map) of each $\mathbf{M}_\alpha$ into $\mathbf{M}_1$.
		It is not difficult to see that $\mathbf{V}_{\infty} := \bigcup_{\alpha \in A} \mathbf{M}_\alpha$, 
		equipped with the metric $\mathbf{g}$ defined by $\mathbf{g} := \mathbf{g}_1\restriction_{\mathbf{M}_\alpha} = \mathbf{g}_\alpha$,
		is a well-defined CGHD of $(\mathbf{M}_1,\mathbf{g}_1)$ and $(\mathbf{M}_2,\mathbf{g}_2)$; see 
		\cite{jSb2016}*{Theorem 3.4} for a detailed proof.
		In addition, by construction, it is maximal and unique in the following sense: 
		$(\mathbf{V}_{\infty},\mathbf{g})$ is an extension of any other CGHD of $(\mathbf{M}_1,\mathbf{g}_1)$ and $(\mathbf{M}_2,\mathbf{g}_2)$. The reader should compare this portion of the proof with Step 1 of Theorem~\,\ref{T:GLOBALUNIQUENESS}. The main new step needed here is to first embed the $\mathbf{M}_\alpha$ into one of the GHDs in question (indeed, this is ``built into'' the definition of a CGHD from Prop.\,\ref{P:LOCALUNIQUENESSEINSTEIN}), 
		so that the union is a well-defined GHD of the data.
		\item[Step 2:] The most technical part of the paper \cite{jSb2016} is proving the analog of Step 2 from our proof sketch of Theorem~\ref{T:GLOBALUNIQUENESS}. We now describe the main ideas. Let $(\mathbf{V},\mathbf{g})$ be a CGHD of two GHDs $(\mathbf{M}_1,\mathbf{g}_1)$ and $(\mathbf{M}_2,\mathbf{g}_2)$, and let
        $\psi: \mathbf{V} \to \mathbf{M}_2$ be the isometric embedding corresponding 
		to the extension into $\mathbf{M}_2$. Recall that $\mathbf{V} \subset \mathbf{M}_1$ by assumption; see Prop.\,\ref{P:LOCALUNIQUENESSEINSTEIN}.
		We say that boundary points $p_1 \in \partial \mathbf{V}$ and $p_2 \in \partial\psi( \mathbf{V})$ are \emph{corresponding boundary points} if and only if given any open neighborhood $\mathbf{V}_1$ of $p_1$ and $\mathbf{V}_2$ of $p_2$, $\mathbf{V}_1 \cap \psi^{-1}(\mathbf{V}_2 \cap \psi(\mathbf{V})) \neq \emptyset$.\footnote{We clarify that this technical definition is was not needed in Step 2 of Theorem~\ref{T:GLOBALUNIQUENESS} because there, $p \in \partial \mathbf{V}_{\infty} \cap \mathbf{M}_1 \cap \mathbf{M}_2$ was a boundary point in a well-defined intersection (since there,
$\mathbf{M}_1$ and $\mathbf{M}_2$ were both subsets of the same ambient space $\R^{1+3}$).} 
Then, similar (but much more technically challenging) methods as in Step 2 of our proof sketch of Theorem~\ref{T:GLOBALUNIQUENESS} can be used to prove the following: 
if a CGHD $(\mathbf{V},\mathbf{g})$ has corresponding boundary points $p_1,p_2$ in both $\mathbf{M}_1$ and $\mathbf{M}_2$, then there exists a strictly larger extension of  $(\mathbf{V},\mathbf{g})$ that is also a CGHD of  $(\mathbf{M}_1,\mathbf{g}_1)$  and  $(\mathbf{M}_2,\mathbf{g}_2)$; see \cite{jSb2016}*{Theorem 3.6}.\footnote{To construct this larger extension, one exploits the existence of corresponding boundary points $p_1,p_2$ to prove that $p_1 \in \partial\mathbf{V} \cap \psi^{-1}(p_2)$ is a ``spacelike'' point in the sense that there is a spacelike hypersurface $\Sigma \subseteq \overline{\mathbf{V}}$ with $p_1 \in \Sigma$. By locally solving Einstein's equations on $\Sigma$, we obtain a strictly larger CGHD.} Equivalently, if a CGHD of $(\mathbf{M}_1,\mathbf{g}_1)$ and $(\mathbf{M}_2,\mathbf{g}_2)$ is maximal, then there are no corresponding boundary points. 
		
		We note that this Step 2, though similar in spirit to Step 2 of Theorem~\ref{T:GLOBALUNIQUENESS}, is distinct precisely because \emph{we are not trying to prove that $\mathbf{V}_{\infty}=  \mathbf{M}_1\cap \mathbf{M}_2$}; as we have stressed, $\mathbf{M}_1\cap \mathbf{M}_2$ is not even
		a priori well-defined. Instead, we are constructing CGHDs via Prop.\,\ref{P:LOCALUNIQUENESSEINSTEIN},
		using the definitions of GHD extensions given in Def.\,\ref{D:EXTENSIONOFEINSTEINGHDS}.
		\item[Step 3:] Roughly speaking, one finishes the proof of Theorem~\ref{T:GLOBALUNIQEUNESSEINSTEIN} by 
		constructing a common extension of both $\mathbf{M}_1$ and $\mathbf{M}_2$, obtained by ``gluing'' them together along their maximal CGHD manifold
		$\mathbf{V}_{\infty}$. We now provide an overview on how this is achieved in \cite{jSb2016}.
		Let $\psi_1: \mathbf{V}_{\infty} \to \mathbf{M}_1$ and $\psi_2: \mathbf{V}_{\infty} \to \mathbf{M}_2$ be the isometric embeddings corresponding 
		to the fact that $(\mathbf{M}_1,\mathbf{g}_1)$ and $(\mathbf{M}_2,\mathbf{g}_2)$
		are both extensions of $(\mathbf{V}_{\infty},\mathbf{g})$.
		Consider the disjoint union $\mathbf{M}_1 \sqcup \mathbf{M}_2$ endowed with the natural topology. Then for $p,q \in \mathbf{M}_1 \sqcup \mathbf{M}_2$, we define $p \sim q$ if and only if $p = \psi_1(\psi_2^{-1}(q))$, or $q = \psi_2(\psi_1^{-1}(p))$, or $p = q$. \cite{jSb2016}*{Theorem 2.7} proves that the quotient manifold $\mathbf{M}:=( \mathbf{M}_1 \sqcup \mathbf{M}_2)/{\sim}$ is the desired GHD which extends both $ \mathbf{M}_1$ and $\mathbf{M}_2$. In particular, the most subtle step is to prove that $\mathbf{M}$ is Hausdorff, which Sbierski proves using Step 2; see \cite{jSb2016}*{Sect.~3.3}. The connectedness of $\mathbf{M}$  
		follows easily from that of the connectedness of the Cauchy hypersurface $\Sigma$; for it is easy to prove that $\mathbf{M}$ is globally hyperbolic, and all globally hyperbolic Lorentzian manifolds with a Cauchy hypersurface $\Sigma$ are diffeomorphic to $\mathbb{R} \times \Sigma$ 
		\cite{aBsM2003}
		and hence are path-connected whenever $\Sigma$ is.
\end{itemize}
\end{proof}

We are now ready to state the main theorem from \cite{jSb2016}.
\begin{theorem}[Existence and uniqueness of an MGHD for Einstein's equations; \cite{jSb2016}*{Theorem 2.8}] \label{T:UNIQUEMGHDEINSTEIN}
	Given initial data for Einstein's equations, there exists a GHD $(\mathbf{M}^{\textsf{Max}},\mathbf{g}^{\textsf{Max}})$
	that is an extension, in the sense of Def.\,\ref{D:EXTENSIONOFEINSTEINGHDS},
	of any other GHD of the same data. This GHD is unique up to isometry,
	and is called the MGHD of the data.
\end{theorem}

\begin{proof}[Discussion of the ideas in the proof]
At a broad level, to prove Theorem~\ref{T:UNIQUEMGHDEINSTEIN},
one adapts the ideas in the paragraph at the start of Sect.\,\ref{SSS:UNIQUEMGHDSWITHUNJUSITFIEDCONNECTEDNESS}. 
The key difference is that the ``maximal element'' is not obtained by unions, but rather by gluing all GHDs together along their maximal CGHD. 
There is one important caveat in this procedure, explained in \cite{jSb2016}*{Sect.~3.3}: the collection of all GHDs is not a set, 
but rather forms a \emph{proper class}. This issue did not arise in Step 2 of our proof sketch of Theorem~\ref{T:GLOBALUNIQEUNESSEINSTEIN} because there, 
the union was taken over the collection of GHDs that embed into the \emph{fixed} manifolds $\mathbf{M}_1$ and $\mathbf{M}_2$, and this collection \emph{is} a set. 
In proving \cite{jSb2016}*{Theorem 2.8}, Sbierski does \emph{not} consider the collection of all GHDs 
of the initial data on the Cauchy hypersurface $\Sigma$, 
but rather, GHDs whose spacetime manifolds are subsets of $\Sigma \times \R$. In \cite{jSb2016}*{Sect.~3.3}, it is shown that 
this collection \emph{is} a set, denoted here by $\mathfrak{X}$, and one can therefore ``glue''  all the maximal CGHDs as in the proof of 
Theorem~\,\ref{T:GLOBALUNIQEUNESSEINSTEIN}, thereby constructing a maximal element $(\mathbf{M}^{\textsf{Max}},\mathbf{g}^{\textsf{Max}}) \in \mathfrak{X}$. 
The proof of maximality and uniqueness -- up to isometry -- 
concludes by proving that any GHD of the same initial data isometrically embeds into the specific one $(\mathbf{M}^{\textsf{Max}},\mathbf{g}^{\textsf{Max}})$.
\end{proof}

We close this section by highlighting two important difficulties that would arise if one attempted to modify the proof of Theorem~\ref{T:GLOBALUNIQEUNESSEINSTEIN} 
to apply to the wave system \eqref{E:WAVESYSTEM} in the case that $\mathbf{M}_1\cap\mathbf{M}_2$ is disconnected in the ambient spacetime $\R^{1+3}$. 
First, one can show \cite{fEhRjS2019} that given two GHDs for equation \eqref{E:WAVESYSTEM} as in the hypotheses of Prop.\,\ref{P:LOCALUNIQUENESS},
the maximal CGHD of $\mathbf{M}_1 \cap \mathbf{M}_2$ is only the connected component that contains the Cauchy hypersurface. 
Second, we expect that if one tried to glue $\mathbf{M}_1$ and $\mathbf{M}_2$ along their maximal CGHD, 
this would yield a solution defined \emph{not} on a subset of $\R^{1+3}$, 
but rather on a manifold which projects down on $\R^{1+3}$ 
and contains the connected components of $\mathbf{M}_1 \cap \mathbf{M}_2$ twice. 
\begin{quote}
	It is an interesting open problem to rigorously construct and make sense of such ``multi-valued'' solutions to \eqref{E:WAVESYSTEM}.
\end{quote}

\subsubsection{The E--R--S example of MGHD non-uniqueness for a quasilinear wave equation}
\label{SSS:EPREXAMPLE}
In this final section, we present a fascinating result from \cite{fEhRjS2019}*{Theorem 3.45}, 
discussed already in the introduction,
which provides a constructive example of a quasilinear wave \eqref{E:WAVESYSTEM} and two GHDs $(\mathbf{M}_1,\Phi_1)$,$(\mathbf{M}_2,\Phi_2)$ of the same smooth data such that $\mathbf{M}_1\cap\mathbf{M}_2$ is disconnected and such that $\Phi_1$ \emph{disagrees} with $\Phi_2$ on part of 
$\mathbf{M}_1 \cap \mathbf{M}_2$.

\begin{figure}  
\centering
\begin{overpic}[scale=.65, grid = false, tics=5, trim=-.5cm -1cm -1cm -.5cm, clip]{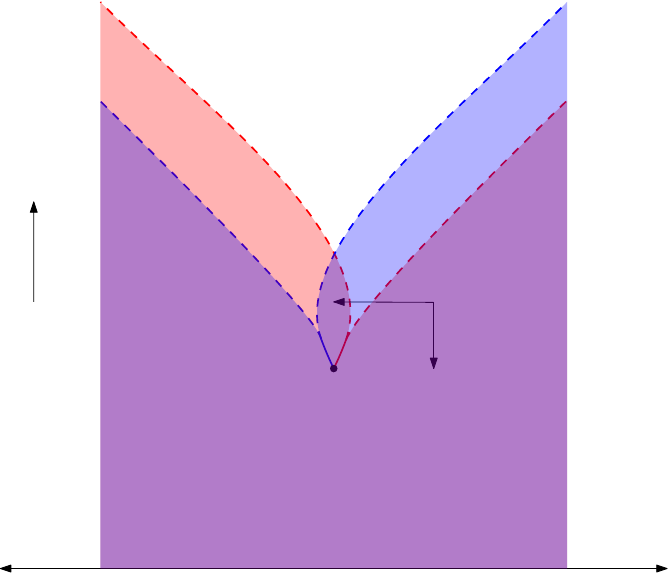}

\put (60,45) {\small $\mathbf{M}_1\cap\mathbf{M}_2$}   
\put (6,44) {\small $t$}
\put (46,20) {\small $\mathbf{V}$}  
\put (48,3) {\small $\Sigma_0$}  

\end{overpic}                                                
\caption{The non-unique GHDs of \cite{fEhRjS2019} for the same data on $\Sigma_0$.
The GHD $\mathbf{M}_1$ is the subset of $\mathbb{R}^{1+1}$ between the dashed red curves and $\Sigma_0$, but does \emph{not} include the solid red curve or black dot. That is, $\mathbf{M}_1$ is the union of the red and purple regions without the solid red curve and black dot. The GHD $\mathbf{M}_2$ is the subset of $\mathbb{R}^{1+1}$ between the dashed blue curves and $\Sigma_0$, but does \emph{not} include solid blue curve or black dot. That is, $\mathbf{M}_2$ is the union of the blue and purple regions without the solid blue curve and black dot.
$\mathbf{M}_1 \cap \mathbf{M}_2$ is disconnected. The connected component of $\mathbf{M}_1 \cap \mathbf{M}_2$ that contains $\Sigma_0$
is denoted by $\mathbf{V}$. $\Phi_1$ and $\Phi_2$ agree on $\mathbf{V}$.
Both solutions are singular at the black dot, but ${\partial}\Phi_1$ blows up precisely on the portion of the solid red curve in $\mathbf{M}_1\cap \mathbf{M}_2$ while $\Phi_2$ is smooth there. In contrast, ${\partial}\Phi_2$ blows up on the portion of the solid blue curve in $\mathbf{M}_1\cap \mathbf{M}_2$ while $\Phi_1$ is smooth there, and both solutions have the same data: $\Phi_1|_{\Sigma_0} = \Phi_2|_{\Sigma_0}$ and $\partial_t \Phi_1|_{\Sigma_0} = \partial_t \Phi_2|_{\Sigma_0}$.} 
\label{F:NONUNIQUEMGHDJAN}        
\end{figure}
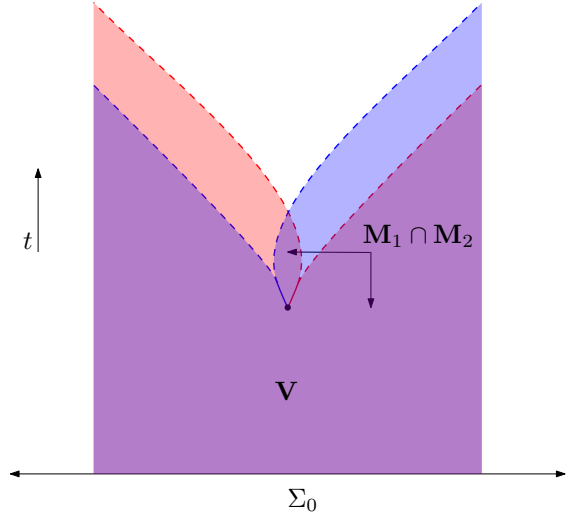

\begin{theorem}[\cite{fEhRjS2019}*{Theorem 3.45}]  \label{T:NONUNIQUEMGHD} 
Consider the following quasilinear wave equation on $\mathbb{R}^{1+1}$: 
    \begin{align} \label{E:BORNINFELDEQUATION}
        (\mathbf{h}^{-1})^{\alpha\beta}[\partial\Phi] \partial_{\alpha}\partial_{\beta}\Phi = 0,
    \end{align}
where relative to the standard Cartesian coordinates $(x^0,x^1 ) := (t,x^1)$,
$(\mathbf{h}^{-1})^{\alpha\beta}$ are the components of the following inverse Lorentzian metric:
		\begin{align}
			\mathbf{h}^{-1} 
			& := \begin{pmatrix}
				(\mathbf{h}^{-1})^{00} & (\mathbf{h}^{-1})^{01} \\ 
				(\mathbf{h}^{-1})^{10} \Phi & (\mathbf{h}^{-1})^{11}
				\end{pmatrix}
			:= \begin{pmatrix}
				-(1+(\partial_1 \Phi)^2) & \partial_t\Phi \partial_1 \Phi \\ 
				\partial_t\Phi \partial_1 \Phi & 1 - (\partial_t \Phi)^2
				\end{pmatrix}. \label{E:BORNINFELDSCALARMETRIC}
		\end{align}
Then there exists $C^{\infty}$ Cauchy initial data $(\Phi|_{t=0},\partial_t \Phi|_{t=0}) = (\mathring{\Phi}_1,\mathring{\Phi}_2)$ for 
\eqref{E:BORNINFELDEQUATION} and two GHDs $(\mathbf{M}_1,\Phi_1)$,$(\mathbf{M}_2,\Phi_2)$ such that $\mathbf{M}_1 \cap \mathbf{M}_2$ is disconnected and such that there exists a $p \in \mathbf{M}_1 \cap \mathbf{M}_2$ with $\Phi_1(p) \neq \Phi_2(p)$; see Fig.\,\ref{F:NONUNIQUEMGHDJAN}.
\end{theorem}

Equation \eqref{E:BORNINFELDEQUATION} is known as the \emph{Born--Infeld equation} (BIE). It appears as a model in various contexts,  including nonlinear electrodynamics and cosmological string theory. It is well known that the wave equation \eqref{E:BORNINFELDEQUATION}
satisfies a very strong ``\emph{double} null condition,'' i.e., there are no resonant interaction terms in the equations,
and in particular, the non-uniqueness from Theorem~\ref{T:NONUNIQUEMGHD}
is not tied to the formation of a shock. We do point out, however, that in the example of 
Theorem~\ref{T:NONUNIQUEMGHD}, each of the two solutions forms a singularity that is more severe than a shock, corresponding to the blowup
of the Lorentzian metric itself, e.g., the gradients $\partial \Phi_1$ and $\partial \Phi_2$ both blow up at the black dot in
Fig.\,\ref{F:NONUNIQUEMGHDJAN}; see \cite{fEhRjS2019}*{Section~3.6} for further details.\footnote{To clarify, the reason we call the gradient blow-up of the non-unique solutions of Theorem~\ref{T:NONUNIQUEMGHD} more ``severe'' than a shock -- even though gradient blow-up also occurs in shock formation for compressible Euler -- is because the components of the metric itself \eqref{E:BORNINFELDSCALARMETRIC} blow up. In contrast, our main result \eqref{T:MAINMGHDEXISTENCETHEOREM} implies that the acoustical metric 
$\mathbf{g}$ \eqref{E:ACOUSTICALMETRIC} extends as a well-defined continuous Lorentzian metric to the entire singular boundary and crease. More generally, ``shock formation'' results for solutions to a quasilinear wave equation of the form $(\mathbf{h}^{-1})^{\alpha\beta}[\partial\Phi]\partial_\alpha\partial_\beta\Phi = \mathcal{N}[\Phi,\partial\phi]$ such as \cites{sA1999a,sA1999b} feature blow-up of $\partial^2 \Phi$, and not $\partial \Phi$.}

We also highlight that $\mathbf{M}_1$ and $\mathbf{M}_2$ \emph{fail} the ``one-sidedness'' hypothesis of Theorem~\ref{T:UNIQUEONESIDEDNESS} for guaranteeing uniqueness because they both lie on ``both sides'' of the solid red and blue curves in Fig.\,\ref{F:NONUNIQUEMGHDJAN}. To provide further insight into how this ``two-sidedness'' is tied to the lack-of-uniqueness, we refer to our discussion of the ideas in the proof of 
Theorem~\ref{T:UNIQUEONESIDEDNESS}. In particular, in the context of the proof of
Theorem~\ref{T:UNIQUEONESIDEDNESS}, $\mathbf{M}_1$ lied on one side of its boundary. Hence, in a proof-by-contradiction-argument, it
was shown to be possible to construct a spacelike hypersurface portion $\Sigma'$ containing some point in $\partial \mathbf{M}_1$,
which allows one to extend the solution to the larger region $\mathbf{M}_1 \cup \mathscr{D}(\Sigma')$,
where $\mathscr{D}(\Sigma')$ is the domain of dependence of $\Sigma'$ in $\mathbf{V}$
($\mathbf{V}$ is the connected component of $\mathbf{M}_1\cap\mathbf{M}_2$ containing $\Sigma_0$).
In the context of Fig.\,\ref{F:NONUNIQUEMGHDJAN}, this extension procedure fails in the sense that \emph{$\mathbf{M}_1$ is already defined on both sides of $\partial \mathbf{M}_1$,} so there is no ``room'' left to extend.

The last sentence of Theorem~\,\ref{T:NONUNIQUEMGHD} is a direct example of a breakdown of classical determinism due to non-uniqueness. By this, we mean that neither the equation nor the dynamics indicate a preferred choice between $\Phi_1$ and $\Phi_2$. In addition, we note that the non-unique GHDs of Theorem~\ref{T:NONUNIQUEMGHD} are 
not maximal. However, by Theorem~\ref{T:EXISTENCEOFMGHDS}, there exists an MGHD 
$(\mathbf{M}_1^{\textsf{Max}},\Phi_1^{\textsf{Max}})$ that is an extension of
$(\mathbf{M}_1,\Phi_1)$. Similarly, there exists an MGHD $(\mathbf{M}_2^{\textsf{Max}},\Phi_2^{\textsf{Max}})$ that extends
$(\mathbf{M}_2,\Phi_2)$ and therefore cannot be the same as $(\mathbf{M}_1^{\textsf{Max}},\Phi_1^{\textsf{Max}})$.
In particular, equation \eqref{E:BORNINFELDEQUATION} admits smooth data that lead to non-unique MGHDs.

We close this section by highlighting the following interesting open problem: 
\begin{quote}
	Given smooth data for equation \eqref{E:WAVESYSTEM}, can one prove the existence of an MGHD \emph{without} invoking Zorn's lemma, 
	just as Sbierski did for Einstein's equations \cite{jSb2016}.
\end{quote}

\subsubsection{There might not even be a single unique GHD}
\label{SSS:NOGHD} 
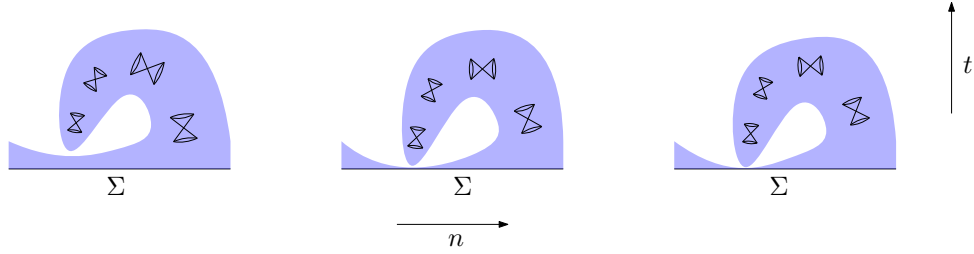
\begin{figure}  
\centering
\begin{overpic}[scale=.65, grid = false, tics=5, trim=-.5cm -1cm -1cm -.5cm, clip]{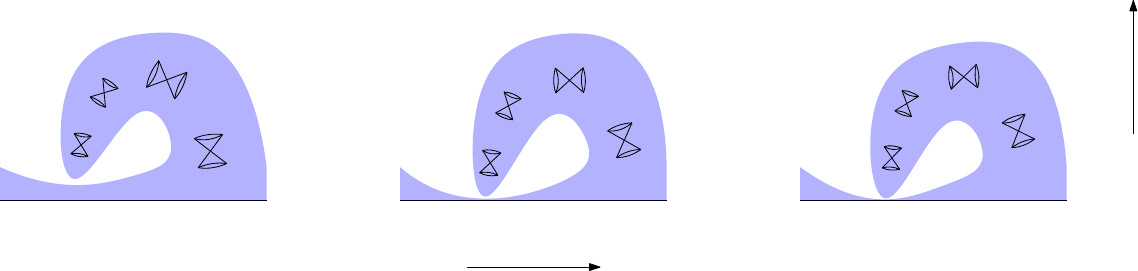}

\put (46,8) {\small{$\Sigma$}}
\put (77,8) {\small{$\Sigma$}}
\put (12,8) {\small{$\Sigma$}}
\put (96,20) {\small{$t$}}
\put (45.5,3) {\small{$n$}}

\end{overpic}                                                
\caption{This figure, which is adapted from \cite{fEhRjS2019}*{Figure 16}, depicts a potential mechanism for resolving Conj.\,\ref{CONJ:NONEXISTENCEOFMGHDS}. The blue regions represent a sequence of GHDs $\{(\mathbf{M}_n,\Phi_n)\}_{n=1}^\infty$ of the same initial data posed on a Cauchy hypersurface $\Sigma$. As $n\to +\infty$, the causal regions of the GHDs get arbitrarily close to the Cauchy hypersurface. Hence, there would not exist \underline{any} open set containing $\Sigma$ on which the data determines a unique solution.
} 
\label{F:NOUNIQUEMGHDATALL}        
\end{figure}

The discussion in Sect.\,\ref{SSS:EPREXAMPLE} was primarily about non-uniqueness of a \underline{maximal} GHD for a wave equation. 
In principle, even worse pathologies could arise. As an example, we highlight the following conjecture,  
which was originally stated as \cite{fEhRjS2019}*{Conjecture 4.81}:

\begin{conjecture} \label{CONJ:NONEXISTENCEOFMGHDS}
There exist quasilinear wave systems of the form \eqref{E:WAVESYSTEM} and smooth initial data for which there does \emph{not} exist any unique GHD
in the sense of Def.\,\ref{D:UNIQUEGHD}; see Fig.\,\ref{F:NOUNIQUEMGHDATALL}. 
\end{conjecture}

Informally speaking, for an example that verifies Conj.\,\ref{CONJ:NONEXISTENCEOFMGHDS}, it would not be possible to \underline{uniquely} classically 
solve the Cauchy problem in \underline{any} neighborhood -- however small -- of the Cauchy hypersurface on which the data are posed.


\bibliography{JBib}

@article {aBmS2006,
    AUTHOR = {Bernal, Antonio N. and S\'anchez, Miguel},
     TITLE = {Further results on the smoothability of {C}auchy hypersurfaces
              and {C}auchy time functions},
   JOURNAL = {Lett. Math. Phys.},
  FJOURNAL = {Letters in Mathematical Physics},
    VOLUME = {77},
      YEAR = {2006},
    NUMBER = {2},
     PAGES = {183--197},
      ISSN = {0377-9017,1573-0530},
   MRCLASS = {53C50 (53C80 81T20)},
  MRNUMBER = {2254187},
MRREVIEWER = {Paul\ E.\ Ehrlich},
       DOI = {10.1007/s11005-006-0091-5},
       URL = {https://doi-org.proxy.library.vanderbilt.edu/10.1007/s11005-006-0091-5},
}

@article {aBsM2003,
    AUTHOR = {Bernal, Antonio N. and S\'anchez, Miguel},
     TITLE = {On smooth {C}auchy hypersurfaces and {G}eroch's splitting
              theorem},
   JOURNAL = {Comm. Math. Phys.},
  FJOURNAL = {Communications in Mathematical Physics},
    VOLUME = {243},
      YEAR = {2003},
    NUMBER = {3},
     PAGES = {461--470},
      ISSN = {0010-3616,1432-0916},
   MRCLASS = {53C50 (53C80 83C05)},
  MRNUMBER = {2029362},
MRREVIEWER = {Paul\ E.\ Ehrlich},
       DOI = {10.1007/s00220-003-0982-6},
       URL = {https://doi-org.proxy.library.vanderbilt.edu/10.1007/s00220-003-0982-6},
}

@article{abbresciaBlueSierskiSpeck2025quasilinear,
  title={A quasilinear wave with a supersonic shock in a weak solution interrupting the classical development},
  author={Abbrescia, Leonardo and Blue, Pieter and Sbierski, Jan and Speck, Jared},
  journal={arXiv preprint arXiv:2511.07594},
  year={2025}
}

@article{shkoller2024geometry,
  title={The geometry of maximal development and shock formation for the Euler equations in multiple space dimensions},
  author={Shkoller, Steve and Vicol, Vlad},
  journal={Inventiones mathematicae},
  volume={237},
  number={3},
  pages={871--1252},
  year={2024},
  publisher={Springer}
}

@article{abbrescia2023relativistic,
  title={The relativistic {E}uler equations: {ESI} notes on their geo-analytic structures and implications for shocks in {$1D$} and multi-dimensions},
  author={Abbrescia, Leonardo and Speck, Jared},
  journal={Classical and Quantum Gravity},
  volume={40},
  number={24},
  pages={243001},
  year={2023},
  publisher={IOP Publishing}
}

@ARTICLE{lAjS2022,
       author = {{Abbrescia}, Leo and {Speck}, Jared},
        title = "{The emergence of the singular boundary from the crease in $3D$ compressible Euler flow}",
      journal = {arXiv e-prints},
         year = 2022,
        month = jul,
          eid = {arXiv:2207.07107},
        pages = {arXiv:2207.07107},
          doi = {10.48550/arXiv.2207.07107},
archivePrefix = {arXiv},
       eprint = {2207.07107},
 primaryClass = {math.AP},
       adsurl = {https://ui.adsabs.harvard.edu/abs/2022arXiv220707107A}
}

@article{abbrescia2025remarkable,
  title={Remarkable Localized Integral Identities for 3 D Compressible Euler Flow and the Double-Null Framework},
  author={Abbrescia, Leonardo and Speck, Jared},
  journal={Archive for Rational Mechanics and Analysis},
  volume={249},
  number={1},
  pages={10},
  year={2025},
  publisher={Springer}
}

@article{abbrescia2026emergencecauchyhorizoncrease,
      title={The emergence of the {C}auchy horizon from the crease for {$3D$} compressible {E}uler flow}, 
      author={Leonardo Abbrescia and Jared Speck},
      year={2026},
      month = sep,
      journal={arXiv preprint arXiv:2609.03101},
      url={https://arxiv.org/abs/2609.03101}, 
}

@preamble{
   "\def\cprime{$'$} "
}

@article {sA1999a,
    AUTHOR = {Alinhac, Serge},
     TITLE = {Blowup of small data solutions for a quasilinear wave equation
              in two space dimensions},
   JOURNAL = {Ann. of Math. (2)},
  FJOURNAL = {Annals of Mathematics. Second Series},
    VOLUME = {149},
      YEAR = {1999},
    NUMBER = {1},
     PAGES = {97--127},
      ISSN = {0003-486X},
     CODEN = {ANMAAH},
   MRCLASS = {35L70 (35B05)},
  MRNUMBER = {1680539 (2000d:35147)},
MRREVIEWER = {Nickolai A. Lar{\cprime}kin},
       DOI = {10.2307/121020},
       URL = {http://dx.doi.org/10.2307/121020},
}

@article {sA1999b,
    AUTHOR = {Alinhac, Serge},
     TITLE = {Blowup of small data solutions for a class of quasilinear wave
              equations in two space dimensions. {II}},
   JOURNAL = {Acta Math.},
  FJOURNAL = {Acta Mathematica},
    VOLUME = {182},
      YEAR = {1999},
    NUMBER = {1},
     PAGES = {1--23},
      ISSN = {0001-5962},
     CODEN = {ACMAA8},
   MRCLASS = {35L70 (35B05)},
  MRNUMBER = {1687180 (2000d:35148)},
MRREVIEWER = {Nickolai A. Lar{\cprime}kin},
       DOI = {10.1007/BF02392822},
       URL = {http://dx.doi.org/10.1007/BF02392822},
}

@article {gCfeKyHyS2025,
    AUTHOR = {Chen, Geng and El-Katri, Faris A. and Hu, Yanbo and Shen,
              Yannan},
     TITLE = {Global solution and singularity formation for the supersonic
              expanding wave of compressible {E}uler equations with radial
              symmetry},
   JOURNAL = {J. Lond. Math. Soc. (2)},
  FJOURNAL = {Journal of the London Mathematical Society. Second Series},
    VOLUME = {111},
      YEAR = {2025},
    NUMBER = {6},
     PAGES = {Paper No. e70208, 34},
      ISSN = {0024-6107,1469-7750},
   MRCLASS = {76N15 (35L65 35L67)},
  MRNUMBER = {4922416},
       DOI = {10.1112/jlms.70208},
       URL = {https://doi-org.proxy.library.vanderbilt.edu/10.1112/jlms.70208},
}

@ARTICLE{CB1952,
    author = "Yvonne Foures (Choquet)-Bruhat",
    title = "{Th\'{e}or\`{e}me d'existence pour certains syst\`{e}mes d'\'{e}quations aux d\'{e}riv\'{e}es partielles non 
    lin\'{e}aires}",
    journal = "Acta Mathematica",
    volume = 88,
   	year = 1952,
   	pages = "141-225",
}

@article {cBgR1969,
    AUTHOR = {Choquet-Bruhat, Yvonne and Geroch, Robert},
     TITLE = {Global aspects of the {C}auchy problem in general relativity},
   JOURNAL = {Comm. Math. Phys.},
  FJOURNAL = {Communications in Mathematical Physics},
    VOLUME = {14},
      YEAR = {1969},
     PAGES = {329--335},
      ISSN = {0010-3616},
   MRCLASS = {83.53},
  MRNUMBER = {MR0250640 (40 \#3872)},
MRREVIEWER = {H. A. Buchdahl},
}

@book{dC2007,
    AUTHOR = {Christodoulou, Demetrios},
     TITLE = {The Formation of Shocks in 3-Dimensional Fluids},
    SERIES = {EMS Monographs in Mathematics},
 PUBLISHER = {European Mathematical Society (EMS), Z\"urich},
      YEAR = {2007},
     PAGES = {viii+992},
      ISBN = {978-3-03719-031-9},
   MRCLASS = {76L05 (35L67 35Q35 35Q75 76Y05 83C55)},
  MRNUMBER = {2284927 (2008e:76104)},
MRREVIEWER = {Philippe G. LeFloch},
       DOI = {10.4171/031},
       URL = {http://dx.doi.org/10.4171/031},
}

@BOOK{dCsK1993,
    AUTHOR = {Christodoulou, Demetrios and Klainerman, Sergiu},
     TITLE = {The Global Nonlinear Stability of the {M}inkowski Space},
    SERIES = {Princeton Mathematical Series},
    VOLUME = {41},
 PUBLISHER = {Princeton University Press},
   ADDRESS = {Princeton, NJ},
      YEAR = {1993},
     PAGES = {x+514},
      ISBN = {0-691-08777-6},
   MRCLASS = {83C05 (35Q75 58G16 83C35)},
  MRNUMBER = {MR1316662 (95k:83006)},
MRREVIEWER = {Alan D. Rendall},
}

@article {dC1986a,
    AUTHOR = {Christodoulou, Demetrios},
     TITLE = {Global solutions of nonlinear hyperbolic equations for small
              initial data},
   JOURNAL = {Comm. Pure Appl. Math.},
  FJOURNAL = {Communications on Pure and Applied Mathematics},
    VOLUME = {39},
      YEAR = {1986},
    NUMBER = {2},
     PAGES = {267--282},
      ISSN = {0010-3640},
     CODEN = {CPAMA},
   MRCLASS = {35L70},
  MRNUMBER = {820070 (87c:35111)},
MRREVIEWER = {R. Glassey},
       DOI = {10.1002/cpa.3160390205},
       URL = {http://dx.doi.org/10.1002/cpa.3160390205},
}

@book {dCsM2014,
    AUTHOR = {Christodoulou, Demetrios and Miao, Shuang},
     TITLE = {Compressible flow and {E}uler's equations},
    SERIES = {Surveys of Modern Mathematics},
    VOLUME = {9},
 PUBLISHER = {International Press, Somerville, MA; Higher Education Press,
              Beijing},
      YEAR = {2014},
     PAGES = {x+iv+583},
      ISBN = {978-1-57146-297-8},
   MRCLASS = {76-02 (35Q31 35Q35 76Nxx)},
  MRNUMBER = {3288725},
}

@article {abbresciaspecknotices,
    AUTHOR = {Abbrescia, Leonardo and Speck, Jared},
     TITLE = {The geometry of maximal globally hyperbolic developments for
              {$3D$} compressible fluids terminating in shocks},
   JOURNAL = {Notices Amer. Math. Soc.},
  FJOURNAL = {Notices of the American Mathematical Society},
    VOLUME = {73},
      YEAR = {2026},
    NUMBER = {3},
     PAGES = {185--199},
      ISSN = {0002-9920,1088-9477},
   MRCLASS = {35Q35 (35L67 76L05)},
  MRNUMBER = {5040784},
       DOI = {10.1090/noti3313},
       URL = {https://doi-org.proxy.library.vanderbilt.edu/10.1090/noti3313},
}

@article {fEhRjS2019,
    AUTHOR = {Eperon, Felicity C. and Reall, Harvey S. and Sbierski, Jan J.},
     TITLE = {Predictability of subluminal and superluminal wave equations},
   JOURNAL = {Comm. Math. Phys.},
  FJOURNAL = {Communications in Mathematical Physics},
    VOLUME = {368},
      YEAR = {2019},
    NUMBER = {2},
     PAGES = {585--626},
      ISSN = {0010-3616},
   MRCLASS = {58J45 (83C05)},
  MRNUMBER = {3949719},
MRREVIEWER = {Vadim Valentinovich Kuzmichev},
       DOI = {10.1007/s00220-019-03428-1},
       URL = {https://doi-org.proxy.library.vanderbilt.edu/10.1007/s00220-019-03428-1},
}

@article {mGra1998,
    AUTHOR = {Grassin, Magali},
     TITLE = {Global smooth solutions to {E}uler equations for a perfect
              gas},
   JOURNAL = {Indiana Univ. Math. J.},
  FJOURNAL = {Indiana University Mathematics Journal},
    VOLUME = {47},
      YEAR = {1998},
    NUMBER = {4},
     PAGES = {1397--1432},
      ISSN = {0022-2518,1943-5258},
   MRCLASS = {76N10 (35Q30)},
  MRNUMBER = {1687130},
MRREVIEWER = {Beno\^it\ P.\ Desjardins},
       DOI = {10.1512/iumj.1998.47.1608},
       URL = {https://doi-org.proxy.library.vanderbilt.edu/10.1512/iumj.1998.47.1608},
}

@article {mHjJ2018b,
    AUTHOR = {{Had\v{z}i\'{c}}, Mahir and Jang, Juhi},
     TITLE = {Expanding large global solutions of the equations of
              compressible fluid mechanics},
   JOURNAL = {Invent. Math.},
  FJOURNAL = {Inventiones Mathematicae},
    VOLUME = {214},
      YEAR = {2018},
    NUMBER = {3},
     PAGES = {1205--1266},
      ISSN = {0020-9910,1432-1297},
   MRCLASS = {76N99 (35Q35 35R35)},
  MRNUMBER = {3878730},
MRREVIEWER = {Alessandro\ Morando},
       DOI = {10.1007/s00222-018-0821-1},
       URL = {https://doi-org.proxy.library.vanderbilt.edu/10.1007/s00222-018-0821-1},
}

@article {mHjJ2018a,
    AUTHOR = {{Had\v{z}i\'{c}}, Mahir and Jang, Juhi},
     TITLE = {Nonlinear stability of expanding star solutions of the
              radially symmetric mass-critical {E}uler-{P}oisson system},
   JOURNAL = {Comm. Pure Appl. Math.},
  FJOURNAL = {Communications on Pure and Applied Mathematics},
    VOLUME = {71},
      YEAR = {2018},
    NUMBER = {5},
     PAGES = {827--891},
      ISSN = {0010-3640,1097-0312},
   MRCLASS = {85A35 (35Q85 76U05 76W05)},
  MRNUMBER = {3794516},
       DOI = {10.1002/cpa.21721},
       URL = {https://doi-org.proxy.library.vanderbilt.edu/10.1002/cpa.21721},
}

@article{gHsKjSwW2016,
author = {Holzegel, Gustav and Klainerman, Sergiu and Speck, Jared and Wong, Willie Wai-Yeung},
title = {Small-data shock formation in solutions to 3D quasilinear wave equations: An overview},
journal = {Journal of Hyperbolic Differential Equations},
volume = {13},
number = {01},
pages = {1-105},
year = {2016},
doi = {10.1142/S0219891616500016},

URL = {http://www.worldscientific.com/doi/abs/10.1142/S0219891616500016},
eprint = {http://www.worldscientific.com/doi/pdf/10.1142/S0219891616500016}
}

@article{fJ1985,
    AUTHOR = {John, Fritz},
     TITLE = {Blow-up of radial solutions of {$u_{tt}=c^2(u_t)\Delta u$} in three space dimensions},
   JOURNAL = {Mat. Apl. Comput.},
  FJOURNAL = {Matem\'atica Aplicada e Computacional},
    VOLUME = {4},
      YEAR = {1985},
    NUMBER = {1},
     PAGES = {3--18},
      ISSN = {0101-8205},
   MRCLASS = {35L70 (35L67)},
  MRNUMBER = {808321 (87c:35114)},
MRREVIEWER = {R. Glassey},
}

@article{fJ1981,
    AUTHOR = {John, Fritz},
     TITLE = {Blow-up for quasilinear wave equations in three space
              dimensions},
   JOURNAL = {Comm. Pure Appl. Math.},
  FJOURNAL = {Communications on Pure and Applied Mathematics},
    VOLUME = {34},
      YEAR = {1981},
    NUMBER = {1},
     PAGES = {29--51},
      ISSN = {0010-3640},
     CODEN = {CPAMA},
   MRCLASS = {35L67},
  MRNUMBER = {600571 (83d:35096)},
MRREVIEWER = {Ronald DiPerna},
       DOI = {10.1002/cpa.3160340103},
       URL = {http://dx.doi.org/10.1002/cpa.3160340103},
}

@article {fJsK1984,
    AUTHOR = {John, Fritz and Klainerman, Sergiu},
     TITLE = {Almost global existence to nonlinear wave equations in three
              space dimensions},
   JOURNAL = {Comm. Pure Appl. Math.},
  FJOURNAL = {Communications on Pure and Applied Mathematics},
    VOLUME = {37},
      YEAR = {1984},
    NUMBER = {4},
     PAGES = {443--455},
      ISSN = {0010-3640},
     CODEN = {CPAMA},
   MRCLASS = {35L70},
  MRNUMBER = {745325 (85k:35147)},
MRREVIEWER = {R. Glassey},
       DOI = {10.1002/cpa.3160370403},
       URL = {http://dx.doi.org/10.1002/cpa.3160370403},
}

@inproceedings {sK1984,
    AUTHOR = {Klainerman, Sergiu},
     TITLE = {Long time behaviour of solutions to nonlinear wave equations},
 BOOKTITLE = {Proceedings of the {I}nternational {C}ongress of
              {M}athematicians, {V}ol.\ 1, 2 ({W}arsaw, 1983)},
     PAGES = {1209--1215},
 PUBLISHER = {PWN, Warsaw},
      YEAR = {1984},
   MRCLASS = {35B40 (35G20)},
  MRNUMBER = {804771},
}

@article{klainerman2026inevitable,
  title={Inevitable shock formation for 3-{D} compressible {E}uler flows},
  author={Klainerman, Sergiu and Wang, Qian and Yang, Shiwu and Yu, Pin},
  journal={arXiv preprint arXiv:2608.09843},
  year={2026}
}

@article{hLiR2003,
    AUTHOR = {Lindblad, Hans and Rodnianski, Igor},
     TITLE = {The weak null condition for {E}instein's equations},
   JOURNAL = {C. R. Math. Acad. Sci. Paris},
  FJOURNAL = {Comptes Rendus Math\'ematique. Acad\'emie des Sciences. Paris},
    VOLUME = {336},
      YEAR = {2003},
    NUMBER = {11},
     PAGES = {901--906},
      ISSN = {1631-073X},
   MRCLASS = {83C05 (35Q75 58J45)},
  MRNUMBER = {1994592 (2004h:83008)},
MRREVIEWER = {Norbert Noutchegueme},
       DOI = {10.1016/S1631-073X(03)00231-0},
       URL = {http://dx.doi.org/10.1016/S1631-073X(03)00231-0},
}

@ARTICLE{hLiR2010,
    AUTHOR = {Lindblad, Hans and Rodnianski, Igor},
     TITLE = {The global stability of {Minkowski} space-time in harmonic gauge},
   JOURNAL = {Annals of Mathematics},
  	VOLUME = {171},
      YEAR = {2010},
    NUMBER = {3},
     PAGES = {1401--1477},
}

@article {jLjS2020a,
    AUTHOR = {Luk, Jonathan and Speck, Jared},
     TITLE = {The hidden null structure of the compressible {E}uler
              equations and a prelude to applications},
   JOURNAL = {J. Hyperbolic Differ. Equ.},
  FJOURNAL = {Journal of Hyperbolic Differential Equations},
    VOLUME = {17},
      YEAR = {2020},
    NUMBER = {1},
     PAGES = {1--60},
      ISSN = {0219-8916},
   MRCLASS = {35Q35 (35L67 35L72 35Q31 76N10)},
  MRNUMBER = {4109292},
       DOI = {10.1142/S0219891620500010},
       URL = {https://doi-org.proxy.library.vanderbilt.edu/10.1142/S0219891620500010},
}

@article {jLjS2018,
    AUTHOR = {Luk, Jonathan and Speck, Jared},
     TITLE = {Shock formation in solutions to the 2{D} compressible {E}uler
              equations in the presence of non-zero vorticity},
   JOURNAL = {Invent. Math.},
  FJOURNAL = {Inventiones Mathematicae},
    VOLUME = {214},
      YEAR = {2018},
    NUMBER = {1},
     PAGES = {1--169},
      ISSN = {0020-9910},
   MRCLASS = {35L67 (35L05 35Q31 76N10)},
  MRNUMBER = {3858399},
       DOI = {10.1007/s00222-018-0799-8},
       URL = {https://doi-org.libproxy.mit.edu/10.1007/s00222-018-0799-8},
}

@article {jLjS2024,
    AUTHOR = {Luk, Jonathan and Speck, Jared},
     TITLE = {The stability of simple plane-symmetric shock formation for
              three-dimensional compressible {E}uler flow with vorticity and
              entropy},
   JOURNAL = {Anal. PDE},
  FJOURNAL = {Analysis \& PDE},
    VOLUME = {17},
      YEAR = {2024},
    NUMBER = {3},
     PAGES = {831--941},
      ISSN = {2157-5045,1948-206X},
   MRCLASS = {35L67 (35L05 35Q31 76L05 76N10)},
  MRNUMBER = {4736521},
MRREVIEWER = {Marta\ Lewicka},
       DOI = {10.2140/apde.2024.17.831},
       URL = {https://doi-org.proxy.library.vanderbilt.edu/10.2140/apde.2024.17.831},
}

@book{bO1983,
    AUTHOR = {O'Neill, Barrett},
     TITLE = {Semi-{R}iemannian geometry},
    SERIES = {Pure and Applied Mathematics},
    VOLUME = {103},
      NOTE = {With applications to relativity},
 PUBLISHER = {Academic Press Inc. [Harcourt Brace Jovanovich Publishers]},
   ADDRESS = {New York},
      YEAR = {1983},
     PAGES = {xiii+468},
      ISBN = {0-12-526740-1},
   MRCLASS = {53-01 (53B30 53C50 83-02)},
  MRNUMBER = {719023 (85f:53002)},
MRREVIEWER = {N. V. Mitskevich},
}

@article {jSb2016,
    AUTHOR = {Sbierski, Jan},
     TITLE = {On the existence of a maximal {C}auchy development for the
              {E}instein equations: a dezornification},
   JOURNAL = {Ann. Henri Poincar\'e},
  FJOURNAL = {Annales Henri Poincar\'e. A Journal of Theoretical and
              Mathematical Physics},
    VOLUME = {17},
      YEAR = {2016},
    NUMBER = {2},
     PAGES = {301--329},
      ISSN = {1424-0637},
   MRCLASS = {83C05},
  MRNUMBER = {3447847},
MRREVIEWER = {Willie W. Wong},
       DOI = {10.1007/s00023-015-0401-5},
       URL = {http://dx.doi.org/10.1007/s00023-015-0401-5},
}

@ARTICLE{tS1985,
    AUTHOR = {Sideris, Thomas C.},
     TITLE = {Formation of singularities in three-dimensional compressible
              fluids},
   JOURNAL = {Comm. Math. Phys.},
  FJOURNAL = {Communications in Mathematical Physics},
    VOLUME = {101},
      YEAR = {1985},
    NUMBER = {4},
     PAGES = {475--485},
      ISSN = {0010-3616},
     CODEN = {CMPHAY},
   MRCLASS = {35Q20 (76N10)},
  MRNUMBER = {MR815196 (87d:35127)},
MRREVIEWER = {Charles J. Amick},
}

@article {jS2019c,
    AUTHOR = {Speck, Jared},
     TITLE = {A {N}ew {F}ormulation of the 3{D} {C}ompressible {E}uler
              {E}quations with {D}ynamic {E}ntropy: {R}emarkable {N}ull
              {S}tructures and {R}egularity {P}roperties},
   JOURNAL = {Arch. Ration. Mech. Anal.},
  FJOURNAL = {Archive for Rational Mechanics and Analysis},
    VOLUME = {234},
      YEAR = {2019},
    NUMBER = {3},
     PAGES = {1223--1279},
      ISSN = {0003-9527},
   MRCLASS = {76N10 (35)},
  MRNUMBER = {4011696},
       DOI = {10.1007/s00205-019-01411-7},
       URL = {https://doi-org.libproxy.mit.edu/10.1007/s00205-019-01411-7},
}

@book{jS2016b,
  title={Shock formation in small-data solutions to {$3D$} quasilinear wave equations},
  author={Speck, Jared},
  isbn={9781470428570},
  series={Mathematical Surveys and Monographs},
  url={https://books.google.com/books?id=0gK5DQAAQBAJ},
  year={2016}
}

@misc{qW2025,
      title={On global dynamics of $3$-D irrotational compressible fluids}, 
      author={Qian Wang},
      year={2025},
      eprint={2407.13649},
      archivePrefix={arXiv},
      primaryClass={math.AP},
      url={https://arxiv.org/abs/2407.13649}, 
}

@article {dY2025,
    AUTHOR = {Yu, Dongxiao},
     TITLE = {Nontrivial global solutions to some quasilinear wave equations
              in three space dimensions},
   JOURNAL = {Ann. PDE},
  FJOURNAL = {Annals of PDE. Journal Dedicated to the Analysis of Problems
              from Physical Sciences},
    VOLUME = {11},
      YEAR = {2025},
    NUMBER = {2},
     PAGES = {Paper No. 25, 159},
      ISSN = {2524-5317,2199-2576},
   MRCLASS = {35L72 (35A01 35L05 35Q31)},
  MRNUMBER = {4949961},
MRREVIEWER = {Dongbing\ Zha},
       DOI = {10.1007/s40818-025-00205-3},
       URL = {https://doi-org.proxy.library.vanderbilt.edu/10.1007/s40818-025-00205-3},
}

@article{tBsSvV2019a,
  author = {Buckmaster, T. and Shkoller, S. and Vicol, V.},
  title = {Formation of {P}oint {S}hocks for 3{D} {C}ompressible {E}uler},
  journal = {Communications on Pure and Applied Mathematics},
  volume = {76},
  number = {9},
  pages = {2073-2191},
  doi = {https://doi.org/10.1002/cpa.22068},
  url = {https://onlinelibrary.wiley.com/doi/abs/10.1002/cpa.22068},
  year = {2023}
}

@article {fMpRiRjS2022b,
    AUTHOR = {Merle, Frank and Rapha\"{e}l, Pierre and Rodnianski, Igor and
              Szeftel, Jeremie},
     TITLE = {On the implosion of a compressible fluid {II}: {S}ingularity
              formation},
   JOURNAL = {Ann. of Math. (2)},
  FJOURNAL = {Annals of Mathematics. Second Series},
    VOLUME = {196},
      YEAR = {2022},
    NUMBER = {2},
     PAGES = {779--889},
      ISSN = {0003-486X},
   MRCLASS = {35Q35 (35B44)},
  MRNUMBER = {4445443},
MRREVIEWER = {Robert Schippa},
       DOI = {10.4007/annals.2022.196.2.4},
       URL = {https://doi-org.proxy.library.vanderbilt.edu/10.4007/annals.2022.196.2.4},
}

@article {fMpRiRjS2022a,
    AUTHOR = {Merle, Frank and Rapha\"{e}l, Pierre and Rodnianski, Igor and
              Szeftel, Jeremie},
     TITLE = {On the implosion of a compressible fluid {I}: {S}mooth
              self-similar inviscid profiles},
   JOURNAL = {Ann. of Math. (2)},
  FJOURNAL = {Annals of Mathematics. Second Series},
    VOLUME = {196},
      YEAR = {2022},
    NUMBER = {2},
     PAGES = {567--778},
      ISSN = {0003-486X},
   MRCLASS = {35Q35 (34C37)},
  MRNUMBER = {4445442},
MRREVIEWER = {Wei Lian},
       DOI = {10.4007/annals.2022.196.2.3},
       URL = {https://doi-org.proxy.library.vanderbilt.edu/10.4007/annals.2022.196.2.3},
}

\end{document}